%% file: memoire.tex
\documentclass[11pt, francais]{smfbook}
\usepackage{amsfonts, amssymb, amsmath, amsthm, latexsym, enumerate, epsfig, color, mathrsfs, fancyhdr,supertabular}
\usepackage[all]{xy}

\usepackage[frenchb]{babel}

\usepackage{mysections}

\usepackage[mac]{inputenc}

\makeindex

\def\N{\ensuremath{\mathbf{N}}}
\def\Z{\ensuremath{\mathbf{Z}}}
\def\Q{\ensuremath{\mathbf{Q}}}
\def\R{\ensuremath{\mathbf{R}}}
\def\C{\ensuremath{\mathbf{C}}}
\def\D{\ensuremath{\mathbf{D}}}
\def\A{\ensuremath{\mathbf{A}}}
\def\P{\ensuremath{\mathbf{P}}}
\def\F{\ensuremath{\mathbf{F}}}

\def\E#1#2{\ensuremath{\mathbf{A}^{#1,\mathrm{an}}_{#2}}}
\def\AZ{\ensuremath{\mathbf{A}^{1,\mathrm{an}}_{\mathbf{Z}}}}
\def\AA{\ensuremath{\mathbf{A}^{1,\mathrm{an}}_{A}}}

\def\Frac{\textrm{Frac}}
\def\Spec{\textrm{Spec}}
\def\Ker{\textrm{Ker}}

\def\p{\ensuremath{\mathfrak{p}}}
\def\m{\ensuremath{\mathfrak{m}}}

\def\As{\ensuremath{\mathscr{A}}}
\def\Bs{\ensuremath{\mathscr{B}}}
\def\Cs{\ensuremath{\mathscr{C}}}
\def\Ds{\ensuremath{\mathscr{D}}}
\def\Fs{\ensuremath{\mathscr{F}}}
\def\Gs{\ensuremath{\mathscr{G}}}
\def\Hs{\ensuremath{\mathscr{H}}}
\def\Is{\ensuremath{\mathscr{I}}}

\def\Ks{\ensuremath{\mathscr{K}}}

\def\Ms{\ensuremath{\mathscr{M}}}
\def\Ns{\ensuremath{\mathscr{N}}}
\def\Os{\ensuremath{\mathscr{O}}}
\def\Ps{\ensuremath{\mathscr{P}}}
\def\Qs{\ensuremath{\mathscr{Q}}}
\def\Rs{\ensuremath{\mathscr{R}}}
\def\Ss{\ensuremath{\mathscr{S}}}
\def\Ts{\ensuremath{\mathscr{T}}}
\def\Us{\ensuremath{\mathscr{U}}}
\def\Vs{\ensuremath{\mathscr{V}}}
\def\Ws{\ensuremath{\mathscr{W}}}

\def\Pc{\ensuremath{\mathcal{P}}}

\def\cn#1#2{\ensuremath{[\![ {#1},{#2}]\!]}}
\def\of#1#2#3{\ensuremath{\mathopen{#1}{#2}\mathclose{#3}}}
\def\surj{\ensuremath{\to \hspace{-0.35cm}\to}}
\def\la{\ensuremath{\langle}}
\def\ra{\ensuremath{\rangle}}

\def\eps{\ensuremath{\varepsilon}}

\def\bk{\ensuremath{{\boldsymbol{k}}}}
\def\br{\ensuremath{{\boldsymbol{r}}}}
\def\bs{\ensuremath{{\boldsymbol{s}}}}
\def\bt{\ensuremath{{\boldsymbol{t}}}}
\def\bu{\ensuremath{{\boldsymbol{u}}}}
\def\bv{\ensuremath{{\boldsymbol{v}}}}
\def\bT{\ensuremath{{\boldsymbol{T}}}}

\def\alphab{\ensuremath{{\boldsymbol{\alpha}}}}
\def\betab{\ensuremath{{\boldsymbol{\beta}}}}

\def\bmax{\ensuremath{\textrm{\bf{max}}}}
\def\bmin{\ensuremath{\textrm{\bf{min}}}}
\def\b0{\ensuremath{{\boldsymbol{0}}}}

\def\disp{\displaystyle}

\def\Zr{\ensuremath{\Z_{r^+}[\![T]\!]}}
\def\Zun{\ensuremath{\Z_{1^-}[\![T]\!]}}

\begin{document}
\setlength{\baselineskip}{0.55cm}


\title{La droite de Berkovich\\ \qquad \qquad \quad \ sur Z}
\author{J\'er\^ome Poineau}
\address{Institut de recherche math\'ematique avanc\'ee\\ 7, rue Ren\'e Descartes, 67084 Strasbourg, France}
\email{jerome.poineau@math.u-strasbg.fr}

\date{\today}


\begin{abstract}
Ce texte est consacr\'e \`a l'\'etude de la droite de Berkovich au-dessus d'un anneau d'entiers de corps de nombres. Cet objet contient naturellement des copies de la droite analytique complexe (ou de son quotient par la conjugaison), associ\'ees aux places infinies, et des droites de Berkovich usuelles au-dessus de corps ultram\'etriques complets, associ\'ees au places finies. Nous montrons qu'il jouit de bonnes propri\'et\'es, topologiques aussi bien qu'alg\'ebriques. Nous exhibons \'egalement quelques espaces de Stein naturels contenus dans cette droite. 

Nous proposons des applications de cette th\'eorie \`a l'\'etude des s\'eries arith\-m\'e\-tiques convergentes : prescription de z\'eros et de p\^oles, noeth\'erianit\'e d'anneaux globaux et probl\`eme inverse de Galois. Des exemples typiques de telles s\'eries sont fournis par les fonctions holomorphes sur le disque unit\'e ouvert complexe dont le d\'eveloppement en~$0$ est \`a coefficients entiers.
\end{abstract}

\maketitle










\pagestyle{empty}

\tableofcontents

\pagestyle{headings}

\pagenumbering{roman}

\include{introduction}

\pagenumbering{arabic}

\include{Banach}

\include{series}

\include{espace}

\include{droite}

\include{fini}

\include{Stein}

\include{applications}

\nocite{}
\bibliographystyle{plain}
\bibliography{biblio}

\input{notations}


\printindex

\end{document}

%% file: introduction.tex

\chapter*{Introduction}


\`A la fin des ann\'ees quatre-vingts, Vladimir G. Berkovich a propos\'e une nouvelle approche de la g\'eom\'etrie analytique $p$-adique. Ses id\'ees, d\'evelopp\'ees dans l'ouvrage \cite{rouge} et approfondies dans l'article \cite{bleu} se sont r\'ev\'el\'ees tr\`es fructueuses ; elles ont permis de d\'emontrer plusieurs conjectures de g\'eom\'etrie arithm\'etique et trouvent maintenant des applications dans des domaines vari\'es : syst\`emes dynamiques, th\'eorie d'Arakelov, dessins d'enfants $p$-adiques, variation de structure de Hodge, \emph{etc.} Pour une introduction au sujet et une pr\'esentation des diff\'erentes applications, nous renvoyons le lecteur int\'eress\'e aux textes de vulgarisation \cite{AntoineBourbaki} et \cite{AntoineGazette}.

Bien que la th\'eorie n'ait \'et\'e v\'eritablement d\'evelopp\'ee que sur les corps ultram\'etriques complets, V. Berkovich propose, dans \cite{rouge}, une d\'efinition d'espace analytique au-dessus de n'importe quel anneau de Banach. Elle s'applique donc lorsque l'on consid\`ere comme anneau de base l'anneau $\Z$ des nombres entiers, muni de la valeur absolue usuelle $|.|_{\infty}$.
Nous nous proposons ici d'entreprendre l'\'etude des espaces analytiques dans ce cas particulier.    

Diff\'erentes valeurs absolues joueront un r\^ole dans notre \'etude. Si $p$ d\'esigne un nombre premier, nous d\'efinissons la valeur absolue $p$-adique $|.|_{p}$ sur $\Z$ de la fa\c{c}on suivante : nous posons $|0|_{p} = 0$ et, pour tout nombre entier $n = p^r\, m\in\Z^*$, o\`u~$m$ est premier \`a~$p$,
$$|n|_{p} = |p^r\, m|_{p} = p^{-r}.$$
Elle se prolonge de fa\c{c}on unique \`a $\Q$. Notons~$\Q_{p}$ le compl\'et\'e de~$\Q$ pour cette valeur absolue et choisissons-en une cl\^oture alg\'ebrique~$\overline{\Q}_{p}$. La valeur absolue~$|.|_{p}$ se prolonge encore de fa\c{c}on unique en une valeur absolue sur~$\overline{\Q}_{p}$. Nous noterons~$\C_{p}$ son compl\'et\'e. Ce corps, qui est alg\'ebriquement clos et complet, est parfois appel\'e corps des {nombres complexes $p$-adiques}. Nous noterons~$|.|_{p}$ l'unique valeur absolue sur~$\C_{p}$ qui prolonge la valeur absolue $p$-adique sur~$\Q$. 

Pour $f\in\Q[\![T]\!]$, notons $R_{\infty}(f)$ le rayon de convergence de la s\'erie $f$ vue comme s\'erie de $\C[\![T]\!]$ et,
pour tout nombre premier $p$, notons $R_{p}(f)$ le rayon de convergence de la s\'erie $f$ vue comme s\'erie de $\C_{p}[\![T]\!]$.
Appelons s\'erie arithm\'etique toute s\'erie de la forme
$$f\in \Z\left[\frac{1}{p_{1}\cdots p_{t}}\right][\![T]\!]$$ 
v\'erifiant des conditions du type
$$R_{\infty}(f) > r_{\infty} \textrm{ et } \forall i\in\cn{1}{t}, \, R_{p_{i}}(f) > r_{i},$$
o\`u $t$ est un nombre entier, $p_{1},\ldots,p_{t}$ des nombres premiers et $r_{1},\ldots,r_{t},r_{\infty}$ des nombres r\'eels strictement positifs. De telles fonctions apparaissent naturellement lorsque l'on \'etudie les anneaux locaux de la droite analytique sur $\Z$ ou certains anneaux de sections globales. L'\'etude g\'eom\'etrique que nous allons mener nous permettra d'obtenir des informations sur certains anneaux de s\'eries de ce type.


\begin{figure}[htb]
\begin{center}
\input{MZ.pstex_t}
\caption{L'espace topologique $\Ms(\Z)$.}
\end{center}
\end{figure}

\bigskip

\pagebreak

{\bf Description des espaces en jeu}

Par d\'efinition, l'espace $\Ms(\Z)$ est l'ensemble des semi-normes multiplicatives sur $\Z$, c'est-\`a-dire des applications de $\Z$ dans $\R_{+}$ qui sont sous-additives, multiplicatives, envoient $0$ sur $0$ et $1$ sur $1$. 
Topologiquement, il est constitu\'e d'une branche, hom\'eomorphe \`a un segment, pour chaque nombre premier~$p$ et d'une branche suppl\'ementaire, associ\'ee \`a la valeur absolue archim\'edienne usuelle. Ces branches se rejoignent en un point, que nous appellerons central, associ\'e \`a la valeur absolue triviale~$|.|_{0}$  (\emph{cf.} figure 
\thefigure
). Signalons que la topologie au voisinage du point central est strictement plus grossi\`ere que la topologie d'arbre.

Soit $n\in\N$. L'espace affine analytique de dimension $n$ au-dessus de $\Z$, que nous noterons $\E{n}{\Z}$, est l'ensemble des semi-normes multiplicatives sur l'anneau de polyn\^omes $\Z[T_{1},\ldots,T_{n}]$. Il est muni d'une projection continue vers la base~$\Ms(\Z)$. Au-dessus des points de la branche archim\'edienne, les fibres de cette projection sont isomorphes \`a l'espace $\C^n$ quotient\'e par l'action de la conjugaison complexe et, au-dessus des points de la branche $p$-adique, ce sont des espaces de Berkovich $p$-adiques de dimension $n$. Il appara\^{\i}t donc clairement que, pour \'etudier cet espace, il nous faudra mettre en {\oe}uvre des techniques pouvant s'appliquer tant dans un cadre archim\'edien qu'ultram\'etrique.


\bigskip

{\bf G\'eom\'etrie analytique complexe}

Dans le cas archim\'edien, la g\'eom\'etrie analytique complexe met \`a notre disposition de nombreux outils.
Les fondations de cette th\'eorie reposent sur une \'etude locale des vari\'et\'es et des fonctions. La com\-pr\'e\-hen\-sion des anneaux locaux des espaces affines y joue donc un r\^ole pr\'epond\'erant. Fixons~$n\in\N$. L'anneau local $\Os_{0}$ de l'espace affine $\C^n$ en $0$ est constitu\'e des s\'eries de la forme
$$\sum_{(k_{1},\ldots,k_{n})\in\N^n} a_{k_{1},\ldots,k_{n}}\, T_{1}^{k_{1}} \cdots T_{n}^{k_{n}}$$
dont le rayon de convergence est strictement positif. Le th\'eor\`eme de division de Weierstra{\ss} nous permet, sous certaines conditions, de diviser une s\'erie de la forme pr\'ec\'edente par une autre et d'obtenir un reste polynomial en la derni\`ere variable. Une fois ce r\'esultat connu, on d\'emontre sans peine que l'anneau $\Os_{0}$ est un anneau local noeth\'erien, r\'egulier et de dimension $n$. Signalons que la d\'emonstration classique du th\'eor\`eme de division de Weierstra{\ss} repose sur le th\'eor\`eme de Rouch\'e et la formule de Cauchy.

\bigskip

\pagebreak

{\bf G\'eom\'etrie analytique $p$-adique}

Bien que le corps des nombres complexes $p$-adiques $\C_{p}$ pr\'esente des analogies avec le corps des nombres complexes $\C$, il en diff\`ere par la topologie. Indiquons, par exemple, que le corps $\C_{p}$ est totalement discontinu (ses composantes connexes sont r\'eduites \`a des points) et n'est pas localement compact. Dans cette situation, il n'est gu\`ere ais\'e de mettre en place une g\'eom\'etrie analytique jouissant de propri\'et\'es raisonnables : il existe bien trop de fonctions localement analytiques. On v\'erifie, par exemple, que la fonction qui vaut $0$ sur le disque ouvert de centre $0$ et de rayon 1 de $\C_{p}$ et $1$ sur son compl\'ementaire est localement d\'eveloppable en s\'erie enti\`ere ! 

Au d\'ebut des ann\'ees soixante, John Tate a apport\'e une solution \`a ce probl\`eme (\emph{cf.} \cite{Tate}). Les espaces qu'il construit, appel\'es {espaces analytiques rigides}, ne sont pas des espaces topologiques, mais des {\it sites} : on distingue certains ouverts et on n'autorise que certains recouvrements. Par exemple, le recouvrement de $\C_{p}$ que nous avons d\'ecrit pr\'ec\'edemment est interdit. Ce formalisme permet de mettre en place, dans le cas $p$-adique, une g\'eom\'etrie analytique fort semblable \`a celle que nous connaissons dans le cas complexe.

Entrons un peu dans les d\'etails. Les objets de base \`a partir desquels est construite la g\'eom\'etrie analytique rigide sont les alg\`ebres que l'on appelle, aujourd'hui, alg\`ebres de Tate. Contrairement \`a ceux de la th\'eorie complexe, ce ne sont pas des anneaux locaux, mais globaux. Soit $n\in\N$. L'alg\`ebre de Tate $\C_p\{T_1,\ldots,T_n\}$ est constitu\'ee des \'el\'ements de la forme
$$\sum_{(k_1,\ldots,k_n)\in\N^n}  a_{k_1,\ldots,k_n}\, T_1^{k_1}\ldots T_n^{k_n} \in \C_p[\![T_1,\ldots,T_n]\!]$$ 
v\'erifiant la condition
$$\lim_{(k_1,\ldots,k_n) \to +\infty} |a_{k_1,\ldots,k_n}|_p = 0.$$
Cet anneau est pr\'ecis\'ement l'anneau des s\'eries convergentes sur le disque ferm\'e de centre~$0$ et de polyrayon $(1,\ldots,1)$ de~$\C_p^n$. C'est le caract\`ere ultram\'etrique de la valeur absolue $p$-adique qui nous permet de donner un sens \`a cette notion de convergence sur un disque ferm\'e. 

Dans cette th\'eorie, il existe \'egalement un th\'eor\`eme de division de Weierstra{\ss} qui rend les m\^emes services que dans le cadre complexe. En l'utilisant, on d\'emontre ais\'ement que l'alg\`ebre de Tate $\C_p\{T_1,\ldots,T_n\}$ est un anneau noeth\'erien et r\'egulier de dimension $n$. Signalons que, cette fois-ci, la d\'e\-mons\-tra\-tion du th\'eo\-r\`e\-me de division de Weierstra{\ss} repose sur des arguments de r\'eduction modulo $p$.


\begin{figure}[htb]
\begin{center}
\input{A1Cp2.pstex_t}
\caption{La droite projective $\P^{1,\mathrm{an}}_{\C_{p}}$. 
}
\end{center}
\end{figure}

\bigskip

{\bf L'approche de Vladimir G. Berkovich}

Les descriptions pr\'ec\'edentes laissent entrevoir les difficult\'es qui se pr\'esentent lorsque l'on cherche \`a r\'eunir les espaces analytiques archim\'ediens et ultra\-m\'e\-tri\-ques dans un formalisme commun. 
L'approche que propose V. Berkovich des espaces analytiques $p$-adiques va permettre d'apporter une solution \`a ce probl\`eme.

Choisissant un point de vue diff\'erent de celui de J. Tate, V. Berkovich ajoute de tr\`es nombreux points aux espaces. \`A titre d'exemple, la droite affine analytique $\E{1}{\C_{p}}$ sur $\C_{p}$ qu'il d\'efinit poss\`ede une structure d'arbre et les points de~$\C_{p}$ sont confin\'es aux extr\'emit\'es de certaines branches. Nous avons esquiss\'e une repr\'esentation de la droite projective analytique $\P^{1,\mathrm{an}}_{\C_{p}}$ \`a la figure \thefigure. On obtient la droite affine $\E{1}{\C_{p}}$ en enlevant le point not\'e $\infty$.

Le proc\'ed\'e de construction qu'utilise V. Berkovich rend la description explicite de ses espaces d\'elicate, mais ils b\'en\'eficient d'autres avantages. Ce sont de v\'eritables espaces topologiques, localement compacts et localement connexes par arcs. Ces propri\'et\'es ouvrent la voie \`a une d\'efinition locale du faisceau structural. Dans le cas de l'espace affine, V. Berkovich propose de le d\'efinir comme le faisceau des fonctions qui sont localement limites uniformes de fractions rationnelles sans p\^oles. Indiquons que l'on retrouve bien ainsi le faisceau construit \`a partir des alg\`ebres de Tate. C'est d'ailleurs v\'eritablement sur la th\'eorie des espaces analytiques rigides que V. Berkovich b\^atit la sienne et il n'utilise gu\`ere la d\'efinition locale du faisceau.

Les d\'efinitions propos\'ees par V.~Berkovich valent \'egalement dans le cas des corps archim\'ediens. Signalons que l'espace de Berkovich affine de dimension~$n$ sur~$\C$ co\"{\i}ncide avec~$\C^n$ et que le faisceau dont il est muni est bien celui des fonctions analytiques.



\bigskip

\bigskip

\pagebreak

Nous venons d'expliquer que les espaces analytiques de Berkovich permettent d'envisager une \'etude locale des espaces analytiques sur $\Z$. Le pr\'esent travail constitue un premier pas dans cette direction. Soulignons que, bien que les id\'ees et d\'efinitions introduites par V. Berkovich invitent \`a adopter ce point de vue, une telle \'etude n'a, \`a notre connaissance, jamais encore \'et\'e entreprise. Indiquons, \`a pr\'esent, le plan que nous allons adopter.

\bigskip 

{\bf Espaces analytiques sur un anneau de Banach}

Le premier chapitre de ce m\'emoire est consacr\'e aux espaces analytiques sur un anneau de Banach quelconque. Nous y rappelons la d\'efinition d'espace analytique au sens de V.~Berkovich ainsi que la construction du faisceau structural qu'il propose. Nous donnons quelques exemples et d\'ecrivons explicitement la droite analytique au-dessus de tout corps valu\'e.

\bigskip

{\bf Alg\`ebres de s\'eries convergentes}

Nous consacrons le deuxi\`eme chapitre \`a l'\'etude d'anneaux de s\'eries convergentes \`a coefficients dans un anneau de Banach. En prenant des limites inductives de tels anneaux, nous obtenons un anneau local sur lequel nous parvenons \`a d\'emontrer un th\'eor\`eme de division de Weierstra{\ss}. Bien entendu, notre preuve ne peut faire appel ni \`a la formule de Cauchy, ni \`a la r\'eduction modulo $p$, faute d'analogue de la premi\`ere dans le cas ultram\'etrique et de la seconde dans le cas archim\'edien. Nous utilisons donc une m\'ethode, inspir\'ee des travaux de H.~Grauert et R.~Remmert, faisant simplement appel \`a des techniques d'alg\`ebres de Banach. \`A l'aide de ce th\'eor\`eme, nous obtenons des r\'esultats de noeth\'erianit\'e et de r\'egularit\'e pour les anneaux locaux consid\'er\'es. 

Afin de pouvoir utiliser ces r\'esultats, nous entreprenons ensuite une \'etude topologique locale aboutissant \`a la d\'emonstration du fait que les anneaux locaux en certains points des espaces de Berkovich sont isomorphes \`a de tels anneaux de s\'eries convergentes.

Nous terminons ce chapitre par la d\'emonstration que les anneaux locaux des espaces de Berkovich sont hens\'eliens. Ce r\'esultat g\'en\'eralise le r\'esultat classique valable pour les espaces au-dessus d'un corps valu\'e et complet, archim\'edien ou non.



\bigskip

{\bf Espace affine analytique au-dessus d'un anneau d'entiers de corps de nombres}

Dans le troisi\`eme chapitre, nous consid\'erons un anneau d'entiers de corps de nombres~$A$ et restreignons notre propos aux espaces analytiques dont la base est le spectre analytique~$\Ms(A)$ de cet anneau. Nous commen\c{c}ons par d\'ecrire cette base elle-m\^eme, aussi bien l'espace topologique sous-jacent, \`a l'aide du th\'eor\`eme d'Ostrowski, que les sections du faisceau structural.

Dans un second temps, nous nous int\'eressons aux espaces affines de dimension quelconque au-dessus de l'anneau~$A$. Nous d\'emontrons quelques r\'esultats concernant la topologie de ces espaces et \'etudions les anneaux locaux en certains points. Nous faisons ici appel aux r\'esultats sur les anneaux de s\'eries convergentes d\'emontr\'es au deuxi\`eme chapitre ainsi qu'\`a la propri\'et\'e d'hens\'elianit\'e, qui nous permet d'\'etablir l'existence d'isomorphismes locaux. Malheureusement, cette \'etude n'est pas compl\`ete et il est vraisemblable qu'il faille introduire de nouvelles techniques afin de la mener \`a terme.

Signalons que nous parvenons \'egalement \`a d\'ecrire explicitement certains anneaux locaux et les anneaux de sections globales sur les disques et les couronnes en termes de s\'eries convergentes. En utilisant le fait que les anneaux locaux sont hens\'eliens, nous obtenons une nouvelle d\'emonstration du th\'eor\`eme classique d'Eisenstein.

\begin{thmi}[Eisenstein]
Soit~$K$ un corps de nombres. Notons~$A$ l'anneau de ses entiers. Soit~$f$ un \'el\'ement de~$K[\![T]\!]$ qui est entier sur~$K[T]$. Alors
\begin{enumerate}[\it i)]
\item il existe un \'el\'ement~$a$ de~$A^*$ tel que la s\'erie $f(aT)$ soit \`a coefficients dans~$A$ ;
\item le rayon de convergence de la s\'erie~$f$ est strictement positif en toute place.
\end{enumerate}
\end{thmi}


\bigskip


{\bf Droite affine analytique au-dessus d'un anneau d'entiers de corps de nombres}

Le quatri\`eme chapitre est consacr\'e sp\'ecifiquement \`a la droite affine analytique au-dessus d'un anneau d'entiers de corps de nombres~$A$. Dans ce cadre, nous parvenons \`a compl\'eter les r\'esultats du chapitre pr\'ec\'edent et \`a \'etudier tous les points. Nous obtenons les r\'esultats suivants, conformes \`a l'intuition.


\begin{thmi}
\begin{enumerate}[\it i)]
\item La droite analytique $\E{1}{A}$ est un espace topologique m\'e\-tri\-sable, localement compact, connexe par arcs et localement connexe par arcs, de dimension topologique 3.
\item Le morphisme de projection $\E{1}{A} \to \Ms(A)$ est ouvert.
\item En tout point $x$ de $\E{1}{A}$, l'anneau local $\Os_{x}$ est hens\'elien, noeth\'erien et r\'egulier.
\item Le principe du prolongement analytique vaut.
\item Le faisceau structural $\Os$ sur $\E{1}{A}$ est coh\'erent.
\end{enumerate}
\end{thmi}

Nous disposons, \`a pr\'esent, de r\'esultats aboutis concernant les propri\'et\'es topologiques et alg\'ebriques de la droite analytique sur l'anneau~$A$. Un peu de travail suppl\'ementaire nous permettra d'en d\'eduire des applications \`a l'\'etude des s\'eries arithm\'etiques, ainsi que nous l'exposerons au chapitre~$7$.

\bigskip

{\bf Morphismes finis}

Le cinqui\`eme chapitre est consacr\'e \`a quelques cas particuliers de morphismes finis entre espaces analytiques. Nous r\'ealisons cette \'etude dans un cadre g\'en\'eral, au-dessus d'un anneau de Banach quelconque. Le r\'esultat principal du chapitre prend la forme suivante.

\begin{thmi}
Soient $(\As,\|.\|)$ un anneau de Banach et~$P$ un polyn\^ome unitaire \`a coefficients dans~$\As$. Sous certaines conditions, l'on peut munir l'anneau quotient $\As[T]/(P(T))$ d'une norme~$\|.\|_{P}$ telle que
\begin{enumerate}[\it i)]
\item le couple $(\As[T]/(P(T)),\|.\|_{P})$ est un anneau de Banach ;
\item le morphisme naturel $\As \to \As[T]/(P(T))$ est born\'e ;
\item le morphisme induit $\varphi : \Ms(\As[T]/(P(T))) \to \Ms(A)$ est ferm\'e et \`a fibres finies;
\item le faisceau $\varphi_{*}\Os$, o\`u~$\Os$ d\'esigne le faisceau structural sur $\Ms(\As[T]/(P(T)))$, est coh\'erent.
\end{enumerate} 
\end{thmi}

Nous d\'emontrons, au passage, un th\'eor\`eme de division de Weierstra{\ss} pour les points rigides des fibres qui nous semble pr\'esenter un int\'er\^et propre.

Nous appliquerons, par la suite, les r\'esultats de ce chapitre \`a certains endomorphismes de la droite analytique au-dessus d'un anneau d'entiers de corps de nombres. Indiquons que nous pensons que les techniques introduites ici permettent d'\'etudier les courbes analytiques au-dessus d'un tel anneau. Dans ce m\'emoire, nous n'en dirons pas plus \`a ce sujet, mais d\'evelopperons ces id\'ees dans un texte \`a venir.

\bigskip

{\bf Espaces de Stein}

Dans le sixi\`eme chapitre, nous reprenons le cadre de la droite affine analytique au-dessus d'un anneau d'entiers de corps de nombres~$A$. Nous cherchons \`a jeter les bases d'une th\'eorie des espaces de Stein pour les parties de cet espace. Les d\'efinitions que nous prenons sont les d\'efinitions cohomologiques habituelles : une partie $P$ de la droite $\E{1}{A}$ est dite de Stein si elle v\'erifie le th\'eor\`eme A :
\begin{quote}
pour tout faisceau coh\'erent $\Fs$ sur $P$ et tout point $x$ de $P$,
la fibre~$\Fs_{x}$ est engendr\'ee par les sections globales $\Fs(P)$
\end{quote}
et le th\'eor\`eme $B$ :
\begin{quote} 
pour tout faisceau coh\'erent $\Fs$ sur $P$ et tout entier $q\in\N^*$,
nous avons $H^q(P,\Fs)=0$.
\end{quote}

L'objet de ce chapitre est de d\'emontrer le th\'eor\`eme suivant.

\begin{thmi}
Soit~$V$ une partie connexe de l'espace~$\Ms(A)$. Soient~$s$ et~$t$ deux nombres r\'eels tels que $0\le s\le t$. Soit~$P$ un polyn\^ome \`a coefficients dans~$\Os(V)$ dont le coefficient dominant est inversible. Les parties suivantes de la droite analytique~$\E{1}{A}$ sont des espaces de Stein :
\begin{enumerate}[\it i)]
\item $\left\{x\in \pi^{-1}(V)\, \big|\, s \le |P(T)(x)| \le t\right\}$ ;
\item $\left\{x\in \pi^{-1}(V)\, \big|\, s \le |P(T)(x)| < t\right\}$ ;
\item $\left\{x\in \pi^{-1}(V)\, \big|\, s < |P(T)(x)| \le t\right\}$ ;
\item $\left\{x\in \pi^{-1}(V)\, \big|\, s < |P(T)(x)| < t\right\}$ ;
\item $\left\{x\in \pi^{-1}(V)\, \big|\, |P(T)(x)| \ge s\right\}$ ;
\item $\left\{x\in \pi^{-1}(V)\, \big|\, |P(T)(x)| > s\right\}$.
\end{enumerate}
\end{thmi}



Nous commen\c{c}ons par traiter le cas des couronnes ferm\'ees. La d\'emonstration que nous proposons reprend la structure de la preuve classique, en g\'eom\'etrie analytique complexe, du fait que les blocs compacts, c'est-\`a-dire les produits de segments r\'eels dans $\C^n$, sont des espaces de Stein. Les ingr\'edients essentiels en sont le lemme de Cousin, qui permet, sous certaines hypoth\`eses, d'\'ecrire une fonction analytique~$f$ d\'efinie sur une intersection de compacts $K^-\cap K^+$ comme diff\'erence d'une fonction analytique~$f^-$ sur~$K^-$ et~$f^+$ sur~$K^+$ et le lemme de Cartan, qui en est la version multiplicative.

La d\'emonstration de ces lemmes met en jeu des outils \`a la fois analytiques et arithm\'etiques. Si les compacts~$K^-$ et~$K^+$ sont d\'efinis, respectivement, par les in\'egalit\'es $|T|\le r$ et $r\le |T|\le s$, il s'agit essentiellement d'\'ecrire une s\'erie de la forme
$$f = \sum_{k\in \Z} a_{k}\, T^k$$
comme diff\'erence $f^--f^+$, avec
$$f^- = \sum_{k\in \N} a_{k}\, T^k \textrm{ et } f^+ = \sum_{k <0} a_{k}\, T^k.$$

Supposons, \`a pr\'esent, que $A=\Z$ et que $K^-$ et $K^+$ sont les compacts de~$\Ms(\Z)$ d\'efinis, respectivement, par les in\'egalit\'es $|p|\le \frac{1}{2}$ et $|p|\ge  \frac{1}{2}$, o\`u $p$ est un nombre premier. Il s'agit alors d'\'ecrire un \'el\'ement de~$\Q_{p}$ comme somme, ou produit, d'un \'el\'ement de~$\Z_{p}$ et d'un \'el\'ement de~$\Z_{(p)}$. Bien entendu, dans un corps de nombres quelconque, ce probl\`eme peut se r\'ev\'eler plus d\'elicat et nous ferons appel au th\'eor\`eme d'approximation forte et \`a la finitude du groupe de Picard.

En ce qui concerne les couronnes ouvertes, le principe de la d\'emonstration consiste \`a construire une exhaustion par des couronnes ferm\'ees. Le fait que les couronnes ouvertes soient de Stein ne d\'ecoule cependant pas formellement de l'existence d'une telle exhaustion. Comme dans le cadre de la g\'eom\'etrie analytique complexe, des propri\'et\'es suppl\'ementaires sont n\'ecessaires et nous sommes amen\'es \`a introduire une notion d'exhaustion de Stein. Signalons que la d\'emonstration des propri\'et\'es requises passe, notamment, par un r\'esultat de fermeture pour les sous-modules d'un module libre, d'int\'er\^et ind\'ependant.

Finalement, le passage du cas des couronnes au cas des parties plus g\'en\'erales qui figurent dans le th\'eor\`eme s'effectue \`a l'aide des r\'esultats sur les morphismes finis d\'emontr\'es au chapitre pr\'ec\'edent.

\bigskip

{\bf Applications}

De m\^eme que la g\'eom\'etrie analytique complexe permet de d\'emontrer des r\'esultats sur les fonctions holomorphes, nous obtenons, \`a l'aide des th\'eor\`emes que nous avons \'etablis concernant la droite affine analytique sur un anneau d'entiers de corps de nombres, des propri\'et\'es des s\'eries arithm\'etiques convergentes (au sens du d\'ebut de l'introduction). C'est l'objet de notre septi\`eme et dernier chapitre. Donnons un exemple de telle propri\'et\'e. Notons $D$ le disque unit\'e ouvert de $\C$.

\begin{thmi}
Soient $E$ et~$F$ deux parties disjointes, ferm\'ees et discr\`etes de $D$ ne contenant pas le point~$0$. Soient~$(n_{a})_{a\in E}$ une famille d'entiers positifs et~$(P_{b})_{b\in F}$ une famille de polyn\^omes \`a coefficients complexes sans terme constant. Nous supposerons que 
\begin{enumerate}
\item quel que soit $a\in E$, $\bar{a}\in E$ et $n_{\bar{a}} = n_{a}$ ;
\item quel que soit $b\in F$, $\bar{b}\in F$ et $P_{\bar{b}} = \overline{P_{b}}$.
\end{enumerate}
Alors il existe $g,h \in \Z[\![T]\!] \cap \Os(D)$, avec~$h\ne 0$, qui v\'erifient les propri\'et\'es suivantes :
\begin{enumerate}[\it i)]
\item la fonction $f=g/h$ est holomorphe sur $D\setminus F$ ;
\item quel que soit $a\in E$, la fonction $f$ s'annule en~$a$ \`a un ordre sup\'erieur \`a~$n_{a}$ ;
\item quel que soit $b\in F$, on a $f(z)- P_{b}\left(\frac{1}{z-b}\right) \in \Os_{b}$ ;
\item on a $f \in\Z[\![T]\!] \cap \Os_{0}$.
\end{enumerate}
\end{thmi}

Ce r\'esultat se d\'emontre par des m\'ethodes cohomologiques. Lorsque la partie~$E$ est vide, nous utilisons la suite exacte courte $0 \to \Os \to \Ms \to \Ms/\Os \to 0$ et le fait que le disque ouvert de centre 0 et de rayon 1 de $\AZ$ est une partie de Stein. Lorsqu'elle ne l'est pas, nous utilisons le m\^eme argument en rempla\c{c}ant le faisceau~$\Os$ par un diviseur de Cartier ad\'equat.

Soit $\Pc$ un ensemble fini de nombres premiers. Notons $N\in\N^*$ leur produit. Il est possible de contr\^oler \'egalement le comportement de~$f$ en tant que fonction m\'eromorphe sur le disque ouvert de centre~$0$ et de rayon~$1$ de~$\C_{p}$, pour tout nombre premier $p\in \Pc$. Il nous faut alors autoriser les coefficients de $g$, de~$h$ et du d\'eveloppement de $f$ en $0$ \`a appartenir \`a $\Z[1/N]$. Bien entendu, nous disposons d'un r\'esultat analogue pour tout corps de nombres.

\bigskip

Nous proposons, ensuite, une application de nos m\'ethodes \`a la noeth\'erianit\'e d'anneaux de s\'eries arithm\'etiques convergentes. Pour l'obtenir, nous nous sommes inspir\'e du th\'eor\`eme suivant de J. Frisch (\emph{cf.} \cite{Frisch}).

\begin{thm*}[J.~Frisch]
Soit $X$ une vari\'et\'e analytique r\'eelle ou complexe. Soit $K$ une partie de~$X$ compacte, semi-analytique et de Stein. Alors l'anneau des fonctions analytiques au voisinage de $K$ est noeth\'erien.
\end{thm*}  

Comme l'ont montr\'e des r\'esultats ult\'erieurs (\emph{cf.} \cite{Siu}, th\'eor\`eme 1), l'hypoth\`ese de semi-analyticit\'e peut \^etre affaiblie. C'est pourquoi nous introduisons ici une notion de partie morcelable. Nous obtenons alors le r\'esultat suivant.

\begin{thmi}
Soit~$A$ un anneau d'entiers de corps de nombres. Soit $L$ une partie de la droite analytique $\E{1}{A}$ compacte, morcelable et de Stein. Alors l'anneau $\Os(L)$ des fonctions analytiques au voisinage de $L$ est noeth\'erien.
\end{thmi}

En appliquant ce th\'eor\`eme aux disques ferm\'es au-dessus des parties semi-analytiques de $\Ms(\Z)$, nous obtenons le r\'esultat suivant.

\begin{cori}
Soient $t$ un entier, $p_{1},\ldots,p_{t}$ des nombres premiers, $r_{1},\ldots,r_{t},r_{\infty}$ des \'el\'ements de l'intervalle~$\of{]}{0,1}{[}$. Alors, l'anneau form\'e des s\'eries 
$$f\in \Z\left[\frac{1}{p_{1}\cdots p_{t}}\right][\![T]\!]$$ 
v\'erifiant les conditions
$$R_{\infty}(f) > r_{\infty} \textrm{ et } \forall i\in\cn{1}{t}, \, R_{p_{i}}(f) > r_{i}$$
est un anneau noeth\'erien.
\end{cori}

Si l'on consid\`ere uniquement des s\'eries \`a coefficients entiers et que l'on n'impose donc des conditions que sur le rayon de convergence complexe, nous retrouvons un r\'esultat de D.~Harbater (\emph{cf.} \cite{Harbater}, th\'eor\`eme 1.8). La preuve qu'il en propose est tr\`es alg\'ebrique : elle consiste \`a d\'ecrire tous les id\'eaux premiers de l'anneau \`a l'aide de manipulations astucieuses sur les s\'eries. Insistons sur le fait que notre d\'emonstration repose sur des arguments g\'eom\'etriques et suit de pr\`es les m\'ethodes de la g\'eom\'etrie analytique complexe. En ce sens, elle nous semble porter des promesses de g\'en\'eralisation. Signalons, enfin, que notre r\'esultat s'\'etend \`a tout anneau d'entiers de corps de nombres.

\bigskip

Pour finir, nous proposons une application au probl\`eme de Galois inverse. L\`a encore, nous proposons une nouvelle d\'emonstration d'un r\'esultat de D.~Harbater (\emph{cf.} \cite{galoiscovers}, corollaire~3.8).

\begin{thmi}
Notons~$\Z_{1^-}[\![T]\!]$ le sous-anneau de~$\Z[\![T]\!]$ form\'e des s\'eries 
$$\sum_{k\ge 0} a_{k}\, T^k$$
qui v\'erifient la condition suivante :
$$\forall r<1,\ \lim_{k \to +\infty}\|a_{k}\|\, r^k = 0.$$

Tout groupe fini est groupe de Galois d'une extension finie et galoisienne du corps $\Frac(\Z_{1^-}[\![T]\!])$.
\end{thmi}

Les m\'ethodes que nous mettons ici en {\oe}uvre nous semblent conceptuellement plus simples et plus g\'eom\'etriques que celles propos\'ees par D.~Harbater. La th\'eorie des espaces de Berkovich nous permet, en effet, d'interpr\'eter l'anneau $\Z_{1^-}[\![T]\!]$ comme un anneau de sections, \`a savoir l'anneau des sections du disque~$\D$, le disque relatif ouvert de rayon~$1$ centr\'e en la section nulle. 

Un groupe fini \'etant donn\'e, nous pouvons alors construire un rev\^etement du disque~$\D$ poss\'edant le groupe de Galois voulu. Nous proc\'edons de fa\c{c}on classique, en exhibant d'abord des rev\^etements cycliques d\'efinis localement, puis en les recollant. La seule \'etape d\'elicate est celle du recollement. C'est le caract\`ere Stein du disque~$\D$, d\'emontr\'e au chapitre pr\'ec\'edent, qui nous permettra de la mener \`a bien.

De nouveau, notre r\'esultat s'\'etend \`a tout anneau d'entiers de corps de nombres. Nous esp\'erons que cette vision tr\`es g\'eom\'etrique du probl\`eme permettra d'y effectuer quelques progr\`es.

\pagebreak

{\bf Remerciements} 

Le pr\'esent m\'emoire pr\'ecise et \'etend les r\'esultats que j'ai obtenus lors de ma th\`ese, effectu\'ee \`a l'institut de recherche math\'ematique de Rennes. Je souhaite remercier ici mes deux directeurs, Antoine Chambert-Loir et Antoine Ducros, pour m'avoir invit\'e \`a me lancer dans ces recherches et m'avoir prodigu\'e sans compter leurs conseils et leurs encouragements. J'ai beaucoup appris \`a leur contact : des math\'ematiques, bien s\^ur, mais \'egalement, je le crois, le m\'etier de math\'ematicien.

J'ai fini de r\'ediger ce texte lors de mon s\'ejour \`a l'universit\'e de Ratisbonne. Je tiens \`a remercier Klaus K\"unnemann, qui m'y a accueilli, pour les excellentes conditions de travail dont j'ai b\'en\'efici\'e et les encouragements qu'il m'a prodigu\'es.

%% file: MZ.pstex_t
\begin{picture}(0,0)%
\includegraphics{MZ.pstex}%
\end{picture}%
\setlength{\unitlength}{3947sp}%
\begingroup\makeatletter\ifx\SetFigFont\undefined%
\gdef\SetFigFont#1#2#3#4#5{%
  \reset@font\fontsize{#1}{#2pt}%
  \fontfamily{#3}\fontseries{#4}\fontshape{#5}%
  \selectfont}%
\fi\endgroup%
\begin{picture}(3008,4390)(4456,-6409)
\put(7209,-2521){\makebox(0,0)[lb]{\smash{{\SetFigFont{12}{14.4}{\rmdefault}{\mddefault}{\updefault}{\color[rgb]{0,0,0}$\Ms(\Z)$}%
}}}}
\put(4471,-5431){\makebox(0,0)[lb]{\smash{{\SetFigFont{12}{14.4}{\rmdefault}{\mddefault}{\updefault}{\color[rgb]{0,0,0}2}%
}}}}
\put(4966,-5626){\makebox(0,0)[lb]{\smash{{\SetFigFont{12}{14.4}{\rmdefault}{\mddefault}{\updefault}{\color[rgb]{0,0,0}3}%
}}}}
\put(7321,-5491){\makebox(0,0)[lb]{\smash{{\SetFigFont{12}{14.4}{\rmdefault}{\mddefault}{\updefault}{\color[rgb]{0,0,0}$p$}%
}}}}
\put(5491,-4051){\makebox(0,0)[lb]{\smash{{\SetFigFont{12}{14.4}{\rmdefault}{\mddefault}{\updefault}{\color[rgb]{0,0,0}$|.|_0$}%
}}}}
\put(6361,-4936){\makebox(0,0)[lb]{\smash{{\SetFigFont{12}{14.4}{\rmdefault}{\mddefault}{\updefault}{\color[rgb]{0,0,0}$|.|_p^\eps$}%
}}}}
\put(5476,-3106){\makebox(0,0)[lb]{\smash{{\SetFigFont{12}{14.4}{\rmdefault}{\mddefault}{\updefault}{\color[rgb]{0,0,0}$|.|_\infty^\eps$}%
}}}}
\put(4741,-4141){\makebox(0,0)[lb]{\smash{{\SetFigFont{12}{14.4}{\rmdefault}{\mddefault}{\updefault}{\color[rgb]{0,0,0}0}%
}}}}
\put(7081,-6331){\makebox(0,0)[lb]{\smash{{\SetFigFont{12}{14.4}{\rmdefault}{\mddefault}{\updefault}{\color[rgb]{0,0,0}$+\infty$}%
}}}}
\put(4711,-3106){\makebox(0,0)[lb]{\smash{{\SetFigFont{12}{14.4}{\rmdefault}{\mddefault}{\updefault}{\color[rgb]{0,0,0}$\eps$}%
}}}}
\put(6421,-5656){\makebox(0,0)[lb]{\smash{{\SetFigFont{12}{14.4}{\rmdefault}{\mddefault}{\updefault}{\color[rgb]{0,0,0}$\eps$}%
}}}}
\put(5881,-5176){\makebox(0,0)[lb]{\smash{{\SetFigFont{12}{14.4}{\rmdefault}{\mddefault}{\updefault}{\color[rgb]{0,0,0}0}%
}}}}
\put(4741,-2251){\makebox(0,0)[lb]{\smash{{\SetFigFont{12}{14.4}{\rmdefault}{\mddefault}{\updefault}{\color[rgb]{0,0,0}1}%
}}}}
\end{picture}%

%% file: A1Cp2.pstex_t
\begin{picture}(0,0)%
\includegraphics{A1Cp2.pstex}%
\end{picture}%
\setlength{\unitlength}{3947sp}%
\begingroup\makeatletter\ifx\SetFigFont\undefined%
\gdef\SetFigFont#1#2#3#4#5{%
  \reset@font\fontsize{#1}{#2pt}%
  \fontfamily{#3}\fontseries{#4}\fontshape{#5}%
  \selectfont}%
\fi\endgroup%
\begin{picture}(6036,5749)(3260,-6870)
\put(4736,-6551){\makebox(0,0)[lb]{\smash{{\SetFigFont{10}{12.0}{\rmdefault}{\mddefault}{\updefault}{\color[rgb]{0,0,0}$-p$}%
}}}}
\put(5861,-1256){\makebox(0,0)[lb]{\smash{{\SetFigFont{10}{12.0}{\rmdefault}{\mddefault}{\updefault}{\color[rgb]{0,0,0}$\infty$}%
}}}}
\put(6056,-6806){\makebox(0,0)[lb]{\smash{{\SetFigFont{10}{12.0}{\rmdefault}{\mddefault}{\updefault}{\color[rgb]{0,0,0}$p-p^2$}%
}}}}
\put(7046,-6581){\makebox(0,0)[lb]{\smash{{\SetFigFont{10}{12.0}{\rmdefault}{\mddefault}{\updefault}{\color[rgb]{0,0,0}$p+p^2$}%
}}}}
\put(6776,-6671){\makebox(0,0)[lb]{\smash{{\SetFigFont{10}{12.0}{\rmdefault}{\mddefault}{\updefault}{\color[rgb]{0,0,0}$p$}%
}}}}
\put(7976,-6491){\makebox(0,0)[lb]{\smash{{\SetFigFont{10}{12.0}{\rmdefault}{\mddefault}{\updefault}{\color[rgb]{0,0,0}$1$}%
}}}}
\put(9281,-5306){\makebox(0,0)[lb]{\smash{{\SetFigFont{10}{12.0}{\rmdefault}{\mddefault}{\updefault}{\color[rgb]{0,0,0}$2$}%
}}}}
\put(3431,-6281){\makebox(0,0)[lb]{\smash{{\SetFigFont{10}{12.0}{\rmdefault}{\mddefault}{\updefault}{\color[rgb]{0,0,0}$-1$}%
}}}}
\put(5741,-6731){\makebox(0,0)[lb]{\smash{{\SetFigFont{10}{12.0}{\rmdefault}{\mddefault}{\updefault}{\color[rgb]{0,0,0}$0$}%
}}}}
\end{picture}%

%% file: Banach.tex
\chapter[Espaces analytiques]{Espaces analytiques sur un anneau de Banach}

Le premier chapitre de ce m\'emoire est consacr\'e aux espaces analytiques sur un anneau de Banach quelconque, au sens de Vladimir G.~Berkovich. Au num\'ero~\ref{defespaceBerko}, nous rappelons les constructions qu'il propose dans l'ouvrage~\cite{rouge}, \`a la fois pour l'espace topologique et le faisceau structural. Nous donnons, en particulier, une description explicite de la droite affine analytique au-dessus d'un corps valu\'e complet quelconque.

Au num\'ero \ref{pcsc}, nous nous int\'eressons \`a certaines parties compactes des espaces analytiques, que nous avons appel\'ees spectralement convexes. Elles pos\-s\`e\-dent notamment la propri\'et\'e d'\^etre hom\'eomorphes \`a des spectres analytiques d'anneaux de Banach que nous pouvons d\'ecrire explicitement. Nous en donnons des exemples et d\'emontrons quelques r\'esultats de permanence \`a leur sujet. Par la suite, les parties spectralement convexes nous seront fort utiles pour mener des raisonnements par r\'ecurrence, puisqu'elles permettent de ramener l'\'etude d'une partie d'un espace affine de dimension~$n$ \`a celle d'un espace de dimension~$0$.

Le num\'ero \ref{parflot} est consacr\'e \`a une application naturelle continue, que nous avons appel\'ee flot, d'une partie de~$\R_{+}$ dans un espace analytique donn\'e. Nous l'\'etudions et comparons les propri\'et\'es des points situ\'es sur une m\^eme trajectoire.



\section{D\'efinitions}\label{defespaceBerko}



\subsection{Spectre analytique d'un anneau de Banach}

Soit $A$ un anneau commutatif unitaire. Par d\'efinition, l'ensemble sous-jacent au spectre Spec($A$) de l'anneau $A$ est l'ensemble des id\'eaux premiers de $A$. D'apr\`es \cite{EGAISpringer}, Introduction, {\bf 13}, il est en bijection avec l'ensemble des classes d'\'equivalence de morphismes unitaires
$$A \to k,$$
o\`u $k$ est un corps. Deux morphismes de $A$ vers des corps $k_{1}$ et $k_{2}$ sont dits \'equivalents s'ils prennent place dans un diagramme commutatif de la forme suivante :
$$\xymatrix{
& & k_{1}\\
A \ar[r] \ar[urr] \ar[drr] & k_{0} \ar[ur] \ar[dr] &\\
&& k_{2}.
}$$

La bijection pr\'ec\'edente peut \^etre d\'ecrite explicitement. Tout d'abord, si \mbox{$A\to k$} est un morphisme unitaire vers un corps, son noyau est un id\'eal premier de $A$ et donc un \'el\'ement de Spec($A$). R\'eciproquement, si $x$ est un point de Spec($A$), il correspond \`a un id\'eal premier $\p_{x}$ de $A$. On construit alors un morphisme de $A$ vers un corps de la fa\c{c}on suivante :
$$A \to A/\p_{x} \to \textrm{Frac}(A/\p_{x}).$$
Le corps $k(x)=\textrm{Frac}(A/\p_{x})$ est appel\'e corps r\'esiduel du point $x$. Par ailleurs, on v\'erifie que tous les morphismes repr\'esentant $x$ se factorisent par le morphisme $A \to k(x)$.

\bigskip

Si nous d\'esirons faire de la g\'eom\'etrie analytique, nous aurons besoin de disposer de notions de normes et de convergence. Nous allons donc consid\'erer non plus un simple anneau, mais un anneau {\it de Banach}. De m\^eme, nous allons remplacer les morphismes vers des corps par des morphismes {\it born\'es}, et donc continus, vers des corps {\it valu\'es}. Rappelons les d\'efinitions de ces notions.

\begin{defi}\index{Anneau de Banach}\index{Morphisme borne@Morphisme born\'e}\index{Norme!d'anneau}
Soient~$\As$ un anneau commutatif unitaire et~$\|.\|$ une application de~$\As$ dans~$\R_{+}$. Nous dirons que l'application~$\|.\|$ est une {\bf norme d'anneau sur l'anneau~$\As$} si elle v\'erifie les propri\'et\'es suivantes :
\begin{enumerate}[\it i)]
\item $(\|f\|=0) \Leftrightarrow (f=0)$ ;
\item $\|1\|=1$ ;
\item $\forall f,g\in\As,\, \|f+g\| \le \|f\| +\|g\|$ ;
\item $\forall f,g\in\As,\, \|fg\| \le \|f\|\,\|g\|$.
\end{enumerate} 

Nous dirons que le couple~$(\As,\|.\|)$ est un {\bf anneau de Banach} si l'application~$\|.\|$ est une norme d'anneau sur l'anneau~$\As$ et si l'espace topologique~$\As$ est complet pour cette norme.
\newcounter{As}\setcounter{As}{\thepage}

Soient~$(\As',\|.\|')$ un anneau de Banach et~$\varphi$ une application de~$\As$ dans~$\As'$. Nous dirons que l'application~$\varphi$ est un {\bf morphisme born\'e d'anneaux de Banach} si l'application~$\varphi$ est un morphisme d'anneaux et s'il existe un nombre r\'eel~$C$ tel que
$$\forall f\in \As,\, \|\varphi(f)\|'\le C\, \|f\|.$$
\end{defi}

\begin{rem}
Cette d\'efinition du caract\`ere born\'e ne co\"{\i}ncide pas avec la d\'efinition habituelle, mais elle est naturelle dans le cadre des morphismes d'anneaux. Nous utiliserons uniquement celle-ci.  
\end{rem}


\begin{defi}\index{Corps!value@valu\'e}\index{Caractere@Caract\`ere}
Nous appellerons {\bf corps valu\'e} tout couple~$(K,|.|)$, o\`u~$K$ est un corps commutatif et~$|.|$ une valeur absolue sur~$K$, c'est-\`a-dire une application de~$K$ dans~$\R_{+}$ qui v\'erifie les propri\'et\'es suivantes :
\begin{enumerate}[\it i)]
\item $(|f|=0) \Leftrightarrow (f=0)$ ;
\item $|1|=1$ ;
\item $\forall f,g\in K,\, |f+g| \le |f\| +|g|$ ;
\item $\forall f,g\in K,\, |fg| = |f|\,|g|$.
\end{enumerate} 

Soit~$(\As,\|.\|)$ un anneau de Banach. Nous appellerons {\bf caract\`ere de l'anneau de Banach $(\As,\|.\|)$} tout morphisme born\'e de la forme
$$\chi : (\As,\|.\|) \to (K,|.|),$$
o\`u $(K,|.|)$ d\'esigne un corps valu\'e complet.
\end{defi} 

\begin{rem}
Dire que le morphisme $\chi : (\As,\|.\|) \to (K,|.|)$ est born\'e signifie, par d\'efinition, qu'il existe un nombre r\'eel $C >0$ tel que, quel que soit $f\in \As$, nous ayons
$$|\chi(f)| \le C\, \|f\|.$$
Soient $f\in \As$ et $n\in\N^*$. Nous avons alors
$$|\chi(f)| = |\chi(f^n)|^{1/n} \le C^{1/n}\, \|f^n\|^{1/n} \le C^{1/n}\, \|f\|.$$
En passant \`a la limite quand $n$ tend vers $+\infty$, nous obtenons
$$|\chi(f)| \le \|f\|.$$
Nous pourrons donc toujours supposer que $C=1$.
\end{rem}

\begin{defi}\index{Caractere@Caract\`ere!equivalence@\'equivalence}
Soit~$(\As,\|.\|)$ un anneau de Banach. Nous dirons que deux {\bf caract\`eres} de~$(\As,\|.\|)$
$$\chi_{1} : (\As,\|.\|) \to (K_{1},|.|_{1}) \textrm{ et } \chi_{2} : (\As,\|.\|) \to (K_{2},|.|_{2})$$ 
sont {\bf \'equivalents} s'il existe un troisi\`eme caract\`ere de~$(\As,\|.\|)$
$$\chi_{0} : (\As,\|.\|) \to (K_{0},|.|_{0})$$ et deux morphismes isom\'etriques 
$$j_{1} : (K_{0},|.|_{0})\to (K_{1},|.|_{1}) \textrm{ et } j_{2} : (K_{0},|.|_{0})\to (K_{2},|.|_{2})$$
qui font commuter le diagramme
$$\xymatrix{
& & & (K_{1},|.|_{1})\\
(\As,\|.\|) \ar[rr]^{\chi_{0}} \ar[urrr]^{\chi_{1}} \ar[drrr]_{\chi_{2}} & & (K_{0},|.|_{0}) \ar[ur]_{j_{1}} \ar[dr]^{j_{2}} &\\
&&& (K_{2},|.|_{2}).
}$$
\end{defi}

Comme dans le cas des sch\'emas, nous pouvons d\'ecrire les classes d'\'equivalence de caract\`eres d'une fa\c{c}on explicite. \`A cet effet, nous aurons besoin de la d\'efinition suivante.

\begin{defi}\index{Semi-norme!multiplicative bornee@multiplicative born\'ee}
Soit~$(\As,\|.\|)$ un anneau de Banach. Une {\bf semi-norme multiplicative born\'ee sur l'anneau de Banach $(\As,\|.\|)$} est une application $|.| : \As \to \R_{+}$ qui v\'erife les propri\'et\'es suivantes :
\begin{enumerate}[\it i)]
\item $|0|=0$ ;
\item $|1|=1$ ;
\item $\forall f,g\in \As,\, |f+g|\le |f|+|g|$ ;
\item $\forall f,g\in \As, \,|fg|= |f||g|$ ;
\item $\exists C >0,\, \forall f\in \As,\, |f|\le C\|f\|$.  
\end{enumerate}
\end{defi}

\begin{rem}
Le m\^eme raisonnement que pour les caract\`eres montre que l'on peut supposer que $C=1$.
\end{rem}

Soit~$(\As,\|.\|)$ un anneau de Banach. L'ensemble des classes d'\'equivalence de caract\`eres sur $(\As,\|.\|)$ est en bijection avec l'ensemble des semi-normes multiplicatives born\'ees sur $(\As,\|.\|)$. Nous pouvons d\'ecrire cette bijection explicitement. \`A tout caract\`ere  
$$\chi : (\As,\|.\|) \to (K,|.|),$$
on associe la semi-norme multiplicative
$$\As \xrightarrow[]{\chi} K \xrightarrow[]{|.|} \R_{+}.$$
Elle est born\'ee car le morphisme $\chi$ est born\'e. On v\'erifie imm\'ediatement que la semi-norme obtenue ne d\'epend que de la classe d'\'equivalence du caract\`ere~$\chi$.

R\'eciproquement, soit $|.|_{x}$ une semi-norme multiplicative born\'ee sur $(\As,\|.\|)$. L'ensemble 
$$\p_{|.|_{x}} = \{f\in \As,\, |f|_{x}=0\}$$   \newcounter{px}\setcounter{px}{\thepage}
est un id\'eal premier de $\As$. Le quotient $A/\p_{|.|_{x}}$ est un anneau int\`egre sur lequel la semi-norme~$|.|_{x}$ induit une valeur absolue. Nous noterons $\Hs(|.|_{x})$ \newcounter{Hsx}\setcounter{Hsx}{\thepage} le compl\'et\'e du corps des fractions de cet anneau pour cette valeur absolue. La construction nous fournit un morphisme 
$$\As \to \Hs(|.|_{x}).$$
On v\'erifie sans peine qu'il est born\'e et donc que c'est un caract\`ere. Comme dans le cas des sch\'emas, tout caract\`ere repr\'esentant la semi-norme multiplicative $|.|_{x}$ se factorise par le caract\`ere $\As \to \Hs(|.|_{x}).$ 

Ces consid\'erations motivent la d\'efinition suivante.

\begin{defi}[V.~Berkovich]\index{Spectre analytique}\newcounter{MAs}\setcounter{MAs}{\thepage}
Soit~$(\As,\|.\|)$ un anneau de Banach. On appelle {\bf spectre analytique de l'anneau de Banach $(\As,\|.\|)$} et l'on note $\Ms(\As,\|.\|)$, ou plus simplement $\Ms(\As)$ si aucune ambigu\"{\i}t\'e n'en r\'esulte, l'ensemble des semi-normes multiplicatives born\'ees sur~$(\As,\|.\|)$. 
\end{defi}

Soit~$(\As,\|.\|)$ un anneau de Banach. Soient $f$ un \'el\'ement de~$\As$ et $x$ un point de $\Ms(\As)$. Notons $|.|_{x}$ la semi-norme multiplicative born\'ee sur~$\As$ associ\'ee au point~$x$. Nous noterons $\p_{x}=\p_{|.|_{x}}$ l'id\'eal premier d\'efini pr\'ec\'edemment. Nous appellerons {\bf corps r\'esiduel compl\'et\'e}\index{Corps!residuel complete@r\'esiduel compl\'et\'e} du point $x$ et noterons $\Hs(x)$ le corps $\Hs(|.|_{x})$ d\'efini pr\'e\-c\'e\-dem\-ment. Nous noterons $f(x)$\newcounter{fdex}\setcounter{fdex}{\thepage} l'image de l'\'el\'ement $f$ de $\As$ par le caract\`ere $\As \to \Hs(x)$. Le corps $\Hs(x)$ est muni d'une valeur absolue, que nous noterons toujours $|.|$. Cela n'entra\^{\i}nera aucune confusion. Avec ces notations, nous avons donc 
$$|f(x)| = |f|_{x}.$$
Comme les notations l'indiquent, nous consid\'ererons d\'esormais les \'el\'ements de~$\As$ comme des fonctions sur l'espace $\Ms(\As)$.

\bigskip

Munissons, \`a pr\'esent, le spectre analytique $\Ms(\As)$ d'une topologie : la topologie la plus grossi\`ere rendant continues les applications d'\'evaluation, c'est-\`a-dire les applications de la forme 
$$\begin{array}{ccc}
\Ms(\As) & \to & \R_{+}\\
x & \mapsto & |f(x)|
\end{array},$$
o\`u~$f$ est un \'el\'ement de~$\As$. Cette topologie est \'egalement celle de la convergence faible, ou encore celle induite par la topologie produit sur~$\R^\As$. Le spectre analytique~$\Ms(\As)$ v\'erifie alors des propri\'et\'es remarquables (\emph{cf.} \cite{rouge}, th\'eor\`eme 1.2.1).

\begin{thm}[V.~Berkovich]
Le spectre analytique $\Ms(\As)$ est un espace topologique compact. Si l'anneau $\As$ n'est pas nul, cet espace n'est pas vide.
\end{thm}

\subsection{Espace affine analytique}\index{Espace affine analytique|(}

Soit $(\As,\|.\|)$ un anneau de Banach. Maintenant que nous avons d\'efini le spectre analytique de cet anneau, nous pouvons d\'efinir ce qu'est l'espace affine au-dessus de celui-ci. Soit $n\in\N$. 

\begin{defi}[V.~Berkovich]
On appelle {\bf espace affine analytique de di\-men\-sion~$n$ sur $(\As,\|.\|)$} l'ensemble des semi-normes multiplicatives sur $\As[T_{1},\ldots,T_{n}]$ dont la restriction \`a $(\As,\|.\|)$ est born\'ee. Nous le noterons $\A^{n,\mathrm{an}}_{\As}$.
\newcounter{AnAs}\setcounter{AnAs}{\thepage}
\end{defi}

En reprenant le raisonnement du paragraphe pr\'ec\'edent, on montre que l'ensemble~$\A^{n,\mathrm{an}}_{\As}$ est en bijection avec l'ensemble des classes d'\'equivalence de morphismes
$$\As[T_{1},\ldots,T_{n}] \to K,$$
o\`u $K$ est un corps valu\'e complet, dont la restriction \`a $\As$ est born\'ee. Comme pr\'e\-c\'e\-dem\-ment, nous associons \`a chaque point $x$ de $\A^{n,\mathrm{an}}_{\As}$ un id\'eal premier~$\p_{x}$ et un corps r\'esiduel compl\'et\'e~$\Hs(x)$. Pour tout \'el\'ement $f$ de $\As[T_{1},\ldots,T_{n}]$, nous d\'esignons par~$f(x)$ l'image de~$f$ par le morphisme $\As[T_{1},\ldots,T_{n}]\to \Hs(x)$.
\index{Corps!residuel complete@r\'esiduel compl\'et\'e}
\newcounter{pxbis}\setcounter{pxbis}{\thepage}
\newcounter{Hsxbis}\setcounter{Hsxbis}{\thepage}
\newcounter{fdexbis}\setcounter{fdexbis}{\thepage}

Nous munissons \'egalement l'espace $\A^{n,\mathrm{an}}_{\As}$ de la topologie la plus grossi\`ere pour laquelle les applications d'\'evaluation sont continues. Il v\'erifie alors encore certaines propri\'et\'es topologiques (\emph{cf.} \cite{rouge}, remarque 1.5.2.(i)). Nous les red\'emontrons ici.

\begin{prop}\label{disquecompact}\index{Compacite@Compacit\'e!disques}
Pour tout nombre r\'eel positif~$r$, la partie de l'espace analytique~$\E{n}{\As}$ d\'efinie par
$$\left\{x\in\E{n}{\As}\, \big|\, \forall i\in\cn{1}{n},\, |T_{i}(x)|\le r\right\}$$
est compacte.
\end{prop}
\begin{proof}
L'application $\|.\|_{r} : \As[T_{1},\ldots,T_{n}] \to \R_{+}$ d\'efinie par
$$
\left\| \sum_{(k_{1},\ldots,k_{n})\in\N^n} a_{k_{1},\ldots,k_{n}}\, T_{1}^{k_{1}}\ldots T_{n}^{k_{n}} \right\|_{r} = \sum_{(k_{1},\ldots,k_{n})\in\N^n} \|a_{k_{1},\ldots,k_{n}}\| \, r^{k_{1}+\cdots+ k_{n}}
$$
est une norme sur la $\As$-alg\`ebre~$\As[T_{1},\ldots,T_{n}]$. Notons~$\Bs$ le compl\'et\'e de l'anneau~$\As$ pour cette norme. L'application naturelle
$$\As[T_{1},\ldots,T_{n}] \to \Bs$$
est born\'ee sur~$\As$. Elle induit donc un morphisme
$$\varphi : \Ms(\Bs) \to \E{n}{\As}.$$


Posons
$$K = \left\{x\in\E{n}{\As}\, \big|\, \forall i\in\cn{1}{n},\, |T_{i}(x)|\le r\right\}.$$
Montrons que l'image de~$\varphi$ contient la partie~$K$. Soit~$x$ un point de~$K$. Il est associ\'e \`a un morphisme
$$\chi_{x} : \As[T_{1},\ldots,T_{n}] \to \Hs(x),$$
qui est born\'e sur~$\As$. Pour tout \'el\'ement~$i$ de~$\cn{1}{n}$, nous avons
$$|T_{i}(x)|\le r =\|T_{i}\|_{r}.$$
On en d\'eduit que le morphisme~$\chi_{x}$ est born\'e lorsque l'on munit l'alg\`ebre $\As[T_{1},\ldots,T_{n}]$ de la norme~$\|.\|_{r}$. Par cons\'equent, le morphisme~$\chi_{x}$ se factorise par un morphisme
$$\Bs \to \Hs(x).$$
On en d\'eduit que le point~$x$ appartient \`a l'image du morphisme~$\varphi$.

Puisque l'espace~$\Ms(\Bs)$ est compact, l'image du morphisme~$\varphi$ l'est \'egalement. Par d\'efinition de la topologie de l'espace~$\E{n}{\As}$, la partie~$K$ est ferm\'ee. Puisqu'elle est contenue dans l'image du morphisme~$\varphi$, elle est compacte. 
\end{proof}

Notons $\pi : \E{n}{\As} \to \Ms(\As)$ l'application de projection induite par le morphisme $\As \to \As[T_{1},\ldots,T_{n}]$.

\begin{cor}\label{partiecompacte}\index{Compacite@Compacit\'e!lemniscates}
Soit~$U$ une partie de~$\Ms(\As)$ et $P_{1}(T_{1}),\ldots,P_{n}(T_{n})$ des polyn\^omes \`a coefficients dans~$\Os(U)$ dont le coefficient dominant est inversible. Pour toute partie compacte~$V$ de~$U$ et tout \'el\'ement~$r$ de~$\R_{+}$, la partie de l'espace analytique~$\E{n}{\As}$ d\'efinie par
$$\left\{x\in \pi^{-1}(V)\, \big|\, \forall i\in\cn{1}{n},\, |P_{i}(T_{i})(x)|\le r\right\}$$
est compacte.
\end{cor}
\begin{proof}
Soit~$i$ un \'el\'ement de~$\cn{1}{n}$. Il existe un entier~$d_{i}$, un \'el\'ement~$a_{i,d_{i}}$ de~$\Os(U)$ inversible et des \'el\'ements $a_{i,d_{i}-1},\ldots,a_{i,0}$ de~$\Os(U)$ tels que
$$P_{i}(T_{i}) = \sum_{k=0}^{d_{i}} a_{i,k}\, T_{i}^k \textrm{ dans } \Os(U)[T_{i}].$$
Puisque la fonction~$a_{i,d_{i}}$ est inversible sur~$U$, la quantit\'e
$$m_{i}=\inf_{v\in V} (|a_{i,d_{i}}(v)|)$$
est strictement positive. Pour tout \'el\'ement~$x$ de~$\pi^{-1}(V)$, nous avons donc
$$\begin{array}{rcl}
|P_{i}(T_{i})(x)| & \ge & \disp |a_{i,d_{i}}(x)|\, |T_{i}(x)|^{d_{i}} - \sum_{k=0}^{d_{i}-1} |a_{i,k}(x)|\, |T_{i}(x)|^k\\
& \ge & \disp m_{i}\, |T_{i}(x)|^{d_{i}} - \sum_{k=0}^{d_{i}-1} \|a_{i,k}\|_{V}\, |T_{i}(x)|^k.
\end{array}$$
La fonction
$$\alpha_{i} : t\in\R \mapsto m_{i}\, t^{d_{i}} - \sum_{k=0}^{d_{i}-1} \|a_{i,k}\|_{V}\, t^k \in\R$$
tend vers~$+\infty$ quand~$t$ tend vers~$+\infty$. Par cons\'equent, il existe un \'el\'ement~$s_{i}$ de~$\R_{+}$ tel que, quel que soit~$t>s_{i}$, on ait~$\alpha_{i}(t)>r$. Pour tout \'el\'ement~$x$ de~$\pi^{-1}(V)$ v\'erifiant $|P_{i}(T_{i})(x)| \le r$, nous avons donc~$|T_{i}(x)|\le s_{i}.$ 

Posons~$s=\max_{1\le i\le n}(s_{i})$. La partie 
$$K = \left\{x\in \pi^{-1}(V)\, \big|\, \forall i\in\cn{1}{n},\, |P_{i}(T_{i})(x)|\le r\right\}$$
est ferm\'ee dans~$\E{n}{\As}$ puisque la partie~$V$ est ferm\'ee. En outre, elle est contenue dans la partie
$$\left\{x\in\E{n}{\As}\, \big|\, \forall i\in\cn{1}{n},\, |T_{i}(x)|\le s\right\},$$
qui est compacte, en vertu du lemme pr\'ec\'edent. On en d\'eduit que la partie~$K$ est compacte.
\end{proof}

\begin{thm}[V.~Berkovich]\label{topoAn}\index{Separation@S\'eparation}\index{Compacite@Compacit\'e!sigma-compacite@$\sigma$-compacit\'e}
L'espace analytique $\A^{n,\mathrm{an}}_{\As}$ est un espace topologique s\'epar\'e, $\sigma$-compact et localement compact.
\end{thm}
\begin{proof}
Soient~$x$ et~$y$ deux points distincts de l'espace~$\E{n}{\As}$. Il existe alors un \'el\'ement~$f$ de $\As[T_{1},\ldots,T_{n}]$ tel que $|f(x)|\ne |f(y)|$. Quitte \`a \'echanger les points~$x$ et~$y$, nous pouvons supposer que $|f(x)|<|f(y)|$. Soit~$r$ un \'el\'ement de l'intervalle~$\of{]}{|f(x)|,|f(y)|}{[}$. Les ouverts
$$\left\{\left. z\in \E{n}{\As}\, \right| |f(z)|<r\right\} \textrm{ et } \left\{\left. z\in \E{n}{\As}\, \right| |f(z)|>r\right\}$$
s\'eparent les points~$x$ et~$y$. Par cons\'equent, l'espace~$\E{n}{\As}$ est s\'epar\'e.

\bigskip

L'espace~$\E{n}{\As}$ est r\'eunion des espaces
$$D_{n} = \left\{x\in\E{n}{\As}\, \big|\, \forall i\in\cn{1}{n},\, |T_{i}(x)|\le n\right\},$$
pour~$n$ d\'ecrivant~$\N$. D'apr\`es la proposition \ref{disquecompact}, ces espaces sont compacts. On en d\'eduit que l'espace~$\E{n}{\As}$ est $\sigma$-compact.

En outre, par d\'efinition de la topologie, tout point poss\`ede un syst\`eme fondamental de voisinages ferm\'es. Puisque tout point est contenu dans l'int\'erieur de l'espace~$D_{n}$, pour un certain entier positif~$n$, et que cet espace est compact, on en d\'eduit tout point poss\`ede un syst\`eme fondamental de voisinages compacts.
\end{proof}

\bigskip

Donnons, \`a pr\'esent, quelques exemples de points d'espaces analytiques. Nous nous restreindrons au cas o\`u l'anneau de Banach $(\As,\|.\|)$ est un corps valu\'e complet $(k,|.|)$. Son spectre analytique~$\Ms(k)$ est alors constitu\'e d'un seul point. Remarquons que l'espace analytique~$\E{n}{k}$ contient l'ensemble~$k^n$. En effet, \`a tout \'el\'ement $\alphab=(\alpha_{1},\ldots,\alpha_{n})$ de~$k^n$, nous pouvons associer le point de~$\E{n}{k}$ d\'efini par
$$\begin{array}{ccc}
k[T_{1},\ldots,T_{n}] & \to & \R_{+}\\
P(T_{1},\ldots,T_{n}) & \mapsto & |P(\alpha_{1},\ldots,\alpha_{n})|
\end{array}.$$
Nous noterons encore~$\alphab$ ce point. Un tel point sera appel\'e point rationnel de l'espace analytique~$\E{n}{k}$. En voici une d\'efinition \'equivalente.  
\newcounter{alphab}\setcounter{alphab}{\thepage}

\begin{defi}\index{Point!rationnel}
Soient~$(k,|.|)$ un corps valu\'e complet. Nous dirons qu'un point~$x$ de l'espace analytique~$\E{n}{k}$ est un {\bf point rationnel} si l'extension de corps $k\to \Hs(x)$ est un isomorphisme.
\end{defi}

En g\'en\'eral, l'espace analytique~$\E{n}{k}$ contient beaucoup plus de points que l'espace~$k^n$. C'est en particulier le cas si le corps~$k$ n'est pas alg\'ebriquement clos et si~$n\ge 1$. Consid\'erons une cl\^oture alg\'ebrique~$\bar{k}$ du corps~$k$. La valeur absolue~$|.|$ sur~$k$ se prolonge de fa\c{c}on unique en une valeur absolue sur~$\bar{k}$, que nous noterons encore~$|.|$. \`A tout \'el\'ement $\betab=(\beta_{1},\ldots,\beta_{n})$ de~$\bar{k}^n$, nous pouvons associer le point de~$\E{n}{k}$ d\'efini par
$$\begin{array}{ccc}
k[T_{1},\ldots,T_{n}] & \to & \R_{+}\\
P(T_{1},\ldots,T_{n}) & \mapsto & |P(\beta_{1},\ldots,\beta_{n})|
\end{array}.$$
Nous noterons encore~$\betab$ ce point. Attention, cependant : si~$\sigma$ d\'esigne un \'el\'ement de Gal($\bar{k}/k$), les points $(\beta_{1},\ldots,\beta_{n})$ et $(\sigma(\beta_{1}),\ldots,\sigma(\beta_{n}))$ co\"incident ! Un tel point sera appel\'e point rigide de l'espace analytique~$\E{n}{k}$. En voici une d\'efinition \'equivalente.  

\begin{defi}\index{Point!rigide}
Soient~$(k,|.|)$ un corps valu\'e complet et~$n$ un nombre entier positif. Nous dirons qu'un point~$x$ de l'espace analytique~$\E{n}{k}$ est un {\bf point rigide} si l'extension de corps $k\to \Hs(x)$ est une extension finie.
\end{defi}

\bigskip

Dans les num\'eros qui suivent, nous d\'ecrivons explicitement l'espace et sa topologie dans quelques cas simples. Si le corps~$k$ est archim\'edien, nous ferons le lien entre l'espace~$\A^{n,\mathrm{an}}_{k}$ et les espaces analytiques r\'eels et complexes usuels. Si le corps~$k$ est ultram\'etrique, nous nous contenterons de d\'ecrire la droite~$\E{1}{k}$. Nous observerons, en particulier, qu'elle contient beaucoup plus de points que~$\bar{k}$.

\index{Espace affine analytique|)}

\subsubsection{Espace affine analytique sur un corps archim\'edien}\index{Espace affine analytique!sur un corps archimedien@sur un corps archim\'edien|(}

Commen\c{c}ons par supposer que le corps $(k,|.|)$ est un corps muni d'une valeur absolue archim\'edienne pour laquelle il est complet. D'apr\`es \cite{algcom56}, VI, \S 6, \no 4, th\'eor\`eme 2, il existe un \'el\'ement~$s$ de l'intervalle $\of{]}{0,1}{]}$ tel que le corps valu\'e $(k,|.|)$ soit isom\'etriquement isomorphe au corps $(\R,|.|_{\infty}^s)$ ou au corps $(\C,|.|_{\infty}^s)$, o\`u $|.|_{\infty}$ d\'esigne la valeur absolue usuelle. 

Supposons que $(k,|.|)=(\C,|.|_{\infty})$. Soit~$n$ un entier positif. Nous savons que les points de l'espace $\E{n}{\C}$ sont en bijection avec les classes d'\'equivalences de caract\`eres de $\C[T_{1},\ldots,T_{n}]$. Soit 
$$\chi : \C[T_{1},\ldots,T_{n}] \to L$$
un tel caract\`ere. D'apr\`es le th\'eor\`eme de Gelfand-Mazur (\emph{cf.} \cite{algcom56}, VI, \S 6, \no 4, th\'eor\`eme~1), le corps~$L$ est isomorphe \`a~$\C$. Posons 
$$\alphab = (\chi(T_{1}),\ldots,\chi(T_{n})) \in \C^n.$$
Alors, le caract\`ere~$\chi$ n'est autre que le morphisme \'evaluation au point~$\alphab$ de~$\C^n$. On en d\'eduit que les ensembles~$\E{n}{\C}$ et~$\C^n$ sont en bijection, autrement dit, tous les points de l'espace analytique~$\E{n}{\C}$ sont rationnels. D'autre part, il est clair que les topologies co\"{\i}ncident. Les espaces~$\E{n}{\C}$ et~$\C^n$ sont donc hom\'eomorphes.

\bigskip

Supposons, \`a pr\'esent, que $(k,|.|)=(\R,|.|_{\infty})$. Soit~$n$ un entier positif. Le m\^eme raisonnement que pr\'ec\'edemment montre que l'espace~$\E{n}{\R}$ est hom\'eomorphe au quotient de l'espace~$\C^n$ par la conjugaison complexe. En particulier, tous les points de l'espace analytique~$\E{n}{\R}$ sont rigides.

\index{Espace affine analytique!sur un corps archimedien@sur un corps archim\'edien|)}





\subsubsection{Droite sur un corps trivialement valu\'e}\label{droitetrivval}\index{Droite affine analytique!sur un corps trivialement value@sur un corps trivialement valu\'e|(}

Dans cette partie, nous supposerons que le corps $k$ est muni de la valeur absolue triviale $|.|_{0}$. Nous nous contenterons de d\'ecrire la droite affine analytique~$\E{1}{k}$. 

Soit~$x$ un point de~$\E{1}{k}$. Il lui correspond une semi-norme multiplicative~$|.|_{x}$ born\'ee sur~$k$. Notons 
$$\p_{x} = \left\{f\in k[T]\, \big|\, |f|_{x}=0\right\}.$$ 
C'est un id\'eal premier de $k[T]$. 

Supposons, tout d'abord, que l'id\'eal~$\p_{x}$ n'est pas l'id\'eal nul. Il existe alors un polyn\^ome irr\'eductible $P(T)$ de $k[T]$ qui engendre l'id\'eal $\p_{x}$. La semi-norme multiplicative $|.|_{x}$ induit une valeur absolue sur le quotient 
$$k[T]/\p_{x} = k[T]/(P(T)),$$
qui est une extension finie du corps $k$. Cette valeur absolue ne peut \^etre que la valeur absolue triviale. Par cons\'equent, nous avons
$$|.|_{x} : 
\begin{array}{ccc}
k[T] & \to & \R_{+}\\
Q(T) & \mapsto & 
\left\{
\begin{array}{cl}
0 & \textrm{si } P\,|\,Q\\
1 & \textrm{sinon}
\end{array}
\right.
\end{array}.$$
Nous noterons $\eta_{P,0}$ le point de $\E{1}{k}$ correspondant. Nous avons
$$\Hs(\eta_{P,0}) = k[T]/(P(T)).$$
\newcounter{etaPzero}\setcounter{etaPzero}{\thepage}

Supposons, \`a pr\'esent, que l'id\'eal~$\p_{x}$ est l'id\'eal nul. La semi-norme multiplicative~$|.|_{x}$ est alors en fait une valeur absolue sur~$k[T]$. Par hypoth\`ese, la restriction de cette valeur absolue \`a~$k$ est born\'ee par la valeur absolue triviale. En particulier, pour tout entier positif~$n$, nous avons $|n.1|_{x}\le 1$. On en d\'eduit que la valeur absolue~$|.|_{x}$ est ultram\'etrique en utilisant le lemme classique suivant.

\begin{lem}\label{vaum}\index{Valeur absolue!ultrametrique@ultram\'etrique}
Soit $(k,|.|)$ un corps valu\'e. La valeur absolue $|.|$ est ul\-tra\-m\'e\-tri\-que si, et seulement si, il existe un nombre r\'eel~$C$ tel que, pour tout entier positif~$n$, nous ayons $|n.1| \le C$.
\end{lem}
\begin{proof}
Supposons que la valeur absolue $|.|$ est ultram\'etrique. En utilisant l'in\'egalit\'e ultram\'etrique et le fait que $|1|=1$, on montre par r\'ecurrence que, pour tout entier positif~$n$, nous avons $|n.1|\le 1$.

Supposons qu'il existe un nombre r\'eel~$C$ tel que, pour tout entier positif~$n$, nous ayons $|n.1| \le C$. Soient $a,b\in k$. Soit $p\in\N^*$. Nous avons
$${\renewcommand{\arraystretch}{2}\begin{array}{rcl}
|a+b|^p &=& |(a+b)^p|\\
&=& \disp\left|\sum_{i=0}^p C_{p}^i\, a^i\, b^{p-i} \right|\\
&\le&\disp\sum_{i=0}^p \left|C_{p}^i\right|\, |a|^i\, |b|^{p-i}\\
&\le& p\, C\, \max(|a|,|b|)^p. 
\end{array}}$$
En \'elevant l'in\'egalit\'e obtenue \`a la puissance $1/p$ et en faisant tendre $p$ vers $+\infty$, on obtient
$$|a+b| \le \max(|a|,|b|).$$
\end{proof}

\begin{figure}[htb]
\begin{center}
\input{valtriv.pstex_t}
\caption{Droite analytique sur un corps trivialement valu\'e.}
\end{center}
\end{figure}
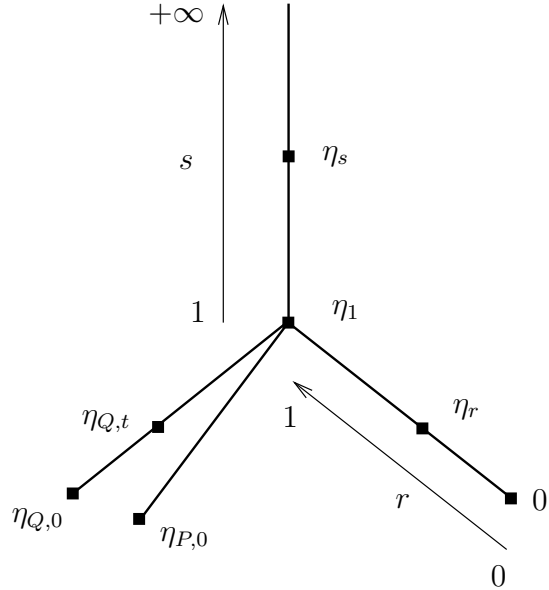

Nous allons distinguer deux cas. Supposons, tout d'abord, que $|T|_{x}\le 1$. On montre alors facilement que, quel que soit $f\in k[T]$, nous avons
$$|f|_{x}\le 1.$$
L'in\'egalit\'e ultram\'etrique assure alors que la partie 
$$\p'_{x} = \left\{ f \in k[T]\, \big|\, |f|_{x} < 1\right\}$$
est un id\'eal premier de $k[T]$. Si cet id\'eal est nul, alors nous avons $|.|_{x} = |.|_{0}$. Nous appellerons {\bf point de Gau{\ss}}\index{Point!de Gau{\ss}} ce point. Nous le noterons~$\eta_{1}$. \newcounter{etauntriv}\setcounter{etauntriv}{\thepage}

Dans les autres cas, l'id\'eal $\p'_{x}$ est engendr\'e par un polyn\^ome irr\'eductible $P$ de $k(T)$. Notons $v_{P}$ la valuation $P$-adique sur $k[T]$. Il existe $r\in\of{]}{0,1}{[}$ tel que $|P|_{x}=r$. Pour tout \'el\'ement~$Q(T)$ de~$k[T]$, nous avons alors 
$$|Q|_{x} = r^{v_{P}(Q)}.$$
Nous noterons $\eta_{P,r}$ le point de $\E{1}{k}$ correspondant. Le corps r\'esiduel compl\'et\'e $\Hs(\eta_{P,r})$ en ce point est le compl\'et\'e du corps $k(T)$ pour la topologie $P$-adique. Si $P(T)=T$, nous noterons~$\eta_{r}$ le point correspondant. Le corps r\'esiduel compl\'et\'e $\Hs(\eta_{r})$ est alors isomorphe au corps des s\'eries de Laurent~$k(\!(T)\!)$. 
\newcounter{etaPr}\setcounter{etaPr}{\thepage}
\newcounter{etar}\setcounter{etar}{\thepage}

Supposons, \`a pr\'esent, que $|T|_{x} >1$. Il existe $r>1$ tel que $|T|_{x}=r$. L'in\'egalit\'e ultram\'etrique montre alors que, quel que soit $Q(T)\in k[T]$, nous avons
$$|Q|_{x} = r^{\deg(Q)}.$$
Nous noterons $\eta_{r}$ le point de $\E{1}{k}$ correspondant. Le corps r\'esiduel compl\'et\'e $\Hs(\eta_{r})$ en ce point est isomorphe au corps $k(\!(T^{-1})\!)$.

Introduisons encore quelques notations. Pour $\alpha\in k$ et $r\in\of{[}{0,1}{]}$, nous noterons
$$\eta_{\alpha,r} = \eta_{T-\alpha,r}.$$
\newcounter{etaalphar}\setcounter{etaalphar}{\thepage}
Si~$r=0$, nous noterons parfois simplement $\alpha$ le point $\eta_{\alpha,0}$. 
\newcounter{alpha}\setcounter{alpha}{\thepage}

\bigskip

Pour finir, d\'ecrivons la topologie de la droite~$\E{1}{k}$. Nous ne d\'emontrerons pas les r\'esultats qui suivent. Pour se faire une id\'ee des preuves, le lecteur int\'eress\'e peut se reporter au num\'ero \ref{descriptionMA}, o\`u nous d\'ecrivons la topologie du spectre d'un anneau d'entiers de corps de nombres. La topologie des branches est particuli\`erement simple. En effet, pour tout polyn\^ome irr\'eductible~$P(T)$ de~$k[T]$, l'application
$$\begin{array}{ccc}
\of{[}{0,1}{]} & \to & \E{1}{k}\\
r & \mapsto & \eta_{P,r}
\end{array}
$$
r\'ealise un hom\'eomorphisme sur son image. De m\^eme, l'application
$$\begin{array}{ccc}
\of{[}{1,+\infty}{[} & \to & \E{1}{k}\\
r & \mapsto & \eta_{r}
\end{array}
$$
r\'ealise un hom\'eomorphisme sur son image. 

Afin d'\^etre complets, il nous reste \`a d\'ecrire un syst\`eme fondamental de voisinages du point de Gau{\ss}~$\eta_{1}$ ; l'ensemble des parties de~$\E{1}{k}$ qui contiennent un voisinage du point~$\eta_{1}$ dans un nombre fini de branches et la totalit\'e des branches restantes en est un. 

\index{Droite affine analytique!sur un corps trivialement value@sur un corps trivialement valu\'e|)}

\subsubsection{Droite sur un corps ultram\'etrique quelconque}\label{droiteum}\index{Droite affine analytique!sur un corps ultrametrique@sur un corps ultram\'etrique|(}

Il est en fait possible de d\'ecrire la droite analytique au-dessus de tout corps ultram\'etrique complet. Nous allons nous limiter au cas des corps qui sont \'egalement alg\'ebriquement clos. Cette restriction ne nuit pas \`a la g\'en\'eralit\'e de notre propos. En effet, d'apr\`es \cite{rouge}, corollaire 1.3.6, si $k$ d\'esigne un corps valu\'e complet, $\bar{k}$ l'une de ses cl\^otures alg\'ebriques et $\hat{\bar{k}}$ le compl\'et\'e de cette derni\`ere, alors le groupe de Galois Gal($\bar{k}/k$) agit sur $\hat{\bar{k}}$ et le morphisme naturel 
$$\E{1}{\hat{\bar{k}}} / \textrm{Gal}(\bar{k}/k) \xrightarrow[]{\sim} \E{1}{k}$$
est un isomorphisme.

\bigskip

Nous supposerons donc, d\'esormais, que $k$ est un corps ultram\'etrique complet alg\'ebriquement clos. Nous reprenons la description donn\'ee par V.~Berkovich  au num\'ero 1.4.4 de l'ouvrage~\cite{rouge}. Il distingue quatre types de points. Soit $\alpha\in k$. L'application d'\'evaluation
$$\begin{array}{ccc}
k[T] & \to & \R_{+}\\
P(T) & \mapsto & |P(\alpha)|
\end{array}$$
d\'efinit une semi-norme multiplicative sur~$k[T]$ born\'ee sur~$k$ et donc un point de~$\E{1}{k}$. Nous noterons~$\alpha$ 
ce point. Un tel point est dit {\bf de type 1}\index{Point!de type 1}. En ce point le corps r\'esiduel compl\'et\'e est simplement 
$$\Hs(\alpha)=k.$$
\newcounter{alphaac}\setcounter{alphaac}{\thepage}

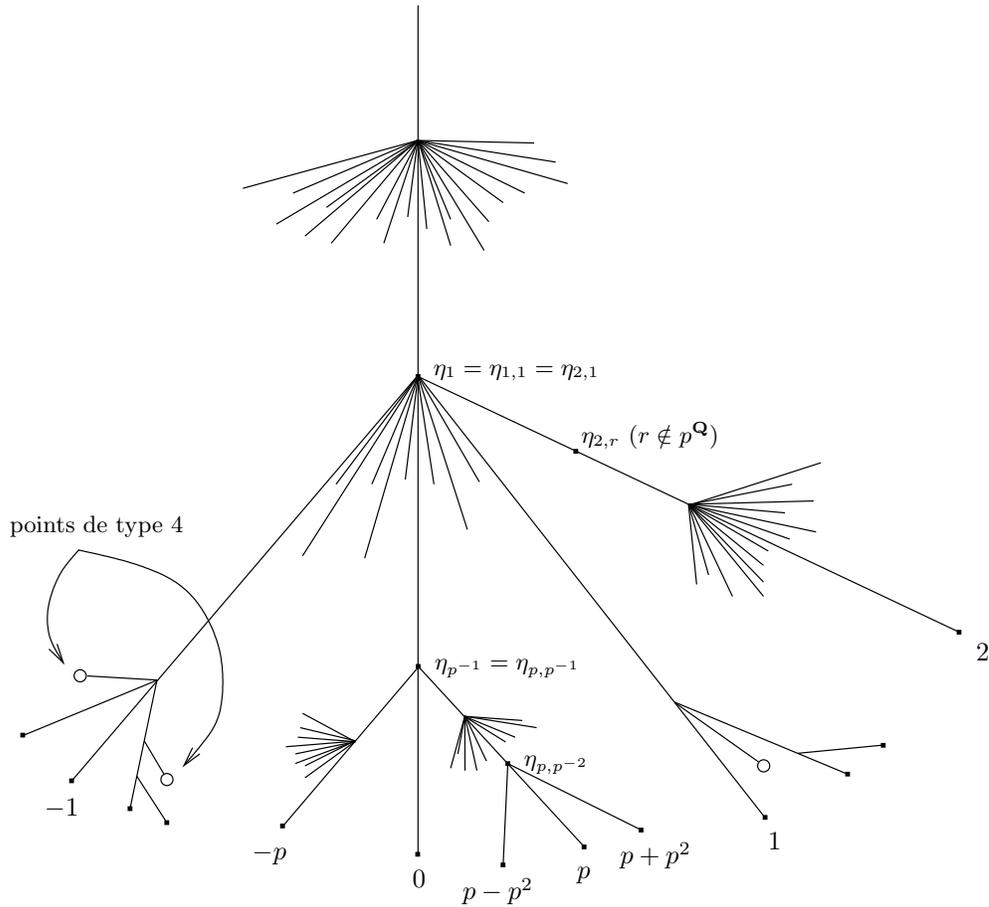
\begin{figure}[htb]
\begin{center}
\input{A1Cppoints.pstex_t}
\caption{Droite analytique sur le corps $\C_{p}$ muni de la valeur absolue $|.|_{p}$.}
\end{center}
\end{figure}

\bigskip

Soient $\alpha\in k$ et $r>0$. L'application
$${\renewcommand{\arraystretch}{1.5}\begin{array}{ccc}
k[T] & \to & \R_{+}\\
\disp \sum_{n\in\N} c_{n}\, (T-\alpha)^n & \mapsto & \disp \max_{n\in\N} (|c_{n}|\, r^n)
\end{array}}$$
d\'efinit encore une semi-norme multiplicative sur $k[T]$ born\'ee sur $k$. Seul le caract\`ere multiplicatif n'est pas imm\'ediat. Il provient en fait de l'in\'egalit\'e ultram\'etrique. Nous noterons $\eta_{\alpha,r}$ le point de la droite $\E{1}{k}$ correspondant.\newcounter{etaalpharac}\setcounter{etaalpharac}{\thepage} Il est remarquable que, contrairement \`a ce que notre notation peut laisser croire, le point $\eta_{\alpha,r}$ ne d\'epend que du disque de centre~$\alpha$ et de rayon~$r$. En particulier, pour $\beta\in k$, nous avons
$$\eta_{\alpha,r} = \eta_{\beta,r} \textrm{ si } |\alpha-\beta|\le r.$$
Les diff\'erents points $\eta_{\alpha,r}$ se comportent diff\'eremment selon que le nombre r\'eel~$r$ appartient ou non au groupe~$|k^*|$. Lorsque~$r$ appartient \`a~$|k^*|$, le point $\eta_{\alpha,r}$ est dit {\bf de type~2}.\index{Point!de type 2} Nous avons alors
$$\widetilde{\Hs(\eta_{\alpha,r})}\simeq \tilde{k}(T) \textrm{ et } |\Hs(\eta_{\alpha,r})^*|=|k^*|.$$
Lorsque~$r$ n'appartient pas \`a~$|k^*|$, le point $\eta_{\alpha,r}$ est dit {\bf de type 3}\index{Point!de type 3}. Nous avons alors
$$\widetilde{\Hs(\eta_{\alpha,r})} = \tilde{k} \textrm{ et le groupe } |\Hs(\eta_{\alpha,r})^*| \textrm{ est engendr\'e par } |k^*| \textrm{ et } r.$$
Lorsque $a=0$, nous noterons simplement $\eta_{r}=\eta_{\alpha,r}$.
\newcounter{etarac}\setcounter{etarac}{\thepage}

\bigskip

Il nous reste un type de points \`a d\'ecrire. Soient~$I$ un ensemble ordonn\'e, \mbox{$\alphab=(\alpha_{i})_{i\in I}$} une famille d'\'el\'ements de~$k$ et $\br=(r_{i})_{i\in I}$ une famille de nombres r\'eels strictement positifs qui v\'erifient les propri\'et\'es suivantes :
\begin{enumerate}[\it i)]
\item $\forall i \le j$, $\overline{D}(\alpha_{i},r_{i}) \subset \overline{D}(\alpha_{j},r_{j})$ ;
\item $\disp \bigcap_{i\in I} \overline{D}(\alpha_{i},r_{i}) = \emptyset$.
\end{enumerate}
De telles familles existent lorsque le corps~$k$ n'est pas maximalement complet\index{Corps!maximalement complet} (\emph{cf.} \cite{Lazard}, d\'efinition 5.2). C'est le cas du corps~$\C_{p}$, pour tout nombre premier~$p$. Remarquons que de telles familles v\'erifient 
$$\inf_{i\in I} (r_{i})>0,$$
sinon le caract\`ere complet du corps~$k$ imposerait \`a l'intersection des disques de contenir un point. L'application
$${\renewcommand{\arraystretch}{1}\begin{array}{ccc}
k[T] & \to & \R_{+}\\
P(T) & \mapsto & \disp \inf_{i\in I} (|P(\eta_{\alpha_{i},r_{i}})|)
\end{array}}$$
d\'efinit une semi-norme multiplicative sur $k[T]$ born\'ee sur $k$. Nous noterons $\eta_{\alphab,\br}$ le point de la droite $\E{1}{k}$ correspondant. Un tel point est dit {\bf de type 4}\index{Point!de type 4}. Le corps r\'esiduel compl\'et\'e en ce point est une extension imm\'ediate\index{Extension immediate@Extension imm\'ediate} du corps~$k$ : il v\'erifie
$$\widetilde{\Hs(\eta_{\alphab,\br})} = \tilde{k} \textrm{ et } |\Hs(\eta_{\alphab,\br})^*|=|k^*|.$$
\newcounter{etaalphabbr}\setcounter{etaalphabbr}{\thepage}

\bigskip

Pour terminer, revenons au cas d'un corps $k$ ultram\'etrique complet quelconque et donc plus n\'ecessairement al\-g\'e\-bri\-que\-ment clos. Consid\'erons le morphisme de changement de base
$$\varphi : \E{1}{\hat{\bar{k}}} \to \E{1}{k}.$$
C'est un morphisme surjectif. Nous dirons qu'un point~$x$ de la droite analytique~$\E{1}{k}$ est de type~$i$\index{Point!de type 1}\index{Point!de type 2}\index{Point!de type 3}\index{Point!de type 4}, pour $i\in\cn{1}{4}$, si l'un des ses ant\'ec\'edents par le morphisme~$\varphi$ est de type~$i$ (et c'est alors le cas pour tous). En outre, pour tous \'el\'ements~$\alpha$ de~$k$ et~$r$ de~$\R_{+}$, nous noterons identiquement les points $\alpha$, $\eta_{\alpha,r}$ et~$\eta_{r}$ de $\E{1}{\hat{\bar{k}}}$ et leur image par le morphisme~$\varphi$. 
\newcounter{alphapac}\setcounter{alphapac}{\thepage}
\newcounter{etaalpharpac}\setcounter{etaalpharpac}{\thepage}
\newcounter{etarpac}\setcounter{etarpac}{\thepage}
De nouveau, nous appellerons {\bf point de Gau{\ss}} le point~$\eta_{1}$.\index{Point!de Gau{\ss}} 

Soit $P(T)$ un polyn\^ome irr\'eductible de $k[T]$. Notons $\alpha_{1},\ldots,\alpha_{d}$, avec $d\in\N^*$, ses racines dans $\bar{k}$. L'application 
$$\begin{array}{ccc}
k[T] & \to & \R_{+}\\
Q(T) & \mapsto & 
\left\{
\begin{array}{cl}
0 & \textrm{si } P\,|\,Q\\
1 & \textrm{sinon}
\end{array}
\right.
\end{array}$$
est une semi-norme multiplicative sur $k[T]$, born\'ee sur $k$. Nous noterons $\eta_{P,0}$ le point de la droite $\E{1}{k}$ correspondant. Nous avons
$$\varphi^{-1}(\eta_{P,0}) = \{\alpha_{1},\ldots,\alpha_{d}\}$$
\newcounter{etaPzeropac}\setcounter{etaPzeropac}{\thepage}
et
$$\Hs(\eta_{P,0}) = k[T]/(P(T)).$$
En particulier, le point~$\eta_{P,0}$ est un point de type~$1$ \index{Point!de type 1|(} et un point rigide \index{Point!rigide|(} de la droite analytique~$\E{1}{k}$. R\'eciproquement, si le corps~$k$ est parfait, le th\'eor\`eme de l'\'el\'ement primitif assure que tout point rigide de cette droite peut s'\'ecrire sous la forme~$\eta_{Q,0}$, o\`u~$Q$ est un polyn\^ome irr\'eductible \`a coefficients dans~$k$.


Les points rigides sont des points de type~$1$ de la droite~$\E{1}{k}$, mais la r\'eciproque n'est en g\'en\'eral pas valable, m\^eme dans le cas des corps parfaits. Consid\'erons, par exemple, le corps~$\Q_{p}$ muni de la valeur absolue $p$-adique usuelle~$|.|_{p}$. Cette valeur absolue se prolonge de fa\c{c}on unique en une valeur absolue sur~$\C_{p}$, que nous noterons identiquement. Soit~$\alpha$ un point de~$\C_{p}$ qui n'est pas alg\'ebrique sur~$\Q_{p}$. L'application
$$\begin{array}{ccc}
\Q_{p}[T] & \to & \R_{+}\\
Q(T) & \mapsto & |Q(\alpha)|_{p}
\end{array}$$
d\'efinit un point de type~$1$ de la droite~$\E{1}{\Q_{p}}$ qui n'est pas un point rigide. En effet, le corps r\'esiduel compl\'et\'e en ce point n'est autre que le corps~$\Q_{p}(\alpha)$, une extension transcendante de~$\Q_{p}$.
\index{Point!de type 1|)}\index{Point!rigide|)} 

\bigskip

La topologie de la droite analytique sur un corps ultram\'etrique complet quelconque est, en g\'en\'eral, assez compliqu\'ee et nous ne la d\'ecrirons pas, mais la figure~\thefigure\ nous semble permettre de se la repr\'esenter assez fid\`element. En particulier, les segments que l'on y voit trac\'es sont hom\'eomorphes \`a des segments. Il faut cependant \^etre prudents en ce qui concerne les voisinages des points de type~$2$, autrement dit, les points de branchement. Soient~$x$ un tel point et~$C_{x}$ l'ensemble des composantes connexes du compl\'ementaire du point~$x$ dans la droite~$\E{1}{k}$ (cet ensemble est naturellement en bijection avec la droite projective sur le corps~$\widetilde{\Hs(x)}$). Alors, pour tout voisinage~$V$ du point~$x$, il n'existe qu'un nombre fini d'\'el\'ements de~$C_{x}$ qui ne soient pas enti\`erement contenus dans~$V$.

\index{Droite affine analytique!sur un corps ultrametrique@sur un corps ultram\'etrique|)}

\subsection{Faisceau structural}

Pour parvenir \`a faire de la g\'eom\'etrie sur les espaces analytiques au sens pr\'ec\'edent, nous devons en faire des espaces localement annel\'es. Nous suivrons la construction d\'evelopp\'ee par V.~Berkovich au num\'ero 1.5 de l'ouvrage~\cite{rouge}. Soient $(\As,\|.\|)$ un anneau de Banach et~$n$ un entier positif. Nous nous restreindrons \`a certains types de normes.

\begin{defi}\label{defsnsp}\index{Norme!spectrale}\index{Semi-norme!spectrale}
On appelle {\bf semi-norme spectrale} sur l'anneau de Banach $(\As,\|.\|)$ la semi-norme d\'efinie par 
$$\forall f\in\As,\, \|f\|_{sp} = \max_{x\in\Ms(\As)} (|f(x)|) = \inf_{k\in\N^*} \left(\left\|f^k\right\|^{\frac{1}{k}}\right).$$
\end{defi}

\begin{rem}
Les deux derni\`eres quantit\'es sont \'egales en vertu de~\cite{rouge}, th\'eor\`eme 1.3.1.
\end{rem}

\begin{defi}\label{defuniforme}\index{Norme!uniforme}\index{Anneau de Banach!uniforme}
Nous dirons que la {\bf norme} $\|.\|$ est {\bf uniforme} si elle est \'equi\-va\-lente \`a la semi-norme spectrale, c'est-\`a-dire s'il existe deux constantes $C_{-}>0$ et $C_{+}>0$ telles que
$$\forall f\in\As,\, C_{-}\, \|f\|_{sp} \le \|f\| \le C_{+}\, \|f\|_{sp}.$$
Dans ce cas, nous dirons que {\bf l'anneau de Banach} $(\As,\|.\|)$ est {\bf uniforme}.
\end{defi}

{\it Dans la suite de ce texte, nous supposerons toujours que la norme $\|.\|$ est uniforme.} Cela impose en particulier \`a la semi-norme spectrale d'\^etre une norme et donc \`a l'anneau~$\As$ d'\^etre r\'eduit. Nous disposons \'egalement d'un hom\'eomorphisme
$$\Ms(\As,\|.\|) \xrightarrow[]{\sim} \Ms(\As,\|.\|_{sp})$$
induit par l'application identit\'e.

\bigskip

D\'efinissons, \`a pr\'esent, le pr\'efaisceau~$\Ks$ des fractions rationnelles sans p\^oles sur~$\E{n}{\As}$ de la fa\c{c}on suivante : pour tout ouvert~$U$ de~$\E{n}{\As}$, l'anneau~$\Ks(U)$ est le localis\'e de $\As[T_{1},\ldots,T_{n}]$ par l'ensemble de ses \'el\'ements qui ne s'annulent en aucun point de $U$. Exprimons cette d\'efinition \`a l'aide de notations math\'ematiques. 

\begin{defi}\label{defKV}
\index{Fractions rationnelles sans poles@Fractions rationnelles sans p\^oles}
\newcounter{Ks}\setcounter{Ks}{\thepage}
Pour tout ouvert~$U$ de l'espace~$\E{n}{\As}$, posons 
$$S_{U} = \left\{ P\in\As[T_{1},\ldots,T_{n}]\, |\, \forall x\in U,\, P(x) \ne 0 \right\}.$$
Nous d\'efinissons {\bf le pr\'efaisceau~$\Ks$ des fractions rationnelles sans p\^oles sur l'espace~$\E{n}{\As}$} comme le foncteur contravariant qui \`a tout ouvert~$U$ de~$\E{n}{\As}$ associ\'e l'anneau
$$\Ks(U) = S_{U}^{-1} \As[T_{1},\ldots,T_{n}].$$

\end{defi}


Nous allons maintenant d\'efinir les fonctions analytiques comme les fonctions qui sont localement limites uniformes de fractions rationnelles sans p\^oles. 

\begin{defi}\index{Faisceau structural}\index{Fonctions analytiques|see{Faisceau structural}}\index{Fonctions holomorphes|see{Faisceau structural}}
\index{Faisceau!structural|see{Faisceau structural}}
\newcounter{Os}\setcounter{Os}{\thepage}
Nous d\'efinissons {\bf le faisceau structural~$\Os$ sur l'espace~$\E{n}{\As}$}, que nous appellerons encore {\bf faisceau des fonctions analytiques sur l'espace~$\E{n}{\As}$}, comme le foncteur contravariant qui \`a tout ouvert~$U$ de~$\E{n}{\As}$ associe l'anneau~$\Os(U)$ constitu\'e de l'ensemble des applications 
$$f : U\to \bigsqcup_{x\in U} \Hs(x)$$
telles que, pour tout \'el\'ement~$x$ de~$U$, on ait
$$f(x)\in\Hs(x)$$
et qui v\'erifient la condition suivante : pour tout \'el\'ement~$x$ de~$U$, il existe un voisinage ouvert~$V$ de~$x$ dans~$U$ et une suite $(R_{i})_{i\in\N}$ d'\'el\'ements de~$\Ks(V)$ telle que, quel que soit $\eps>0$, il existe un entier positif~$j$ pour lequel on ait
$$\forall i\ge j,\, \forall y\in V,\, |f(y) - R_{i}(y)| \le \eps.$$ 
\end{defi}

\begin{rem}
Cette d\'efinition locale assure que~$\Os$ est bien un faisceau d'anneaux sur~$\E{n}{\As}$. On v\'erifie qu'en tout point~$x$ de l'espace~$\E{n}{\As}$, le germe~$\Os_{x}$ est un anneau local dont l'id\'eal maximal est l'ensemble des germes de fonctions qui s'annulent au point~$x$.
\end{rem}

\begin{defi}
Soit~$x$ un point de l'espace~$\E{n}{\As}$. Nous noterons~$\m_{x}$ l'id\'eal maximal de l'anneau local~$\Os_{x}$. Nous appellerons {\bf corps r\'esiduel} du point~$x$ le corps
$$\kappa(x)=\Os_{x}/\m_{x}.$$ 
\newcounter{mx}\setcounter{mx}{\thepage}
\newcounter{kappax}\setcounter{kappax}{\thepage}
\end{defi}

\begin{rem}\index{Espace affine analytique!sur un corps archimedien@sur un corps archim\'edien}
Si l'anneau de Banach consid\'er\'e est l'anneau $\C$ muni de la valeur absolue usuelle, nous retrouvons la notion habituelle de fonction holomorphe. En effet, toutes les fractions rationnelles sans p\^oles sur un ouvert de $\C^n$ sont holomorphes sur cet ouvert et il est bien connu qu'une limite uniforme de fonctions holomorphes reste holomorphe.

R\'eciproquement, toute fonction holomorphe sur un ouvert $U$ de $\C^n$ est localement limite uniforme de polyn\^omes. Il suffit, par exemple, de recouvrir l'ouvert~$U$ par des disques ouverts dont l'adh\'erence est contenue dans $U$.
\end{rem}

Le r\'esultat qui suit justifie le fait que nous ayons choisi de munir l'anneau~$\As$ d'une norme uniforme.

\begin{lem}
Le morphisme d'anneaux naturel
$$\As[T_{1},\ldots,T_{n}] \to \Os(\E{n}{\As})$$
est injectif.
\end{lem}
\begin{proof}
Soit~$P$ un \'el\'ement de~$\As[T_{1},\ldots,T_{n}] $ dont l'image dans~$\Os(\E{n}{\As})$ est nulle. Cela signifie qu'en tout point~$x$ de l'espace~$\E{n}{\As}$, nous avons
$$P(T_{1},\ldots,T_{n})(x) = 0.$$
Il existe une famille presque nulle $(a_{k_{1},\ldots,k_{n}})_{(k_{1},\ldots,k_{n})\in\N^n}$ d'\'el\'ements de~$\As$ telle que l'on ait 
$$P(T_{1},\ldots,T_{n})= \sum_{(k_{1},\ldots,k_{n})\in\N^n} a_{k_{1},\ldots,k_{n}}\, T_{1}^{k_{1}}\cdots T_{n}^{k_{n}} \textrm{ dans } \As[T_{1},\ldots,T_{n}].$$

Soit~$b$ un point de~$\Ms(\As)$. Si le polyn\^ome
$$P_{b}(T_{1},\ldots,T_{n})= \sum_{(k_{1},\ldots,k_{n})\in\N^n} a_{k_{1},\ldots,k_{n}}(b)\, T_{1}^{k_{1}}\cdots T_{n}^{k_{n}} \textrm{ de } \Hs(b)[T_{1},\ldots,T_{n}]$$
n'est pas nul, il existe une extension finie~$L_{b}$ de~$\Hs(b)$ et un \'el\'ement~$\alpha_{b}$ de~$L_{b}^n$ tel que l'on ait
$$P_{b}(\alpha_{b})\ne 0 \textrm{ dans } L_{b}.$$
L'\'el\'ement~$\alpha_{b}$ de~$L_{b}^n$ d\'efinit alors un point (rigide)~$\alpha'_{b}$ de l'espace~$\Hs(b)[T_{1},\ldots,T_{n}]$ en lequel nous avons
$$P(T_{1},\ldots,T_{n})(\alpha'_{b}) = P_{b}(T_{1},\ldots,T_{n})(\alpha'_{b}) \ne 0.$$
C'est impossible. 

Soit $(k_{1},\ldots,k_{n})$ un \'el\'ement de~$\N^n$. Nous avons montr\'e que, pour tout point~$b$ de~$\Ms(\As)$, nous avons
$$a_{k_{1},\ldots,k_{n}}(b)=0.$$
On en d\'eduit que
$$\|a_{k_{1},\ldots,k_{n}}\|_{sp} = 0$$
et donc que
$$\|a_{k_{1},\ldots,k_{n}}\| = 0,$$
puisque la semi-norme~$\|.\|_{sp}$ et la norme~$\|.\|$ sont \'equivalentes. Par cons\'equent, nous avons
$a_{k_{1},\ldots,k_{n}}=0$ dans~$\As$.
On en d\'eduit que le polyn\^ome~$P$ est nul.
\end{proof}



\begin{rem}
L'application identit\'e de $(\As,\|.\|)$ vers $(\As,\|.\|_{sp})$ induit un isomorphisme d'espaces annel\'es
$$\E{n}{\As,\|.\|} \xrightarrow[]{\sim} \E{n}{\As,\|.\|_{sp}}.$$
Pour de nombreuses questions, nous pourrons donc supposer que la norme~$\|.\|$ {\it est} la norme spectrale.
\end{rem}

Nous disposons, \`a pr\'esent, d'une notion de fonction analytique sur les ouverts de l'espace $\E{n}{\As}$. Nous pouvons en d\'eduire une d\'efinition g\'en\'erale d'espace analytique. Nous la donnons ci-dessous dans un souci d'exhaustivit\'e, mais ne l'utiliserons pas. Dans le cas complexe, un espace est dit analytique s'il est localement isomorphe \`a un ferm\'e analytique d'un ouvert d'un espace affine. La d\'efinition suivante s'impose donc naturellement.

\begin{defi}[V.~Berkovich]\label{espan}\index{Espace analytique}
On dit qu'un espace localement annel\'e \mbox{$(V,\Os_{V})$} est un {\bf mod\`ele local d'un espace analytique} sur~$\As$ s'il existe un entier positif~$n$, un ouvert~$U$ de~$\E{n}{\As}$ et un faisceau~$\Is$ d'id\'eaux de type fini de~$\Os_{U}$ tels que $(V,\Os_{V})$ soit isomorphe au support du faisceau $\Os_{U}/\Is$, muni du faisceau $\Os_{U}/\Is$.

On appelle {\bf espace analytique} sur $\As$ tout espace localement annel\'e qui est localement isomorphe \`a un mod\`ele local d'un espace analytique sur $\As$.
\end{defi}


\bigskip

\`A titre d'exemple, donnons, sans d\'emonstration, quelques propri\'et\'es des anneaux locaux en les points de la droite analytique sur un corps ultram\'etrique. 

\begin{prop}\label{anneauxlocauxdroite}\index{Droite affine analytique!sur un corps ultrametrique@sur un corps ultram\'etrique!anneaux locaux}
Soit $(k,|.|)$ un corps ultram\'etrique complet. Notons $X=\E{1}{k}$ la droite analytique sur le corps~$k$. 
\begin{enumerate}[\it i)]
\item Soit~$x$ un point rigide \index{Point!rigide}de l'espace~$X$. Alors, l'anneau local~$\Os_{X,x}$ est un anneau de valuation discr\`ete. S'il existe un polyn\^ome~$P$ irr\'eductible \`a coefficients dans~$k$ tel que le point~$x$ soit le point~$\eta_{P,0}$ (c'est toujours le cas si le corps~$k$ est parfait), alors l'id\'eal maximal de~$\Os_{X,x}$ est engendr\'e par~$P$.
\item Soit~$x$ un point de type~$1$\index{Point!de type 1} de l'espace~$X$ qui n'est pas un point rigide. Alors, l'anneau local~$\Os_{X,x}$ est un corps.
\item Soit~$x$ un point de type~$2$, $3$ ou~$4$ \index{Point!de type 2}\index{Point!de type 3}\index{Point!de type 4}de l'espace~$X$. Alors, l'anneau local~$\Os_{X,x}$ est un corps.
\end{enumerate}
\end{prop}

\bigskip

\bigskip

Dans la suite de ce texte, nous consid\'ererons souvent les sections d'un faisceau \index{Faisceau!sections sur une partie quelconque} au-dessus d'une partie qui n'est pas ouverte. Voici quelques rappels sur cette notion. Soit~$Y$ un espace topologique et~$\Fs$ un faisceau d'ensembles sur~$Y$. \`A ce faisceau est associ\'e un espace \'etal\'e $(\tilde{\Fs},p)$, o\`u~$\tilde{\Fs}$ est un espace topologique et $p : \tilde{\Fs} \to Y$ un hom\'eomorphisme local. Pour tout partie~$V$ de~$Y$, notons $\tilde{\Fs}(V)$ l'ensemble des sections continues de l'application~$p$ au-dessus de~$V$. Pour toute partie {\it ouverte}~$U$ de~$Y$, il existe alors une bijection canonique
$$\Fs(U) \xrightarrow[]{\sim} \tilde{\Fs}(U).$$
Pour des pr\'ecisions sur cette construction, on se reportera \`a~\cite{TF}, II, 1.2. 

\begin{defi}\label{defiFV}
Pour toute partie~$V$ de~$Y$, on pose 
$$\Fs(V)=\tilde{\Fs}(V).$$
\end{defi}

Sous certaines conditions, il est possible de d\'ecrire l'ensemble~$\Fs(V)$ directement en termes des ensembles de sections du faisceau~$\Fs$ sur les ouverts de l'espace~$Y$. Citons, \`a ce propos, le corollaire~1 au th\'eor\`eme 3.3.1 du chapitre~II de l'ouvrage~\cite{TF}.

\begin{thm}\label{thFV}
Soient~$V$ une partie de~$Y$ qui poss\`ede un syst\`eme fondamental de voisinages paracompacts. Alors l'application canonique
$$\varinjlim \Fs(U) \to \Fs(V),$$ 
o\`u la limite inductive est prise sur l'ensemble des voisinages ouverts~$U$ de~$V$ dans~$Y$, est bijective.
\end{thm}

Nous n'utiliserons l'ensemble~$\Fs(V)$ que dans les cas o\`u les hypoth\`eses du th\'eor\`eme sont satisfaites. C'est en particulier le cas lorsque
\begin{enumerate}
\item la partie~$V$ est ferm\'ee et l'espace~$Y$ paracompact (par exemple, si~$Y$ est une partie ferm\'ee d'un espace affine analytique au-dessus d'un anneau de Banach, d'apr\`es le th\'eor\`eme \ref{topoAn}) ;
\item la partie~$V$ est quelconque et l'espace~$Y$ est m\'etrisable (par exemple, si~$Y$ est une partie d'un espace affine analytique au-dessus d'un anneau d'entiers de corps de nombres, comme nous le verrons au th\'eor\`eme \ref{metrisable}).
\end{enumerate}

\bigskip

Signalons que cette notation peut malheureusement pr\^eter \`a confusion lorsque l'on consid\`ere un espace analytique au-dessus d'un corps ultram\'etrique complet. Soient $(k,|.|)$ un tel corps et $n$ un entier positif. Notons~$\overline{D}$ le disque unit\'e ferm\'e centr\'e en~$0$ de l'espace~$\E{n}{k}$. L'alg\`ebre~$\Os(\overline{D})$ n'est alors pas l'alg\`ebre de Tate, form\'ee des s\'eries qui convergent sur~$\overline{D}$, mais l'alg\`ebre de Washnitzer, constitu\'ee des s\'eries qui convergent au voisinage de~$\D$ (\emph{cf.}~\cite{GK}, 1.2). 

\bigskip

Ajoutons quelques mots au sujet de la restriction des faisceaux.\index{Faisceau!restriction \`a une partie quelconque} 

\begin{defi}
Pour toute partie~$V$ de~$Y$, on d\'efinit un faisceau d'ensembles~$\Fs_{|V}$ sur~$V$ comme le foncteur contravariant qui \`a tout ouvert~$U$ de~$V$ associe l'ensemble $\Fs(V)$.
\end{defi}

Ces restrictions jouissent de bonnes propri\'et\'es, comme le montre le lemme qui suit.

\begin{lem}
Soient~$V$ et~$W$ deux parties de l'espace~$Y$. Pour toute partie~$U$ de~$V\cap W$, nous avons une bijection
$$\Fs_{|V}(U) \simeq \Fs_{|W}(U).$$
En particulier, pour tout point~$x$ de~$V\cap W$, nous avons une bijection entre les germes
$$\left(\Fs_{|V}\right)_{x} \simeq \left(\Fs_{|W}\right)_{x}.$$
\end{lem}

Pour finir, signalons que les constructions et r\'esultats qui pr\'ecd\`ent restent \'evidemment valables {\it mutatis mutandis} pour les faisceaux \`a valeurs dans n'importe quelle cat\'egorie.

\section{Parties compactes spectralement convexes}\label{pcsc}

Soient~$(\As,\|.\|)$ un anneau de Banach uniforme et~$n$ un entier positif. Introduisons deux nouvelles d\'efinitions. Rappelons (\emph{cf.} d\'efinition \ref{defKV}) que, pour toute partie~$V$ de l'espace analytique~$\E{n}{\As}$, nous d\'efinissons l'anneau~$\Ks(V)$ comme le localis\'e de l'anneau~$\As[T_{1},\ldots,T_{n}]$ par la partie multiplicative form\'ee des \'el\'ements qui ne s'annulent pas au voisinage de~$V$.

\begin{defi}\label{defBV}\newcounter{BsV}\setcounter{BsV}{\thepage}
Soit $V$ une partie compacte de l'espace analytique~$\E{n}{\As}$. Nous notons~$\Bs(V)$ le compl\'et\'e de l'anneau~$\Ks(V)$ pour la norme uniforme~$\|.\|_{V}$ sur~$V$. 
\end{defi}

\begin{rem}
Quel que soient~$P\in\Ks(V)$ et~$k\in\N$, nous avons l'\'egalit\'e 
$$\|P^k\|_{V}=\|P\|_{V}^k$$ 
et cette propri\'et\'e s'\'etend \`a~$\Bs(V)$. On en d\'eduit que la norme~$\|.\|_{V}$ sur~$\Bs(V)$ est la norme spectrale. En particulier, le couple $(\Bs(V),\|.\|_{V})$ est un anneau de Banach uniforme.
\end{rem}


Soit $V$ une partie compacte de l'espace analytique~$\E{n}{\As}$. Le morphisme naturel 
$$f : \As[T_{1},\ldots,T_{n}] \to \Bs(V)$$
est born\'e sur~$\As$. Il induit donc un morphisme entre espaces localement annel\'es 
$$\varphi : \Ms(\Bs(V)) \to \E{n}{\As}.$$ 
Nous allons chercher ici \`a d\'ecrire l'image de ce morphisme et, plus g\'en\'eralement, \`a comprendre ses propri\'et\'es.

\bigskip

Commen\c{c}ons par une propri\'et\'e topologique simple.

\begin{lem}\label{homeoimage}
Le morphisme~$\varphi$ r\'ealise un hom\'eomorphisme sur son image.
\end{lem}
\begin{proof}
Puisque l'espace~$\Ms(\Bs(V))$ est compact, il suffit de montrer que le morphisme~$\varphi$ est injectif. Soient~$x$ et~$y$ deux points distincts de~$\Ms(\Bs(V))$. Notons~$|.|_{x}$ et~$|.|_{y}$ les semi-normes multiplicatives born\'ees sur~$\Bs(V)$ associ\'ees. Par hypoth\`ese, il existe un \'el\'ement~$P$ de~$\Bs(V)$ tel que
$$|P|_{x} \ne |P|_{y}.$$
La densit\'e de $\Ks(V)$ dans $\Bs(V)$ nous permet d'en d\'eduire qu'il existe~\mbox{$Q\in\Ks(V)$} tel que 
$$|Q|_{x} \ne |Q|_{y}.$$
En \'ecrivant $Q$ comme \'el\'ement du localis\'e de $\As[T_{1},\ldots,T_{n}]$, on montre alors qu'il existe un polyn\^ome~$P\in\As[T_{1},\ldots,T_{n}]$ tel que
$$|f(P)|_{x} \ne |f(P)|_{y}.$$
Par cons\'equent, les points $\varphi(x)$ et $\varphi(y)$ de $\E{n}{\As}$ sont distincts.
\end{proof}

En fait, nous disposons m\^eme d'un isomorphisme d'espaces annel\'es si l'on s'autorise \`a restreindre le morphisme \`a la source et au but. 
\begin{lem}\label{isointerieur}
Notons $U$ l'int\'erieur de l'image de $\varphi$ dans $\E{n}{\As}$. Le morphisme  
$$\psi : \varphi^{-1}(U) \to U$$
induit par $\varphi$ est un isomorphisme d'espaces annel\'es.
\end{lem}
\begin{proof}
Soit $x\in\varphi^{-1}(U)$. Notons $y=\psi(x)=\varphi(x)$. Il nous suffit de montrer que le morphisme induit
$$\psi_{x}^* : \Os_{U,y} \to \Os_{\varphi^{-1}(U),x}$$
est un isomorphisme. L'injectivit\'e provient directement du fait que~$\varphi$ est un hom\'eomorphisme.

Montrons que ce morphisme est surjectif. Soit $g\in \Os_{\varphi^{-1}(U),x}$. Notons~$\Ks'$ le pr\'efaisceau des fractions rationnelles sur~$\Ms(\Bs(V))$. Il existe un voisinage compact~$W$ de~$x$ dans~$\varphi^{-1}(U)$ et une suite~$(R_{k})_{k\in\N}$ d'\'el\'ements de~$\Ks'(W)$ qui converge uniform\'ement vers~$g$ sur~$W$. Soit~$k\in\N$. Par d\'efinition de~$\Ks'(W)$, il existe un \'el\'ement~$S_{k}$ de $\Ks(\psi(W))$ tel que
$$\|\psi_{|W}^*(S_{k}) - R_{k}\|_{W} \le \frac{1}{2^k}.$$
La suite~$(S_{k})_{k\in\N}$ \'etant de Cauchy uniforme sur~$\psi(W)$, elle converge vers un \'el\'ement de~$\Bs(\psi(W))$. Son image dans l'anneau local~$\Os_{U,y}$ est envoy\'ee sur~$g$ par~$\psi^*_{x}$. 
\end{proof}

\begin{rem}
Le r\'esultat est, en g\'en\'eral, faux si l'on ne restreint pas le morphisme. Supposons que le compact~$V$ soit r\'eduit \`a un point~$x$. Par d\'efinition, nous avons alors $\Bs(V)=\Hs(x)$. L'hom\'eomorphisme induit par~$\varphi$ est donc
$$\Ms(\Hs(x)) \xrightarrow[]{\sim} \{x\}$$
et le morphisme induit entre les anneaux locaux est
$$\Os_{X,x} \to \Hs(x).$$
Ce n'est pas, en g\'en\'eral, un isomorphisme.
\end{rem}

D\'emontrons, \`a pr\'esent, un premier r\'esultat sur l'image de $\varphi$.

\begin{lem}\label{contientV}
L'image du morphisme $\varphi$ contient le compact $V$.
\end{lem}
\begin{proof}
Soit $x$ un point de $V$. Il lui correspond un caract\`ere 
$$\chi_{x} : \As[T_{1},\ldots,T_{n}] \to \Hs(x).$$
Puisque $x\in V$, un \'el\'ement $P$ de $\As[T_{1},\ldots,T_{n}]$ qui ne s'annule pas au voisinage de $V$ ne s'annule pas en $x$. Son image est donc inversible dans $\Hs(x)$. Par cons\'equent, le morphisme $\chi_{x}$ induit, par localisation, un morphisme
$$\Ks(V) \to \Hs(x).$$
Puisque $x$ appartient \`a $V$, ce morphisme est born\'e. Il induit donc un morphisme entre les compl\'et\'es
$$\Bs(V) \to \Hs(x),$$
ce qu'il fallait d\'emontrer.
\end{proof}

\begin{rem}
La r\'eciproque de ce r\'esultat n'est pas vraie en g\'en\'eral. Montrons-le sur un exemple. Choisissons pour alg\`ebre de Banach $\As$ un corps alg\'ebriquement clos~$k$ que nous munissons de la valeur absolue triviale~$|.|_{0}$. Notons~$D$ le disque ferm\'e de centre~$0$ et de rayon~$1$ de~$X = \E{n}{k}$ :
$$D = \bigcap_{1\le i\le n} \left\{ x\in X\, \big|\, |T_{i}(x)|\le 1\right\}.$$
Consid\'erons la partie compacte $V$ de $X$ d\'efinie par 
$$V = \bigcup_{1\le i\le n} \left\{ x \in D\, \big|\, |T_{i}(x)|=1\right\}.$$
Supposons que $n\ge 2$. Tout polyn\^ome non constant $P$ de $k[T_{1},\ldots,T_{n}]$ s'annule alors sur $V$, puisqu'il s'annule en un point non nul de~$k^n$. Par cons\'equent, nous avons
$$\Ks(V) = k[T_{1},\ldots,T_{n}].$$
La norme uniforme sur la partie $V$ n'est autre que la norme triviale. On en d\'eduit que~$\Bs(V)$ est l'alg\`ebre~$k[T_{1},\ldots,T_{n}]$ munie de la norme triviale, autrement dit l'alg\`ebre~$k\{T_{1},\ldots,T_{n}\}$ munie de la norme de {Gau\ss}. Par cons\'equent, l'image de~$\Ms(\Bs(V))$ dans~$X$ est le disque~$D$ tout entier.
\end{rem}

\bigskip


Dans certains cas, nous pouvons cependant affirmer que l'image du morphisme~$\varphi$ co\"incide avec la partie compacte~$V$.

\begin{defi}\label{defrationnel}\index{Compact!rationnel}\index{Compact!pro-rationnel}
Nous dirons que la partie compacte~$V$ de l'espace affine~$\E{n}{\As}$ est {\bf rationnelle} s'il existe $p\in\N$, des polyn\^omes $P_{1},\ldots,P_{p},Q$ de $\As[T_{1},\ldots,T_{n}]$ ne s'annulant pas simultan\'ement sur~$V$ et des nombres r\'eels \mbox{$r_{1},\ldots,r_{p}>0$} tels que
$$V = \bigcap_{1\le i\le p} \left\{ x\in X\, \big|\, |P_{i}(x)|\le r_{i}\, |Q(x)|\right\}.$$

Une partie compacte~$V$ de l'espace affine~$\E{n}{\As}$ est dite {\bf pro-rationnelle} si elle est intersection de parties compactes rationnelles.
\end{defi}

\begin{rem}
Soient $p\in\N$, $P_{1},\ldots,P_{p}$ des \'el\'ements de $\As[T_{1},\ldots,T_{n}]$ et \mbox{$s_{1},\ldots,s_{p},t_{1},\ldots,t_{p}$} des nombres r\'eels positifs. Alors la partie de $\E{n}{\As}$ d\'efinie par
$$\bigcap_{1\le i\le p} \left\{ x\in X\, \big|\, s_{i} \le |P_{i}(x)| \le t_{i}\right\}$$
est une partie compacte rationnelle de~$\E{n}{\As}$, d\`es qu'elle est compacte. Rappelons que nous avons donn\'e des exemples de parties compactes \`a la proposition \ref{disquecompact} et au corollaire \ref{partiecompacte}. On en d\'eduit ais\'ement que tout point de $\E{n}{\As}$ poss\`ede un syst\`eme fondamental de voisinages constitu\'e de parties compactes rationnelles. 
\end{rem}

\begin{lem}
Si le compact~$V$ est pro-rationnel, alors l'image du morphisme~$\varphi$ est contenue dans~$V$.
\end{lem}
\begin{proof}
Supposons qu'il existe un ensemble $J$ et une famille $(V_{j})_{j\in J}$ de parties compactes rationnelles telles que
$$V = \bigcap_{j\in J} V_{j}.$$
Soit $j\in J$. Il existe un entier $p\in\N$, des polyn\^omes $P_{1},\ldots,P_{p},Q$ de $\As[T_{1},\ldots,T_{n}]$ ne s'annulant pas simultan\'ement sur~$V$ et des nombres r\'eels $r_{1},\ldots,r_{p}>0$ tels que
$$V_{j} = \bigcap_{1\le i\le p} \left\{ x\in X\, \big|\, |P_{i}(x)|\le r_{i}\, |Q(x)|\right\}.$$
Soit~$x$ un point de~$\Ms(\Bs(V))$. Il est associ\'e \`a une semi-norme multiplicative~$|.|_{x}$ born\'ee sur~$\Bs(V)$. Rappelons que nous notons~$f$ le morphisme naturel de~$\As[T_{1},\ldots,T_{n}]$ dans $\Bs(V)$. Le point $y=\varphi(x)$ est alors associ\'e \`a la semi-norme multiplicative born\'ee sur $\As[T_{1},\ldots,T_{n}]$ d\'efinie par $|f(.)|_{x}$. Par hypoth\`ese, le polyn\^ome~$Q$ ne s'annule pas sur~$V_{j}$ et donc sur~$V$. On en d\'eduit que l'\'el\'ement~$f(Q)$ est inversible dans~$\Bs(V)$. Par cons\'equent, nous avons~$|f(Q)|_{x} \ne 0$. Soit $i\in\cn{1}{p}$. Nous avons
$$\frac{|f(P_{i})|_{x}}{|f(Q)|_{x}} = \left| \frac{f(P_{i})}{f(Q)} \right|_{x} \le  \left\| \frac{f(P_{i})}{f(Q)} \right\|_{V}.$$
Or, par d\'efinition de $V_{j}$, quel que soit $z\in V$, nous avons $|P_{i}(z)| \le r_{i}\, |Q(z)|$. On en d\'eduit que
$$ \left\| \frac{f(P_{i})}{f(Q)} \right\|_{V} = \sup_{z\in V} \left(\frac{|P_{i}(z)|}{|Q(z)|}\right) \le \sup_{z\in V_{j}} \left(\frac{|P_{i}(z)|}{|Q(z)|}\right) \le r_{i}.$$
Par cons\'equent, nous avons
$$|f(P_{i})|_{x} \le r_{i}\, |f(Q)|_{x}.$$
Cette in\'egalit\'e \'etant v\'erifi\'ee quel que soit $i\in\cn{1}{p}$, la semi-norme multiplicative~$|f(.)|_{x}$ correspond bien \`a un \'el\'ement de $V_{j}$. 

Finalement, nous avons montr\'e que
$$y \in \bigcap_{j\in J} V_{j} = V.$$
L'image du morphisme $\varphi$ est donc contenue dans $V$. 
\end{proof}

Regroupons dans un m\^eme \'enonc\'e les r\'esultats que nous avons d\'emontr\'es dans le cas des parties compactes pro-rationnelles.

\begin{thm}\label{compactrationnel}\index{Compact!rationnel}\index{Compact!pro-rationnel}
Si le compact~$V$ est pro-rationnel, alors le morphisme
$$\varphi : \Ms(\Bs(V)) \to \E{n}{\As}$$ 
induit par le morphisme naturel
$$\As[T_{1},\ldots,T_{n}] \to \Bs(V)$$
r\'ealise un hom\'eomorphisme de~$\Ms(\Bs(V))$ sur son image, qui est \'egale \`a $V$. En outre, le morphisme  
$$\varphi^{-1}\big(\mathring{V}\big) \to \mathring{V}$$
induit par $\varphi$ est un isomorphisme d'espace annel\'es.
\end{thm}

\bigskip

Afin d'y faire r\'ef\'erence par la suite, nous donnons un nom aux parties compactes qui poss\`edent des propri\'et\'es analogues \`a celles des parties compactes pro-rationnelles.

\begin{defi}\label{defspconvexe}\index{Compact!spectralement convexe}
Nous dirons que la partie compacte~$V$ de l'espace analytique~$\E{n}{\As}$ est {\bf spectralement convexe} si le morphisme naturel
$$\varphi : \Ms(\Bs(V)) \to \E{n}{\As}$$
induit un hom\'eomorphisme entre $\Ms(\Bs(V))$ et~$V$ et si le morphisme induit
$$\varphi^{-1}\big(\mathring{V}\big) \to \mathring{V}$$
est un isomorphisme d'espace annel\'es.
\end{defi}

\begin{rem}\label{remdefspconvexe}
D'apr\`es les lemmes~\ref{homeoimage}, \ref{isointerieur} et~\ref{contientV}, une partie compacte~$V$ est spectralement convexe si, et seulement si, l'image du morphisme~$\varphi$ est contenue dans~$V$. 
\end{rem} 

\`A partir d'une partie compacte spectralement convexe donn\'ee, il est facile d'en construire d'autres, ainsi que le montrent les r\'esultats qui suivent. Introduisons des notations suppl\'ementaires. Soit~$m\in\N$. Le morphisme 
$$\As[T_{1},\ldots,T_{n}] \to \Bs(V)$$ 
induit un morphisme
$$\As[T_{1},\ldots,T_{n},S_{1},\ldots,S_{m}] \to \Bs(V)[S_{1},\ldots,S_{m}]$$ 
et un diagramme commutatif
$$\xymatrix{
\E{m}{\Bs(V)} \ar[r]^\psi \ar[d]_{\pi'} & \E{n+m}{\As} \ar[d]_{\pi''}\\
\Ms(\Bs(V)) \ar[r]^\varphi & \E{n}{\As}
}.$$
Nous noterons~$\Ks'$ et~$\Bs'$ (respectivement~$\Ks''$ et~$\Bs''$) le pr\'efaisceau des fractions rationnelles sur~$\E{m}{\Bs(V)}$ (respectivement~$\E{n+m}{\As}$) et celui que l'on obtient en compl\'etant les anneaux de sections pour la norme uniforme.

\begin{lem}
Soit~$r\in\R_{+}$. Notons~$D''$ la partie compacte de~$\E{n+m}{\As}$ d\'efinie par
$$D'' = \left\{ x\in {\pi''}^{-1}(V)\, \big|\, \forall j\in\cn{1}{m},\, |S_{j}(x)|\le r \right\}.$$
Si le compact~$V$ est spectralement convexe, alors le compact~$D''$ l'est \'egalement.
\end{lem}
\begin{proof}
Supposons que le compact~$V$ est spectralement convexe. D'apr\`es la remarque \ref{remdefspconvexe}, il suffit de montrer que l'image~$Z$ du morphisme naturel
$$\Ms(\Bs''(D'')) \to \E{n+m}{\As}$$
est contenue dans~$D''$. Remarquons, tout d'abord, que, pour tout \'el\'ement~$j$ de~$\cn{1}{m}$, nous avons $\|S_{j}\|_{D''}\le r$. On en d\'eduit que
$$Z \subset  \left\{ x\in\E{n+m}{\As}\, \big|\, \forall j\in\cn{1}{m},\, |S_{j}(x)|\le r \right\}.$$

\medskip

Consid\'erons, maintenant, le morphisme
$$\As[T_{1},\ldots,T_{n}] \to \As[T_{1},\ldots,T_{n},S_{1},\ldots,S_{m}].$$
Pour tout \'el\'ement~$P$ de $\As[T_{1},\ldots,T_{n}$ et tout \'el\'ement~$x$ de~$\E{n+m}{\As}$, nous avons l'\'egalit\'e $|P(x)|=|P(\pi''(x))|$. On en d\'eduit que le morphisme pr\'ec\'edent induit un morphisme
$$\Ks(V) \to \Ks''(W),$$
puis un morphisme born\'e
$$\Bs(V) \to \Bs''(W).$$
Nous obtenons alors le diagramme commutatif suivant :
$$\xymatrix{
\Ms(\Bs''(W)) \ar[r] \ar[d] & \E{n+m}{\As} \ar[d]_{\pi''}\\
\Ms(\Bs(V)) \ar[r]^\varphi & \E{n}{\As}
}.$$
Puisque le compact~$V$ est spectralement convexe, l'image du morphisme~$\varphi$ est contenue dans~$V$. On en d\'eduit que l'image~$Z$ de~$\Ms(\Bs''(W))$ est contenue dans~${\pi''}^{-1}(V)$ et, finalement, dans~$D''$.
\end{proof}

\begin{prop}\label{compactrationnelrelatif}
Si le compact~$V$ est spectralement convexe, alors le morphisme
$$\psi : \E{m}{\Bs(V)} \to \E{n+m}{\As}$$
induit un hom\'eomorphime sur son image, qui est \'egale \`a~${\pi''}^{-1}(V)$. En outre, le morphisme induit au-dessus de~$\mathring{V}$ est un isomorphisme d'espaces annel\'es.
\end{prop}
\begin{proof}
Supposons que le compact~$V$ est spectralement convexe. Soit~$r>0$. Posons
$$D' = \left\{ x\in \E{m}{\Bs(V)}\, \big|\, \forall j\in\cn{1}{s},\, |S_{j}(x)|\le r \right\}$$
et
$$D'' = \left\{ x\in {\pi''}^{-1}(V)\, \big|\, \forall j\in\cn{1}{s},\, |S_{j}(x)|\le r \right\}.$$
Puisque~$V$ est spectralement convexe, le morphisme $D'\to D''$ induit par~$\psi$ est bijectif. En particulier, un \'el\'ement~$P$ de $\As[T_{1},\ldots,T_{n},S_{1},\ldots,S_{m}]$ s'annule sur~$D'$ si, et seulement si, il s'annule sur~$D''$ et satisfait l'\'egalit\'e $\|P\|_{D'}=\|P\|_{D''}$. On d\'eduit de ces propri\'et\'es que le morphisme naturel
$$\Bs''(D'')\to\Bs'(D')$$ 
est un isomorphisme. 

Consid\'erons le diagramme commutatif
$$\xymatrix{
\Ms(\Bs'(D')) \ar[r]^\sim \ar[d]_{\alpha} & \Ms(\Bs''(D'')) \ar[d]_{\beta}\\
\E{m}{\Bs(V)} \ar[r]^\psi& \E{n+m}{\As}
}.$$
Puisque~$D'$ est un compact rationnel, le morphisme~$\alpha$ induit un hom\'eomorphisme sur son image, qui est \'egale \`a~$D'$, et un isomorphisme d'espaces annel\'es sur l'int\'erieur. D'apr\`es le lemme pr\'ec\'edent, le compact~$D''$ est spectralement convexe et le morphisme~$\beta$ induit un hom\'eomorphisme sur son image, qui est \'egale \`a~$D''$, et un isomorphisme d'espaces annel\'es sur l'int\'erieur. On en d\'eduit que le morphisme~$\psi$ induit un hom\'eomorphisme entre les espaces~$D'$ et~$D''$ et un isomorphisme d'espaces annel\'es entre leur int\'erieur. 

On obtient finalement le r\'esultat voulu en remarquant que les espaces~$\E{m}{\Bs(V)}$ et~$\E{n+m}{\As}$ sont obtenus comme r\'eunion des espaces~$D'$ et~$D''$ pour~$r\in\R_{+}$.

\end{proof}



\begin{prop}\label{stabilitespconvexe}
Supposons que le compact~$V$ est spectralement convexe. Soit~$W$ une partie compacte et spectralement convexe de~$\E{m}{\Bs(V)}$. Alors, la partie compacte~$\psi(W)$ de~$\E{n+m}{\As}$ est encore spectralement convexe.
\end{prop}
\begin{proof}
D'apr\`es la proposition pr\'ec\'edente, le morphisme
$$\psi : \E{m}{\Bs(V)} \to \E{n+m}{\As}$$
induit un hom\'eomorphime sur son image, qui est \'egale \`a~${\pi''}^{-1}(V)$. En raisonnant comme dans la d\'emonstration pr\'ec\'edente, on en deduit que le morphisme naturel
$$\Bs''(\psi(W)) \to \Bs'(W)$$
est un isomorphisme. 

Consid\'erons, \`a pr\'esent, le diagramme commutatif
$$\xymatrix{
\Ms(\Bs'(W)) \ar[r]^\sim \ar[d]_{\alpha} & \Ms(\Bs''(\psi(W))) \ar[d]_{\beta}\\
\E{m}{\Bs(V)} \ar[r]^\psi_{\sim}& \E{n+m}{\As}
}.$$
Puisque le compact~$W$ est spectralement convexe, l'image du morphisme~$\alpha$ est \'egale \`a~$W$. On en d\'eduit que l'image du morphisme~$\beta$ est contenue dans~$\psi(W)$. On conclut alors \`a l'aide de la remarque \ref{remdefspconvexe}.
\end{proof}


%

Pour finir, nous montrons l'existence de l'enveloppe spectralement convexe d'une partie compacte. 

\begin{prop}\label{envspconvexe}
Notons~$W$ l'image du morphisme naturel
$$\varphi : \Ms(\Bs(V)) \to \E{n}{\As}.$$
C'est une partie compacte et spectralement convexe de~$\E{n}{\As}$.
\end{prop}
\begin{proof}
D'apr\`es les lemmes \ref{homeoimage} et \ref{contientV}, le morphisme~$\varphi$ r\'ealise un hom\'eomorphisme sur le compact~$W$ et ce dernier contient~$V$. Soit~$f$ un \'el\'ement de $\As[T_{1},\ldots,T_{n}]$ qui ne s'annule pas au voisinage du compact~$V$. Il poss\`ede alors un inverse dans $\Ks(V)$ et donc dans $\Bs(V)$. On en d\'eduit que, pour tout \'el\'ement $y$ de $\Ms(\Bs(V))$, nous avons $f(y)\ne 0$. La fonction~$f$ est donc minor\'ee par une constante strictement positive sur le compact~$W$. Elle ne s'annule donc pas au voisinage de~$W$. On en d\'eduit que le morphisme 
$$\Ks(W)\to\Ks(V)$$
induit par l'inclusion $V\subset W$ est un isomorphisme. Puisque le morphisme~$\varphi$ a pour image~$W$, la norme uniforme sur~$W$ co\"incide avec la norme sur~$\Bs(V)$, qui n'est autre que la norme uniforme sur~$V$. On en d\'eduit que le morphisme naturel
$$\Bs(W)\to\Bs(V)$$
est un isomorphisme. Il en est donc de m\^eme pour le morphisme
$$\Ms(\Bs(W)) \to \Ms(\Bs(V)) \xrightarrow[]{\sim} W.$$
\end{proof}




%

%

\section{Flot}\label{parflot}

Nous consacrons cette partie \`a la d\'emonstration de quelques propri\'et\'es des semi-normes multiplicatives. Nous nous int\'eresserons notamment \`a l'application qui consiste \`a \'elever une semi-norme multiplicative \`a une certaine puissance. 

Commen\c{c}ons par rappeler un r\'esultat classique permettant de d\'emontrer qu'une application est une valeur absolue (\emph{cf.} \cite{algcom56}, VI, \S 6, \no 1, proposition~2).

\begin{prop}\label{vabourbaki}
Soit $k$ un corps. Soit $f$ une application de $k$ dans $\R_{+}$ v\'erifiant les propri\'et\'es suivantes :
\begin{enumerate}[\it i)]
\item $(f(x)=0) \Longleftrightarrow (x=0)$ ;
\item $\forall x,y \in K,\, f(xy)=f(x)f(y)$ ;
\item $\exists A>0,\, \forall x,y \in K,\, f(x+y) \le A\max(f(x),f(y))$ ;
\item $\exists C>0,\, \forall n\in\N^*,\, f(n.1)\le Cn$.
\end{enumerate}
Alors l'application $f$ est une valeur absolue sur $k$.
\end{prop}

\begin{lem}\label{pum}\index{Valeur absolue!inegalite ultrametrique relachee@in\'egalit\'e ultram\'etrique rel\^ach\'ee}
Soit $k$ un corps muni d'une valeur absolue $|.|$. Supposons qu'il existe $\lambda\in \of{[}{0,1}{]}$ tel que, quel que soit $n\in\N$, on ait $$|n.1|\le n^\lambda.$$
Alors, quels que soient les \'el\'ements $x$ et $y$ de $k$, on a 
$$|x+y|\le 2^\lambda \max\{|x|,|y|\}.$$
\end{lem}
\begin{proof}
Soient $x,y\in k$. Soit $r\in\N^*$. On a 
$${\renewcommand{\arraystretch}{2}\begin{array}{rcl}
|x+y|^r &=& |(x+y)^r|\\
&\le& \disp \sum_{i=0}^r |C_{r}^i|\, |x|^i\, |y|^{r-i}\\
&\le & \disp(r+1) \max_{0\le i\le r} ((C_{r}^i)^\lambda\, |x|^{i}\, |y|^{r-i})\\
&\le & \disp(r+1) \left(\max_{0\le i\le r} \left(C_{r}^i\, |x|^{i/\lambda}\, |y|^{(r-i)/\lambda}\right)\right)^\lambda\\   
&\le & \disp(r+1) \left(\sum_{i=0}^r C_{r}^i\, |x|^{i/\lambda}\, |y|^{(r-i)/\lambda} \right)^\lambda\\   
&\le & \disp(r+1) \left(|x|^{1/\lambda} + |y|^{1/\lambda} \right)^{r\lambda}\\   
&\le & \disp(r+1) \left( 2\, \max(|x|, |y|)^{1/\lambda} \right)^{r\lambda}\\
&\le & \disp(r+1)\, 2^{r\lambda}\,   \max(|x|, |y|)^r.
\end{array}}$$
En \'elevant cette in\'egalit\'e \`a la puissance $1/r$ et en faisant tendre $r$ vers l'infini, on obtient le r\'esultat annonc\'e.
\end{proof}

Soient $x$ un point de $\E{n}{\As}$ et $b$ son projet\'e sur $\Ms(\As)$. Le point $b$ est associ\'e \`a une semi-norme multiplicative $|.|_{b}$ sur $\As$. Un calcul \'el\'ementaire montre que l'ensemble
$$\{\eps\in \R_{+}^*\, |\, \forall f\in\As,\, |f|_{b}^\eps \le \|f\| \}$$
est un intervalle. Nous le noterons indiff\'eremment $I_{x}$ ou $I_{b}$.
\newcounter{Ib}\setcounter{Ib}{\thepage}


Soit $\eps\in I_{b}$. Notons $|.|_{x}$ la semi-norme multiplicative sur $\As[T_{1},\ldots,T_{n}]$ associ\'ee au point $x$ de $\E{n}{\As}$. L'application 
$$|.|_{x}^\eps :
\newcounter{vaeps}\setcounter{vaeps}{\thepage}
\begin{array}{ccc}
\As[T_{1},\ldots,T_{n}] & \to & \R_{+}\\
P & \mapsto & |P|_{x}^\eps
\end{array}$$ 
est multiplicative, envoie $0$ sur $0$ et $1$ sur $1$. 

Montrons, \`a pr\'esent, que c'est une semi-norme. Consid\'erons le corps r\'esiduel compl\'et\'e $(\Hs(x),|.|)$ du point $x$. Quel que soient $f,g\in \Hs(x)$, nous avons
$$|f+g| \le |f| + |g| \le 2 \max(|f|,|g|)$$
et donc
$$|f+g|^\eps \le 2^\eps \max(|f|^\eps,|g|^\eps).$$
En outre, quel que soit $n\in\N$, nous avons
$$|n|^\eps = |n|_{x}^\eps = |n|_{b}^\eps \le \|n\| \le n.$$
D'apr\`es la proposition \ref{vabourbaki}, l'application $|.|^\eps$ est une valeur absolue sur le corps~$\Hs(x)$. On en d\'eduit que l'application $|.|_{x}^\eps$ est une semi-norme multiplicative sur $\As[T_{1},\ldots,T_{n}]$. Elle est born\'ee sur $\As$, par d\'efinition de $I_{b}$, et d\'efinit donc un point de $\E{n}{\As}$. Nous le noterons $x^\eps$. Remarquons que les corps $\Hs(x)$ et $\Hs(x^\eps)$ sont canoniquement isomorphes. Seule la valeur absolue change.

Nous avons volontairement exclu la valeur~$0$ de notre d\'efinition de~$I_{b}$. Il est cependant possible de d\'efinir \'egalement le point~$x^0$, comme nous le montrons ici. Pour cela, il nous faut supposer que l'intervalle~$I_{b}$ a pour borne inf\'erieure~$0$. L'application 
$$|.|_{x}^0 : 
\newcounter{vazero}\setcounter{vazero}{\thepage}
\begin{array}{ccc}
\As[T_{1},\ldots,T_{n}] & \to & \R_{+}\\
P & \mapsto & 
\left\{ \begin{array}{cl}
0 & \textrm{si } |P(x)|=0\\
1 & \textrm{si } |P(x)|\ne 0
\end{array}\right.
\end{array}$$ 
est multiplicative, envoie $0$ sur $0$ et $1$ sur $1$. Le m\^eme raisonnement que pr\'e\-c\'e\-dem\-ment montre que c'est une semi-norme multiplicative sur $\As[T_{1},\ldots,T_{n}]$ qui est born\'ee sur~$\As$. Nous noterons~$x^0$ le point de l'espace $\E{n}{\As}$ qui lui est associ\'e. Contrairement au cas pr\'ec\'edent, les corps~$\Hs(x)$ et~$\Hs(x^0)$ ne sont, en g\'en\'eral, pas isomorphes.

\bigskip

Dans la suite de cette partie, nous noterons $X=\E{n}{\As}$.

\begin{defi}\index{Flot}
D\'efinissons une partie de $X \times \R_{+}^*$ par
$$D = \left\{ (x,\eps),\, x\in X,\, \eps\in I_{x}\right\}.$$
Nous appellerons {\bf flot} l'application 
$$\begin{array}{ccc}
D & \to & X\\
(x,\eps) & \mapsto & x^\eps
\end{array}.$$
\end{defi}


\begin{prop}
Le flot est une application continue. 
\end{prop}
\begin{proof}
Rappelons que la topologie de $X=\E{n}{\As}$ est, par d\'efinition, la topologie la plus grossi\`ere qui rend continues les applications de la forme
$$\begin{array}{ccc}
X & \to & \R_{+}\\
x & \mapsto & |P(x)|
\end{array},$$
avec $P\in A[T_{1},\ldots,T_{n}]$. Pour montrer que le flot est continu, il suffit donc de montrer que, quel que soit $P\in A[T_{1},\ldots,T_{n}]$, l'application compos\'ee
$$\begin{array}{ccc}
D & \to & \R_{+}\\
(x,\eps) & \mapsto & |P(x^\eps)| = |P(x)|^\eps
\end{array}$$
est continue. Cette propri\'et\'e est bien v\'erifi\'ee car l'application pr\'ec\'edente est obtenue en composant deux applications continues : l'application d'\'evaluation de $P$ et l'\'el\'evation \`a la puissance $\eps$.
\end{proof}

Le flot peut parfois se prolonger \`a une partie de~$X\times\R_{+}$, mais il n'est alors, en g\'en\'eral, plus continu. Nous disposons cependant du r\'esultat, plus faible, suivant.

\begin{lem}\label{flot0}
Soit~$x$ un point de~$X$ tel que l'intervalle~$I_{x}$ ait pour borne inf\'erieure~$0$. Alors l'application
$$\begin{array}{ccc}
I_{x}\cup\{0\} & \to & X\\
\eps & \mapsto & x^\eps
\end{array}$$
est continue.
\end{lem}
\begin{proof}
Par d\'efinition de la topologie de~$X$, il suffit de montrer que, quel que soit~$P\in\As[T_{1},\ldots,T_{n}]$, l'application
$$\begin{array}{ccc}
I_{x}\cup\{0\} & \to & \R_{+}\\
\eps & \mapsto & |P(x^\eps)| = |P(x)|^\eps
\end{array}$$
est continue. Ce r\'esultat est imm\'ediat.
\end{proof}




En pratique, il est plus facile d'utiliser le flot en se restreignant \`a certaines parties de l'espace~$X$. Introduisons des notations adapt\'ees. Soit~$Y$ une partie ouverte de~$X$. Posons
$$D_{Y} = \{(z,\lambda)\in D\, |\, z\in Y,\, z^\lambda\in Y\}.$$
\newcounter{DY}\setcounter{DY}{\thepage}
Soit~$x$ un point de~$Y$. Nous notons
$$I_{Y}(x) = \left\{\eps\in I_{x}\, |\, x^\eps\in Y\right\},$$
\newcounter{IYx}\setcounter{IYx}{\thepage}
$$T_{Y}(x) = \left\{x^\eps,\, \eps \in I_{Y}(x)\right\}$$
\newcounter{TYx}\setcounter{TYx}{\thepage}
et
$$D_{Y}(x) = \left\{ (z,\lambda),\, z\in T_{Y}(x),\, \lambda\in I_{Y}(z) \right\}.$$
\newcounter{DYx}\setcounter{DYx}{\thepage}


\begin{defi}\label{voisflot}\index{Voisinages flottants}\index{Flot!voisinages flottants|see{Voisinages flottants}}
Nous dirons que le point~$x$ de~$Y$ a des {\bf voisinages flottants} dans~$Y$ si le flot est une application ouverte en chaque point de $D_{Y}(x)$. 
\end{defi}
\begin{rem}
\begin{enumerate}[a)]
\item Cette d\'efinition ne d\'epend que de la partie $T_{Y}(x)$ et pas du point $x$ lui-m\^eme.
\item Pour tout point~$p$ de~$D_{Y}$, il est \'equivalent de demander que le flot soit ouvert au point~$p$ ou que sa restriction \`a~$D_{Y}$ soit ouverte au point~$p$.
\end{enumerate}
\end{rem}

Lorsque le flot est d\'efini sur une partie suffisamment grande, par exemple lorsque la partie~$D_{Y}$ est un voisinage de~$D_{Y}(x)$ dans~$Y\times\R_{+}^*$, tous les points ont des voisinages flottants. Le lemme qui suit pr\'ecise cet \'enonc\'e. Nous n'avons donc introduit cette notion que pour prendre en compte les effets de bord qui peuvent appara\^{\i}tre.
\begin{lem}\label{criterevoisflot}
Supposons que, quel que soit~$(z,\lambda)\in D_{Y}(x)$, il existe un voisinage~$U$ de~$z$ dans~$Y$ tel que 
$$U \times \{\lambda\} \subset D_{Y}.$$
Alors, le point~$x$ a des voisinages flottants dans~$Y$.
\end{lem}
\begin{proof}
Soit $(z,\lambda) \in D_{Y}(x)$. Puisque $D_{Y}(z) = D_{Y}(x)$, nous pouvons supposer que $z=x$. Soit~$U$ un voisinage du point~$x$ dans~$Y$.
Quitte \`a restreindre~$U$, nous pouvons supposer qu'il est de la forme
$$U = \bigcap_{1\le i\le r} \left\{ z\in Y\, \big|\, \alpha_{i} < |f_{i}(z)|< \beta_{i} \right\},$$
avec $r\in\N$, $f_{1},\ldots,f_{r} \in \As[T_{1},\ldots,T_{n}]$, $\alpha_{1},\ldots,\alpha_{r},\beta_{1},\ldots,\beta_{r}\in\R_{+}$.

L'\'el\'ement $(x^\lambda,1/\lambda)$ appartient \`a $D_{Y}(x)$. Par cons\'equent, il existe un voisinage~$V$ de~$x^\lambda$ dans~$Y$ tel que 
$$V \times \{\eps\} \subset D_{Y}.$$
Consid\'erons la partie~$W$ de~$Y$ d\'efinie par 
$$U = \bigcap_{1\le i\le r} \left\{ z\in Y\, \big|\, \alpha_{i}^\eps < |f_{i}(z)|< \beta_{i}^\eps \right\}.$$
C'est une partie ouverte de $Y$ qui contient le point~$x^\lambda$. Par cons\'equent, la partie~$V\cap W$ de~$Y$ est un voisinage du point~$x^\lambda$ dans~$Y$. Or, quel que soit~$y\in V\cap W$, il existe~$z\in U$ tel que~$y=z^\lambda$. On en d\'eduit que le flot est une application ouverte au point~$(x,\lambda)$.  
\end{proof}

\bigskip

Int\'eressons-nous, \`a pr\'esent, aux propri\'et\'es du flot. Nous allons montrer que, sous certaines hypoth\`eses, il suffit de conna\^{\i}tre les fonctions au voisinage d'un point~$x$ pour les conna\^{\i}tre au voisinage de toute la trajectoire~$T_{Y}(x)$.

\begin{lem}\label{leminj}
Supposons que le point $x$ de $Y$ a des voisinages flottants dans~$Y$. Soit~$U$ un voisinage ouvert de~$x$ dans~$Y$. Soit $(R_{n})_{n\in\N}$ une suite de~$\Ks(U)$ qui converge uniform\'ement sur~$U$. Notons $f \in \Os(U)$ sa limite. Supposons que la fonction~$f$ soit nulle au voisinage du point~$x$. Alors la fonction~$f$ est nulle au voisinage de~$T_{Y}(x) \cap U$.
\end{lem}
\begin{proof}
Il existe un voisinage $U'$ de $x$ dans $U$ tel que, quel que soit $z\in U'$, nous ayons
$$\lim_{n\to +\infty} R_{n}(z) = 0 \textrm{ dans } \Hs(z),$$
c'est-\`a-dire
$$\lim_{n\to +\infty} |R_{n}(z)| = 0.$$
Soit $y\in T_{Y}(x)\cap U$. Il existe $\eps\in I_{Y}(x)$ tel que $y=x^\eps$. Soit~$J$ un voisinage de~$\eps$ dans~$\R_{+}^*$. Alors la partie~$V=D_{Y}\cap (U'\times J)$ est un voisinage de~$(x,\eps)$ dans~$D_{Y}$. Puisque le flot est ouvert au voisinage de $(x,\eps)$, la partie 
$$\left\{ z^\lambda,\, (z,\lambda) \in V \right\}$$
est un voisinage de $y$ dans $Y$. Soit $(z,\lambda)\in V$. Nous avons
$$\lim_{n\to +\infty} |R_{n}(z^\lambda)| = \lim_{n\to +\infty} |R_{n}(z)|^{\lambda} = 0.$$
Par cons\'equent, $f(z^\lambda)=0$ et la fonction $f$ est nulle au voisinage de $y$ dans~$Y$. 
\end{proof}

\begin{prop}\label{flot}
Supposons que le point $x$ de $Y$ a des voisinages flottants dans~$Y$ et que l'ensemble~$I_{Y}(x)$ est un intervalle. Alors le morphisme de restriction 
$$\Os_{Y}(T_{Y}(x)) \to \Os_{Y,x}$$
est un isomorphisme.

Soit $f$ une fonction d\'efinie sur un voisinage de $y$ dans $Y$. Alors la fonction~$f$ poss\`ede un et un seul prolongement au voisinage de $T_{Y}(x)$, que nous noterons encore $f$. Nous avons alors
$$\forall \eps\in I_{Y}(x),\, |f(x^\eps)|=|f(x)|^\eps.$$
En outre, si l'intervalle~$I_{Y}(x)$ a pour borne inf\'erieure~$0$, si le point~$x^0$ appartient \`a~$Y$ et si la fonction~$f$ est \'egalement d\'efinie au point~$x^0$, alors nous avons
$$|f(x^0)|=|f(x)|^0.$$

\end{prop}
\begin{proof}
Commen\c{c}ons par montrer l'injectivit\'e du morphisme. Soit \mbox{$f\in \Os_{Y}(T_{Y}(x))$} telle que $f$ soit nulle au voisinage de $x$. Notons $V$ l'ensemble des points de $T_{Y}(x)$ au voisinage desquels la fonction $f$ est nulle. Il est clair que $V$ est une partie ouverte de $T_{Y}(x)$. Par hypoth\`ese, elle n'est pas vide. Montrons, \`a pr\'esent, que $V$ est une partie ferm\'ee de $T_{Y}(x)$. Soit $y$ un point de $T_{Y}(x)$ adh\'erent \`a $V$. Il existe un voisinage $U$ de $y$ dans $Y$ et une suite $(R_{n})_{n\in\N}$ de~$\Ks(U)$ qui converge uniform\'ement vers $f$ sur $U$. Puisque $y$ est adh\'erent \`a $V$, il existe un point $z$ appartenant \`a $V\cap U$, c'est-\`a-dire un point de $T_{Y}(x) \cap U$ au voisinage duquel la fonction $f$ est nulle. D'apr\`es le lemme \ref{leminj}, la fonction $f$ est nulle au voisinage de $T_{Y}(z) \cap U$ et, en particulier, au voisinage de $y$. On en d\'eduit que la partie $V$ est ferm\'ee. Puisque $I_{Y}(x)$ est un intervalle, l'image $T_{Y}(x)$ de $\{x\}\times I_{Y}(x)$ par le flot est connexe. On en d\'eduit que $V=T_{Y}(x)$ et donc que la fonction $f$ est nulle au voisinage de $T_{Y}(x)$.

Montrons, \`a pr\'esent, que le morphisme est surjectif. Soit $f\in \Os_{Y,x}$. Il existe un voisinage $U$ de $x$ dans $Y$ et une suite $(R_{n})_{n\in\N}$ de $\Ks(U)$ qui converge uniform\'ement vers $f$ sur $U$. Soit $\eps\in I_{Y}(x)$. Nous allons construire une fonction $g_{y}$ au voisinage de $y=x^\eps$. Soit~$J$ un voisinage compact de~$\eps$ dans~$\R_{+}^*$. Alors la partie~$V=D_{Y}\cap (U\times J)$ est un voisinage de~$(x,\eps)$ dans~$D_{Y}$. Puisque le flot est ouvert au voisinage du point~$(x,\eps)$, la partie 
$$\left\{ z^\lambda,\, (z,\lambda) \in V \right\}$$
est un voisinage $V_{y}$ de $y$ dans $Y$. Soit $(z,\lambda)\in V$. Posons 
$$g_{y}(z^\lambda) = f(z) \textrm{ dans } \Hs(z^\lambda).$$ 
Quel que soit $n\in\N$, nous avons encore \mbox{$R_{n} \in \Ks(V_{y})$}. Montrons que la suite $(R_{n})_{n\in\N}$ converge uniform\'ement vers $g_{y}$ sur $V_{y}$. Soit $\eta \in\of{]}{0,1}{]}$. Il existe $N\in\N$ tel que, quels que soient $n\ge N$ et $z\in U$, on ait
$$|R_{n}(z)-f(z)| \le \eta.$$ 
Soient $z\in U'$, $\lambda\in J$ et $n\ge N$. Nous avons alors
$$|R_{n}(z^\lambda)-g_{y}(z^\lambda)| = |R_{n}(z)-f(z)|^{\lambda} \le \eta^{\lambda} \le \eta^{\alpha},$$ 
o\`u $\alpha>0$ d\'esigne la borne inf\'erieure de $J$. Par cons\'equent, la suite $(R_{n})_{n\in\N}$ de~$\Ks(V_{y})$ converge uniform\'ement vers $g_{y}$ sur $V_{y}$. 

Quel que soient $y_{1},y_{2} \in T_{Y}(x)$ et $z\in V_{y_{1}} \cap V_{y_{2}}$, nous avons
$$g_{y_{1}}(z) = \lim_{n\to +\infty} R_{n}(z) = g_{y_{2}}(z) \textrm{ dans } \Hs(z).$$
De m\^eme, quel que soient $y\in T_{Y}(x)$ et $z\in U\cap V_{y}$, nous avons 
$$f(z) = \lim_{n\to +\infty} R_{n}(z) = g_{y_{2}}(z) \textrm{ dans } \Hs(z).$$
Toutes les fonctions que nous avons construites co\"{\i}ncident donc sur les domaines de d\'efinition communs. Par cons\'equent, la fonction $f$ se prolonge bien au voisinage de $T_{Y}(x)$. 

Les r\'esultats sur la valeur absolue des fonctions proviennent directement de la construction du prolongement de $f$ \`a $T_{Y}(x)$. Le r\'esultat pour $x^0$ s'obtient, quant \`a lui, en utilisant le lemme~\ref{flot0} et la continuit\'e de $f$.


\end{proof}

Nous aurons parfois besoin de montrer qu'une fonction d\'efinie au voisinage du point~$x$ se prolonge sur un voisinage connexe de sa trajectoire~$T_{Y}(x)$. Sous certaines hypoth\`eses, le lemme suivant nous permet d'\'etablir un tel r\'esultat. 

\begin{lem}\label{connexeflot}
Supposons que le point~$x$ poss\`ede un syst\`eme fondamental de voisinages connexes (respectivement connexes par arcs) dans~$Y$. Supposons \'egalement que la partie~$D_{Y}$ est un voisinage de~$D_{Y}(x)$ dans~$Y\times\R_{+}^*$. Alors, tout point de~$T_{Y}(x)$ poss\`ede un syst\`eme fondamental de voisinages connexes (respectivement connexes par arcs) dans~$Y$. 
\end{lem}
\begin{proof}
Commen\c{c}ons par remarquer que la seconde hypoth\`ese impose au point~$x$ d'avoir des voisinages flottants dans~$Y$, en vertu du lemme~\ref{criterevoisflot}.

Soient~$y$ un point de~$T_{Y}(x)$ et~$V$ un voisinage de~$y$ dans~$Y$. Il existe~$\eps\in I_{Y}(x)$ tel que~$x^\eps = y$. Notons~$W$ l'image r\'eciproque de~$V$ par le flot. C'est un voisinage du point~$(x,\eps)$ de~$D_{Y}(x)$ dans~$D_{Y}$. Il existe donc un voisinage~$U$ de~$x$ dans~$Y$ et un intervalle ouvert~$J$ contenant~$\eps$ tels que la partie~$U\times J$ soit contenue dans~$W$. Les hypoth\`eses nous permettent de supposer que la partie~$U$ est connexe (respectivement connexe par arcs). Dans ce cas, la partie~$U\times J$ est encore connexe (respectivement connexe par arcs) et il en est de m\^eme pour son image par le flot. Puisque le point~$x$ poss\`ede des voisinages flottants dans~$Y$, cette image est un voisinage du point~$Y$ dans~$V$.
\end{proof}

%% file: valtriv.pstex_t
\begin{picture}(0,0)%
\includegraphics{valtriv.pstex}%
\end{picture}%
\setlength{\unitlength}{3947sp}%
\begingroup\makeatletter\ifx\SetFigFont\undefined%
\gdef\SetFigFont#1#2#3#4#5{%
  \reset@font\fontsize{#1}{#2pt}%
  \fontfamily{#3}\fontseries{#4}\fontshape{#5}%
  \selectfont}%
\fi\endgroup%
\begin{picture}(3300,3798)(4246,-5816)
\put(4261,-5347){\makebox(0,0)[lb]{\smash{{\SetFigFont{12}{14.4}{\rmdefault}{\mddefault}{\updefault}{\color[rgb]{0,0,0}$\eta_{Q,0}$}%
}}}}
\put(6271,-4037){\makebox(0,0)[lb]{\smash{{\SetFigFont{12}{14.4}{\rmdefault}{\mddefault}{\updefault}{\color[rgb]{0,0,0}$\eta_1$}%
}}}}
\put(6211,-3067){\makebox(0,0)[lb]{\smash{{\SetFigFont{12}{14.4}{\rmdefault}{\mddefault}{\updefault}{\color[rgb]{0,0,0}$\eta_s$}%
}}}}
\put(7031,-4657){\makebox(0,0)[lb]{\smash{{\SetFigFont{12}{14.4}{\rmdefault}{\mddefault}{\updefault}{\color[rgb]{0,0,0}$\eta_r$}%
}}}}
\put(4701,-4727){\makebox(0,0)[lb]{\smash{{\SetFigFont{12}{14.4}{\rmdefault}{\mddefault}{\updefault}{\color[rgb]{0,0,0}$\eta_{Q,t}$}%
}}}}
\put(5961,-4747){\makebox(0,0)[lb]{\smash{{\SetFigFont{12}{14.4}{\rmdefault}{\mddefault}{\updefault}{\color[rgb]{0,0,0}$1$}%
}}}}
\put(6671,-5257){\makebox(0,0)[lb]{\smash{{\SetFigFont{12}{14.4}{\rmdefault}{\mddefault}{\updefault}{\color[rgb]{0,0,0}$r$}%
}}}}
\put(7271,-5747){\makebox(0,0)[lb]{\smash{{\SetFigFont{12}{14.4}{\rmdefault}{\mddefault}{\updefault}{\color[rgb]{0,0,0}$0$}%
}}}}
\put(5191,-5467){\makebox(0,0)[lb]{\smash{{\SetFigFont{12}{14.4}{\rmdefault}{\mddefault}{\updefault}{\color[rgb]{0,0,0}$\eta_{P,0}$}%
}}}}
\put(7531,-5267){\makebox(0,0)[lb]{\smash{{\SetFigFont{12}{14.4}{\rmdefault}{\mddefault}{\updefault}{\color[rgb]{0,0,0}$0$}%
}}}}
\put(5381,-4091){\makebox(0,0)[lb]{\smash{{\SetFigFont{12}{14.4}{\rmdefault}{\mddefault}{\updefault}{\color[rgb]{0,0,0}1}%
}}}}
\put(5311,-3126){\makebox(0,0)[lb]{\smash{{\SetFigFont{12}{14.4}{\rmdefault}{\mddefault}{\updefault}{\color[rgb]{0,0,0}$s$}%
}}}}
\put(5121,-2201){\makebox(0,0)[lb]{\smash{{\SetFigFont{12}{14.4}{\rmdefault}{\mddefault}{\updefault}{\color[rgb]{0,0,0}$+\infty$}%
}}}}
\end{picture}%

%% file: A1Cppoints.pstex_t
\begin{picture}(0,0)%
\includegraphics{A1Cppoints.pstex}%
\end{picture}%
\setlength{\unitlength}{3947sp}%
\begingroup\makeatletter\ifx\SetFigFont\undefined%
\gdef\SetFigFont#1#2#3#4#5{%
  \reset@font\fontsize{#1}{#2pt}%
  \fontfamily{#3}\fontseries{#4}\fontshape{#5}%
  \selectfont}%
\fi\endgroup%
\begin{picture}(6100,5696)(3196,-6870)
\put(3211,-4497){\makebox(0,0)[lb]{\smash{{\SetFigFont{9}{10.8}{\rmdefault}{\mddefault}{\updefault}{\color[rgb]{0,0,0}points de type 4}%
}}}}
\put(6056,-6806){\makebox(0,0)[lb]{\smash{{\SetFigFont{10}{12.0}{\rmdefault}{\mddefault}{\updefault}{\color[rgb]{0,0,0}$p-p^2$}%
}}}}
\put(7046,-6581){\makebox(0,0)[lb]{\smash{{\SetFigFont{10}{12.0}{\rmdefault}{\mddefault}{\updefault}{\color[rgb]{0,0,0}$p+p^2$}%
}}}}
\put(6776,-6671){\makebox(0,0)[lb]{\smash{{\SetFigFont{10}{12.0}{\rmdefault}{\mddefault}{\updefault}{\color[rgb]{0,0,0}$p$}%
}}}}
\put(7976,-6491){\makebox(0,0)[lb]{\smash{{\SetFigFont{10}{12.0}{\rmdefault}{\mddefault}{\updefault}{\color[rgb]{0,0,0}$1$}%
}}}}
\put(9281,-5306){\makebox(0,0)[lb]{\smash{{\SetFigFont{10}{12.0}{\rmdefault}{\mddefault}{\updefault}{\color[rgb]{0,0,0}$2$}%
}}}}
\put(3431,-6281){\makebox(0,0)[lb]{\smash{{\SetFigFont{10}{12.0}{\rmdefault}{\mddefault}{\updefault}{\color[rgb]{0,0,0}$-1$}%
}}}}
\put(5741,-6731){\makebox(0,0)[lb]{\smash{{\SetFigFont{10}{12.0}{\rmdefault}{\mddefault}{\updefault}{\color[rgb]{0,0,0}$0$}%
}}}}
\put(4736,-6551){\makebox(0,0)[lb]{\smash{{\SetFigFont{10}{12.0}{\rmdefault}{\mddefault}{\updefault}{\color[rgb]{0,0,0}$-p$}%
}}}}
\put(6441,-5967){\makebox(0,0)[lb]{\smash{{\SetFigFont{9}{10.8}{\rmdefault}{\mddefault}{\updefault}{\color[rgb]{0,0,0}$\eta_{p,p^{-2}}$}%
}}}}
\put(5881,-5357){\makebox(0,0)[lb]{\smash{{\SetFigFont{9}{10.8}{\rmdefault}{\mddefault}{\updefault}{\color[rgb]{0,0,0}$\eta_{p^{-1}}=\eta_{p,p^{-1}}$}%
}}}}
\put(5871,-3507){\makebox(0,0)[lb]{\smash{{\SetFigFont{9}{10.8}{\rmdefault}{\mddefault}{\updefault}{\color[rgb]{0,0,0}$\eta_1=\eta_{1,1}=\eta_{2,1}$}%
}}}}
\put(6801,-3937){\makebox(0,0)[lb]{\smash{{\SetFigFont{9}{10.8}{\rmdefault}{\mddefault}{\updefault}{\color[rgb]{0,0,0}$\eta_{2,r}\ (r\notin p^{\Q})$}%
}}}}
\end{picture}%

%% file: series.tex
\chapter[S\'eries convergentes]{Alg\`ebres de s\'eries convergentes}

Nous consacrons ce chapitre \`a l'\'etude de certains anneaux de s\'eries convergentes \`a coefficients dans un anneau de Banach. Au num\'ero \ref{algglob}, nous nous int\'eressons \`a des alg\`ebres globales, dans la lign\'ee des alg\`ebres de Tate. Nous \'etudions leur spectre analytique et comparons leur norme en tant qu'alg\`ebre de s\'eries \`a la semi-norme uniforme sur leur spectre.

Au num\'ero \ref{algloc}, nous \'etudions certaines limites inductives d'alg\`ebres globales de disques, en un sens que pr\'ecisons. Ce sont des anneaux locaux dont nous montrons qu'il satisfont des th\'eor\`emes de division et de pr\'eparation de Weierstra{\ss}, comme les anneaux locaux des espaces analytiques complexes. Nous en d\'eduisons plusieurs propri\'et\'es, telles la noeth\'erianit\'e ou la r\'egularit\'e. Le num\'ero suivant est consacr\'e \`a l'\'etude de limites inductives d'alg\`ebres globales de couronnes. Nous d\'emontrons, de nouveau, quelques propri\'et\'es alg\'ebriques de ces anneaux, mais, cette fois-ci, de fa\c{c}on directe, sans avoir recours aux th\'eor\`emes de Weierstra{\ss}.

Nous entreprenons ensuite, au num\'ero \ref{exannloc}, une br\`eve \'etude de la topologie des espaces affines analytiques au voisinage de certains points. Nous en d\'eduisons une description explicite de certains anneaux locaux en termes d'alg\`ebres de s\'eries convergentes. 

Pour finir, le num\'ero \ref{henselianite} est consacr\'e \`a l'hens\'elianit\'e des anneaux locaux des espaces analytiques. Nous expliquons comment cette propri\'et\'e peut \^etre utilis\'ee pour d\'emontrer l'existence d'isomorphismes locaux entre espaces analytiques.


\pagebreak

\bigskip

Dans tout ce chapitre, nous fixons un anneau de Banach uniforme $(\As,\|.\|)$ et un entier positif~$n$. Nous noterons
$$B=\Ms(\As,\|.\|) \textrm{ et } X=\E{n}{(\As,\|.\|)}.$$ 
Les faisceaux structuraux sur ces espaces seront respectivement not\'es~$\Os_{B}$ et~$\Os_{X}$.

Nous noterons encore 
$$\pi : X \to B$$ 
le morphisme de projection induit par le morphisme naturel $\As\to \As[T_{1},\ldots,T_{n}]$. Pour toute partie~$V$ de~$B$, nous posons
$$X_{V} = \pi^{-1}(V)$$
et, pour tout point~$b$ de~$B$,
$$X_{b} = \pi^{-1}(b).$$

\section{Alg\`ebres globales de disques et de couronnes}\label{algglob}

Nous commen\c{c}ons par introduire quelques notations. Pour des \'el\'ements \mbox{$\bk=(k_{1},\ldots,k_{n})$} de~$\Z^n$ et $\bs = (t_{1},\ldots,t_{n})$ de~$(\R_{+}^*)^n$, posons
$$\bs^\bk = \prod_{i=1}^n s_{i}^{k_{i}}.$$
\newcounter{bs}\setcounter{bs}{\thepage}
\newcounter{bk}\setcounter{bk}{\thepage}
\newcounter{bsbk}\setcounter{bsbk}{\thepage}
D\'efinissons encore
$$\bT = (T_{1},\ldots,T_{n})$$\newcounter{bT}\setcounter{bT}{\thepage}
et, quel que soit $\bk=(k_{1},\ldots,k_{n})\in\Z^n$,
$$\bT^\bk = \prod_{i=1}^n T_{i}^{k_{i}}.$$
\newcounter{bTbk}\setcounter{bTbk}{\thepage}

Soit $\bt = (t_{1},\ldots,t_{n})\in (\R_{+}^*)^n$. Nous noterons 
$$\As\of{\la}{|\bT| \le \bt}{\ra}$$
\index{Disque|see{Couronne}}\index{Disque!algebre associee@alg\`ebre associ\'ee}
\newcounter{algd}\setcounter{algd}{\thepage}
l'alg\`ebre constitu\'ee des s\'eries de la forme 
$$\sum_{\bk \in \N^n} a_\bk\, \bT^\bk,$$ 
o\`u $(a_{\bk})_{\bk\in \Z^n}$ d\'esigne une famille de $\As$ v\'erifiant la condition suivante : 
$$\textrm{la famille } \left(\|a_\bk\|\, \bt^\bk\right)_{\bk \in \N^n} \textrm{ est sommable.}$$ 
Cette alg\`ebre est compl\`ete pour la norme d\'efinie par 
$$\left\|\sum_{\bk \in \N^n} a_\bk\, \bT^\bk \right\|_{\bt} = \sum_{\bk\in\N^n} \|a_\bk\|\, \bt^\bk.$$
\newcounter{normed}\setcounter{normed}{\thepage}
Comme nous l'expliquerons plus loin, elle est li\'ee \`a l'alg\`ebre des fonctions sur le disque de rayon~$\bt$ : 
$$\overline{D}(\bt) = \left\{x\in X\, \big|\, \forall i\in\cn{1}{n},\, |T_{i}(x)| \le t_{i} \right\}.$$
\newcounter{Dbt}\setcounter{Dbt}{\thepage}

D\'efinissons, \`a pr\'esent, deux relations, $\le$ et $<$, sur $\R^n$ de la fa\c{c}on suivante : pour deux \'el\'ements $\bs = (s_{1},\ldots,s_{n})$ et $\bt = (t_{1},\ldots,t_{n})$ de $\R^n$, nous posons
$$\bs \le \bt \textrm{ si } \forall i\in\cn{1}{n},\, s_{i} \le t_{i}$$
et
$$\bs < \bt \textrm{ si } \forall i\in\cn{1}{n},\, s_{i} < t_{i}.$$
\newcounter{relo}\setcounter{relo}{\thepage}

D\'efinissons \'egalement une relation~$\prec$ sur~$\R_{+}^n$ : pour deux \'el\'ements $\bs = (s_{1},\ldots,s_{n})$ et $\bt = (t_{1},\ldots,t_{n})$ de $\R_{+}^n$, nous posons
$$\bs \prec \bt \textrm{ si } \forall i\in\cn{1}{n},\, s_{i}<t_{i} \textrm{ ou } s_{i}=0.$$

\bigskip

Soient~$\bs$ et~$\bt$ dans~$(\R_{+}^*)^n$ v\'erifiant $\bs\le \bt$. Nous allons d\'efinir, sur le mod\`ele pr\'ec\'edent, une alg\`ebre associ\'ee \`a la couronne de rayon int\'erieur~$\bs$ et de rayon ext\'erieur~$\bt$ : 
$$\overline{C}(\bs,\bt) = \{x\in X\, |\, \forall i\in\cn{1}{n},\, s_{i}\le |T_{i}(x)| \le t_{i}\}.$$
\newcounter{Cbsbt}\setcounter{Cbsbt}{\thepage}
Pour $\bk=(k_{1},\ldots,k_{n}) \in \Z^n$, nous posons
$$\bmax(\bs^\bk,\bt^\bk) = \prod_{i=1}^n \max(s_{i}^{k_{i}},t_{i}^{k_{i}}) \in \of{]}{0,+\infty}{[}.$$
\newcounter{bmax}\setcounter{bmax}{\thepage}
Cette notation a \'et\'e choisie pour son caract\`ere naturel. Elle peut malheureusement pr\^eter \`a confusion :  attention \`a ne pas confondre la quantit\'e pr\'ec\'edente avec
$$\max(\bs^\bk,\bt^\bk) = \max\left(\prod_{i=1}^n s_{i}^{k_{i}}, \prod_{i=1}^n t_{i}^{k_{i}}\right).$$
Nous d\'efinissons l'alg\`ebre
$$\As\of{\la}{\bs \le |\bT| \le \bt}{\ra}$$
\index{Couronne!algebre associee@alg\`ebre associ\'ee}
\newcounter{algc}\setcounter{algc}{\thepage}
comme l'alg\`ebre constitu\'ee des s\'eries de la forme 
$$\sum_{\bk \in \Z^n} a_\bk\, \bT^\bk,$$ 
o\`u $(a_{\bk})_{\bk\in \Z^n}$ d\'esigne une famille de $\As$ v\'erifiant la condition suivante : 
$$\textrm{la famille } \left(\|a_\bk\|\, \bmax(\bs^\bk,\bt^\bk)\right)_{\bk \in \Z^n} \textrm{ est sommable.}$$ 
Cette alg\`ebre est compl\`ete pour la norme d\'efinie par 
$$\left\|\sum_{\bk \in \Z^n} a_\bk\, \bT^\bk \right\|_{\bs,\bt} = \sum_{\bk\in\Z^n} \|a_\bk\|\, \bmax(\bs^\bk,\bt^\bk).$$
\newcounter{normec}\setcounter{normec}{\thepage}

Afin de pouvoir traiter simultan\'ement les deux types d'alg\`ebres pr\'esent\'es ci-dessus, ainsi que celui associ\'e aux produits de disques et de couronnes, nous introduisons de nouvelles notations. Pour $\bk=(k_{1},\ldots,k_{n}) \in \Z^n$ et \mbox{$\bs = (s_{1},\ldots,s_{n})\in \R_{+}^n$} v\'erifiant la condition 
$$\forall i\in\cn{1}{n},\, (k_{i} <0 \implies s_{i} >0),$$
nous posons
$$\bs^\bk = \prod_{i=1}^n s_{i}^{k_{i}}.$$
\newcounter{bsbkbis}\setcounter{bsbkbis}{\thepage}
Pour~$k<0$, nous posons~$0^k = +\infty$. Pour $\bk=(k_{1},\ldots,k_{n}) \in \Z^n$,  \mbox{$\bs = (s_{1},\ldots,s_{n})\in \R_{+}^n$} et \mbox{$\bt = (t_{1},\ldots,t_{n})\in (\R_{+}^*)^n$}, nous posons
$$\bmax(\bs^\bk,\bt^\bk) = \prod_{i=1}^n \max(s_{i}^{k_{i}},t_{i}^{k_{i}}) \in \of{]}{0,+\infty}{]}.$$
\newcounter{bmaxbis}\setcounter{bmaxbis}{\thepage}
Si $\bs$ appartient \`a $(\R_{+}^*)^n$, nous posons
$$\bmin(\bs^\bk,\bt^\bk) = \prod_{i=1}^n \min(s_{i}^{k_{i}},t_{i}^{k_{i}}) \in \of{]}{0,+\infty}{[}.$$
\newcounter{bmin}\setcounter{bmin}{\thepage}

Soient $\bs = (s_{1},\ldots,s_{n})\in \R_{+}^n$ et $\bt = (t_{1},\ldots,t_{n})\in (\R_{+}^*)^n$  tels que $\bs\le \bt$. Dans la suite de ce paragraphe, nous nous int\'eresserons \`a l'alg\`ebre 
$$\As\of{\la}{\bs \le |\bT| \le \bt}{\ra}$$\newcounter{algcbis}\setcounter{algcbis}{\thepage}
constitu\'ee des s\'eries de la forme 
$$\sum_{\bk \in \Z^n} a_\bk\, \bT^\bk,$$ 
o\`u $(a_{\bk})_{\bk\in \Z^n}$ d\'esigne une famille de $\As$ v\'erifiant la condition suivante : 
$$\textrm{la famille } \left(\|a_\bk\|\, \bmax(\bs^\bk,\bt^\bk)\right)_{\bk \in \Z^n} \textrm{ est sommable.}$$ 
Remarquons, que s'il existe un indice $i\in\cn{1}{n}$ tel que $s_{i}=0$, alors, quel que soit $\bk\in\Z^n$ avec $k_{i}<0$, nous avons \mbox{$\bmax(\bs^\bk,\bt^\bk) = +\infty$}. La condition de sommabilit\'e impose alors que $a_{\bk}=0$.

L'alg\`ebre~$\As\of{\la}{\bs \le |\bT| \le \bt}{\ra}$ est compl\`ete pour la norme d\'efinie par 
$$\left\|\sum_{\bk \in \Z^n} a_\bk\, \bT^\bk \right\|_{\bs,\bt} = \sum_{\bk\in\Z^n} \|a_\bk\|\, \bmax(\bs^\bk,\bt^\bk).$$\newcounter{normecbis}\setcounter{normecbis}{\thepage}

L'alg\`ebre $\As\of{\la}{\bs \le |\bT| \le \bt}{\ra}$ est li\'ee \`a l'anneau des fonctions sur la couronne de rayon int\'erieur~$\bs$ et de rayon ext\'erieur~$\bt$ : 
$$\overline{C}(\bs,\bt) = \{x\in X\, |\, \forall i\in\cn{1}{n},\, s_{i}\le |T_{i}(x)| \le t_{i}\}.$$
\newcounter{Cbsbtbis}\setcounter{Cbsbtbis}{\thepage}
Pr\'ecisons ce r\'esultat.

\begin{prop}\label{spectrecouronne}
Le morphisme 
$$\Ms\left(\As\of{\la}{\bs \le |\bT| \le \bt}{\ra}\right) \to \E{n}{\As}$$
induit par l'injection naturelle
$$\As[\bT] \to \As\of{\la}{\bs \le |\bT| \le \bt}{\ra}$$
r\'ealise un hom\'eomorphisme sur son image $\overline{C}(\bs,\bt)$. En particulier, pour tout \'el\'ement~$f$ de~$\As\of{\la}{\bs \le |\bT| \le \bt}{\ra}$, nous avons
$$\|f\|_{\overline{C}(\bs,\bt)} = \inf_{j \ge 1}\left( \|f^j\|_{\bs,\bt}^{1/j}\right).$$
\end{prop}
\begin{proof}
Posons 
$$\Bs=\left\{ \sum_{\bk \in I} a_{\bk}\, \bT^\bk,\, I\subset J_{\bs} \right\},$$
o\`u $J_{\bs}$ d\'esigne l'ensemble des parties finies de l'ensemble
$$\{\bk=(k_{1},\ldots,k_{n})\in\Z^n\, |\, k_{i}\ge 0 \textrm{ si } s_{i}=0 \}.$$
Par exemple, si $\bs=0$, nous avons $\Bs=\As[\bT]$. L'anneau $\Bs$, qui est contenu dans l'anneau total des fractions de $\As[\bT]$, est dense dans $\As\of{\la}{\bs \le |\bT| \le \bt}{\ra}$ pour la norme $\|.\|_{\bs,\bt}$. On en d\'eduit que le morphisme 
$$\varphi : \Ms\left(\As\of{\la}{\bs \le |\bT| \le \bt}{\ra}\right) \to \E{n}{\As}$$
est injectif. Puisque l'espace $\Ms\left(\As\of{\la}{\bs \le |\bT| \le \bt}{\ra}\right)$ est compact, le morphisme $\varphi$ r\'ealise un hom\'eomorphisme sur son image. 

Il nous reste \`a montrer que l'image du morphisme $\varphi$ est \'egale \`a $\overline{C}(\bs,\bt)$. Soit $x\in \Ms\left(\As\of{\la}{\bs \le |\bT| \le \bt}{\ra}\right)$. Quel que soit $i\in\cn{1}{n}$, nous avons 
$$|T_{i}(x)|\le \|T_{i}\|_{\bs,\bt} = t_{i}.$$
Quel que soit $i\in\cn{1}{n}$, avec $s_{i}>0$, nous avons 
$$|T_{i}^{-1}(x)|\le \|T_{i}^{-1}\|_{\bs,\bt} = s_{i}^{-1}$$
et donc
$$|T_{i}(x)| \ge s_{i}.$$
On en d\'eduit que 
$$\varphi\left(\Ms\left(\As\of{\la}{\bs \le |\bT| \le \bt}{\ra}\right)\right) \subset \overline{C}(\bs,\bt).$$

R\'eciproquement, soit $x\in\overline{C}(\bs,\bt)$. Pour montrer que $x\in \Ms\left(\As\of{\la}{\bs \le |\bT| \le \bt}{\ra}\right)$, nous devons montrer que la semi-norme multiplicative $|.|_{x}$ sur $\As[\bT]$, born\'ee sur $\As$, associ\'ee \`a $x$ se prolonge en une semi-norme multiplicative born\'ee sur $(\As\of{\la}{\bs \le |\bT| \le \bt}{\ra},\|.\|_{\bs,\bt})$. Soit $i\in\cn{1}{n}$ tel que $s_{i}>0$. La fonction~$T_{i}$ ne s'annule pas sur la couronne~$\overline{C}(\bs,\bt)$. On en d\'eduit que la semi-norme multiplicative~$|.|_{x}$ se prolonge \`a~$\Bs$. Expliquons-en la raison. Pour $i\in\cn{1}{n}$, posons $r_{i}=0$ si $s_{i}=0$ et $r_{i}=1$ si $s_{i}>0$. Posons encore $\br=(r_{1},\ldots,r_{n})$. Tout \'el\'ement $Q$ de~$\Bs$ poss\`ede une \'ecriture sous la forme $\left(\bT^{-\br}\right)^l\, P$, avec $l\in\N$ et $P\in\As[\bT]$, et nous pouvons alors poser
$$|Q|_{x} = |\bT^\br|_{x}^{-l}\, |P|_{x}.$$
Cette quantit\'e ne d\'epend pas de l'\'ecriture de $Q$ choisie. On v\'erifie que l'application prolong\'ee, que nous notons encore $|.|_{x}$, d\'efinit bien une semi-norme multiplicative sur~$\Bs$.

Soit $i\in\cn{1}{n}$. Nous avons 
$$|T_{i}(x)| \le \max_{y\in \overline{C}(\bs,\bt)} (|T_{i}(y)|) = t_{i}.$$
Si $s_{i}>0$, nous avons \'egalement
$$|T_{i}^{-1}(x)| = |T_{i}(x)|^{-1} \le \min_{y\in \overline{C}(\bs,\bt)} (|T_{i}(y)|^{-1}) = s_{i}^{-1}.$$
Soit $Q(\bT)=\sum_{k\in\Z^n} a_{\bk}\, \bT^\bk\in \Bs$. Notons $b=\pi(x)$. Nous avons alors
$$|Q(\bT)|_{x} \le \disp \sum_{k\in\Z^n} |a_{\bk}(b)|\, \bmax(\bs^\bk,\bt^\bk)
\le   \disp \sum_{k\in\Z^n} \|a_{\bk}\|\, \bmax(\bs^\bk,\bt^\bk)
= \|P\|_{\bs,\bt}.$$
Le r\'esultat de densit\'e mentionn\'e plus haut montre finalement que~$|.|_{x}$ se prolonge en une semi-norme multiplicative born\'ee sur $\As\of{\la}{\bs \le |\bT| \le \bt}{\ra}$.\\
\end{proof}



Les r\'esultats qui suivent ont pour objet de comparer la norme $\|.\|_{\bs,\bt}$ et la norme uniforme $\|.\|_{\overline{C}(\bs,\bt)}$ sur la couronne $\overline{C}(\bs,\bt)$. Rappelons que nous avons suppos\'e que la norme $\|.\|$ d\'efinie sur l'anneau $\As$ est \'equivalente \`a la norme spectrale : il existe deux constantes $C_{-},C_{+}>0$ telles que
$$\forall f\in\As,\, C_{-}\, \|f\|_{sp} \le \|f\| \le C_{+}\, \|f\|_{sp}.$$

\begin{lem}\label{inegcoeff}
Soit $R = \sum_{k\in \Z^n} a_{\bk}\, \bT^\bk \in \As[\bT,\bT^{-1}]$. Quel que soit $\bk \in\Z^n$, nous avons 
$$\|a_{\bk}\|\, \bmax(\bs^\bk,\bt^\bk) \le C_{+}\, \|R\|_{\overline{C}(\bs,\bt)}.$$
\end{lem}
\begin{proof}
Commen\c{c}ons par remarquer que ce r\'esultat est bien connu lorsque l'anneau de Banach $(\As,\|.\|)$ est un corps valu\'e. En effet, lorsque le corps est ultram\'etrique, cela d\'ecoule imm\'ediatement de la description de la norme $\|R\|_{\overline{C}(\bs,\bt)}$ que l'on sait justement \^etre \'egale \`a
$$\max_{\bk \in\Z^n} \left(\|a_{\bk}\|\, \bmax(\bs^\bk,\bt^\bk) \right).$$
Lorsque le corps est archim\'edien, l'in\'egalit\'e provient de la formule de Cauchy.

Revenons au cas g\'en\'eral. Soit $\bk \in\Z^n$. Consid\'erons un point $z$ de $B$ en lequel l'\'egalit\'e $|a_\bk(z)|=\|a_\bk\|_{sp}$ a lieu. Il en existe car la partie $B$ est compacte. Le raisonnement pr\'ec\'edent assure que 
$$\|a_{\bk}\|\, \bmax(\bs^\bk,\bt^\bk) \le C_{+}\, |a_\bk(z)|\, \bmax(\bs^\bk,\bt^\bk) \le C_{+}\, \|R\|_{\pi^{-1}(z)\cap \overline{C}(\bs,\bt)}.$$
On en d\'eduit imm\'ediatement l'in\'egalit\'e demand\'ee.
\end{proof}

\begin{prop} \label{comparaisoncouronne}\index{Couronne!comparaison des normes}\index{Disque!comparaison des normes}
Soient $\bu=(u_{1},\ldots,u_{n})$ un \'el\'ement de~$\R_{+}^n$ et $\bv=(v_{1},\ldots,v_{n})$ un \'el\'ement de~$\left(\R_{+}^*\right)^n$ tels que~$\bs \prec \bu\le\bv < \bt$. Alors, pour tout \'el\'ement~$R$ de $\As\of{\la}{\bs\le |\bT|\le \bt}{\ra}$, on a l'in\'egalit\'e 
$$\|R\|_{\bu,\bv} \le C_{+}\,\left( \prod_{i=1}^n \frac{s_{i}}{u_{i}-s_{i}}+\frac{t_{i}}{t_{i}-v_{i}} \right) \, \|R\|_{\overline{C}(\bs,\bt)},$$
o\`u, pour tout \'el\'ement~$i$ de~$\cn{1}{n}$, nous posons $s_{i}/(u_{i}-s_{i})=0$ si~$s_{i}=0$.
\end{prop}
\begin{proof}
Il suffit de reprendre la preuve du lemme pr\'ec\'edent en rempla\c{c}ant l'anneau $\As[\bT]$ par l'anneau $\Bs$ introduit dans la d\'emonstration du lemme \ref{spectrecouronne}.
\end{proof}
\begin{proof}
Comme dans la preuve de la proposition \ref{spectrecouronne}, posons 
$$\Bs=\left\{ \sum_{\bk \in I} a_{\bk}\, \bT^\bk,\, I\subset J_{\bs} \right\},$$
o\`u $J_{\bs}$ d\'esigne l'ensemble des parties finies de l'ensemble
$$\{\bk=(k_{1},\ldots,k_{n})\in\Z^n\, |\, k_{i}\ge 0 \textrm{ si } s_{i}=0 \}.$$
Soit 
$$P=\sum_{\bk \in \N^n}  a_\bk\, \bT^\bk$$
un \'el\'ement de~$\Bs$. D'apr\`es le lemme pr\'ec\'edent, quel que soit \mbox{$\bk \in\Z^n$}, nous avons
$$\|a_\bk\|\, \bmax(\bs^\bk,\bt^\bk) \le C_{+}\, \|P\|_{\overline{C}(\bs,\bt)}.$$ 
On en d\'eduit que 
$${\renewcommand{\arraystretch}{2.3}\begin{array}{rcl}
\|P\|_{\bu,\bv} &=& \disp \sum_{\bk \in \Z^n} \|a_\bk\|\, \bmax(\bu^\bk,\bv^\bk)\\ 
&=& \disp C_{+}\, \|P\|_{\overline{C}(\bs,\bt)} \sum\limits_{\bk \in\Z^n}\left( \prod_{i=1}^n \frac{\max(u_{i}^{k_{i}},v_{i}^{k_{i}})}{\max(s_{i}^{k_{i}},t_{i}^{k_{i}})} \right)\\ 
& \le & \disp C_{+}\, \|P\|_{\overline{C}(\bs,\bt)}  \prod_{i=1}^n  \left( \sum_{k_{i}<0} \left(\frac{u_{i}}{s_{i}}\right)^{k_{i}} + \sum_{k_{i}\ge 0} \left(\frac{v_{i}}{t_{i}}\right)^{k_{i}} \right)\\ 
& \le &  \disp C_{+}\, \|P\|_{\overline{C}(\bs,\bt)}\, \left( \prod_{i=1}^n \frac{s_{i}}{u_{i}-s_{i}}+\frac{t_{i}}{t_{i}-v_{i}} \right).
\end{array}}$$
On conclut par densit\'e de $\Bs$ dans $\As\of{\la}{\bs \le |\bT|\le \bt}{\ra}$ pour la norme~$\|.\|_{\bs,\bt}$ et donc la norme~$\|.\|_{\bu,\bv}$.
\end{proof}

\section{Limites d'alg\`ebres de disques}\label{algloc}

Soit $V$ une partie compacte de $B$. Rappelons (\emph{cf.} d\'efinitions \ref{defKV} et \ref{defBV}) que $\Ks(V)$ d\'esigne le localis\'e de l'anneau $\As$ par l'ensemble des \'el\'ements qui ne s'annulent pas au voisinage de~$V$ et $\Bs(V)$ le compl\'et\'e de l'anneau $\Ks(V)$ pour la norme uniforme $\|.\|_{V}$ sur~$V$. Pour $\bt\in (\R_{+}^*)^n$, nous noterons $\|.\|_{V,\bt}$ la norme sur l'anneau $\Bs(V)\of{\la}{|\bT|\le \bt}{\ra}$ d\'efinie au paragraphe pr\'ec\'edent.
\newcounter{nVt}\setcounter{nVt}{\thepage}

Soit $b$ un point de $B$. Rappelons que nous notons $\m_{b}$ l'id\'eal maximal de l'anneau local $\Os_{B,b}$ et $\kappa(b)$ son corps r\'esiduel. Nous noterons 
$$L_{b} = \varinjlim_{V,\bt}\, \Bs(V)\of{\la}{|\bT|\le \bt}{\ra},$$
\newcounter{Lb}\setcounter{Lb}{\thepage}
o\`u $V$ parcourt l'ensemble des voisinages compacts du point $b$ dans $B$ et $\bt$ parcourt~$(\R_{+}^*)^n$. Pour commencer, \'enon\c{c}ons un lemme qui assure que certaines d\'ecompositions formelles, comme somme ou produit, d'\'el\'ements de $L_{b}$ existent dans $L_{b}$.

\begin{lem}\label{formelLb} 
Soit 
$$G = \sum_{\bk\ge 0} a_{\bk}\, \bT^\bk \in L_{b}.$$
Soit $E$ une partie de $\N^n$. Alors les s\'eries 
$$G_{1} = \sum_{\bk\in E} a_{\bk}\, \bT^\bk \textrm{ et } G_{2} = \sum_{\bk\notin E} a_{\bk}\, \bT^\bk$$
appartiennent \`a $L_{b}$ et v\'erifient
$$G = G_{1} + G_{2}.$$

Soit $i\in\cn{1}{n}$. Supposons qu'il existe $H \in \Os_{B,b}[\![\bT]\!]$ telle que $G = T_{i}\, H$. Alors~$H$ appartient \`a $L_{b}$ et l'\'egalit\'e $G = T_{i}\, H$ vaut dans $L_{b}$.  
\end{lem}
\begin{proof}
Il suffit de revenir \`a la d\'efinition des \'el\'ements de $L_{b}$ et de prendre garde \`a ce que les conditions de convergence restent v\'erifi\'ees.
\end{proof}

\begin{lem}\label{idealmax}
L'anneau $L_{b}$ est un anneau local d'id\'eal maximal 
$$\m = (\m_{b},T_{1},\ldots,T_{n})$$ 
et de corps r\'esiduel~$\kappa(b)$.
\end{lem}
\begin{proof}
En utilisant le lemme~\ref{formelLb}, on montre que le morphisme naturel
$$\Os_{B,b} \to L_{b}/(T_{1},\ldots,T_{n})$$
est un isomorphisme. On en d\'eduit un isomorphisme naturel 
$$\kappa(b) = \Os_{B,b}/\m_{b} \xrightarrow[]{\sim} L_{b}/\m.$$
Par cons\'equent, l'id\'eal $\m$ est maximal. 

Pour montrer que l'anneau $L_{b}$ est un anneau local d'id\'eal $\m$, il suffit de montrer que tout \'el\'ement de $L_{b}$ qui n'appartient pas \`a $\m$ est inversible. Soit $F\in L_{b}\setminus\m$. Il existe $V$ un voisinage compact du point $b$ dans $B$ et $t\in (\R_{+}^*)^n$ tels que $F \in \Bs(V)\of{\la}{|\bT|\le \bt}{\ra}$. Nous pouvons \'ecrire $F$ sous la forme
$$F = a_{0} + \sum_{i=1}^n T_{i}\, G_{i}(\bT),$$
avec $a_{0}\in \Bs(V)$ et, quel que soit $i\in\cn{1}{n}$, $G_{i}\in \Bs(V)\of{\la}{|\bT|\le \bt}{\ra}$. Puisque $F$ n'appartient pas \`a $\m$, son premier coefficient $a_{0}$ n'appartient pas \`a $\m_{b}$. On en d\'eduit que $a_{0}$ est inversible au voisinage de $b$ dans $B$. Quitte \`a restreindre $V$ et \`a multiplier $F$ par $a_{0}^{-1}$, nous pouvons supposer que $a_{0}=1$. Notons 
$$M = \max_{1\le i\le n} (\|G_{i}\|_{V,\bt}).$$
Soit $\bs=(s_{1},\ldots,s_{n})\in (\R_{+}^*)^n$ tel que 
$$\sum_{i=1}^n s_{i}\, M <1.$$
Nous avons alors
$$\left\| \sum_{i=1}^n T_{i}\, G_{i}(\bT) \right\|_{V,\bs} <1.$$
On en d\'eduit que la fonction 
$$F = 1 + \sum_{i=1}^n T_{i}\, G_{i}(\bT) $$
est inversible dans l'anneau de Banach $\Bs(V)\of{\la}{|\bT|\le \bs}{\ra}$ et donc dans $L_{b}$.
\end{proof}

\subsection{Th\'eor\`emes de Weierstra\ss}

Dans ce paragraphe, nous montrerons que l'anneau $L_{b}$ satisfait les conclusions des th\'eor\`emes de division et de pr\'eparation de Weierstra\ss. Notre preuve est calqu\'ee sur celle que mettent en {\oe}uvre  H. Grauert et R. Remmert dans le cadre de la g\'eom\'etrie analytique complexe. 

Nous noterons $\bT'=(T_{1},\ldots,T_{n-1})$ et 
$$L'_{b} = \varinjlim_{V,\bt'}\, \Bs(V)\of{\la}{|\bT'|\le \bt'}{\ra},$$
o\`u $V$ parcourt l'ensemble des voisinages compacts du point $b$ dans $B$ et $\bt'$ parcourt l'ensemble $(\R_{+}^*)^{n-1}$. 

\begin{thm}[Th\'eor\`eme de division de Weierstra\ss]\label{division}\index{Theoreme de Weierstrass@Th\'eor\`eme de Weierstra{\ss}!division locale}\index{Theoreme@Th\'eor\`eme!de Weierstrass@de Weierstra{\ss}|see{Th\'eor\`eme de Weierstra{\ss}}}\index{Weierstrass@Weierstra{\ss}|see{Th\'eor\`eme de Weierstra{\ss}}}
Soit $G \in L_{b}$ une s\'erie telle que 
$G(0,T_n)(b)\ne 0 \textrm{ dans } \Hs(b)[\![T_{n}]\!].$
Notons $p$ la valuation en $T_{n}$ de la s\'erie $G(0,T_n)(b)$. Soit $F\in L_{b}$. Alors il existe un unique couple $(Q,R) \in (L_{b})^2$ tel que 
\begin{enumerate}[\it i)]
\item $R\in L'_{b}[T_n]$ est un polyn\^ome de degr\'e strictement inf\'erieur \`a $p$ ;
\item $F=QG+R.$
\end{enumerate}
\end{thm}
\begin{proof}
Notons $G=\sum_{k\in\N} g_k(\bT')\, T_n^k$ o\`u, quel que soit $k\in\N$, \mbox{$g_k\in L'_{b}$}, $g_0(0)(b)=\cdots=g_{p-1}(0)(b)=0$ et $g_p(0)(b)\ne 0$. Quitte \`a choisir un voisinage compact assez petit $V$ du point $b$ et un r\'eel strictement positif $r$ assez petit \'egalement, nous pouvons supposer que $G\in \Bs(V)\of{\la}{|\bT|\le \br}{\ra}$, o\`u $\br=(r,\ldots,r)\in (\R_{+}^*)^n$, et que $g_p(\bT')$ est inversible dans $\Bs(V)\of{\la}{|\bT'|\le \br'}{\ra}$, o\`u $\br'=(r,\ldots,r)\in (\R_{+}^*)^{n-1}$. Quitte \`a multiplier alors $G$ par $g_p^{-1}$, nous pouvons supposer que $g_p=1$.

Soient $\bs' \in (\R_{+}^*)^{n-1}$, avec $\bs'\le \br'$, et $s\in\of{]}{0,r}{]}$. Posons $\bs=(\bs',s)\in (\R_{+}^*)^{n}$. 
Tout \'el\'ement $\varphi$ de $\Bs(V)\of{\la}{|\bT|\le \bs}{\ra}$ peut s'\'ecrire de fa\c{c}on unique sous la forme 
$$\varphi = \alpha(\varphi)\, T_n^p + \beta(\varphi),$$
o\`u $\alpha(\varphi)$ d\'esigne un \'el\'ement de $\Bs(V)\of{\la}{|\bT|\le \bs}{\ra}$ et $\beta(\varphi)$ un \'el\'ement de $\Bs(V)\of{\la}{|\bT'|\le \bs'}{\ra}[T_n]$ de degr\'e strictement inf\'erieur \`a $p$. Remarquons, d\`es \`a pr\'esent, que, quel que soit $\varphi \in \Bs(V)\of{\la}{|\bT|\le \bs}{\ra}$, on a 
$$\|\varphi\|_{V,\bs} = \|\alpha(\varphi)\|_{V,\bs}\, s^p + \|\beta(\varphi)\|_{V,\bs}.$$

Consid\'erons, \`a pr\'esent, l'endomorphisme 
$$ A_{\bs} :
{\renewcommand{\arraystretch}{1.2}\begin{array}{ccc}
\Bs(V)\of{\la}{|\bT|\le \bs}{\ra} & \to &\Bs(V)\of{\la}{|\bT|\le \bs}{\ra} \\
\varphi & \mapsto & \alpha(\varphi)\, G + \beta(\varphi)
\end{array}}.$$
Il nous suffit de trouver un $n$-uplet $\bs$ assez petit pour lequel l'endomorphisme~$A_{\bs}$ soit bijectif. Remarquons que, quel que soit $\varphi \in \Bs(V)\of{\la}{|\bT|\le \bs}{\ra}$, on a 
$${\renewcommand{\arraystretch}{1.3}\begin{array}{rcl}
\|A_{\bs}(\varphi)-\varphi\|_{V,\bs} &=& \|\alpha(\varphi)\, (G-T_n^p)\|_{V,\bs}\\
&\le&  \|\alpha(\varphi)\|_{V,\bs}\, \|G-T_n^p\|_{V,\bs}\\
&\le& s^{-p}\,\|\varphi\|_{V,\bs}\, \|G-T_n^p\|_{V,\bs}.
\end{array}}$$

Soient $u,v\in\of{]}{0,\min(r,1)}{[}$. Nous noterons $(\bu,v)$ le $n$-uplet $(u,\ldots,u,v)$. Soit \mbox{$k\in\cn{0}{p-1}$}. Il existe une constante $M_{k}\in\R$, ind\'ependante de $u$ et de $v$, telle que l'on ait 
$$\|g_k\|_{V,\bu}\le \|g_{k}(0)\|_{V} + M_{k}\, u.$$ 
Il existe \'egalement une constante $N\in\R$, encore ind\'ependante de $u$ et de $v$, telle que l'on ait 
$$\left\|\sum_{k\ge p+1} g_k(\bT')\, T_n^k\right\|_{V,(\bu,v)}\le N v^{p+1}.$$ 
Par cons\'equent, il existe une constante $M\in\R$, ind\'ependante de $u$ et de $v$, telle que
$$\|G-T_n^p\|_{V,(\bu,v)} \le \sum_{k=0}^{p-1} \|g_{k}(0)\|_{V} + M (u+v^{p+1}).$$

Soit $\eps\in\of{]}{0,1}{[}$. Quitte \`a choisir judicieusement $v$ puis $u$, nous pouvons supposer que $M(u+v^{p+1})\le \eps v^{p}/2$. Quel que soit $k\in\cn{0}{p-1}$, nous avons $g_{k}(0)(b)=0$, par hypoth\`ese. Par cons\'equent, quitte \`a restreindre le voisinage $V$ de $b$, nous pouvons supposer que 
$$\sum_{k=0}^{p-1} \|g_{k}(0)\|_{V}\le \eps v^p/2.$$ 
On dispose alors de l'in\'egalit\'e $$\|A_{(\bu,v)}-I\|_{V,(\bu,v)}\le \eps<1$$
et on en d\'eduit que l'endomorphisme $A_{(\bu,v)}=I+(A_{(\bu,v)}-I)$ est inversible.
\end{proof}  

Nous pouvons obtenir une version plus pr\'ecise du th\'eor\`eme de Weierstra{\ss} lorsque l'on divise par des s\'eries d'un type particulier. 

\begin{defi}\label{distingue}\index{Polynome distingue@Polyn\^ome distingu\'e}
Soit $p\in\N$. Nous dirons qu'un polyn\^ome $h\in L'_{b}[T_n]$ est \textbf{distingu\'e} de degr\'e $p$ s'il est unitaire, de degr\'e $p$ et v\'erifie 
$$h(0,T_n)(b)=T_n^p \textrm{ dans } \Hs(b)[\![T_{n}]\!].$$
\end{defi}

\begin{thm}[Th\'eor\`eme de division de Weierstra{\ss} semi-local]\label{global}\index{Theoreme de Weierstrass@Th\'eor\`eme de Weierstra{\ss}!division semi-locale}
Soient $p\in\N$ et $G \in L'_{b}[T_n]$ un polyn\^ome distingu\'e de degr\'e $p$. Soient $V$ un voisinage compact de $b$ dans $B$ et $\br' \in (\R_{+}^*)^{n-1}$ tel que $G\in\Bs(V)\of{\la}{|\bT'|\le \br'}{\ra}[T_{n}]$. Soient~$v_{-}$ et~$v_{+}$ deux nombres r\'eels v\'erifiant $0<v_{-}\le v_{+}$. Alors il existe un voisinage compact $W$ de $b$ dans $V$ et un $(n-1)$-uplet $\bs' \in (\R_{+}^*)^{n-1}$, avec $\bs'\le \br'$, v\'erifiant la propri\'et\'e suivante : pour tout voisinage compact $U$ de $b$ dans $W$, tout $(n-1)$-uplet $\bt' \in (\R_{+}^*)^{n-1}$  v\'erifiant $\bt' \le \bs'$, tout nombre r\'eel $w\in\of{[}{v_{-},v_{+}}{]}$ et tout \'el\'ement $F$ de $\Bs(U)\of{\la}{|\bT|\le (\bt',w)}{\ra}$, il existe un unique couple $(Q,R) \in (\Bs(U)\of{\la}{|\bT|\le (\bt',w)}{\ra})^2$ tel que 
\begin{enumerate}[\it i)]
\item $R$ soit un polyn\^ome de degr\'e strictement inf\'erieur \`a $p$ ;
\item $F=QG+R$.
\end{enumerate}
En outre, il existe une constante $C\in\R_{+}^*$, ind\'ependante de $U$, $\bt'$, $w$ et $F$, telle que l'on ait les in\'egalit\'es
\begin{enumerate}[a)]
\item $\|Q\|_{U,(\bt',w)} \le C\,\|F\|_{U,(\bt',w)}$ ;
\item $\|R\|_{U,(\bt',w)} \le C\,\|F\|_{U,(\bt',w)}$.
\end{enumerate}
\end{thm}
\begin{proof}
Notons $$G= T_{n}^p + \sum_{k=0}^{p-1} g_k(\bT')\, T_n^k$$ o\`u, quel que soit \mbox{$k\in\cn{0}{p-1}$}, \mbox{$g_k\in \Bs(V)$} et $g_k(0)(b)=0$. Soient $\bs' \in (\R_{+}^*)^{n-1}$, avec $\bs'\le \br'$, $u\in\of{]}{0,v_{+}}{]}$ et $W$ un voisinage compact de $b$ dans $V$. Tout \'el\'ement $\varphi$ de $\Bs(W)\of{\la}{|\bT|\le (\bs',u)}{\ra}$ peut s'\'ecrire de fa\c{c}on unique sous la forme 
$$\varphi = \alpha(\varphi)\, T_n^p + \beta(\varphi),$$
o\`u $\alpha(\varphi)$ d\'esigne un \'el\'ement de $\Bs(W)\of{\la}{|\bT|\le (\bs',u)}{\ra}$ et $\beta(\varphi)$ un \'el\'ement de $\Bs(W)\of{\la}{|\bT'|\le \bs'}{\ra}[T_n]$ de degr\'e strictement inf\'erieur \`a $p$. Remarquons, d\`es \`a pr\'esent, que, quel que soit $\varphi \in \Bs(W)\of{\la}{|\bT|\le (\bs',u)}{\ra}$, nous avons 
$$\|\varphi\|_{W,(\bs',u)} = \|\alpha(\varphi)\|_{W,(\bs',u)}\, u^p + \|\beta(\varphi)\|_{W,(\bs',u)}.$$

Consid\'erons, \`a pr\'esent, l'endomorphisme 
$$ A_{W,(\bs',u)} :
{\renewcommand{\arraystretch}{1.2}\begin{array}{ccc}
\Bs(W)\of{\la}{|\bT|\le (\bs',u)}{\ra} & \to &\Bs(W)\of{\la}{|\bT|\le (\bs',u)}{\ra} \\
\varphi & \mapsto & \alpha(\varphi)\, G + \beta(\varphi)
\end{array}}.$$
Remarquons que, quel que soit $\varphi \in \Bs(W)\of{\la}{|\bT|\le (\bs',u)}{\ra}$, nous avons
$${\renewcommand{\arraystretch}{1.3}\begin{array}{rcl}
\|A_{W,(\bs',u)}(\varphi)-\varphi\|_{W,(\bs',u)} &=& \|\alpha(\varphi)\, (G-T_n^p)\|_{W,(\bs',u)}\\
&\le&  \|\alpha(\varphi)\|_{W,(\bs',u)}\, \|G-T_n^p\|_{W,(\bs',u)}\\
&\le& u^{-p}\,\|\varphi\|_{W,(\bs',u)}\, \|G-T_n^p\|_{W,(\bs',u)}.
\end{array}}$$

Si $\bs'=(s_{1},\ldots,s_{n-1})$, nous noterons $\max(\bs') = \max(s_{1},\ldots,s_{n-1})$. Soit \mbox{$k\in\cn{0}{p-1}$}. Il existe une constante $M_k\in\R$, ind\'ependante de $\bs'$, telle que l'on ait 
$$\|g_k\|_{W,\bs'}\le \|g_{k}(0)\|_{W} + M_{k}\, \max(\bs').$$ 
Par cons\'equent, il existe une constante $M\in\R$, ind\'ependante de $\bs'$, telle que l'on ait
$${\renewcommand{\arraystretch}{1.5}\begin{array}{rcl}
\|G-T_n^p\|_{W, (\bs',u)} &\le& \disp \sum_{k=0}^{p-1} \|g_{k}(0)\|_{W}\, u^k + M  \max(\bs')\\
&\le& \disp \sum_{k=0}^{p-1} \|g_{k}(0)\|_{W}\, v_{+}^k + M  \max(\bs').
\end{array}}$$

Soit $\eps\in\of{]}{0,1}{[}$. Quel que soit $k\in\cn{0}{p-1}$, nous avons $g_{k}(0)(b)=0$, par hypoth\`ese. Par cons\'equent, il existe un voisinage $W$ de $b$ dans $V$ tel que l'on ait 
$$\disp \sum_{k=0}^{p-1} \|g_{k}(0)\|_{W}\, v_{+}^k \le \eps\, \frac{v_{-}^p}{2}.$$
Il existe \'egalement $\bs'\le \br'$ tel que
$$M \max(\bs') \le \eps\, \frac{v_{-}^p}{2}.$$ 

Soient $U$ un voisinage compact de $b$ dans $W$, $\bt'\le\bs'$ et $w\in\of{[}{v_{-},v_{+}}{]}$. On dispose alors de l'in\'egalit\'e 
$${\renewcommand{\arraystretch}{1.5}\begin{array}{rcl}
\|G-T_{n}^p\|_{U,(\bt',w)}\, w^{-p} &\le & \|G-T_{n}^p\|_{W,(\bs',w)}\, v_{-}^{-p}\\
&\le& \left(\sum\limits_{k=0}^{p-1} \|g_{k}(0)\|_{W}\, v_{+}^k + M \max(\bs')\right) v_{-}^{-p}\\
&\le & \eps\, v_{-}^p\, v_{-}^{-p} \le\eps.
\end{array}}$$
Nous avons donc 
$$\|A_{U,(\bt',w)}-I\|_{U,(\bt',w)}\le \eps<1.$$
Par cons\'equent, l'endomorphisme $A_{U,(\bt',w)}=I+(A_{U,(\bt',w)}-I)$ est inversible.

Soit $F\in \Bs(U)\of{\la}{|\bT|\le (\bt',w)}{\ra}$. Il existe un unique couple $(Q,R)$, avec \mbox{$Q\in\Bs(U)\of{\la}{|\bT|\le (\bt',w)}{\ra}$} et $R\in\Bs(U)\of{\la}{|\bT'|\le \bt'}{\ra}[T_n]$ de degr\'e strictement inf\'erieur \`a $p$, tel que 
$$F=QG+R.$$
Avec les notations pr\'ec\'edentes, nous avons $Q=\alpha(A_{U,(\bt',w)}^{-1}(F))$ et $R=\beta(A_{U,(\bt',w)}^{-1}(F))$. Puisque $\|A_{U,(\bt',w)}-I\|_{U,|\bT|\le (\bt',w)}\le \eps$, nous avons 
$$\|A_{U,(\bt',w)}^{-1}\|_{U,(\bt',w)}\le \sum_{i=0}^{+\infty} \eps^i = \frac{1}{1-\eps}.$$ 
On en d\'eduit que 
$$\|Q\|_{U,(\bt',w)}\le \frac{v_{-}^{-p}}{1-\eps}\, \|F\|_{U,(\bt',w)}$$
et que 
$$\|R\|_{U,(\bt',w)}\le \frac{1}{1-\eps}\, \|F\|_{U,(\bt',w)}.$$
\end{proof}

\begin{thm}[Th\'eor\`eme de pr\'eparation de Weierstra\ss]\label{preparation}\index{Theoreme de Weierstrass@Th\'eor\`eme de Weierstra{\ss}!preparation locale@pr\'eparation locale}
Soit $G \in L_{b}$ une s\'erie telle que $G(0,T_n)(b)\ne 0$ dans $\Hs(b)[\![T_{n}]\!]$. Notons $p$ la valuation en $T_{n}$ de la s\'erie $G(0,T_n)(b)$. Alors il existe un unique couple \mbox{$(\Omega,E) \in (L_{b})^2$} v\'erifiant les conditions suivantes :
\begin{enumerate}[\it i)]
\item $\Omega \in L'_{b}[T_n]$ est un polyn\^ome distingu\'e de degr\'e $p$ ;
\item $E$ est inversible dans $L_{b}$ ;
\item $G=E\, \Omega.$
\end{enumerate}
\end{thm}
\begin{proof}
Supposons que des s\'eries $\Omega$ et $E$ v\'erifiant les conditions requises existent. Alors $\Omega$ s'\'ecrit sous la forme $T_n^p + S$, o\`u $S\in L'_{b}[T_n]$ d\'esigne un polyn\^ome de degr\'e strictement inf\'erieur \`a $p$. Les s\'eries $S$ et $E$ sont alors reli\'ees par l'\'egalit\'e $T_n^p = E^{-1}\, G -S$. Le th\'eor\`eme de division de {Weierstra\ss} \ref{division} nous assure l'unicit\'e des s\'eries $E^{-1}$ et $S$. On en d\'eduit l'unicit\'e des s\'eries $\Omega$ et $E$.

D\'emontrons, \`a pr\'esent, l'existence de ces s\'eries. Le th\'eor\`eme \ref{division} appliqu\'e avec $T_n^p$ et $G$ nous assure qu'il existe $Q \in L_{b}$ et $R \in L'_{b}[T_n]$ de degr\'e strictement inf\'erieur \`a $p$ tels que
$$T_n^p = QG+R.$$
Montrons, tout d'abord, que $R(0,T_n)(b)=0$. Si $H$ d\'esigne un \'el\'ement de $L_{b}$, nous noterons $v_{b}(H)$ la valuation en $T_{n}$ de la s\'erie $H(0,T_n)(b)$ dans $\Hs(b)[\![T_{n}]\!]$. 

Nous avons alors 
$${\renewcommand{\arraystretch}{1.3}\begin{array}{rcl}
v_{b}(R) &=& v_{b}(T_{n}^p-QG)\\
&\ge & \min\left(v_{b}(T_{n}^p), v_{b}(Q) + v_{b}(G)\right)\\
&\ge & p.
\end{array}}$$
Puisque $R(0,T_{n})$ est suppos\'e de degr\'e strictement inf\'erieur \`a $p$, nous avons donc $R(0,T_n)(b)=0$. On en d\'eduit que $v_{b}(T_n^p-R)=p$ et donc que 
$$v_{b}(Q) = v_{b}(QG) - v_{b}(G) = p -p =  0.$$ 
Par cons\'equent, $Q$ est inversible dans $L_{b}$. Les s\'eries $E=Q^{-1}$ et $\Omega = T_{n}^p - R$ conviennent.
\end{proof}

Par la suite, nous aurons \'egalement besoin du lemme suivant, fort utile pour nous ramener \`a une situation dans laquelle on peut utiliser les th\'eor\`emes de {Weierstra\ss}.

\begin{lem}\label{auto}
Soit $G \in L_{b}$ tel que $G(b) \ne 0$ dans $\Hs(b)[\![\bT]\!]$. Il existe un automorphisme $\sigma$ de $L_{b}$ tel que l'on ait $\sigma(G)(0,T_{n})(b)\ne 0$ dans $\Hs(b)[\![T_{n}]\!]$.
\end{lem}
\begin{proof}
D'apr\`es \cite{diviseurs}, \S 3, \no 7, lemme 3, il existe \mbox{$u(1),\ldots,u(n-1)\in\N^*$} tels que l'automorphisme $\tau$ de $\Hs(b)[\![\bT]\!]$ d\'efini par 
\begin{enumerate}[\it i)]
\item $\forall i\in\cn{1}{n-1},\, \tau(T_{i}) = T_{i} + T_{n}^{u(i)}$ ;
\item $\tau(T_{n}) = T_{n}$
\end{enumerate}
envoie $G$ sur un \'el\'ement $\tau(G)$ qui v\'erifie $\tau(G)(0,T_{n})(b) \ne 0$. 

Montrons que l'application $\tau$ peut \^etre d\'efinie sur $L_{b}$. Soient $U$ un voisinage compact de $b$ dans $B$ et $\br=(r_{1},\ldots,r_{n}) \in (\R_{+}^*)^n$. Quel que soit $i\in\cn{0}{n-1}$, il existe $s_{i},s_{n,i} \in\R_{+}^*$ tels que $s_{i} + s_{n,i}^{u(i)} \le r_{i}$. Posons $s_{n} = \min(s_{n,1},\ldots,s_{n,n-1},r_{n})$ et $\bs=(s_{1},\ldots,s_{n})$. D\'efinissons alors un endomorphisme $\tau_{U}$ de $\Bs(U)[\![\bT]\!]$ par les m\^emes formules que $\tau$. On v\'erifie alors que, quel que soit $F\in\Bs(U)\of{\la}{|\bT|\le \br}{\ra}$, on a 
$$\tau_{U}(F)\in\As(U)\of{\la}{|\bT|\le \bs}{\ra}.$$ 
On en d\'eduit un morphisme $\sigma_{U} : \Bs(U)\of{\la}{|\bT|\le \br}{\ra} \to L_{b}$. On v\'erifie sans peine que tous ces morphismes sont compatibles et d\'efinissent donc un endomorphisme $\sigma$ de $L_{b}$. En outre, l'endomorphisme $\sigma$ induit l'endomorphisme $\tau$ sur $\Os_{B,b}[\![\bT]\!]$. On en d\'eduit, en particulier, que $\sigma(G)(0,T_{n})(b)\ne 0$.

En appliquant le m\^eme proc\'ed\'e \`a partir de $\tau^{-1}$, on construit un endomorphisme~$\sigma^{-1}$ de $L_{b}$ qui est l'inverse de $\sigma$. Par cons\'equent, $\sigma$ est un automorphisme de $L_{b}$.
\end{proof}

\subsection{Propri\'et\'es}

Nous consacrerons cette partie \`a d\'emontrer quelques propri\'et\'es de l'anneau local~$L_{b}$. 

\begin{thm}\label{noetheriencorps}
Supposons que l'anneau local $\Os_{B,b}$ est un corps. Alors l'anneau local~$L_{b}$ est noeth\'erien.
\end{thm}
\begin{proof}
Nous allons proc\'eder par r\'ecurrence. Si $n=0$, l'isomorphisme \mbox{$L_{b}\simeq \Os_{B,b}$} nous montre que le r\'esultat est vrai.

Supposons, \`a pr\'esent, que le r\'esultat soit vrai pour $L'_{b}$. Soit $I$ un id\'eal de $L_{b}$. L'id\'eal nul \'etant \'evidemment de type fini, nous pouvons supposer que $I\ne (0)$. Choisissons un \'el\'ement non nul $G$ de $I$. Puisque $\Os_{B,b}$ est un corps, il s'injecte dans $\Hs(b)$ et nous avons donc $G(b)\ne 0$. D'apr\`es le lemme \ref{auto}, quitte \`a appliquer un automorphisme de $L_{b}$, nous pouvons donc supposer que $G(0,T_n)(b)\ne 0$. D'apr\`es le th\'eor\`eme de division de {Weierstra\ss} \ref{division}, l'id\'eal $I$ est engendr\'e par $G$ et par la partie $I\cap L'_{b}[T_n]$. Or l'anneau $L'_{b}[T_n]$ est noeth\'erien, puisque $L'_{b}$ l'est, donc l'id\'eal $I\cap L'_{b}[T_n]$ est engendr\'e par un nombre fini d'\'el\'ements, ce qui suffit pour conclure.
\end{proof}

Nous souhaitons, maintenant, traiter le cas o\`u l'anneau local $\Os_{B,b}$ est un anneau de valuation discr\`ete. Nous aurons besoin d'une hypoth\`ese suppl\'ementaire. 

\begin{defi}\label{condU}\index{Condition (U)}
Soit~$b$ un point de~$B$ en lequel l'anneau local~$\Os_{B,b}$ est de valuation discr\`ete. Choisissons une uniformisante~$\pi$ de cet anneau et~$V$ un voisinage de~$b$ dans~$B$ sur lequel elle est d\'efinie. Nous dirons que l'uniformisante~$\pi$ v\'erifie la {\bf condition ($\textrm{U}_{V}$)} s'il existe une constante $C_{V}>0$ telle que pour toute fonction $f\in\Bs(V)$ v\'erifiant $f(b)=0$, il existe une fonction $g\in\Bs(V)$ v\'erifiant les propri\'et\'es suivantes :
\begin{enumerate}[\it i)]
\item $f = \pi\, g$ dans $\Bs(V)$ ;
\item $\|g\|_{V} \le C_{V}\, \|f\|_{V}$.
\end{enumerate}

Nous dirons que l'anneau de valuation discr\`ete~$\Os_{B,b}$ v\'erifie la {\bf condition (U)} s'il existe une uniformisante~$\pi$ de~$\Os_{B,b}$ d\'efinie sur un voisinage~$V$ du point~$b$ dans~$B$ et un syst\`eme fondamental~$\Ws$ de voisinages compacts du point~$b$ dans $V$ tel que, pour tout \'el\'ement~$W$ de~$\Ws$, l'uniformisante~$\pi$ v\'erifie la condition ($\textrm{U}_{W}$).

\end{defi}

\begin{rem}
Il est clair que la condition~(U) ne d\'epend pas de l'ouvert de d\'efinition~$V$ de~$\pi$ que nous avons choisi. En outre, si $\pi'$ d\'esigne une uniformisante de $\Os_{B,b}$, il existe une fonction $\alpha$ inversible dans $\Os_{B,b}$ telle que $\pi = \alpha\, \pi'$ dans $\Os_{B,b}$. Si les propri\'et\'es pr\'ec\'edentes sont v\'erifi\'ees pour l'uniformisante $\pi$, elles le sont donc encore pour l'uniformisante $\pi'$. Par cons\'equent, la condition~(U) porte bien sur l'anneau local lui-m\^eme et ne d\'epend pas des choix de $\pi$ et de $V$ effectu\'es.
\end{rem}

Nous utiliserons la condition~(U) sous la forme du lemme suivant.
\begin{lem}\label{lemU}
Supposons que l'anneau local $\Os_{B,b}$ est un anneau de valuation discr\`ete v\'erifiant la condition (U). Soit $\pi$ une uniformisante de $\Os_{B,b}$ et notons $v_{\pi}$ la valuation $\pi$-adique sur cet anneau. Soit $G\in L_{b}\setminus\{0\}$. Notons $\sum_{\bk\ge 0} a_{\bk} \bT^\bk$ son image dans $\Os_{B,b}[\![ \bT]\!]$. Posons 
$$v(G) = \min\{ v_{\pi}(a_{\bk}),\, \bk \ge 0 \} \in \N.$$
Alors, il existe une fonction $H$ de $L_{b}$ v\'erifiant les propri\'et\'es suivantes :
\begin{enumerate}[\it i)]
\item $H(b) \ne 0$ dans $\Hs(b)[\![\bT]\!]$ ; 
\item $G = \pi^{v(G)} H$ dans $L_{b}$.
\end{enumerate}
\end{lem}
\begin{proof}
Soit $V$ un voisinage de $b$ dans $B$ sur lequel $\pi$ est d\'efinie. Par hypoth\`ese, il existe un syst\`eme fondamental $\Ws$ de voisinages de $b$ dans $V$ tel que, quel que soit $W\in\Ws$, l'uniformisante~$\pi$ v\'erifie la condition ($\textrm{U}_{W}$), avec une certaine constante $C_{W}>0$. 
Il existe un voisinage compact $U$ de $b$ dans $B$ et $\bt\in (\R_{+}^*)^n$ tels que la s\'erie $G$ soit un \'el\'ement de $\Bs(U)\of{\la}{|\bT|\le \bt}{\ra}$. Par cons\'equent, il existe une famille $(a_{\bk})_{\bk\ge 0}$ d'\'el\'ements de $\Bs(U)$ telle que
$$G = \sum_{\bk \ge 0} a_{\bk}\, \bT^\bk$$
et
$$\sum_{\bk\ge 0} \|a_{\bk}\|_{U}\, \bt^\bk < +\infty.$$

Soit $W$ un \'el\'ement de $\Ws$ contenu dans $U$. Soit $\bk\ge 0$. Par hypoth\`ese, $\pi^{v(G)}$ divise $a_{\bk}$ dans $\Os_{B,b}$. La condition~($\textrm{U}_{W}$) nous assure qu'il existe $b_{\bk} \in \Bs(W)$ v\'erifiant les propri\'et\'es suivantes :
\begin{enumerate}[\it i)]
\item $a_{\bk} = \pi^{v(G)}\, b_{\bk}$ dans $\Bs(W)$ ;
\item $\|b_{\bk}\|_{W} \le C_{W}^{v(G)}\, \|a_{\bk}\|_{W}$.
\end{enumerate}
Nous avons 
$$\sum_{\bk\ge 0} \|b_{\bk}\|_{W}\, \bt^\bk \le C_{W}^{v(G)}\, \sum_{\bk\ge 0} \|a_{\bk}\|_{U}\, \bt^\bk < +\infty.$$
Par cons\'equent, la s\'erie $\sum_{\bk\ge 0} b_{\bk}\, \bT^\bk$ d\'efinit un \'el\'ement $H$ de $\Bs(W)\of{\la}{|\bT|\le \bt}{\ra}$. Il v\'erifie bien $G = \pi^{v(G)} H$ et $H(b)\ne 0$.
\end{proof}

\begin{thm}\label{noetherienavd}
Supposons que l'anneau local $\Os_{B,b}$ est un anneau de valuation discr\`ete v\'erifiant la condition (U). Alors, l'anneau local $L_{b}$ est noeth\'erien.
\end{thm}
\begin{proof}
Nous allons proc\'eder par r\'ecurrence sur $n$. Si $n=0$, nous avons $L_{b}\simeq \Os_{B,b}$ et le r\'esultat est vrai.

Supposons, \`a pr\'esent, que le r\'esultat soit vrai pour $L'_{b}$. Soit $I$ un id\'eal de $L_{b}$. L'id\'eal nul \'etant de type fini, nous pouvons supposer que $I\ne (0)$. Notons 
$$v(I)=\min\{v(G),\, G\in I\}.$$
D'apr\`es le lemme \ref{lemU}, il existe un id\'eal $J$ de $L_{b}$ v\'erifiant les propri\'et\'es suivantes :
\begin{enumerate}[\it i)]
\item $I = \pi^{v(I)} J$ ;
\item l'id\'eal $J$ contient un \'el\'ement $G$ v\'erifiant $G(b) \ne 0$ dans $\Hs(b)[\![\bT]\!]$.
\end{enumerate} 
Nous pouvons alors utiliser le m\^eme raisonnement que dans la preuve du th\'eor\`eme \ref{noetheriencorps} pour montrer que l'id\'eal $J$ est de type fini. Il en est donc de m\^eme pour l'id\'eal $I$.
\end{proof}

\begin{thm}\label{factoriel}
Supposons que l'anneau local $\Os_{B,b}$ est un corps ou un anneau de valuation discr\`ete v\'erifiant la condition~(U). Alors, l'anneau local $L_{b}$ est factoriel.
\end{thm}
\begin{proof}
Il nous suffit de reprendre la structure des raisonnements pr\'ec\'edents en utilisant, cette fois-ci, le th\'eor\`eme de pr\'eparation de {Weierstra\ss} \ref{preparation}, joint au lemme \ref{auto}, et le th\'eor\`eme de {Gau\ss}. 
\end{proof}

Nous pouvons, en fait, obtenir un r\'esultat plus fort et d\'emontrer, sous les m\^emes hypoth\`eses, que l'anneau local $L_{b}$ est r\'egulier.

\begin{thm}\label{regulier}
Supposons que l'anneau local $\Os_{B,b}$ est un corps ou un anneau de valuation discr\`ete v\'erifiant la condition~(U). Alors, l'anneau $L_{b}$ est un anneau local r\'egulier de dimension \'egale \`a $\dim(\Os_{B,b})+n$.
\end{thm}
\begin{proof}
Rappelons que nous notons $\m = (\m_{b},T_{1},\ldots,T_{n})$ l'id\'eal maximal de $L_{b}$ et que nous avons 
$$\kappa(b) = \Os_{B,b}/\m_{b} \simeq L_{b}/\m.$$

Supposons, tout d'abord, que $\Os_{B,b}$ est un corps. Nous avons \mbox{$\m=(T_{1},\ldots,T_{n})$}, $\Os_{B,b}=\kappa(b)$ et $\dim(\Os_{B,b})=0$. La suite 
$$(0) \subset (T_{1}) \subset \cdots \subset (T_{1},\ldots,T_{n})$$
est une suite strictement croissante d'id\'eaux premiers de $L_{b}$. On en d\'eduit que 
$$\dim(L_{b})\ge n.$$

Montrons, \`a pr\'esent, que la famille $(T_{1},\ldots,T_{n})$ engendre le $\kappa(b)$-espace vectoriel $\m/\m^2$. Soit $G\in\m$. Par d\'efinition de $\m$, il existe $G_{1},\ldots,G_{n} \in L_{b}$ tels que 
$$G = \sum_{i=1}^n T_{i}\, G_{i} \textrm{ dans } L_{b}.$$
Quel que soit $i\in\cn{1}{n}$, il existe $h_{i}\in\Os_{B,b}$, $H_{i,1},\ldots,H_{i,n} \in \Os_{B,b}[\![\bT]\!]$ tels que 
$$G_{i} = h_{i} + \sum_{j=1}^n T_{j}\, H_{i,j} \textrm{ dans } \Os_{B,b}[\![\bT]\!].$$
D'apr\`es le lemme \ref{formelLb}, cette d\'ecomposition vaut encore dans $L_{b}$. Par cons\'equent, nous avons
$$G =  \sum_{i=1}^n h_{i}\, T_{i} + \sum_{1\le i,j\le n} T_{i}\,T_{j}\,H_{i,j} \textrm{ dans } L_{b}.$$
Or, quels que soient $i,j\in\cn{1}{n}$, nous avons $T_{i}\,T_{j}\in \m^2$. On en d\'eduit que 
$$G =  \sum_{i=1}^n h_{i}\, T_{i} \textrm{ dans } \m/\m^2.$$
Nous avons bien montr\'e que la famille $(T_{1},\ldots,T_{n})$ engendre le $\kappa(b)$-espace vectoriel $\m/\m^2$.

Comme tout anneau local noeth\'erien, l'anneau $L_{b}$ v\'erifie
$$\dim(L_{b}) \le \dim_{\kappa(b)}(\m/\m^2)\le n.$$ 
Finalement, nous avons donc
$$\dim(L_{b}) = \dim_{\kappa(b)}(\m/\m^2)=n.$$
On en d\'eduit que l'anneau $L_{b}$ est un anneau local r\'egulier de dimension $n$.\\

Supposons, \`a pr\'esent, que $\Os_{B,b}$ est un anneau de valuation discr\`ete v\'erifiant la condition U. Nous avons alors $\dim(\Os_{B,b})=1$. Soit $\pi$ une uniformisante de~$\Os_{B,b}$.  La suite 
$$(0) \subset (\pi) \subset (\pi,T_{1}) \subset \cdots \subset (\pi,T_{1},\ldots,T_{n})$$
est une suite strictement croissante d'id\'eaux premiers de $L_{b}$. Observons que pour montrer que ce sont des id\'eaux premiers, il faut faire appel \`a la condition U et, plus pr\'ecis\'ement, au lemme \ref{lemU}. Nous avons montr\'e que 
$$\dim(L_{b})\ge n+1.$$

Montrons, \`a pr\'esent, que la famille $(\pi,T_{1},\ldots,T_{n})$ engendre le $\kappa(b)$-espace vectoriel $\m/\m^2$. Soit $G\in\m$. Par d\'efinition de $\m$, il existe $G_{0},\ldots,G_{n} \in L_{b}$ tels que 
$$G = \pi\, G_{0} + \sum_{i=1}^n T_{i}\, G_{i} \textrm{ dans } L_{b}.$$
Par le m\^eme raisonnement que dans le cas des corps, on montre qu'il existe $h_{1},\ldots,h_{n} \in \Os_{X,x}$ tels que 
$$\sum_{i=1}^n T_{i}\, G_{i} =  \sum_{i=1}^n h_{i}\, T_{i} \textrm{ dans } \m/\m^2.$$
En utilisant de nouveau le lemme \ref{formelLb}, on montre qu'il existe $h_{0}\in\Os_{B,b}$, $H_{0,1},\ldots,H_{0,n} \in L_{b}$ tels que 
$$G_{0} = h_{0} + \sum_{j=1}^n T_{j}\, H_{0,j} \textrm{ dans } L_{b}.$$
Par cons\'equent, nous avons
$$\pi\, G_{0} = \pi\, h_{0} +  \sum_{j=1}^n \pi\, T_{j}\,H_{0,j} \textrm{ dans } L_{b}.$$
Or, quel que soit $j\in\cn{1}{n}$, nous avons $\pi\,T_{j}\in \m^2$. On en d\'eduit que 
$$G =  h_{0}\, \pi + \sum_{i=1}^n h_{i}\, T_{i} \textrm{ dans } \m/\m^2.$$
Nous avons bien montr\'e que la famille $(\pi,T_{1},\ldots,T_{n})$ engendre le $\kappa(b)$-espace vectoriel $\m/\m^2$.

L'anneau local noeth\'erien $L_{b}$ v\'erifie donc 
$$\dim(L_{b}) \le \dim_{\kappa(b)}(\m/\m^2) \le n+1.$$ 
On en d\'eduit que
$$\dim(L_{b}) = \dim_{\kappa(b)}(\m/\m^2)=n+1.$$
Finalement, l'anneau $L_{b}$ est un anneau local r\'egulier de dimension $n+1$.
\end{proof}


\section{Limites d'alg\`ebres de couronnes}

Soit $V$ une partie compacte de $B$. Pour $\bs\in\R_{+}^n$ et $\bt\in (\R_{+}^*)^n$, nous noterons~$\|.\|_{V,\bs,\bt}$ la norme sur l'anneau $\Bs(V)\of{\la}{\bs \le |\bT|\le \bt}{\ra}$ d\'efinie au num\'ero \ref{algglob}.
\newcounter{nVst}\setcounter{nVst}{\thepage}

Soit $b$ un point de $B$. Soit $\br=(r_{1},\ldots,r_{n}) \in  (\R_{+}^*)^n$ tel que la famille~$(r_{1},\ldots,r_{n})$ soit libre dans l'espace vectoriel~$\Q \otimes_{\Z} (\R_{+}^*/|\Hs(b)^*|)$. Nous noterons 
$$L_{b,\br} = \varinjlim_{V,\bs,\bt}\, \Bs(V)\of{\la}{\bs \le |\bT|\le \bt}{\ra},$$
\newcounter{Lbbr}\setcounter{Lbbr}{\thepage}
o\`u $V$ parcourt l'ensemble des voisinages compacts du point $b$ dans $B$, $\bs$ parcourt~$\prod_{i=1}^n \of{]}{0,r_{i}}{[}$ et $\bt$ parcourt~$\prod_{i=1}^n \of{]}{r_{i},+\infty}{[}$.

Comme pr\'ec\'edemment, lorsque l'anneau local~$\Os_{B,b}$ est un corps ou un anneau de valuation discr\`ete soumis \`a la condition~(U), nous pouvons mener une \'etude pr\'ecise de l'anneau~$L_{b,\br}$. Signalons que les r\'esultats s'obtiennent bien plus facilement que pr\'ec\'edemment. En particulier, nous n'aurons pas besoin de faire appel aux th\'eor\`emes de division et de pr\'eparation de Weierstra{\ss}. Nous commen\c{c}ons par \'enoncer un lemme qui g\'en\'eralise, en un certain sens, l'in\'egalit\'e ultram\'etrique.

\begin{lem}\label{inegpum}
Soit $k$ un corps muni d'une valeur absolue $|.|$ v\'erifiant l'in\'egalit\'e suivante : quels que soient les \'el\'ements $x$ et $y$ de $k$, on a 
$$|x+y|\le 2^\lambda \max(|x|,|y|).$$
Soient $n\in\N$ et $x_{0},\ldots,x_{n}\in k$. Alors on a 
$$\left|\sum_{i=0}^n x_{i}\right| \le 2^{n\lambda}\max_{0\le i\le n} (|x_{i}|).$$
Si l'on suppose que, quel que soit $i\in\cn{1}{n}$, on a $|x_{i}|<2^{-n\lambda} |x_{0}|$, alors on a
$$\left|\sum_{i=0}^n x_{i}\right| \ge 2^{-n\lambda} |x_{0}|.$$
\end{lem}
\begin{proof}
La premi\`ere in\'egalit\'e s'obtient facilement par r\'ecurrence. D\'emontrons la seconde. Supposons donc que, quel que soit $i\in\cn{1}{n}$, on a $|x_{i}|<2^{-n\lambda} |x_{0}|$. Alors 
$${\renewcommand{\arraystretch}{1.5}\begin{array}{rcl}
|x_{0}| & = & \disp \left| \sum_{i=0}^n x_{i} - x_{n} - \cdots - x_{1} \right|\\
& \le & \disp 2^{n\lambda}\, \max\left( \left| \sum_{i=0}^n x_{i} \right|, |x_{n}|,\ldots, |x_{1}| \right),
\end{array}}$$
d'apr\`es la premi\`ere in\'egalit\'e. Supposons, par l'absurde, qu'il existe $i\in\cn{1}{n}$ tel que 
$$ \max\left( \left| \sum_{i=0}^n x_{i} \right|, |x_{n}|,\ldots, |x_{1}| \right) = |x_{i}|.$$
Nous obtenons alors
$$|x_{0}| \le 2^{n\lambda} \, |x_{i}| < |x_{0}|,$$
ce qui est impossible. Par cons\'equent, nous avons 
$$ \max\left( \left| \sum_{i=0}^n x_{i} \right|, |x_{n}|,\ldots, |x_{1}| \right) =  \left| \sum_{i=0}^n x_{i} \right|.$$
On en d\'eduit la seconde in\'egalit\'e.
\end{proof}

\begin{thm}\label{corpscorps}
Supposons que l'anneau local $\Os_{B,b}$ est un corps. Alors l'anneau~$L_{b,\br}$ est un corps.
\end{thm}
\begin{proof}

Soit~$f$ un \'el\'ement non nul de l'anneau~$L_{b,\br}$. Il nous suffit de montrer que cet \'el\'ement est inversible. Il existe un voisinage compact~$V$ de~$b$ dans~$B$, des \'el\'ements~$\bs$ et~$\bt$ de~$\R_{+}^n$ v\'erifiant~$\bs < \br$ et~$\bt > \br$ tels que
$$f \in \Bs(V)\of{\la}{\bs \le |\bT|\le \bt}{\ra}.$$
Dans ce dernier anneau, la fonction~$f$ poss\`ede une \'ecriture sous la forme
$$f = \sum_{\bk\in\Z^n} a_{\bk}\, \bT^{\bk},$$
o\`u, quel que soit~$\bk\in\Z^n$, nous avons~$a_{\bk} \in \Bs(V)$ et la famille~$(\|a_{\bk}\|_{V}\, \bmax(\bs^\bk,\bt^\bk))_{\bk\in\Z^n}$ est sommable.

Les conditions impos\'ees au $n$-uplet~$\br$ nous assurent qu'il existe un \'el\'ement~$\bk_{0}$ de~$\Z^n$ tel que, quel que soit~$\bk\ne\bk_{0}$, on ait
$$|a_{\bk_{0}}(b)|\, \bmax(\bs^{\bk_{0}},\bt^{\bk_{0}}) > |a_{\bk}(b)|\, \bmax(\bs^\bk,\bt^\bk).$$
En utilisant le fait que la famille $(\|a_{\bk}\|_{V}\, \bmax(\bs^\bk,\bt^\bk))_{\bk\in\Z^n}$ est sommable, on en d\'eduit qu'il existe~$u,v\in\R$ tels que, quel que soit~$\bk\ne\bk_{0}$, on ait m\^eme
$$|a_{\bk_{0}}(b)|\, \bmin(\bs^{\bk_{0}},\bt^{\bk_{0}}) > v > u > |a_{\bk}(b)|\, \bmax(\bs^\bk,\bt^\bk).$$
Il existe un voisinage $E$ de~$-\infty$ dans~$\Z^n\setminus\{\bk_{0}\}$ tel que
$$\sum_{\bk\in E} \|a_{\bk}\|_{V}\, \bmax(\bs^\bk,\bt^\bk) \le u.$$
De m\^eme, il existe un voisinage $F$ de~$+\infty$ dans~$\Z^n\setminus (E\cup \{\bk_{0}\})$ tel que
$$\sum_{\bk\in F} \|a_{\bk}\|_{V}\, \bmax(\bs^\bk,\bt^\bk) \le u.$$
La partie $G=\Z^n\setminus(E\cup F \cup\{\bk_{0}\})$ ne contient qu'un nombre fini de termes. On en d\'eduit qu'il existe deux \'el\'ements~$\bs_{0}$ et~$\bt_{0}$ de~$(\R_{+}^*)^n$ v\'erifiant~$\bs\le\bs_{0}<\br$ et~$\br<\bt_{0}\le\bt$ tels que l'on ait
$$|a_{\bk_{0}}(b)|\, \bmin(\bs_{0}^\bk,\bt_{0}^\bk) > v$$
et, quel que soit~$\bk\in G$, 
$$|a_{\bk}(b)|\, \bmax(\bs_{0}^\bk,\bt_{0}^\bk) < u.$$
D\'efinissons deux voisinages compacts du point~$b$ dans~$V$ par
$$W_{0} = \left\{ c\in V\, \big|\, \forall \bk\in G,\, |a_{\bk}(b)|\, \bmax(\bs_{0}^\bk,\bt_{0}^\bk) \le u \right\}$$
et
$$W_{1} = \left\{ c\in V\, \big|\, |a_{\bk_{0}}(c)|\, \bmin(\bs_{0}^{\bk_{0}},\bt_{0}^{\bk_{0}}) \ge v \right\}.$$
Il existe un \'el\'ement~$\lambda$ de l'intervalle~$\of{]}{0,1}{]}$ v\'erifiant
$$2^{(c+2)\lambda}\, u <v.$$
Les conditions que nous avons impos\'ees sur~$\br$ imposent au corps valu\'e~$\Hs(b)$ d'\^etre ultram\'etrique. En particulier, nous avons~$|2(b)|\le 1$. Par cons\'equent, la partie
$$W_{2} = \left\{ c\in V\, \big|\, |2(c)| \le 2^{\lambda} \right\}$$
est un voisinage compact de~$b$ dans~$V$.
Choisissons un voisinage compact rationnel~$W$ de~$b$ contenu dans~$W_{0}\cap W_{1} \cap W_{2}$. Nous allons montrer que la fonction~$f$ est inversible dans l'anneau~$\Bs(W)\of{\la}{\bs_{0} \le |\bT| \le \bt_{0}}{\ra}$. Notons 
$$D = \pi^{-1}(W)\cap \overline{C}(\bs_{0},\bt_{0}).$$
En utilisant le fait que~$\Bs(W)=W$ et le lemme~$\ref{spectrecouronne}$, on montre que
$$\Ms(\Bs(W)\of{\la}{\bs_{0} \le |\bT| \le \bt_{0}}{\ra}) = D.$$
D'apr\`es \cite{rouge}, corollaire 1.2.4, pour montrer que la fonction~$f$ est inversible dans l'anneau~$\Bs(W)\of{\la}{\bs_{0} \le |\bT| \le \bt_{0}}{\ra}$, il suffit de montrer qu'elle ne s'annule par sur son spectre analytique~$D$. Soit~$y$ un point de~$D$. Notons~$c$ son projet\'e sur~$B$. C'est un \'el\'ement de~$W$. Nous avons
$${\renewcommand{\arraystretch}{2}
\begin{array}{rcl}
|f(y)| & = & \disp \left| \sum_{\bk\in\Z^n} a_{\bk}(c)\, \bT(y)^\bk \right|\\
&=& \disp \left| a_{\bk_{0}}(c)\, \bT(y)^{\bk_{0}} + \sum_{\bk\in E} a_{\bk}(c)\, \bT(y)^\bk +  \sum_{\bk\in F} a_{\bk}(c)\, \bT(y)^\bk \right.\\
&& \disp + \left. \sum_{\bk\in G} a_{\bk}(c)\, \bT(y)^\bk \right|.
\end{array}
}$$
\'Ecrivons l'expression \`a l'int\'erieur de la valeur absolue comme une somme de~\mbox{$3+\sharp G$} termes. \`A l'exception du premier, chacun de ces termes $g$ v\'erifie
$$|g| \le u < 2^{-(\sharp G +2)\lambda}\, v \le |a_{\bk_{0}}(c)|\, |\bT(c)^{\bk_{0}}|.$$
D'apr\`es le lemme~\ref{inegpum}, nous avons donc
$$|f(y)|\ge 2^{-(\sharp G +2)\lambda}\,  |a_{\bk_{0}}(c)|\, |\bT(c)^{\bk_{0}}| >0.$$
On en d\'eduit le r\'esultat.  
\end{proof}

Venons-en, \`a pr\'esent, au cas o\`u l'anneau local $\Os_{B,b}$ est un anneau de valuation discr\`ete v\'erifiant la condition~(U) de la d\'efinition~\ref{condU}. Soit~$\pi$ une uniformisante de~$\Os_{B,b}$ et~$v_{\pi}$ la valuation associ\'ee. Nous disposons d'un r\'esultat analogue \`a celui du lemme~\ref{lemU}. Avant de l'\'enoncer, d\'efinissons une application~$v$ de~$L_{b,\br}$ dans~$\N\cup\{+\infty\}$. Soit~$f$ un \'el\'ement de~$L_{b,\br}$. Il existe un voisinage compact~$V$ de~$b$ dans~$B$, des \'el\'ements~$\bs$ et~$\bt$ de~$\R_{+}^n$ v\'erifiant~$\bs < \br$ et~$\bt > \br$ tels que
$$f \in \Bs(V)\of{\la}{\bs \le |\bT|\le \bt}{\ra}.$$
Dans ce dernier anneau, la fonction~$f$ poss\`ede une \'ecriture sous la forme
$$f = \sum_{\bk\in\Z^n} a_{\bk}\, \bT^{\bk},$$
o\`u, quel que soit~$\bk\in\Z^n$, nous avons~$a_{\bk} \in \Bs(V)$ et la famille~$(\|a_{\bk}\|_{V}\, \bmax(\bs^\bk,\bt^\bk))_{\bk\in\Z^n}$ est sommable. Posons
$$v(f)=  \min\{ v_{\pi}(a_{\bk}),\, \bk \in \Z^n \} \in \N\cup\{+\infty\}.$$
Cette quantit\'e ne d\'epend pas du repr\'esentant de~$f$ choisi.

\begin{lem}\label{lemU3}
Supposons que l'anneau local $\Os_{B,b}$ est un anneau de valuation discr\`ete v\'erifiant la condition (U). Soit $\pi$ une uniformisante de $\Os_{B,b}$ et notons $v_{\pi}$ la valuation associ\'ee. Soit $f$ un \'el\'ement non nul de $L_{b,\br}\setminus\{0\}$. Alors, il existe une fonction $g$ de $L_{b,\br}$ v\'erifiant les propri\'et\'es suivantes :
\begin{enumerate}[\it i)]
\item $v(g)=0$ ; 
\item $f = \pi^{v(f)} g$ dans $L_{b,\br}$.
\end{enumerate}
\end{lem}

Nous en d\'eduisons le th\'eor\`eme suivant.

\begin{thm}\label{avdavd}
Supposons que l'anneau local $\Os_{B,b}$ est un anneau de valuation discr\`ete v\'erifiant la condition~(U). Alors l'anneau~$L_{b,\br}$ est un anneau de valuation discr\`ete, de valuation~$v$ et d'id\'eal maximal~$\m_{b}\,L_{b,\br}$.
\end{thm}
\begin{proof}
On v\'erifie directement sur la d\'efinition de l'application~$v$ que les deux propri\'et\'es suivantes sont v\'erifi\'ees : quels que soient~$f$ et~$g$ dans~$L_{b,\br}$, nous avons
\begin{enumerate}[\it i)]
\item $v(f+g) \ge \min(v(f),v(g))$ ;
\item $v(fg) = v(f) + v(g)$.
\end{enumerate}

En outre, la condition~(U) assure que nous avons~$v(f)=+\infty$ si, et seulement si, la fonction~$f$ est nulle. De cette propri\'et\'e, jointe \`a la propri\'et\'e~{\it ii)}, on d\'eduit que l'anneau~$L_{b,\br}$ est int\`egre. Notons~$F$ son corps des fractions. L'application~$v$ se prolonge en un morphisme surjectif de~$F^*$ dans~$\Z$ qui v\'erifie encore la propri\'et\'e~{\it i)}. C'est donc une valuation discr\`ete. 

Pour conclure, il nous reste \`a montrer que nous avons les deux \'egalit\'es suivantes :
\begin{enumerate}[a)]
\item $L_{b,\br} = \{f\in F\, |\, v(f) \ge 0\}$ ;
\item $\m_{b}\,L_{b,\br} = \{f\in F\, |\, v(f) > 0\}.$
\end{enumerate}
L'\'egalit\'e b) se d\'eduit de l'\'egalit\'e a) en utilisant la condition~(U). En outre, en utilisant le lemme~\ref{lemU3}, on se ram\`ene \`a montrer que tout \'el\'ement de~$L_{b,\br}$ de valuation nulle est inversible dans~$L_{b,\br}$. Soit~$f$ un \'el\'ement de~$L_{b,\br}$ tel que~\mbox{$v(f)=0$}. Il existe un voisinage compact~$V$ de~$b$ dans~$B$, des \'el\'ements~$\bs$ et~$\bt$ de~$\R_{+}^n$ v\'erifiant~$\bs < \br$ et~$\bt > \br$ tels que
$$f \in \Bs(V)\of{\la}{\bs \le |\bT|\le \bt}{\ra}.$$
Dans ce dernier anneau, la fonction~$f$ poss\`ede une \'ecriture sous la forme
$$f = \sum_{\bk\in\Z^n} a_{\bk}\, \bT^{\bk},$$
o\`u, quel que soit~$\bk\in\Z^n$, nous avons~$a_{\bk} \in \Bs(V)$ et la famille~$(\|a_{\bk}\|_{V}\, \bmax(\bs^\bk,\bt^\bk))_{\bk\in\Z^n}$ est sommable. Puisque~$v(f)=0$, la famille~$(|a_{\bk}(b)|)_{\bk\in\Z^n}$ n'est pas nulle. Les conditions impos\'ees au $n$-uplet~$\br$ nous assurent alors qu'il existe un \'el\'ement~$\bk_{0}$ de~$\Z^n$ tel que, quel que soit~$\bk\ne\bk_{0}$, on ait
$$|a_{\bk_{0}}(b)|\, \bmax(\bs^{\bk_{0}},\bt^{\bk_{0}}) > |a_{\bk}(b)|\, \bmax(\bs^\bk,\bt^\bk).$$
On en utilisant le m\^eme raisonnement que dans la preuve du th\'eor\`eme~\ref{corpscorps}, on montre que la fonction~$f$ est inversible dans l'anneau~$L_{b,\br}$.
\end{proof}

\section{Exemples d'anneaux locaux}\label{exannloc}

Il est possible d'exhiber des bases de voisinages explicites de certains points de l'espace affine. Ces r\'esultats nous seront, par la suite, tr\`es utiles pour \'etudier les anneaux locaux en ces points. Commen\c{c}ons par nous int\'eresser \`a des parties compactes plus g\'en\'erales.

\begin{lem}\label{lemvois}\index{Voisinages d'un compact}
Soient~$U$ un ouvert de~$B$, $Y$ un ouvert de~$X_{U}$, $p$ un entier et $f_{1},\ldots,f_{p}$ des \'el\'ements de~$\Os_{X}(Y)$. Pour toute partie compacte~$V$ de~$U$ et tous \'el\'ements $\bs=(s_{1},\ldots,s_{p})$ et $\bt=(t_{1},\ldots,t_{p})$ de~$\R_{+}^p$, nous posons
$$M_{V}(\bs,\bt) = \{y\in Y\cap X_{V}\, |\, \forall i\in\cn{1}{p},\, s_{i} \le |f_{i}(y)|\le t_{i}\}.$$
Nous supposerons que toutes ces parties sont compactes.

Soient~$V$ une partie compacte de~$U$ et~$\bs$ et~$\bt$ deux \'el\'ements de~$\R^p$. Soit~$N$ un voisinage du compact $M_{V}(\bs,\bt)$ dans~$Y$. Il existe un voisinage compact~$V'$ de~$V$ dans~$U$ et deux \'el\'ements~$\bs'$ et~$\bt'$ de~$\R_{+}^p$ v\'erifiant les in\'egalit\'es $\bs'\prec \bs$ et $\bt' > \bt$ tels que l'on ait l'inclusion
$$M_{V'}(\bs',\bt') \subset N.$$
\end{lem}
\begin{proof}
Posons $M=M_{V}(\bs,\bt)$. Soient~$V_{0}$ un voisinage compact de~$V$ dans~$U$ et~$\bs_{0}$ et~$\bt_{0}$ deux \'el\'ements de~$\R^p$ v\'erifiant les in\'egalit\'es $\bs'\prec \bs$ et $\bt' > \bt$. La partie compacte~$M_{0} = M_{V_{0}}(\bs_{0},\bt_{0})$ est alors un voisinage compact de~$M$ dans~$Y$. Sans perdre en g\'en\'eralit\'e, nous pouvons supposer que~$N$ est un voisinage ouvert de~$M$ dans~$M_{0}$. 

Posons $M_{1}=M_{V_{0}}(\bs,\bt)$. La partie $N\cap M_{1}$ est un voisinage ouvert de~$M$ dans~$M_{1}$. Son compl\'ementaire~$S_{1}$ est une partie compacte. Puisque~$M_{1}\cap X_{V}=M$, le compact~$S_{1}$ ne coupe pas~$X_{V}$. Par cons\'equent, le compact~$\pi(S_{1})$ ne coupe pas~$V$. Choisissons un voisinage compact~$V'$ de~$V$ dans~$V_{0}$ contenu dans le compl\'ementaire de~$\pi(S_{1})$. Nous avons alors 
$$M_{V'}(\bs,\bt) = M_{1} \cap X_{V'} \subset M_{1}\cap N \subset N.$$ 

Posons $M_{2}=M_{V'}(\bt,\bt_{0})$. La partie $N\cap M_{2}$ est un voisinage ouvert de~$M$ dans~$M_{2}$. Son compl\'ementaire~$S_{2}$ est une partie compacte. La fonction
$$\max_{1\le i\le p} (|f_{i}|-t_{i})$$
atteint son minimum~$m$ sur~$S_{2}$. Puisque~$S_{2}$ est disjoint de~$M$, le nombre r\'eel~$m$ est strictement positif. Pour tout \'el\'ement~$i$ de~$\cn{1}{p}$, choisissons un \'el\'ement~$t'_{i}$ de l'intervalle $\of{]}{t_{i},t_{i}+m}{[}$. Posons $\bt'=(t'_{1},\ldots,t'_{p})$. Nous avons alors $\bt'>\bt$ et
$$M_{V'}(\bt,\bt') = \subset M_{2}\cap N \subset N.$$ 

Nous montrons de m\^eme qu'il existe un \'el\'ement~$\bs'$ de~$\R_{+}^p$ v\'erifiant $\bs'\prec\bs$ tel que 
$$M_{V'}(\bs',\bs) = \subset N.$$

On en d\'eduit que
$$M_{V'}(\bs',\bt')\subset N,$$
ce qui d\'emontre le r\'esultat.
\end{proof}

Nous allons maintenant appliquer ce r\'esultat afin d'obtenir une description explicite de syst\`emes fondamentaux de voisinages pour certains points.

\begin{defi}\label{defdeploye}\index{Point!deploye@d\'eploy\'e}
Soient~$b$ un point de~$B$,  $\alpha_{1},\ldots,\alpha_{n}$ des \'el\'ements de~$\Os_{B,b}$ et $r_{1},\ldots,r_{n}$ des \'el\'ements de~$\R_{+}$. Notons $I$ l'ensemble des \'el\'ements~$i$ de~$\cn{1}{n}$ tels que $r_{i}\ne 0$. Supposons que la famille $(r_{i})_{i\in I}$ est libre dans l'espace vectoriel $\Q\otimes_{\Z} (\R_{+}^*/|\Hs(b)^*|)$. Il existe alors un unique point~$x$ de la fibre~$X_{b}$ qui v\'erifie les in\'egalit\'es suivantes :
$$\forall i\in\cn{1}{n}, |(T_{i}-\alpha_{i})(x)|=r_{i}.$$
Un tel point est dit {\bf d\'eploy\'e}.
\end{defi}

Soient $b$ un point de $B$ et ${\alphab} =(\alpha_{1},\ldots,\alpha_{n})$ un \'el\'ement de~$\Os_{B,b}^n$. Soit $B_{0}$ un voisinage de $b$ dans $B$ sur lequel les fonctions $\alpha_{1},\ldots,\alpha_{n}$ sont d\'efinies.

Soient~$I$ une partie de~$\cn{1}{n}$ et $(r_{i})_{i\in I}$ une famille de $\R_{+}^*$ dont l'image dans l'espace vectoriel $\Q\otimes_{\Z} (\R_{+}^*/|\Hs(b)^*|)$ est libre. Notons $J=\cn{1}{n}\setminus I$ et, pour $i\in J$, posons $r_{i}=0$. Posons encore $\br=(r_{1},\ldots,r_{n})$. Notons~$x$ l'unique point de la fibre~$X_{b}$ qui v\'erifie 
$$\forall i\in\cn{1}{n}, |(T_{i}-\alpha_{i})(x)|=r_{i}.$$

\begin{prop}\label{voisdep}\index{Voisinages d'un point!deploye@d\'eploy\'e}
Soit~$U$ un voisinage du point~$x$ dans~$X$. Pour tout \'el\'ement~$i$ de~$J$, posons $s_{i}=0$. Il existe un voisinage~$V$ du point~$b$ dans~$B_{0}$, pour tout \'el\'ement~$i$ de~$I$, un \'el\'ement~$s_{i}$ de $\of{]}{0,r_{i}}{[}$ et, pour tout \'el\'ement~$i$ de $\cn{1}{n}$, un \'el\'ement~$t_{i}$ de $\of{]}{r_{i},+\infty}{[}$ tels que l'on ait l'inclusion
$$\left\{y\in X_{V}\, \big|\, \forall i\in\cn{1}{n},\, s_{i}\le |(T_{i}-\alpha_{i})(y)| \le t_{i}\right\} \subset U.$$

\end{prop}
\begin{proof}
D'apr\`es le corollaire \ref{partiecompacte}, pour toute partie compacte~$V$ de~$B_{0}$ et tous \'el\'ements $s_{1},\ldots,s_{n},t_{1},\ldots,t_{n}$ de~$\R_{+}$, la partie de~$X$ d\'efinie par
$$\{y\in X_{V}\, |\, \forall i\in\cn{1}{p},\, s_{i} \le |(T_{i}-\alpha)(y)|\le t_{i}\}$$
est compacte. Le r\'esultat d\'ecoule alors du lemme \ref{lemvois}.
\end{proof}

Nous allons, \`a pr\'esent, pr\'eciser ce r\'esultat. \`A cet effet, nous allons construire une application~$\sigma_{\alphab,\br}$ de~$B_{0}$ dans
$$W = \left\{y\in X_{B_{0}}\, \big|\, \forall i\in\cn{1}{n},\,  |(T_{i}-\alpha_{i})(y)| = r_{i}\right\}$$
qui soit une section du morphisme~$\pi$ au-dessus de~$B_{0}$. 

Soit~$c$ un point de~$B_{0}$. Si le point~$c$ est associ\'e \`a une valeur absolue ultram\'etrique, nous d\'efinissons~$\sigma_{\alphab,\br}(c)$ comme le point associ\'e \`a la semi-norme multiplicative
$${\renewcommand{\arraystretch}{1.5}\begin{array}{ccc}
\As[T_{1},\ldots,T_{n}] & \to & \R_{+}\\
\disp \sum_{\bk\ge 0} a_{\bk}\, \prod_{i=1}^n (T_{i}-\alpha_{i})^{k_{i}} & \mapsto & \disp \max_{\bk\ge 0} \left(|a_{\bk}(c)| \prod_{i=1}^n r_{i}^{k_{i}}\right)
\end{array}}.$$
Si le point~$c$ est associ\'e \`a une valeur absolue archim\'edienne, alors le corps r\'esiduel compl\'et\'e~$\Hs(c)$ est~$\R$ ou~$\C$ muni de la valeur absolue~$|.|_{\infty}^\eps$, avec~$\eps\in\of{]}{0,1}{]}$. Nous d\'efinissons~$\sigma_{\alphab,\br}(c)$ comme le point $(\alpha_{1}+r_{1}^{1/\eps},\ldots,\alpha_{n}+r_{n}^{1/\eps})$ de la fibre~$X_{c}$, autrement dit, comme le point associ\'e \`a la semi-norme multiplicative
$$\begin{array}{ccc}
\As[T_{1},\ldots,T_{n}] & \to & \R_{+}\\
\disp \sum_{\bk\ge 0} a_{\bk}\, \prod_{i=1}^n (T_{i}-\alpha_{i})^{k_{i}} & \mapsto & \disp \left|\sum_{\bk\ge 0} a_{\bk}(c) \prod_{i=1}^n r_{i}^{k_{i}/\eps} \right|_{\infty}^\eps
\end{array}.$$

\begin{lem}
L'application
$$\sigma_{\alphab,\br} : B_{0} \to W$$
est une section continue du morphisme~$\pi$ au-dessus de~$B_{0}$.
\end{lem}
\begin{proof}
Le fait que l'application~$\sigma_{\alphab,\br}$ prenne ses valeurs dans~$W$ et soit une section de~$\pi$ est imm\'ediat. Int\'eressons-nous, maintenant, \`a sa continuit\'e. Rappelons que, par d\'efinition de la topologie de~$X$, l'application~$\sigma_{\alphab,\br}$ est continue si, et seulement si, pour tout \'el\'ement~$P$ de $\As[T_{1},\ldots,T_{n}]$, l'application
$$|P(.)| \circ \sigma_{\alphab,\br} : \begin{array}{ccc}
B_{0} & \to & \R_{+}\\
c & \mapsto & |P(\sigma_{\alphab,\br}(c))|
\end{array}$$
est continue.

Consid\'erons l'ouvert de~$B_{0}$ d\'efini par
$$B_{1}=\left\{c\in B_{0}\, \big |2(c)|<1 \right\}.$$
Chacun des points de cet ouvert est associ\'e \`a une valeur absolue ultram\'etrique. Par cons\'equent, pour tout \'el\'ement $P=\sum_{\bk\ge 0} a_{\bk}\, \prod_{i=1}^n (T_{i}-\alpha_{i})^{k_{i}}$ de $\As[T_{1},\ldots,T_{n}]$, nous avons
$$|P(\sigma_{\alphab,\br}(c))| = \max_{\bk\ge 0} \left(|a_{\bk}(c)| \prod_{i=1}^n r_{i}^{k_{i}}\right).$$
On en d\'eduit que l'application~$\sigma_{\alphab,\br}$ est continue sur~$B_{1}$.

Consid\'erons de m\^eme l'ouvert de~$B_{0}$ d\'efini par
$$B_{2}=\left\{c\in B_{0}\, \big |2(c)|>1 \right\}.$$
Chacun des points de cet ouvert est associ\'e \`a une valeur absolue archim\'edienne. Par cons\'equent, pour tout \'el\'ement $P=\sum_{\bk\ge 0} a_{\bk}\, \prod_{i=1}^n (T_{i}-\alpha_{i})^{k_{i}}$ de $\As[T_{1},\ldots,T_{n}]$, nous avons
$$|P(\sigma_{\alphab,\br}(c))| = \left|P\left(\alpha+r_{1}^{1/\eps},\ldots,\alpha+r_{n}^{1\eps}\right)(c)\right|.$$
On en d\'eduit que l'application~$\sigma_{\alphab,\br}$ est continue sur~$B_{2}$.

Si le point central~$a_{0}$ de~$B$ n'appartient pas \`a~$B_{0}$, alors $B_{0}=B_{1}\cup B_{2}$ et nous avons montr\'e que l'application~$\sigma_{\alphab,\br}$ est continue. Supposons, \`a pr\'esent, que le point~$a_{0}$ appartienne \`a~$B_{0}$. Par hypoth\`ese, l'image dans l'espace vectoriel $\Q\otimes_{\Z} (\R_{+}^*/|\Hs(b)^*|)$  de la famille $(r_{i})_{i\in I}$ de $\R_{+}^*$ est libre. Puisque $|\Hs(a_{0})^*|=\{1\}$ est contenu dans~$|\Hs(b)^*|$ son image est encore libre dans l'espace vectoriel $\Q\otimes_{\Z} (\R_{+}^*/|\Hs(a_{0})^*|)$. On en d\'eduit que le point~$\sigma_{\alphab,\br}(a_{0})$ est l'unique point du compact
$$\{y\in X_{0}\, \big|\, \forall i\in\cn{1}{n},\, |(T_{i}-\alpha_{i})(y)|=r_{i}\}.$$  
Soit~$U$ un voisinage du point~$\sigma_{\alphab,\br}(a_{0})$ dans~$X$. D'apr\`es la proposition \ref{voisdep}, il contient une partie de la forme
$$\left\{y\in X_{V}\, \big|\, \forall i\in\cn{1}{n},\, s_{i}\le |(T_{i}-\alpha_{i})(y)| \le t_{i}\right\},$$
o\`u~$V$ est un voisinage du point~$a_{0}$ dans~$B_{0}$, pour tout \'el\'ement~$i$ de~$J$, $s_{i}=0$, pour tout \'el\'ement~$i$ de~$I$, $s_{i}$ appartient \`a $\of{]}{0,r_{i}}{[}$ et, pour tout \'el\'ement~$i$ de $\cn{1}{n}$, $t_{i}$ appartient \`a $\of{]}{r_{i},+\infty}{[}$. En particulier, il contient la partie~$W\cap X_{V}$. Par cons\'equent, la partie $\sigma_{\alphab,\br}^{-1}(U)$ contient le voisinage~$V$ de~$a_{0}$ dans~$B_{0}$. On en d\'eduit que l'application~$\sigma_{\alphab,\br}$ est continue au voisinage du point~$a_{0}$. 

Nous avons, \`a pr\'esent, trait\'e le cas de tous les points de~$B_{0}$. Nous avons donc bien montr\'e que l'application~$\sigma_{\alphab,\br}$ est continue.
\end{proof}

Cette section nous permet d'obtenir des informations suppl\'ementaires sur les voisinages des points d\'eploy\'es des fibres.

\begin{cor}\label{voissectiondep}\index{Voisinages d'un point!deploye@d\'eploy\'e}
Soient~$b$ un point de~$B$, $x$ un point d\'eploy\'e de la fibre~$X_{b}$ et~$U$ un voisinage du point~$x$ dans~$X$. Il existe un voisinage~$V$ du point~$x$ dans~$U$ v\'erifiant les propri\'et\'es suivantes :
\begin{enumerate}[\it i)]
\item la projection $\pi(V)$ est un voisinage du point~$\pi(x)=b$ dans~$B$ ;
\item il existe une section continue~$\sigma$ du morphisme de projection $V \to \pi(V)$ ;
\item pour tout point~$b$ de~$\pi(V)$, la trace de la fibre~$X_{b}$ sur~$V$ est connexe par arcs.
\end{enumerate} 
\end{cor}
\begin{proof}
Ce r\'esultat d\'ecoule directement de la proposition et du lemme qui pr\'ec\`edent. Le point {\it iii)} est vrai car pour tout corps valu\'e complet~$k$, tous \'el\'ements $\alpha_{1},\ldots,\alpha_{n}$ de~$k$ et $s_{1},\ldots,s_{n},t_{1},\ldots,t_{n}$ de~$\R_{+}$, la partie de l'espace analytique~$\E{n}{k}$ d\'efinie par 
$$\left\{y\in\E{n}{k}\, \big|\, \forall i\in\cn{1}{n},\, s_{i}\le |(T-\alpha_{i})(y)|\le t_{i}\right\}$$
est connexe par arcs. 
\end{proof}

\begin{cor}\label{ouvertdep}\index{Ouverture au voisinage d'un point!deploye@d\'eploy\'e}
Soient~$b$ un point de~$B$ et~$x$ un point d\'eploy\'e de la fibre~$X_{b}$. Le morphisme $\pi$ est ouvert au point~$x$.
\end{cor}

\begin{cor}\label{cpadep}\index{Connexite par arcs au voisinage d'un point@Connexit\'e par arcs au voisinage d'un point!deploye@d\'eploy\'e}
Soient~$b$ un point de~$B$ et~$x$ un point d\'eploy\'e de la fibre~$X_{b}$. Si le point~$b$ de~$B$ poss\`ede un syst\`eme fondamental de voisinages connexes par arcs, alors il en est de m\^eme pour le point~$x$ de~$X$.
\end{cor}

Nous pouvons, \`a pr\'esent, d\'ecrire explicitement les anneaux locaux aux points d\'eploy\'es des fibres. Reprenons les notations du d\'ebut de ce num\'ero. Soient $b$ un point de $B$ et ${\alphab} =(\alpha_{1},\ldots,\alpha_{n})$ un \'el\'ement de~$\Os_{B,b}^n$. Soit $B_{0}$ un voisinage de $b$ dans $B$ sur lequel les fonctions $\alpha_{1},\ldots,\alpha_{n}$ sont d\'efinies.

Soient~$I$ une partie de~$\cn{1}{n}$ et $(r_{i})_{i\in I}$ une famille de $\R_{+}^*$ dont l'image dans l'espace vectoriel $\Q\otimes_{\Z} (\R_{+}^*/|\Hs(b)^*|)$ est libre. Notons $J=\cn{1}{n}\setminus I$ et, pour $i\in J$, posons $r_{i}=0$. Posons encore $\br=(r_{1},\ldots,r_{n})$. Notons~$x$ l'unique point de la fibre~$X_{b}$ qui v\'erifie 
$$\forall i\in\cn{1}{n}, |(T_{i}-\alpha_{i})(x)|=r_{i}.$$

\begin{thm}\label{anneaulocal}\index{Anneau local en un point!deploye@d\'eploy\'e}
Le morphisme $\As[\bT] \to \Os_{X,x}$ induit un isomorphisme
$$\varinjlim_{V,\bs,\bt}\, \Bs(V)\of{\la}{\bs \le |\bT-\alphab|\le \bt}{\ra} \xrightarrow{\sim} \Os_{X,x},$$
o\`u $V$ parcourt l'ensemble des voisinages de $b$ dans $B_{0}$, quel que soit $i\in J$, $s_{i}=0$ et $t_{i}$ parcourt $\R_{+}^*$, quel que soit $i\in I$, $s_{i}$ et $t_{i}$ parcourent respectivement $\of{]}{0,r_{i}}{[}$ et $\of{]}{r_{i},+\infty}{[}$.
\end{thm}
\begin{proof}
Quitte \`a remplacer l'anneau~$\As$ par~$\Bs(U)$, o\`u~$U$ d\'esigne un voisinage compact rationnel de~$b$ assez petit, nous pouvons supposer que~\mbox{$\alphab \in\As^n$}. Cette op\'eration est licite d'apr\`es le th\'eor\`eme~\ref{compactrationnel}. Quitte \`a appliquer la translation par le vecteur~$-\alphab$, qui est un automorphisme, nous pouvons supposer que~$\alphab =0$.

Soit $V$ un voisinage compact du point~$b$ dans~$\Ms(\As)$, $\bs$ un \'el\'ement de~$\R_{+}^n$ et~$\bt$ un \'el\'ement de $\left(\R_{+}^*\right)^n$ tels que $\bs\le\bt$. D'apr\`es la proposition \ref{spectrecouronne}, le morphisme naturel 
$\As[T] \to \Bs(\overline{C}_{V}(\bs,\bt))$
se prolonge en un morphisme 
$$\Bs(V)\of{\la}{\bs\le |\bT|\le \bt}{\ra} \to \Bs(\overline{C}_{V}(\bs,\bt)).$$ 
La proposition \ref{comparaisoncouronne} assure que ce morphisme est injectif. En utilisant la proposition~\ref{voisdep}, on en d\'eduit qu'il existe un morphisme injectif
$$\varphi : \varinjlim_{V,\bs,\bt}\, \Bs(V)\of{\la}{\bs \le|\bT|\le \bt}{\ra} \hookrightarrow \Os_{X,x},$$
o\`u~$V$ parcourt l'ensemble des voisinages compacts du point~$b$ dans~$B$ et~$\bs$ et~$\bt$ l'ensemble des \'el\'ements de~$\R_{+}^n$ qui v\'erifient $\bs \prec\br< \bt$.

Il nous reste \`a montrer que ce morphisme est surjectif. Soit~$f$ un \'el\'ement de~$\Os_{X,x}$. Par d\'efinition du faisceau structural, il existe un voisinage~$U$ du point~$x$ dans~$X$ sur lequel la fonction~$f$ est la limite uniforme d'une suite de fractions rationnelles~$(R_j)_{j\ge 0}$ \`a coefficients dans~$\As$ sans p\^oles sur~$U$. D'apr\`es la proposition~\ref{voisdep}, nous pouvons supposer que le voisinage~$U$ est de la forme 
$$U = \overline{C}_{V}(\bs,\bt),$$
o\`u~$V$ d\'esigne un voisinage compact rationnel du point~$b$ dans~$B$, et~$\bs$ et~$\bt$ deux \'el\'ements de~$\R_{+}^n$ qui v\'erifient $\bs \prec\br< \bt$. Le morphisme naturel 
$$\As[\bT] \to \Bs(V)\of{\la}{\bs \le |\bT|\le \bt}{\ra}$$
est injectif. D'apr\`es les propositions \ref{compactrationnelrelatif} et \ref{spectrecouronne}, ce morphisme induit un hom\'eomorphisme 
$$\Ms(\Bs(V)\of{\la}{\bs\le |\bT|\le \bt}{\ra}) \xrightarrow[]{\sim} U.$$
Soit~$P$ un \'el\'ement de $\As[\bT]$ qui ne s'annule en aucun point de~$U$. D'apr\`es \cite{rouge}, corollaire 1.2.4, l'image de~$P$ est inversible dans l'anneau $\Bs(V)\of{\la}{\bs\le|\bT|\le \bt}{\ra}$. On en d\'eduit que l'anneau~$\Ks(U)$ s'injecte dans $\Bs(V)\of{\la}{\bs\le |\bT|\le \bt}{\ra}$.

Soient~$\bu$ un \'el\'ement de~$\R_{+}^n$ tel que $\bs\prec \bu\prec \br$ et~$\bv$ un \'el\'ement de~$\left(\R_{+}^*\right)^n$ tel que $\br<\bv<\bt$. L'anneau~$\Ks(U)$ s'injecte encore dans l'anneau $\Bs(V)\of{\la}{\bu \le |\bT|\le \bv}{\ra}$. L'in\'egalit\'e sur les normes d\'emontr\'ee dans la proposition \ref{comparaisoncouronne} assure que la suite~$(R_j)_{j\ge 0}$ est une suite de Cauchy dans l'anneau $\Bs(V)\of{\la}{\bu \le |\bT|\le \bv}{\ra}$. Puisque ce dernier anneau est complet, la suite~$(R_j)_{j\ge 0}$ y converge et sa limite est envoy\'ee sur la fonction~$f$ par le morphisme~$\varphi$.

\end{proof}

\section{Hens\'elianit\'e}\label{henselianite}

Nous commen\c{c}ons par montrer que les anneaux locaux de l'espace affine analytique~$X$ au-dessus de $B$ sont hens\'eliens. Nous d\'ecrivons ensuite un cadre dans lequel cette propri\'et\'e peut d\'eboucher sur l'existence d'un isomorphisme local entre espaces analytiques.

\subsection{D\'emonstration}

\index{Henselianite@Hens\'elianit\'e|(}

\begin{prop}\label{Hensel}
Soit $x$ un point de $X$. L'anneau local $\Os_{X,x}$ est hens\'elien.
\end{prop}
\begin{proof}

Rappelons que nous notons $\kappa(x) = \Os_{X,x}/\m_{x}$. Soit $P(T)$ un polyn\^ome unitaire de $\Os_{X,x}[T]$ dont l'image dans $\kappa(x)[T]$ poss\`ede une racine simple $\alpha$. D'apr\`es \cite{alh}, chapitre VII, proposition 3, il nous suffit de montrer que $\alpha$ se rel\`eve en une racine de $P(T)$ dans $\Os_{X,x}$. 

Choisissons un \'el\'ement $f$ de $\Os_{X,x}$ relevant $\alpha$. Nous pouvons alors retraduire les hypoth\`eses sous la forme $P(f)(x)=0$ et $P'(f)(x)\ne 0$.

Soit $U$ un voisinage compact de $x$ dans $X$ tel que les coefficients du polyn\^ome~$P$ et l'\'el\'ement $f$ appartiennent \`a $\Bs(U)$. Quitte \`a restreindre $U$, nous pouvons supposer que la fonction $P'(f)$ y est inversible. Il existe un polyn\^ome \mbox{$Q(T_{1},T_{2})\in\Bs(U)[T_{1},T_{2}]$}, ind\'ependant de $f$, tel que, quel que soit $g\in\Bs(U)$, on ait
$${\renewcommand{\arraystretch}{1.5}\begin{array}{rcl}
P(f+P(f)g) &=& P(f) + P'(f)P(f)g + P(f)^2 g^2 Q(f,g)\\
&=& P'(f)P(f) \left( \frac{1}{P'(f)} + g + \frac{P(f)}{P'(f)}\, g^2 Q(f,g)\right).
\end{array}}$$

Notons $d\in\N$ le degr\'e du polyn\^ome $Q(f,T)$. Soit $t\in\of{]}{0,1}{[}$. Quitte \`a restreindre encore le voisinage $U$ de $x$, nous pouvons supposer que $t/(d+1)$ majore la norme uniforme sur $U$ de tous les coefficients du polyn\^ome 
$$R(T)= - \frac{P(f)}{P'(f)}\, T^2 Q(f,T).$$ 
On a alors
$$\begin{array}{rcl}
\forall g\in\Bs(U), \, \|R(g)\|_{U} &\le& \sum\limits_{i=2}^{d+2} \frac{t}{d+1}\,\|g\|_{U}^{i}\\
&\le& t\,  \max\left(\|g\|_{U}^{2},\|g\|_{U}^{d+2}\right).
\end{array}$$ 
En particulier, si $g\in\Bs(U)$ v\'erifie $\|g\|_{U}\le 1$, alors nous avons encore \mbox{$\|R(g)\|_{U} \le 1$}.

Quitte \`a diminuer $t$, nous pouvons supposer que 
$$t\, \max\left(  \left\| \frac{-1}{P'(f)} \right\|_{U}^{2} , \left\| \frac{-1}{P'(f)} \right\|_{U}^{d+2} \right) \le 1.$$
Nous avons alors 
$$\left\|  R\left( \frac{-1}{P'(f)}\right) \right\|_{U} \le 1.$$
On en d\'eduit que, quel que soit $n\in\N^*$, nous avons
$$\left\| R^{\circ n}\left( \frac{-1}{P'(f)}\right) \right\|_{U} \le 1,$$
o\`u $R^{\circ n}$ d\'esigne l'application $R$ \'elev\'ee \`a la puissance $n$ pour la loi de composition. 

En utilisant le fait que, si un \'el\'ement $b$ de $\Bs(U)$ v\'erifie $\|b\|_{U} \le 1$, alors 
$$\|R(b)\|\le t\, \|b\|_{U}^2,$$ 
on montre, \`a l'aide d'une r\'ecurrence, que, quel que soit $n\in\N^*$, nous avons 
$$\left\| R^{\circ n}\left( \frac{-1}{P'(f)}\right) \right\|_{U} \le t^{2^{n-1} -1}.$$
En particulier, la s\'erie 
$$\sum_{n\in\N} R^{\circ n}\left( \frac{-1}{P'(f)}\right)$$
converge dans $\Bs(U)$. Notons $s$ sa somme. Elle v\'erifie l'\'equation 
$$s - R(s) =  -\frac{1}{P'(f)}.$$
On en d\'eduit que $P(f+P(f)s)=0$. Puisque $P(f)$ est nul dans $\kappa(x)$, l'\'el\'ement $f+P(f)s$ de $\Os_{X,x}$ rel\`eve bien $\alpha$.
\end{proof}

\begin{cor}\label{Henselgeneral}
Soit $(Z,\Os_{Z})$ un espace analytique sur $\As$ (au sens de la d\'efinition \ref{espan}). Pour tout point $z$ de $Z$, l'anneau local $\Os_{Z,z}$ est hens\'elien.
\end{cor}
\begin{proof}
Par d\'efinition, l'anneau local $\Os_{Z,z}$ est le quotient de l'anneau local en un point d'un espace affine analytique sur $\As$. Ce dernier anneau est hens\'elien, d'apr\`es la proposition pr\'ec\'edente. Cela suffit pour conclure car tout quotient d'un anneau hens\'elien est hens\'elien. 
\end{proof}

\index{Henselianite@Hens\'elianit\'e|)}

\subsection{Isomorphismes locaux}

Le caract\`ere hens\'elien d'un anneau local peut
\^etre interpr\'et\'e comme une sorte de th\'eor\`eme des fonctions
implicites. Par la suite, nous utiliserons effectivement cette
propri\'et\'e pour d\'emontrer des r\'esultats d'isomorphie. La
proposition qui suit donne un exemple d'application.

Soit $P(S)$ un polyn\^ome unitaire \`a coefficients dans $\As$. Notons $d\in\N$ son degr\'e. Nous nous int\'eresserons \`a l'alg\`ebre
$$\As' = \As[S]/(P(S)).$$
Puisque le polyn\^ome est unitaire, le morphisme 
$$n : \begin{array}{ccc}
\As^d & \to & \As'\\
(a_{0},\ldots,a_{d-1}) & \mapsto & \disp \sum_{i=0}^{d-1} a_{i}\, S^i
\end{array}$$
est un isomorphisme. Munissons l'alg\`ebre $\As^d$ de la norme~$\|.\|_{\infty}$ donn\'ee par le maximum des normes des coefficients. 
On d\'efinit alors une norme, not\'ee~$\|.\|_{\textrm{div}}$, sur~$\As'$ de la fa\c{c}on suivante :
$$\forall f\in \As',\, \|f\|_{\textrm{div}} = \|n^{-1}(f)\|_{\infty}.$$
Cette norme n'est pas, \emph{a priori}, une norme d'alg\`ebre. Nous supposerons donc que l'alg\`ebre~$\As'$ est munie d'une norme d'alg\`ebre~$\|.\|'$ \'equivalente \`a la norme~$\|.\|_{\textrm{div}}$ : il existe deux constantes $D_{-},D_{+}>0$ telles que 
$$\forall f\in\As',\, D_{-}\, \|f\|_{\textrm{div}} \le \|f\|' \le D_{+}\, \|f\|_{\textrm{div}}.$$
Munie de la norme~$\|.\|'$, l'alg\`ebre $\As'$ est une alg\`ebre de Banach. En outre, le morphisme $(\As,\|.\|)\to (\As',\|.\|')$ est born\'e. Nous noterons
$$\varphi : X'=\E{n}{\As'} \to \E{n}{\As}=X$$
le morphisme induit entre les espaces analytiques.

Soit $U$ une partie ouverte de $X$ et supposons qu'il existe une fonction $R$ d\'efinie sur $U$ v\'erifiant $P(R)=0$. Signalons qu'en pratique, nous d\'eduirons l'existence d'une telle fonction du caract\`ere hens\'elien d'un certain anneau local.


Nous pouvons alors d\'efinir une application~$\sigma$ de $U\subset X$ vers $X'$. Soit $x$ un point de $U$. Soit $p(\bT) = \sum_{\bk \ge 0} p_{\bk}\, \bT^\bk$, o\`u la famille $(p_{\bk})_{\bk\ge 0}$ est une famille presque nulle d'\'el\'ements de $\As'$. Quel que soit $\bk\in\N^n$, relevons l'\'el\'ement $p_{\bk}$ de $\As'$ en un \'el\'ement $q_{\bk}(S)$ de $\As[S]$. Consid\'erons l'application
$$\chi_{\sigma(x)} : {\renewcommand{\arraystretch}{1.5}\begin{array}{ccc}
\As'[\bT] & \to & \Hs(x)\\
p(\bT) & \mapsto & \disp \sum_{\bk\ge 0} q_{\bk}(R(x))\,\bT^\bk(x)
\end{array}}.$$ 
Puisque $P(R(x))=0$, cette application ne d\'epend pas du choix des diff\'erents relev\'es. On en d\'eduit aussit\^ot que $\chi_{\sigma(x)}$ est un morphisme de $\As$-alg\`ebres. Montrons que ce morphisme est born\'e sur $\As'$. Soit $f\in\As'$. Il existe $a_{0},\ldots,a_{d-1}\in\As$ tels que 
$$f = \sum_{i=0}^{d-1} a_{i}\, S^i \textrm{ dans } \As'.$$ 
Nous avons alors 
$${\renewcommand{\arraystretch}{2}\begin{array}{rcl}
\left|\chi_{\sigma(x)}(f)\right| &=& \disp \left|  \sum_{i=0}^{d-1} a_{i}(x)\, R(x)^i \right|\\
& \le &  \disp \left( \sum_{i=0}^{d-1} |R(x)^i| \right)\, \max_{0\le i\le d-1}(|a_{i}(x)|)\\
& \le &  \disp   \left( \sum_{i=0}^{d-1} |R(x)^i| \right)\, \max_{0\le i\le d-1}(\|a_{i}\|)\\
& \le &  \disp  \left( \sum_{i=0}^{d-1} |R(x)^i| \right)\, D_{-}^{-1}\,\|f\|'.\\ 
\end{array}}$$
Par cons\'equent, le morphisme $\chi_{\sigma(x)}$ est born\'e sur $\As'$. C'est donc un caract\`ere de~$\As'[\bT]$. Nous noterons $\sigma(x)$ le point de $X'$ associ\'e. L'application $\sigma$ ainsi construite est une section continue de $\varphi$ au-dessus de $U$. Sous certaines hypoth\`eses, nous pouvons obtenir un r\'esultat bien plus fort. Nous noterons $\alpha$ l'image de $S$ dans $\As'$.

\begin{prop}\label{isolocal}\index{Isomorphisme local}
Supposons que
\begin{enumerate}[\it i)]
\item la norme $\|.\|'$ sur $\As'$ est uniforme et \'equivalente \`a la norme~$\|.\|_{\textrm{div}}$ ;
\item l'ouvert $U$ est connexe ;
\item la fonction~$P'(\alpha)$ est inversible sur $\varphi^{-1}(U)$ ;
\item il existe un point \mbox{$x_{0} \in U$} tel que $R(\sigma(x_{0})) = \alpha$ dans $\Hs(\sigma(x_{0}))$.
\end{enumerate} 
Alors la partie~$\sigma(U)$ est ouverte dans~$X'$ et la section $\sigma$ induit un isomorphisme entre les espaces~$U$ et~$\sigma(U)$, munis des structures d'espaces localement annel\'es induites. 
\end{prop}
\begin{proof}
Le polyn\^ome $P(T)$ poss\`ede une unique factorisation dans $\As'[T]$ sous la forme $P(T) = (T-\alpha)Q(T)$, avec $Q(T)\in \As'[T]$. Quel que soit le point $x'$ de $\varphi^{-1}(U)$, nous avons $P(R(x'))=0$, d'o\`u l'on tire soit $R(x')=\alpha$, soit $Q(R(x'))=0$. Ces deux conditions ne peuvent valoir simultan\'ement, puisque, par hypoth\`ese, nous avons $P'(\alpha)(x')\ne 0$. Par cons\'equent, la partie de~$X'$ d\'efinie par 
$$U'=\{x'\in \varphi^{-1}(U)\,|\, R(x')=\alpha\}$$
est ouverte.

Montrons, \`a pr\'esent, que $\sigma(U)=U'$. Par hypoth\`ese, nous avons $R(\sigma(x_{0}))=\alpha$, autrement dit, le point~$\sigma(x_{0})$ appartient \`a~$U'$. Puisque l'ouvert~$U$ est connexe, la partie~$\sigma(U)$ l'est encore. Nous en d\'eduisons l'inclusion~$\sigma(U)\subset U'$.

R\'eciproquement, soit $x'$ un point de $U'$. Par d\'efinition de
$U'$, nous avons~$R(x')=\alpha$. Notons $x\in U$ son image par le
morphisme~$\varphi$. Soit $p(\bT) = \sum_{\bk \ge 0} p_{\bk}\,
\bT^\bk$, o\`u la famille~$(p_{\bk})_{\bk\ge 0}$ est une famille
presque nulle d'\'el\'ements de~$\As[S]/(P(S))$. Quel que soit
$\bk\in\N^n$, relevons l'\'el\'ement~$p_{\bk}$ en un
\'el\'ement~$q_{\bk}$ de~$\As[S]$. Le caract\`ere~$\chi_{\sigma(x)}$
envoie le polyn\^ome~$p_{\bk}$ sur l'\'el\'ement
$$\sum_{\bk\ge 0} q_{\bk}(R(x))\,\bT^\bk(x) \textrm{ de } \Hs(x).$$
L'image de cet \'el\'ement par l'injection $\Hs(x)\hookrightarrow
\Hs(x')$ n'est autre que
$$\sum_{\bk\ge 0} q_{\bk}(R(x'))\,\bT^\bk(x') = \sum_{\bk\ge 0} q_{\bk}(\alpha)\,\bT^\bk(x) = p(\bT(x')) \textrm{ dans } \Hs(x').$$
On en d\'eduit que $\sigma(x)=x'$.

Nous venons de d\'emontrer que le morphisme~$\varphi$ r\'ealise un
hom\'eomorphisme de l'ouvert~$U'$ de~$X'$ sur l'ouvert~$U$ de~$X$.
Nous allons prouver qu'il induit m\^eme un isomorphisme entre les
espaces annel\'es. Soit~$x'$ un point de~$U'$. Notons $x\in U$ son
image par le morphisme~$\varphi$. Il suffit de montrer que le
morphisme
$$\Os_{X,x} \to \Os_{X',x'}$$
induit par~$\varphi$ est un isomorphisme. Montrons, tout d'abord,
qu'il est injectif. Soit~$f$ une fonction analytique d\'efinie sur un
voisinage~$V$ de~$x$ dans~$U$ dont l'image dans l'anneau
local~$\Os_{X',x'}$ est nulle. Il existe alors un voisinage~$W'$ du
point~$x'$ dans~$\varphi^{-1}(V)$ tel que, quel que soit~$y'$
dans~$W'$, nous ayons
$$\varphi^*(f)(y')=0 \textrm{ dans } \Hs(y').$$
On en d\'eduit que, quel que soit~$y$ dans~$\varphi(W')$, nous avons
$$f(y)=0 \textrm{ dans } \Hs(y).$$
La partie~$W=\varphi(W')$ est un voisinage de~$x$ dans~$X$,
car~$\varphi$ est un hom\'eomorphisme sur~$U'$, et la fonction~$f$ est
nulle en tout point de ce voisinage. Cette condition impose que \`a la
fonction~$f$ d'\^etre nulle en tant qu'\'el\'ement de~$\Os_{X}(W)$, et donc
dans l'anneau local~$\Os_{X,x}$, car l'alg\`ebre~$\As$ est munie d'une
norme uniforme.

Montrons, \`a pr\'esent, que le morphisme entre les anneaux locaux est
surjectif. Soit $f'\in\Os_{X',x'}$. Il existe un voisinage
ouvert~$V'$ de~$x'$ dans~$U'$ et une suite~$(R'_{n})_{n\in\N}$
d'\'el\'ements de~$\Ks(V')$ qui converge vers~$f'$ uniform\'ement sur~$V'$.

Soit~$n\in\N$. Il existe deux \'el\'ements~$p'$ et~$q'$
de~$\As'[\bT]$, $q'$ ne s'annulant pas sur~$V'$, tels que~$R'_n=p'/q'$
dans~$\Ks(V')$. Il existe une famille presque
nulle~$(p'_{\bk})_{\bk\ge 0}$ d'\'el\'ements de~$\As'$ telle que
$$p'(\bT)  = \sum_{\bk \ge 0} p'_{\bk}\, \bT^\bk.$$
Quel que soit~$\bk$ dans~$\N^n$, relevons l'\'el\'ement~$p'_{\bk}$
en un \'el\'ement~$p_{\bk}$ de~$\As[S]$. Puisque $V=\varphi(V')$
est contenu dans~$U$, la fonction~$R$ y est d\'efinie. Il en est
donc de m\^eme pour la fonction
$$p(\bT) = \sum_{\bk \ge 0} p_{\bk}(R)\, \bT^\bk \textrm{ de } \Os_{X}(V).$$
Par d\'efinition de $U'$, au-dessus de~$U'$, nous avons~$R=\alpha$.
On en d\'eduit que
$$\varphi^*(p)=\sum_{\bk \ge 0} p_{\bk}(\alpha)\, \bT^\bk = p' \textrm{ dans } \Os_{X}(V').$$

En proc\'edant
comme pr\'ec\'edemment, on montre qu'il existe un \'el\'ement~$q$
de~$\Os_{X}(V)$ tel que
$$\varphi^*(q)=q' \textrm{ dans } \Os_{X}(V').$$ 
Puisque la fonction~$q'$ ne s'annule pas sur~$V'$, la fonction~$q$ ne s'annule pas sur~$V$
et elle est donc inversible dans~$\Os_{X}(V)$. L'\'el\'ement~$R_{n} = pq^{-1}$ de~$\Os_{X}(V)$ v\'erifie l'\'egalit\'e
$$\varphi^*(R_{n})=R_{n}' \textrm{ dans } \Os_{X}(V').$$ 

Puisque la suite~$(R'_{n})_{n\in\N}$ converge uniform\'ement sur~$V'$, la suite~$(R_{n})_{n\in\N}$ est une suite de Cauchy uniforme sur toute partie compacte de~$V$. Elle converge donc vers une fonction~$f$ de~$\Os_{X}(V)$. Cette fonction v\'erifie
$$\varphi^*(f)=f'$$
dans~$\Os_{X}(V')$ et donc dans~$\Os_{X',x'}$.
C'est ce que nous voulions d\'emontrer.

\end{proof}

\begin{rem}

En g\'en\'eral, il n'est pas ais\'e de montrer que l'hypoth\`ese~{\it i)} de la proposition pr\'ec\'edente est satisfaite. Nous \'etablirons, dans un chapitre ult\'erieur, des crit\`eres permettant de l'assurer (\emph{cf.} lemme~\ref{equivdiv} et proposition~\ref{finiuniforme}). Nous prouverons \'egalement qu'elle est v\'erifi\'ee dans trois cas particuliers, dans la preuve des propositions~\ref{isoext}, \ref{chgtinterne} et~\ref{isocentral}.
\end{rem}

%% file: espace.tex
\chapter[Espace affine sur un corps de nombres]{Espace affine analytique au-dessus d'un anneau d'entiers de corps de nombres}\label{chapitreespace}

Ce chapitre est consacr\'e \`a l'\'etude des espaces analytiques au-dessus d'un anneau d'entiers de corps de nombres. 
Dans ce cadre, nous allons pouvoir pr\'eciser et g\'en\'eraliser les r\'esultat obtenus au chapitre pr\'ec\'edent.

Dans le num\'ero \ref{espacedebase}, nous nous int\'eressons au spectre analytique de l'anneau d'entiers de corps de nombres~$A$. Nous commen\c{c}ons par le d\'e\-cri\-re ensemblistement et poursuivons en \'etablissant ses propri\'et\'es topologiques. Pour finir, nous d\'ecrivons les sections du faisceau structural au-dessus des ouverts de cet espace et en d\'eduisons notamment l'expression des anneaux locaux.

Dans la suite du chapitre, nous passons \`a l'\'etude des espaces affines de dimension quelconque. Au num\'ero \ref{fssea}, nous commen\c{c}ons par reprendre, pour les pr\'eciser dans le cadre que nous avons choisi, les descriptions des anneaux locaux que nous avons d\'ej\`a obtenues. \`A titre d'application, nous utilisons le caract\`ere hens\'elien d'un certain anneau local pour donner une nouvelle d\'emonstration du th\'eor\`eme classique d'Eisenstein (\emph{cf.} th\'eor\`eme \ref{Eisenstein}). \`A la fin de ce num\'ero, nous nous int\'eressons aux anneaux de sections globales sur les disques et couronnes relatifs et en proposons une description explicite. 

Les num\'eros \ref{prdf} et \ref{fibresinternes} sont consacr\'es \`a l'\'etude de certains types de points : points rigides des fibres, puis points internes. Nous d\'ecrivons des syst\`emes fondamentaux de voisinages et d\'emontrons quelques propri\'et\'es alg\'ebriques des anneaux locaux.

Nous parvenons \`a d\'ecrire et \'etudier les anneaux locaux en plusieurs types de points, au nombre desquels les points rigides des fibres. Les r\'esultats que nous obtenons ne sont cependant pas complets : certains probl\`emes ont, jusqu'ici, r\'esist\'e \`a nos tentatives et requi\`erent vraisemblablement une approche nouvelle. Nous proposons \'egalement une description explicite des anneaux de sections globales sur les disques et couronnes relatifs.

Au num\'ero \ref{dimensiontopologique}, nous nous int\'eressons \`a la dimension topologique des espaces affines au-dessus d'un anneau d'entiers de corps de nombres.

Finalement, le num\'ero \ref{sectionprolan} est consacr\'e au prolongement analytique. Il ne contient presqu'aucun r\'esultat et nous nous contentons d'y \'enoncer quelques d\'efinitions et propri\'et\'es li\'ees \`a cette question, en vue d'une utilisation ult\'erieure.


\bigskip

\index{Corps de nombres!notations|(}

Dans ce chapitre, nous fixons un corps de nombres~$K$. Nous noterons~$A$ l'anneau de ses entiers. 
\newcounter{K}\setcounter{K}{\thepage}


%

\section[Spectre analytique]{Spectre analytique d'un anneau d'entiers de corps de nombres}\label{espacedebase}

\index{Spectre analytique!de A@de $A$|(}

Dans cette partie, nous allons \'etudier le spectre analytique de l'anneau d'entiers de corps de nombres~$A$. Pour ce faire, nous devons le munir d'une norme qui en fasse un anneau de Banach. Plusieurs choix d'offrent \`a nous : norme triviale, restriction de la valeur absolue complexe, etc. Nous choisirons la norme~$\|.\|$ d\'efinie de la fa\c{c}on suivante :
$$\forall f\in A,\, \|f\| = \max_{\sigma : K\hookrightarrow \C}(|\sigma(f)|_{\infty}),$$\index{Norme!sur un anneau d'entiers de corps de nombres}
o\`u le maximum est pris sur l'ensemble des plongements $\sigma$ du corps $K$ dans $\C$. Par exemple, lorsque $K=\Q$, cette norme est simplement la valeur absolue usuelle~$|.|_{\infty}$.\newcounter{vainfinie}\setcounter{vainfinie}{\thepage} Notre choix est guid\'e par le fait que cette norme est plus grande que toutes les semi-normes multiplicatives que l'on peut d\'efinir sur l'anneau $A$. Le spectre $\Ms(A,\|.\|)$ contiendra donc tous les points possibles. 

Remarquons que l'anneau $A$ muni de la norme $\|.\|$ est bien un anneau de Banach. En effet, quel que soit $f\in A\setminus\{0\}$, nous avons $\|f\|\ge 1$. Cette in\'egalit\'e d\'ecoule simplement de la formule du produit. Par cons\'equent, la topologie induite sur~$A$ par la norme~$\|.\|$ est discr\`ete. 

Dans la suite de ce texte, nous supposerons toujours que l'anneau~$A$ est muni de la norme~$\|.\|$. Nous \'ecrirons donc~$\Ms(A)$ et~$\E{n}{A}$, pour tout entier~$n$,\newcounter{MA}\setcounter{MA}{\thepage} sans plus de pr\'ecisions. Nous noterons simplement~$\Os$ le faisceau structural sur l'espace $\Ms(A)$. 

\subsection{Description ensembliste et topologique}\label{descriptionMA}

Le th\'eor\`eme d'Ostrowski nous permet de d\'ecrire explicitement toutes les semi-normes multiplicatives sur~$A$, autrement dit l'ensemble~$\Ms(A)$.

\bigskip

Nous avons, tout d'abord, la valeur absolue triviale 
\newcounter{vatriviale}\setcounter{vatriviale}{\thepage}
$$|.|_0 : \begin{array}{ccc}
K & \to & \R_+\\
f & \mapsto & \left\{
\begin{array}{cl}
0 & \textrm{si } f = 0\\
1 & \textrm{sinon}
\end{array}\right.
\end{array}.$$
Nous noterons $a_{0}$\newcounter{azero}\setcounter{azero}{\thepage} le point de $\Ms(A)$ correspondant. Le corps r\'esiduel en ce point est 
$$(\Hs(a_{0}),|.|)= (K,|.|_{0}).$$

\bigskip

Soit $p$ un nombre premier. Nous noterons $v_{p}$ la valuation $p$-adique sur~$\Q$ et~$|.|_{p}$\newcounter{vapadique}\setcounter{vapadique}{\thepage} la valeur absolue $p$-adique d\'efinie par
$$|.|_p : \begin{array}{ccc}
\Q & \to & \R_+\\
f & \mapsto & p^{-v_{p}(f)}
\end{array}.$$

Soit~$\m$\newcounter{m}\setcounter{m}{\thepage} un id\'eal maximal de~$A$. L'anneau local~$A_{\m}$ est un anneau de valuation discr\`ete. Notons $k_{\m}=A/\m$ son corps r\'esiduel. Choisissons une uniformisante~$\pi_{\m}$ de~$A_{\m}$. Nous noterons encore~$\hat{A}_{\m}$ le compl\'et\'e de~$A_{\m}$ pour la topologie $\m$-adique et~$\hat{K}_{\m}$ son corps des fractions. Notons~$p_{\m}$ le nombre premier tel que $\m \cap \Z = p_{\m}\,\Z$. Le corps~$\hat{K}_{\m}$ est alors une extension finie du corps~$\Q_{p_{\m}}$, dont nous noterons~$n_{\m}$ le degr\'e. Nous noterons~$|.|_{\m}$ l'unique valeur absolue sur~$K$ qui prolonge la valeur absolue~$|.|_{p_{\m}}$ sur~$\Q$. Pour tout \'el\'ement~$f$ de~$K$, nous avons
$$|f|_{\m} = \left| N_{\hat{K}_{\m}/\Q_{p_{\m}}}(f)\right|_{p_{\m}}^{1/n_{\m}}.$$
Nous noterons~$a_{\m}$ le point de~$\Ms(A)$ correspondant \`a la valeur absolue~$|.|_{\m}$. 


\`A chaque nombre r\'eel strictement positif~$\eps$, on associe alors la valeur absolue~$|.|_{\m}^\eps$ sur~$K$. Nous noterons~$a_{\m}^{\eps}$ le point de~$\Ms(A)$ correspondant. Le corps r\'esiduel en ce point est 
$$(\Hs(a_{\m}^{\eps}),|.|)= (\hat{K}_{\m},|.|_{\m,\eps}).$$
\newcounter{ame}\setcounter{ame}{\thepage}

Lorsque nous faisons tendre $\eps$ vers 0 dans la formule pr\'ec\'edente, nous retrouvons la valeur absolue triviale. Nous noterons donc 
$$a_{\m}^0 = a_{0}.$$ 

Lorsque nous faisons tendre $\eps$ vers $+\infty$, nous obtenons la semi-norme multiplicative induite par la valeur absolue triviale sur le corps fini $k_{\m}$ :
$$|.|_{\m,\infty} : \begin{array}{ccc}
A & \to & \R_+\\
f & \mapsto & \left\{
\begin{array}{cl}
0 & \textrm{si } f\in \m\\
1 & \textrm{sinon}
\end{array}\right.
\end{array}.$$
Nous noterons $\tilde{a}_{\m}$, ou encore $a_{\m}^{+\infty}$, le point de $\Ms(A)$ correspondant. Le corps r\'esiduel en ce point est 
$$(\Hs(\tilde{a}_{\m}),|.|)= (\tilde{k}_{\m},|.|_{0}).$$
\newcounter{tam}\setcounter{tam}{\thepage}

\bigskip

\newcounter{sigma}\setcounter{sigma}{\thepage}
Soit $\sigma$ un plongement du corps~$K$ dans~$\C$. Nous poserons $\hat{K}_{\sigma}=\R$ si le plongement est r\'eel, c'est-\`a-dire si son image est contenue dans~$\R$, et $\hat{K}_{\sigma}=\C$ dans les autres cas. Nous noterons~$|.|_{\sigma}$ la valeur absolue sur~$K$ d\'efinie par
$$|.|_\sigma : \begin{array}{ccc}
K & \to & \R_+\\
f & \mapsto & |\sigma(f)|_{\infty}
\end{array},$$
o\`u $|.|_\infty$ d\'esigne la valeur absolue usuelle sur~$\C$. Nous noterons~$a_{\sigma}$ le point de~$\Ms(A)$ correspondant. Remarquons que deux plongements complexes conjugu\'es d\'efinissent la m\^eme valeur absolue et donc le m\^eme point de $\Ms(A)$. Nous noterons~$n_{\sigma}$ le degr\'e de l'extension~$\hat{K}_{\sigma}/\R$.

\`A chaque nombre r\'eel $\eps\in\mathopen{[}0,1\mathclose{]}$, on associe la valeur absolue~$|.|_{\sigma}^{\eps}$ 
sur~$K$. Nous noterons~$a_{\sigma}^{\eps}$ le point de~$\Ms(A)$ correspondant. Le corps r\'esiduel en ce point est 
$$(\Hs(a_{\sigma}^{\eps}),|.|)= (\hat{K}_{\sigma},|.|_{\sigma,\eps}).$$

\begin{rem}
Pour $\eps>1$, l'application $|.|_{\sigma}^\eps$ ne d\'efinit plus une norme, car elle ne satisfait plus l'in\'egalit\'e triangulaire.
\end{rem}

Comme pr\'ec\'edemment, lorsque nous faisons tendre $\eps$ vers 0, nous retrouvons la valeur absolue triviale. Nous noterons donc 
$$a_{\sigma}^0 = a_0.$$ 

\bigskip

Adoptons quelques notations suppl\'ementaires. Nous noterons $\Sigma_{f}=\textrm{Max}(A)$ l'ensemble des id\'eaux maximaux de~$A$ et $\Sigma_{\infty}$ l'ensemble des plongements du corps~$K$ dans le corps~$\C$, \`a conjugaison pr\`es. D\'esignons par~$r_{1}$\newcounter{run}\setcounter{run}{\thepage} le nombre de plongements r\'eels de~$K$ et par~$2 r_{2}$ 
son nombre de plongements complexes non r\'eels. Nous avons alors
$$\sharp\left(\Sigma_{\infty}\right) = r_{1} + r_{2}.$$
Rappelons que l'on a l'\'egalit\'e $r_{1} + 2r_{2} = [K : \Q]$.

Pour finir, nous notons $\Sigma = \Sigma_{f} \cup \Sigma_{\infty}$ et posons
$$l(\sigma) = \left\{ 
\begin{array}{cl}
+\infty & \textrm{si } \sigma\in\Sigma_{f}\ ;\\
1 &  \textrm{si } \sigma\in\Sigma_{\infty}.
\end{array}
\right.$$

\index{Corps de nombres!notations|)}

\begin{prop}[formule du produit]\label{formuleduproduit}\index{Corps de nombres!formule du produit}\index{Formule du produit|see{Corps de nombres}}
Pour tout \'el\'ement non nul~$f$ de~$K$, nous avons l'\'egalit\'e
$$\prod_{\sigma\in\Sigma} |f|_{\sigma}^{n_{\sigma}} = 1.$$
\end{prop}



\begin{thm}[Ostrowski]
L'ensemble $\Ms(A)$ est constitu\'e exactement des points d\'ecrits pr\'ec\'edemment.
\end{thm}
\begin{proof}
Soit~$b$ un point de l'espace~$\Ms(A)$. Notons
$$\p_{b} = \left\{f\in A\, \big|\, |f(b)|=0\right\}.$$
C'est un id\'eal premier de l'anneau~$A$. Puisque l'anneau~$A$ est un anneau de Dedekind, l'id\'eal~$\p_{b}$ est soit l'id\'eal nul, soit un id\'eal maximal.

Supposons, tout d'abord, que~$\p_{b}$ est un id\'eal maximal~$\m$ de~$A$. Dans ce cas, la semi-norme multiplicative~$|.|_{b}$ associ\'ee au point~$b$ induit une valeur absolue sur le quotient~$A/\m$. Or, ce quotient est un corps fini. Il ne peut donc \^etre muni que de la valeur absolue triviale. On en d\'eduit que le point~$b$ n'est autre que le point~$\tilde{a}_{\m}$.

Supposons, maintenant, que~$\p_{b}$ est l'd\'eal nul. Dans ce cas, la semi-norme multiplicative~$|.|_{b}$ associ\'ee au point~$b$ est une valeur absolue sur l'anneau~$A$. La version habituelle du th\'eor\`eme d'Ostrowski entra\^{\i}ne alors le r\'esultat.
\end{proof}

\bigskip

La description explicite des points nous permet de d\'ecrire, de fa\c{c}on tout aussi explicite, la topologie de l'espace $\Ms(A)$.
\begin{lem}\label{topobranche}
Soit $\sigma\in\Sigma$. L'application
$$a_{\sigma}^. : \begin{array}{ccc}
\of{[}{0,l(\sigma)}{]} & \to & \Ms(A)\\
\eps & \mapsto & a_{\sigma}^{\eps}
\end{array}$$
induit un hom\'eomorphisme sur son image. 
\end{lem}
\begin{proof}
Par d\'efinition de la topologie de~$\Ms(A)$, pour montrer que l'application~$a_{\sigma}^.$ est continue, il suffit de montrer que, quel que soit~$f\in A$, l'application compos\'ee
$$\begin{array}{ccccc}
\of{[}{0,l(\sigma)}{]} & \to & \Ms(A) & \to &\R_{+}\\
\eps & \mapsto & a_{\sigma}^{\eps} & \mapsto & |f(a_{\sigma}^\eps)| = |f|_{\sigma}^\eps
\end{array}$$
est continue. Ce r\'esultat est imm\'ediat. Puisque l'espace~$\of{[}{0,l(\sigma)}{]}$ est compact et que l'espace~$\Ms(A)$ est s\'epar\'e, l'application~$a_{\sigma}^.$ induit un hom\'eomorphisme sur son image. 
\end{proof}

\begin{defi}\label{branches}\index{Branche sigma-adique@Branche $\sigma$-adique}\index{Branche sigma-adique@Branche $\sigma$-adique!ouverte}\index{Branche sigma-adique@Branche $\sigma$-adique!semi-ouverte}
\newcounter{branche}\setcounter{branche}{\thepage}
Soit~$\sigma\in\Sigma$. Nous appellerons {\bf branche~$\sigma$-adique} l'image de l'application pr\'ec\'edente et la noterons $\Ms(A)_{\sigma}$. Nous appellerons {\bf branche $\sigma$-adique ouverte}, et noterons~$\Ms(A)_{\sigma}'$, la branche~$\sigma$-adique priv\'ee des points associ\'es \`a une valeur absolue triviale. Nous \^oterons donc deux points si $\sigma\in\Sigma_{f}$, mais un seul point si $\sigma\in\Sigma_{\infty}$. Signalons que ces branches ouvertes sont les trajectoires du flot, au sens du num\'ero \ref{parflot}. Pr\'ecis\'ement, quel que soit~\mbox{$\eps\in\of{]}{0,l(\sigma)}{]}$} tel que~$a_{\sigma}^\eps \in \Ms(A)_{\sigma}'$, nous avons
$$T_{\Ms(A)}(a_{\sigma}^\eps) \simeq \Ms(A)_{\sigma}'.$$
Nous appellerons {\bf branche $\sigma$-adique semi-ouverte}, et noterons~$\Ms(A)_{\sigma}''$, la branche~$\sigma$-adique priv\'ee du point associ\'e \`a la valeur absolue triviale sur~$A$. Cette d\'efinition co\"{\i}ncide avec la pr\'ec\'edente dans le cas des \'el\'ements de~$\Sigma_{\infty}$.

Nous appellerons {\bf point central}\index{Point!central} de $\Ms(A)$ le point $a_{0}$. Nous appellerons {\bf point extr\^eme}\index{Point!extreme@extr\^eme} de $\Ms(A)$ tout point de la forme $\tilde{a}_{m}$, o\`u~$\m$ est un \'el\'ement de~$\Sigma_{f}$. Enfin, nous appellerons {\bf point interne}\index{Point!interne} de $\Ms(A)$ tout autre point. En particulier, quel que soit $\sigma\in\Sigma_{\infty}$, le point $a_{\sigma}=a_{\sigma}^1$ est un point interne.
\end{defi}



\bigskip

Afin de d\'ecrire plus pr\'ecis\'ement la topologie de l'espace~$\Ms(A)$, nous aurons besoin de quelques r\'esultats de th\'eorie des nombres.

\begin{lem}\label{diviseur}
Soit $\m\in\Sigma_{f}$. Alors il existe un \'el\'ement~$f$ de~$A$ qui v\'erifie les propri\'et\'es suivantes :
\begin{enumerate}[\it i)]
\item $|f|_{\m}<1$ ;
\item $\forall \m'\in\Sigma_{f}\setminus\{\m\}$, $|f|_{\m'}=1$.
\end{enumerate}
\end{lem}
\begin{proof}
Notons~$P$ le point de~$\Spec(A)$ associ\'e \`a l'id\'eal maximal~$\m$. Puisque le groupe de Picard de $\Spec(A)$ est fini, il existe $N\in\N^*$ tel que le diviseur $N[P]$ soit principal. Tout \'el\'ement~$f$ de~$A$ dont~$N[P]$ est le diviseur convient. 
\end{proof}

\begin{lem}\label{pasQquad}
Supposons que le corps $K$ ne soit ni $\Q$, ni un corps quadratique imaginaire. Alors, quel que soit~$\sigma\in\Sigma$, il existe un \'el\'ement~$f$ de~$A$ qui v\'erifie les conditions suivantes :
\begin{enumerate}[\it i)]
\item $|f|_{\sigma} < 1$ ;
\item $\forall \sigma'\in\Sigma_{f}\setminus\{\sigma\}$, $|f|_{\sigma'}=1$ ;
\item $\forall \sigma'\in \Sigma_{\infty}\setminus\{\sigma\}$, $|f|_{\sigma'} > 1$.
\end{enumerate}
\end{lem}
\begin{proof}
Notons~$\sigma_{1},\ldots,\sigma_{r_{1}}$, avec~$r_{1}\in\N$, les plongements r\'eels du corps~$K$ et \mbox{$\sigma_{r_{1}+1},\ldots,\sigma_{r_{1}+r_{2}}$}, avec~$r_{2}\in\N$, ses plongements complexes non r\'eels \`a conjugaison pr\`es. Par hypoth\`ese, nous avons~$r_{1}+r_{2}\ge 2$. Rappelons que, d'apr\`es le th\'eor\`eme des unit\'es de Dirichlet, le morphisme de groupes $L$ qui \`a toute unit\'e $f\in A^\times$ associe l'\'el\'ement 
$$\big(\log(|\sigma_{1}(g)|),\ldots,\log(|\sigma_{r_{1}}(g)|),2\log(|\sigma_{r_{1}+1}(g)|),\ldots, 2\log(|\sigma_{r_{1}+r_{2}}(g)|)\big)$$
de $\R^{r_{1}+r_{2}}$ a pour image un r\'eseau de l'hyperplan~$H$ de~$\R^{r_{1}+r_{2}}$ d\'efini par l'\'equation
$$H\ :\ x_{1}+\cdots+x_{r_{1}+r_{2}} = 0.$$

Supposons, tout d'abord, que $\sigma\in\Sigma_{\infty}$. Il existe alors $i\in\cn{1}{r_{1}+r_{2}}$ tel que~$\sigma=\sigma_{i}$. Consid\'erons le quadrant de~$\R^{r_{1}+r_{2}}$ d\'efini par 
$$Q = \{(x_{1},\ldots,x_{r_{1}+r_{2}}) \in \R^{r_{1}+r_{2}}\, |\, x_{i}<0,\, \forall j\ne i,\, x_{j}>0\}.$$
Le r\'esultat rappel\'e ci-dessus assure qu'il existe une unit\'e $f\in A^\times$ telle que 
$$L(f) \in Q.$$
Nous avons alors~$|f|_{\sigma_{i}}<1$, quel que soit~$j\ne i$, $|f|_{\sigma_{j}}>1$ et, quel que soit~$\m\in\Sigma_{f}$, $|f|_{\m}=1$.


D'apr\`es le lemme \ref{diviseur}, il existe un \'el\'ement~$f$ de~$A$ qui v\'erifie $|f|_{\m}<1$ et, pour tout \'el\'ement~$\m$ de $\Sigma_{f}\setminus\{\m\}$, $|f|_{\m'}=1$. La formule du produit assure alors que
$$\prod_{i=1}^{r_{1}} |f|_{\sigma_{i}}\, \prod_{i=r_{1}+1}^{r_{1}+r_{2}} |f|_{\sigma_{i}}^2 >1.$$
Notons $L(f) = (y_{1},\ldots,y_{r_{1}+r_{2}})\in R^{r_{1}+r_{2}}$. Nous avons alors
$$S = \sum_{i=1}^{r_{1}+r_{2}} y_{i} >0.$$
Soit $\eps>0$ tel que $S > (r_{1}+r_{2}-1)\eps$. Posons
$$z_{0}= \big(-y_{1}+\eps,\ldots,-y_{r_{1}+r_{2}-1}+\eps, -y_{r_{1}+r_{2}} +S -(r_{1}+r_{2}-1)\eps\big) \in H.$$
Nous avons $L(f)+z_{0} \in (\R_{+}^*)^{r_{1}+r_{2}}$. Par cons\'equent, il existe un voisinage ouvert~$U$ de~$z_{0}$ dans~$H$ de volume~$v$ strictement positif tel que
$$L(f) + U \subset (\R_{+}^*)^{r_{1}+r_{2}}.$$
Soit $n\in\N^*$ tel $nv$ soit strictement plus grand que le volume d'une maille du r\'eseau $L(A^\times)$. La partie
$$nL(f) + nU \subset (\R_{+}^*)^{r_{1}+r_{2}}$$
contient alors un \'el\'ement $z$ du r\'eseau $L(A^\times)$. Il existe $g\in A^\times$ tel que $L(g)=z$. Posons $h = f^n\, g$. Nous avons toujours~$|h|_{\m}<1$ et, quel que soit~$\m'\in\Sigma_{f}\setminus\{\m\}$, $|h|_{\m'}=1$. En outre, nous avons
$$L(h) \in (\R_{+}^*)^{r_{1}+r_{2}},$$
autrement dit, quel que soit $i\in\cn{1}{r_{1}+r_{2}}$, $|h|_{\sigma_{i}}>1$.
\end{proof}

\begin{lem}\label{Qquad}
Supposons que le corps $K$ soit $\Q$ ou un corps quadratique imaginaire. Dans ce cas, $\Sigma_{\infty}$ est r\'eduit \`a un \'el\'ement que nous noterons $\sigma_{\infty}$. Alors, pour tout \'el\'ement~$\m$ de~$\Sigma_{f}$, il existe un \'el\'ement~$f$ de~$A$ qui v\'erifie les conditions suivantes :
\begin{enumerate}[\it i)]
\item $|f|_{\m} < 1$ ;
\item $|f|_{\sigma_{\infty}} > 1$ ;
\item $\forall \m'\in \Sigma_{f}\setminus\{\m\}$, $|f|_{\m'} = 1$.
\end{enumerate}
\end{lem}
\begin{proof}
D'apr\`es le lemme \ref{diviseur}, il existe un \'el\'ement~$f$ de~$A$ v\'erifiant $|f|_{\m} < 1$ et, pour tout \'el\'ement~$\m'$ de $\Sigma_{f}\setminus\{\m\}$, $|f|_{\m'} = 1$. La formule du produit (\emph{cf.} proposition \ref{formuleduproduit}) assure alors que $|f|_{\sigma_{\infty}} >1$.
\end{proof}

\begin{cor}
Soit~$\sigma\in\Sigma$. L'ensemble
$$\Ms(A)''_{\sigma} = \{a_{\sigma}^\eps,\, \eps\in\of{]}{0,l(\sigma)}{]}\}$$
est un ouvert de l'espace~$\Ms(A)$.
\end{cor}

Ce corollaire joint au lemme \ref{topobranche} permet de d\'ecrire la topologie au voisinage de tout point de l'espace~$\Ms(A)$ diff\'erent du point central. Int\'eressons-nous, \`a pr\'esent, \`a ce dernier.

\begin{lem}
Soit~$V$ un voisinage du point~$a_{0}$ dans~$\Ms(A)$. Il existe un sous-ensemble fini~$\Sigma_{V}$ de~$\Sigma$ tel que, pour tout \'el\'ement~$\sigma$ de $\Sigma\setminus\Sigma_{V}$, la branche~$\Ms(A)_{\sigma}$ soit contenue dans~$V$.
\end{lem}
\begin{proof}
Par d\'efinition de la topologie, il existe $r\in\N$, $f_{1},\ldots,f_{r} \in A$, $s_{1},\ldots,s_{r},t_{1},\ldots,t_{r}\in\R$ tels que la partie
$$W = \bigcap_{1\le i\le r} \left\{ b\in \Ms(A)\, |\, s_{i} < |f_{i}(b)| < t_{i} \right\}$$
soit un voisinage du point~$a_{0}$ dans~$V$. 

Soit~$i\in\cn{1}{r}$. Supposons, tout d'abord, que $f_{i}=0$. Nous avons alors \mbox{$f_{i}(a_{0})=0$} et donc $s_{i}<0$ et $t_{i}>0$. Posons~$\Sigma_{i}=\emptyset$. Nous avons alors
$$\left\{ b\in \Ms(A)\, |\, s_{i} < |f_{i}(b)| < t_{i} \right\} = \Ms(A) = \bigcup_{\sigma\notin\Sigma_{i}} \Ms(A)_{\sigma}.$$
Supposons, \`a pr\'esent, que~$f_{i}\ne 0$. Nous avons alors $|f_{i}(a_{0})|=1$ et donc $s_{i}<1$ et $t_{i}>1$. Posons
$$\Sigma_{i} = \{ \m \in \Sigma_{f}\, |\, f_{i}\in \m\} \cup \Sigma_{\infty}.$$
C'est un sous-ensemble fini de~$\Sigma$ qui v\'erifie
$$\left\{ b\in \Ms(A)\, |\, s_{i} < |f_{i}(b)| < t_{i} \right\} \supset \bigcup_{\sigma\notin\Sigma_{i}} \Ms(A)_{\sigma}.$$
Le sous-ensemble fini $\Sigma_{V}=\bigcup_{1\le i\le r} \Sigma_{i}$ satisfait alors la condition voulue. 
\end{proof}

\begin{lem}\label{voisa0}
Notons~$\Vs_{0}$ l'ensemble des parties de~$\Ms(A)$ qui v\'erifient les propri\'et\'es suivantes :  pour tout \'el\'ement~$V$ de~$\Vs_{0}$, il existe un sous-ensemble fini~$\Sigma_{V}$ de~$\Sigma$ et, pour tout \'el\'ement~$\sigma$ de~$\Sigma_{V}$, il existe un \'el\'ement~$\eps_{\sigma}$ de $\of{]}{0,l(\sigma)}{]}$ tels que
$$V = \bigcup_{\sigma\in\Sigma_{0}} \of{[}{a_{0},a_{\sigma}^{\eps_{\sigma}}}{[}\, \cup \bigcup_{\sigma\notin\Sigma_{0}} \Ms(A)_{\sigma}.$$
L'ensemble~$\Vs_{0}$ est un syst\`eme fondamental de voisinages ouverts du point~$a_{0}$ dans~$\Ms(A)$.
\end{lem}
\begin{proof}
Le fait que les \'el\'ements de~$\Vs_{0}$ soient des ouverts de~$\Ms(A)$ d\'ecoule des lemmes \ref{pasQquad} et \ref{Qquad}. Il nous suffit donc de montrer que tout voisinage du point central~$a_{0}$ contient un \'el\'ement de~$\Vs_{0}$.

Soit~$U$ un voisinage du point~$a_{0}$ dans~$\Ms(A)$. D'apr\`es le lemme pr\'ec\'edent, il existe un sous-ensemble fini~$\Sigma_{U}$ de~$\Sigma$ tel que, pour tout \'el\'ement~$\sigma$ de $\Sigma\setminus\Sigma_{U}$, la branche~$\Ms(A)_{\sigma}$ soit contenue dans~$U$. Pour tout \'el\'ement~$\sigma$ de~$\Sigma_{U}$, la partie $U\cap \Ms(A)_{\sigma}$ est un voisinage du point~$a_{0}$ dans~$\Ms(A)_{\sigma}$. Le lemme \ref{topobranche} assure alors qu'il existe un \'el\'ement~$\eps_{\sigma}$ de $\of{]}{0,l(\sigma)}{]}$ tel que la partie $U\cap \Ms(A)_{\sigma}$ contienne $\of{[}{a_{0},a_{\sigma}^{\eps_{\sigma}}}{[}$. On en d\'eduit que le voisinage~$U$ du point central~$a_{0}$ contient l'\'el\'ement de~$\Vs_{0}$ d\'efini par
$$\bigcup_{\sigma\in\Sigma_{0}} \of{[}{a_{0},a_{\sigma}^{\eps_{\sigma}}}{[}\, \cup \bigcup_{\sigma\notin\Sigma_{0}} \Ms(A)_{\sigma}.$$
\end{proof}

\begin{figure}[htb]
\begin{center}
\input{MZa0.pstex_t}
\caption{Un voisinage du point central $a_{0}$.}
\end{center}
\end{figure}
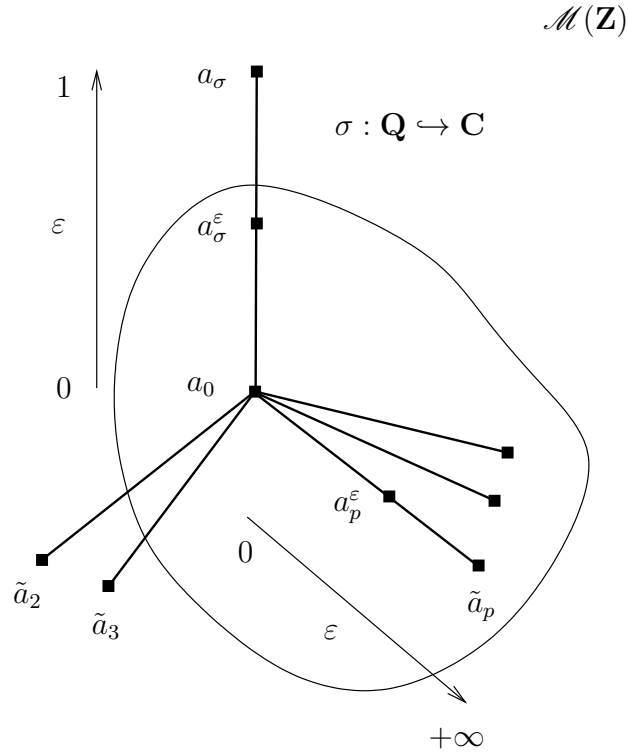

Regroupons, finalement, les r\'esultats obtenus.

\begin{cor}
Consid\'erons l'espace topologique
$$P=\bigsqcup_{\sigma\in\Sigma} \of{]}{0,l(\sigma)}{]}.$$
Notons~$\hat{P} = \hat{P}\cup\{\infty\}$ son compactifi\'e d'Alexandrov. L'application
$$\begin{array}{ccc}
P & \to & \Ms(A)\\
\eps_{\sigma}\in  \of{]}{0,l(\sigma)}{]} & \mapsto & a_{\sigma}^{\eps_{\sigma}}
\end{array}$$
se prolonge en un hom\'eomorphime
$$\hat{P} \xrightarrow[]{\sim} \Ms(A)$$
qui envoie le point~$\infty$ de~$\hat{P}$ sur le point central~$a_{0}$ de~$\Ms(A)$.
\end{cor}


Remarquons qu'\`a partir de la description de la topologie que nous venons de donner, on red\'emontre facilement la compacit\'e de l'espace $\Ms(A)$. D'autres propri\'et\'es sont v\'erifi\'ees. Nous les r\'esumons dans le th\'eor\`eme suivant.

\begin{thm}\index{Connexite par arcs au voisinage d'un point@Connexit\'e par arcs au voisinage d'un point!de MA@de $\Ms(A)$}
L'espace $\Ms(A)$ est compact, connexe par arcs et localement connexe par arcs.
\end{thm}

Remarquons que nous pouvons d\'ecrire facilement les parties connexes de l'espace~$\Ms(A)$. Deux cas se pr\'esentent. Si une partie connexe de~$\Ms(A)$ \'evite le point central~$a_{0}$, alors elle est contenue dans l'une des branches et est donc hom\'eomorphe \`a un intervalle. Si une partie connexe de~$\Ms(A)$ contient le point central~$a_{0}$, alors sa trace sur toute branche est une partie connexe, et donc hom\'eomorphe \`a un intervalle, contenant le point~$a_{0}$. On en d\'eduit le r\'esultat suivant.

\begin{prop}\label{connexeMA}
Une intersection de parties connexes de~$\Ms(A)$ est connexe.
\end{prop}

Indiquons pour finir un r\'esultat concernant les morphismes de changement de base.

\begin{thm}\label{chgtbase}
Soit $K'$ une extension finie de $K$. Notons $A'$ l'anneau des entiers de $K'$. Alors le morphisme 
$$\Ms(A') \to \Ms(A)$$
induit par l'injection $A \to A'$ est continu, ouvert, propre, surjectif et \`a fibres finies.
\end{thm}

\subsection{Faisceau structural}

Nous allons d\'ecrire les sections du faisceau structural~$\Os$ sur plusieurs types d'ouverts connexes de $\Ms(A)$. Auparavant, il est utile de calculer explicitement la norme uniforme sur certains compacts et le compl\'et\'e pour cette norme de l'anneau des fractions rationnelles sans p\^oles au voisinage du compact.

\subsubsection{Parties compactes}\label{partiescompactes}

Nous allons d\'ecrire ici toutes les parties compactes, connexes et non vides de~$\Ms(A)$. Soit~$L$ une telle partie. Nous allons distinguer plusieurs cas.

\begin{enumerate}
\item Il existe~$\sigma\in\Sigma_{\infty}$ tel que~$L$ soit contenue dans la branche $\sigma$-adique de~$\Ms(A)$.

\begin{enumerate}
\item La partie $L$ \'evite le point central $a_{0}$.

Dans ce cas, il existe $u,v\in\of{]}{0,1}{]}$, avec $u\le v$, tels que 
$$L=\of{[}{a_{\sigma}^u,a_{\sigma}^v}{]} = \{a_{\sigma}^{\eps},\, u\le \eps\le v\}.$$
Les fonctions rationnelles d\'efinies au voisinage de ce compact sont \mbox{$\Ks(L) = K$} et la norme uniforme est $\|.\|_{L} = \max(|.|_{\sigma}^u,|.|_{\sigma}^v)$. On en d\'eduit que $\Bs(L)\simeq \hat{K}_{\sigma}$. Attirons l'attention du lecteur sur le fait que l'isomorphisme pr\'ec\'edent est un isomorphisme de corps topologiques mais pas de corps norm\'es (sauf dans le cas o\`u $u=v$) !

\item La partie $L$ contient le point central $a_{0}$.

Il existe alors $v\in\of{[}{0,1}{]}$ tel que 
$$L=\of{[}{a_{0},a_{\sigma}^v}{]}.$$
Les fonctions rationnelles d\'efinies au voisinage de ce compact sont \mbox{$\Ks(L) = K$} et la norme uniforme est $\|.\|_{L} = \max(|.|_{0},|.|_{\sigma}^v)$. On en d\'eduit que $\Bs(L)\simeq K$.
\end{enumerate}

\item Il existe~$\m\in\Sigma_{f}$ tel que~$L$ soit contenue dans la branche $\m$-adique de~$\Ms(A)$.

\begin{enumerate}
\item La partie $L$ \'evite le point central $a_{0}$ et le point extr\^eme $\tilde{a}_{\m}$.

Il existe alors $u,v\in\of{]}{0,+\infty}{[}$, avec $u\le v$, tels que 
$$L=\of{[}{a_{\m}^u,a_{\m}^v}{]}.$$
Nous avons~$\Ks(L)=K$, $\|.\|_{L} = \max(|.|_{\m}^u,|.|_{\m}^v)$ et~\mbox{$\Bs(L)\simeq \hat{K}_{\m}$}.

\item La partie $L$ \'evite le point central $a_{0}$ et contient le point extr\^eme $\tilde{a}_{\m}$.

Il existe alors $u\in\of{]}{0,+\infty}{]}$ tel que 
$$L=\of{[}{a_{\m}^u,\tilde{a}_{\m}}{]}.$$
Dans ce cas, les \'elements de~$K$ peuvent avoir un p\^ole au point~$\tilde{a}_{\m}$ et nous avons donc~$\Ks(L)=A_{\m}$, $\|.\|_{L} = |.|_{\m}^u$ et~$\Bs(L)\simeq \hat{A}_{\m}$.

\item La partie $L$ contient le point central $a_{0}$ et \'evite le point extr\^eme $\tilde{a}_{\m}$.

Il existe alors $v\in\of{[}{0,+\infty}{[}$ tel que 
$$L=\of{[}{a_{0},a_{\m}^v}{]}.$$
Nous avons~$\Ks(L)=K$, $\|.\|_{L} = \max(|.|_{0},|.|_{\m}^u)$ et~\mbox{$\Bs(L)\simeq K$}.

\item La partie $L$ contient le point central $a_{0}$ et le point extr\^eme $\tilde{a}_{\m}$.

Dans ce cas, la partie $L$ est la branche $\m$-adique tout enti\`ere :
$$L = \Ms(A)_{\m}.$$
Nous avons~$\Ks(L)=A_{\m}$, $\|.\|_{L} =|.|_{0}$ et~$\Bs(L)\simeq A_{\m}$.
\end{enumerate}

\item La partie $L$ n'est contenue dans aucune branche de $\Ms(A)$.

D'apr\`es le raisonnement pr\'ec\'edant la proposition~\ref{connexeMA},
quel que soit \mbox{$\sigma\in\Sigma$}, il existe~$v_{\sigma}\in\of{[}{0,l(\sigma)}{]}$ tel que 
$$L = \bigcup_{\sigma\in\Sigma} \of{[}{a_{0},a_{\sigma}^{v_{\sigma}}}{]}.$$
Notons $\Sigma'= \{\m \in\Sigma_{f}\, |\, v_{\sigma}=l(\sigma)\}$. Nous avons alors
$$\Ks(L) = \bigcap_{\sigma\in\Sigma} \Ks(L_{\sigma}) = \bigcap_{\m\in\Sigma'} A_{\m}.$$
La norme uniforme sur cet anneau est
$$\|.\|_{L} = \max_{\sigma\in\Sigma} (\|.\|_{L_{\sigma}}) = \max \left( \max_{\sigma\in\Sigma} (|.|_{\sigma}^{v_{\sigma}}), |.|_{0}\right)$$
et nous avons donc
$$\Bs(L)=\Ks(L)=\bigcap_{\m\in\Sigma'} A_{\m}.$$
\end{enumerate}


Nous venons de d\'ecrire toutes les parties compactes et connexes de l'espace~$\Ms(A)$. Nous allons montrer qu'elles sont pro-rationnelles, au sens de la d\'efinition~\ref{defrationnel}. 

\begin{prop}\label{prorat}\index{Compact!de la base}
Toute partie compacte et connexe~$L$ de l'espace~$\Ms(A)$ est pro-rationnelle et donc spectralement convexe. En particulier, le morphisme naturel
$$\Ms(\Bs(L)) \to \Ms(\As)$$
induit un hom\'eomorphisme entre les espaces $\Ms(\Bs(L))$ et~$L$.
\end{prop}
\begin{proof}
Commen\c{c}ons par d\'emontrer le r\'esultat pour certaines parties compactes simples. Soient $\sigma\in\Sigma$ et $\eps\in\of{]}{0,l(\sigma)}{[}$. Consid\'erons le compact 
$$L = \Ms(A)\setminus \of{]}{a_{\sigma}^\eps,a_{\sigma}^{l(\sigma)}}{]}.$$
Supposons, tout d'abord, que~$\sigma\in\Sigma_{f}$ ou que $\sigma\in\Sigma_{\infty}$ et que le corps~$K$ n'est ni~$\Q$, ni un corps quadratique imaginaire. D'apr\`es le lemme \ref{pasQquad}, il existe alors un \'el\'ement~$f$ de~$A$ qui v\'erifie les conditions suivantes :
\begin{enumerate}[\it i)]
\item $|f|_{\sigma} < 1$ ;
\item $\forall \sigma'\ne \sigma$, $|f|_{\sigma'}\ge 1$.
\end{enumerate}
Nous avons alors
$$\left\{b\in \Ms(A)\, \big|\, |f(b)|\ge |f|_{\sigma}^\eps \right\} = L.$$
Le compact $L$ est rationnel.

Supposons, \`a pr\'esent, que le corps~$K$ est soit~$\Q$, soit un corps quadratique imaginaire et que~$\sigma=\sigma_{\infty}$. D'apr\`es le lemme~\ref{Qquad}, il existe alors un \'el\'ement~$f$ de~$A$ qui v\'erifie les conditions suivantes :
\begin{enumerate}[\it i)]
\item $|f|_{\sigma_{\infty}} > 1$ ;
\item $\forall \sigma'\ne \sigma$, $|f|_{\sigma'}\le 1$.
\end{enumerate}
Nous avons alors
$$\left\{b\in \Ms(A)\, \big|\, |f(b)|\le |f|_{\sigma_{\infty}}^\eps \right\} = L.$$
De nouveau, le compact $L$ est donc rationnel.

Consid\'erons, \`a pr\'esent, le compact 
$$M = \of{[}{a_{\sigma}^\eps,a_{\sigma}^{l(\sigma)}}{]}.$$
En utilisant la m\^eme fonction $f$ que pr\'ec\'edemment, nous pouvons \'ecrire, dans le premier cas, 
$$\left\{b\in \Ms(A)\, \big|\, |f(b)|\le |f|_{\sigma}^\eps \right\} = M,$$
et, dans le second,
$$\left\{b\in \Ms(A)\, \big|\, |f(b)|\ge |f|_{\sigma_{\infty}}^\eps \right\} = M.$$
Le compact $M$ est donc rationnel.

Puisque toutes les parties compactes et connexes de $\Ms(A)$ s'obtiennent comme intersection de compacts de l'un des deux types pr\'ec\'edents, la premi\`ere partie du r\'esultat est d\'emontr\'ee. Nous d\'eduisons la seconde partie du th\'eor\`eme \ref{compactrationnel}.
\end{proof}

\subsubsection{Parties ouvertes}\label{partiesouvertes}\index{Faisceau structural!sections sur MA@sections sur $\Ms(A)$|(}

Pour d\'eterminer les sections globales sur les ouverts de la base, il suffit \`a pr\'esent de recoller les compl\'et\'es pr\'ec\'edents. Introduisons tout d'abord une notation. 

\begin{defi}
Pour tout sous-ensemble~$\Sigma_{0}$ de~$\Sigma$, nous posons
$$A\left[\frac{1}{\Sigma_{0}}\right] = \left\{\left. \frac{a}{b}\in K\, \right|\, a\in A,\, b\in A^*,\, \exists \m\in\Sigma_{f}\cap\Sigma_{0}, b\in\m \right\}.$$
\end{defi}

\begin{rem}
Supposons que l'ensemble~$\Sigma_{0}$ pr\'ec\'edent est fini. Le localis\'e $A[1/\Sigma_{0}]$ poss\`ede alors une expression plus simple. Nous pouvons alors, en effet, consid\'erer le diviseur $\sum_{\m\in\Sigma_{f}\cap\Sigma_{0}} (\m)$ sur $\Spec(A)$. Puisque le groupe de Picard de $\Spec(A)$ est fini, ce diviseur est de torsion. Il existe donc $n\in\N^*$ et $f\in A$ tels que 
$$n \sum_{\m\in\Sigma_{f}\cap\Sigma_{0}} (\m) = (f).$$ 
Nous avons donc 
$$A\left[\frac{1}{\Sigma_{0}}\right]  = A\left[\frac{1}{f}\right].$$
\end{rem}

\bigskip

Soit $U$ un ouvert connexe et non vide de $\Ms(A)$. Comme pr\'ec\'edemment, nous allons distinguer plusieurs cas. 

\begin{enumerate}
\item Il existe~$\sigma\in\Sigma_{\infty}$ tel que~$U$ soit contenu dans la branche $\sigma$-adique de~$\Ms(A)$.

Alors, il existe $u,v\in\of{[}{0,1}{]}$, avec $u < v$, tels que 
$$U = \of{]}{a_{\sigma}^u,a_{\sigma}^v}{[} \textrm{ ou } \of{]}{a_{\sigma}^u,a_{\sigma}}{]}.$$ 
Dans les deux cas, nous avons $\Os(U) = \hat{K}_{\sigma}.$

\item Il existe~$\m\in\Sigma_{f}$ tel que~$U$ soit contenu dans la branche $\m$-adique de~$\Ms(A)$.

\begin{enumerate}
\item L'ouvert $U$ \'evite le point extr\^eme $\tilde{a}_{\m}$. 

Alors, il existe $u,v\in\of{[}{0,+\infty}{]}$, avec $u < v$, tels que 
$$U = \of{]}{a_{\m}^u,a_{\m}^v}{[}.$$ 
Comme pr\'ec\'edemment, nous avons $\Os(U)=\hat{K}_{\m}.$

\item L'ouvert $U$ contient le point extr\^eme $\tilde{a}_{\m}$.

Alors, il existe $u\in\of{[}{0,+\infty}{[}$ tel que 
$$U = \of{]}{a_{\m}^u,\tilde{a}_{\m}}{]}.$$ 
Dans ce cas, nous avons  $\Os(U) = \hat{A}_{\m}.$
\end{enumerate}

\item L'ouvert $U$ n'est contenu dans aucune branche de $\Ms(A)$.

Dans ce cas, c'est un voisinage du point central $a_{0}$ et il poss\`ede une \'ecriture de la forme
$$\begin{array}{rcl}
U &=&\disp \Ms(A) \setminus  \left(\bigcup_{\sigma\in\Sigma_{0}} \of{[}{a_{\sigma_{i}}^{u_{i}},a_{\sigma_{i}}^{l(\sigma_{i})}}{]}\right)\\
&=&\disp \left(\bigcup_{\sigma\in\Sigma_{0}} \of{[}{a_{0},a_{\sigma_{i}}^{u_{i}}}{[}\right) \cup  \left(\bigcup_{\sigma\notin \Sigma_{0}} \Ms(A)_{\sigma} \right),
\end{array}$$
o\`u~$\Sigma_{0}$ est un sous-ensemble fini de~$\Sigma$ et, pour tout \'el\'ement~$\sigma$ de~$\Sigma_{0}$, $u_{\sigma}$ est un \'el\'ement de $\of{]}{0,l(\sigma)}{]}$. Nous avons alors 
$$\Os(U) = A\left[\frac{1}{\Sigma_{0}}\right].$$

\end{enumerate}

\begin{figure}
\begin{center}
\input{sheaf.pstex_t}
\caption{Anneaux de sections globales.}
\end{center}
\end{figure}
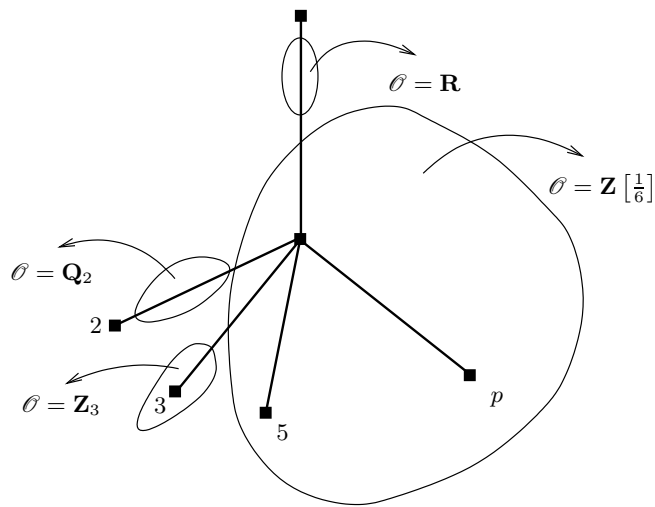

\index{Faisceau structural!sections sur MA@sections sur $\Ms(A)$|)}

Nous pouvons, \`a pr\'esent, d\'ecrire les anneaux locaux en les points de la base. Soit $b$ un point de $\Ms(A)$. Nous allons, de nouveau, distinguer plusieurs cas.\index{Anneau local en un point!de MA@de $\Ms(A)$}

\begin{enumerate}
\item Il existe $\sigma\in\Sigma$ tel que le point $b$ est un point interne de la branche $\sigma$-adique.

Dans ce cas, nous avons
$$\Os_{b} \simeq \hat{K}_{\sigma}.$$

\item Il existe $\m\in\Sigma_{f}$ tel que le point $b$ est le point extr\^eme $\tilde{a}_{\m}$.

Dans ce cas, nous avons
$$\Os_{\tilde{a}_{\m}} \simeq \hat{A}_{\m}.$$

\item Le point $b$ est le point central $a_{0}$ de $\Ms(A)$.

Nous avons alors
$$\Os_{a_{0}} \simeq K.$$
\end{enumerate}

\begin{rem}
La topologie de l'espace~$\Ms(A)$ laisse penser que c'est, en quelque sorte, un espace ad\'elique. La connaissance du faisceau structural permet de pr\'eciser cette id\'ee. Consid\'erons le morphisme d'inclusion
$$j : \Ms(A)\setminus\{a_{0}\} \to \Ms(A).$$ 
Le germe $(j_{*}\Os)_{a_{0}}$ est isomorphe \`a l'anneau des ad\`eles.\index{Adeles@Ad\`eles|see{Corps de nombres}}\index{Corps de nombres!adeles@ad\`eles}
\end{rem}

\bigskip

Gr\^ace aux descriptions explicites que nous avons obtenues, il est d\'esormais facile de montrer que les anneaux locaux de l'espace~$\Ms(A)$ qui sont des anneaux de valuation discr\`ete -- ce sont exactement les anneaux locaux en les points extr\^emes -- satisfont la condition~(U), au sens de la d\'efinition \ref{condU}. 

\begin{lem}\index{Condition (U)!pour les points de MA@pour les points de $\Ms(A)$}
Soit $\m\in\Sigma_{f}$. L'anneau de valuation discr\`ete~$\Os_{\tilde{a}_{\m}}$ satisfait la condition~(U).
\end{lem}
\begin{proof}
Consid\'erons l'uniformisante~$\pi_{\m}$ de l'anneau de valuation discr\`ete $\Os_{\tilde{a}_{\m}} = \hat{A}_{\m}$. Elle est d\'efinie sur l'ouvert $V=\of{]}{a_{0},\tilde{a}_{\m}}{]}$. L'ensemble
$$\Ws = \left\{\of{[}{a_{\m}^\eps,\tilde{a}_{\m}}{]},\, \eps>0\right\}$$
est un syst\`eme fondamental de voisinages compacts du point~$b$ dans~$V$.

Soit~$\eps>0$. Posons $W=\of{[}{a_{\m}^\eps,\tilde{a}_{\m}}{]}$. Les descriptions pr\'ec\'edentes montrent que nous avons
$$\Bs(W)=\hat{A}_{\m}$$
et
$$\forall f\in \hat{A}_{\m},\, \|f\|_{W}=|f|_{\m}^\eps.$$
Soit~$f$ un \'el\'ement de~$\Bs(W)$ tel que $f(\tilde{a}_{\m})=0$. Cela signifie que~$f$ est divisible par~$\pi_{\m}$ dans $\Os_{\tilde{a}_{\m}}$, c'est-\`a-dire dans~$\hat{A}_{\m}$, et donc dans~$\Bs(W)$. Il existe donc un \'el\'ement~$g$ de~$\hat{A}_{\m}$ tel que $f=\pi_{\m}\, g$. En outre, nous avons
$$\|f\|_{W} = |f|_{\m}^\eps = |\pi_{\m}|_{\m}^\eps\, |g|_{\m}^\eps = |\pi_{\m}|_{\m}^\eps\, \|g\|_{W}.$$
Par cons\'equent, l'uniformisante~$\pi_{\m}$ v\'erifie la condition~$\textrm{U}_{W}$. On en d\'eduit le r\'esultat attendu.
\end{proof}

\bigskip

Les r\'esultats qui pr\'ec\`edent permettent \'egalement de d\'ecrire explicitement les anneaux de fonctions d\'efinies au voisinages des parties compactes de~$\Ms(A)$. Nous obtenons en particulier le r\'esultat suivant.

\begin{prop}\label{isoBO}
Soit~$V$ une partie compacte et connexe de l'espace~$\Ms(A)$ qui n'est pas r\'eduite \`a un point extr\^eme. Alors le morphisme naturel 
$$\Ks(V)\to \Os(V)$$ 
se prolonge en un isomorphisme
$$\Bs(V)\xrightarrow[]{\sim} \Os(V).$$ 
En particulier, pour tout point~$b$ de l'espace~$\Ms(A)$ qui n'est pas extr\^eme, le morphisme naturel
$$\Os_{b} \to \Hs(b)$$
est un isomorphisme.
\end{prop}

De la description explicite des anneaux locaux d\'ecoule un autre fait important, qui permet de ramener l'\'etude de n'importe quel point rationnel d'une fibre d'un espace analytique affine au-dessus de~$\As$ \`a celle du point~$0$ de cette m\^eme fibre.

\begin{lem}\label{translation}
Soit~$b$ un point de l'espace~$\Ms(A)$. Le morphisme naturel
$$\Os_{b} \to \Hs(b)$$
est surjectif.

Soient $n$ un entier, $\E{n}{A}$ l'espace analytique de dimension~$n$ au-dessus de~$A$ et $\pi : \E{n}{A} \to \Ms(A)$ le morphisme naturel de projection. Soit~$x$ un point rationnel de la fibre~$\pi^{-1}(b)$. Il existe un voisinage ouvert~$U$ du point~$b$ dans~$\Ms(A)$ et un automorphisme~$\varphi$ de~$\pi^{-1}(U)$ qui fait commuter le diagramme
$$\xymatrix{
\pi^{-1}(U) \ar[rr]^{\varphi}_{\sim} \ar[rd]_\pi & & \pi^{-1}(U) \ar[ld]^{\pi}\\
& U &
}$$
et envoie le point~$x$ de la fibre~$\pi^{-1}(b)$ sur le point~$0$ de cette m\^eme fibre.
\end{lem}
\begin{proof}
La premi\`ere partie du r\'esultat provient de la description explicite des anneaux locaux et des corps r\'esiduels compl\'et\'es.

Passons \`a la deuxi\`eme partie. Il existe des \'el\'ements $\alpha_{1},\ldots,\alpha_{n}$ de~$\Hs(b)$ tels que le point~$x$ de la fibre~$\pi^{-1}(b)$ soit d\'efini par les \'equations
$$(T_{1}-\alpha_{1})(b) = \cdots = (T_{n}-\alpha_{n})(b)=0.$$ 
D'apr\`es la premi\`ere partie, pour tout \'el\'ement~$i$ de~$\cn{1}{n}$, il existe un \'el\'ement~$\beta_{i}$ de~$\Os_{b}$ dont l'image dans le corps r\'esiduel compl\'et\'e~$\Hs(b)$ est \'egal \`a~$\alpha_{i}$. Choisissons un voisinage ouvert~$U$ du point~$b$ dans~$\Ms(A)$ sur lequel les fonctions $\beta_{1},\ldots,\beta_{n}$ sont d\'efinies. Nous pouvons alors choisir comme automorphisme la translation de vecteur $(\beta_{1},\ldots,\beta_{n})$ au-dessus de~$\pi^{-1}(U)$.
\end{proof}

\subsubsection{Bord de Shilov}\label{borddeShilovbase}

Commencer par rappeler la notion de bord de Shilov et celles qui lui sont li\'ees. 

\begin{defi}\index{Bord analytique}\index{Bord de Shilov|see{Bord analytique}}
Soit~$(\As,\|.\|)$ un anneau de Banach. Nous dirons qu'une partie ferm\'ee~$\Gamma$ de~$\Ms(\As,\|.\|)$ est un {\bf bord analytique de l'anneau norm\'e $(\As,\|.\|)$} si elle v\'erifie la condition suivante :
$$\forall f\in\As,\, \|f\|_{\Ms(\As,\|.\|)} = \|f\|_{\Gamma}.$$

Nous appellerons {\bf bord de Shilov} de l'anneau de Banach~$(\As,\|.\|)$ le plus petit bord, pour la relation d'inclusion, de l'anneau de Banach~$(\As,\|.\|)$, s'il existe.

Soient~$n$ un nombre entier positif et~$V$ une partie compacte et spectralement convexe de l'espace analytique $\E{n}{\As}$. Par d\'efinition (\emph{cf.}~\ref{defspconvexe}), le morphisme naturel
$$\Ms(\Bs(V)) \to \E{n}{\As}$$
induit un hom\'eomorphisme entre les espaces $\Ms(\Bs(V))$ et~$V$. Nous appellerons {\bf bord analytique (respectivement bord de Shilov) du compact~$V$} l'image par cet hom\'eomorphisme d'un bord analytique (respectivement du bord de Shilov, s'il existe) de l'anneau de Banach~$(\Bs(V),\|.\|_{V})$.  



\end{defi}

\begin{rem}
Le lemme de Zorn assure que tout anneau de Banach poss\`ede un bord analytique minimal.
\end{rem}

Signalons qu'A.~Escassut et~N.~Ma\"inetti ont prouv\'e l'existence du bord de Shilov dans de nombreux cas (\emph{cf.}~\cite{Shilov}, th\'eor\`eme~C).

\begin{thm}[Escassut, Ma\"inetti]\label{escassutmainetti}
Soit~$(k,|.|)$ un corps valu\'e et complet dont la valuation n'est pas triviale. Soit~$\As$ une $k$-alg\`ebre de Banach munie d'une norme d'alg\`ebre~$\|.\|$ qui induit la valeur absolue~$|.|$ sur~$k$. Alors l'alg\`ebre $(\As,\|.\|)$ poss\`ede un bord de Shilov.
\end{thm}

\medskip

Int\'eressons-nous, \`a pr\'esent, au bord de Shilov des parties compactes et con\-nexes de l'espace~$\Ms(A)$. Cela a un sens puisque ces parties sont pro-rationnelles (\emph{cf.}~proposition~\ref{prorat}) et donc spectralement convexes (\emph{cf.}~th\'eor\`eme~\ref{compactrationnel}). En reprenant les r\'esultats des paragraphes~\ref{partiescompactes} et~\ref{partiesouvertes}, l'on montre simplement que les parties compactes et connexes de~$\Ms(A)$ poss\`edent un bord de Shilov. Nous pouvons en donner une description explicite.\index{Bord analytique!des compacts de MA@des compacts de $\Ms(A)$|(}

\begin{enumerate}
\item Pour tout \'el\'ement~$\m$ de~$\Sigma_{f}$ et tous \'el\'ements~$u$ et~$v$ de~$\R_{+}$ v\'erifiant l'in\'egalit\'e $u\le v$, la partie compacte $\of{[}{a_{\m}^u,a_{\m}^v}{]}$ poss\`ede un bord de Shilov \'egal \`a l'ensemble $\{a_{\m}^u,a_{\m}^v\}$.
\item Pour tout \'el\'ement~$\m$ de~$\Sigma_{f}$ et tout \'el\'ement~$u$ de~$\R_{+}$, la partie compacte $\of{[}{a_{\m}^u,\tilde{a}_{\m}}{]}$ poss\`ede un bord de Shilov \'egal au singleton $\{a_{\m}^u\}$.
\item Pour tout \'el\'ement~$\m$ de~$\Sigma_{f}$, la partie compacte~$\{\tilde{a}_{\m}\}$ poss\`ede un bord de Shilov \'egal au singleton~$\{\tilde{a}_{\m}\}$.
\item Pour tout \'el\'ement~$\sigma$ de~$\Sigma_{\infty}$ et tous \'el\'ements~$u$ et~$v$ de~$\of{[}{0,1}{]}$ v\'erifiant l'in\'egalit\'e $u\le v$, la partie compacte $\of{[}{a_{\m}^u,a_{\m}^v}{]}$ poss\`ede un bord de Shilov \'egal \`a l'ensemble $\{a_{\m}^u,a_{\m}^v\}$.
\item Lorsque la partie compacte et connexe n'est pas contenue dans une branche, le r\'esultat est plus difficile \`a \'etablir. Les lemmes~\ref{pasQquad} et~\ref{Qquad} permettent cependant d'y parvenir rapidement. Soit~$L$ une partie compacte et connexe de~$B$ qui n'est pas contenue dans une branche. Alors, il existe un \'el\'ement~$(v_{\sigma})_{\sigma\in\Sigma}$ de~$\prod_{\sigma\in\Sigma} \of{[}{0,l(\sigma)}{]}$ tel que l'on ait l'\'egalit\'e
$$L = \bigcup_{\sigma\in\Sigma} \of{[}{a_{0},a_{\sigma}^{v_{\sigma}}}{]}.$$  
Posons
$$\Sigma_{0} = \{\sigma\in\Sigma_{f}\, |\, 0 < v_{\sigma} < +\infty \}.$$
Alors, la partie compacte~$L$ de~$B$ poss\`ede un bord de Shilov \'egal \`a l'ensemble
$$\bigcup_{\sigma\in\Sigma_{0}} \{a_{\sigma}^{v_{\sigma}}\} \cup \bigcup_{\sigma\in\Sigma_{\infty}} \{a_{\sigma}^{v_{\sigma}}\}.$$
\end{enumerate}

Ces descriptions explicites permettent d'obtenir le r\'esultat suivant.

\begin{prop}\label{Shilovbase}
Soit~$V$ une partie compacte et connexe de l'espace~$B$. Cette partie poss\`ede un bord de Shilov~$\Gamma_{V}$. C'est un ensemble fini. 

En outre, pour tout point~$\gamma$ de~$\Gamma_{V}$, il existe un \'el\'ement~$f$ de~$\Ks(V)$ v\'erifiant la propri\'et\'e suivante :
$$|f(\gamma)|=\|f\|_{V} \textrm{ et } \forall b\in V\setminus\{\gamma\},\, |f(b)|< \|f\|_{V}.$$

Si la partie compacte et connexe~$V$ n'est pas r\'eduite \`a un point extr\^eme, alors, en tout point~$\gamma$ de~$\Gamma_{V}$, l'anneau local~$\Os_{\gamma}$ est un corps.

\end{prop}

\index{Bord analytique!des compacts de MA@des compacts de $\Ms(A)$|)}




Introduisons une nouvelle d\'efinition.

\begin{defi}\index{Algebriquement trivial@Alg\'ebriquement trivial}
Soient~$(\As,\|.\|)$ un anneau de Banach et~$n$ un nombre entier positif. Posons~$X=\E{n}{(\As,\|.\|)}$ et notons~$\Os_{X}$ le faisceau structural sur cet espace. Nous dirons qu'une partie~$S$ de l'espace analytique~$X$ est {\bf alg\'ebriquement triviale} si, pour tout point~$x$ de~$S$, l'anneau local~$\Os_{X,x}$ est un corps.
\end{defi}


\begin{cor}\index{Bord analytique!au voisinage d'un point de MA@au voisinage d'un point de $\Ms(A)$}
Tout point de l'espace~$\Ms(A)$ poss\`ede un syst\`eme fondamental de voisinages compacts et connexes dont le bord de Shilov est une partie finie et alg\'ebriquement triviale.
\end{cor}

\index{Spectre analytique!de A@de $A$|)}

\section{Faisceau structural sur l'espace affine}\label{fssea}

\index{Espace affine analytique!sur A@sur $A$|(}

Dans la suite de ce chapitre, nous fixons un entier positif~$n$. Nous posons
$$B=\Ms(A) \textrm{ et } X=\E{n}{A}.$$ 
Les faisceaux structuraux sur ces espaces seront respectivement not\'es~$\Os_{B}$ et~$\Os_{X}$. Lorsqu'aucune confusion ne peut en d\'ecouler, nous les noterons simplement~$\Os$.

Nous noterons encore 
$$\pi : X \to B$$ 
\newcounter{pi}\setcounter{pi}{\thepage}
le morphisme de projection induit par le morphisme naturel $A \to A[T_{1},\ldots,T_{n}]$. Pour toute partie~$V$ de~$B$, nous posons
$$X_{V} = \pi^{-1}(V)$$ 
et, pour tout point~$b$ de~$B$,
$$X_{b} = \pi^{-1}(b).$$ 

Introduisons encore quelques \'el\'ements de terminologie pour l'espace affine de dimension~$n$ au-dessus de~$\Ms(A)$, dans la lign\'ee de la d\'efinition \ref{branches}.

\begin{defi}\index{Partie sigma-adique@Partie $\sigma$-adique}\index{Partie sigma-adique@Partie $\sigma$-adique!ouverte}\index{Partie sigma-adique@Partie $\sigma$-adique!semi-ouverte}
\newcounter{partie}\setcounter{partie}{\thepage}
Pour $\sigma\in\Sigma$, nous appellerons {\bf partie $\sigma$-adique} de~$X$ (respectivement {\bf partie $\sigma$-adique ouverte} de~$X$, {\bf partie $\sigma$-adique semi-ou\-ver\-te} de~$X$), et noterons~$X_{\sigma}$ (respectivement~$X'_{\sigma}$, $X''_{\sigma}$), l'image r\'eciproque par la projection~$\pi$ de la branche $\sigma$-adique (respectivement branche $\sigma$-adique ouverte, branche $\sigma$-adique semi-ouverte) de~$\Ms(A)$.

Nous appellerons {\bf fibre centrale}\index{Fibre!centrale} de $X$, et noterons $X_{0}$, la fibre de~$\pi$ au-dessus du point central de $\Ms(A)$. Nous appellerons {\bf fibre extr\^eme}\index{Fibre!extreme@extr\^eme} de~$X$ toute fibre de~$\pi$ au-dessus d'un point extr\^eme de~$\Ms(A)$. Pour $\m\in\Sigma_{f}$, nous noterons $\tilde{X}_{\m} = \pi^{-1}(\tilde{a}_{\m})$. Finalement, nous appellerons {\bf fibre interne}\index{Fibre!interne} de~$X$ toute fibre de~$\pi$ au-dessus d'un point interne de~$\Ms(A)$. Nous appellerons {\bf point interne} de~$X$\index{Point!interne} tout point d'une telle fibre.
\end{defi}

\subsection{Anneaux locaux}\index{Anneau local en un point!deploye@d\'eploy\'e!de An@de \E{n}{A}|(}

Au th\'eor\`eme \ref{anneaulocal}, nous avons d\'ecrit les anneaux locaux en les points d\'eploy\'es en fonction d'anneaux de sections sur la base. Gr\^ace aux r\'esultats \'etablis \`a la section pr\'ec\'edente, nous pouvons pr\'eciser cette description dans le cas o\`u la base est le spectre d'un anneau d'entiers de corps de nombres. Soit~$b$ un point de~$B$. Soient $\alpha_{1},\ldots,\alpha_{n}$ des \'el\'ements de~$\Os_{B,b}$. Soient~$I$ une partie de~$\cn{1}{n}$ et $(r_{i})_{i\in I}$ une famille de $\R_{+}^*$ dont l'image dans l'espace vectoriel $\Q\otimes_{\Z} (\R_{+}^*/|\Hs(b)^*|)$ est libre. Notons $J=\cn{1}{n}\setminus I$ et, pour $i\in J$, posons $r_{i}=0$. Notons~$x$ l'unique point de la fibre~$X_{b}$ qui v\'erifie 
$$\forall i\in\cn{1}{n}, |(T_{i}-\alpha_{i})(x)|=r_{i}.$$
Notons
$$S = \{(s_{1},\ldots,s_{n}) \in\R_{+}^n,\, |\, \forall i\in I,\, s_{i}\in\of{]}{0,r_{i}}{[},\ \forall i\in J,\, s_{i}=0\}$$
et
$$T = \{(t_{1},\ldots,t_{n}) \in\R_{+}^n,\, |\, \forall i\in I,\, t_{i}>r_{i},\ \forall i\in J,\, t_{i}>0\}.$$

\begin{prop}
Supposons que le point~$b$ est un point interne de l'espace~$B$. Il existe un \'el\'ement~$\sigma$ de~$\Sigma$ et un nombre r\'eel~$\eps>0$ tels que $b=a_{\sigma}^\eps$. Notons~$L_{i}$ l'anneau compos\'e des s\'eries \`a coefficients dans~$\hat{K}_{\sigma}$ de la forme
$$\sum_{\bk\in\Z^n} a_{\bk}\, (\bT-\alphab)^\bk$$
qui v\'erifient la condition suivante : il existe des \'el\'ements~$\bs$ de~$S$ et~$\bt$ de~$T$ tels que la famille
$$\left(|a_{\bk}|_{\sigma}^\eps\, \bmax(\bs^\bk,\bt^\bk)\right)_{\bk\in\Z^n}$$
est sommable. Une telle famille v\'erifie en particulier la condition suivante : pour tout \'el\'ement~$i$ de~$J$ et tout \'el\'ement~$\bk$ de~$\Z^n$ v\'erifiant~$k_{i}<0$, nous avons~$a_{\bk}=0$.

Le morphisme naturel $A[\bT] \to \Os_{X,x}$ induit un isomorphisme
$$L_{i} \xrightarrow[]{\sim} \Os_{X,x}.$$
\end{prop}
\begin{proof}
Nous supposerons que le nombre r\'eel~$\eps$ appartient \`a l'intervalle~$\of{]}{0,l(\sigma)}{[}$. Le cas o\`u~$\eps=l(\sigma)$, et donc $\sigma\in\Sigma_{\infty}$, ne pr\'esente pas de difficult\'e suppl\'ementaire et nous ne le traiterons pas.

La famille
$$\Vs = (V_{\alpha,\beta} = \of{[}{a_{\sigma}^\alpha,a_{\sigma}^\beta}{]})_{0<\alpha<\eps<\beta< l(\sigma)}$$
est un syst\`eme fondamental de voisinages du point~$a_{\sigma}^\eps$ dans~$B$. En outre, quel que soient les nombres r\'eels~$\alpha$ et~$\beta$ v\'erifiant $0<\alpha<\eps<\beta< l(\sigma)$, nous avons
$$(\Bs(V_{\alpha,\beta}),\|.\|_{V_{\alpha,\beta}}) =  \left(\hat{K}_{\sigma},\max(|.|_{\sigma}^\alpha,|.|_{\sigma}^\beta)\right).$$
D'apr\`es le th\'eor\`eme~\ref{anneaulocal}, le morphisme $\As[\bT] \to \Os_{X,x}$ induit un isomorphisme
$$\varinjlim_{V,\bs,\bt}\, \Bs(V)\of{\la}{\bs \le |\bT-\alphab|\le \bt}{\ra} \xrightarrow{\sim} \Os_{X,x},$$
o\`u $V$ parcourt la famille~$\Vs$, $\bs$ l'ensemble~$S$ et~$\bt$ l'ensemble~$T$. Soient~$\alpha,\beta\in\R$ v\'erifiant $0<\alpha<\eps<\beta<l(\sigma)$, $\bs\in S$ et~$\bt\in T$. Soit~$f$ un \'el\'ement de l'anneau $\Bs(V_{\alpha,\beta})\of{\la}{\bs \le |\bT-\alphab|\le \bt}{\ra}$. Il existe alors une famille $(a_{\bk})_{\bk\in\Z^n}$ d'\'el\'ements de~$\hat{K_{\sigma}}$ telle que 
$$\textrm{la famille } \left(\max(|a_{\bk}|_{\sigma}^\alpha,|a_{\bk}|_{\sigma}^\beta)\, \bmax(\bs^\bk,\bt^\bk)\right)_{\bk\in\Z^n} \textrm{ est sommable}$$
et telle que l'on ait l'\'egalit\'e
$$f = \sum_{\bk\in\Z^n} a_{\bk}\, \bT^\bk.$$
Pour conclure, il reste \`a constater que l'ensemble des familles $(a_{\bk})_{\bk\in\Z^n}$ d'\'el\'ements de~$\hat{K_{\sigma}}$ pour lesquelles il existe des \'el\'ements~$\bs$ de~$S$ et~$\bt$ de~$T$ tels que
$$\textrm{la famille } \left(\max(|a_{\bk}|_{\sigma}^\alpha,|a_{\bk}|_{\sigma}^\beta)\, \bmax(\bs^\bk,\bt^\bk)\right)_{\bk\in\Z^n} \textrm{ est sommable}$$
est identique \`a l'ensemble des familles $(a_{\bk})_{\bk\in\Z^n}$ d'\'el\'ements de~$\hat{K_{\sigma}}$ pour lesquelles il existe des \'el\'ements~$\bs$ de~$S$ et~$\bt$ de~$T$ tels que
$$\textrm{la famille } \left(|a_{\bk}|_{\sigma}^\eps\, \bmax(\bs^\bk,\bt^\bk)\right)_{\bk\in\Z^n} \textrm{ est sommable}.$$
\end{proof}

\begin{rem}
Supposons que le point~$b$ est un point interne de l'espace~$B$. La description explicite que nous venons d'obtenir montre que le morphisme naturel
$$\Os_{X,x} \to \Os_{X_{b},x}$$
est un isomorphisme. Ce r\'esultat vaut, en fait, pour tous les points des fibres internes, ainsi que nous le d\'emontrerons plus tard (\emph{cf.}~proposition~\ref{isointerne}). 
\end{rem}

\begin{prop}
Supposons que le point~$b$ est un point extr\^eme de l'espace~$B$. Il existe un \'el\'ement~$\m$ de~$\Sigma_{f}$ tel que $b=\tilde{a}_{\m}$. Notons~$L_{e}$ l'anneau compos\'e des s\'eries \`a coefficients dans~$\hat{A}_{\m}$ de la forme
$$\sum_{\bk\in\Z^n} a_{\bk}\, (\bT-\alphab)^\bk$$
qui v\'erifient la condition suivante : il existe des \'el\'ements~$\eps$ de~$\R_{+}^*$, $\bs$ de~$S$ et~$\bt$ de~$T$ tels que la famille
$$\left(|a_{\bk}|_{\m}^\eps\, \bmax(\bs^\bk,\bt^\bk)\right)_{\bk\in\Z^n}$$
est sommable. Une telle famille v\'erifie en particulier la condition suivante : pour tout \'el\'ement~$i$ de~$J$ et tout \'el\'ement~$\bk$ de~$\Z^n$ v\'erifiant~$k_{i}<0$, nous avons~$a_{\bk}=0$.

Le morphisme naturel $A[\bT] \to \Os_{X,x}$ induit un isomorphisme
$$L_{e} \xrightarrow[]{\sim} \Os_{X,x}.$$
\end{prop}
\begin{proof}
La famille
$$\Vs = (V_{\eps} = \of{[}{a_{\m}^\eps,\tilde{a}_{\m}}{]})_{\eps>0}$$
est un syst\`eme fondamental de voisinages du point~$\tilde{a}_{\m}$ dans~$B$. En outre, pour tout \'el\'ement~$\eps$ de~$\R_{+}^*$, nous avons
$$(\Bs(V_{\eps}),\|.\|_{V_{\eps}}) =  \left(\hat{A}_{\m},|.|_{\m}^\eps\right).$$
D'apr\`es le th\'eor\`eme~\ref{anneaulocal}, le morphisme $\As[\bT] \to \Os_{X,x}$ induit un isomorphisme
$$\varinjlim_{V,\bs,\bt}\, \Bs(V)\of{\la}{\bs \le |\bT-\alphab|\le \bt}{\ra} \xrightarrow{\sim} \Os_{X,x},$$
o\`u $V$ parcourt la famille~$\Vs$, $\bs$ l'ensemble~$S$ et~$\bt$ l'ensemble~$T$. On en d\'eduit le r\'esultat annonc\'e.
\end{proof}

\begin{cor}\label{descratext}
Supposons qu'il existe un \'el\'ement~$\m$ de~$\Sigma_{f}$ tel que le point~$b=\tilde{a}_{\m}$ et que~$I=\emptyset$. Le point~$x$ est alors un point rationnel de la fibre extr\^eme~$\tilde{X}_{\m}$. Le morphisme naturel $A[\bT] \to \Os_{X,x}$ induit un isomorphisme
$$\hat{A}_{\m}[\![\bT-\alphab]\!] \xrightarrow[]{\sim} \Os_{X,x}.$$
\end{cor}
\begin{proof}
Reprenons les notations de la proposition pr\'ec\'edente. Nous souhaitons montrer que l'anneau~$L_{e}$ n'est autre que l'anneau~$\hat{A}_{\m}[\![\bT-\alphab]\!]$. Tout d'abord, puisque~$I$ est vide, nous disposons de l'inclusion
$$L_{e} \subset \hat{A}_{\m}[\![\bT-\alphab]\!].$$
R\'eciproquement, soit
$$f =  \sum_{\bk\in\Z^n} a_{\bk}\, (\bT-\alphab)^\bk$$
un \'el\'ement de~$\hat{A}_{\m}[\![\bT-\alphab]\!]$. Soient~$\eps>0$ et $t_{1},\ldots,t_{n}\in\of{]}{0,1}{[}$. Le $n$-uplet~$(t_{1},\ldots,t_{n})$ appartient \`a~$T$. Puisque~$I$ est vide, l'ensemble~$S$ est r\'eduit au $n$-uplet nul. Remarquons finalement que, pour tout \'el\'ement~$\bk$ de~$\N^n$, nous avons~$|a_{\bk}|_{\m}^\eps\le 1$. On en d\'eduit que la famille
$$\left(|a_{\bk}|_{\m}^\eps\, \bmax(\bs^\bk,\bt^\bk)\right)_{\bk\in\Z^n} = \left(|a_{\bk}|_{\m}^\eps\, \bt^\bk\right)_{\bk\in\Z^n}$$
est sommable et donc que l'\'el\'ement~$f$ appartient \`a~$L_{e}$.
\end{proof}

Dans le cas de la droite, nous pouvons simplifier la description. Pour traiter ce cas, nous supposerons que~$n=1$ et supprimerons les indices des notations. 

\begin{cor}\label{descriptionextreme3}
Supposons que~$n=1$, que~$r<1$ et que le point~$b$ est un point extr\^eme de l'espace~$B$. Il existe un \'el\'ement~$\m$ de~$\Sigma_{f}$ tel que $b=\tilde{a}_{\m}$. Notons~$L_{e}^{(1)}$ l'anneau compos\'e des s\'eries \`a coefficients dans~$\hat{A}_{\m}$ de la forme
$$\sum_{k\in\Z} a_{k}\, (T-\alpha)^k$$
telles que le rayon de convergence de la s\'erie
$$\sum_{k<0} a_{k}\, U^k$$
soit strictement sup\'erieur \`a~$1$. C'est un anneau de valuation discr\`ete d'id\'eal maximal~$(\pi_{\m})$ et de corps r\'esiduel $\tilde{k}_{\m}(\!(T)\!)$. Le morphisme naturel $A[T] \to \Os_{X,x}$ induit des isomorphismes
$$L_{e}^{(1)} \xrightarrow[]{\sim} \Os_{X,x}$$
et
$$\tilde{k}_{\m}(\!(T)\!)  \xrightarrow[]{\sim} \kappa(x) \xrightarrow[]{\sim} \Hs(x).$$
\end{cor}
\begin{proof}
Commen\c{c}ons par nous int\'eresser \`a l'anneau local~$\Os_{X,x}$. Nous savons qu'il est isomorphe \`a l'anneau compos\'e des s\'eries \`a coefficients dans~$\hat{A}_{\m}$ de la forme
$$\sum_{k\in\Z} a_{k}\, (T-\alpha)^k$$
qui v\'erifient la condition suivante : il existe des \'el\'ements~$\eps$ de~$\R_{+}^*$, $s$ de~$\of{]}{0,r}{[}$ et~$t$ de~$\of{]}{r,+\infty}{[}$ tels que la famille
$$\left(|a_{k}|_{\m}^\eps\, \max(s^k,t^k)\right)_{k\in\Z}$$
soit sommable. Cette condition est \'equivalente \`a la conjonction des deux conditions suivantes :
\begin{enumerate}[a)]
\item il existe~$\eps>0$ et~$t>r$ tel que la famille $(|a_{k}|_{\m}^\eps\, t^k)_{k\ge 0}$ est sommable ;
\item il existe~$\eps>0$ et~$s\in\of{]}{0,r}{[}$ tel que la famille $(|a_{k}|_{\m}^\eps\, s^k)_{k< 0}$ est sommable.
\end{enumerate}
La condition a) est toujours satisfaite. En effet, la suite $(|a_{k}|_{\m})_{k\ge 0}$ est born\'ee. Le rayon de convergence de la s\'erie $\sum_{k\ge 0} a_{k}\, U^k$ est donc sup\'erieur \`a~$1$. On v\'erifie ensuite sans peine que la condition b) est \'equivalente \`a celle de l'\'enonc\'e du corollaire.

Pour d\'emontrer l'assertion finale, il suffit de remarquer que le corps $\kappa(x)\simeq \tilde{k}_{\m}(\!(T)\!)$ est complet pour la valuation $T$-adique et donc pour la valeur absolue associ\'ee au point~$x$. On en d\'eduit que le morphisme naturel 
$$\kappa(x) \to \Hs(x)$$
est un isomorphisme.
\end{proof}

\begin{prop}
Supposons que le point~$b$ est le point central~$a_{0}$ de l'espace~$B$. Notons~$L_{c}$ l'anneau compos\'e des s\'eries \`a coefficients dans~$K$ de la forme
$$\sum_{\bk\in\Z^n} a_{\bk}\, (\bT-\alphab)^\bk$$
qui v\'erifient les conditions suivantes : 
\begin{enumerate}[\it i)]
\item il existe un sous-ensemble fini~$\Sigma_{0}$ de~$\Sigma$ contenant~$\Sigma_{\infty}$ tel que, quel que soit~$\bk$ dans~$\Z^n$, l'\'el\'ement~$a_{\bk}$ appartient \`a~$A[1/\Sigma_{0}]$ ;
\item quel que soit~$\sigma$ dans~$\Sigma_{0}$, il existe des \'el\'ements~$\eps$ de~$\of{]}{0,l(\sigma)}{[}$, $\bs$ de~$S$ et~$\bt$ de~$T$ tels que la famille
$$\left(|a_{\bk}|_{\sigma}^\eps\, \bmax(\bs^\bk,\bt^\bk)\right)_{\bk\in\Z^n}$$
est sommable. 
\end{enumerate}
Une telle famille v\'erifie en particulier la condition suivante : pour tout \'el\'ement~$i$ de~$J$ et tout \'el\'ement~$\bk$ de~$\Z^n$ v\'erifiant~$k_{i}<0$, nous avons~$a_{\bk}=0$. Pour~$i$ dans~$\cn{1}{n}$, posons $\eps_{i}=1$, si~$r_{i}>1$, et $\eps_{i}=-1$, si $r_{i}<1$. La famille pr\'ec\'edente v\'erifie \'egalement la condition suivante : l'ensemble
$$\{\bk\in\Z^n\, |\, \forall i\in\cn{1}{n},\, \eps_{i} k_{i} >0 \textrm{ et } a_{\bk}\}$$
est fini.

Le morphisme naturel $A[\bT] \to \Os_{X,x}$ induit un isomorphisme
$$L_{c} \xrightarrow[]{\sim} \Os_{X,x}.$$
\end{prop}
\begin{proof}
Soit~$\Sigma_{0}$ une partie finie de~$\Sigma$ qui contient~$\Sigma_{\infty}$. Soit $(\eps_{\sigma})_{\sigma\in\Sigma_{0}}$ un \'el\'ement de $\prod_{\sigma\in\Sigma_{0}} \of{]}{0,l(\sigma)}{[}$. Posons
$$M = B\setminus \bigcup_{\sigma\in\Sigma_{0}} \of{]}{a_{\sigma}^{\eps_{\sigma}},a_{\sigma}^{l(\sigma)}}{]}.$$  
C'est un voisinage compact du point~$a_{0}$ dans~$B$ et l'ensemble des parties construites de cette mani\`ere est un syst\`eme fondamental du point~$a_{0}$ dans~$B$.

Nous avons
$$\Bs(M) = A\left[\frac{1}{\Sigma_{0}}\right]$$
et, pour tout \'el\'ement~$f$ de~$\Bs(M)$,
$$\|f\|_{M} = \max_{\sigma\in\Sigma_{0}} (|f|_{\sigma}^{\eps_{\sigma}}).$$
Nous d\'eduisons alors le r\'esultat attendu du th\'eor\`eme~\ref{anneaulocal}. 

\`A l'aide de la formule du produit, l'on d\'emontre que tout \'el\'ement non nul~$a$ de~$\Bs(M)$ satisfait l'in\'egalit\'e~$\|a\|_{M}\ge 1$. Le r\'esultat concernant la forme des s\'eries en d\'ecoule aussit\^ot. 
\end{proof}

Dans le cas de la droite, nous pouvons simplifier la description. Pour traiter ce cas, nous supposerons que~$n=1$ et supprimerons les indices des notations. Nous adopterons les notations suivantes. Si $f=\sum_{k\in\N} a_{k}\, T^k$ est une s\'erie \`a coefficients dans~$K$ et~$\sigma$ un \'el\'ement de~$\Sigma$, nous noterons~$R_{\sigma}(f)$ le rayon de convergence de la s\'erie $\sum_{k\in\N} |a_{k}|_{\sigma}\, T^k$. 

\begin{cor}\label{descriptioncentral13}
Supposons que~$n=1$ et que le point~$b$ est le point central~$a_{0}$ de l'espace~$B$. Notons~$E$ l'anneau compos\'e des s\'eries \`a coefficients dans~$K$ de la forme
$$f = \sum_{k\in\N} a_{k}\, (T-\alpha)^k$$
qui v\'erifient les conditions suivantes : 
\begin{enumerate}[\it i)]
\item il existe un \'el\'ement~$N$ de~$A^*$ tel que, quel que soit~$k$ dans~$\N$, l'\'el\'ement~$a_{k}$ appartient \`a~$A[1/N]$ ;
\item quel que soit~$\sigma$ dans~$\Sigma$, nous avons $R_{\sigma}(f)>0$.
\end{enumerate}
C'est un anneau de valuation discr\`ete d'id\'eal maximal~$(T)$ et de corps r\'esiduel~$K$. Si~$r=0$, le morphisme naturel $A[T] \to \Os_{X,x}$ induit un isomorphisme
$$E \xrightarrow[]{\sim} \Os_{X,x}.$$
Si~$r\in\of{]}{0,1}{[}$, le morphisme naturel $A[T] \to \Os_{X,x}$ induit un isomorphisme
$$\textrm{Frac}(E) = E\left[\frac{1}{T-\alpha}\right] \xrightarrow[]{\sim} \Os_{X,x}.$$
L'anneau local~$\Os_{X,x}$ est alors un corps hens\'elien.
\end{cor}
\begin{proof}
Supposons, tout d'abord, que~$r=0$. Reprenons les notations de la proposition pr\'ec\'edente. Soit~$f= \sum_{k\in\N} a_{k}\, (T-\alpha)^k$ un \'el\'ement de~$L_{c}$. Il existe un sous-ensemble fini~$\Sigma_{0}$ de~$\Sigma$ contenant~$\Sigma_{\infty}$ tel que, quel que soit~$k$ dans~$\N$, l'\'el\'ement~$a_{k}$ appartient \`a~$A[1/\Sigma_{0}]$. En utilisant la finitude du groupe des classes de l'anneau~$A$, on montre qu'il existe un \'el\'ement~$N$ de~$A^*$ tel que
$$A\left[\frac{1}{\Sigma_{0}}\right] = A\left[\frac{1}{N}\right].$$
Soit~$\sigma$ dans~$\Sigma_{0}$. Il existe des \'el\'ements~$\eps$ de~$\of{]}{0,l(\sigma)}{[}$ et~$t$ de~$T$ tels que la s\'erie
$$\sum_{k\in\N} |a_{k}|_{\sigma}^\eps\, t^k$$
converge. On en d\'eduit que~$R_{\sigma}(f)\ge t^{1/\eps} >0$. 

Soit~$\sigma\in\Sigma\setminus\Sigma_{0}$. Pour tout \'el\'ement~$k$ de~$\N$, nous avons $|a_{k}|_{\sigma}\le 1$. On en d\'eduit que $R_{\sigma}(f)\ge 1 >0$. Par cons\'equent, l'\'el\'ement~$f$ appartient \`a~$E$.

R\'eciproquement, soit~$f= \sum_{k\in\N} a_{k}\, (T-\alpha)^k$ un \'el\'ement de~$E$. Il existe un \'el\'ement~$N$ de~$A^*$ tel que la s\'erie~$f$ appartienne \`a~$A[1/N][\![T]\!]$. Posons
$$\Sigma_{0} = \{\m\in\Sigma_{f}\, |\, N\in \m\} \cup \Sigma_{\infty}.$$
C'est une partie finie de~$\Sigma$ qui v\'erifie
$$A\left[\frac{1}{\Sigma_{0}}\right] = A\left[\frac{1}{N}\right].$$
Choisissons un \'el\'ement~$t>0$ qui satisfait la condition suivante :
$$\forall \sigma\in\Sigma_{0},\, t < R_{\sigma}(f).$$ 
Alors, pour tout \'el\'ement~$\sigma$ de~$\Sigma_{0}$, la famille
$$(|a_{k}|_{\sigma}\, t^k)_{k\in\N}$$
est sommable. On en d\'eduit que l'\'el\'ement~$f$ appartient \`a~$L_{c}$.\\

Le cas o\`u le nombre r\'eel~$r$ appartient \`a l'intervalle~$\of{]}{0,1}{[}$ se traite de la m\^eme mani\`ere. Remarquons que la proposition pr\'ec\'edente assure d\'ej\`a que n'interviennent dans le d\'eveloppement en s\'erie d'un \'el\'ement de~$L_{c}$ qu'un nombre fini de termes non nuls d'indice n\'egatif. Montrons, \`a pr\'esent, que le corps Frac($E$) est hens\'elien. D'apr\`es~\cite{rouge}, lemme~2.3.2\,\footnote{V. Berkovich \'enonce, en fait, ce r\'esultat pour des corps suppos\'es \og quasi-complete \fg. La d\'efinition 2.3.1 montre que cette notion co\"{\i}ncide avec celle de corps hens\'elien.}, il suffit de montrer que l'anneau local~$E$ est hens\'elien. La proposition~\ref{Hensel} assure que tel est bien le cas.
\end{proof}

\begin{rem}
Soient~$N$ un \'el\'ement de~$A^*$ et $f = \sum_{k\in\N} a_{k}\, (T-\alpha)^k$ une s\'erie \`a coefficients dans~$A[1/N]$. Posons
$$\Sigma_{0} = \{\m\in\Sigma_{f}\, |\, N\in \m\} \cup \Sigma_{\infty}.$$
C'est une partie finie de~$\Sigma$. Pour tout \'el\'ement~$\m$ de~$\Sigma\setminus\Sigma_{0}$, la s\'erie~$f$ est \`a coefficients dans~$\hat{A}_{\m}$ et nous avons donc $R_{\m}(f)\ge 1$. Par cons\'equent, pour assurer que la s\'erie ~$f$ appartient \`a l'anneau~$E$, il suffit de tester un nombre fini de conditions.
\end{rem}

\index{Anneau local en un point!deploye@d\'eploy\'e!de An@de \E{n}{A}|)}

Donnons, \`a pr\'esent, un exemple d'application de ces descriptions explicites. Nous nous pla\c{c}erons de nouveau dans le cadre de la droite et consid\'ererons le point~$x$ d\'efini comme \'etant le point~$0$ de la fibre centrale. Reprenons les notations du corollaire pr\'ec\'edent. Nous identifierons l'anneau local~$\Os_{X,x}$ avec l'anneau~$E$. Nous noterons~$F$ son corps des fractions. Nous avons d\'emontr\'e que c'est un corps hens\'elien. Observons que cette propri\'et\'e permet de retrouver le th\'eor\`eme d'Eisenstein.

\begin{thm}[Eisenstein]\label{Eisenstein}\index{Eisenstein|see{Th\'eor\`eme d'Eisenstein}}\index{Theoreme@Th\'eor\`eme!d'Eisenstein}
Tout \'el\'ement de~$K[\![T]\!]$ entier sur~$K[T]$ appartient \`a~$E$.
\end{thm}
\begin{proof}
Soit~$f$ un \'el\'ement de~$K[\![T]\!]$ entier sur~$K[T]$. Il est encore entier sur l'anneau local~$E$. Par cons\'equent, il existe un polyn\^ome~$P\in E[U]$ unitaire qui annule~$f$. Puisque l'anneau local~$E$ est factoriel, l'anneau~$E[U]$ l'est \'egalement. Il existe donc un entier~$r$, des polyn\^omes~$P_{1},\ldots,P_{r}$ \`a coefficients dans~$E$, irr\'eductibles et unitaires et des entiers $n_{1},\ldots,n_{r}$ tels que l'on ait l'\'egalit\'e
$$P = \prod_{i=1}^r P_{i}^{n_{i}} \textrm{ dans } E[U].$$

Soit~$i\in\cn{1}{r}$. Puisque la caract\'eristique du corps~$F$ est nulle, le polyn\^ome~$P_{i}$ est s\'eparable. D'apr\`es~\cite{bleu}, proposition~2.4.1, la ca\-t\'e\-go\-rie des extensions s\'eparables finies du corps~$F$ est \'equivalente \`a celle des extensions s\'eparables finies de son compl\'et\'e~$\hat{F}$. On en d\'eduit que le polyn\^ome~$P_{i}$ est encore irr\'eductible dans~$\hat{F}[U]$. 

Remarquons, \`a pr\'esent, que le corps~$\hat{F}$ n'est autre que le corps des s\'eries de Laurent~$K(\!(T)\!)$. L'\'ecriture
$$P = \prod_{i=1}^r P_{i}^{n_{i}}$$
est encore la d\'ecomposition du polyn\^ome~$P$ en produits de facteurs irr\'eductibles et unitaires dans~$K[\![T]\!][U]$. Par cons\'equent, il existe~$i\in\cn{1}{r}$ tel que $P_{i} = U-f$. On en d\'eduit que la s\'erie~$f$ est un \'el\'ement de~$E$.
\end{proof}

\subsection{Anneaux de sections globales}

Dans cette partie, nous voulons d\'ecrire les anneaux de sections globales de certaines parties de l'espace affine~$X$. Plus pr\'ecis\'ement, nous allons nous int\'eresser aux disques et couronnes compacts au-dessus de parties compactes et connexes de l'espace~$B$. 

Introduisons quelques notations. Pour une partie~$V$ de~$B$ et des $n$-uplets $\bs=(s_{1},\ldots,s_{n})$ et $\bt=(t_{1},\ldots,t_{n})$ dans $\R_{+}^n$, nous posons
$$\mathring{D}_{V}(\bt) = \{x\in X\, |\, \pi(x)\in V,\, \forall i\in\cn{1}{n},\, |T_{i}(x)| < t_{i} \},$$
$$\overline{D}_{V}(\bt) = \{x\in X\, |\, \pi(x)\in V,\, \forall i\in\cn{1}{n},\, |T_{i}(x)| \le t_{i} \},$$
$$\mathring{C}_{V}(\bs,\bt) = \{x\in X\, |\, \pi(x)\in V,\, \forall i\in\cn{1}{n},\, s_{i} < |T_{i}(x)| < t_{i} \}$$
et
$$\overline{C}_{V}(\bs,\bt) = \{x\in X\, |\, \pi(x)\in V,\, \forall i\in\cn{1}{n},\, s_{i} \le |T_{i}(x)| \le t_{i} \}.$$
\newcounter{crel}\setcounter{crel}{\thepage}

Toutes ces parties sont compactes, en vertu de la proposition \ref{disquecompact}.

\begin{defi}\index{Disque!relatif}\index{Couronne!relative}
Nous appellerons {\bf disque relatif} de~$X$ toute partie de la forme $\mathring{D}_{V}(\bt)$ ou $\overline{D}_{V}(\bt)$, o\`u~$V$ d\'esigne une partie de~$B$ et~$\bt$ un \'el\'ement de~$\R_{+}^n$.

Nous appellerons {\bf couronne relative} de~$X$ toute partie de la forme $\mathring{C}_{V}(\bs,\bt)$ ou $\overline{C}_{V}(\bs,\bt)$, o\`u~$V$ d\'esigne une partie de~$B$ et~$\bs$ et~$\bt$ deux \'el\'ements de~$\R_{+}^n$.
\end{defi}

Rappelons que, d'apr\`es la d\'efinition \ref{defiFV} et les remarques qui la suivent, si $K$ est une partie compacte de~$X$, la notation~$\Os(K)$ d\'esigne l'anneau des fonctions qui sont d\'efinies au voisinage de~$K$. En particulier, si~$(k,|.|)$ est un corps ultram\'etrique complet et~$\overline{D}$ le disque unit\'e de~$\E{n}{k}$, l'alg\`ebre~$\Os(\overline{D})$ n'est pas l'alg\`ebre affino\"ide~$k\{\bT\}$, mais l'alg\`ebre des s\'eries surconvergentes, constitu\'ee de l'ensemble des s\'eries de~$k[\![\bT]\!]$ dont le rayon de convergence est strictement sup\'erieur \`a~$1$.

Commen\c{c}ons par \'enoncer un r\'esultat topologique. C'est un cas particulier du lemme \ref{lemvois}.

\begin{lem}\label{lemvoisdisquedroite}\index{Voisinages d'un compact!disque}\index{Voisinages d'un compact!couronne}\index{Disque!voisinages}\index{Couronne!voisinages}
Soient~$V$ une partie compacte de~$B$ et~$\bt$ un \'el\'ement de~$\R_{+}^n$. Tout voisinage du disque~$\overline{D}_{V}(\bt)$ contient un disque de la forme~$\overline{D}_{V}(\bt')$, o\`u~$\bt'$ est un \'el\'ement de~$\R_{+}^n$ qui v\'erifie l'in\'egalit\'e~$\bt' > \bt$.

Soit~$\bs$ un \'el\'ement de~$\R_{+}^n$ tel que~$\bs\le\bt$. Tout voisinage de la couronne~$\overline{C}_{V}(\bs,\bt)$ contient une couronne de la forme~$\overline{C}_{V}(\bs',\bt')$, o\`u~$\bs'$ et~$\bt'$ sont deux \'el\'ements de~$\R_{+}^n$ qui v\'erifient les in\'egalit\'es~$\bs'\prec \bs$ et~$\bt' > \bt$.
\end{lem}



Consacrons-nous, \`a pr\'esent, \`a l'\'etude des fonctions d\'efinies au voisinage de disques compacts. Nous commen\c{c}ons par montrer que ces fonctions admettent un d\'eveloppement en s\'erie.

\begin{prop}\label{dvptserie}
Soit $V$ une partie compacte de $B$. Soit $\bt \in \R_{+}^n$. Alors le morphisme naturel 
$$\Os(V)[\bT] \to \Os(V)[\![\bT]\!]$$ 
se prolonge en un morphisme injectif
$$\varphi_{V,\bt} : \Os\left(\overline{D}_{V}(\bt)\right) \hookrightarrow \Os(V)[\![\bT]\!].$$
\end{prop}
\begin{proof}
Soit $f \in \Os\left(\overline{D}_{V}(\bt)\right)$. D'apr\`es le lemme~\ref{lemvoisdisquedroite}, il existe un polyrayon $\br >\bt$ telle que la fonction $f$ soit d\'efinie sur $\mathring{D}_{V}(\br)$.

Soit $b$ un point de~$V$. La fonction~$f$ est d\'efinie au voisinage du point~$0$ de la fibre~$X_{b}$. D'apr\`es le th\'eor\`eme~\ref{anneaulocal}, il existe un voisinage compact~$V^b$ du point~$b$ dans~$B$ et un nombre r\'eel~$r_b>0$ tels qu'au voisinage de la partie compacte~$\overline{D}_{V^{b}}(r_b)$ de~$X$, la fonction~$f$ poss\`ede une expression de la forme 
$$f = \sum_{\bk\in\N^n} a_\bk \bT^\bk,$$
o\`u, quel que soit $\bk\in\N^n$, $a_\bk \in \Bs(V^b)$.

En identifiant localement les diff\'erents d\'eveloppements en s\'erie, on montre que, quel que soit $\bk\in\N^n$, l'\'el\'ement $a_{\bk}$ appartient \`a $\Os(V)$. Nous avons donc construit un morphisme 
$$\varphi_{V,\bt} : \Os\left(\overline{D}_{V}(\bt)\right) \to \Os(V)[\![\bT]\!]$$ 
qui co\"{\i}ncide avec le morphisme naturel $\Os(V)[\bT] \to \Os(V)[\![\bT]\!]$ sur~$\Os(V)[\bT]$.

Montrons que le morphisme $\varphi_{V,\bt}$ est injectif. Supposons que deux fonctions~$f$ et~$g$ de~$\Os\left(\overline{D}_{V}(\bt)\right)$ aient la m\^eme image. Soit $b\in V$. Notons $x$ le point $0$ de la fibre $X_{b}$. Les fonctions $f$ et $g$ ont m\^eme d\'eveloppement dans $L_{b} \simeq \Os_{X,x}$. On en d\'eduit que les fonctions $f$ et $g$ co\"{\i}ncident sur un voisinage de $x$ dans la fibre $X_{b}$. Puisque cette fibre est un espace irr\'eductible, les fonctions $f$ et $g$ co\"{\i}ncident n\'ecessairement sur toute la fibre. On en d\'eduit finalement que $f=g$.
\end{proof}

Afin de d\'ecrire explicitement l'image du morphisme pr\'ec\'edent, introduisons une notation. Pour toute partie compacte~$V$ de~$B$ et tout \'el\'ement~$\bt$ de~$\R_{+}^n$, nous noterons
$$\Os(V)\of{\la}{|\bT| \le \bt}{\ra}^\dag$$
l'anneau des s\'eries \`a coefficients dans~$\Os(V)$ de la forme
$$\sum_{\bk \in\N^n} a_{\bk} \bT^\bk$$
qui v\'erifient la condition suivante :
$$\exists \br > \bt,\, \lim_{\bk\to +\infty} \|a_\bk\|_{V}\, \br^\bk = 0.$$
\newcounter{surcvd}\setcounter{surcvd}{\thepage}

\begin{prop}\label{imagedisque}
Soit $V$ une partie compacte de $B$. Soit $\bt \in\R_{+}^n$. L'image du morphisme $\varphi_{V,\bt}$ est contenue dans $\Os(V)\of{\la}{|\bT| \le \bt}{\ra}^\dag$.
\end{prop}
\begin{proof}
Soit $f \in \Os(\overline{D}_{V}(\bt))$. D'apr\`es le lemme~\ref{lemvoisdisquedroite}, il existe un polyrayon $\bv >\bt$ telle que la fonction $f$ soit d\'efinie sur~$\mathring{D}_{V}(\bv)$. La proposition pr\'ec\'edente nous montre que la fonction~$f$ poss\`ede un d\'eveloppement en s\'erie de la forme 
$$f = \sum_{\bk \in\N^n} a_{\bk} \bT^\bk \in \Os(V)[\![\bT]\!].$$

Soit $b\in V$. Puisque le groupe $|\Hs(b)^*|$ est discret dans $\R_{+}^*$, il existe une famille \mbox{$\bu=(u_{1},\ldots,u_{n})$} de $\R_{+}^*$ qui v\'erifie $\bt < \bu < \bv$ et dont l'image est libre dans le \mbox{$\Q$-espace} vectoriel \mbox{$\Q\times_{\Z} \left(\R_{+}^*/|\Hs(b)^*|\right)$}. Notons $x$ l'unique point de la fibre~$X_{b}$ qui v\'erifie
$$\forall i\in\cn{1}{n},\, |T_{i}(x)|=u_{i}.$$
La description de l'anneau local au point $x$ obtenue au th\'eor\`eme \ref{anneaulocal} nous assure qu'il existe un voisinage $V^{b}$ de $b$ dans $B$ et $\br_{b} > \bv > \bt$ tels que
$$\lim_{\bk \to +\infty} \|a_{\bk}\|_{V^{b}}\, \br_{b}^\bk = 0.$$

Par compacit\'e, nous pouvons recouvrir la partie $V$ par un nombre fini de compacts $V^{b_{1}},\ldots,V^{b_{p}}$, avec $p\in\N$ et $b_{1},\ldots,b_{p} \in V$. On en d\'eduit qu'il existe~\mbox{$\br > \bt$}  tel que  
$$\lim_{\bk\to +\infty} \|a_\bk\|_{V}\, \br^\bk = 0.$$
\end{proof}

\begin{rem}
Ce r\'esultat cache un principe du prolongement analytique. Nous n'insisterons pas ici sur ce point, mais consacrerons la section~\ref{sectionprolan} \`a ce propos.
\end{rem}

Int\'eressons-nous, \`a pr\'esent, \`a la r\'eciproque de ce r\'esultat. Nous n'allons consid\'erer que certaines parties compactes de la base. 

\begin{thm}\label{isodisque}\index{Disque!sections globales}\index{Faisceau structural!sections sur un disque de An@sections sur un disque de $\E{n}{A}$}
Soit $V$ une partie compacte et connexe de $B$. Supposons que le point central de $B$ n'appartienne pas au bord du compact $V$. Soit \mbox{$\bt \in \R_{+}^n$}. Alors le morphisme 
$$\varphi_{V,\bt} : \Os\left(\overline{D}_{V}(\bt)\right) \hookrightarrow \Os(V)[\![\bT]\!]$$
r\'ealise un isomorphisme sur l'anneau~$\Os(V)\of{\la}{|\bT| \le \bt}{\ra}^\dag$.
\end{thm}
\begin{proof}
D'apr\`es les propositions qui pr\'ec\`edent, il nous suffit de montrer que toute s\'erie de la forme donn\'ee appartient \`a l'image de $\varphi_{V,\bt}$. Nous allons distinguer plusieurs cas, en fonction du compact $V$.

Commen\c{c}ons par consid\'erer un compact de la forme 
$$V = \of{[}{a_{\sigma,\alpha},a_{\sigma,l(\sigma)}}{]},$$
avec $\sigma\in\Sigma$ et $\alpha\in\of{]}{0,l(\sigma)}{[}$. 

Soit $\br'\in\R_{+}^n$ tel que $\bt < \br' < \br$. Soit $\mu>1$ tel que $\bt < (\br')^\mu < \br$. Soit~\mbox{$\bk\in\N^n$}. Nous avons 
$$\lim_{\bk\to +\infty} |a_\bk|_{\sigma}^\alpha\, (\br')^\bk = 0$$
et l'on en d\'eduit que 
$$\lim_{\bk\to +\infty} |a_\bk|_{\sigma}^{\alpha \mu}\, \left((\br')^\mu\right)^\bk = 0.$$

Remarquons, \`a pr\'esent, que, quel que soit $\bk\in\N^n$, l'\'el\'ement $a_{\bk}$ de $\Os(V)=\hat{A}_{\sigma}$ se prolonge \`a l'ouvert $U = \of{]}{a_{\sigma,\alpha\mu},a_{\sigma,l(\sigma)}}{]}$ et v\'erifie
$$\|a_{\bk}\|_{U} = |a_{\bk}|_{\sigma}^{\alpha\mu}.$$
On en d\'eduit que la s\'erie $f$ d\'efinit un \'el\'ement de $\Os\big(\mathring{D}_{U}(\br')\big)$ et donc de~$\Os\left(\overline{D}_{V}(\bt)\right)$.\\

Ce raisonnement met en \'evidence le fait que la difficult\'e du probl\`eme r\'eside dans l'\'etude du comportement au bord du compact~$V$. Remarquons que ce bord ne peut contenir qu'un nombre fini de points. En effet, si le compact~$V$ ne contient pas le point central de~$B$, sa connexit\'e lui impose d'\^etre contenu dans une branche de~$B$. Il est donc de la forme 
$$V = \of{[}{a_{\sigma}^u,a_{\sigma}^v}{]},$$ 
avec $\sigma\in\Sigma$, $u,v\in\of{]}{0,l(\sigma)}{]}$ et $u\le v$. Son bord contient alors au plus deux points. Si le compact~$V$ contient le point central~$a_{0}$ de~$B$, alors, par hypoth\`ese, il contient un voisinage de ce point et il n'existe donc qu'un nombre fini de branches de~$B$ que~$V$ ne contient pas enti\`erement. On en d\'eduit que le bord du compact~$V$ n'est constitu\'e que d'un nombre fini de points. En reprenant le raisonnement pr\'ec\'edent en chaque point du bord du compact~$V$, on obtient le r\'esultat annonc\'e. 
\end{proof}

\begin{rem}
\'Enonc\'ee de la m\^eme fa\c{c}on, la proposition pr\'ec\'edente est fausse si le point central de~$B$ se situe sur le bord du compact~$V$. Consid\'erons, par exemple, la partie compacte constitu\'ee du seul point central de~$\Ms(\Z)$,
$$V = \{a_{0}\}.$$
L'anneau $\Os(V)$ est alors l'anneau~$\Q$ et la norme $\|.\|_{V}$ est la norme triviale. 

Pla\c{c}ons-nous sur la droite $\E{1}{\Z}$. Soit $t\in\of{[}{0,1}{[}$. L'anneau $\Os(V)\of{\la}{|T| \le t}{\ra}^\dag$ n'est autre que l'anneau 
$\Q[\![T]\!]$.
Consid\'erons la s\'erie 
$$f = \sum_{k\in\N} k!\, T^k.$$
Elle appartient bien \`a l'anneau pr\'ec\'edent, mais ne peut se prolonger \`a aucun disque de centre $0$ et de rayon strictement positif de la branche archim\'edienne de $\Ms(\Z)$.

De m\^eme, pour tout nombre premier $p$, la s\'erie 
$$f = \sum_{k\in\N} q^{-k^2}\, T^k \in \Q[\![T]\!]$$
ne peut se prolonger \`a aucun disque de centre $0$ et de rayon strictement positif de la branche $p$-adique de $\Ms(\Z)$.



\end{rem}

Le cas des couronnes se traite de fa\c{c}on analogue \`a celui des disques. Introduisons de nouveau une notation. Soient~$V$ une partie compacte de~$B$ et~$\bs$ et~$\bt$ deux \'el\'ements de~$\R_{+}^n$. Posons $I=\{i\in\cn{1}{n}\, |\, s_{i}> 0\}$ et $J=\cn{1}{n}\setminus I$. Nous noterons
$$\Os(V)\of{\la}{\bs \le |\bT| \le \bt}{\ra}^\dag$$
l'anneau constitu\'e des s\'eries \`a coefficients dans~$\Os(V)$ de la forme
$$\sum_{\bk \in\Z^n} a_{\bk} \bT^\bk,$$
qui v\'erifient les trois conditions suivantes :
$$\forall \bk\in\Z^n\setminus \left(\prod_{i\in I} \R \times \prod_{i\in J} \R_{+}\right),\, a_{\bk}=0,$$
$$\exists \br > \bt,\, \lim_{\bk\to +\infty} \|a_\bk\|_{V}\, \br^\bk = 0$$
et
$$\exists \br \prec \bs,\, \lim_{\bk\to -\infty} \|a_\bk\|_{V}\, \br^\bk = 0.$$
\newcounter{surcvc}\setcounter{surcvc}{\thepage}

En particulier, si~$\bs=\boldsymbol{0}$, alors cet anneau est contenu dans $\Os(V)[\![\bT]\!]$ et nous avons l'\'egalit\'e
$$\Os(V)\of{\la}{\boldsymbol{0} \le |\bT| \le \bt}{\ra}^\dag = \Os(V)\of{\la}{|\bT| \le \bt}{\ra}^\dag.$$

\begin{prop}\label{imagecouronne}
Soit $V$ une partie compacte de $B$. Soient~$\bs$ et~$\bt$ deux \'el\'ements de~$\R_{+}^n$ v\'erifiant l'in\'egalit\'e $\bs \prec \bt$. Alors le morphisme naturel 
$$\Os(V)[\bT] \to \Os(V)[\![\bT]\!]$$ 
se prolonge en un morphisme injectif 
$$\varphi_{V,\bs,\bt} : \Os\left(\overline{C}_{V}(\bs,\bt)\right) \hookrightarrow \Os(V)\of{\la}{\bs \le \bT \le \bt}{\ra}^\dag.$$
\end{prop}
\begin{proof}
Il suffit de reprendre la preuve des propositions~\ref{dvptserie} et~\ref{imagedisque}. Il faut cependant prendre garde au fait que nous ne pouvons plus consid\'erer un voisinage du point~$0$ d'une fibre. Il est cependant possible de remplacer ce point par un point de type~$3$ d\'eploy\'e, c'est-\`a-dire un point~$x$ d\'efini par des \'equations du type
$$\forall i\in\cn{1}{n},\, |T_{i}(x)|=r_{i},$$
o\`u $r_{1},\ldots,r_{n}$ sont des \'el\'ements de $\R_{+}^*$ tels que l'image de la famille $(r_{1},\ldots,r_{n})$ dans le $\Q$-espace vectoriel $\Q\otimes_{\Z} \left(\R_{+}^*/|\Hs(b)^*|\right)$ est libre. Un tel choix est possible car le groupe $|\Hs(b)^*|$ est discret dans $\R_{+}^*$. Dans ce cas, nous disposons encore d'une description de l'anneau local en termes de s\'eries, par le th\'eor\`eme \ref{anneaulocal}.
\end{proof}

Comme dans le cas des disques, nous pouvons raffiner cette proposition pour obtenir, dans certains cas, un r\'esultat d'isomorphie similaire \`a celui de la proposition~\ref{isodisque}. La d\'emonstration en \'etant compl\`etement analogue, nous ne la r\'edigerons pas.

\begin{thm}\label{isocouronne}\index{Couronne!sections globales}\index{Faisceau structural!sections sur une couronne de An@sections sur une couronne de $\E{n}{A}$}
Soit $V$ une partie compacte et connexe de $B$. Supposons que le point central de $B$ n'appartienne pas au bord du compact $V$. Soient~$\bs$ et~$\bt$ deux \'el\'ements de~$\R_{+}^n$ v\'erifiant l'in\'egalit\'e $\bs \prec \bt$. Alors, le morphisme 
$$\varphi_{V,\bs,\bt} : \Os\left(\overline{C}_{V}(\bs,\bt)\right) \xrightarrow[]{\sim} \Os(V)\of{\la}{\bs \le |\bT| \le \bt}{\ra}^\dag$$
est un isomorphisme.
\end{thm}

\bigskip

Int\'eressons-nous, \`a pr\'esent, au bord analytique des couronnes. Dans le cas d'espaces d\'efinis au-dessus d'un corps ultram\'etrique, nous disposons d'une description explicite.


\begin{lem}
Soit~$(k,|.|)$ un corps ultram\'etrique complet. Soient~$\bs$ et~$\bt$ deux \'el\'ement de~$\R_{+}^n$ v\'erifiant l'in\'egalit\'e $\bs\prec\bt$. Consid\'erons la couronne~$C$ de~$\E{n}{k}$ de rayon int\'erieur~$\bs$ et de rayon ext\'erieur~$\bt$. Pour tout \'el\'ement~$i$ de~$\cn{1}{n}$, notons
$$R_{i} = \{s_{i},t_{i}\}\cap\R_{+}^*.$$
La couronne~$C$ poss\`ede un bord de Shilov. C'est l'ensemble fini et simple
$$\Gamma_{C} = \{\eta_{r_{1},\ldots,r_{n}}\, |\, \forall i\in\cn{1}{n},\, r_{i}\in R_{i}\}.$$
\end{lem}
\begin{proof}
La description explicite des fonctions d\'efinies au voisinage de la couronne~$C$ et de la norme uniforme sur~$C$ montre que, pour tout \'el\'ement~$f$ de~$\Os(C)$, nous avons
$$\|f\|_{C} = \max_{z\in\Gamma_{C}} (|f(x)|).$$
Puisque~$\Os(C)$ est dense dans~$\Bs(C)$ pour la norme~$\|.\|_{C}$, ce r\'esultat vaut encore pour les \'el\'ements de~$\Bs(C)$. On en d\'eduit que la partie~$\Gamma_{C}$ est un bord analytique du compact~$C$. 

En outre, pour tout point~$z$ de~$\Gamma_{C}$, il existe un \'el\'ement~$\bk\in\Z^n$ tel que la fonction $\bT^\bk$ appartienne \`a~$\Ks(C)$ et atteigne son maximum en valeur absolue au point~$z$ et uniquement en ce point. Par cons\'equent, tout bord analytique du compact~$C$ contient la partie~$\Gamma_{C}$. Cette derni\`ere est donc bien le bord de Shilov du compact~$C$.
\end{proof}

Dans le cas archim\'edien, nous disposons \'egalement de r\'esultats.

\begin{lem}
Soit~$(k,|.|)$ un corps archim\'edien complet. Soient~$\bs$ et~$\bt$ deux \'el\'ements de~$\R_{+}^n$ v\'erifiant l'in\'egalit\'e $\bs\prec\bt$. Consid\'erons la couronne~$C$ de~$\E{n}{k}$ de rayon int\'erieur~$\bs$ et de rayon ext\'erieur~$\bt$. Pour tout \'el\'ement~$i$ de~$\cn{1}{n}$, notons
$$R_{i} = \{s_{i},t_{i}\}\cap\R_{+}^*.$$
La couronne~$C$ poss\`ede un bord de Shilov. Il est contenu dans l'ensemble compact
$$\Gamma_{C} = \{x\in\E{n}{\C}\, |\, \forall i\in\cn{1}{n},\, \exists r_{i}\in R_{i},\, |T_{i}(x)|=r_{i}\}.$$
En outre, si le corps~$k$ est le corps des nombres complexes~$\C$, l'\'egalit\'e vaut.
\end{lem}
\begin{proof}
L'existence du bord de Shilov d\'ecoule du r\'esultat d'A.~Escassut et N.~Ma\"inetti d\'ej\`a cit\'e (\emph{cf.} th\'eor\`eme~\ref{escassutmainetti}).

Le principe du maximum assure que le bord de Shilov de la couronne~$C$ est contenu dans son bord topologique, qui n'est autre que la partie~$\Gamma_{C}$. Cette remarque permet de d\'emontrer le premier point.

Supposons, \`a pr\'esent, que le corps~$k$ est le corps~$\C$. Nous avons alors
$$\Gamma_{C} = \{(z_{1},\ldots,z_{n}) \in \C^n\, |\, \forall i\in\cn{1}{n},\, \exists r_{i}\in R_{i},\, |z_{i}|=r_{i}\}.$$
Pour tout point~$z=(z_{1},\ldots,z_{n})$ de~$\Gamma_{C}$, il existe un \'el\'ement~$(\alpha_{1},\ldots,\alpha_{n})$ de~$\C^n$ et un \'el\'ement$(k_{1},\ldots,k_{n})$ de~$\{-1,1\}^n$ tels que la fonction $\prod_{1\le i\le n} (z_{i}-\alpha_{i})^{k_{i}}$ soit d\'efinie au voisinage de la couronne~$C$ et atteigne son maximum en valeur absolue au point~$z$ et uniquement en ce point. Par cons\'equent, tout bord analytique du compact~$C$ contient la partie~$\Gamma_{C}$. Cette derni\`ere est donc bien le bord de Shilov du compact~$C$.
\end{proof}

Ces rappels nous permettent de d\'ecrire un bord analytique non trivial des couronnes relatives \`a l'aide du lemme suivant. Remarquons que toute couronne compacte au-dessus d'une partie compacte et connexe de~$B$ (et donc pro-rationnelle, d'apr\`es la proposition~\ref{prorat}) est pro-rationnelle et donc spectralement convexe, d'apr\`es le th\'eor\`eme~\ref{compactrationnel}.

\begin{lem}\label{Shilovcouronne}
Soit $V$ une partie compacte et connexe de $B$ et~$C$ une couronne compacte au-dessus de~$V$. Pour tout point~$v$ de~$V$, notons~$\gamma_{v}$ le bord de Shilov du compact~$C\cap X_{v}$ dans~$X_{v}$. Alors, la partie
$$\Gamma = \bigcup_{v\in V} \gamma_{v}$$
est un bord analytique de la couronne~$C$.
\end{lem}
\begin{proof}
Puisque~$\Ks(C)$ est dense dans~$\Bs(C)$ pour la norme~$\|.\|_{C}$, il suffit de d\'emontrer que, pour tout \'el\'ement~$f$ de~$\Ks(C)$, nous avons
$$\|f\|_{C} = \|f\|_{\Gamma}.$$
Soit~$f$ un \'el\'ement de~$\Ks(C)$. Il existe un \'el\'ement~$v$ de~$V$ tel que l'on ait
$$\|f\|_{C} = \|f\|_{C\cap X_{v}}.$$
La fonction~$f$ induit par restriction une section sur $C\cap X_{v}$ du pr\'efaisceau des fonctions rationnelles sur~$X_{v}$. Nous avons donc
$$\|f\|_{C\cap X_{v}} = \|f\|_{\gamma_{v}}.$$
On en d\'eduit le r\'esultat attendu.
\end{proof}

La description des fonctions au voisinage des couronnes obtenue plus haut permet de pr\'eciser ce r\'esultat dans le cas ultram\'etrique.

\begin{prop}\label{Shilovcouronneum}\index{Couronne!bord de Shilov}\index{Bord analytique!d'une couronne}
Soit $V$ une partie compacte et connexe de $B_{\textrm{um}}$ et~$C$ une couronne compacte au-dessus de~$V$. Notons~$\Gamma_{V}$ le bord de Shilov du compact~$V$ dans~$B$. Pour tout point~$v$ de~$V$, notons~$\Gamma_{v}$ le bord de Shilov du compact~$C\cap X_{v}$ dans~$X_{v}$. La couronne~$C$ poss\`ede un bord de Shilov. C'est l'ensemble fini
$$\Gamma = \bigcup_{v\in \Gamma_{V}} \gamma_{v}.$$
\end{prop}

\begin{proof}
Dans le cas o\`u la couronne est vide, le r\'esultat est imm\'ediat. Dans le cas contraire, il existe deux \'el\'ements~$\bs$ et~$\bt$ de~$\R_{+}^n$ v\'erifiant $\bs\prec\bt$ tels que $C=\overline{C}_{V}(\bs,\bt)$. D'apr\`es la proposition~\ref{imagecouronne}, le morphisme naturel $\Os(V)[\bT] \to \Os(V)[\![\bT]\!]$ 
se prolonge en un morphisme injectif 
$$\Os(C) \hookrightarrow \Os(V)\of{\la}{\bs \le \bT \le \bt}{\ra}^\dag.$$
Commen\c{c}ons par montrer que, pour tout \'el\'ement~$f=\sum_{\bk\in\Z^n} a_{\bk}\, \bT^{\bk}$ de~$\Os(C)$, nous avons
$$\|f\|_{C} = \max_{\bk\in\Z^n}\left(\|a_{\bk}\|_{V}\, \bmax(\bs^\bk,\bt^\bk) \right).$$
Puisque la couronne~$C$ est compacte, il existe un \'el\'ement~$z$ de~$C$ en lequel nous avons l'\'egalit\'e
$$\|f\|_{C} = |f(z)|.$$
Nous avons alors
$$|f(z)|=\|f\|_{C\cap X_{\pi(z)}} = \max_{\bk\in\Z^n}\left(|a_{\bk}(\pi(z))|\, \bmax(\bs^\bk,\bt^\bk) \right),$$
puisque le point~$\pi(z)$ appartient \`a la partie ultram\'etrique~$B_{\textrm{um}}$ de~$B$. On en d\'eduit l'\'egalit\'e annonc\'ee.

De cette description explicite de la norme, on d\'eduit que tout \'el\'ement de~$\Os(C)$, et donc tout \'el\'ement de~$\Bs(C)$ atteint son maximum sur~$\Gamma$, autrement dit que~$\Gamma$ est un bord analytique de~$C$. En outre, pour tout point~$z$ de~$\Gamma$, il existe un \'el\'ement~$a$ de~$\Ks(V)$ (son existence est assur\'ee par la proposition~\ref{Shilovbase}) et un \'el\'ement~$\bk$ de~$\Z^n$ tels que la fonction $a\, \bT^\bk$ appartienne \`a~$\Ks(C)$ et atteigne son maximum en valeur absolue au point~$z$ et uniquement en ce point. Par cons\'equent, la partie~$\Gamma$ est le bord de Shilov de la couronne compacte~$C$.
\end{proof}

\bigskip

Pour finir, calculons explicitement ces anneaux globaux dans un cas particulier, celui des couronnes au-dessus de voisinages compacts du point central.

\begin{prop}\label{descriptionsectionscouronnes}\index{Disque!sections globales}\index{Couronne!sections globales}
Soit~$\Sigma'$ une partie finie de~$\Sigma$ contenant~$\Sigma_{\infty}$. Pour $\sigma\in\Sigma'$, choisissons un \'el\'ement~$\eps_{\sigma}\in\of{]}{0,1}{]}$. Consid\'erons la partie compacte~$V$ de~$B$ d\'efinie par 
$$V = \left( \bigcup_{\sigma\in\Sigma'} \of{[}{a_{0},a_{\sigma}^{\eps_{\sigma}}}{]} \right) \cup \left( \bigcup_{\sigma\notin \Sigma'} B_{\sigma} \right).$$
Soient~$\bs$ et~$\bt$ deux \'el\'ements de~$\R_{+}^n$. Posons $I=\{i\in\cn{1}{n}\, |\, s_{i}> 0\}$ et $J=\cn{1}{n}\setminus I$. L'anneau $\Os(V)\of{\la}{\bs\le |\bT|\le\bt}{\ra}^\dag$ est constitu\'e des s\'eries \`a coefficients dans~$K$ de la forme 
$$\sum_{\bk\in\Z^n} a_{\bk}\, \bT^\bk$$ 
v\'erifiant les conditions suivantes :
\begin{enumerate}[\it i)]
\item $\forall \bk\ge 0$, $\disp a_{k}\in A\left[\frac{1}{\Sigma'}\right]$ ;
\item $\forall \bk\in\Z^n\setminus \left(\prod_{i\in I} \R \times \prod_{i\in J} \R_{+}\right),\, a_{\bk}=0$ ;
\item $\forall \sigma\in\Sigma'$, $\exists \br \prec \bs^{\eps_{\sigma}}$, $\disp \lim_{\bk\to -\infty} |a_{\bk}|_{\sigma}\, \br^\bk =0$ ;
\item $\forall \sigma\in\Sigma'$, $\exists \br > \bt^{\eps_{\sigma}}$, $\disp \lim_{\bk\to +\infty} |a_{\bk}|_{\sigma}\, \br^\bk =0$.
\end{enumerate}
Si~$\bt\ge 1$, pour toute s\'erie du type pr\'ec\'edent, l'ensemble
$$\{\bk\in\N^n\, |\, a_{\bk}\ne 0\}$$
est fini. Si~$\bs\le 1$, pour toute s\'erie du type pr\'ec\'edent, l'ensemble 
$$\{\bk\in \Z^n\cap \of{]}{-\infty,0}{]}\, |\, a_{\bk}\ne 0\}$$
est fini. En particulier, si~$\bs=\b0$ et~$\bt\ge 1$, alors l'anneau $\Os(V)\of{\la}{|\bT|\le\bt}{\ra}^\dag$ n'est autre que l'anneau de polyn\^omes~$ A[1/(\Sigma'\cap\Sigma_{f})][\bT]$. 
\end{prop}
\begin{proof}
Les r\'esultats d\'emontr\'es aux num\'eros~\ref{partiesouvertes} et~\ref{borddeShilovbase} permettent de d\'emontrer que nous avons
$$\Os(V) = A\left[\frac{1}{\Sigma'\cap\Sigma_{f}}\right] \textrm{ et } \|.\|_{V} = \max_{\sigma\in\Sigma'}(|.|_{\sigma}^{\eps_{\sigma}}).$$
La premi\`ere partie du r\'esultat d\'ecoule alors imm\'ediatement de la d\'efinition de l'anneau $\Os(V)\of{\la}{|\bT|\le\bt}{\ra}^\dag$.

D'apr\`es la formule du produit, pour tout \'el\'ement non nul~$a$ de $A[1/\Sigma']$, nous avons $\prod_{\sigma\in\Sigma'} |a|_{\sigma}\ge 1$ et donc
$$\|a\|_{V} = \max_{\sigma\in\Sigma'}(|a|_{\sigma}^{\eps_{\sigma}})\ge 1.$$
On en d\'eduit la seconde partie du r\'esultat.
\end{proof}

\section{Points rigides des fibres}\label{prdf}

\index{Point!rigide!d'une fibre de An@d'une fibre de \E{n}{A}|(}

Soit $b$ un point de~$B$. La proposition \ref{voisdep} nous permet de d\'ecrire un syst\`eme fondamental de voisinages explicite d'un point $x$ de la fibre $X_{b}$ d\'efini par des \'equations du type
$$(T_{1}-\alpha_{1})(x)=\cdots = (T_{n}-\alpha_{n})(x)=0,$$
avec $\alpha_{1},\ldots,\alpha_{n} \in \Os_{B,b}$. Remarquons que, lorsque l'espace de base est le spectre d'un anneau d'entiers de corps de nombres, tous les points rationnels de la fibre~$X_{b}$ sont de ce type. En effet, d'apr\`es le lemme \ref{translation}, le morphisme naturel $\Os_{B,b} \surj \Hs(b)$ est surjectif.


Dans ce num\'ero, nous montrons qu'il est possible de ramener l'\'etude de certains points de l'espace~$X$, \`a savoir les points rigides des fibres, \`a celle des points rationnels par le biais d'un isomorphisme local (\emph{cf.} proposition \ref{isolocal}).



\subsection{Isomorphismes locaux}

Nous montrons ici que nous nous trouvons bien dans le cadre d'application de la proposition~\ref{isolocal} et en pr\'ecisons les conclusions. Nous distinguerons selon le type de la fibre dans laquelle se situe le point rigide consid\'er\'e. Commen\c{c}ons par le cas le plus simple : celui des fibres extr\^emes.

\begin{prop}\label{isoext}\index{Isomorphisme local!au voisinage d'un point rigide!d'une fibre extreme@d'une fibre extr\^eme}
Soient $\m$ un \'el\'ement de~$\Sigma_{f}$ et $x$ un point rigide de la fibre extr\^eme~$\tilde{X}_{\m}$. Supposons que le point $x$ poss\`ede un syst\`eme fondamental de voisinages connexes. Alors, il existe une extension finie $K'$ de $K$, un point $x'$ de $\E{n}{A'}$, o\`u $A'$ d\'esigne l'anneau des entiers de $K'$, rationnel dans sa fibre, tel que le morphisme naturel
$$\E{n}{A'} \to \E{n}{A}$$
envoie le point $x'$ sur le point $x$ et induise un isomorphisme d'un voisinage de~$x'$ sur un voisinage de $x$. 

\end{prop}
\begin{proof}
L'extension de corps $k_{\m}\to \Hs(x)$ est une extension finie et s\'eparable, puisque le corps $k_{\m}$ est fini. D'apr\`es le th\'eor\`eme de l'\'el\'ement primitif, il existe un \'el\'ement~$\tilde{\alpha}$ de~$\Hs(x)$ tel que $k_{\m}[\tilde{\alpha}]=\Hs(x)$. Notons $\tilde{P}(S) \in k_{\m}[S]$ le polyn\^ome minimal unitaire de $\tilde{\alpha}$ sur $k_{\m} = A/\m$. Choisissons un relev\'e unitaire~$P(S)$ de~$\tilde{P}(S)$ dans~$A[S]$. Ce polyn\^ome est encore irr\'eductible. Consid\'erons l'extension finie~$K' = K[S]/(P(S))$ de~$K$. C'est un corps de nombres dont nous noterons~$A'$ l'anneau des entiers. 


Posons $V=\of{[}{a_{\m},\tilde{a}_{\m}}{]}$. L'anneau de Banach $(\Bs(V),\|.\|_{V})$ n'est autre que l'anneau $(\hat{A}_{\m},|.|_{\m})$. Puisque le polyn\^ome~$P(S)$ est irr\'eductible dans~$k_{\m}[S]$, l'id\'eal maximal~$\m$ de~$A$ est divis\'e par un unique id\'eal maximal~$\m'=\m A'$ de~$A'$ et nous disposons d'un isomorphisme
$$u : \hat{A}_{\m}[S]/(P(S)) \xrightarrow[]{\sim} \hat{A'}_{\m'}.$$
Munissons l'anneau $\hat{A}_{\m}[S]/(P(S))$ de la norme
$$\|.\|'=|u(.)|_{\m'}.$$
C'est alors un anneau de Banach muni d'une norme uniforme. Notons~$W$ le segment~$\of{[}{a_{\m'},\tilde{a}_{\m'}}{]}$ de~$\Ms(A')$. L'isomorphisme~$u$ identifie alors les alg\`ebres norm\'ees $(\Bs(V)[S]/(P(S)),\|.\|')$ et $(\Bs(W),\|.\|_{W})$.

Puisque le polyn\^ome~$P$ est unitaire, le morphisme de~$\Bs(V)$-modules
$$n : \begin{array}{ccc}
\Bs(V)^d & \to & \Bs(V)[S]/(P(S))\\
(a_{0},\ldots,a_{d-1}) & \mapsto & \disp \sum_{i=0}^{d-1} a_{i}\, S^i
\end{array}$$
est un isomorphisme. Munissons l'alg\`ebre $\Bs(V)^d$ de la norme~$\|.\|_{\infty}$ donn\'ee par le maximum des normes des coefficients. Nous d\'efinissons alors une norme, not\'ee~$\|.\|_{V,\textrm{div}}$, sur~$\Bs(V)[S]/(P(S))$ de la fa\c{c}on suivante :
$$\forall f\in \Bs(V)[S]/(P(S)),\, \|f\|_{V,\textrm{div}} = \|n^{-1}(f)\|_{\infty}.$$

Pour appliquer la proposition~\ref{isolocal}, nous devons d\'emontrer que les normes~$\|.\|'$ et~$\|.\|_{V,\textrm{div}}$, d\'efinies sur~$\Bs(V)$, sont \'equivalentes. Tel est bien le cas car ce sont deux normes sur un m\^eme $\hat{K}_{\m}$-espace vectoriel de dimension finie qui induisent la m\^eme valeur absolue sur~$\hat{K}_{\m}$, \`a savoir~$|.|_{\m}$.

Notons
$$Y=\E{n}{\Bs(V)} \textrm{ et } Y'=\E{n}{\Bs(V)[S]/(P(S))}.$$
Notons encore 
$$\varphi : Y' \to Y \textrm{ et } \psi : \E{n}{A'} \to \E{n}{A}$$
les morphismes naturels. La partie~$V$ est une partie compacte et connexe de~$\Ms(A)$. Notons~$L_{V}$ son image r\'eciproque dans~$\E{n}{A}$. La partie~$W$ est une partie compacte et connexe de~$\Ms(A')$. Notons~$L'_{W}$ son image r\'eciproque dans~$\E{n}{A'}$. Consid\'erons, \`a pr\'esent, le diagramme commutatif suivant :
$$\xymatrix{
Y' \ar[r]^\varphi \ar[d]_{\chi'} &  Y \ar[d]_{\chi}\\
\E{n}{A'}  \ar[r]^\psi & \E{n}{A}}.$$
D'apr\`es la proposition~\ref{prorat}, les parties compactes~$V$ et~$W$ sont pro-rationnelles. D'apr\`es la proposition \ref{compactrationnelrelatif}, les morphismes~$\chi$ et~$\chi'$ sont des isomorphismes d'espaces annel\'es au-dessus, respectivement, de l'int\'erieur de~$L_{V}$ et de l'int\'erieur de~$L'_{W}$. Remarquons que le point~$x$ appartient \`a la fibre extr\^eme~$\pi^{-1}(\tilde{a}_{\m})$, situ\'ee \`a l'int\'erieur de~$L_{V}$. En outre, tout ant\'ec\'edent de~$x$ par le morphisme~$\psi$ appartient \`a la fibre extr\^eme situ\'ee au-dessus de~$\tilde{a}_{\m'}$ et donc \`a l'int\'erieur de~$L'_{W}$. Nous noterons encore~$x$ l'ant\'ec\'edent du point~$x$ par le morphisme~$\chi$. Pour conclure, il nous suffit de trouver un point~$x'$ de~$Y'$, rationnel dans sa fibre, tel que le morphisme~$\varphi$ induise un isomorphisme d'un voisinage de~$x'$ sur un voisinage de~$x$.

Notons~$\alpha$ l'image de~$S$ dans l'anneau~$\Bs(V)[S]/(P(S))$. D'apr\`es la proposition~\ref{Hensel}, il existe une fonction~$R$ d\'efinie sur un voisinage~$U$ de~$x$ dans~$Y$ telle que~$P(R)=0$ et $R(x)=\tilde{\alpha}$ dans~$\Hs(x)$. Construisons alors une section~$\sigma$ du morphisme~$\varphi$ au-dessus de~$U$, par le proc\'ed\'e d\'ecrit imm\'ediatement avant la proposition~\ref{isolocal}. Par sa d\'efinition m\^eme, nous avons 
$$S(\sigma(x))=R(x) \textrm{ dans } \Hs(x),$$
autrement dit,
$$R(\sigma(x))=\alpha \textrm{ dans } \Hs(\sigma(x)).$$  

Soit $b$ un point de $\Ms(\Bs(V))$. Le corps~$\Hs(b)$ est \'egal au corps~$k_{\m}$ ou au corps~$\hat{K}_{\m}$. Dans tous les cas, l'image du polyn\^ome $P(T)$ est irr\'eductible dans $\Hs(b)[T]$. Puisque le corps~$\Hs(b)$ est parfait, elle est \'egalement s\'eparable. Soit $c$ un point de $\Ms(\Bs(V)[S]/(P(S)))$ au-dessus du point~$b$. L'\'el\'ement~$\alpha$ de l'anneau $\Bs(V)[S]/(P(S))$ s'envoie sur une racine du polyn\^ome~$P(T)$ dans~$\Hs(c)$. Puisque le polyn\^ome~$P$ est s\'eparable, nous avons $P'(\alpha)=0$.

Pour finir, d'apr\`es le corollaire~\ref{cpainterne}, le point~$x$ poss\`ede, dans~$X$, et donc dans~$Y$, un syst\`eme fondamental de voisinages connexes. Nous pouvons donc appliquer la proposition~\ref{isolocal}. Nous obtenons, au voisinage du point~$x$, une section du morphisme~$\varphi$ qui est un isomorphisme local.

Pour conclure, il nous reste \`a montrer que le point $x'=\sigma(x)$ est rationnel dans sa fibre. Consid\'erons la projection~$b'$ de ce point sur~$\Ms(\Bs(V)[S]/(P(S)))$. Par d\'efinition, le caract\`ere associ\'e est
$$\begin{array}{ccc}
\Bs(V)[S]/(P(S)) & \to & \Hs(x)\\
Q(S) & \mapsto & Q(R(x))=Q(\tilde{\alpha})
\end{array}.$$
L'image de ce morphisme est le corps~$k_{\m}[\tilde{\alpha}]=\Hs(x)=\Hs(x')$. On en d\'eduit que le morphisme $\Hs(b) \to \Hs(x')$ est un isomorphisme et donc que le point~$x'$ est rationnel dans sa fibre.
\end{proof}

Int\'eressons-nous, \`a pr\'esent, aux fibres internes.

\begin{prop}\label{chgtinterne}\index{Isomorphisme local!au voisinage d'un point rigide!d'une fibre interne}
Soient $\tau$ un \'el\'ement de~$\Sigma$, $\eps$ un \'el\'ement de $\of{]}{0,l(\tau)}{[}$ et $x$ un point rigide de la fibre interne $X_{a_{\tau}^\eps}$. Supposons que le point $x$ poss\`ede un syst\`eme fondamental de voisinages connexes par arcs. Alors, il existe une extension finie~$K'$ de~$K$, un point~$x'$ de~$\E{n}{A'}$, o\`u~$A'$ d\'esigne l'anneau des entiers de~$K'$, rationnel dans sa fibre, tel que le morphisme naturel
$$\E{n}{A'} \to \E{n}{A}$$
envoie le point~$x'$ sur le point~$x$ et induise un isomorphisme d'un voisinage de~$x'$ sur un voisinage de~$x$. 
\end{prop}
\begin{proof}
L'extension de corps $\hat{K}_{\tau} \to \Hs(x)$ est une extension finie et s\'eparable, puisque la caract\'eristique du corps~$\hat{K}_{\tau}$ est nulle. D'apr\`es le th\'eor\`eme de l'\'el\'ement primitif, il existe un \'el\'ement $\alpha$ de $\Hs(x)$ tel que $\hat{K}_{\tau}[\alpha]=\Hs(x)$. Si le corps~$\hat{K}_{\tau}$ est ultram\'etrique, le lemme de Krasner assure que nous pouvons supposer que l'\'el\'ement~$\alpha$ est alg\'ebrique sur le corps~$K$. Si le corps~$\hat{K}_{\tau}$ est archim\'edien, nous pouvons encore supposer que~$\alpha$ est alg\'ebrique sur le corps~$K$, et m\^eme que c'est une racine carr\'ee de~$-1$. Notons $P(S) \in K[S]$ le polyn\^ome minimal unitaire de $\alpha$ sur $K$. Ce polyn\^ome est encore irr\'eductible sur le corps~$\hat{K}_{\tau}$. Consid\'erons l'extension finie~$K' = K[S]/(P(S))$ de~$K$. C'est un corps de nombres dont nous noterons~$A'$ l'anneau des entiers.


Soient~$\lambda\in\of{]}{0,\eps}{[}$ et~$\mu\in\of{]}{\eps,l(\tau)}{[}$. Posons $V=\of{[}{a_{\tau}^\lambda,a_{\tau}^\mu}{]}$. L'anneau de Banach $(\Bs(V),\|.\|_{V})$ n'est autre que l'anneau $(\hat{K}_{\tau},\max(|.|_{\tau}^\lambda,|.|_{\tau}^\mu))$. Puisque le polyn\^ome~$P(S)$ est irr\'eductible dans~$\hat{K}_{\tau}[S]$, la place~$\tau$ de~$K$ se prolonge en une unique place~$\tau'$ de~$K'$ et nous disposons d'un isomorphisme
$$u : \hat{K}_{\tau}[S]/(P(S)) \xrightarrow[]{\sim} \hat{K'}_{\tau'}.$$
Munissons l'anneau $\hat{K}_{\tau}[S]/(P(S))$ de la norme
$$\|.\|'=\max(|u(.)|_{\tau'}^\lambda,|u(.)|_{\tau'}^\mu).$$
C'est alors un anneau de Banach muni d'une norme uniforme. Notons~$W$ le segment~$\of{[}{a_{\tau'}^\lambda,a_{\tau'}^\mu}{]}$ de~$\Ms(A')$. L'isomorphisme~$u$ identifie alors les alg\`ebres norm\'ees $(\Bs(V)[S]/(P(S)),\|.\|')$ et $(\Bs(W),\|.\|_{W})$.

Introduisons une notation. Soient~$L$ une~$K$-alg\`ebre et~$\|.\|$ une norme sur~$L$. Puisque le polyn\^ome~$P$ est unitaire, le morphisme de~$L$-modules
$$n_{L} : \begin{array}{ccc}
L^d & \to & L[S]/(P(S))\\
(a_{0},\ldots,a_{d-1}) & \mapsto & \disp \sum_{i=0}^{d-1} a_{i}\, S^i
\end{array}$$
est un isomorphisme. Munissons l'alg\`ebre $L^d$ de la norme~$\|.\|_{\infty}$ donn\'ee par le maximum des normes des coefficients. Nous d\'efinissons alors une norme, not\'ee~$\|.\|_{\textrm{div}}$, sur~$L[S]/(P(S))$ de la fa\c{c}on suivante :
$$\forall f\in L[S]/(P(S)),\, \|f\|_{\textrm{div}} = \|n_{L}^{-1}(f)\|_{\infty}.$$

Pour appliquer la proposition~\ref{isolocal}, nous devons d\'emontrer que les normes~$\|.\|'$ et~$\|.\|_{V,\textrm{div}}$, d\'efinies sur~$\Bs(V)$, sont \'equivalentes. Or la norme~$\|.\|_{V,\textrm{div}}$ est \'e\-qui\-va\-len\-te \`a la norme $\max(|.|_{\tau,\textrm{div}}^\lambda,|.|_{\tau,\textrm{div}}^\mu)$. Il nous suffit, \`a pr\'esent, de remarquer que, quel que soit~$\nu\in\{\lambda,\mu\}$, les normes~$|.|_{\tau,\textrm{div}}^nu$ et~$|.|_{\tau'}^\nu$ sont \'equivalentes. En effet, ce sont deux normes sur un m\^eme $\hat{K}_{\tau}$-espace vectoriel de dimension finie qui induisent la m\^eme valeur absolue sur~$\hat{K}_{\tau}$, \`a savoir~$|.|_{\tau}^\nu$.

Le reste du raisonnement se d\'eroule exactement comme dans la preuve pr\'e\-c\'e\-den\-te.
\end{proof}

Pour terminer, traitons le cas de la fibre centrale.

\begin{prop}\label{isocentral}\index{Isomorphisme local!au voisinage d'un point rigide!de la fibre centrale}
Soit $x$ un point rigide de la fibre centrale $X_{0}$. Supposons que le point $x$ poss\`ede un syst\`eme fondamental de voisinages connexes par arcs. Alors, il existe une extension finie $K'$ de $K$, un point $x'$ de $\E{n}{A'}$, o\`u $A'$ d\'esigne l'anneau des entiers de $K'$, rationnel dans sa fibre, tel que le morphisme naturel
$$\E{n}{A'} \to \E{n}{A}$$
envoie le point $x'$ sur le point $x$ et induise un isomorphisme d'un voisinage de~$x'$ sur un voisinage de $x$. 
\end{prop}
\begin{proof}
L'extension de corps $\Hs(a_{0}) = K \to \Hs(x)$ est une extension finie et s\'eparable, puisque la caract\'eristique du corps~$K$ est nulle. D'apr\`es le th\'eor\`eme de l'\'el\'ement primitif, il existe un \'el\'ement~$\alpha$ de~$\Hs(x)$ tel que~$K[\alpha]=\Hs(x)$. Notons~\mbox{$P(S)\in K[S]$} le polyn\^ome minimal unitaire de~$\alpha$ sur~$K$. Il existe un unique isomorphisme 
$$K[S]/(P(S)) \xrightarrow[]{\sim} \Hs(x)$$ 
envoyant $S$ sur $\alpha$.

Le caract\`ere s\'eparable de l'extension $\Hs(x)/K$ assure \'egalement que l'anneau des entiers $A'$ de $\Hs(x)$ est un anneau de Dedekind de type fini sur $A$. Par cons\'equent, il existe un \'el\'ement $f$ de $K$ tel que 
$$A[f,\alpha]=A'[f].$$
Choisissons un sous-ensemble fini~$\Sigma_{0}$ de~$\Sigma_{f}$ de sorte que la fonction~$f$ soit d\'efinie et inversible sur l'ouvert de~$B$ d\'efini par
$$U = B \setminus \bigcup_{\m\in\Sigma_{0}} \{\tilde{a}_{\m}\}.$$ 
Quitte \`a augmenter l'ensemble~$\Sigma_{0}$, nous pouvons supposer que les coefficients du polyn\^ome~$P(S)$ sont d\'efinis en tout point de~$U$ et que, quel que soit~$b$ dans~$U$, l'image du polyn\^ome~$P(S)$ est s\'eparable sur~$\Hs(b)$. Pour~$\m$ dans~$\Sigma_{f}$, notons~$r(\m)$ l'ensemble des id\'eaux maximaux de~$A'$ divisant l'id\'eal maximal~$\m$ de~$A$. Pour~$\sigma$ dans~$\Sigma_{\infty}$, notons~$r(\sigma)$ l'ensemble des plongements complexes de~$\Hs(x)$ \`a conjugaison pr\`es qui prolongent $\sigma$. Notons
$$\Sigma'_{0}=\bigcup_{\m\in\Sigma_{0}} r(\m) \textrm{ et } \Sigma'_{\infty}=\bigcup_{\sigma\in\Sigma_{\infty}}.$$ 

Consid\'erons la partie compacte contenue dans~$U$ d\'efinie par
$$M = B \setminus \bigcup_{\m\in\Sigma_{0}} \of{]}{a_{\m},\tilde{a}_{\m}}{]}.$$
Consid\'erons l'alg\`ebre de Banach $(\Bs(M),\|.\|_{M})$. Nous avons
$$\Bs(M) = A\left[\frac{1}{\Sigma_{0}}\right] = \left\{ \frac{a}{b} \in K,\, a,b\in A,\, b\ne 0,\, \exists \m\in\Sigma_{0},\, b\in\m \right\}$$
et
$$\|.\|_{M} = \max_{\sigma\in\Sigma_{0}\cup\Sigma_{\infty}}(|.|_{\sigma}).$$

Le compact~$M$ \'etant contenu dans~$U$, l'anneau~$\As$ est un localis\'e de l'anneau~$A[u]$. On en d\'eduit que le morphisme
$$A\left[\frac{1}{\Sigma_{0}}\right][S]/(P(S)) \xrightarrow[]{\sim} A'\left[\frac{1}{\Sigma'_{0}}\right]$$
est un isomorphime. Munissons l'anneau $\Bs(M)[S]/(P(S))$ de la norme
$$\|.\|'=\max_{\sigma\in\Sigma'_{0}\cup\Sigma'_{\infty}}(|.|_{\sigma}).$$
C'est alors un anneau de Banach muni d'une norme uniforme. 

Introduisons une notation. Soient~$L$ une~$K$-alg\`ebre et~$\|.\|$ une norme sur~$L$. Puisque le polyn\^ome~$P$ est unitaire, le morphisme de~$K$-espaces vectoriels
$$n_{L} : \begin{array}{ccc}
L^d & \to & L[S]/(P(S))\\
(a_{0},\ldots,a_{d-1}) & \mapsto & \disp \sum_{i=0}^{d-1} a_{i}\, S^i
\end{array}$$
est un isomorphisme. Munissons l'alg\`ebre $L^d$ de la norme~$\|.\|_{\infty}$ donn\'ee par le maximum des normes des coefficients. Nous d\'efinissons alors une norme, not\'ee~$\|.\|_{\textrm{div}}$, sur~$L[S]/(P(S))$ de la fa\c{c}on suivante :
$$\forall f\in L[S]/(P(S)),\, \|f\|_{\textrm{div}} = \|n_{L}^{-1}(f)\|_{\infty}.$$

Pour appliquer la proposition~\ref{isolocal}, nous devons d\'emontrer que les normes~$\|.\|'$ et~$\|.\|_{M,\textrm{div}}$, d\'efinies sur~$\Bs(M)$, sont \'equivalentes. Or la norme~$\|.\|_{M,\textrm{div}}$ est \'e\-qui\-va\-len\-te \`a la norme 
$$\max_{\sigma\in\Sigma_{0}\cup\Sigma_{\infty}}(|.|_{\sigma,\textrm{div}}).$$
Soit~$\m\in\Sigma_{0}$. Nous disposons alors d'un isomorphisme
$$\hat{K}_{\m}[S]/(P(S)) \xrightarrow[]{\sim} \prod_{\m'\in r(\m)} \hat{\Hs(x)}_{\m'}.$$
On en d\'eduit que les normes~$|.|_{\m,\textrm{div}}$ et~$\max_{\m'\in r(\m)} (|.|_{\m'})$ sont \'equivalentes car ce sont deux normes sur un m\^eme $\hat{K}_{\m}$-espace vectoriel de dimension finie qui induisent la m\^eme valeur absolue sur~$\hat{K}_{\m}$, \`a savoir~$|.|_{\m}$.

On raisonne de m\^eme pour les \'el\'ements de~$\Sigma_{\infty}$ en prenant garde au fait que l'isomorphisme ne vaut que si l'on consid\`ere tous les plongements complexes et pas seulement les classes de conjugaison.

Le reste du raisonnement se d\'eroule exactement comme dans la preuve de la premi\`ere proposition.

\end{proof}

Pour plus de clart\'e, nous regroupons les trois r\'esultats obtenus dans la proposition suivante.

\begin{prop}\label{isorigcpa}
Soit $x$ un point rigide de l'une des fibres de l'espace~$X$. Supposons que le point $x$ poss\`ede un syst\`eme fondamental de voisinages connexes par arcs. Alors, il existe une extension finie~$K'$ de~$K$, un point~$x'$ de~$\E{n}{A'}$, o\`u~$A'$ d\'esigne l'anneau des entiers de~$K'$, rationnel dans sa fibre, tel que le morphisme naturel
$$\E{n}{A'} \to \E{n}{A}$$
envoie le point~$x'$ sur le point~$x$ et induise un isomorphisme d'un voisinage de~$x'$ sur un voisinage de~$x$. 
\end{prop}

\subsection{Voisinages sur la droite}

Pour utiliser la proposition qui pr\'ec\`ede, il est n\'ecessaire de disposer d'un r\'esultat de connexit\'e locale. Nous consacrons donc une section \`a l'\'etude de la topologie au voisinage des points rigides des fibres dans le cas le plus simple : celui de la droite. Dans les propositions qui suivent, nous supposerons donc que $n=1$ et que $X=\E{1}{A}$. 

\begin{prop}\label{voisrig1}\index{Voisinages d'un point!rigide de A1@rigide de \E{1}{A}}
Soient $b$ un point de~$B$ et~$P(T)$ un polyn\^ome unitaire \`a coefficients dans~$\Os_{B,b}$ dont l'image dans~$\Hs(b)[T]$ est irr\'eductible. Soit~$x$ le point de la fibre~$X_{b}$ d\'efini par l'\'equation $P(T)(x)=0$. Soient $B_{0}$ un voisinage de $b$ dans $B$ sur lequel les coefficients du polyn\^ome~$P$ sont d\'efinis.

Soit~$U$ un voisinage du point~$x$ dans~$X$. Il existe un voisinage~$V$ du point~$b$ dans~$B_{0}$ et un nombre r\'eel $t>0$ tels que l'on ait l'inclusion
$$\left\{y\in X_{V}\, \big|\, |P(T)(y)| < t\right\} \subset U.$$
\end{prop}
\begin{proof}
D'apr\`es le corollaire \ref{partiecompacte}, pour toute partie compacte~$V$ de~$B_{0}$ et tout \'el\'ement $s$ de~$\R_{+}$, la partie de~$X$ d\'efinie par
$$\{y\in X_{V}\, |\, |P(T)(y)| \le s\}$$
est compacte. Le r\'esultat d\'ecoule alors du lemme \ref{lemvois}.
\end{proof}

Nous souhaitons, \`a pr\'esent, montrer que les voisinages qui figurent dans l'\'enonc\'e de la proposition sont connexes par arcs lorsque leur projection sur la base l'est. \`A cet effet, nous commencerons par d\'emontrer quelques r\'esultats sur la topologie des fibres.

\begin{lem}\label{ccfibresum}
Soit $(k,|.|)$ un corps valu\'e, ultram\'etrique, maximalement complet et alg\'ebriquement clos. Soient $d\in\N$, $\alpha_{1},\ldots,\alpha_{d}\in k$ et $t\in\R_{+}^*$. Posons 
$$P(T) = \prod_{i=1}^d (T-\alpha_{i})$$
et 
$$U = \left\{ x\in \E{1}{k}\, \big|\, |P(T)(x)|<t\right\}.$$
Alors, pour tout point $y$ de $U$, il existe un chemin trac\'e sur~$U$ qui joint le point~$y$ \`a l'un des points~$\alpha_{i}$, avec~$i\in\cn{1}{d}$.  
\end{lem}
\begin{proof}
Soit $y$ un point de~$U$. Puisque le corps~$k$ est maximalement complet, il existe $\beta\in k$ et $r\in\R_{+}$ tels que~$y = \eta_{\beta,r}$ dans~$\E{1}{k}$. Supposons, tout d'abord, qu'il existe~$i\in\cn{1}{d}$ tel que~$\beta=\alpha_{i}$. Consid\'erons alors le chemin
$$l : \begin{array}{ccc}
\of{[}{0,1}{]} & \to & \disp \E{1}{k}\\
u & \mapsto & \disp \eta_{\alpha_{i},(1-u)r}
\end{array}.$$
Il joint le point~$y$ au point~$\alpha_{i}$ et tout polyn\^ome d\'ecro\^{\i}t le long de ce chemin. En particulier, il est \`a valeurs dans~$U$.

Revenons, \`a pr\'esent, au cas g\'en\'eral. Nous distinguerons deux cas. Dans un premier temps, supposons, qu'il existe~\mbox{$i\in\cn{1}{d}$} tel que~$|\beta-\alpha_{i}|\le r$. Alors le point~$y=\eta_{\beta,r}$ n'est autre que le point~$\eta_{\alpha_{i},r}$ et nous sommes ramen\'es au cas pr\'ec\'edent. Il nous reste \`a traiter le cas o\`u, quel que soit~$i\in\cn{1}{d}$, nous avons~$|\beta-\alpha_{i}| > r$. Dans ce cas, nous avons
$$|P(T)(\eta_{\beta,r})| = \prod_{i=1}^d |(T-\alpha_{i})(\eta_{\beta,r})| = \prod_{i=1}^d |\beta-\alpha_{i}|.$$
Notons $s = \min_{1\le i\le d}(|\beta-\alpha_{i}|)$. Consid\'erons le chemin
$$l' : \begin{array}{ccc}
\of{[}{0,1}{]} & \to & \disp \E{1}{s}\\
u & \mapsto & \disp \eta_{\beta,(1-u)r+us}
\end{array}.$$
Il joint le point $y$ au point~$\eta_{\beta,s}$, qui est du type consid\'er\'e pr\'ec\'edemment. En outre, la fonction~$P$ est constante le long du chemin~$l'$. Il est donc bien \`a valeurs dans~$U$. On en d\'eduit le r\'esultat annonc\'e.
\end{proof}

\begin{lem}\label{ccfibresarc}
Soient $d\in\N$, $\alpha_{1},\ldots,\alpha_{d}\in \C$ et $t\in\R_{+}^*$. Posons 
$$P(T) = \prod_{i=1}^d (T-\alpha_{i})$$
et 
$$U = \left\{ z\in \C\, \big|\, |P(z)|_{\infty}<t\right\}.$$
Alors, pour tout point $y$ de $U$, il existe un chemin trac\'e sur~$U$ qui joint le point~$y$ \`a l'un des points~$\alpha_{i}$, avec~$i\in\cn{1}{d}$.  
\end{lem}
\begin{proof}
Consid\'erons l'application continue
$$\begin{array}{ccc}
\C & \to & \C\\
z & \mapsto & P(z)
\end{array}.$$
C'est un rev\^etement ramifi\'e. Consid\'erons le chemin trac\'e sur la base
$$\begin{array}{ccc}
\of{[}{0,1}{]} & \to & \C\\
u & \mapsto & (1-u)\,P(y)
\end{array}.$$
En relevant ce chemin \`a partir du point~$y$, on obtient un chemin trac\'e sur~$U$ qui aboutit \`a l'un des racines du polyn\^ome~$P$. 
\end{proof}

\begin{cor}\label{ccfibres}
Soit $(k,|.|)$ un corps valu\'e complet. Soient $d$ un entier, \mbox{$Q_{1}(T),\ldots,Q_{d}(T)\in k[T]$} des polyn\^omes irr\'eductibles et $t\in\R_{+}^*$ un nombre r\'eel strictement positif. Pour~$i\in\cn{1}{d}$, notons~$x_{i}$ le point de la droite~$\E{1}{k}$ d\'efini par l'\'equation $Q_{i}(T)(x_{i})=0$. Posons 
$$P(T) = \prod_{i=1}^d Q_{i}$$
et 
$$U = \left\{ x\in \E{1}{k}\, \big|\, |P(T)(x)|<t\right\}.$$
Alors, pour tout point $y$ de $U$, il existe un chemin trac\'e sur~$U$ qui joint le point~$y$ \`a l'un des points~$x_{i}$, avec~$i\in\cn{1}{d}$.  

En particulier, si le polyn\^ome~$P(T)$ est une puissance d'un polyn\^ome irr\'eductible, alors la partie~$U$ est connexe par arcs.
\end{cor}
\begin{proof}
Soit $(L,|.|)$ une extension du corps valu\'e $(k,|.|)$. Le morphisme induit
$$\E{1}{L} \to \E{1}{k}$$
est continu et surjectif. On en d\'eduit qu'il suffit de d\'emontrer le r\'esultat pour une extension de $k$. Nous pouvons donc utiliser nous ramener \`a la situation du lemme \ref{ccfibresum}, si la valeur absolue $|.|$ est ultram\'etrique, ou du lemme \ref{ccfibresarc}, si elle est archim\'edienne.
\end{proof}

Revenons, \`a pr\'esent, aux voisinages des points rigides dans l'espace total.

\begin{prop}\label{voiscparig1}
Soient $b$ un point de $B$ et $V$ un voisinage connexe par arcs de $b$ dans $B$. Soit $P(T)\in\Os(V)[T]$ un polyn\^ome unitaire dont l'image dans~$\Hs(b)[T]$ est irr\'eductible. Soit~$t\in\R_{+}^*$ un nombre r\'eel strictement positif. Alors, la partie~$U$ de~$X=\E{1}{A}$ d\'efinie par
$$U = \left\{ y\in X\, \big|\, \pi(y)\in V,\, |P(T)(y)|<t \right\}$$
est connexe par arcs.
\end{prop}
\begin{proof}
Nous noterons~$x$ l'unique point de la fibre~$X_{b}$ qui v\'erifie
$$P(T)(x)=0.$$
Nous allons montrer que tout point de $U$ peut \^etre joint au point~$x$ par un chemin trac\'e sur~$U$. Nous allons distinguer plusieurs cas selon le type du point~$b$.\\ 

Supposons, tout d'abord, que le point~$b$ est le point central~$a_{0}$ de~$B$. Soit~$y$ un point de~$U$. Posons~$c=\pi(y)$. D\'ecomposons le polyn\^ome $P(T)$ en produit de facteurs irr\'eductibles et unitaires dans $\Hs(c)[T]$ : il existe~$d\in\N^*$, $Q_{1}(T),\ldots,Q_{d}(T)$ des polyn\^omes irr\'eductibles distincts et~$n_{1},\ldots,n_{d}\in\N^*$ tels que 
$$P(T) = \prod_{i=1}^d Q_{i}(T)^{n_{i}} \textrm{ dans } \Hs(c)[T].$$
Quel que soit~$i\in\cn{1}{d}$, notons~$y_{i}$ le point de la fibre~$X_{c}$ d\'efini par l'\'equation $Q_{i}(T)(y_{i})=0$. D'apr\`es le lemme~\ref{ccfibres}, il existe un indice~$i\in\cn{1}{d}$ et un chemin trac\'e sur~$X_{c}\cap U$ qui joint le point $y$ au point $y_{i}$.  Il nous reste \`a montrer que l'on peut joindre le point~$y_{i}$ au point~$x$ par un chemin trac\'e sur~$U$. Si le point~$c$ est le point~$a_{0}$, c'est \'evident. 

Supposons, que le point~$c$ est un point interne de~$B$. Il existe alors~$\sigma\in\Sigma$ et~\mbox{$\eps>0$} tels que~$c=a_{\sigma}^\eps$. Puisque la partie~$V$ est suppos\'ee connexe par arcs, elle contient le segment~$W=\of{[}{a_{0},a_{\sigma}^\eps}{]}$. Remarquons que, quel que soit~$\lambda\in\of{]}{0,\eps}{]}$, le polyn\^ome~$Q_{i}(T)$ est encore irr\'eductible dans~$\Hs(a_{\sigma}^\lambda)[T]$. D\'efinissons une section~$\varphi$ de~$\pi$ au-dessus de~$W$ de la fa\c{c}on suivante : au point~$a_{\sigma}^\lambda$, avec~$\lambda\in\of{]}{0,\eps}{]}$, nous associons le point~$\varphi(a_{\sigma}^\lambda)$ de la fibre~$X_{a_{\sigma}^\lambda}$ d\'efini par l'\'equation $Q_{i}(T)(\varphi(a_{\sigma}^\lambda))=0$ et au point~$a_{0}$, nous associons le point~$\varphi(a_{0})=x$. L'application~$\varphi$ est une section continue de~$\pi$ au-dessus de~$W$ \`a valeurs dans~$U$ et son image est un chemin joignant le point~$y_{i}$ au point~$x$.

Pour finir, supposons que point~$c$ est un point extr\^eme de~$B$. Il existe alors \mbox{$\m\in\Sigma_{f}$} tel que~$c=\tilde{a}_{\m}$. La d\'ecomposition $P(T) = \prod_{i=1}^d Q_{i}(T)^{n_{i}}$ vaut donc dans l'anneau de polyn\^omes~$k_{\m}[T]$. Le lemme de Hensel nous assure qu'il existe des polyn\^omes~$R_{1}(T),\ldots,R_{d}(T)\in \hat{A}_{\m}$ unitaires tels que l'on ait la d\'ecomposition
$$P(T) = \prod_{i=1}^d R_{i}(T) \textrm{ dans } \hat{A}_{\m}[T]$$
et, quel que soit~$i\in\cn{1}{d}$, 
$$R_{i}(T) = Q_{i}(T)^{n_{i}} \mod \m.$$
Choisissons un facteur irr\'eductible $S_{i}(T)$ du polyn\^ome~$R_{i}(T)$ dans~$\hat{A}_{\m}[T]$. Puisque la partie~$V$ est suppos\'ee connexe par arcs, elle contient le segment~$W=\of{[}{a_{0},\tilde{a}_{\m}}{]}$. Nous d\'efinissons alors une section~$\varphi$ de~$\pi$ au-dessus de~$W$ de la fa\c{c}on suivante : au point~$a_{\m}^\lambda$, avec~$\lambda\in\of{]}{0,\infty}{[}$, nous associons le point~$\varphi(a_{\sigma}^\lambda)$ de la fibre~$X_{a_{\sigma}^\lambda}$ d\'efini par l'\'equation~$S_{i}(T)(\varphi(a_{\sigma}^\lambda))=0$, au point~$a_{0}$ nous associons le point~$\varphi(a_{0})=x$ et au point~$\tilde{a}_{\m}$, nous associons le point~$y_{i}$. Comme pr\'ec\'edemment, l'application~$\varphi$ est une section continue de~$\pi$ au-dessus de~$W$ \`a valeurs dans~$U$ et son image est un chemin joignant le point~$y_{i}$ au point~$x$.\\ 

Supposons, \`a pr\'esent, que le point~$b$ est un point extr\^eme de~$B$. Il existe alors \mbox{$\m\in\Sigma_{f}$} tel que~$b=\tilde{a}_{\m}$. Supposons, dans un premier temps que $a_{0}\in V$. Alors le polyn\^ome~$P(T)$ est \`a coefficients dans $A_{\m}$ et il est irr\'eductible dans~$A_{\m}[T]$ puisqu'il est unitaire et que sa r\'eduction modulo~$\m$ est irr\'eductible. Nous sommes donc ramen\'es au cas pr\'ec\'edent. 

Supposons, \`a pr\'esent, que le point central~$a_{0}$ n'appartient pas \`a~$V$. Si la partie~$V$ est r\'eduite au point extr\^eme~$\tilde{a}_{\m}$, le r\'esultat provient directement du lemme~\ref{ccfibres}. Dans les autres cas, la partie~$V$ est un intervalle contenu dans~$\of{]}{a_{0},\tilde{a}_{\m}}{]}$. Le polyn\^ome~$P(T)$ est alors \`a coefficients dans~$\hat{A}_{\m}$. Puisqu'il est unitaire et que son image modulo~$\m$ est irr\'eductible, il est irr\'eductible dans~$\hat{A}_{\m}[T]$ et donc dans~$\hat{K}_{\m}[T]$. Soit~$y$ un point de~$U$. Il existe alors~\mbox{$\eps\in\of{]}{0,+\infty}{]}$} tel que~$\pi(y)=a_{\m}^\eps$. D'apr\`es le lemme~\ref{ccfibres}, il existe un chemin trac\'e sur \mbox{$X_{a_{\m}^\eps} \cap U$} joignant le point~$y$ au point~$z$ d\'efini par l'\'equation~$P(T)(z)=0$. Il nous suffit, \`a pr\'esent, de montrer que l'on peut joindre le point~$z$ au point~$x$ par un chemin trac\'e sur~$U$. Puisque la partie~$V$ est suppos\'ee connexe par arcs, elle contient le segment~$W=\of{[}{a_{\m}^\eps,\tilde{a}_{\m}}{]}$. D\'efinissons une section~$\varphi$ de~$\pi$ au-dessus de~$W$ de la fa\c{c}on suivante : \`a tout point~$c$ de~$W$ nous associons le point~$\varphi(c)$ de la fibre~$X_{c}$ d\'efini par l'\'equation $P(T)(\varphi(c))=0$. L'application~$\varphi$ est une section continue de~$\pi$ au-dessus de~$W$ \`a valeurs dans~$U$ et son image est un chemin joignant le point~$z$ au point~$x$.\\

Il nous reste \`a traiter le cas o\`u le point~$b$ est un point interne de~$B$ : il existe~$\sigma\in\Sigma$ et~$\eps>0$ tel que~$b=a_{\sigma}^\eps$. Si la partie~$V$ contient un point extr\^eme ou le point central de~$B$, nous sommes ramen\'es \`a l'un des cas pr\'ec\'edents. Nous supposerons donc que la partie~$V$ est contenue dans~$B'_{\sigma}$. Dans ce cas, pour tout point~$c$ de~$V$, le corps~$\Hs(c)$ est isomorphe au corps~$\hat{K}_{\sigma}$ et le polyn\^ome~$P(T)$ est irr\'eductible dans~$\Hs(c)[T]$. Soit~$y$ un point de~$U$. D'apr\`es le lemme~\ref{ccfibres}, il existe un chemin trac\'e sur~$X_{\pi(y)} \cap U$ joignant le point~$y$ au point~$z$ d\'efini par l'\'equation~$P(T)(z)=0$. Il nous suffit, \`a pr\'esent, de montrer que l'on peut joindre le point~$z$ au point~$x$ par un chemin trac\'e sur~$U$. D\'efinissons une section~$\varphi$ de~$\pi$ au-dessus de~$V$ de la fa\c{c}on suivante : \`a tout point~$c$ de~$V$ nous associons le point~$\varphi(c)$ de la fibre~$X_{c}$ d\'efini par l'\'equation $P(T)(\varphi(c))=0$. L'application~$\varphi$ est une section continue de~$\pi$ au-dessus de~$V$ \`a valeurs dans~$U$ et son image est un chemin passant par les points~$z$ et~$x$.
\end{proof}

\begin{cor}\label{cparig1}\index{Connexite par arcs au voisinage d'un point@Connexit\'e par arcs au voisinage d'un point!rigide! de A1@de \E{1}{A}}
Soient $b$ un point de $B$ et $x$ un point rigide de la fibre~$X_{b}$. Alors, le point~$x$ poss\`ede un syst\`eme fondamental de voisinages connexes par arcs.
\end{cor}
\begin{proof}
D'apr\`es le lemme \ref{translation}, le morphisme naturel $\Os_{B,b} \to \Hs(b)$ est surjectif. Nous pouvons donc supposer que le polyn\^ome~$P(T)$ d\'efinissant le point~$x$ est \`a coefficients dans~$\Os_{B,b}$. Il nous suffit alors d'appliquer les propositions~\ref{voisrig1} et~\ref{voiscparig1}.
\end{proof}

\subsection{\'Etude topologique locale}

Revenons, \`a pr\'esent, au cas d'un espace affine de dimension quelconque : $X=\E{n}{A}$, avec $n\in\N$. Les r\'esultats obtenus sur la topologie de la droite nous permettent de mettre en {\oe}uvre un raisonnement par r\'ecurrence. 


\begin{prop}\label{cparig}\index{Connexite par arcs au voisinage d'un point@Connexit\'e par arcs au voisinage d'un point!rigide! de An@de \E{n}{A}}
Soit $b$ un point de $B$. Soit $x$ un point rigide de la fibre~$X_{b}$. Alors, le point~$x$ poss\`ede un syst\`eme fondamental de voisinages connexes par arcs.
\end{prop}
\begin{proof}
Nous allons d\'emontrer le r\'esultat attendu par r\'ecurrence sur l'entier $n\in\N$. Le cas $n=0$ est imm\'ediat. 

Soit $n\in\N^*$ tel que le r\'esultat soit vrai pour $n-1$. Notons 
$$\varphi_{1} : \E{n}{A} \to \E{1}{A}$$
le morphisme induit par l'injection $i_{1} : A[T_{1}] \to A[T_{1},\ldots,T_{n}]$ et 
$$\varphi_{0} : \E{1}{A} \to \Ms(A)$$
celui induit par l'injection $i_{0} : A \to A[T_{1}]$. Posons $y=\varphi_{1}(x)$. 

D'apr\`es la proposition \ref{cparig1}, le point $y$ de $\E{1}{A}$ poss\`ede un syst\`eme fondamental de voisinages connexes par arcs. Nous pouvons donc appliquer la proposition~\ref{isorigcpa}. Elle assure qu'il existe une extension finie $K'$ de $K$, un point $y'$ de $\E{1}{A'}$, o\`u $A'$ d\'esigne l'anneau des entiers de $K'$, rationnel dans sa fibre, tel que le morphisme naturel
$$\alpha : \E{1}{A'} \to \E{1}{A}$$
envoie le point $y'$ sur le point $y$ et induise un isomorphisme 
$$\beta : U' \to U$$
d'un voisinage~$U'$ de~$y'$ dans~$\E{1}{A'}$ sur un voisinage~$U$ de~$y$ dans~$\E{1}{A}$. Consid\'erons le diagramme commutatif suivant
$$\xymatrix{
\E{n}{A'} \ar^{\alpha_{n}}[r] \ar^{\varphi_{1}'}[d]  & \E{n}{A}  \ar^{\varphi_{1}}[d]\\
\E{1}{A'} \ar^{\varphi'_{0}}[d] \ar^{\alpha}[r]& \E{1}{A}  \ar^{\varphi_{0}}[d]\\
\Ms(A') \ar^{\alpha_{0}}[r]& \Ms(A)
}.$$

Quitte \`a restreindre le voisinage~$U$ de~$y$, nous pouvons supposer qu'il est compact et rationnel. Le voisinage~$U'$ l'est alors \'egalement. Nous pouvons donc appliquer le th\'eor\`eme \ref{compactrationnel} et la proposition \ref{compactrationnelrelatif}. On en d\'eduit un isomorphisme
$$\gamma : \Ms(\Bs(U')) \xrightarrow[]{\sim} \Ms(\Bs(U))$$
qui co\"{\i}ncide avec $\beta$ en tant qu'application et m\^eme en tant que morphisme d'espace annel\'es si l'on se restreint \`a l'int\'erieur des espaces consid\'er\'es. On en d\'eduit un diagramme commutatif
$$\xymatrix{
\E{n-1}{\Bs(U')}\ar^{\delta}_{\sim}[r]\ar^{\psi'}[d] & \E{n-1}{\Bs(U)}\ar^{\psi}[d]\\
 \Ms(\Bs(U'))\ar^{\gamma}_{\sim}[r] & \Ms(\Bs(U))
}.$$  
En tant que morphisme d'espaces topologiques, le morphisme~$\psi$ n'est autre que le morphisme~$\varphi_{1}$ restreint \`a~$\varphi_{1}^{-1}(U)$ \`a la source et~$U$ au but. De m\^eme, le morphisme~$\psi'$ co\"{\i}ncide avec le morphisme~$\varphi'_{1}$ restreint \`a~${\varphi'_{1}}^{-1}(U')$ \`a la source et~$U'$ au but. Par cons\'equent, il suffit de montrer que le point~$x$ poss\`ede un syst\`eme fondamental de voisinages connexes par arcs dans~$\E{n-1}{\Bs(U)}$. Puisque~$\delta$ est un hom\'eomorphisme, il suffit de montrer que le point~$\delta^{-1}(x)$ poss\`ede un syst\`eme fondamental de voisinages connexes par arcs dans~$\E{n-1}{\Bs(U')}$. Or le point~$\delta^{-1}(x)$ est envoy\'e sur le point~$\gamma^{-1}(y)=y'$ dans~$\Ms(\Bs(U'))$. Ce dernier point est rationnel dans sa fibre~${\varphi'_{0}}^{-1}(\varphi'_{0}(y'))$. Par cons\'equent, quitte \`a changer~$x$ en~$\delta^{-1}(x)$, nous pouvons supposer que le point~$y$ est rationnel dans sa fibre, autrement dit que le morphisme
$$\Hs(b) \xrightarrow[]{\sim} \Hs(y)$$  
est un isomorphisme.

Notons 
$$\lambda_{n-1} : \E{n}{A} \to \E{n-1}{A}$$
le morphisme induit par l'injection $j_{n-1} : A[T_{2},\ldots,T_{n}] \to A[T_{1},\ldots,T_{n}]$ et 
$$\lambda_{0} : \E{n-1}{A} \to \Ms(A)$$
celui induit par l'injection $j_{0} : A \to A[T_{2},\ldots,T_{n-1}]$. Posons $z=\lambda_{n-1}(x)$. De l'isomorphisme $\Hs(b) \xrightarrow[]{\sim} \Hs(y)$, on d\'eduit un isomorphisme
$$\Hs(z) \xrightarrow[]{\sim} \Hs(x).$$  
D'apr\`es l'hypoth\`ese de r\'ecurrence, le point $z$ de $\E{n-1}{A}$ poss\`ede un syst\`eme fondamental de voisinages connexes par arcs. Nous pouvons donc appliquer la proposition~\ref{isorigcpa}. Par le m\^eme raisonnement que pr\'ec\'edemment, nous en d\'eduisons que nous pouvons supposer que le point~$z$ est rationnel dans la fibre~${\lambda_{0}}^{-1}(b)$. Autrement dit, le morphisme  
$$\Hs(b) \xrightarrow[]{\sim} \Hs(z)$$  
est un isomorphisme. Nous nous sommes finalement ramen\'es au cas d'un point~$x$ rationnel dans sa fibre~$X_{b}$, puisque le morphisme~$\Hs(b)\to \Hs(x)$ est un isomorphisme. Or d'apr\`es le lemme \ref{translation}, le morphisme canonique $\Os_{B,b}\to \Hs(b)$ est surjectif. Nous pouvons donc appliquer le corollaire \ref{cpadep}. On en d\'eduit le r\'esultat attendu. 
\end{proof}

En utilisant cette proposition, nous pouvons rel\^acher les hypoth\`eses de la proposition~\ref{isorigcpa}.

\begin{prop}\label{isorig}\index{Isomorphisme local!au voisinage d'un point rigide!zzz@de \E{n}{A}}
Soit $x$ un point rigide de l'une des fibres de l'espace~$X$. Alors, il existe une extension finie~$K'$ de~$K$, un point~$x'$ de~$\E{n}{A'}$, o\`u~$A'$ d\'esigne l'anneau des entiers de~$K'$, rationnel dans sa fibre, tel que le morphisme naturel
$$\E{n}{A'} \to \E{n}{A}$$
envoie le point~$x'$ sur le point~$x$ et induise un isomorphisme d'un voisinage de~$x'$ sur un voisinage de~$x$. 
\end{prop}

\begin{cor}\label{ouvertrig}\index{Ouverture au voisinage d'un point!rigide de An@rigide de \E{n}{A}}
Soit $b$ un point de $B$. Soit $x$ un point rigide de la fibre~$X_{b}$. Alors, le morphisme $\pi$ est ouvert au point $x$.
\end{cor}
\begin{proof}
La proposition~\ref{isorig} assure qu'il existe une extension finie $K'$ de $K$, un point $x'$ de $\E{n}{A'}$, o\`u $A'$ d\'esigne l'anneau des entiers de $K'$, rationnel dans sa fibre, tel que le morphisme naturel
$$\alpha : \E{n}{A'} \to \E{n}{A}$$
envoie le point $x'$ sur le point $x$ et induise un isomorphisme 
$$\beta : U' \to U$$
d'un voisinage~$U'$ de~$x'$ dans~$\E{n}{A'}$ sur un voisinage~$U$ de~$x$ dans~$\E{n}{A}$. Consid\'erons le diagramme commutatif suivant :
$$\xymatrix{
U' \ar^{\beta}_{\sim}[r] \ar^{\pi'}[d]& U\ar^{\pi}[d]\\
\Ms(A') \ar^{\gamma}[r]& \Ms(A)
}.$$
Soit $V$ un voisinage du point~$x$ dans~$X$. Nous pouvons supposer qu'il est contenu dans~$U$. Nous avons alors 
$$\pi(V) = \gamma(\pi'(\beta^{-1}(V))).$$
Le morphisme~$\beta^{-1}$ \'etant un hom\'eomorphisme, il envoie le voisinage~$V$ du point~$x$ sur un voisinage~$\beta^{-1}(V)$ du point~$x'$. Puisque le point~$x'$ est rationnel dans sa fibre, le corollaire~\ref{ouvertdep} nous assure que la partie~$\pi'(\beta^{-1}(V))$ est un voisinage du point~$\pi'(x')$ dans~$\Ms(A')$. D'apr\`es le th\'eor\`eme~\ref{chgtbase}, le morphisme~$\gamma$ est ouvert. On en d\'eduit que la partie~$\pi(V) = \gamma(\pi'(\beta^{-1}(V)))$ est un voisinage du point~$b = \gamma(\pi'(\beta^{-1}(x)))$ dans~$\Ms(A)$.
\end{proof}

\subsection{\'Etude alg\'ebrique locale}

Nous disposons dor\'enavant de la connexit\'e locale au voisinage des points rigides des fibres et pouvons donc appliquer le r\'esultat d'isomorphie locale que nous avons d\'emontr\'e dans la proposition~\ref{isorig}. Cela va nous permettre d'\'etudier les anneaux locaux en ces points. Commen\c{c}ons par le cas des fibres extr\^emes.


\begin{thm}\label{noetherienrigext}\index{Anneau local en un point!rigide!d'une fibre extreme de An@d'une fibre extr\^eme de $\E{n}{A}$}
Soient~$\m$ un \'el\'ement de~$\Sigma_{f}$ et~$x$ un point rigide de la fibre extr\^eme~$\tilde{X}_{\m}$. Alors, l'anneau~$\Os_{X,x}$ est un anneau local noeth\'erien, r\'egulier, de dimension~$n+1$. Son corps r\'esiduel~$\kappa(x)$ est complet, et donc isomorphe \`a~$\Hs(x)$.
\end{thm}
\begin{proof}
D'apr\`es la proposition \ref{cparig}, le point~$x$ poss\`ede un syst\`eme fondamental de voisinages connexes par arcs. Nous pouvons donc utiliser la proposition \ref{isorig} et nous ramener au cas d'un point $x$ rationnel. Il existe alors des \'el\'ements~\mbox{$\alpha_{1},\ldots,\alpha_{n}$} de $k_{\m}$ tels que le point $x$ soit l'unique point de la fibre~$\tilde{X_{\m}}$ v\'erifiant 
$$(T_{1}-\alpha_{1})(x)=\cdots=(T_{n}-\alpha_{n})(x)=0.$$
Bien entendu, quel que soit $i\in\cn{1}{n}$, l'\'el\'ement $\alpha_{i}$ de $k_{\m}$ se rel\`eve en un \'el\'ement de $\hat{A}_{\m}$ et donc en un \'el\'ement de $\Os_{B,b}$. Nous pouvons donc appliquer le th\'eor\`eme~\ref{anneaulocal}. Il nous assure qu'il existe un isomorphisme 
$$\Os_{X,x} \simeq \varinjlim_{V,\bt}\, \Bs(V)\of{\la}{|\bT|\le \bt}{\ra},$$
o\`u $V$ d\'ecrit l'ensemble des voisinages compacts du point $\tilde{a}_{\m}$ de~$B$ et $\bt$ l'ensemble~$(\R_{+}^*)^n$. Il ne nous reste plus, \`a pr\'esent, qu'\`a appliquer les th\'eor\`emes~\ref{noetherienavd} et~\ref{regulier} pour d\'emontrer la premi\`ere partie du th\'eor\`eme. La seconde d\'ecoule du lemme~\ref{idealmax} et de la description des corps r\'esiduels aux points de l'espace~$B$.

\end{proof}

\begin{rem}
Signalons que, dans ce cas particulier, nous pouvons conclure sans l'aide des th\'eor\`emes~\ref{noetherienavd} et~\ref{regulier}. En effet, nous connaissons un syst\`eme fondamental de voisinages compacts explicite du point $\tilde{a}_{\m}$ de $B$ : il s'agit de l'ensemble des intervalles $\of{[}{a_{\m}^{\eps},\tilde{a}_{\m}}{]}$, avec~$\eps\in\of{]}{0,+\infty}{[}$. Quel que soit $\eps\in\of{]}{0,+\infty}{[}$, l'alg\`ebre $\Bs(\of{[}{a_{\m}^{\eps},\tilde{a}_{\m}}{]})$ n'est autre que l'alg\`ebre $\hat{A}_{\m}$. Elle est munie de la norme $\|.\|_{\of{[}{a_{\m}^{\eps},\tilde{a}_{\m}}{]}}=|.|_{\m}^\eps$. On en d\'eduit imm\'ediatement un isomorphisme 
$$\Os_{X,x} \simeq \hat{A}_{\m}[\![\bT]\!].$$
\end{rem}

Le cas des points rigides des fibres internes et centrale se traite de mani\`ere identique. Il suffit de remplacer, dans la d\'emonstration ci-dessus, le th\'eor\`eme~\ref{noetherienavd} par le th\'eor\`eme~\ref{noetheriencorps}. Nous obtenons le r\'esultat suivant.

\begin{thm}\label{noetherienrigcentralinterne}\index{Anneau local en un point!rigide!d'une fibre interne de An@d'une fibre interne de $\E{n}{A}$}\index{Anneau local en un point!rigide!de la fibre centrale de An@de la fibre centrale de $\E{n}{A}$}
Soit~$x$ un point rigide d'une fibre interne ou centrale de l'espace~$X$. Alors, l'anneau $\Os_{X,x}$ est un anneau local noeth\'erien, r\'egulier, de dimension $n$. Son corps r\'esiduel~$\kappa(x)$ est complet, et donc isomorphe \`a~$\Hs(x)$.
\end{thm}

\index{Point!rigide!d'une fibre de An@d'une fibre de \E{n}{A}|)}

\section{Fibres internes}\label{fibresinternes}

\index{Fibre!interne|(}\index{Point!interne|(}

Nous \'etudions ici les points des fibres internes de l'espace~$X$, en utilisant les propri\'et\'es du flot. Nous retrouverons, en particulier, par ce biais, les r\'esultats sur les points rigides des fibres internes obtenus \`a la section pr\'ec\'edente.

Nous reprenons, ici, les notations du paragraphe \ref{parflot}, consacr\'e au flot. Soit \mbox{$\m\in\Sigma_{f}$}. Rappelons que la fibre~$X_{a_{\m}}$ est l'espace affine de dimension~$n$ au-dessus du corps $\hat{K}_{\m}$. D'apr\`es le lemme \ref{topobranche}, l'application
$$\psi_{\m} : {\renewcommand{\arraystretch}{1.2}\begin{array}{ccc}
\of{]}{0,+\infty}{[} & \to & B'_{\m}\\
\eps & \mapsto & a_{\m}^\eps
\end{array}}$$
est un hom\'eomorphisme. L'intervalle de d\'efinition de la trajectoire de tout point de la fibre~$X_{a_{\m}}$ est $\of{]}{0,+\infty}{[}$. Par cons\'equent, nous disposons d'une application
$$\varphi_{\m} : {\renewcommand{\arraystretch}{1.2}\begin{array}{ccc}
X_{a_{\m}} \times \of{]}{0,+\infty}{[} & \to & X'_{\m}\\
(x,\eps) & \mapsto & x^\eps 
\end{array}}.$$
Notons $p_{2} : X_{a_{\m}} \times \of{]}{0,+\infty}{[} \to \of{]}{0,+\infty}{[}$ l'application de projection sur le second facteur. Ces diff\'erentes applications s'inscrivent dans le diagramme commutatif qui suit :
$$\xymatrix{
X_{a_{\m}} \times \of{]}{0,+\infty}{[} \ar[r]^{\qquad \ \ \varphi_{\m}} \ar[d]_{p_{2}} & X'_{\m} \ar[d]_{\pi}\\
\of{]}{0,+\infty}{[} \ar[r]^{\quad \psi_{\m}} & B'_{\m}
}.$$

\begin{prop}\label{produitum}\index{Fibre!interne!structure de produit}\index{Voisinages d'un point!interne!de An@de $\E{n}{A}$}
L'application~$\varphi_{\m}$ est un hom\'eomorphisme.
\end{prop}
\begin{proof}
Pour $x\in X'_{\m}$, notons 
$$\lambda(x) = \frac{\log(|\pi(x)|)}{\log(|\pi_{\m}|_{\m})}.$$
L'application $\lambda$ est continue et, quel que soit $x\in X'_{\m}$, nous avons
$$\pi(x) = a_{\m}^{\lambda(x)}.$$

Il est clair que l'application $\varphi_{\m}$ est bijective d'inverse 
$$\varphi_{\m}^{-1} : 
{\renewcommand{\arraystretch}{1.3}\begin{array}{ccc}
X'_{\m} & \to & X_{a_{\m}} \times \of{]}{0,+\infty}{[}\\ 
x & \mapsto & \left(x^{1/\lambda(x)},\pi(x)\right)
\end{array}}.$$

Montrons que l'application $\varphi_{\m}$ est un hom\'eomorphisme. Rappelons que la topologie de $X'_{\m}$ est, par d\'efinition, la topologie la plus grossi\`ere qui rend continues les applications de la forme
$$|P| : \begin{array}{ccc}
X'_{\m} & \to & \R_{+}\\
x & \mapsto & |P(x)|
\end{array},$$
avec $P\in A[\bT]$. Pour montrer que l'application $\varphi_{\m}$ est continue, il suffit donc de montrer que, quel que soit $P\in A[\bT]$, l'application 
$$|P| \circ \varphi_{\m} : {\renewcommand{\arraystretch}{1.2}\begin{array}{ccc}
X_{a_{\m}} \times \of{]}{0,+\infty}{[} & \to & \R_{+}\\
(x,\eps) & \mapsto & |P(x^\eps)| = |P(x)|^\eps
\end{array}}$$
est continue. Cette propri\'et\'e est bien v\'erifi\'ee.

De m\^eme, la topologie sur $X_{a_{\m}} \times \of{]}{0,+\infty}{[}$ est, par d\'efinition, la topologie la plus grossi\`ere qui rend continues la projection $p_{1}$ vers $X_{a_{\m}}$ et la projection~$p_{2}$ vers $\of{]}{0,+\infty}{[}$. Il nous suffit donc de montrer la compos\'ee de $\varphi_{\m}^{-1}$ avec chacune de ces deux applications est continue. C'est imm\'ediat pour l'application 
$$p_{2} \circ \varphi_{\m}^{-1} = \psi_{\m}^{-1} \circ \pi.$$ 
Consid\'erons donc l'application 
$$p_{1} \circ \varphi_{\m}^{-1} : 
{\renewcommand{\arraystretch}{1.2}\begin{array}{ccc}
X'_{\m} & \to & X_{a_{\m}}\\ 
x & \mapsto & x^{1/\lambda(x)}
\end{array}}.$$
Pour montrer que cette application est continue, il suffit de montrer que, quel que soit $P\in \hat{K}_{\m}[\bT]$, l'application 
$$|P| \circ p_{1} \circ \varphi_{\m}^{-1} : 
{\renewcommand{\arraystretch}{1.2}\begin{array}{ccc}
X'_{\m} & \to & \R_{+}\\ 
x & \mapsto & \left|P(x^{1/\lambda(x)})\right| = |P(x)|^{1/\lambda(x)}
\end{array}}$$
est continue. Puisqu'une fonction qui est limite uniforme, sur tout compact, d'applications continues est encore continue, il suffit de montrer que les applications de la forme $|P| \circ p_{1} \circ \varphi_{\m}^{-1}$, avec $P\in A[\bT]$ sont continues. Cela d\'ecoule alors directement de la d\'efinition de la topologie de $X'_{\m}$ et de la continuit\'e de la projection.
\end{proof}

Un r\'esultat similaire est valable pour la partie archim\'edienne de l'espace~$X$. La preuve en est compl\`etement analogue et nous ne la d\'etaillerons pas. Soit \mbox{$\sigma\in\Sigma_{\infty}$}. Rappelons que la fibre~$X_{a_{\sigma}}$ est isomorphe \`a l'espace~$\C^n$ si~$\hat{K}_{\sigma}=\C$ et \`a son quotient par la conjugaison complexe si~$\hat{K}_{\sigma}=\R$. D'apr\`es le lemme \ref{topobranche}, l'application
$$\psi_{\sigma} : {\renewcommand{\arraystretch}{1.2}\begin{array}{ccc}
\of{]}{0,1}{]} & \to & B'_{\sigma}\\
\eps & \mapsto & a_{\sigma}^\eps
\end{array}}$$
est un hom\'eomorphisme. L'intervalle de d\'efinition de la trajectoire de tout point de la fibre~$X_{a_{\sigma}}$ est $\of{]}{0,1}{]}$. Par cons\'equent, nous disposons d'une application
$$\varphi_{\sigma} : {\renewcommand{\arraystretch}{1.2}\begin{array}{ccc}
X_{a_{\sigma}} \times \of{]}{0,1}{]} & \to & X'_{\sigma}\\
(x,\eps) & \mapsto & x^\eps 
\end{array}}.$$
Notons $p_{2} : X_{a_{\sigma}} \times \of{]}{0,1}{]} \to \of{]}{0,1}{]}$ l'application de projection sur le second facteur. Ces diff\'erentes applications s'inscrivent dans le diagramme commutatif qui suit :
$$\xymatrix{
X_{a_{\sigma}} \times \of{]}{0,1}{]} \ar[r]^{\qquad  \varphi_{\sigma}} \ar[d]_{p_{2}} & X'_{\m} \ar[d]_{\pi}\\
\of{]}{0,1}{]} \ar[r]^{\ \psi_{\sigma}} & B'_{\sigma}
}.$$

\begin{prop}\label{produitarc}\index{Fibre!interne!structure de produit}\index{Voisinages d'un point!interne!de An@de $\E{n}{A}$}
L'application~$\varphi_{\sigma}$ est un hom\'eomorphisme.
\end{prop}

Nous d\'eduisons de ces r\'esultats deux corollaires topologiques.

\begin{cor}\label{ouvertinterne}\index{Ouverture au voisinage d'un point!interne}
Le morphisme $\pi$ est ouvert en tout point d'une fibre interne de $X$.
\end{cor}

\begin{cor}\label{cpainterne}\index{Connexite par arcs au voisinage d'un point@Connexit\'e par arcs au voisinage d'un point!interne}
Tout point d'une fibre interne de l'espace $X$ poss\`ede un syst\`eme fondamental de voisinages connexes par arcs. 
\end{cor}

\begin{cor}\label{flotinterne}\index{Voisinages flottants!pour un point interne}
Tout point interne de $X$ poss\`ede des voisinages flottants, au sens de la d\'efinition \ref{voisflot}.
\end{cor}
\begin{proof}
Soient $\sigma\in\Sigma$ et $x$ un point de $X'_{\sigma}$. Reprenons les notations du paragraphe~\ref{parflot}. Nous avons $D=X'_{\sigma}$ et la structure de produit dont les propositions pr\'ec\'edentes d\'emontrent l'existence assurent que le flot est une application ouverte.
\end{proof}


\begin{prop}\label{isointerne}\index{Fibre!interne!structure de produit}
Soit $b$ un point interne de $B$. Alors l'inclusion 
$$j_{b} : X_{b} \hookrightarrow X$$ 
de la fibre dans l'espace total induit un isomorphisme entre les espaces annel\'es
$$(X_{b},j_{b}^{-1}\Os_{X}) \xrightarrow[]{\sim} (X_{b},\Os_{X_{b}}).$$
\end{prop}
\begin{proof}
Signalons tout d'abord qu'en d\'epit de ce que les notations utilis\'ees peuvent laisser penser les espaces topologiques sous-jacents sont, \emph{a priori}, diff\'erents. En effet, sur l'un ce sont les valeurs absolues de polyn\^omes \`a coefficients dans $A$ qui doivent \^etre continues, et, sur l'autre, ce sont celles des polyn\^omes \`a coefficients dans $\hat{K}_{\sigma}$. Cependant, la continuit\'e \'etant une propri\'et\'e stable par limite uniforme sur tout compact, les topologies sont bien identiques. L'application identit\'e d\'efinit donc bien un hom\'eomorphisme.

Int\'eressons-nous, \`a pr\'esent, aux faisceaux structuraux. Soit~ $x\in X_{b}$. Il nous suffit de montrer que le morphisme naturel
$$\Os_{X,x} \to \Os_{X_{b},x}$$
est un isomorphisme. Commen\c{c}ons par montrer qu'il est injectif. Soit~$f$ un \'el\'ement de~$\Os_{X,x}$ nul dans~$\Os_{X_{b},x}$. Il existe un voisinage~$V$ de~$x$ dans~$X_{b}$ sur lequel la fonction~$f$ est nulle. D'apr\`es les propositions~\ref{produitum} et~\ref{produitarc}, la fonction~$f$ est d\'efinie sur un voisinage~$U$ de~$x$ dans~$X$ de la forme
$$U = \{y^\eps,\, y\in W,\, \alpha < \eps < \beta\},$$  
o\`u~$W$ est un voisinage de~$x$ dans~$V$, $\alpha$ un \'el\'ement de~$\of{]}{0,1}{[}$ et~$\beta$ un \'el\'ement de~$\of{]}{1,+\infty}{[}$. Soit~$z\in U$. Il existe un \'el\'ement~$y$ de~$W$ et un nombre r\'eel~$\eps\in\of{]}{\alpha,\beta}{[}$ tels que~$z=y^\eps$. D'apr\`es le corollaire~\ref{flotinterne}, le point~$y$ poss\`ede des voisinages flottants. D'apr\`es la proposition~\ref{flot}, nous avons donc
$$|f(z)| = |f(y)|^\eps = 0.$$
On en d\'eduit que la fonction~$f$ est nulle sur~$U$ et donc dans l'anneau local~$\Os_{X,x}$.

Montrons, \`a pr\'esent, que le morphisme entre les anneaux locaux est surjectif. Soit~$f\in\Os_{X_{b},s}$. Il existe un voisinage compact~$V$ de~$x$ dans~$X_{b}$ et une suite~$(R_{k})_{k\in\N}$ d'\'el\'ements de~$\hat{K}_{\sigma}(\bT)$, sans p\^oles sur~$V$, qui converge vers la fonction~$f$ sur~$V$. Soit $k\in\N$. Il existe un \'el\'ement $S_{k}$ de~Frac($A[\bT]$) sans p\^oles sur~$V$ qui v\'erifie
$$\|S_{k}-R_{k}\|_{V} \le 2^{-k}.$$
Consid\'erons le voisinage~$U$ du point~$x$ de~$X$ d\'efini par
$$U = \left\{y^\eps,\, y\in V,\, \frac{1}{2} \le \eps \le \frac{3}{2}\right\}.$$
Quel que soit $k\in\N$, la fonction~$S_{k}$ n'a pas de p\^oles sur la partie compacte~$U$. 

Soit $\eta>0$. Il existe un entier $p\in\N$ tel que, quels que soient $k,l\ge p$, nous ayons 
$$\|R_{k}-R_{l}\|_{V} \le \eta.$$
Quitte \`a augmenter~$p$, nous pouvons supposer que~$2^{-p}\le \eta$. Soit~$z\in U$. Il existe un \'el\'ement~$y$ de~$V$ et un nombre r\'eel~$\eps\in\of{[}{1/2,3/2}{]}$ tels que~$z=y^\eps$. Quel que soient $k,l\ge p$, nous avons alors
$${\renewcommand{\arraystretch}{1.3}\begin{array}{rcl}
|(S_{k}-S_{l})(y)| &=& |(S_{k}-S_{l})|^\eps\\
&\le& \left(\|R_{k}-R_{l}\|_{V} + 2^{-k} + 2^{-l}\right)^\eps\\
&\le& (3\eta)^\eps\\
&\le& \max\left((3\eta)^{1/2},(3\eta)^{3/2}\right).
\end{array}}$$
Par cons\'equent, la suite $(S_{k})_{k\in\N}$ converge uniform\'ement sur~$U$ vers un \'el\'ement~$g$ de~$\Bs(U)$ et donc de~$\Os_{X,x}$. L'image de cet \'el\'ement dans l'anneau local~$\Os_{X_{b},x}$ n'est autre que l'\'el\'ement~$f$.
\end{proof}

\begin{thm}\label{noetherieninterne}\index{Anneau local en un point!interne!de An@de $\E{n}{A}$}
Soit $x$ un point interne de $X$. L'anneau local $\Os_{X,x}$ est hens\'elien, noeth\'erien, r\'egulier, excellent et de dimension inf\'erieure \`a~$n$. Le corps~$\kappa(x)$ est hens\'elien.
\end{thm}
\begin{proof}
La proposition qui pr\'ec\`ede nous permet de nous ramener au cas o\`u l'espace de base est le spectre d'un corps, cadre dans lequel ces r\'esultats sont connus. 
\end{proof}

\bigskip

Pour finir, d\'emontrons des r\'esultats indiquant que l'on peut contr\^oler le bord de Shilov des voisinages de certains points. Commen\c{c}ons par rappeler quelques propri\'et\'es du le bord de Shilov dans le cadre des espaces analytiques sur un corps ultram\'etrique complet. 

\begin{prop}[V.~Berkovich]\label{ShilovBerko}\index{Bord analytique!de MA@de $\Ms(\As)$}
Soient~$(k,|.|)$ un corps ultram\'etrique complet et~$(\As,\|.\|)$ une alg\`ebre $k$-affino\"ide. L'anneau de Banach $(\As,\|.\|)$ poss\`ede un bord de Shilov~$\Gamma$. C'est un ensemble fini.

Soient~$(k,|.|)$ un corps ultram\'etrique complet et~$m$ un entier positif. Le bord de Shilov de tout domaine affino\"ide de~$\E{m}{k}$ est contenu dans son int\'erieur topologique.
\end{prop}
\begin{proof}
La premi\`ere partie de la proposition provient du corollaire 2.4.5 de~\cite{rouge}. La seconde provient du corollaire 2.5.13 (ii) et de la proposition 2.5.20 (que l'on applique, par exemple, en prenant comme espace affino\"ide~$X$ un disque de rayon assez grand et comme domaine affino\"ide~$V$ le domaine affino\"ide en question). 
\end{proof}

Apportons une pr\'ecision gr\^ace \`a la proposition suivante.

\begin{prop}\label{ShilovAntoine}
Soient~$(k,|.|)$ un corps ultram\'etrique complet et~$m$ un entier positif. Soit~$V$ un domaine strictement affino\"ide irr\'eductible de l'espace affine~$Y=\E{m}{k}$. Notons~$\Gamma$ son bord de Shilov. En tout point~$\gamma$ de~$\Gamma$, le corps r\'esiduel~$\widetilde{\Hs(\gamma)}$ du corps r\'esiduel compl\'et\'e~$\Hs(\gamma)$ est de degr\'e de transcendance~$m$. En particulier, en tout point~$\gamma$ de~$\Gamma$, l'anneau local~$\Os_{Y,\gamma}$ est un corps. 
\end{prop}
\begin{proof}
La premi\`ere partie de la proposition d\'ecoule de la proposition 2.4.4. (ii) de~\cite{rouge}. Puisque le corps $\widetilde{\Hs(\gamma)}$ a pour degr\'e de transcendance~$m$, le point~$\gamma$ ne peut se trouver localement sur aucun ferm\'e de Zariski de dimension strictement inf\'erieure \`a~$m$. L'espace~$Y$ \'etant r\'eduit, on en d\'eduit que l'anneau local~$\Os_{Y,\gamma}$ est un corps.
\end{proof}

\begin{rem}
Ainsi que nous l'a fait remarquer A.~Ducros, le r\'esultat pr\'ec\'edent vaut pour tout domaine affino\"ide de tout espace de Berkovich bon et r\'eduit.
\end{rem}

\begin{cor}\label{Shilovrationnelinterne}
Soient~$\m\in\Sigma_{f}$ et~$b\in B'_{\m}$. Soit~$V$ une partie compacte et spectralement convexe de~$X$ contenue dans la fibre~$X_{b}$ qui est un domaine strictement affino\"ide irr\'eductible de cette fibre, vue comme espace analytique sur~$\Hs(b)$. Alors la partie~$V$ poss\`ede un bord analytique fini et alg\'ebriquement trivial.
\end{cor}
\begin{proof}
Le bord de Shilov~$\Gamma$ de l'affino\"ide~$V$ est un bord analytique de~$V$ dans~$X$. En effet, l'alg\`ebre affino\"ide de~$V$ contient~$\Bs(V)$. Il suffit ensuite de combiner les deux propositions pr\'ec\'edentes avec la proposition \ref{isointerne}.
\end{proof}

\begin{lem}
Soit $(k,|.|)$ un corps ultram\'etrique complet dont la valuation n'est pas triviale. Soient~$m\in\N$, $y$ un point de l'espace~$Y=\E{m}{k}$ et~$U$ un voisinage de ce point. Il existe $r\in\N$ et $P_{1},\ldots,P_{r},Q_{1},\ldots,Q_{r}\in k[\bT]$ tels que la partie de~$Y$ d\'efinie par
$$\bigcap_{1\le i\le r} \left\{z\in Y\, \big|\, |P_{i}(z)|\le |Q_{i}(z)|\right\}$$
soit un voisinage strictement affino\"ide irr\'eductible du point~$y$ dans~$U$.
\end{lem}
\begin{proof}
Soit~$\alpha\in k$ tel que~$|\alpha|\in\of{]}{0,1}{[}$. L'ensemble
$$E=\left\{|\alpha|^{\frac{p}{q}},\, p\in\Z, q\in\N^*\right\}$$
est alors dense dans~$\R_{+}$. Par d\'efinition de la topologie, il existe $r,s\in\N$, $G_{1},\ldots,G_{r},H_{1},\ldots,H_{s}\in k[\bT]$, $u_{1},\ldots,u_{r},v_{1},\ldots,v_{s}$ tels que la partie
$$V = \bigcap_{1\le i\le r} \left\{z\in Y\, \big|\, |G_{i}(z)|\le u_{i}\right\} \cap  \bigcap_{1\le i\le s} \left\{z\in Y\, \big|\, |H_{i}(z)|\ge v_{i}\right\}$$
soit un voisinage compact du point~$y$ dans~$U$. 

Soit $i\in\cn{1}{r}$. Il existe~$p\in\Z$ et~$q\in\N^*$ tels que $u_{i}=|\alpha|^{p/q}$. Remarquons que
$$\left\{z\in Y\, \big|\, |G_{i}(z)|\le u_{i}\right\} = \left\{z\in Y\, \big|\, |(\alpha^{-p}\, G_{i}^q)(z)|\le 1\right\}.$$
Par cons\'equent, nous pouvons supposer que, quel que soit $i\in\cn{1}{r}$, nous avons $u_{i}=1$. De m\^eme, nous pouvons supposer que, quel que soit $i\in\cn{1}{s}$, nous avons $v_{i}=1$. 

Consid\'erons un disque compact~$D$ de~$Y$ qui contient le compact~$V$. Notons~$\As_{D}$ l'alg\`ebre $k$-affino\"ide associ\'ee. La partie~$V$ est un domaine rationnel de~$D$. Notons~$\As_{V}$ l'alg\`ebre $k$-affino\"ide associ\'ee. Consid\'erons la composante connexe~$W$ de~$V$ qui contient le point~$x$. Il existe un \'el\'ement~$f$ de~$\As_{V}$ qui est nul sur~$W$ et vaut identiquement~$1$ sur la r\'eunion~$R$ des autres composantes connexes de~$V$. Puisque~$V$ est un domaine rationnel de~$D$, il existe des \'el\'ements~$g$ et~$h$ de~$\As_{D}$ tels que la fonction~$h$ ne s'annule pas sur~$V$ et tels que l'on ait
$$\left\{z\in V\, \left|\, \left|\frac{g}{h}(z)\right| < |\alpha| \right.\right\} =W \textrm{ et } \left\{z\in V\, \left|\, \left|\frac{g}{h}(z)\right| > |\alpha| \right.\right\} =R.$$
Puisque~$D$ est un disque, les polyn\^omes sont denses dans~$\As_{D}$. Il existe donc des \'el\'ements~$G$ et~$H$ de~$k[\bT]$ tels que la fonction~$H$ ne s'annule pas sur~$V$ et tels que l'on ait
$$\left\{z\in V\, \left|\, \left|\frac{G}{H}(z)\right| \le |\alpha| \right.\right\} =W.$$
Pour conclure, il suffit d'\'ecrire le voisinage compact et connexe~$W$ du point~$y$ dans~$U$ sous la forme
$$W = V \cap \left\{z\in Y\, \big|\, |G(z)|\le |\alpha H(z)|\right\}.$$
C'est bien un domaine strictement affino\"ide de~$Y$. Il est irr\'eductible car il est connexe et que l'espace analytique~$Y$ est normal.
\end{proof}


%

\begin{prop}\label{shilovfiniinterne}\index{Bord analytique!au voisinage d'un point interne!de An@de $\E{n}{A}$}
Soit~$\m\in\Sigma_{f}$. Tout point de~$X'_{\m}$ poss\`ede un syst\`eme fondamental de voisinages compacts, connexes et spectralement convexes qui poss\`edent un bord analytique fini et alg\'ebriquement trivial.
\end{prop}
\begin{proof}
Soient~$b$ un point de~$B'_{\m}$ et~$x$ un point de la fibre~$X_{b}$. Soit~$U$ un voisinage du point~$x$ dans~$X$. D'apr\`es le lemme pr\'ec\'edent, il existe $r\in\N$ et $P_{1},\ldots,P_{r},Q_{1},\ldots,Q_{r}\in \hat{K}_{\m}[\bT]$ tels que la partie de~$Y$ d\'efinie par
$$V_{1} = \bigcap_{1\le i\le r} \left\{z\in Y\, \big|\, |P_{i}(z)|\le |Q_{i}(z)|\right\}$$
soit un voisinage strictement affino\"ide irr\'eductible du point~$y$ contenu dans l'int\'erieur de~$U\cap X_{b}$ dans~$X_{b}$.

D'apr\`es la proposition \ref{produitum}, il existe $\alpha,\beta\in I_{x}$ v\'erifiant $0<\alpha<1<\beta$ et tels que la partie
$$V = \{y^\eps,\, y\in V_{1}, \eps\in\of{[}{\alpha,\beta}{]}\}$$
soit un voisinage compact et connexe du point~$x$ dans~$U$. Remarquons que
$$V =  \bigcap_{1\le i\le r} \left\{z\in \pi^{-1}(\of{[}{b^\alpha,b^\beta}{]})\, \big|\, |P_{i}(z)|\le |Q_{i}(z)|\right\}.$$
On d\'eduit alors du th\'eor\`eme \ref{compactrationnel} et des propositions \ref{prorat} et \ref{stabilitespconvexe} que le compact~$V$ est spectralement convexe.

Notons~$\Gamma_{1}$ le bord analytique fini et alg\'ebriquement trivial du compact~$V_{1}$ dont l'existence nous est assur\'ee par le corollaire \ref{Shilovrationnelinterne}. Posons
$$\Gamma = \{x^\alpha,\, x\in\Gamma_{1}\} \cup \{x^\beta,\, x\in\Gamma_{1}\}.$$
On d\'eduit du corollaire \ref{flotinterne} et de la proposition \ref{flot} que, pour tout \'el\'ement~$f$ de~$\Os(V)$, nous avons
$$\|f\|_{V} = \|f\|_{\Gamma}$$
et que la partie~$\Gamma$ est encore finie et alg\'ebriquement triviale. La partie~$\Gamma$ est donc un bord analytique du compact~$V$ qui satisfait les propri\'et\'es voulues.
\end{proof}

\begin{prop}\label{Shilovrigum}\index{Bord analytique!au voisinage d'un point rigide!de An@de $\E{n}{A}$}
Soient~$\sigma$ un \'el\'ement de~$\Sigma_{f}$ et~$b$ un point de $B_{\sigma}\setminus\{a_{0}\}$. Tout point rigide de la fibre~$X_{b}$ poss\`ede un syst\`eme fondamental de voisinages compacts, connexes et spectralement convexes qui poss\`edent un bord de Shilov fini et alg\'ebriquement trivial.
\end{prop}
\begin{proof}
Soit~$x$ un point rigide de la fibre~$X_{b}$. D'apr\`es la proposition \ref{isorig} et le lemme \ref{translation}, nous pouvons supposer que le point~$x$ est le point~$0$ de la fibre~$X_{b}$. La proposition~\ref{voisdep} assure alors que ce point poss\`ede un syt\`eme fondamental de voisinages qui sont des disques compacts au-dessus de parties compactes et connexes de~$B_{\sigma}\setminus\{a_{0}\}$. La discussion men\'ee au num\'ero~\ref{borddeShilovbase} et la proposition~\ref{Shilovcouronneum} montrent qu'une telle partie poss\`ede un bord de Shilov et en fournissent une description explicite. C'est en particulier un ensemble fini compos\'e de points internes. On d\'eduit de la proposition~\ref{isointerne} qu'il est alg\'ebriquement trivial. 
\end{proof}

\begin{prop}\label{Shilov3depum}\index{Bord analytique!au voisinage d'un point deploye@au voisinage d'un point d\'eploy\'e}
Soient~$\sigma$ un \'el\'ement de~$\Sigma_{f}$ et~$b$ un point de $B_{\sigma}\setminus\{a_{0}\}$. Tout point d\'eploy\'e de la fibre~$X_{b}$ poss\`ede un syst\`eme fondamental de voisinages compacts, connexes et spectralement convexes qui poss\`edent un bord de Shilov fini et alg\'ebriquement trivial.
\end{prop}
\begin{proof}
La proposition~\ref{voisdep} assure qu'un point d\'eploy\'e poss\`ede un syt\`eme fondamental de voisinages qui sont des couronnes compactes au-dessus de parties compactes et connexes de~$B_{\sigma}\setminus\{a_{0}\}$. On conclut alors comme dans la preuve pr\'ec\'edente.
\end{proof}

\index{Fibre!interne|)}\index{Point!interne|)}

\section{Dimension topologique}\label{dimensiontopologique}

Nous consacrons cette partie \`a l'\'etude de la dimension topologique de l'espace affine analytique~$X=\E{n}{A}$ d\'efini au-dessus de l'anneau d'entiers de corps de nombres~$A$. La notion de dimension topologique n'est agr\'eable que lorsque l'espace consid\'er\'e est m\'etrisable. Dans ce cas, la dimension de recouvrement (\emph{cf.} \cite{dimension}, d\'efinition I.4) et la dimension inductive forte (\emph{cf.} \cite{dimension}, d\'efinition I.5) co\"{\i}ncident (\emph{cf.} \cite{dimension}, th\'eor\`eme II.7). Commen\c{c}ons par v\'erifier que nous nous trouvons bien dans cette situation.

\begin{thm}\label{metrisable}\index{Espace affine analytique!sur A@sur $A$!metrisabilit\'e@m\'etrisabilit\'e}\index{Metrisabilite@M\'etrisabilit\'e|see{Espace affine analytique}}
L'espace analytique~$X=\E{n}{A}$ est m\'etrisable.
\end{thm}
\begin{proof}
Soient~$x$ un point de~$X$ et~$U$ un voisinage du point~$x$ dans~$X$. Par d\'efinition de la topologie, il existe $r\in\N$, $P_{1},\ldots,P_{r}\in A[T_{1},\ldots,T_{n}]$ et $u_{1},\ldots,u_{r},v_{1},\ldots,v_{r}\in\R$
tels que la partie
$$V = \bigcap_{1\le i\le r} \{y\in X\, |\, u_{i} < |P_{i}(y)| < v_{i}\}$$
soit un voisinage du point~$x$ contenu dans~$U$. Nous pouvons supposer que les nombres $u_{1},\ldots,u_{r},v_{1},\ldots,v_{r}$ sont rationnels. Puisque l'ensemble~$A$ est d\'e\-nom-bra\-ble, l'ensemble des voisinages de la forme pr\'ec\'edente est alors d\'enombrable. On en d\'eduit que l'espace~$X$ est s\'eparable.

D'apr\`es le th\'eor\`eme \ref{topoAn}, l'espace~$X$ est localement compact et donc r\'egulier. Le th\'eor\`eme d'Urysohn (\emph{cf.} \cite{dimension}, corollaire du th\'eor\`eme I.3) assure alors qu'il est m\'etrisable.
\end{proof}

Nous pouvons, \`a pr\'esent, calculer la dimension topologique de l'espace~$\E{n}{A}$. Commen\c{c}ons par l'espace de base~$B=\Ms(A).$

\begin{prop}\index{Spectre analytique!de A@de $A$!dimension topologique}\index{Dimension topologique!de MA@de $\Ms(A)$}
La dimension topologique de l'espace~$B$ est \'egale \`a~$1$.
\end{prop}
\begin{proof}
Soit~$\sigma\in\Sigma$. La branche~$\sigma$-adique~$B_{\sigma}$ est hom\'eomorphe au segment~$\of{[}{0,1}{]}$. Elle est donc de dimension~$1$. D'apr\`es~\cite{dimension}, th\'eor\`eme~II.3, nous avons donc
$$\dim(B)\ge 1.$$
En outre, nous avons
$$B = \bigcup_{\sigma\in\Sigma} B_{\sigma}$$
et ce recouvrement est d\'enombrable. D'apr\`es~\cite{dimension}, th\'eor\`eme~II.1, nous avons donc
$$\dim(B) \le 1.$$
On en d\'eduit le r\'esultat voulu.
\end{proof}

Traitons, maintenant, le cas g\'en\'eral.

\begin{prop}\label{dim}\index{Espace affine analytique!sur A@sur $A$!dimension topologique}\index{Dimension topologique!de An@de $\E{n}{A}$}
La dimension topologique de l'espace~$\E{n}{A}$ est \'egale \`a~$2n+1$.
\end{prop}
\begin{proof}
Commen\c{c}ons par minorer la dimension. Soit~$\sigma\in\Sigma_{\infty}$. D'apr\`es la proposition~\ref{produitarc}, la partie~$X'_{\sigma}$ de~$X$ est hom\'eomorphe \`a~$X_{a_{\sigma}} \times \of{]}{0,1}{]}$. Si~$\sigma$ est un plongement r\'eel, la fibre~$X_{a_{\sigma}}$ est hom\'eomorphe au quotient de l'espace~$\C^n$ par l'action de la conjugaison complexe. Elle est donc de dimension \'egale \`a~$2n$. Si~$\sigma$ est un plongement complexe non r\'eel, la fibre~$X_{a_{\sigma}}$ est hom\'eomorphe \`a l'espace~$\C^n$ lui-m\^eme et est donc encore de dimension \'egale \`a~$2n$. Dans tous les cas, la dimension de~$X'_{\sigma}$ est \'egale \`a~$2n+1$. D'apr\`es~\cite{dimension}, th\'eor\`eme~II.3, nous avons donc
$$\dim(X)\ge 2n+1.$$

Soit~$k\in\N^*$. Consid\'erons le disque de centre~$0$ et de rayon~$k$ de~$X$ :
$$D(k) = \left\{x\in X\, \big|\, |T(x)|\le k\right\}.$$
C'est une partie compacte de~$X$. L'application de projection
$$\pi_{k} : D(k) \to B$$ 
est continue et ferm\'ee. Soit~$b$ un point de~$B$. Si la valeur absolue sur le corps r\'esiduel compl\'et\'e~$\Hs(b)$ est archim\'edienne, la dimension de la fibre~$\pi_{k}^{-1}(b)$ est \'egale \`a~$2n$. Si elle est ultram\'etrique, la fibre~$\pi_{k}^{-1}(b)$ est le disque de centre~$0$ et de rayon~$k$ de l'espace affine de Berkovich de dimension~$n$ au-dessus du corps~$\Hs(b)$. D'apr\`es~\cite{bleu}, proposition~1.2.18, sa dimension est inf\'erieure \`a~$n$. D'apr\`es~\cite{dimension}, th\'eor\`eme~III.6, nous avons
$$\dim(D(k)) \le \dim(B) + 2n \le 2n+1.$$

Bien entendu, nous avons
$$X = \bigcup_{k\in\N^*} D(k).$$
D'apr\`es~\cite{dimension}, th\'eor\`eme~II.1, nous avons donc
$$\dim(X) \le 2n+1.$$
On en d\'eduit le r\'esultat annonc\'e.
\end{proof}

\section{Prolongement analytique}\label{sectionprolan}

Int\'eressons-nous, \`a pr\'esent, au probl\`eme du prolongement analytique. Commen\c{c}ons par pr\'eciser ce que nous entendons par ce terme.

\begin{defi}\index{Prolongement analytique}
Soit~$(S,\Os_{S})$ un espace localement annel\'e. Nous dirons que {\bf le principe du prolongement analytique vaut sur l'espace~$(S,\Os_{S})$} si, pour tout point~$s$ de~$S$, le morphisme naturel
$$\Os_{S}(S) \to \Os_{S,s}$$
est injectif.

Soit~$T$ une partie de l'espace topologique~$S$. Notons $j_{T} : T \hookrightarrow S$ le morphisme d'inclusion. Nous dirons que {\bf le principe du prolongement analytique vaut sur la partie~$T$ de l'espace~$S$} s'il vaut sur l'espace $(T,j_{T}^{-1}\Os_{S})$. 
\end{defi}

Introduisons \'egalement une version locale.

\begin{defi}
Soit~$(S,\Os_{S})$ un espace localement annel\'e. Soit~$s$ un point de~$S$. Nous dirons que {\bf le principe du prolongement analytique vaut au voisinage du point~$s$} si, pour tout voisinage~$U$ du point~$s$ dans~$S$ et tout \'el\'ement~$f$ de~$\Os_{S}(U)$ dont l'image n'est pas nulle dans~$\Os_{S,s}$, il existe un voisinage~$V$ de~$s$ dans~$U$ tel que l'image de la fonction~$f$ ne soit nulle dans aucun des anneaux locaux~$\Os_{S,t}$, pour~$t$ appartenant \`a~$V$.
\end{defi}

Donnons un exemple de point v\'erifiant cette propri\'et\'e.

\begin{lem}
Soit~$(S,\Os_{S})$ un espace analytique (au sens de la d\'efinition \ref{espan}). Soit~$s$ un point de~$S$ en lequel l'anneau local est un corps. Alors le principe du prolongement analytique vaut au voisinage du point~$s$.
\end{lem}

Le lemme qui suit, de d\'emonstration imm\'ediate, relie les d\'efinitions locale et globale de prolongement analytique.

\begin{lem}\label{lemprolan}
Soit~$(S,\Os_{S})$ un espace localement annel\'e. Soit~$T$ une partie connexe de l'espace topologique~$S$. Supposons que le principe du prolongement analytique vaut au voisinage de tout point de~$T$. Alors, il vaut sur~$T$.

Soit~$s$ un point de~$S$. Supposons que le point~$S$ poss\`ede un syst\`eme fondamental de voisinages sur lesquels vaut le principe du prolongement analytique. Alors il vaut au voisinage du point~$s$.
\end{lem}

Commen\c{c}ons par nous int\'eresser au cas de l'espace de base~$B$.

\begin{prop}\label{prolanB}\index{Prolongement analytique!au voisinage d'un point de MA@au voisinage d'un point de $\Ms(A)$}
Le principe du prolongement analytique vaut au voisinage de tout point~$b$ de~$B$. En particulier, il vaut sur tout ouvert connexe de l'espace~$B$. 
\end{prop}

Consid\'erons, \`a pr\'esent, le cas de l'espace affine de dimension~$n$, $X=\E{n}{A}$. Commen\c{c}ons par nous int\'eresser aux points internes de cet espace. L'utilisation du flot permet d'obtenir facilement des r\'esultats.

\begin{prop}\label{prolinterne}\index{Prolongement analytique!au voisinage d'un point interne!de An@de $\E{n}{A}$}
Le principe du prolongement analytique vaut au voisinage de tout point interne de l'espace~$X$. En particulier, pour tout \'el\'ement~$\sigma$ de~$\Sigma$, le principe du prolongement analytique vaut sur tout ouvert connexe de l'espace~$X'_{\sigma}$. 
\end{prop}
\begin{proof}
Soient~$b$ un point interne de l'espace~$B$ et~$x$ un point de la fibre~$X_{b}$. Soient~$U$ un voisinage du point~$x$ dans~$X$ et~$f$ un \'el\'ement de~$\Os_{X}(U)$ dont l'image dans l'anneau local~$\Os_{X,x}$ n'est pas nul. La proposition~\ref{isointerne} nous assure que l'image de~$f$ dans l'anneau local~$\Os_{X_{b},x}$ diff\`ere encore de~$0$. Soit~$V_{0}$ un voisinage connexe du point~$x$ dans la fibre~$X_{b}$. C'est un espace analytique normal et connexe d\'efini sur un corps valu\'e complet. Le principe du prolongement analytique y vaut donc. Par cons\'equent, pour tout \'el\'ement~$y$ de~$V_{0}$, l'image de la fonction~$f$ dans l'anneau local~$\Os_{X_{b},y}$, et donc dans l'anneau local~$\Os_{X,y}$, diff\`ere de~$0$. Les propositions~\ref{produitum} et~\ref{produitarc} assurent que le point~$x$ poss\`ede un voisinage~$V$ dans~$U$ form\'e de trajectoires d'\'el\'ements de~$V_{0}$. Le corollaire~\ref{flotinterne} et la proposition~\ref{flot} assurent alors que, pour tout point~$y$ de~$V$, l'image de la fonction~$f$ diff\`ere de~$0$ dans l'anneau local~$\Os_{X,y}$.
\end{proof}

Nous n'irons, pour le moment, gu\`ere plus loin dans cette direction. Mentionnons cependant quelques r\'esultats partiels. 


%

\begin{lem}
Soient~$V$ une partie ouverte de l'espace~$B$ et~$Y$ une couronne ouverte au-dessus de~$V$. Soit~$f$ un \'el\'ement de~$\Os_{X}(Y)$. Notons~$C$ l'ensemble des points de~$V$ qui poss\`edent un voisinage~$W$ v\'erifiant la propri\'et\'e suivante : 
$$\forall y \in X_{W}\cap Y,\, f(y)=0.$$
La partie~$C$ est ouverte et ferm\'ee dans~$V$.
\end{lem}
\begin{proof}
Par d\'efinition, la partie~$C$ est ouverte dans~$V$. Il nous reste \`a montrer qu'elle est ferm\'ee dans~$V$.

Soit~$c$ un point de~$V\setminus C$. Soit~$y$ un point d\'eploy\'e (\emph{cf.} d\'efinition \ref{defdeploye}) contenu dans~$X_{c}\cap Y$. Supposons, par l'absurde, que l'image de~$f$ dans l'anneau local~$\Os_{X,y}$ soit nulle. D'apr\`es la proposition~\ref{voisdep}, il existe un voisinage~$W$ de~$c$ dans~$V$ tel que, pour tout point~$d$ de~$W$, la fonction~$f$ est nulle sur une partie ouverte de la fibre~$X_{d}\cap Y$. Soit~$d$ un point de~$W$. Puisque l'espace analytique~$X_{d}\cap Y$ est normal et connexe, la fonction~$f$ y est identiquement nulle. On en d\'eduit que le point~$c$ appartient \`a~$C$, ce qui contredit l'hypoth\`ese.

Nous avons donc montr\'e que l'image de~$f$ dans l'anneau local~$\Os_{X,y}$ n'est pas nulle. Supposons, tout d'abord, que le point~$c$ est un point interne ou central. La description explicite de l'anneau local~$\Os_{X,y}$ nous permet d'affirmer que le morphisme naturel
$$\Os_{X,y} \to \Os_{X_{c},y}$$
est injectif. Par cons\'equent, l'image de~$f$ dans l'anneau local~$\Os_{X_{c},y}$ n'est pas nulle. Puisque l'espace analytique~$X_{c}\cap Y$ est normal et connexe, il poss\`ede un point~$z$ en lequel nous avons~$|f(z)|> 0$. En outre, nous pouvons supposer que le point~$z$ est d\'eploy\'e, car l'ensemble de ces points est dense. D'apr\`es le corollaire~\ref{ouvertdep}, le morphisme~$\pi$ est ouvert au voisinage du point~$z$. En outre, il existe un voisinage du point~$z$ dans~$Y$ sur lequel la fonction~$f$ ne s'annule pas. On en d\'eduit que la partie~$V\setminus C$ est un voisinage du point~$c$ dans~$V$.

Supposons, \`a pr\'esent, que le point~$c$ est un point extr\^eme : il existe un \'el\'ement~$\m$ de~$\Sigma_{f}$ tel que~$c=\tilde{a}_{\m}$. Il existe alors un nombre r\'eel~$\eps>0$ tel que l'intervalle~$\of{]}{a_{\m}^\eps,\tilde{a}_{\m}}{]}$ soit contenu dans~$V$. Notons
$$U=\pi^{-1}(\of{]}{a_{\m}^\eps,\tilde{a}_{\m}}{]}) \cap Y.$$
D'apr\`es la proposition~\ref{voisdep}, il existe un point de $U$ au voisinage duquel la fonction~$f$ n'est pas nulle. Puisque l'ouvert~$U$ est connexe, la proposition~\ref{prolinterne} nous assure que, pour tout point~$z$ de~$U$, l'image de la fonction~$f$ dans l'anneau local~$\Os_{X,z}$ n'est pas nulle. On en d\'eduit que, pour tout \'el\'ement~$\delta$ de~$\of{]}{\eps,+\infty}{[}$, il existe un point de $X_{a_{\m}^\eps}\cap Y$ en lequel la fonction~$f$ n'est pas nulle. En particulier, l'intervalle $\of{]}{a_{\m}^\eps,\tilde{a}_{\m}}{]}$ est contenu dans~$V\setminus C$.  On en d\'eduit que la partie~$V\setminus C$ est un voisinage du point~$c$ dans~$V$.
\end{proof}

\begin{cor}\label{corprolan}
Soient~$V$ une partie ouverte et connexe de l'espace~$B$ et~$Y$ une couronne ouverte au-dessus de~$V$. Soit~$x$ un point de~$Y$ en lequel le morphisme~$\pi$ est ouvert. Alors le morphisme naturel
$$\Os_{X}(Y) \to \Os_{X,x}$$
est injectif.
\end{cor}
\begin{proof}
Soit~$f$ un \'el\'ement de~$\Os_{X}(Y)$ dont l'image dans l'anneau local~$\Os_{X,x}$ est nulle. Puisque le morphisme~$\pi$ est ouvert en~$x$, il existe un voisinage~$W$ de~$\pi(x)$ dans~$V$ tel que, pour tout point~$b$ de~$W$, la fonction~$f$ est nulle sur une partie ouverte de la fibre~$X_{b}\cap Y$. Soit~$b$ un point de~$W$. Puisque l'espace analytique~$X_{b}\cap Y$ est normal et connexe, la fonction~$f$ y est identiquement nulle. 

D\'efinissons la partie~$C$ de~$V$ de la m\^eme fa\c{c}on que dans le lemme qui pr\'ec\`ede. Nous venons de montrer qu'elle n'est pas vide. Puisque la partie~$V$ est suppos\'ee connexe, nous avons n\'ecessairement l'\'egalit\'e~$C=V$. En d'autres termes, la fonction~$f$ est nulle en tout point de la couronne~$Y$ et donc dans~$\Os_{X}(Y)$.
\end{proof}

\index{Espace affine analytique!sur A@sur $A$|)}

%% file: MZa0.pstex_t
\begin{picture}(0,0)%
\includegraphics{MZa0.pstex}%
\end{picture}%
\setlength{\unitlength}{3947sp}%
\begingroup\makeatletter\ifx\SetFigFont\undefined%
\gdef\SetFigFont#1#2#3#4#5{%
  \reset@font\fontsize{#1}{#2pt}%
  \fontfamily{#3}\fontseries{#4}\fontshape{#5}%
  \selectfont}%
\fi\endgroup%
\begin{picture}(3643,4771)(4456,-6409)
\put(7799,-1821){\makebox(0,0)[lb]{\smash{{\SetFigFont{12}{14.4}{\rmdefault}{\mddefault}{\updefault}{\color[rgb]{0,0,0}$\Ms(\Z)$}%
}}}}
\put(5626,-2161){\makebox(0,0)[lb]{\smash{{\SetFigFont{12}{14.4}{\rmdefault}{\mddefault}{\updefault}{\color[rgb]{0,0,0}$a_\sigma$}%
}}}}
\put(6491,-2491){\makebox(0,0)[lb]{\smash{{\SetFigFont{12}{14.4}{\rmdefault}{\mddefault}{\updefault}{\color[rgb]{0,0,0}$\sigma : \Q \hookrightarrow \C$}%
}}}}
\put(6481,-4856){\makebox(0,0)[lb]{\smash{{\SetFigFont{12}{14.4}{\rmdefault}{\mddefault}{\updefault}{\color[rgb]{0,0,0}$a_p^\eps$}%
}}}}
\put(5626,-3116){\makebox(0,0)[lb]{\smash{{\SetFigFont{12}{14.4}{\rmdefault}{\mddefault}{\updefault}{\color[rgb]{0,0,0}$a_\sigma^\eps$}%
}}}}
\put(5561,-4101){\makebox(0,0)[lb]{\smash{{\SetFigFont{12}{14.4}{\rmdefault}{\mddefault}{\updefault}{\color[rgb]{0,0,0}$a_0$}%
}}}}
\put(4471,-5431){\makebox(0,0)[lb]{\smash{{\SetFigFont{12}{14.4}{\rmdefault}{\mddefault}{\updefault}{\color[rgb]{0,0,0}$\tilde{a}_2$}%
}}}}
\put(4966,-5626){\makebox(0,0)[lb]{\smash{{\SetFigFont{12}{14.4}{\rmdefault}{\mddefault}{\updefault}{\color[rgb]{0,0,0}$\tilde{a}_3$}%
}}}}
\put(7321,-5491){\makebox(0,0)[lb]{\smash{{\SetFigFont{12}{14.4}{\rmdefault}{\mddefault}{\updefault}{\color[rgb]{0,0,0}$\tilde{a}_p$}%
}}}}
\put(4741,-4141){\makebox(0,0)[lb]{\smash{{\SetFigFont{12}{14.4}{\rmdefault}{\mddefault}{\updefault}{\color[rgb]{0,0,0}0}%
}}}}
\put(4711,-3106){\makebox(0,0)[lb]{\smash{{\SetFigFont{12}{14.4}{\rmdefault}{\mddefault}{\updefault}{\color[rgb]{0,0,0}$\eps$}%
}}}}
\put(6421,-5656){\makebox(0,0)[lb]{\smash{{\SetFigFont{12}{14.4}{\rmdefault}{\mddefault}{\updefault}{\color[rgb]{0,0,0}$\eps$}%
}}}}
\put(5881,-5176){\makebox(0,0)[lb]{\smash{{\SetFigFont{12}{14.4}{\rmdefault}{\mddefault}{\updefault}{\color[rgb]{0,0,0}0}%
}}}}
\put(4741,-2251){\makebox(0,0)[lb]{\smash{{\SetFigFont{12}{14.4}{\rmdefault}{\mddefault}{\updefault}{\color[rgb]{0,0,0}1}%
}}}}
\put(7081,-6331){\makebox(0,0)[lb]{\smash{{\SetFigFont{12}{14.4}{\rmdefault}{\mddefault}{\updefault}{\color[rgb]{0,0,0}$+\infty$}%
}}}}
\end{picture}%

%% file: sheaf.pstex_t
\begin{picture}(0,0)%
\includegraphics{sheaf.pstex}%
\end{picture}%
\setlength{\unitlength}{3947sp}%
\begingroup\makeatletter\ifx\SetFigFont\undefined%
\gdef\SetFigFont#1#2#3#4#5{%
  \reset@font\fontsize{#1}{#2pt}%
  \fontfamily{#3}\fontseries{#4}\fontshape{#5}%
  \selectfont}%
\fi\endgroup%
\begin{picture}(3654,3149)(5177,-5064)
\put(6102,-4476){\makebox(0,0)[lb]{\smash{{\SetFigFont{9}{10.8}{\rmdefault}{\mddefault}{\updefault}{\color[rgb]{0,0,0}3}%
}}}}
\put(8216,-4411){\makebox(0,0)[lb]{\smash{{\SetFigFont{9}{10.8}{\rmdefault}{\mddefault}{\updefault}{\color[rgb]{0,0,0}$p$}%
}}}}
\put(6869,-4648){\makebox(0,0)[lb]{\smash{{\SetFigFont{9}{10.8}{\rmdefault}{\mddefault}{\updefault}{\color[rgb]{0,0,0}5}%
}}}}
\put(8571,-3100){\makebox(0,0)[lb]{\smash{{\SetFigFont{9}{10.8}{\rmdefault}{\mddefault}{\updefault}{\color[rgb]{0,0,0}$\Os=\Z\left[\frac{1}{6}\right]$}%
}}}}
\put(7566,-2450){\makebox(0,0)[lb]{\smash{{\SetFigFont{9}{10.8}{\rmdefault}{\mddefault}{\updefault}{\color[rgb]{0,0,0}$\Os=\R$}%
}}}}
\put(5263,-4470){\makebox(0,0)[lb]{\smash{{\SetFigFont{9}{10.8}{\rmdefault}{\mddefault}{\updefault}{\color[rgb]{0,0,0}$\Os=\Z_3$}%
}}}}
\put(5192,-3643){\makebox(0,0)[lb]{\smash{{\SetFigFont{9}{10.8}{\rmdefault}{\mddefault}{\updefault}{\color[rgb]{0,0,0}$\Os=\Q_2$}%
}}}}
\put(5698,-3976){\makebox(0,0)[lb]{\smash{{\SetFigFont{9}{10.8}{\rmdefault}{\mddefault}{\updefault}{\color[rgb]{0,0,0}2}%
}}}}
\end{picture}%

%% file: droite.tex
\chapter[Droite affine sur un corps de nombres]{Droite affine analytique au-dessus d'un anneau d'entiers de corps de nombres}\label{chapitredroite}

\index{Droite affine analytique!sur A@sur $A$|(}

Dans le chapitre pr\'ec\'edent, nous sommes parvenu \`a exhiber des syst\`emes fondamentaux de voisinages pour certains points de l'espace affine au-dessus d'un anneau d'entiers de corps de nombres et \`a \'etablir certaines propri\'et\'es des anneaux locaux en ces points. Notre \'etude reste cependant incompl\`ete ; nous allons la mener \`a terme dans le cadre de la droite affine. 

Nous commencerons, au num\'ero \ref{recapitulatif}, par rappeler les r\'esultats dont nous disposons d\'ej\`a et les appliquer au cas de la droite. Nous observerons notamment que, dans ce cadre, n'\'echappent \`a  notre \'etude que certains points des fibres centrale et extr\^emes, \`a savoir les points de type~$3$ et~$2$, auxquels nous consacrerons respectivement les num\'eros \ref{chgt} et \ref{pdt2}. 

Nous regroupons au num\'ero \ref{sectionresume} les r\'esultats d\'emontr\'es jusqu'alors et prouvons, en outre, la validit\'e du principe du prolongement analytique.

Finalement, nous montrons au num\'ero \ref{parcoherence} que le faisceau structural sur la droite affine analytique au-dessus d'un anneau d'entiers de corps de nombres est coh\'erent. L'on sait l'importance que rev\^et cette propri\'et\'e en g\'eom\'etrie alg\'ebrique et en g\'eom\'etrie analytique complexe. Elle se r\'ev\'elera, pour nous, capitale au chapitre \ref{chapitreapplications}, puisqu'elle nous permettra d'utiliser les r\'esultats sur les espaces de Stein d\'emontr\'es au chapitre \ref{chapitreStein}.

\bigskip

Dans ce chapitre, comme dans le pr\'ec\'edent, nous fixons un corps de nombres~$K$ et notons~$A$ l'anneau de ses entiers. Nous posons
$$B=\Ms(A).$$ 
Puisque nous nous int\'eressons ici \`a la droite affine analytique, nous posons
$$X=\E{1}{A}.$$
Les faisceaux structuraux sur ces espaces seront respectivement not\'es~$\Os_{B}$ et~$\Os_{X}$. Lorsqu'aucune confusion ne peut en d\'ecouler, nous les noterons simplement~$\Os$.
\newcounter{nod}\setcounter{nod}{\thepage}

Nous noterons~$T$ la variable sur l'espace~$X$. Nous d\'esignerons finalement par
$$\pi : X \to B$$ 
le morphisme de projection induit par le morphisme naturel $A\to A[T]$. Pour toute partie~$V$ de~$B$, nous posons
$$X_{V} = \pi^{-1}(V)$$
et, pour tout point~$b$ de~$B$,
$$X_{b} = \pi^{-1}(b).$$




\section{R\'ecapitulatif}\label{recapitulatif}

Commen\c{c}ons par appliquer au cas de la droite les r\'esultats que nous avons d\'emontr\'es pour les espaces affines. Commen\c{c}ons par les points rigides des fibres.

\begin{thm}\label{toporig}\index{Connexite par arcs au voisinage d'un point@Connexit\'e par arcs au voisinage d'un point!rigide de A1@rigide de \AA}\index{Ouverture au voisinage d'un point!rigide de A1@rigide de \AA}
Soient~$b$ un point de l'espace~$B$ et~$x$ un point rigide de la fibre~$X_{b}$. Le point~$x$ poss\`ede un syst\`eme fondamental de voisinages connexes par arcs dans~$X$ et le morphisme de projection~$\pi$ est ouvert au point~$x$.
\end{thm}
\begin{proof}
Ce r\'esultat est une cons\'equence des propositions~\ref{cparig} et \ref{ouvertrig}.
\end{proof}

En ce qui concerne les propri\'et\'es de l'anneau local, nous distinguerons deux cas.

\begin{thm}\label{recrigavd}\index{Anneau local en un point!rigide!d'une fibre interne de A1@d'une fibre interne de $\AA$}\index{Anneau local en un point!rigide!de la fibre centrale de A1@de la fibre centrale de $\AA$}
Soient~$b$ un point de l'espace~$B$ qui n'est pas un point extr\^eme et~$x$ un point rigide de la fibre~$X_{b}$. L'anneau local~$\Os_{X,x}$ est un anneau de valuation discr\`ete hens\'elien. Son corps r\'esiduel~$\kappa(x)$ est complet, et donc isomorphe \`a~$\Hs(x)$.
\end{thm}
\begin{proof}
Remarquons tout d'abord que l'anneau local~$\Os_{B,b}$ est un corps. D'apr\`es la proposition~\ref{isorig}, nous pouvons supposer que le point~$x$ est rationnel dans sa fibre. D'apr\`es le lemme \ref{translation}, nous pouvons supposer que c'est le point~$0$ de cette fibre. Le th\'eor\`eme~\ref{anneaulocal} permet alors de ramener l'\'etude \`a celle de l'anneau local~$L_{b}$. D'apr\`es les th\'eor\`emes~\ref{noetheriencorps} et~\ref{factoriel}, ce dernier anneau est noeth\'erien et factoriel. D'apr\`es le lemme~\ref{idealmax}, son id\'eal maximal est engendr\'e par l'\'el\'ement~$T$, qui n'est pas nilpotent. La proposition~$2$ de~\cite{CL} assure alors que l'anneau~$L_{b}$ est de valuation discr\`ete. Le caract\`ere hens\'elien, quant \`a lui, d\'ecoule de la proposition~\ref{Hensel}. Le r\'esultat concernant le corps r\'esiduel~$\kappa(x)$ d\'ecoule du th\'eor\`eme~\ref{noetherienrigcentralinterne}.
\end{proof}

\begin{thm}\label{recrigext}\index{Anneau local en un point!rigide!d'une fibre extreme de A1@d'une fibre extr\^eme de $\AA$}
Soient~$b$ un point extr\^eme de l'espace~$B$ et~$x$ un point rigide de la fibre~$X_{b}$. L'anneau local~$\Os_{X,x}$ est un anneau hens\'elien, noeth\'erien et r\'egulier de dimension~$2$. Son corps r\'esiduel~$\kappa(x)$ est complet, et donc isomorphe \`a~$\Hs(x)$.
\end{thm}
\begin{proof}
Ce r\'esultat d\'ecoule de la proposition~\ref{Hensel} et du th\'eor\`eme \ref{noetherienrigext}.
\end{proof}

Pour les points internes, nous disposons de r\'esultats complets.

\begin{thm}\label{topointerne}\index{Anneau local en un point!interne!non-rigide de A1@non-rigide de $\AA$}\index{Connexite par arcs au voisinage d'un point@Connexit\'e par arcs au voisinage d'un point!interne!non-rigide de A1@non-rigide de \AA}\index{Ouverture au voisinage d'un point!interne!non-rigide de A1@non-rigide de \AA}
Soient~$b$ un point interne de l'espace~$B$ et~$x$ un point de la fibre~$X_{b}$ qui n'est pas un point rigide. Le point~$x$ poss\`ede un syst\`eme fondamental de voisinages connexes par arcs dans~$X$ et le morphisme de projection~$\pi$ est ouvert au point~$x$. L'anneau local~$\Os_{X,x}$ est isomorphe au corps~$\kappa(x)$, lequel est hens\'elien.
\end{thm}
\begin{proof}
La premi\`ere partie du r\'esultat d\'ecoule directement des corollaires~\ref{ouvertinterne} et~\ref{cpainterne}. La seconde d\'ecoule de la proposition~\ref{isointerne} et du r\'esultat correspondant pour la droite analytique sur un corps valu\'e complet (qui est alors n\'ecessairement ultram\'etrique).
\end{proof}

Il nous reste donc \`a \'etudier les points des fibres extr\^emes et centrale qui ne sont pas rigides. Rappelons que nous avons \'egalement d\'emontr\'e des r\'esultats pour certains points de type~$3$ de ces fibres.

\begin{thm}\label{topo3}\index{Connexite par arcs au voisinage d'un point@Connexit\'e par arcs au voisinage d'un point!de type 3 d'une fibre extreme@de type $3$ d'une fibre extr\^eme}\index{Ouverture au voisinage d'un point!de type 3 d'une fibre extreme@de type $3$ d'une fibre extr\^eme}\index{Connexite par arcs au voisinage d'un point@Connexit\'e par arcs au voisinage d'un point!de type 3 de la fibre centrale}\index{Ouverture au voisinage d'un point!de type 3 de la fibre centrale}
Soient~$b$ un point extr\^eme ou central de l'espace~$B$, $\alpha$ un \'el\'ement de~$\Hs(b)$ et~$r$ un \'el\'ement de $\R_{+}^*\setminus\{1\}$. Notons~$x$ le point~$\eta_{\alpha,r}$ de la fibre~$X_{b}$. Le point~$x$ poss\`ede un syst\`eme fondamental de voisinages connexes par arcs dans~$X$ et le morphisme de projection~$\pi$ est ouvert au point~$x$.
\end{thm}
\begin{proof}
D'apr\`es le lemme \ref{translation}, nous pouvons supposer que l'\'el\'ement~$\alpha$ est nul. Le r\'esultat d\'ecoule alors des corollaires~\ref{ouvertdep} et~\ref{cpadep}.
\end{proof}

Pour d\'ecrire les propri\'et\'es de l'anneau local, nous distinguerons deux cas.

\begin{thm}\index{Anneau local en un point!de type 3!de la fibre centrale}
Soient~$\alpha$ un \'el\'ement de~$K$ et~$r$ un \'el\'ement de $\R_{+}^*\setminus\{1\}$. Notons~$x$ le point~$\eta_{\alpha,r}$ de la fibre centrale~$X_{0}$. L'anneau local~$\Os_{X,x}$ est isomorphe au corps~$\kappa(x)$, lequel est hens\'elien.
\end{thm}
\begin{proof}
Le r\'esultat d\'ecoule du corollaire~\ref{descriptioncentral13}.

\end{proof}

\begin{thm}\label{recanneau3}\index{Anneau local en un point!de type 3!d'une fibre extreme@d'une fibre extr\^eme}
Soient~$\m$ un \'el\'ement de~$\Sigma_{f}$, $\alpha$ un \'el\'ement de~$\tilde{k}_{\m}$ et~$r$ un \'el\'ement de \mbox{$\R_{+}^*\setminus\{1\}$}. Notons~$x$ le point~$\eta_{\alpha,r}$ de la fibre extr\^eme~$\tilde{X}_{\m}$. L'anneau local~$\Os_{X,x}$ est un anneau de valuation discr\`ete d'uniformisante~$\pi_{\m}$. Son corps r\'esiduel~$\kappa(x)$ est complet, et donc isomorphe \`a~$\Hs(x)$.
\end{thm}
\begin{proof}
D'apr\`es le lemme \ref{translation}, nous pouvons supposer que l'\'el\'ement~$\alpha$ est nul. Quitte \`a changer~$T$ en~$T^{-1}$, nous pouvons supposer que~$r<1$. Le r\'esultat d\'ecoule alors du corollaire~\ref{descriptionextreme3}.

\end{proof}

\bigskip

Lorsque les anneaux locaux sont des anneaux de valuation discr\`ete, nous pouvons obtenir des informations suppl\'ementaires. \`A cet effet, nous introduisons une nouvelle d\'efinition.

\begin{defi}\index{Uniformisante forte}
Soient~$(Y,\Os_{Y})$ un espace analytique et~$y$ un point de~$Y$. Supposons que l'anneau local~$\Os_{Y,y}$ est un anneau de valuation discr\`ete. Soit~$V$ un voisinage du point~$y$ dans~$Y$ et~$\pi$ un \'el\'ement de~$\Os_{Y}(V)$. Nous dirons que la fonction~$\pi$ est une {\bf uniformisante forte de l'anneau~$\Os_{Y,y}$ sur~$V$} s'il existe un nombre r\'eel~$C$ v\'erifiant la propri\'et\'e suivante : pour tout \'el\'ement~$f$ de~$\Os_{Y}(V)$ dont l'image~$f(y)$ dans~$\Hs(y)$ est nulle, il existe un \'el\'ement~$g$ de~$\Os_{Y}(V)$ tel que
\begin{enumerate}[\it i)]
\item $f = \pi g$ dans~$\Os_{Y}(V)$ ;
\item $\|g\|_{V} \le C\, \|f\|_{V}$.
\end{enumerate}
\end{defi}

\begin{rem}
L'image dans l'anneau de valuation discr\`ete~$\Os_{Y,y}$ d'une uniformisante forte est une uniformisante.
\end{rem}

\begin{lem}\label{uniforteB}\index{Uniformisante forte!en un point extreme de MA@en un point extr\^eme de $\Ms(A)$}
Soit~$b$ un point de l'espace~$B$ tel que l'anneau local~$\Os_{B,b}$ soit un anneau de valuation discr\`ete. Soit~$\pi$ une uniformisante de l'anneau~$\Os_{B,b}$ et~$U$ un voisinage du point~$b$ dans~$B$ sur lequel elle est d\'efinie. Alors il existe un syst\`eme fondamental~$\Vs$ de voisinages compacts et connexes du point~$b$ dans~$U$ tel que, pour tout \'el\'ement~$V$ de~$\Vs$, la fonction~$\pi$ est une uniformisante forte de l'anneau~$\Os_{B,b}$ sur~$V$.
\end{lem}
\begin{proof}
Il existe un \'el\'ement~$\m$ de~$\Sigma_{f}$ tel que le point~$b$ soit le point~$\tilde{a}_{\m}$. Les descriptions explicites du num\'ero \ref{partiesouvertes} permettent de montrer que, pour tout nombre r\'eel $\eps>0$, la fonction~$\pi_{\m}$ est une uniformisante forte de l'anneau $\Os_{B,\tilde{a}_{\m}}\simeq \hat{A}_{\m}$ sur $\of{[}{a_{\m}^\eps,\tilde{a}_{\m}}{]}$. Le r\'esultat pour toute autre uniformisante s'en d\'eduit.
\end{proof}

\begin{prop}
Soit~$b$ un point de~$B$ qui n'est pas un point extr\^eme. Notons~$x$ le point~$0$ de la fibre~$X_{b}$. Soit~$V$ un voisinage compact et connexe du point~$b$ dans~$B$ dont le bord ne contient pas le point central~$a_{0}$. Soit~$t$ un nombre r\'eel strictement positif. La fonction~$T$ est une uniformisante forte de l'anneau de valuation discr\`ete~$\Os_{X,x}$ sur le disque~$\overline{D}_{V}(t)$.
\end{prop}
\begin{proof}
Soit~$f$ un \'el\'ement de~$\Os(\overline{D}_{V}(t))$ dont l'image dans~$\Hs(x)$ est nulle. D'apr\`es la proposition~\ref{imagedisque}, il existe un nombre r\'eel~$r > t$ et une suite~$(f_{k})_{k\ge 0}$ d'\'el\'ements de~$\Os(V)$ v\'erifiant la condition $\lim_{k\to +\infty} \|f_k\|_{V}\, r^k = 0$ tels que l'on ait l'\'egalit\'e
$$f = \sum_{k \ge 0} f_{k} T^k.$$
Par hypoth\`ese, nous avons~$f(x)=0$ et donc~$f_{0}(x)=f_{0}(b)=0$. Puisque le point~$b$ n'est pas extr\^eme, l'anneau local~$\Os_{B,b}$ est un corps. Par cons\'equent, la fonction~$f_{0}$ est nulle au voisinage du point~$b$ dans~$V$. D'apr\`es le principe du prolongement analytique, elle est nulle dans~$\Os(V)$. Maintenant, le th\'eor\`eme~\ref{isodisque} assure que la s\'erie
$$g = \sum_{k \ge 0} f_{k+1} T^k$$
d\'efinit un \'el\'ement de~$\Os(\overline{D}_{V}(t))$. Par cons\'equent, nous avons l'\'egalit\'e
$$f = Tg \textrm{ dans } \Os(\overline{D}_{V}(t)).$$
D'apr\`es le lemme~\ref{Shilovcouronne}, le disque~$\overline{D}_{V}(t)$ poss\`ede un bord analytique~$\Gamma$ qui v\'erifie la propri\'et\'e suivante :
$$\forall y\in\Gamma,\, |T(y)|=t.$$
Soit~$y$ un point de~$\Gamma$ en lequel la fonction~$g$ atteint son maximum. Nous avons alors
$${\renewcommand{\arraystretch}{1.5}\begin{array}{rcl}
\|g\|_{\overline{D}_{V}(t)} &=& |g(y)|\\
&=& |T(y)|^{-1}\, |f(y)|\\
&\le & t^{-1}\,\|f\|_{\overline{D}_{V}(t)}.
\end{array}}$$
\end{proof}

\begin{cor}\label{uniforterig}\index{Uniformisante forte!en un point rigide!d'une fibre interne de An@d'une fibre interne de $\E{n}{A}$}\index{Uniformisante forte!en un point rigide!de la fibre centrale de An@de la fibre centrale de $\E{n}{A}$}
Soient~$b$ un point de~$B$ qui n'est pas un point extr\^eme et~$x$ un point rigide de la fibre~$X_{b}$. Soient~$\pi$ une uniformisante de l'anneau de valuation discr\`ete~$\Os_{X,x}$ et~$U$ un voisinage du point~$x$ dans~$X$ sur lequel elle est d\'efinie. Alors il existe un syst\`eme fondamental~$\Vs$ de voisinages compacts et connexes du point~$x$ dans~$U$ tel que, pour tout \'el\'ement~$V$ de~$\Vs$, la fonction~$\pi$ est une uniformisante forte de l'anneau~$\Os_{X,x}$ sur~$V$.
\end{cor}
\begin{proof}
D'apr\`es la proposition~\ref{isorig}, nous pouvons supposer que le point~$x$ est rationnel dans sa fibre. D'apr\`es le lemme \ref{translation}, nous pouvons supposer que c'est le point~$0$ de cette fibre. La proposition pr\'ec\'edente jointe \`a la proposition~\ref{voisdep} nous permet alors de conclure lorsque l'uniformisante consid\'er\'ee est~$T$. Le r\'esultat pour toute autre uniformisante s'en d\'eduit.


\end{proof}

\begin{prop}
Soient~$\m$ un \'el\'ement de~$\Sigma_{f}$ et~$r$ un \'el\'ement de \mbox{$\R_{+}^*\setminus\{1\}$}. Soient~$s$ et~$t$ deux \'el\'ements de~$\R_{+}^*$ qui v\'erifient~$s<r<t$. Notons~$x$ le point~$\eta_{r}$ de la fibre extr\^eme~$\tilde{X}_{\m}$. Soit~$\eps$ un \'el\'ement de~$\R_{+}^*$. Consid\'erons la couronne
$$C = \left\{y\in\pi^{-1}(\of{[}{a_{\m}^\eps,\tilde{a}_{\m}}{]})\, \big|\, s<|T(y)|<t\right\}.$$
La fonction~$\pi_{\m}$ est une uniformisante forte de l'anneau de valuation discr\`ete~$\Os_{X,x}$ sur la couronne~$C$.
\end{prop}
\begin{proof}
Soit~$f$ un \'el\'ement de~$\Os(C)$ dont l'image dans~$\Hs(x)$ est nulle. Remarquons que l'anneau norm\'e $(\Os(V),\|.\|_{V})$ n'est autre que l'anneau $(\hat{A}_{\m},|.|_{\m}^\eps)$. D'apr\`es la proposition~\ref{imagecouronne}, il existe deux nombres r\'eels~$s_{0}$ et~$t_{0}$ v\'erifiant $0<s_{0}<s<t<t_{0}$ et une suite~$(f_{k})_{k\ge 0}$ d'\'el\'ements de~$\hat{A}_{\m}$ v\'erifiant la condition $\lim_{k\to +\infty} |f_k|_{\m}^\eps\, r^k = 0$ tels que l'on ait l'\'egalit\'e
$$f = \sum_{k \ge 0} f_{k} T^k.$$
Par hypoth\`ese, nous avons~$f(x)=0$ et donc
$$\forall k\in\N, f_{k}(x)=f_{k}(\tilde{a}_{\m})=0.$$ 
On en d\'eduit que, pour tout \'el\'ement~$k$ de~$\N$, il existe un \'el\'ement~$g_{k}$ de~$\hat{A}_{\m}$ tel que l'on ait l'\'egalit\'e
$$f_{k} = \pi_{\m}\, g_{k}.$$
En outre, pour tout \'el\'ement~$k$ de~$\N$, nous avons
$$|g_{k}|_{\m}^\eps = |\pi_{\m}|_{\m}^{-\eps}\, |f_{k}|_{\m}^\eps.$$
Par cons\'equent, la s\'erie
$$g = \sum_{k \ge 0} g_{k} T^k$$
d\'efinit un \'el\'ement de l'anneau $\Os(V)\of{\la}{s\le |T| \le t}{\ra}^\dag$ et donc de l'anneau~$\Os(C)$, d'apr\`es le th\'eor\`eme~\ref{isocouronne}. Nous avons alors l'\'egalit\'e
$$f = \pi_{\m}\, g \textrm{ dans } \Os(C).$$
D'apr\`es le lemme~\ref{Shilovcouronneum}, la couronne~$C$ poss\`ede un bord analytique~$\Gamma$ qui v\'erifie la propri\'et\'e suivante :
$$\forall y\in\Gamma,\, |\pi_{\m}(y)|=|\pi_{\m}|_{\m}^\eps.$$
Soit~$y$ un point de~$\Gamma$ en lequel la fonction~$g$ atteint son maximum. Nous avons alors
$${\renewcommand{\arraystretch}{1.5}\begin{array}{rcl}
\|g\|_{C} &=& |g(y)|\\
&=& |\pi_{\m}(y)|^{-1}\, |f(y)|\\
&\le & |\pi_{\m}|_{\m}^{-\eps}\,\|f\|_{C}.
\end{array}}$$
\end{proof}

\begin{cor}\label{recuniforte3depext}\index{Uniformisante forte!en un point de type 3 d'une fibre extreme@en un point de type $3$ d'une fibre extr\^eme}
Soit~$\m$ un \'el\'ement de~$\Sigma_{f}$. Soient~$\alpha$ un \'el\'ement de~$\tilde{k}_{\m}$ et~$r$ un \'el\'ement de~$\R_{+}^*\setminus\{1\}$. Notons~$x$ le point~$\eta_{\alpha,r}$ de la fibre extr\^eme~$\tilde{X}_{\m}$. Soient~$\pi$ une uniformisante de l'anneau de valuation discr\`ete~$\Os_{X,x}$ et~$U$ un voisinage du point~$x$ dans~$X$ sur lequel elle est d\'efinie. Il existe un syst\`eme fondamental~$\Vs$ de voisinages compacts et connexes du point~$x$ dans~$U$ tel que, pour tout \'el\'ement~$V$ de~$\Vs$, la fonction~$\pi$ est une uniformisante forte de l'anneau~$\Os_{X,x}$ sur~$V$.
\end{cor}
\begin{proof}
D'apr\`es le lemme \ref{translation}, nous pouvons supposer que le point~$x$ est le point~$\eta_{r}$ de la fibre~$\tilde{X}_{\m}$. La proposition pr\'ec\'edente jointe \`a la proposition~\ref{voisdep} nous permet alors de conclure lorsque l'uniformisante consid\'er\'ee est~$\pi_{\m}$. Le r\'esultat pour toute autre uniformisante s'en d\'eduit im\-m\'e\-dia\-te\-ment.
\end{proof}

\bigskip

Int\'eressons-nous, maintenant, au bord analytique de voisinages des points. Nous nous contentons de rappeler ici les r\'esultats des propositions \ref{shilovfiniinterne}, \ref{Shilovrigum} et \ref{Shilov3depum}.

\begin{prop}\label{recShilovinterneum}\index{Bord analytique!au voisinage d'un point interne!de A1@de $\AA$}
Soit~$\sigma\in\Sigma_{f}$. Tout point de~$X'_{\sigma}$ poss\`ede un syst\`eme fondamental de voisinages compacts, connexes et spectralement convexes qui poss\`edent un bord analytique fini et al\-g\'e\-bri\-que\-ment trivial.
\end{prop}

\begin{prop}\label{recShilovrigideum}\index{Bord analytique!au voisinage d'un point rigide!d'une fibre extreme de A1@d'une fibre extr\^eme de $\AA$}
Soit~$b$ un point de $B_{\textrm{um}}\setminus\{a_{0}\}$. Tout point rigide de la fibre~$X_{b}$ poss\`ede un syst\`eme fondamental de voisinages compacts, connexes et spectralement convexes qui poss\`edent un bord de Shilov fini et alg\'ebriquement trivial.
\end{prop}

\begin{prop}\label{recShilov3depum}\index{Bord analytique!au voisinage d'un point de type 3 d'une fibre extreme@au voisinage d'un point de type $3$ d'une fibre extr\^eme}
Soit~$b$ un point de $B_{\textrm{um}}\setminus\{a_{0}\}$. Tout point de type~$3$ d\'eploy\'e de la fibre~$X_{b}$ poss\`ede un syst\`eme fondamental de voisinages compacts, connexes et spectralement convexes qui poss\`edent un bord de Shilov fini et al\-g\'e\-bri\-que\-ment trivial.
\end{prop}

\bigskip

Pour finir, int\'eressons-nous au principe du prolongement analytique. 

\begin{prop}
Soit~$b$ un point de l'espace~$B$. Soit~$V$ un voisinage ouvert et connexe du point~$b$ dans~$B$. Soit~$r$ un \'el\'ement de l'intervalle~$\of{]}{0,1}{[}$. Le principe du prolongement analytique vaut sur le disque~$\mathring{D}_{V}(r)$.
\end{prop}
\begin{proof}
D'apr\`es le corollaire~\ref{corprolan}, il suffit de montrer que le morphisme de projection~$\pi$ est ouvert au voisinage de tout point du disque~$\mathring{D}_{V}(r)$. Ce r\'esultat d\'ecoule des th\'eor\`emes~\ref{toporig}, \ref{topointerne} et~\ref{topo3}.
\end{proof}

\begin{cor}\label{prolanrig}\index{Prolongement analytique!au voisinage d'un point rigide!de A1@de $\AA$}
Le principe du prolongement analytique vaut au voisinage des points rigides des fibres de l'espace~$X$.
\end{cor}
\begin{proof}
Soient~$b$ un point de l'espace~$B$ et~$x$ un point rigide de la fibre~$X_{b}$. D'apr\`es la proposition~\ref{isorig} et le lemme \ref{translation}, nous pouvons supposer que le point~$x$ est le point~$0$ de la fibre~$X_{b}$. D'apr\`es la proposition~\ref{voisdep}, la famille des disques ouverts~$\mathring{D}_{V}(r)$, o\`u~$V$ parcourt l'ensemble des voisinages ouverts et connexes de~$b$ dans~$B$ et~$r$ l'intervalle~$\of{]}{0,1}{[}$, est un syst\`eme fondamental de voisinages du point~$x$ dans~$X$. Nous concluons alors en utilisant le lemme~\ref{lemprolan} et la proposition pr\'ec\'edente.
\end{proof}

\begin{prop}
Soit~$b$ un point de l'espace~$B$. Soit~$V$ un voisinage ouvert et connexe du point~$b$ dans~$B$. Soient~$s$ et~$t$ deux nombres r\'eels qui v\'erifient la condition $0<s<t<1$. Le principe du prolongement analytique vaut sur la couronne~$\mathring{C}_{V}(s,t)$.
\end{prop}
\begin{proof}
D'apr\`es le corollaire~\ref{corprolan}, il suffit de montrer que le morphisme de projection~$\pi$ est ouvert au voisinage de tout point de la couronne ouverte~$\mathring{C}_{V}(s,t)$. Ce r\'esultat d\'ecoule des th\'eor\`emes~\ref{topointerne} et~\ref{topo3}.
\end{proof}

\begin{cor}\label{prolan3dep}\index{Prolongement analytique!au voisinage d'un point de type 3!d'une fibre extreme@d'une fibre extr\^eme}\index{Prolongement analytique!au voisinage d'un point de type 3!de la fibre centrale}
Soient~$b$ un point extr\^eme ou central de l'espace~$B$, $\alpha$ un \'el\'ement de~$\Hs(b)$ et~$r$ un \'el\'ement de~$\R_{+}^*\setminus\{1\}$. Notons~$x$ le point~$\eta_{\alpha,r}$ de la fibre~$X_{b}$. Le principe du prolongement analytique vaut au voisinage du point~$x$ de l'espace~$X$.
\end{cor}
\begin{proof}
D'apr\`es le lemme \ref{translation}, nous pouvons supposer que l'\'el\'ement~$\alpha$ de~$\Hs(b)$ est nul. Quitte \`a changer~$T$ en~$T^{-1}$, nous pouvons supposer que l'\'el\'ement~$r$ appartient \`a l'intervalle~$\of{]}{0,1}{[}$. D'apr\`es la proposition~\ref{voisdep}, la famille des couronnes ouvertes~$\mathring{C}_{V}(s,t)$, o\`u~$V$ parcourt l'ensemble des voisinages ouverts et connexes de~$b$ dans~$B$, $s$ l'intervalle~$\of{]}{0,r}{[}$ et~$t$ l'intervalle~$\of{]}{r,1}{[}$, est un syst\`eme fondamental de voisinages du point~$x$ dans~$X$. Nous concluons alors en utilisant le lemme~\ref{lemprolan} et la proposition pr\'ec\'edente.
\end{proof}

Concluons en rappelant le r\'esultat de la proposition~\ref{prolinterne}.

\begin{prop}\label{recprolinterne}\index{Prolongement analytique!au voisinage d'un point interne!de A1@de $\AA$}
Le principe du prolongement analytique vaut au voisinage des points internes de l'espace~$X$.
\end{prop}

\section{Points de type~$3$}\label{chgt}

Nous nous int\'eressons ici aux points de type~$3$ des fibres extr\^emes et centrale. Un changement de base va nous permettre de nous ramener au cas de points de type~$3$ d\'eploy\'es. \`A cet effet, nous allons \'etendre le r\'esultat des propositions~\ref{isoext} et~\ref{isocentral}.

\subsection{Fibres extr\^emes}\index{Point!de type 3!d'une fibre extreme@d'une fibre extr\^eme|(}

Traitons, tout d'abord, le cas des fibres extr\^emes. Nous commencerons par montrer que l'on peut pr\'eciser le r\'esultat de changement de base obtenu \`a la proposition \ref{isoext}. Soit $\m\in\Sigma_{f}$. Soit $P(T)$ un polyn\^ome irr\'eductible \`a coefficients dans~$k_{\m}$. Rappelons que, quel que soit $r\in\of{[}{0,1}{]}$, nous notons~$\eta_{P,r}$ le point de la fibre~$\tilde{X}_{\m}$ associ\'e \`a la valeur absolue 
$$\begin{array}{ccc}
A[T] & \to & \R_{+}\\
F(T) & \mapsto & r^{v_{P(T)}(F(T))}
\end{array},$$
o\`u $v_{P(T)}$ d\'esigne la valuation $P(T)$-adique de $k_{\m}[T]$. Pour $\alpha\in k_{\m}$ et $r\in\of{[}{0,1}{]}$, nous notons 
$$\eta_{\alpha,r} = \eta_{T-\alpha,r}.$$
Pour~$r\in\of{[}{0,1}{]}$, nous notons encore.
$$\eta_{r}=\eta_{0,r}=\eta_{T,r}.$$
Finalement, pour $r\in\of{[}{1,+\infty}{[}$, nous notons~$\eta_{r}$ le point de la fibre~$\tilde{X}_{\m}$ associ\'e \`a la valeur absolue 
$$\begin{array}{ccc}
A[T] & \to & \R_{+}\\
F(T) & \mapsto & r^{-\deg(\tilde{F}(T))}
\end{array},$$
o\`u~$\tilde{F}(T)$ d\'esigne l'image du polyn\^ome~$F(T)$ dans~$k_{\m}[T]$. Nous avons ainsi d\'ecrit tous les points de la fibre extr\^eme~$\tilde{X}_{\m}$ (\emph{cf.}~\ref{droitetrivval} pour la classification, avec d\'emonstration, des points de la droite analytique sur un corps trivialement valu\'e quelconque). Les points de type~$3$ sont ceux pour lesquels le nombre r\'eel~$r$ est diff\'erent de~$0$ et de~$1$.

Nous noterons $x=\eta_{P,0}$ le point rigide de la fibre $\tilde{X}_{\m}$ d\'efini par l'\'equation
$$P(T)(x)=0.$$
D'apr\`es la proposition \ref{isoext}, il existe une extension finie $K'$ de $K$, un point~$x'$ de $X' = \E{1}{A'}$, o\`u $A'$ d\'esigne l'anneau des entiers de $K'$, rationnel dans sa fibre, tel que le morphisme naturel
$$\varphi : \E{1}{A'} \to \E{1}{A}$$
envoie le point $x'$ sur le point $x$ et induise un isomorphisme d'un voisinage de $x'$ sur un voisinage de $x$. Notons $\m'$ l'id\'eal maximal de $A'$ correspondant au point~$\pi(x')$ et $\alpha$ l'\'el\'ement de $k_{\m'}$ qui correspond au point $x'$. Un calcul direct utilisant la s\'eparabilit\'e du polyn\^ome~$P(T)$ montre que, pour tout \'el\'ement~$r$ de l'intervalle $\of{[}{0,1}{]}$, nous avons
$$\varphi(\eta_{\alpha,r}) = \eta_{P,r}.$$

Nous devons reprendre et pr\'eciser ici les arguments de la proposition~\ref{isoext}. Nous aurons besoin d'utiliser certaines propri\'et\'es du flot et commen\c{c}ons donc par montrer l'existence de voisinages flottants. Posons
$$Y_{\m} = X_{\m} \setminus X_{0} = \pi^{-1}(\of{]}{a_{0},\tilde{a}_{\m}}{]}).$$

\begin{lem}\label{voisflotext}\index{Voisinages flottants!pour un point de type 3!d'une fibre extreme@d'une fibre extr\^eme}
Soient~$x\in Y_{\m}$ et~$\eps\in I_{Y_{\m}}(x)$. Alors, la partie~$D_{Y_{\m}}$ est un voisinage de~$(x,\eps)$ dans~$Y_{\m}\times\R_{+}^*$. 

En particulier, tous les points de~$Y_{\m}$ ont des voisinages flottants dans~$Y_{\m}$.
\end{lem} 
\begin{proof}
Ce r\'esultat d\'ecoule directement de l'\'egalit\'e
$$D_{Y_{\m}} = Y_{\m}\times\R_{+}.$$
La cons\'equence suit, par le lemme~\ref{criterevoisflot}.
\end{proof}

\begin{prop}\label{isoext3}\index{Isomorphisme local!au voisinage d'un point de type 3!d'une fibre extreme@d'une fibre extr\^eme}
Le morphisme $\varphi$ induit un isomorphisme d'espaces annel\'es d'un voisinage de 
$$\{\eta_{\alpha,r},\, r\in\of{[}{0,1}{[}\} \textrm{ dans } X'$$ 
sur un voisinage de 
$$\{\eta_{P,r},\, r\in\of{[}{0,1}{[}\} \textrm{ dans } X.$$ 
\end{prop}
\begin{proof}
Consid\'erons le voisinage~$U$ de~$x$ dans~$X$, la fonction~$R$ d\'efinie sur~$U$ v\'erifiant~$P(R)=0$ et la section~$\sigma$ du morphisme~$\varphi$ au-dessus de~$U$ consid\'er\'es dans la preuve de la proposition~\ref{isoext}. Soit~$V$ un voisinage du point~$x$ dans~$X$ v\'erifiant les propri\'et\'es suivantes :  
\begin{enumerate}[\it i)]
\item $V$ est connexe ;
\item la fonction $R$ se prolonge \`a $V$ et la fonction $P(R)$ est nulle sur $V$ ;
\item la fonction $P'(\alpha)$ est inversible sur $\varphi^{-1}(V)$.
\end{enumerate}  
D'apr\`es la proposition \ref{isolocal}, la section $\sigma$ se prolonge alors \`a~$V$ et induit un isomorphisme d'espaces annel\'es sur son image. Il nous suffit donc de montrer qu'il existe un voisinage~$V$ de la partie 
$\{\eta_{P,r},\, r\in\of{[}{0,1}{[}\}$ dans~$X$ qui v\'erifie les propri\'et\'es demand\'ees.

Commen\c{c}ons par la derni\`ere propri\'et\'e. Quel que soit~$b\in B_{\m}\setminus X_{0}$, le polyn\^ome~$P(T)$ est irr\'eductible et s\'eparable sur le corps~$\Hs(b)$. Par cons\'equent, tout voisinage~$V$ contenu dans~$B_{\m}\setminus X_{0}$ satisfait cette propri\'et\'e. 

Passons aux deux propri\'et\'es suivantes. Il existe $r_{0}\in\of{]}{0,1}{[}$ tel que le point~$\eta_{P,r_{0}}$ appartienne \`a~$U$. En utilisant l'isomorphisme~$\sigma$ et le corollaire~\ref{cpadep}, on montre que le point~$\eta_{P,r_{0}} = \sigma^{-1}(\eta_{\alpha,r_{0}})$ de~$X$ poss\`ede un syst\`eme fondamental de voisinages connexes par arcs. Le lemme~\ref{voisflotext} assure que nous sommes dans les conditions d'utilisation de la proposition~\ref{flot} et du lemme~\ref{connexeflot}. On en d\'eduit que la fonction~$R$ se prolonge sur un voisinage connexe~$V$ de l'ensemble
$$T_{Y_{\m}}(\eta_{P,r_{0}}) = \{\eta_{P,r_{0}}^\eps,\, \eps\in\of{]}{0,+\infty}{[}\} = \{\eta_{P,r},\, r\in\of{]}{0,1}{[}\}.$$
En outre, nous avons encore $P(R)=0$ sur $V$, toujours d'apr\`es la proposition~\ref{flot}. On en d\'eduit le r\'esultat annonc\'e.
\end{proof}

Cet \'enonc\'e nous permet de ramener l'\'etude des points de type~$3$ de la fibre extr\^eme \`a celle des points de type~$3$ d\'eploy\'es. Nous en tirons plusieurs cons\'equences.

\begin{cor}\label{cpa3ext}\index{Connexite par arcs au voisinage d'un point@Connexit\'e par arcs au voisinage d'un point!de type 3 d'une fibre extreme@de type $3$ d'une fibre extr\^eme}
Tout point de type~$3$ d'une fibre extr\^eme poss\`ede un syst\`eme fondamental de voisinages connexes par arcs. 
\end{cor}
\begin{proof}
Soient~$\m\in\Sigma_{f}$ et~$x$ un point de type~$3$ de la fibre extr\^eme~$\tilde{X}_{\m}$. Supposons, tout d'abord, qu'il existe un \'el\'ement~$r>1$ tel que le point~$x$ soit le point~$\eta_{r}$. Le r\'esultat d\'ecoule alors du corollaire~\ref{cpadep}.
 
Dans les autres cas, il existe un polyn\^ome irr\'eductible~$P$ \`a coefficients dans~$k_{\m}$ et un \'el\'ement~$r$ de l'intervalle $\of{]}{0,1}{[}$ tel que le point~$x$ soit le point~$\eta_{P,r}$. Dans ce cas, la proposition~\ref{isoext3} nous montre que, quitte \`a remplacer l'anneau~$A$ par l'anneau des entiers d'une extension du corps~$K$, nous pouvons supposer que le polyn\^ome~$P$ est de degr\'e~$1$. Le r\'esultat d\'ecoule alors du corollaire~\ref{cpadep}.
\end{proof}

De m\^eme, en utilisant le corollaire~\ref{ouvertdep}, on d\'emontre le r\'esultat suivant.

\begin{cor}\label{ouvert3ext}\index{Ouverture au voisinage d'un point!de type 3 d'une fibre extreme@de type $3$ d'une fibre extr\^eme}
Le morphisme $\pi$ est ouvert en tout point de type~$3$ d'une fibre extr\^eme. 
\end{cor}

Venons-en, \`a pr\'esent, aux propri\'et\'es des anneaux locaux.

\begin{cor}\label{avd3ext}\index{Anneau local en un point!de type 3!d'une fibre extreme@d'une fibre extr\^eme}
Soient~$\m$ un \'el\'ement de~$\Sigma_{f}$ et~$x$ un point de type~$3$ de la fibre extr\^eme~$\tilde{X}_{\m}$. L'anneau local~$\Os_{X,x}$ est un anneau de valuation discr\`ete d'id\'eal maximal~$\m\, \Os_{X,x}$. Son corps r\'esiduel~$\kappa(x)$ est complet, et donc isomorphe \`a~$\Hs(x)$.
\end{cor}
\begin{proof}
Supposons, tout d'abord, qu'il existe un \'el\'ement~$r>1$ tel que le point~$x$ soit le point~$\eta_{r}$. 

Dans les autres cas, il existe un un polyn\^ome irr\'eductible~$P$ \`a coefficients dans~$k_{\m}$ et un \'el\'ement~$r$ de l'intervalle $\of{]}{0,1}{[}$ tels que le point~$x$ soit le point~$\eta_{P,r}$. La proposition~\ref{isoext3} assure que, quitte \`a remplacer l'anneau~$A$ par l'anneau des entiers d'une extension du corps~$K$, nous pouvons supposer que le polyn\^ome~$P$ est de degr\'e~$1$. La conclusion d\'ecoule alors du th\'eor\`eme~\ref{recanneau3}.
\end{proof}

\bigskip

En proc\'edant de m\^eme, nous pouvons \'etendre les r\'esultats dont nous disposons concernant les uniformisantes forte, le bord analytique des voisinages ou le prolongement analytique. Ces r\'esultats d\'ecoulent du corollaire \ref{recuniforte3depext}, de la proposition \ref{recShilov3depum} et du corollaire \ref{prolan3dep}.

\begin{cor}\label{uniforte3ext}\index{Uniformisante forte!en un point de type 3 d'une fibre extreme@en un point de type $3$ d'une fibre extr\^eme}
Soient~$\m$ un \'el\'ement de~$\Sigma_{f}$ et~$x$ un point de type~$3$ de la fibre extr\^eme~$\tilde{X}_{\m}$. Soient~$\pi$ une uniformisante de l'anneau de valuation discr\`ete~$\Os_{X,x}$ et~$U$ un voisinage du point~$x$ dans~$X$ sur lequel elle est d\'efinie. Il existe un syst\`eme fondamental~$\Vs$ de voisinages compacts et connexes du point~$x$ dans~$U$ tel que, pour tout \'el\'ement~$V$ de~$\Vs$, la fonction~$\pi$ est une uniformisante forte de l'anneau~$\Os_{X,x}$ sur~$V$.
\end{cor}

\begin{cor}\label{Shilov3ext}\index{Bord analytique!au voisinage d'un point de type 3 d'une fibre extreme@au voisinage d'un point de type $3$ d'une fibre extr\^eme}
Soit~$\m$ un \'el\'ement de~$\Sigma_{f}$. Tout point de type~$3$ de la fibre extr\^eme~$\tilde{X}_{\m}$ poss\`ede un syst\`eme fondamental de voisinages compacts, connexes et spectralement convexes qui poss\`edent un bord de Shilov fini et alg\'ebriquement trivial.
\end{cor}

\begin{cor}\label{prolan3ext}\index{Prolongement analytique!au voisinage d'un point de type 3!d'une fibre extreme@d'une fibre extr\^eme}
Soient~$\m$ un \'el\'ement de~$\Sigma_{f}$. Le principe du prolongement analytique vaut au voisinage de tout point de type~$3$ de la fibre extr\^eme~$\tilde{X}_{\m}$.
\end{cor}
\index{Point!de type 3!d'une fibre extreme@d'une fibre extr\^eme|)}

\subsection{Fibre centrale}\index{Point!de type 3!de la fibre centrale|(}

\'Etudions, maintenant, les points de type~$3$ de la fibre centrale. Nous m\`enerons le raisonnement en suivant les m\^emes \'etapes que dans le cas des fibres extr\^emes. Nous commencerons donc par pr\'eciser le r\'esultat de changement de bases obtenu \`a la proposition \ref{isocentral}. Soit $Q(T)$ un polyn\^ome irr\'eductible de $K[T]$. Quel que soit $r\in\of{[}{0,1}{]}$, notons~$\eta_{Q,r}$ le point de la fibre~$X_{0}$ associ\'e \`a la valeur absolue 
$$\begin{array}{ccc}
A[T] & \to & \R_{+}\\
F(T) & \mapsto & r^{v_{Q(T)}(F(T))}
\end{array},$$
o\`u $v_{Q(T)}$ d\'esigne la valuation $Q(T)$-adique de $K[T]$. Pour $\alpha\in K$ et $r\in\of{[}{0,1}{]}$, nous notons 
$$\eta_{\alpha,r} = \eta_{T-\alpha,r}.$$
Pour~$r\in\of{[}{0,1}{]}$, nous notons encore.
$$\eta_{r}=\eta_{0,r}=\eta_{T,r}.$$
Finalement, pour $r\in\of{[}{1,+\infty}{[}$, nous notons~$\eta_{r}$ le point de la fibre~$X_{0}$ associ\'e \`a la valeur absolue 
$$\begin{array}{ccc}
A[T] & \to & \R_{+}\\
F(T) & \mapsto & r^{-\deg(F(T))}
\end{array}.$$
Nous avons ainsi d\'ecrit tous les points de la fibre extr\^eme~$X_{0}$ (\emph{cf.}~\ref{droitetrivval} pour la classification, avec d\'emonstration, des points de la droite analytique sur un corps trivialement valu\'e quelconque). Les points de type~$3$ sont ceux pour lesquels le nombre r\'eel~$r$ est diff\'erent de~$0$ et de~$1$.

Nous noterons $x=\eta_{Q,0}$ le point rigide de la fibre $X_{0}$ d\'efini par l'\'equation
$$Q(T)(x)=0.$$
D'apr\`es la proposition \ref{isocentral}, il existe une extension finie $K'$ de $K$, un point~$x'$ de $X' = \E{1}{A'}$, o\`u $A'$ d\'esigne l'anneau des entiers de $K'$, rationnel dans sa fibre, tel que le morphisme 
$$\psi : \E{1}{A'} \to \E{1}{A}$$
envoie le point $x'$ sur le point $x$ et induise un isomorphisme d'un voisinage de $x'$ sur un voisinage de $x$. Notons $\beta$ l'\'el\'ement de $K'$ qui correspond au point $x'$. Remarquons que, pour tout \'el\'ement~$r$ de l'intervalle $\of{[}{0,1}{]}$, nous avons
$$\psi(\eta_{\beta,r}) = \eta_{Q,r}.$$

Comme pr\'ec\'edemment, \'enon\c{c}ons un r\'esultat assurant l'existence de voisinages flottants. Consid\'erons la partie ouverte~$Y$ de~$X$ obtenue en enlevant les extr\'emit\'es des branches archim\'ediennes :
$$Y = X\setminus\left( \bigcup_{\sigma\in\Sigma_{\infty}} X_{a_{\sigma}}\right).$$

\begin{lem}\label{voisflotcentral}\index{Voisinages flottants!pour un point de type 3!de la fibre centrale}
Soient~$x\in Y$ et~$\eps\in I_{Y}(x)$. Alors, la partie~$D_{Y}$ est un voisinage de~$(x,\eps)$ dans~$Y\times\R_{+}^*$. 

En particulier, tous les points de~$Y$ ont des voisinages flottants dans~$Y$.
\end{lem} 
\begin{proof}
Puisque~$\eps\in I_{Y}(x)$, le point~$x^\eps$ est un \'el\'ement de~$Y$. Nous avons donc~$|2(x)|^{\eps}<2$. Il existe $\lambda>\eps$ tel que l'on ait $|2(x)|^{\eps}< |2(x)|^{\lambda}<2$. La partie
$$\{y\in Y\, \big|\, |2(y)| < 2^{1/\lambda}\} \times \of{]}{0,\lambda}{[}$$
est alors un voisinage de~$(x,\eps)$ dans~$Y\times\R_{+}^*$. 
\end{proof}

Nous tirons de ce r\'esultat les m\^emes cons\'equences que dans le cas des fibres extr\^emes. Les preuves \'etant similaires, nous ne les d\'etaillerons pas.

\begin{prop}\label{isocentral3}\index{Isomorphisme local!au voisinage d'un point de type 3!de la fibre centrale}
Le morphisme $\psi$ induit un isomorphisme d'un voisinage de 
$$\{\eta_{\beta,r},\, r\in\of{[}{0,1}{[}\} \textrm{ dans } X'$$ 
sur un voisinage de 
$$\{\eta_{Q,r},\, r\in\of{[}{0,1}{[}\} \textrm{ dans } X.$$ 
\end{prop}

\begin{cor}\label{cpa3central}\index{Connexite par arcs au voisinage d'un point@Connexit\'e par arcs au voisinage d'un point!de type 3 de la fibre centrale}
Tout point de type~$3$ de la fibre centrale poss\`ede un syst\`eme fondamental de voisinages connexes par arcs. 
\end{cor}

\begin{cor}\label{ouvert3central}\index{Ouverture au voisinage d'un point!de type 3 de la fibre centrale}
Le morphisme $\pi$ est ouvert en tout point de type~$3$ de la fibre centrale. 
\end{cor}

\begin{cor}\label{corps3central}\index{Anneau local en un point!de type 3!de la fibre centrale}
Soit~$x$ un point de type~$3$ de la fibre centrale. En ce point, l'anneau local~$\Os_{X,x}$ co\"incide avec le corps~$\kappa(x)$, lequel est hens\'elien. 
\end{cor}

\begin{cor}\label{prolan3central}\index{Prolongement analytique!au voisinage d'un point de type 3!de la fibre centrale}
Le principe du prolongement analytique vaut au voisinage de tout point de type~$3$ de la fibre centrale de l'espace~$X$.
\end{cor}

\index{Point!de type 3!de la fibre centrale|)}

\section{Points de type~$2$}\label{pdt2}

Pour compl\'eter notre \'etude de la droite analytique sur un corps de nombres, il nous reste \`a \'etudier les points de type 2 des fibres centrale et extr\^emes. Sur ces fibres, et, de fa\c{c}on g\'en\'erale, sur la droite analytique au-dessus de tout corps trivialement valu\'e, il n'existe qu'un point de type~$2$ : le point de Gau{\ss}.

\subsection{Fibres extr\^emes}\index{Point!de type 2!d'une fibre extreme@d'une fibre extr\^eme|(}

Commen\c{c}ons notre \'etude par les fibres extr\^emes. Soit $\m\in\Sigma_{f}$. Notons $x$ le point de Gau{\ss} de la fibre extr\^eme $\tilde{X}_{\m}$. Nous nous int\'eressons, tout d'abord, aux voisinages du point~$x$. Nous notons~$\hat{A}_{\m}^\times$ l'ensemble des \'el\'ements inversibles de l'anneau~$\hat{A}_{\m}$.

\begin{lem}
Soit $U$ un voisinage de $x$ dans $X$. Alors, il existe un entier~$d\in\N$, des polyn\^omes~$P_{1},\ldots,P_{d}\in \hat{A}_{\m}^\times [T]$ et deux nombres r\'eels~$\alpha,\eps>0$ tels que l'on ait
$$U \supset \bigcap_{1\le i\le d} \left\{ y\in \pi^{-1}(\of{]}{a_{\m}^\alpha,\tilde{a}_{\m}}{]})\, \big|\, 1-\eps < |P_{i}(y)| < 1+\eps \right\}.$$
\end{lem}
\begin{proof}
Remarquons que si le r\'esultat vaut pour un nombre fini de voisinages, il vaut encore pour leur intersection. Par cons\'equent, nous pouvons supposer que le voisinage~$U$ est de la forme
$$U = \left\{ y\in X\, \big|\, s < |P(y)| < t\right\},$$
avec $P\in A[T]$ et $s,t\in\R$. En effet, par d\'efinition de la topologie, tout voisinage du point~$x$ contient une intersection finie de voisinages de cette forme.
 
Supposons, tout d'abord, que $P\ne 0 \mod \m$. Il existe alors~$Q\in\hat{A}_{\m}^\times [T]$, $R\in\hat{A}_{\m}[T]$, avec~$R\ne 0 \mod\m$, et~$p\in\N^*$ tels que
$$P = Q + \pi_{\m}^p\, R.$$
Puisque le point~$x$ appartient \`a~$U$ et que~$P(x)=1$, nous avons~$s<1<t$. Par cons\'equent, il existe~$\eps\in\of{]}{0,1}{[}$ tel que~$s<1-\eps$ et~$t>1+\eps$. Soit~$\alpha>0$ tel que 
$$2|\pi_{\m}|_{\m}^{p\alpha} \le 1-\eps.$$ 
Nous avons alors
$$U \supset \left\{ y\in X\, \big|\, 1-\eps < |Q(y)| < 1+\eps\right\} \cap \left\{ y\in \pi^{-1}(\of{]}{a_{\m}^\alpha,\tilde{a}_{\m}}{]})\, \big|\, 0 < |R(y)| < 2\right\}.$$

Supposons, \`a pr\'esent, que $P= 0 \mod \m$. Il existe alors un polyn\^ome~$Q$ de~$\hat{A}_{\m}[T]$, avec~\mbox{$Q\ne 0 \mod\m$}, et~$p\in\N^*$ tels que
$$P = \pi_{\m}^p\, Q.$$
Puisque le point~$x$ appartient \`a~$U$ et que~$P(x)=0$, nous avons~$s<0<t$ et donc
$$U = \left\{ y\in X\, \big|\, |P(y)| < t\right\}.$$
Soit $\alpha>0$ tel que $2|\pi_{\m}|_{\m}^{p\alpha} \le t$. Nous avons alors
$$U \supset \left\{ y\in \pi^{-1}(\of{]}{a_{\m}^\alpha,\tilde{a}_{\m}}{]})\, \big|\, 0 < |Q(y)| < 2\right\}.$$

On d\'emontre finalement le r\'esultat \`a l'aide d'une r\'eccurence sur le nombre de coefficients non nuls du polyn\^ome~$P$ et en utilisant, \`a chaque \'etape, l'un ou l'autre des r\'esultats pr\'ec\'edents.
\end{proof}

\begin{lem}\label{voisGaussext}
Soit $U$ un voisinage de $x$ dans $X$. Alors, il existe deux entiers~$d,e\in\N$, des polyn\^omes~$P_{1},\ldots,P_{d}$ de~$\hat{A}_{\m}^\times [T]$, deux \`a deux distincts, unitaires, irr\'eductibles et dont l'image dans~$k_{\m}[T]$ est une puissance d'un polyn\^ome irr\'eductible, des polyn\^omes~$Q_{1},\ldots,Q_{e}$ de~$\hat{A}_{\m}^\times [T]$, deux \`a deux distincts, unitaires, irr\'eductibles et dont l'image dans~$k_{\m}[T]$ est une puissance d'un polyn\^ome irr\'eductible et deux nombres r\'eels~$\alpha,\eps>0$ tels que l'on ait
$${\renewcommand{\arraystretch}{1.5}
\begin{array}{cccl}
U & \supset & & \disp \bigcap_{1\le i\le d} \left\{ y\in \pi^{-1}(\of{]}{a_{\m}^\alpha,\tilde{a}_{\m}}{]})\, \big|\, |P_{i}(y)| < 1+\eps \right\}\\
&& \cap & \disp \bigcap_{1\le j\le e} \left\{ y\in \pi^{-1}(\of{]}{a_{\m}^\alpha,\tilde{a}_{\m}}{]})\, \big|\, |Q_{j}(y)| > 1-\eps \right\}.
\end{array}
}$$
\end{lem}
\begin{proof}
Comme pr\'ec\'edemment, si le r\'esultat vaut pour un nombre fini de voisinages, il vaut encore pour leur intersection. D'apr\`es le lemme pr\'ec\'edent, nous pouvons donc supposer que le voisinage~$U$ est de la forme
$$U =  \left\{ y\in \pi^{-1}(\of{]}{a_{\m}^\alpha,\tilde{a}_{\m}}{]})\, \big|\, |P(y)| < 1+\eps \right\}$$
ou
$$U =  \left\{ y\in \pi^{-1}(\of{]}{a_{\m}^\alpha,\tilde{a}_{\m}}{]})\, \big|\, |P(y)| > 1-\eps \right\},$$
o\`u $P$ est un polyn\^ome unitaire \`a coefficients dans~$\hat{A}_{\m}$ et~$\alpha$ et~$\eps$ deux nombres r\'eels strictement positifs. Nous supposerons que nous nous trouvons dans le premier cas. Le second se traite de m\^eme. \'Ecrivons le polyn\^ome~$P$ sous la forme
$$P = P_{1}\cdots P_{d},$$
o\`u~$d\in\N$ et~$P_{1},\ldots,P_{d}$ sont des polyn\^omes \`a coefficients dans~$\hat{A}_{\m}$ unitaires, irr\'eductibles et dont l'image dans~$k_{\m}[T]$ est une puissance d'un polyn\^ome ir\-r\'e\-duc\-ti\-ble. La factorialit\'e de l'anneau $\hat{K}_{\m}[T]$ et le lemme de Hensel assurent l'existence d'une telle d\'ecomposition existe. Soit~$i\in\cn{1}{d}$. Puisque le polyn\^ome~$P_{i}$ est unitaire, il v\'erifie~\mbox{$|P_{i}(x)|=1$}. Par cons\'equent, la partie
$$U_{i}=\left\{ y\in \pi^{-1}(\of{]}{a_{\m}^\alpha,\tilde{a}_{\m}}{]})\, \big|\, |P_{i}(y)| < (1+\eps)^{1/d} \right\}$$
est un voisinage du point~$x$ dans~$X$. L'intersection
$$\bigcap_{1\le i\le d} U_{i}$$
est alors un voisinage de~$x$ dans~$U$ de la forme voulue.
\end{proof}

\begin{prop}\label{sectionGaussext}\index{Voisinages d'un point!de type 2!d'une fibre extreme@d'une fibre extr\^eme}
Soit $U$ un voisinage du point $x$ dans $X$. Alors il existe un voisinage ouvert $W$ de $x$ dans $U$ v\'erifiant les propri\'et\'es suivantes :
\begin{enumerate}[\it i)]
\item la projection $\pi(W)$ est un voisinage connexe par arcs de $\pi(x)=\tilde{a}_{\m}$ dans~$B$ ;
\item la section de Gau{\ss} $\sigma_{G}$ restreinte \`a $\pi(W)$ prend ses valeurs dans $W$ ;
\item pour tout point $b$ de $\pi(W)$, la trace de la fibre $X_{b}$ sur $W$ est connexe par arcs.
\end{enumerate} 
\end{prop}
\begin{proof}
Appliquons le lemme pr\'ec\'edent. Le voisinage~$W$ que l'on obtient v\'erifie les propri\'et\'es demand\'ees. Les deux premi\`eres sont imm\'ediates. Int\'eressons-nous \`a la troisi\`eme. Nous conservons les notations du lemme pr\'ec\'edent. Soit~$\beta$ un \'el\'ement de~$\of{]}{\alpha,+\infty}{]}$. Nous voulons montrer que la trace de la fibre $X_{a_{\m}^\beta}$ sur~$W$ est connexe par arcs. Soit~$i\in\cn{1}{d}$. Par d\'efinition, le polyn\^ome~$P_{i}$ est une puissance d'un polyn\^ome irr\'eductible dans~$\Hs(a_{\m}^\beta)[T]$. On en d\'eduit que la partie
$$\left\{ y\in X_{a_{\m}^\beta}\, \big|\, |P_{i}(y)| < 1+\eps \right\}$$
est connexe par arcs. On l'obtient en effet \`a partir de la droite~$\E{1}{\Hs(a_{\m}^\beta)}$ en coupant l'une des branches partant du point de Gau{\ss}. De m\^eme, quel que soit~\mbox{$j\in\cn{1}{e}$}, la partie
$$\left\{ y\in X_{a_{\m}^\beta}\, \big|\, |P_{j}(y)| > 1-\eps \right\}$$
est connexe par arcs. Puisque la droite analytique~$\E{1}{\Hs(a_{\m}^\beta)}$ a une structure d'arbre, une intersection de parties connexes par arcs l'est encore. On en d\'eduit que la partie~$W\cap X_{a_{\m}^\beta}$ est connexe par arcs.
\end{proof}

Quatre corollaires suivent.

\begin{cor}\label{cpaGaussext}\index{Connexite par arcs au voisinage d'un point@Connexit\'e par arcs au voisinage d'un point!de type 2 d'une fibre extreme@de type $2$ d'une fibre extr\^eme}
Le point de Gau{\ss} de la fibre extr\^eme $\tilde{X}_{\m}$ poss\`ede un syst\`eme fondamental de voisinages connexes par arcs.
\end{cor}

\begin{cor}\label{ouvertGaussext}\index{Ouverture au voisinage d'un point!de type 2 d'une fibre extreme@de type $2$ d'une fibre extr\^eme}
Le morphisme $\pi$ est ouvert au voisinage du point de Gau{\ss} de la fibre extr\^eme $\tilde{X}_{\m}$.
\end{cor}

\begin{cor}\label{prolanGaussext}\index{Prolongement analytique!au voisinage d'un point de type 2!d'une fibre extreme@d'une fibre extr\^eme}
Le principe du prolongement analytique vaut au voisinage du point de Gau{\ss} de la fibre extr\^eme $\tilde{X}_{\m}$.
\end{cor}
\begin{proof}
Soient~$U$ un voisinage du point~$x$ dans~$X$ et~$f$ un \'el\'ement de~$\Os(U)$ dont l'image dans l'anneau local~$\Os_{X,x}$ n'est pas nulle. Consid\'erons alors un voisinage ouvert~$W$ du point~$x$ contenu dans~$U\cap\pi^{-1}(\of{]}{a_{0},\tilde{a}_{\m}}{]})$ et v\'erifiant les propri\'et\'es de la proposition \ref{sectionGaussext}. 

Posons $W_{+}=W\cap \pi^{-1}(\of{]}{a_{0},\tilde{a}_{\m}}{[})$. Puisque le morphisme~$\pi$ est ouvert au voisinage du point~$x$, il existe un point interne~$y$ de~$W_{+}$ tel que l'image de la fonction~$f$ n'est pas nulle dans l'anneau local~$\Os_{X,y}$. Par choix de~$W$, l'ouvert~$W_{+}$  est connexe et le principe du prolongement analytique y vaut donc, d'apr\`es la proposition \ref{prolinterne}. On en d\'eduit que, pour tout point~$z$ de~$W_{+}$, l'image de la fonction~$f$ dans l'anneau local~$\Os_{X,z}$ n'est pas nulle.

Posons $W_{0}=W\cap \tilde{X}_{\m}$. Soit~$z$ un point de~$W_{0}\setminus\{x\}$. D'apr\`es le th\'eor\`eme \ref{toporig} et le corollaire \ref{ouvert3ext}, le morphisme~$\pi$ est ouvert au voisinage du point~$z$. Par cons\'equent, tout voisinage du point~$z$ contient un \'el\'ement de~$W_{+}$ et l'image de la fonction~$f$ dans l'anneau local~$\Os_{X,z}$ ne peut pas \^etre nulle. Ceci conclut la preuve.
\end{proof}

\begin{cor}\label{prolanbrancheext}
Le principe du prolongement analytique vaut sur tout ouvert connexe de l'espace~$\pi^{-1}(\of{]}{a_{0},\tilde{a}_{\m}}{]})$.
\end{cor}
\begin{proof}
Ce r\'esultat d\'ecoule des corollaires \ref{prolanrig} et \ref{prolan3ext}, de la proposition \ref{recprolinterne} et du lemme \ref{lemprolan}.
\end{proof}

Int\'eressons-nous, \`a pr\'esent, \`a l'anneau local.

\begin{prop}\label{avd2ext}\index{Anneau local en un point!de type 2!d'une fibre extreme@d'une fibre extr\^eme}
Soit~$\m\in\Sigma_{f}$. Notons~$x$ le point de Gau{\ss} de la fibre extr\^eme~$\tilde{X}_{\m}$. L'anneau local $\Os_{X,x}$ est un anneau de valuation discr\`ete d'id\'eal maximal~$\m\,\Os_{X,x}$. Son corps r\'esiduel~$\kappa(x)$ est complet, et donc isomorphe \`a $\Hs(x) = \tilde{k}_{\m}(T)$.
\end{prop} 
\begin{proof}
Nous allons d\'efinir une valuation discr\`ete~$v$ sur l'anneau local~$\Os_{X,x}$. Soit~$f$ un \'el\'ement de~$\Os_{X,x}$. Il existe un voisinage~$U$ de~$x$ dans~$X$ sur lequel la fonction~$f$ est d\'efinie. Pour~$r\in\of{[}{0,1}{]}$, nous noterons simplement~$\eta_{r}$ le point~$\eta_{r}$ de la fibre extr\^eme~$\tilde{X}_{\m}$. La trace de la partie~$U$ sur la fibre extr\^eme~$\tilde{X}_{\m}$ est un voisinage du point~$x=\eta_{1}$ dans cette fibre. Par cons\'equent, il existe~$R\in\of{]}{0,1}{[}$ tel que, quel que soit~$r\in\of{[}{R,1}{]}$, on ait $\eta_{r} \in U.$
D'apr\`es la proposition~\ref{avd3ext}, l'anneau local~$\Os_{X,\eta_{R}}$ est un anneau de valuation discr\`ete. Notons~$v_{R}$ la valuation sur cet anneau. Nous posons alors
$$v(f)=v_{R}(f) \in \N\cup\{+\infty\}.$$
La proposition~\ref{flot} nous assure que cette quantit\'e ne d\'epend pas du nombre r\'eel~$R$ choisi.

Les deux propri\'et\'es suivantes sont imm\'ediates : quels que soient~$f$ et~$g$ dans~$\Os_{X,x}$, nous avons
\begin{enumerate}
\item $v(f+g)\ge \min(v(f),v(g))$ ;
\item $v(fg) = v(f)+v(g)$.
\end{enumerate}

Nous avons \'egalement $v(0)=+\infty$. Montrons que seule la fonction nulle satisfait cette \'egalit\'e. Soit~$f\in\Os_{X,x}$ telle que~$v(f)=+\infty$. Soit~$U$ un voisinage ouvert de~$x$ dans~$X$ sur lequel la fonction~$f$ est d\'efinie. D'apr\`es la proposition~\ref{sectionGaussext}, quitte \`a restreindre~$U$, nous pouvons supposer qu'il v\'erifie les propri\'et\'es suivantes :
\begin{enumerate}[\it i)]
\item la projection $\pi(U)$ est un voisinage connexe par arcs de $\pi(x)=\tilde{a}_{\m}$ dans~$B$ ;
\item la section de Gau{\ss} $\sigma_{G}$ restreinte \`a $\pi(U)$ prend ses valeurs dans $U$ ;
\item pour tout point $b$ de $\pi(U)$, la trace de la fibre $X_{b}$ sur $U$ est connexe par arcs.
\end{enumerate} 
Soit~$R\in\of{]}{0,1}{[}$ tel que, quel que soit~$r\in\of{[}{R,1}{]}$, on ait $\eta_{r} \in U$. Par d\'efinition de~$v$, nous avons~$v_{R}(f)=+\infty$. Par cons\'equent, l'image de la fonction~$f$ dans l'anneau local~$\Os_{X,\eta_{R}}$ est nulle. Il existe donc un voisinage ouvert~$V$ du point~$\eta_{R}$ dans~$U$ tel que la fonction~$f$ soit nulle sur~$V$. D'apr\`es le corollaire~\ref{ouvert3ext}, la partie~$V_{0}=\pi(V)$ est un voisinage du point extr\^eme~$\tilde{a}_{\m}$ dans~$B$. Soit~$c\in V_{0}$. La fonction~$f$ est nulle sur un l'ouvert $X_{c}\cap V$ de $X_{c}\cap U$. Comme ce dernier espace est normal et connexe, la fonction~$f$ y est identiquement nulle. Finalement, la fonction~$f$ est nulle sur~$U\cap X_{V_{0}}$ et donc dans l'anneau local~$\Os_{X,x}$.

La propri\'et\'e que nous venons de d\'emontrer jointe \`a la propri\'et\'e 2 impose \`a l'anneau local~$\Os_{X,x}$ d'\^etre int\`egre. Consid\'erons son corps des fractions~$L$. L'application~$v$ se prolonge alors en une valuation discr\`ete sur le corps~$L$. Pour parvenir \`a nos fins, il nous reste \`a montrer les deux \'egalit\'es
$$\Os_{X,x} = \{f\in L\, |\, v(f)\ge 0\}$$
et
$$\m\,\Os_{X,x} = \{f\in L\, |\, v(f)> 0\}.$$
Remarquons que la seconde \'egalit\'e d\'ecoule de la premi\`ere et du fait que le g\'en\'erateur~$\pi_{\m}$ de l'id\'eal maximal~$\m$ de~$A$ a pour valuation~$v(\pi_{\m})=1$. D'autre part, pour d\'emontrer la premi\`ere \'egalit\'e, il nous suffit de montrer que tout \'el\'ement~$f$ de~$\Os_{X,x}$ v\'erifiant~$v(f)=0$ est inversible dans l'anneau~$\Os_{X,x}$. Ce r\'esultat se d\'emontre facilement en utilisant les propri\'et\'es du flot (\emph{cf.} proposition~\ref{flot}). En effet, soit~$f$ un \'el\'ement de~$\Os_{X,x}$ v\'erifiant $v(f)=0$. Il existe un nombre r\'eel~$R\in\of{]}{0,1}{[}$ v\'erifiant les propri\'et\'es habituelles tel que l'on ait~$v_{R}(f)=0$. On en d\'eduit que la fonction~$f$ est inversible dans l'anneau local~$\Os_{X,\eta_{R}}$ et donc que~$|f(\eta_{R})|\ne 0$. La proposition~\ref{flot} nous assure alors que l'on a
$$|f(x)|=|f(\eta_{R})|^0=1.$$ 
On en d\'eduit que la fonction~$f$ est inversible dans l'anneau local~$\Os_{X,x}$.

Le corps valu\'e~$\Hs(x)$ est isomorphe au corps~$\tilde{k}_{m}(T)$ muni de la valuation triviale. Le corps valu\'e~$\kappa(x)$ qui en est un sous-corps est donc complet.
\end{proof}

\bigskip




\begin{lem}\label{bordP}
Soit~$\m$ un \'el\'ement de~$\Sigma_{f}$. Soient $r\in\R_{+}^*\setminus\{1\}$ et~$P(T)$ un polyn\^ome unitaire, non constant, \`a coefficients dans~$\hat{A}_{\m}$. Quel que soit~$\eps>0$, posons
$$W_{\eps} = \left\{y\in \pi^{-1}(\of{[}{a_{\m}^{\eps},\tilde{a}_{\m}}{]})\, \big|\, |P(T)(y)|=r\right\}.$$ 
Il existe $\eps_{0}>0$ tel que, pour tout $\eps\ge\eps_{0}$, le compact rationnel~$W_{\eps}$ poss\`ede un bord analytique fini, alg\'ebriquement trivial et contenu dans la fibre~$X_{a_{\m}^\eps}$.
\end{lem}
\begin{proof}
Soient~$K'$ une extension finie de~$K$, $A'$ l'anneau de ses entiers et~$\m'$ un id\'eal maximal de~$A'$ divisant l'id\'eal maximal~$\m$ de~$A$. En utilisant la surjectivit\'e du morphisme $\E{1}{A'} \to \E{1}{A}$, on montre facilement que si le r\'esultat \'enonc\'e vaut en rempla\c{c}ant~$A$ par~$A'$ et~$\m$ par~$\m'$, alors il vaut dans la form\'e originale. Par cons\'equent, quitte \`a remplacer~$K$ par une extension finie~$K'$ bien choisie, nous pouvons supposer que le polyn\^ome~$P(T)$ est scind\'e dans~$\hat{K}_{\m}$. Puisqu'il est unitaire et \`a coefficients entiers, ses racines sont enti\`eres et il existe donc $t\in\N^*$ et $\alpha_{1},\ldots,\alpha_{t}\in\hat{A}_{\m}$ tels que
$$P(T) = \prod_{i=1}^t (T-\alpha_{i}).$$

Remarquons \'egalement qu'il suffit de montrer que le compact rationnel~$W_{\eps}$ poss\`ede un bord analytique qui est contenu dans la fibre~$X_{a_{\m}^\eps}$. Les autres conditions d\'ecoulent alors du corollaire \ref{Shilovrationnelinterne}.

\bigskip
 
Supposons, tout d'abord, que $r>1$. Soit~$\eps>0$. Soit~$y$ un point de~$W_{\eps}$. Puisque $|P(T)(y)=r$, il existe un \'el\'ement~$i$ de~$\cn{1}{t}$ tel que $|(T-\alpha_{i})(y)|\ge u$. Quel que soit~$j\ne i$, nous avons $|(\alpha_{i}-\alpha_{j})(y)|\le 1<u$ et donc $|(T-\alpha_{j})(y)|=u$. Par cons\'equent, nous avons \'egalement 
$$|(T-\alpha_{i})(y)|=u.$$
R\'eciproquement, l'on montre que
$$W_{\eps} = \left\{y\in \pi^{-1}(\of{[}{a_{\m}^{\eps},\tilde{a}_{\m}}{]})\, \big|\, |(T-\alpha_{i})(y)|=u\right\}.$$ 
Le r\'esultat d\'ecoule alors de la proposition \ref{Shilovcouronneum} et des descriptions explicites \'etablies au num\'ero \ref{borddeShilovbase}.

\bigskip

Supposons, \`a pr\'esent, que~$r<1$. Posons
$$D=\left\{(i,j)\in\cn{1}{t}^2\, \big|\, |\alpha_{i}-\alpha_{j}|_{\m}<1\right\}.$$
Il existe $\eps_{0}>0$ tel que, pour tout couple~$(i,j)$ de~$D$, nous ayons
$$|\alpha_{i}-\alpha_{j}|_{\m}^{\eps_{0}}<r.$$
Soit~$\eps\ge\eps_{0}$. Soit~$y$ un point de~$W_{\eps}$. Puisque $|P(T)(y)=r<1$, il existe un \'el\'ement~$i$ de~$\cn{1}{t}$ tel que $|(T-\alpha_{i})(y)|<1$. Posons
$$g_{i} = \left\{j\ne i\, \big|\, |\alpha_{i}-\alpha_{j}|_{\m}=1\right\}$$ 
et
$$p_{i} = \left\{j\ne i\, \big|\, |\alpha_{i}-\alpha_{j}|_{\m}<1\right\}.$$ 
Remarquons que, par d\'efinition de~$\eps_{0}$, pour tout \'el\'ement~$i$ de~$p_{i}$, nous avons
$$ |\alpha_{i}-\alpha_{j}|_{\m}^\eps<r.$$

Supposons, par l'absurde, que $|(T-\alpha_{i})(y)|<r$. Alors, quel que soit~$j\in g_{i}$, nous avons $|(T-\alpha_{j})(y)|=1$ et, quel que soit~$j\in p_{i}$, nous avons $|(T-\alpha_{j})(y)|<r$. On en d\'eduit que
$$|P(T)(y)|<r^{\sharp g_{i} +1} < r,$$
ce qui est impossible. 

Par cons\'equent, nous avons $|(T-\alpha_{i})(y)|\ge r$. On en d\'eduit que, quel que soit~$j\in g_{i}$, nous avons $|(T-\alpha_{j})(y)|=1$ et, quel que soit~$j\in p_{i}$, $|(T-\alpha_{j})(y)|=|(T-\alpha_{i})(y)|$. Par cons\'equent, nous avons
$$|(T-\alpha_{i})(y)|=r^{1/(\sharp g_{i} +1)} \in\of{[}{r,1}{[}.$$ 
R\'eciproquement, l'on montre que, si~$y$ est un point de $\pi^{-1}(\of{[}{a_{\m}^{\eps},\tilde{a}_{\m}}{]})$ tel que $|(T-\alpha_{i})(y)|=r^{1/(\sharp g_{i} +1)}$, alors~$y$ appartient \`a~$W_{\eps}$.

Finalement, nous avons montr\'e que
$$W_{\eps} = \bigcup_{1\le i\le t} \left\{y\in \pi^{-1}(\of{[}{a_{\m}^{\eps},\tilde{a}_{\m}}{]})\, \big|\, |(T-\alpha_{i})(y)|=r^{1/(\sharp g_{i} +1)}\right\}.$$
Le r\'esultat d\'ecoule alors de la proposition \ref{Shilovcouronneum} et des descriptions explicites \'etablies au num\'ero \ref{borddeShilovbase}.
\end{proof}

\begin{cor}\label{Shilov2ext}\index{Bord analytique!au voisinage d'un point de type 2 d'une fibre extreme@au voisinage d'un point de type $2$ d'une fibre extr\^eme}
Soit~$\m$ un \'el\'ement de~$\Sigma_{f}$. Le point de {Gau\ss} de la fibre extr\^eme~$\tilde{X}_{\m}$ poss\`ede un syst\`eme fondamental de voisinages compacts, connexes et spectralement convexes qui poss\`edent un bord analytique fini et alg\'ebriquement trivial.
\end{cor}
\begin{proof}
Notons~$x$ le point de {Gau\ss} de la fibre extr\^eme~$\tilde{X}_{\m}$. Soit~$U$ un voisinage du point~$x$ dans~$X$. D'apr\`es le lemme \ref{voisGaussext}, il existe deux entiers~$d,e\in\N$, des polyn\^omes~$P_{1},\ldots,P_{d}$ de~$\hat{A}_{\m}^\times [T]$, deux \`a deux distincts, unitaires, irr\'eductibles et dont l'image dans~$k_{\m}[T]$ est une puissance d'un polyn\^ome irr\'eductible, des polyn\^omes~$Q_{1},\ldots,Q_{e}$ de~$\hat{A}_{\m}^\times [T]$, deux \`a deux distincts, unitaires, irr\'eductibles et dont l'image dans~$k_{\m}[T]$ est une puissance d'un polyn\^ome irr\'eductible et deux nombres r\'eels~$\alpha>0$ et $\eps\in\of{]}{0,1}{[}$ tels que le voisinage du point~$x$ d\'efini par
$${\renewcommand{\arraystretch}{1.5}
\begin{array}{cccl}
V & = & & \disp \bigcap_{1\le i\le d} \left\{ y\in \pi^{-1}(\of{[}{a_{\m}^\alpha,\tilde{a}_{\m}}{]})\, \big|\, |P_{i}(y)| \le 1+\eps \right\}\\
&& \cap & \disp \bigcap_{1\le j\le e} \left\{ y\in \pi^{-1}(\of{[}{a_{\m}^\alpha,\tilde{a}_{\m}}{]})\, \big|\, |Q_{j}(y)| \ge 1-\eps \right\}
\end{array}}$$
soit contenu dans~$U$. Le voisinage~$V$ est compact, connexe (par le m\^eme raisonnement que dans la preuve de la proposition \ref{sectionGaussext}) et spectralement convexe (d'apr\`es la proposition \ref{stabilitespconvexe}). Notons
$${\renewcommand{\arraystretch}{1.5}
\begin{array}{cccl}
W & = & & \disp \bigcap_{1\le i\le d} \left\{ y\in \pi^{-1}(\of{[}{a_{\m}^\alpha,\tilde{a}_{\m}}{]})\, \big|\, |P_{i}(y)| = 1+\eps \right\}\\
&& \cap & \disp \bigcap_{1\le j\le e} \left\{ y\in \pi^{-1}(\of{[}{a_{\m}^\alpha,\tilde{a}_{\m}}{]})\, \big|\, |Q_{j}(y)| = 1-\eps \right\}.
\end{array}}$$
D'apr\`es la proposition \ref{ShilovBerko}, cette partie compacte contient le bord de Shilov de l'intersection de~$V$ avec chaque fibre au-dessus de $\of{[}{a_{\m}^\alpha,\tilde{a}_{\m}}{]}$. C'est donc un bord analytique de~$V$. Quitte \`a augmenter~$\alpha$, d'apr\`es le lemme \ref{bordP}, le compact rationnel~$V$ poss\`ede un bord analytique fini et alg\'ebriquement trivial.
\end{proof}

\begin{cor}\label{uniforte2ext}\index{Uniformisante forte!en un point de type 2 d'une fibre extreme@en un point de type $2$ d'une fibre extr\^eme}
Soit~$\m$ un \'el\'ement de~$\Sigma_{f}$. Notons~$x$ le point de {Gau\ss} de la fibre extr\^eme~$\tilde{X}_{\m}$. Soient~$\pi$ une uniformisante de l'anneau de valuation discr\`ete~$\Os_{X,x}$ et~$U$ un voisinage du point~$x$ dans~$X$ sur lequel elle est d\'efinie. Il existe un syst\`eme fondamental~$\Vs$ de voisinages compacts et connexes du point~$x$ dans~$U$ tel que, pour tout \'el\'ement~$V$ de~$\Vs$, la fonction~$\pi$ est une uniformisante forte de l'anneau~$\Os_{X,x}$ sur~$V$.
\end{cor}
\begin{proof}
Soient $d,e\in\N$, $P_{1},\ldots,P_{d}\in\hat{A}_{\m}^\times [T]$, deux \`a deux distincts, unitaires, irr\'eductibles et dont l'image dans~$k_{\m}[T]$ est une puissance d'un polyn\^ome irr\'eductible, $Q_{1},\ldots,Q_{e}\in\hat{A}_{\m}^\times [T]$, deux \`a deux distincts, unitaires, irr\'eductibles et dont l'image dans~$k_{\m}[T]$ est une puissance d'un polyn\^ome irr\'eductible, $\alpha>0$ et $\eps\in\of{]}{0,1}{[}$. D\'efinissons un voisinage du point~$x$ dans~$X$ par
$${\renewcommand{\arraystretch}{1.5}
\begin{array}{cccl}
V & = & & \disp \bigcap_{1\le i\le d} \left\{ y\in \pi^{-1}(\of{[}{a_{\m}^\alpha,\tilde{a}_{\m}}{]})\, \big|\, |P_{i}(y)| \le 1+\eps \right\}\\
&& \cap & \disp \bigcap_{1\le j\le e} \left\{ y\in \pi^{-1}(\of{[}{a_{\m}^\alpha,\tilde{a}_{\m}}{]})\, \big|\, |Q_{j}(y)| \ge 1-\eps \right\}.
\end{array}}$$
Nous avons montr\'e dans le corollaire pr\'ec\'edent que, si~$\alpha$ est assez grand, ce que nous supposerons d\'esormais, la partie~$V$ est compacte et connexe et poss\`ede un bord analytique~$\Gamma$ fini et alg\'ebriquement trivial. D'apr\`es le lemme \ref{voisGaussext}, il suffit de montrer que la fonction~$\pi_{\m}$ est une uniformisante forte de l'anneau~$\Os_{X,x}$ sur~$V$. Remarquons que la fonction~$\pi_{\m}$ ne s'annule pas sur l'ensemble~$\Gamma$. Posons
$$C=\|\pi_{\m}^{-1}\|_{\Gamma}.$$

Soit~$f$ un \'el\'ement de~$\Os(V)$ dont l'image dans~$\Hs(x)$ est nulle. Puisque l'espace analytique~$\tilde{X}_{\m}$ est normal, que la partie~$V\cap\tilde{X}_{\m}$ est connexe et que l'anneau local~$\Os_{\tilde{X}_{\m},x}$ est un corps, nous avons
$$\forall y\in V\cap\tilde{X}_{\m},\, f(y)=0.$$
D'apr\`es la proposition \ref{avd2ext}, la fonction~$f$ est multiple de~$\pi_{\m}$ au voisinage du point~$x$. D'apr\`es le corollaire \ref{avd3ext}, elle l'est \'egalement au voisinage de tout point de type~$3$ de $V\cap \tilde{X}_{\m}$.

Soit~$y$ un point de $V\cap \tilde{X}_{\m}$ qui n'est pas de type~$2$ ou~$3$. C'est alors un point rigide de~$\tilde{X}_{\m}$. La proposition \ref{isorig} nous permet de supposer que c'est un point rationnel. En utilisant le d\'eveloppement en s\'erie de la fonction~$f$ donn\'e par le corollaire \ref{descratext} et le r\'esultat concernant les points de type~$3$ voisins, l'on montre que la fonction~$f$ est multiple de~$\pi_{\m}$ au voisinage du point~$y$. 

Soit~$y$ un point de~$V$ qui n'appartient pas \`a la fibre extr\^eme~$\tilde{X}_{\m}$. La fonction~$f$ est multiple de~$\pi_{\m}$ au voisinage du point~$y$, puisque la fonction~$\pi_{\m}$ est inversible au voisinage de ce point.

La connexit\'e de~$V$ et le principe du prolongement analytique (\emph{cf.} corollaire \ref{prolanbrancheext}) assurent qu'il existe un \'el\'ement~$g$ de~$\Os(V)$ tel que l'on ait l'\'egalit\'e
$$f = \pi_{\m}\, g \textrm{ dans } \Os(V).$$ 
En outre, nous avons
$$\|g\|_{V} = \|g\|_{\Gamma} = \max_{\gamma\in\Gamma} (|(\pi_{\m}^{-1}\, f)(\gamma)|) \le C\, \|f\|_{V}.$$
\end{proof}

\index{Point!de type 2!d'une fibre extreme@d'une fibre extr\^eme|)}

\subsection{Fibre centrale}\index{Point!de type 2!de la fibre centrale|(}

Int\'eressons-nous, \`a pr\'esent, au point de Gau{\ss} de la fibre centrale. Comme pr\'ec\'edemment, nous commen\c{c}ons par \'etudier ses voisinages. C'est un probl\`eme bien plus d\'elicat que pour les fibres extr\^emes.

\begin{lem}\label{cparho}
Soit $(k,|.|)$ un corps ultram\'etrique complet. Soit un polyn\^ome $P(T) = \sum_{i=0}^d a_{i}\, T^i \in k[T]$, avec $d\in\N^*$, quel que soit $i\in\cn{0}{d-1}$, $a_{i}\in k$ et $a_{d}\in k^*$. Posons
$$\rho = \max_{0\le i\le d-1} \left( \left|\frac{a_{i}}{a_{d}}\right|^{\frac{1}{d-i}} \right).$$
Soient $\lambda,\mu \in\R$ v\'erifiant la condition $\mu > |a_{d}|\, \rho^d$.
Alors la partie de $\E{1}{k}$ d\'efinie par 
$$U =  \left\{x\in \E{1}{k}\,\big|\, \lambda < |P(x)| < \mu \right\}$$
est connexe par arcs.
\end{lem}
\begin{proof}
Soit $k'$ un corps alg\'ebriquement clos et maximalement complet contenant $k$. Puisque le morphisme de changement de bases $\E{1}{k'} \to \E{1}{k}$ est continu et surjectif, quitte \`a remplacer $k$ par $k'$, nous pouvons supposer que le corps $k$ est alg\'ebriquement clos et maximalement complet. Il existe alors $\alpha_{1},\ldots,\alpha_{d} \in k$ tels que 
$$P(T) = a_{d}\, \prod_{i=1}^d (T-\alpha_{i}).$$

D'apr\`es \cite{BGR}, proposition 3.1.2.1, quel que soit $i\in\cn{1}{d}$, nous avons 
$$|\alpha_{i}| \le \rho.$$
Soit $r \ge \rho$ v\'erifiant la condition $ \lambda < |a_{d}|\, r^d < \mu$. Alors, nous avons 
$$|P(\eta_{r})| = |a_{d}|\, \prod_{i=1}^d |(T-\alpha_{i})(\eta_{r})| = |a_{d}|\, r^d.$$
Par cons\'equent, le point $\eta_{r}$ appartient \`a $U$.

Soit $x$ un point de $U$. Puisque $k$ est maximalement complet, il existe $\beta\in k$ et $s\in\R_{+}$ tels que $x=\eta_{\beta,s}$. Soit $i\in\cn{1}{d}$. Nous avons $T-\alpha_{i} = (T-\beta)+(\beta-\alpha_{i})$ et donc
$$|(T-\alpha_{i})(\eta_{\beta,s})| = \max(s,|\beta-\alpha_{i}|).$$

Supposons que $|\beta|\le r$. Consid\'erons le chemin injectif $l$ trac\'e sur $\E{1}{k}$ d\'efini par
$$\begin{array}{ccc}
\of{[}{0,1}{]} & \to & \disp \E{1}{k}\\
t & \mapsto & \eta_{\beta, tr + (1-t)s}
\end{array}.$$
Il joint le point $\eta_{\beta,s}$ au point $\eta_{\beta,r}=\eta_{r}$. Si $s$ est inf\'erieur \`a $r$, alors, lorsque l'on parcourt~$l$, la fonction $|P|$ cro\^{\i}t de $|P(\eta_{\beta,s})|$ \`a $|P(\eta_{r})|$. En particulier, le chemin reste dans $U$. Il en est de m\^eme si $s > r$.

Supposons, \`a pr\'esent, que $|\beta|>r$. Si $s\ge |\beta|$, alors $\eta_{\beta,s} = \eta_{0,s}$ et nous sommes ramen\'es au cas pr\'ec\'edent. Supposons donc que $s < |\beta|$. Quel que soit $i\in\cn{1}{d}$, nous avons
$$|(T-\alpha_{i})(\eta_{\beta,s})| = \max(s,|\beta-\alpha_{i}|) = \max(s,|\beta|) = |\beta|.$$
Le long du chemin $l'$, joignant le point $\eta_{\beta,s}$ au point $\eta_{\beta,|\beta|}$, d\'efini par 
$$\begin{array}{ccc}
\of{[}{0,1}{]} & \to & \E{1}{k}\\
t & \mapsto & \eta_{\beta, t|\beta| + (1-t)s}
\end{array},$$
la fonction $|P|$ est constante. Le chemin $l'$ est donc trac\'e sur $U$. Nous sommes donc ramen\'es au cas du point $\eta_{\beta,|\beta|} = \eta_{0,|\beta|}$ que nous avons trait\'e pr\'ec\'edemment. Nous pouvons donc joindre le point $\eta_{\beta,s}$ au point $\eta_{r}$ par un chemin trac\'e sur~$U$.
\end{proof}

\begin{lem}\label{cpaHardt}
Soit $(k,|.|)$ un corps archim\'edien complet. Soient~$d\in\N$ et~$P_{1},\ldots,P_{d}$ des polyn\^omes \`a coefficients dans~$k$. Alors, il existe~$S,T\in\R$ tels que, quels que soient~\mbox{$s_{1},\ldots,s_{d}\in\of{[}{0,S}{[}$} et~\mbox{$t_{1},\ldots,t_{d}\in\of{]}{T,+\infty}{[}$}, la partie de~$\E{1}{k}$ d\'efinie par
$$\bigcap_{1\le j\le d} \left\{ z\in\E{1}{k}\, \big|\, s_{j} < |P_{j}(z)| < t_{j}\right\}$$
est connexe par arcs.  
\end{lem}
\begin{proof}
Consid\'erons un plongement du corps~$k$ dans le corps~$\C$. Nous munissons~$\C$ de l'unique valeur absolue qui \'etend celle de~$k$. Le morphisme $\E{1}{\C} \to \E{1}{k}$ induit par le plongement pr\'ec\'edent \'etant continu et surjectif, nous pouvons supposer que $k=\C$.

Nous pouvons supposer qu'aucun des polyn\^omes~$P_{i}$, avec~$i\in\cn{1}{d}$, n'est nul. Notons $E$ l'ensemble des \'el\'ements~$(x,y,s_{1},\ldots,s_{d},t_{1},\ldots,t_{d})$ de $\R^2\times \R_{+}^{2d}$ qui v\'erifient la condition suivante :
$$\forall j\in\cn{1}{d},\, s_{j}<|P_{j}(x+iy)|^2<t_{j}.$$
C'est un ensemble semi-alg\'ebrique r\'eel. Consid\'erons \'egalement l'application 
$$p : E \to \of{[}{0,1}{]}^{2d}$$ 
qui \`a tout \'el\'ement $u=(x,y,s_{1},\ldots,s_{d},t_{1},\ldots,t_{d})$ de $E$ associe
$$p(u) = \left(s_{1},\ldots,s_{d},\frac{t_{1}}{1+t_{1}},\ldots,\frac{t_{d}}{1+t_{d}}\right).$$
Cette application est semi-alg\'ebrique r\'eelle et continue. D'apr\`es le th\'eor\`eme de Hardt (\emph{cf.}\cite{BCR}, th\'eor\`eme 9.3.1), il existe une partition $(T_{1},\ldots,T_{r})$, avec~$r\in\N$, de~$\of{[}{0,1}{]}^{2d}$ en parties semi-alg\'ebrique telle que, quel que soit~$k\in\cn{1}{r}$, il existe un ensemble semi-alg\'ebrique~$F_{k}$ et un hom\'eomorphisme semi-alg\'ebrique
$$\theta_{k} : T_{k}\times F_{k} \xrightarrow[]{\sim} p^{-1}(T_{k})$$ 
tel que l'application~$p\circ\theta_{j}$ soit la projection~$T_{k}\times F_{k} \to T_{k}$. Notons~$v$ le point $(0,\ldots,0,1,\ldots,1)$ de~$\of{[}{0,1}{]}^{2d}$. Pour parvenir au r\'esulat souhait\'e, il suffit de montrer que le point~$v$ poss\`ede un voisinage dans~$\of{[}{0,1}{]}^{2d}$ au-dessus duquel les fibres de l'application~$p$ sont connexes. Autrement dit, il suffit de montrer que pour tout indice~$k\in\cn{1}{r}$ tel que le point~$v$ soit adh\'erent \`a la partie~$T_{k}$, la partie~$F_{k}$ est connexe.

Soit~$k\in\cn{1}{r}$ tel que le point~$v$ soit adh\'erent \`a la partie~$T_{k}$. D'apr\`es le lemme de s\'election des courbes (\emph{cf.}~\cite{BCR}, th\'eor\`eme~2.5.5), il existe une fonction semi-alg\'ebrique continue
$$f : \of{[}{0,1}{]} \to T_{k}$$
telle que~$f(\of{[}{0,1}{[})\subset T_{k}$ et~$f(1)=v$. Puisque la fonction~$f$ est semi-alg\'ebrique, quitte \`a restreindre son intervalle de d\'efinition puis effectuer un changement d'\'echelle pour se ramener \`a~$\of{[}{0,1}{]}$, nous pouvons supposer que les~$d$ premi\`eres fonctions coordonn\'ees de~$f$ sont d\'ecroissantes et que les~$d$ derni\`eres sont croissantes. Soit~$(x,y)$ un point de~$\R^2$ tel que~$(x,y,f(0))\in E$. Quel que soit~$u\in\of{[}{0,1}{[}$, nous avons alors encore~$(x,y,f(u))\in E$. 

Soient~$z_{1},z_{2}$ des \'el\'ements de~$\R^2$ tels que~$(z_{1},f(0))$ et~$(z_{2},f(0))$ appartiennent \`a~$E$. Quand les nombres~$s_{1},\ldots,s_{d}$ sont assez petits et les nombres~$t_{1},\ldots,t_{d}$ assez grands, les points~$z_{1}$ et~$z_{2}$ appartiennent \`a la m\^eme composante connexe de
$$\bigcap_{1\le j\le d} \left\{ (x,y)\in\R^2\, \big|\, s_{j}<|P_{j}(x+iy)|^2<t_{j}\right\}.$$
On en d\'eduit qu'il existe~$u\in\of{[}{0,1}{[}$ tels que les points~$(z_{1},f(u))$ et~$(z_{2},f(u))$ appartiennent \`a la m\^eme composante connexe de~$p^{-1}(f(u))$. Le morphisme~$p$ \'etant semi-alg\'ebriquement trivial au-dessus de~$T_{k}$, les points~$(z_{1},f(0))$ et~$(z_{2},f(0))$ doivent \'egalement appartenir \`a la m\^eme composante connexe de~$p^{-1}(f(0))$. On en d\'eduit que la partie~$F_{k}$ est connexe, ce qui conclut la preuve.
\end{proof}

\begin{prop}\label{sectionGausscentral}\index{Voisinages d'un point!de type 2!de la fibre centrale}
Notons $x$ le point de Gau{\ss} de la fibre centrale. Soit $U$ un voisinage de $x$ dans $X$. Alors il existe un voisinage ouvert $W$ de $x$ dans $U$ v\'erifiant les propri\'et\'es suivantes :
\begin{enumerate}[\it i)]
\item la projection $\pi(W)$ est un voisinage ouvert et connexe par arcs de $\pi(x)=a_{0}$ dans~$B$ ;
\item il existe une section topologique de $\pi$ au-dessus de $\pi(W)$ \`a valeurs dans $W$ ;
\item pour tout point $b$ de $\pi(W)$, la trace de la fibre $X_{b}$ sur $W$ est connexe par arcs ;
\item quels que soient $x\in W$ et $\eps\in\of{[}{0,1}{]}$, le point $x^\eps$ appartient \`a $W$.
\end{enumerate} 
\end{prop}
\begin{proof}
Par d\'efinition de la topologie de $X$, il existe un entier \mbox{$r\in\N^*$}, des polyn\^omes $f_{1},\ldots,f_{r}\in A[T]$ et un nombre r\'eel $\lambda >0$ tels que $U$ contienne une partie de la forme
$$V = \bigcap_{1\le i\le r} \left\{ y\in X\,\big|\, |f_{i}(x)|-\lambda < |f_{i}(y)| < |f_{i}(x)| + \lambda \right\}.$$
Nous pouvons supposer que, quel que soit $i\in\cn{1}{d}$, nous avons $f_{i}\ne 0$. Alors 
$$V = \bigcap_{1\le i\le r} \left\{ y\in X\,\big|\, 1-\lambda < |f_{i}(y)| < 1+ \lambda \right\}.$$
Nous allons montrer qu'il existe un voisinage $E$ de $a_{0}$ dans $B$ tel que le voisinage \mbox{$W=V\cap X_{E}$} de $x$ dans $X$ v\'erifie les propri\'et\'es requises. Nous allons proc\'eder en plusieurs \'etapes en prouvant tout d'abord le r\'esultat au-dessus de la partie ultram\'etrique de $B$, puis au-dessus de chacune des branches archim\'ediennes. Le r\'esultat global en d\'ecoulera pourvu que les sections que nous aurons alors construites se recollent sur la fibre centrale. De fa\c{c}on \`a en \^etre certain, nous imposerons \`a toutes les sections d'envoyer le point central~$a_{0}$ sur le point de Gau{\ss}~$\eta_{1}$ de la fibre centrale.\\

Notons 
$$B_{um} = \bigcup_{\m\in\Sigma_{f}} B_{\m}$$
la partie ultram\'etrique de $B$. On d\'efinit une section topologique $\sigma_{G}$ de la projection $\pi$ au-dessus de $B_{um}$ en associant \`a tout point $b$ de $B_{um}$ le point de Gau{\ss} de la fibre $X_{b}$.

Soit $i\in\cn{1}{r}$. Notons 
$$V_{i} =  \left\{ y\in X\,\big|\, 1-\lambda < |f_{i}(y)| < 1+\lambda \right\}.$$
Remarquons que, quels que soient $x\in V_{i}$ et $\eps\in\of{[}{0,1}{]}$, nous avons $x^\eps\in V_{i}$.

Il existe $d_{i}\in\N^*$ et $f_{i,0},\ldots,f_{i,d_{i}}\in A$, avec $f_{i,d_{i}} \ne 0$, tels que 
$$f_{i}(T) = \sum_{j=0}^{d_{i}} f_{i,j}\, T^j.$$
Posons 
$$C_{i} = \bigcap_{0\le j\le d_{i}} \left\{b\in B_{um}\, \big|\, |f_{i,j}(a_{0})|-\lambda < |f_{i,j}(b)| < |f_{i,j}(a_{0})|+\lambda \right\}.$$
C'est un voisinage du point central $a_{0}$ de $B_{um}$. La section topologique $\sigma_{G}$ de $\pi$ restreinte \`a $C_{i}$ prend ses valeurs dans $V_{i}$. 

Notons $D_{i}$ l'ensemble des points de $B_{um}$ o\`u la fonction $f_{i,d_{i}}$ est inversible. D\'efinissons alors une fonction continue $\rho_{i}$ de $D_{i}$ dans $\R_{+}$ en associant \`a tout point $b$ de $D_{i}$ le nombre r\'eel 
$$\rho_{i}(b) = \max_{0\le j\le d_{i}-1} \left( \left|\frac{f_{i,j}(b)}{f_{i,d_{i}}(b)}\right|^{\frac{1}{d_{i}-j}} \right).$$
Notons $D'_{i}$ le voisinage ouvert de $a_{0}$ dans $D_{i}$ d\'efini par
$$D'_{i} = \left\{b\in D_{i}\, \big|\, |\rho_{i}(b)|< 1+\lambda\right\}.$$
Finalement, choisissons $E_{i}$ un voisinage ouvert et connexe par arcs de $a_{0}$ dans \mbox{$C_{i}\cap D'_{i}$}. Quels que soient $x\in E_{i}$ et $\eps\in\of{[}{0,1}{]}$, nous avons alors $x^\eps\in E_{i}$.

Posons 
$$E=\bigcap_{1\le i\le r} E_{i}$$
et
$$W = V \cap X_{E}.$$
Les premi\`ere, troisi\`eme et quatri\`eme propri\'et\'es de l'\'enonc\'e sont alors clairement v\'erifi\'ees. Soit $b\in E=\pi(W)$. Quel que soit $i\in\cn{1}{r}$, d'apr\`es le lemme~\ref{cparho} et puisque~$b\in D'_{i}$, la partie $V_{i}\cap X_{b}$ est connexe par arcs. Puisque la fibre~$X_{b}$ est un arbre, l'intersection $V\cap X_{b}$ de toutes ces parties est donc connexe par arcs.\\ 

Passons maintenant aux branches archim\'ediennes de $B$. Soit $\sigma\in\Sigma_{\infty}$. Nous avons
$$\lim_{\eps \xrightarrow[>]{} 0} (1-\lambda)^{1/\eps}=0 \textrm{ et } \lim_{\eps \xrightarrow[>]{} 0} (1+\lambda)^{1/\eps}=+\infty.$$
Par cons\'equent, d'apr\`es le lemme~\ref{cpaHardt}, il existe~$\eta>0$ tel que, quel que soit~\mbox{$\eps\in\of{]}{0,\eta}{[}$}, la partie
$$\bigcap_{1\le i\le r} \left\{ y\in \E{1}{\hat{K}_{\sigma}}\, \big|\, (1-\lambda)^{1/\eps} < |f_{i}(y)|_{\sigma} < (1+\lambda)^{1/\eps} \right\}$$ 
est connexe par arcs. En d'autres termes, quel que soit~$\eps\in\of{]}{0,\eta}{[}$, la trace de la fibre~$X_{a_{\sigma}^\eps}$ sur~$V$ est connexe par arcs. Le lemme~\ref{cparho} nous montre que la trace de la fibre centrale $X_{0}=X_{a_{\sigma}^0}$ sur~$V$ est \'egalement connexe par arcs.

Soit~$\alpha$ un nombre r\'eel transcendant. Consid\'erons l'application~$\sigma_{G}$ qui au point~$a_{\sigma}^\eps$ de~$B'_{\sigma}$, avec~$\eps\in\of{]}{0,1}{]}$, associe le point de~$X$ associ\'e \`a la semi-norme multiplicative sur~$A[T]$, born\'ee sur~$A$, d\'efinie par
$$\begin{array}{ccc}
A[T] & \to & \R_{+}\\
P(T) & \mapsto & |P(\alpha)|_{\infty}^\eps
\end{array}
$$
et au point~$a_{0}$ associe le point de Gau{\ss}~$\eta_{1}$ de la fibre centrale~$X_{0}$. Cette application~$\sigma_{G}$ d\'efinit une section topologique continue de la projection~$\pi$ au-dessus de~$B_{\sigma}$.

Soit~$i\in\cn{1}{d}$. Puisque~$\alpha$ est transcendant, nous avons~$f_{i}(\alpha)\ne 0$ dans~$\hat{K}_{\sigma}$. Par cons\'equent, il existe~$\eta_{i}>0$ tel que, quel que soit~\mbox{$\eps\in\of{]}{0,\eta_{i}}{[}$}, on ait
$$ (1-\lambda)^{1/\eps} < |f_{i}(\alpha)|_{\sigma} < (1+\lambda)^{1/\eps}.$$ 
Posons $\zeta = \min_{1\le i\le d} (\eta_{i})$. Au-dessus du voisinage~$\of{[}{a_{0},a_{\sigma}^\zeta}{[}$ de~$a_{0}$ dans~$B_{\sigma}$, la restriction de la section~$\sigma_{G}$ est \`a valeurs dans~$V$. On en d\'eduit le r\'esultat annonc\'e.
\end{proof}

Nous obtenons imm\'ediatement les deux corollaires suivants.

\begin{cor}\label{cpaGausscentral}\index{Connexite par arcs au voisinage d'un point@Connexit\'e par arcs au voisinage d'un point!de type 2 de la fibre centrale}
Le point de Gau{\ss} de la fibre centrale poss\`ede un syst\`eme fondamental de voisinages connexes par arcs.
\end{cor}

\begin{cor}\label{ouvertGausscentral}\index{Ouverture au voisinage d'un point!de type 2 de la fibre centrale}
Le morphisme $\pi$ est ouvert au voisinage du point de Gau{\ss} de la fibre centrale.
\end{cor}

Int\'eressons-nous, \`a pr\'esent, \`a l'anneau local.

\begin{prop}\label{corps2central}\index{Anneau local en un point!de type 2!de la fibre centrale}
Notons~$x$ le point de Gau{\ss} de la fibre centrale. L'anneau local $\Os_{X,x}$ est un corps, canoniquement isomorphe au corps~$K(T)$. Il est complet pour la valeur absolue associ\'ee au point~$x$, qui n'est autre que la valeur absolue triviale.
\end{prop}
\begin{proof}
Commen\c{c}ons par prouver que l'anneau local~$\Os_{X,x}$ est un corps. Il suffit de montrer que son id\'eal maximal est r\'eduit \`a~$(0)$. Soit~$f$ une fonction d\'efinie sur un voisinage~$U$ de~$x$ dans~$X$ et s'annulant en~$x$. Nous voulons montrer que~$f$ s'annule encore au voisinage de~$x$ dans~$X$. 

D'apr\`es la proposition \ref{sectionGausscentral}, il existe un voisinage ouvert $W$ de $x$ dans $U$ v\'erifiant les propri\'et\'es suivantes :
\begin{enumerate}[\it i)]
\item la projection $\pi(W)$ est un voisinage ouvert connexe par arcs de $\pi(x)=a_{0}$ dans~$B$ ;
\item il existe une section topologique de $\pi$ au-dessus de $\pi(W)$ \`a valeurs dans $W$ ;
\item pour tout point $b$ de $\pi(W)$, la trace de la fibre $X_{b}$ sur $W$ est connexe par arcs ;
\item quel que soient $x\in W$ et $\eps\in\of{[}{0,1}{]}$, le point $x^\eps$ appartient \`a $W$.
\end{enumerate} 

Soit $\sigma\in\Sigma$. Notons $W'_{\sigma} =W\cap X'_{\sigma}$. C'est la trace de $W$ sur la branche $\sigma$-adique ouverte. Soit $b\in \pi(W'_{\sigma})$. Soit $u$ un point rigide de $W\cap X_{b}$ tel que l'extension $\hat{K}_{\sigma}=\Hs(b) \to \Hs(u)$ soit transcendante. Consid\'erons l'application suivante, induite par le flot : 
$$\begin{array}{ccc}
[0,1] &\to & X\\
\theta & \mapsto & u^{1-\theta}
\end{array}.$$
Son image d\'efinit un chemin continu trac\'e sur $W$ et joignant le point $u$ au point $u^0$ de la fibre centrale. Puisque l'extension $\hat{K}_{\sigma}=\Hs(b) \to \Hs(u)$ est transcendante, le point~$u^0$ n'est autre que le point~$x$, le point de Gau{\ss} de la fibre centrale. D'apr\`es le lemme~\ref{voisflotcentral} et la proposition~\ref{flot}, quel que soit~$\theta\in\of{[}{0,1}{]}$, nous avons
$$|f(\eta_{1})|=|f(u^0)|=|f(u^\theta)|^0.$$ 
On en d\'eduit que $|f(u)|=0$. La fonction $f$ s'annule donc sur tous les points transcendants de $W\cap X_{b}$. Puisque $W\cap X_{b}$ est normal et connexe, la fonction $f$ y est identiquement nulle. 
Nous avons donc montr\'e que la fonction~$f$ est identiquement nulle sur $W'_{\sigma}$. La continuit\'e de~$f$ nous permet de montrer qu'elle est encore nulle sur~$W\cap X_{\sigma}$. On en d\'eduit finalement que la fonction $f$ est nulle sur $W$.\\

D\'emontrons, \`a pr\'esent, la derni\`ere partie de la proposition. Puisque l'anneau local~$\Os_{X,x}$ est un corps, le morphisme $\Os_{X,x} \to \Hs(x)$ est injectif. L'\'egalit\'e \mbox{$\Hs(x)=K(T)$} nous montre qu'il est \'egalement surjectif.
\end{proof}

\begin{cor}\label{prolanGausscentral}\index{Prolongement analytique!au voisinage d'un point de type 2!de la fibre centrale}
Le principe du prolongement analytique vaut au voisinage du point de {Gau\ss} de la fibre centrale de l'espace~$X$.
\end{cor}

\index{Point!de type 2!de la fibre centrale|)}

\section{R\'esum\'e}\label{sectionresume}

Dans cette partie, nous regroupons les r\'esultats que nous avons obtenu concernant la droite affine analytique sur un anneau d'entiers de corps de nombres. Rappelons que~$A$ d\'esigne un anneau d'entiers de corps de nombres, $B=\Ms(A)$ son spectre analytique, $X=\E{1}{A}$ la droite affine analytique au-dessus de~$A$ et~$\pi : X\to B$ le morphisme de projection.

\begin{thm}\label{resume}\index{Droite affine analytique!sur A@sur $A$!proprietes@propri\'et\'es}\index{Droite affine analytique!sur A@sur $A$!metrisabilite@m\'etrisabilit\'e}\index{Droite affine analytique!sur A@sur $A$!dimension topologique}\index{Dimension topologique! de A1@de $\E{1}{A}$}\index{Connexite par arcs au voisinage d'un point@Connexit\'e par arcs au voisinage d'un point!de A1@de $\AA$}\index{Ouverture au voisinage d'un point!de A1@de $\AA$}\index{Anneau local en un point!de A1@de $\AA$}
\begin{enumerate}[\it i)]
\item L'espace $X$ est localement compact, m\'etrisable et de dimension topologique~$3$.
\item L'espace $X$ est localement connexe par arcs.
\item Le morphisme de projection $\pi : X\to B$ est ouvert.
\item En tout point~$x$ de~$X$, l'anneau local~$\Os_{X,x}$ est hens\'elien, noeth\'erien, r\'egulier, de dimension inf\'erieure \`a $2$ et le corps r\'esiduel~$\kappa(x)$ est hens\'elien.
\end{enumerate}
\end{thm}
\begin{proof}
Le point~{\it i)} provient des th\'eor\`emes \ref{topoAn} et \ref{dim}. Le point~{\it ii)} est obtenu en regroupant les r\'esultats des th\'eor\`emes \ref{toporig} et \ref{topointerne} et des corollaires \ref{cpa3ext}, \ref{cpa3central}, \ref{cpaGaussext} et~\ref{cpaGausscentral}. Le point~{\it iii)} est obtenu en regroupant les r\'esultats des th\'eor\`emes \ref{toporig} et \ref{topointerne} et des corollaires \ref{ouvert3ext}, \ref{ouvert3central}, \ref{ouvertGaussext} et~\ref{ouvertGausscentral}. Le point~{\it iv)} est obtenu en regroupant les r\'esultats des th\'eor\`emes \ref{recrigavd}, \ref{recrigext} et \ref{topointerne}, des corollaires~\ref{avd3ext} et \ref{corps3central} et des propositions~\ref{corps2central} et~\ref{avd2ext}. 

\end{proof}

\bigskip

\begin{thm}\label{prolan}\index{Prolongement analytique!au voisinage d'un point de A1@au voisinage d'un point de $\AA$}
Le principe du prolongement analytique vaut au voisinage de tout point de l'espace~$X$. Par cons\'equent, il vaut sur tout ouvert connexe de l'espace~$X$.
\end{thm}
\begin{proof}
Il suffit de regrouper les r\'esultats de la proposition \ref{prolinterne}, des corollaires \ref{prolanrig}, \ref{prolan3ext}, \ref{prolan3central}, \ref{prolanGaussext} et \ref{prolanGausscentral}. La deuxi\`eme partie provient du lemme \ref{lemprolan}.
\end{proof}

De ce th\'eor\`eme d\'ecoulent plusieurs r\'esultats concernant les anneaux globaux de fonctions holomorphes et m\'eromorphes sur les parties connexes de la droite analytique~$X$.

\begin{defi}\label{defimero}\index{Fonctions meromorphes sur A1@Fonctions m\'eromorphes sur $\AA$|(}\index{Faisceau!des fonctions meromorphes@des fonctions m\'eromorphes|see{Fonctions m\'eromorphes}}
Nous appellerons {\bf faisceau des fonctions m\'eromorphes} et noterons~$\Ms$ le faisceau associ\'e au pr\'efaisceau qui envoie tout ouvert~$U$ de~$X$ sur l'anneau total des fractions de l'anneau~$\Os(U)$. 
\end{defi} 

\begin{cor}\label{prolmero}\index{Prolongement analytique!pour les fonctions meromorphes sur A1@pour les fonctions m\'eromorphes sur $\AA$}
Soient~$U$ une partie connexe de~$X$ et~$x$ un point de~$U$. Le morphisme de restriction
$$\Ms(U) \to \Ms_{x}$$
est injectif.
\end{cor}
\begin{proof}
Soit~$s$ un \'el\'ement de~$\Ms(U)$ dont l'image dans~$\Ms_{x}$ est nulle. Notons
$$V = \{y\in U\, |\, s_{y}=0 \textrm{ dans } \Ms_{y}\}.$$
C'est une partie non vide et ouverte de~$U$. 

Soit~$y$ un point de~$U\setminus V$. Il existe un voisinage~$W$ du point~$y$ dans~$U$ et deux \'el\'ements~$f$ et~$g$ de~$\Os(W)$, $g$ ne divisant pas~$0$, tels que 
$$s =\frac{f}{g} \textrm{ dans } \Ms(W).$$
Par hypoth\`ese, le germe~$s_{y}$ n'est pas nul dans l'anneau local~$\Os_{X,y}$. D'apr\`es le th\'eor\`eme \ref{prolan}, il existe un voisinage~$W'$ du point~$y$ dans~$W$ tel que l'image de la fonction~$f$ ne soit nulle dans aucun des anneaux locaux~$\Os_{X,z}$, pour~$z$ appartenant \`a~$W'$. Soit~$z$ un \'el\'ement de~$W'$. D'apr\`es le th\'eor\`eme \ref{resume}, {\it iv}), l'anneau local~$\Os_{X,z}$ est int\`egre. L'\'el\'ement~$s_{z}$ de $\Ms_{z}=\Frac(\Os_{X,z})$ n'est donc pas nul. On en d\'eduit que le voisinage~$W'$ du point~$y$ est contenu dans~$U\setminus V$. Par cons\'equent, la partie~$V$ est ferm\'ee dans~$U$. La connexit\'e de~$U$ assure qu'elle est \'egale \`a la partie~$U$ tout enti\`ere.
\end{proof}

\begin{cor}\label{anneauintegre}\index{Faisceau structural!sections sur une partie connexe de A1@sections sur une partie connexe de $\AA$}
Soit~$U$ une partie connexe de l'espace~$X$. L'anneau~$\Os(U)$ est int\`egre et l'anneau~$\Ms(U)$ est un corps.
\end{cor}
\begin{proof}
Soit~$x$ un point de~$U$. D'apr\`es le th\'eor\`eme \ref{prolan}, le morphisme naturel
$$\Os(U) \to \Os_{X,x}$$
est injectif. D'apr\`es le th\'eor\`eme \ref{resume}, {\it iv}), l'anneau local~$\Os_{X,x}$ est r\'egulier et donc int\`egre. On en d\'eduit que l'anneau~$\Os(U)$ est int\`egre.

Soit~$s$ un \'el\'ement non nul de~$\Ms(U)$. D'apr\`es le corollaire \ref{prolmero}, en tout point~$x$ de~$U$ l'\'el\'ement~$s_{x}$ de $\Ms_{x}=\Frac(\Os_{X,x})$ est non nul et donc inversible. On en d\'eduit que l'\'el\'ement~$s$ lui-m\^eme est inversible et donc que l'anneau~$\Ms(U)$ est un corps.
\end{proof}

\begin{cor}\label{fnmero}
Soit~$U$ une partie connexe de l'espace~$X$ contenant le point de Gau{\ss} de la fibre centrale. Alors l'anneau des fonctions m\'e\-ro\-mor\-phes sur~$U$ est l'anneau des fractions rationnelles~$K(T)$.
\end{cor}
\begin{proof}
Notons~$x$ le point de Gau{\ss} de la fibre centrale. D'apr\`es la proposition \ref{corps2central}, l'anneau local~$\Os_{X,x}$ est canoniquement isomorphe au corps~$K(T)$. D'apr\`es le corollaire \ref{prolmero}, le morphisme de restriction
$$\Ms(U) \to \Os_{X,x}=K(T)$$
est injectif. Il est \'evident qu'il est \'egalement surjectif, ce qui conclut la preuve.
\end{proof}

\index{Fonctions meromorphes sur A1@Fonctions m\'eromorphes sur $\AA$|)}

\bigskip

\begin{thm}\label{uniforte}\index{Uniformisante forte!en un point de A1@en un point de $\AA$}
Soit~$x$ un point de l'espace~$X$ en lequel l'anneau local~$\Os_{X,x}$ est un anneau de valuation discr\`ete. Soient~$\pi$ une uniformisante de~$\Os_{X,x}$ et~$U$ un voisinage du point~$x$ dans~$X$ sur lequel elle est d\'efinie. Alors il existe un syst\`eme fondamental~$\Vs$ de voisinages compacts et connexes du point~$x$ dans~$U$ tel que, pour tout \'el\'ement~$V$ de~$\Vs$, la fonction~$\pi$ est une uniformisante forte de l'anneau~$\Os_{X,x}$ sur~$V$.
\end{thm}
\begin{proof}
Nous pouvons d\'ecrire exactement l'ensemble des points en lequel l'anneau local est un anneau de valuation discr\`ete : il s'agit des points rigides des fibres internes et centrale et des points de type~$2$ et~$3$ des fibres extr\^emes. Le r\'esultat attendu se d\'eduit alors des corollaires \ref{uniforterig}, \ref{uniforte3ext} et \ref{uniforte2ext}.
\end{proof}

\begin{thm}\label{Shilovum}\index{Bord analytique!au voisinage du point de A1@au voisinage d'un point de $\AA$}
Tout point de $X_{\textrm{um}}\setminus X_{0}$ poss\`ede un syst\`eme fondamental de voisinages compacts, connexes et spectralement convexes qui poss\`edent un bord analytique fini et alg\'ebriquement trivial.
\end{thm}
\begin{proof}
On regroupe les r\'esultats des propositions \ref{recShilovinterneum} et \ref{recShilovrigideum} et des corollaires \ref{Shilov3ext} et \ref{Shilov2ext}.
\end{proof}

\section{Coh\'erence}\label{parcoherence}

Dans cette partie, nous montrons que le faisceau structural $\Os_{X}$ de la droite analytique $X$ est coh\'erent. Rappelons, auparavant, quelques d\'efinitions et notations. Fixons un espace localement annel\'e~$(Y,\Os_{Y})$.

\begin{defi}\label{defitf}\index{Faisceau!de type fini}
Un faisceau de~$\Os_{Y}$-modules~$\Fs$ est dit {\bf de type fini} si, pour tout point~$y$ de~$Y$, il existe un voisinage~$V$ de~$y$ dans~$Y$, un entier~$p$ et des \'el\'ements $F_{1},\ldots,F_{p}$ de~$\Fs(V)$ tels que, pour tout point~$z$ de~$V$, le $\Os_{Y,z}$-module~$\Fs_{z}$ soit engendr\'e par les germes $(F_{1})_{z},\ldots,(F_{p})_{z}$.
\end{defi}

\begin{defi}\index{Faisceau!des relations}
Soient~$V$ une partie ouverte de~$Y$, $\Fs$ un faisceau de~$\Os_{Y}$-modules, $q\in\N$ et~$F_{1},\ldots,F_{q}\in\Fs(V)$. On appelle {\bf faisceau des relations} entre~$F_{1},\ldots,F_{q}$, et on note~$\Rs(F_{1},\ldots,F_{q})$, le noyau du morphisme de faisceau suivant
$$\begin{array}{ccc}
\Os_{V}^q & \to & \Fs_{V}\\
(a_{1},\ldots,a_{q}) & \mapsto & \disp \sum_{i=1}^q a_{i}\, F_{i}
\end{array}.$$
\end{defi}

\begin{defi}\index{Faisceau!coherent@coh\'erent}
Un faisceau de $\Os_{Y}$-modules~$\Fs$ est dit {\bf coh\'erent} s'il v\'erifie les deux propri\'et\'es suivantes :
\begin{enumerate}[\it i)]
\item le faisceau~$\Fs$ est localement de type fini ;
\item quels que soient l'ouvert~$V$ de~$Y$, l'entier~$q$ et les \'el\'ements~\mbox{$F_{1},\ldots,F_{q}$} de~$\Fs(V)$, le faisceau~$\Rs(F_{1},\ldots,F_{q})$ des relations entre~$F_{1},\ldots,F_{q}$ est localement de type fini.
\end{enumerate} 
\end{defi}

Venons-en, \`a pr\'esent, \`a la preuve de la coh\'erence du faisceau~$\Os_{X}$. Il est \'evidemment localement de type fini. Il nous reste \`a \'etudier les faisceaux de relations. Commen\c{c}ons par un lemme.

\begin{lem}
Soit $x$ un point de $X$. Soient $U$ un voisinage ouvert de $x$ dans $X$, $p\in\N^*$ et $f_{1},\ldots,f_{p} \in \Os(U)$. Notons $(e_{1},\ldots,e_{p})$ la base canonique de $\Os_{X}^p$.  Supposons qu'il existe $l\in\cn{1}{p}$ tel que $f_{l}\ne 0$ dans $\Os_{X,x}$. Si l'anneau local $\Os_{X,x}$ est un anneau de valuation discr\`ete ou un corps, alors il existe un voisinage ouvert $V$ de $x$ dans $U$ tel que, quel que soit $y\in V$, la famille 
$$(f_{j} e_{i} - f_{i} e_{j})_{1\le i < j\le p}$$
de~$\Os_{X,y}^p$ engendre le germe~$\Rs(f_{1},\ldots,f_{p})_{y}$. 
\end{lem}
\begin{proof}
Supposons que l'anneau local $\Os_{X,x}$ est un anneau de valuation discr\`ete. Choisissons-en une uniformisante $\tau$. Quitte \`a restreindre $U$, nous pouvons supposer que $\tau$ est d\'efinie sur $U$. Notons $m$ le minimum des valuations des \'el\'ements $f_{1},\ldots,f_{p}$ de $\Os_{X,x}$. Puisque l'un de ces \'el\'ements n'est pas nul, nous avons $m \in\N$. Remarquons que, quel que soit $i\in\cn{1}{p}$, nous avons $\tau^{-m} f_{i} \in\Os_{X,x}$. Par choix de $m$, il existe $j\in\cn{1}{p}$ tel que la fonction $\tau^{-m}f_{j}$ soit inversible dans $\Os_{X,x}$. Il existe donc un voisinage ouvert~$V$ de~$x$ dans~$U$ sur lequel les fonctions $\tau^{-m} f_{1},\ldots,\tau^{-m}f_{p}$ sont d\'efinies et la fonction $\tau^{-m}f_{j}$ inversible. D'apr\`es le th\'eor\`eme~\ref{resume}, nous pouvons supposer que la partie~$V$ est connexe.

Nous disposons de l'inclusion suivante entre faisceaux de~$\Os_{V}$-modules :
$$\Rs(\tau^{-m} f_{1},\ldots,\tau^{-m}f_{p}) \subset \Rs(f_{1},\ldots,f_{p}).$$
Montrons que c'est une \'egalit\'e. Il suffit pour cela de montrer que l'inclusion induit une \'egalit\'e entre les germes. Soit~$y$ un point de~$V$. Remarquons tout d'abord que l'image de~$\tau$ dans l'anneau local~$\Os_{X,y}$ n'est pas nulle. Dans le cas contraire, le principe du prolongement analytique (\emph{cf.} th\'eor\`eme~\ref{prolan}) imposerait en effet \`a la fonction~$\tau$ d'\^etre nulle sur l'ouvert connexe~$V$ tout entier, mais nous savons qu'elle n'est pas nulle au voisinage du point~$x$. Soit~$(a_{1},\ldots,a_{p})\in \Rs(f_{1},\ldots,f_{p})_{y}$. Nous avons alors
$$\sum_{i=1}^p a_{i}f_{i} = \tau^m\, \left(\sum_{i=1}^p a_{i}\tau^{-m}f_{i}\right)=0 \textrm{ dans } \Os_{X,y}.$$
D'apr\`es le th\'eor\`eme~\ref{resume}, l'anneau local~$\Os_{X,y}$ est int\`egre. On en d\'eduit que 
$$(a_{1},\ldots,a_{p})\in \Rs(\tau^{-m}f_{1},\ldots,\tau^{-m}f_{p})_{y}.$$

Par cons\'equent, nous pouvons supposer qu'il existe~$j\in\cn{1}{p}$ tel que la fonction~$f_{j}$ est inversible sur~$V$. Soient $y\in V$ et $(a_{1},\ldots,a_{p})\in \Rs(f_{1},\ldots,f_{p})_{y}$. Nous avons alors 
$$\sum_{i=1}^p a_{i} f_{i} = 0 \textrm{ dans } \Os_{X,y}.$$
Pour conclure, il nous suffit de remarquer que, dans $\Os_{X,y}^p$, nous avons
$$\begin{array}{rcl}
\disp \sum_{i\ne j} a_{i} f_{j}^{-1} (f_{j}e_{i} - f_{i} e_{j}) &=&\disp \sum_{i\ne j} a_{i}e_{i} -  \left(\sum_{i\ne j} a_{i}f_{i}\right)  f_{j}^{-1} e_{j}\\
&=&\disp  \sum_{i\ne j} a_{i}e_{i} - (-a_{j}f_{j}) f_{j}^{-1} e_{j}\\
&=&\disp \sum_{i=1}^p a_{i} e_{i}.
\end{array}$$

On d\'emontre le r\'esultat par la m\^eme m\'ethode lorsque l'anneau local $\Os_{X,x}$ est un corps. Dans ce cas, les r\'eductions pr\'eliminaires sont inutiles et l'on passer directement \`a la derni\`ere \'etape.
\end{proof}

D\'emontrons, finalement, le r\'esultat attendu.

\begin{thm}\label{coherence}\index{Faisceau structural!coherence sur A1@coh\'erence sur $\AA$}\index{Droite affine analytique!sur A@sur $A$!coherence du faisceau structural@coh\'erence du faisceau structural}
Le faisceau structural $\Os_{X}$ est coh\'erent.
\end{thm}
\begin{proof}
Soit $x$ un point de $X$. Soient $U$ un voisinage ouvert de $x$ dans $X$, $p\in\N^*$ et $f_{1},\ldots,f_{p} \in \Os(U)$. Il nous suffit de montrer que le faisceau des relations $\Rs(f_{1},\ldots,f_{p})$ est de type fini au voisinage du point $x$.

Si les fonctions $f_{1},\ldots,f_{p}$ sont nulles dans $\Os_{X,x}$, alors, par le principe du prolongement analytique, elles sont nulles au voisinage du point $x$ et le r\'esultat est imm\'ediat. Par cons\'equent, nous pouvons supposer qu'il existe $l\in\cn{1}{p}$ tel que $f_{l}\ne 0$ dans $\Os_{X,x}$. 

Si l'anneau local $\Os_{X,x}$ est un anneau de valuation discr\`ete ou un corps, alors le lemme pr\'ec\'edent nous permet de conclure.

Il nous reste, \`a pr\'esent, \`a traiter le cas o\`u l'anneau local $\Os_{X,x}$ n'est ni un anneau de valuation discr\`ete, ni un corps. Cela impose au point $x$ d'\^etre un point rigide d'une fibre extr\^eme. 

D'apr\`es le th\'eor\`eme~\ref{noetherienrigext}, l'anneau local~$\Os_{X,x}$ est noeth\'erien. Par cons\'equent, le $\Os_{X,x}$-module $\Rs(f_{1},\ldots,f_{p})_{x}$ est de type fini. Il existe donc un entier $q\in\N^*$, un voisinage ouvert $V$ de $x$ et des fonctions $g_{1},\ldots,g_{q} \in\Os(V)^p$ tels que le $\Os_{X,x}$-module $\Rs(f_{1},\ldots,f_{p})_{x}$ soit engendr\'e par $((g_{1})_{x},\ldots,(g_{q})_{x})$. 
 
Puisque les fibres extr\^emes sont des droites analytiques sur des corps trivialement valu\'es, l'ensemble de leurs points rigides est discret. Par cons\'equent, l'ensemble des points de~$X$ en lequel l'anneau local est de dimension~$2$ forme une partie discr\`ete de l'espace~$X$.  Quitte \`a restreindre~$V$, nous pouvons donc supposer que~$x$ est le seul point de~$V$ en lequel l'anneau local n'est ni un anneau de valuation discr\`ete, ni un corps. Alors, d'apr\`es le lemme pr\'ec\'edent, quel que soit $y\in V\setminus\{x\}$, le $\Os_{X,y}$-module $\Rs(f_{1},\ldots,f_{p})_{y}$ est engendr\'e par la famille $(f_{j} e_{i} - f_{i} e_{j})_{1\le i < j\le p}$ de $\Os_{X,y}^p$. Par cons\'equent, le faisceau $\Rs(f_{1},\ldots,f_{p})$ est de type fini au voisinage du point $x$.
\end{proof}

\index{Droite affine analytique!sur A@sur $A$|)}

%% file: fini.tex
\chapter{Morphismes finis}\label{chapitrefini}

Nous \'etudions, dans ce chapitre, quelques cas particuliers de morphismes finis entre espaces analytiques au sens de V.~Berkovich. Exception faite du dernier num\'ero, nous quittons, ici, les espaces analytiques au-dessus d'un anneau d'entiers de corps de nombres pour revenir au cadre g\'en\'eral, au-dessus d'un anneau de Banach muni d'une norme uniforme.

Le num\'ero \ref{mtf} est consacr\'e aux morphismes finis au sens topologique. Nous nous contentons d'y rappeler les r\'esultats classiques dont nous aurons besoin par la suite.

Au num\'ero \ref{Wglobal}, nous d\'emontrons un th\'eor\`eme de division de Weierstra{\ss} que nous qualifions de global. Il permet en effet de diviser une fonction d\'efinie sur un disque de dimension~$1$ par un polyn\^ome donn\'e, pourvu que le rayon du disque soit assez grand. Si nous disposons d'un anneau de Banach uniforme $(\As,\|.\|)$, ce th\'eor\`eme nous permettra de munir de normes uniformes certaines extensions finies de l'anneau~$\As$.

Au num\'ero \ref{unexemple}, nous nous int\'eresserons \`a un cas particulier de morphisme fini. Un anneau de Banach $(\As,\|.\|)$ muni d'une norme uniforme \'etant fix\'e, nous consid\'ererons un morphisme d'une hypersurface de la droite~$\E{1}{\As}$ au-dessus de~$\As$ vers le spectre analytique~$\Ms(\As)$ de~$\As$. Nous accorderons une attention particuli\`ere \`a l'image directe du faisceau structural de l'hypersurface. Nous d\'ecrirons ses fibres et donnerons des conditions n\'ecessaires pour qu'il soit coh\'erent.

Le num\'ero \ref{Wrigide} est de nouveau consacr\'e \`a la d\'emonstration d'un th\'eor\`eme de division de Weierstra{\ss}. Il s'agit, cette fois-ci, d'un th\'eor\`eme de nature locale, mais qui permet d'effectuer une division au voisinage des points rigides des fibres, et non plus seulement rationnels. \`A l'aide de ce r\'esultat, nous \'etudier, au num\'ero \ref{endodroite}, les endomorphismes de la droite au-dessus d'un anneau de Banach muni d'une norme uniforme donn\'es par un polyn\^ome. L\`a encore, nous nous int\'eresserons particuli\`erement \`a l'image directe, par ce morphisme, du faisceau structural de la droite.

\`A ce stade du chapitre, nous aurons introduit plusieurs conditions assurant que les morphismes \'etudi\'es jouissent de bonnes propri\'et\'es. Dans le num\'ero \ref{adcdn}, nous montrons qu'elles sont satisfaites lorsque l'anneau de Banach consid\'er\'e est un anneau d'entiers de corps de nombres. 

Signalons, pour finir, que nous sommes convaincu que les techniques introduites ici permettent de ramener l'\'etude des courbes analytiques au-dessus d'un anneau d'entiers de corps de nombres \`a celle de la droite, au moins lorsque les courbes en question proviennent de courbes alg\'ebriques.

\section{Morphismes topologiques finis}\label{mtf}

\index{Morphisme fini!au sens topologique|(}
 
Avant d'\'etudier les morphismes du point de vue alg\'ebrique, nous allons les consid\'erer du point de vue topologique. Nous obtiendrons d\'ej\`a ainsi plusieurs r\'esultats dignes d'int\'er\^et. Nous les \'enon\c{c}ons sans d\'emonstration et renvoyons le lecteur int\'eress\'e \`a~\cite{GR}, I, \S 1.

Dans toute cette section, nous fixons deux espaces topologiques s\'epar\'es~$X$ et~$Y$ et une application $\varphi : X \to Y$.

\begin{defi}\label{topofini}
Nous dirons que l'application~$\varphi : X \to Y$ est un {\bf morphisme topologique fini} si elle est continue, ferm\'ee et \`a fibres finies.
\end{defi}

La propri\'et\'e suivante des applications ferm\'ees est imm\'ediate. Elle nous sera utile \`a de nombreuses reprises.

\begin{lem}\label{voisfini}
Supposons que l'application~$\varphi$ est ferm\'ee. Alors, pour toute partie~$V$ de~$Y$, l'ensemble
$$\{\varphi^{-1}(W),\, W \textrm{ voisinage de } V \textrm{ dans } Y\}$$
est un syst\`eme fondamental de voisinages de~$\varphi^{-1}(V)$ dans~$X$.
\end{lem}


\begin{cor}\label{corvoisfini}
Soit~$V$ une partie de~$Y$. Notons
$$\varphi_{V} : \varphi^{-1}(V) \to V$$
le morphisme d\'eduit de~$\varphi$ par restriction et corestriction. Soit~$\Fs$ un faisceau sur~$X$. Si l'application~$\varphi$ est ferm\'ee, alors le morphisme naturel
$$(\varphi_{*}\Fs)_{V} \to (\varphi_{V})_{*}\Fs_{\varphi^{-1}(V)}$$
est un isomorphisme.
\end{cor}

Venons-en, maintenant, aux propri\'et\'es des morphismes topologiques finis.

\begin{thm}\label{finifibres}
Supposons que l'application~$\varphi$ est un morphisme topologique fini. Soient~$y$ un point de~$Y$ et~$x_{1},\ldots,x_{r}$, avec~$r\in\N^*$, ses ant\'ec\'edents par le morphisme~$\varphi$. Soit~$\Ss$ un faisceau en groupes ab\'eliens sur~$X$. Alors le morphisme naturel
$$(\varphi_{*}\Ss)_{y} \to \prod_{i=1}^r \Ss_{x_{i}}$$
est un isomorphisme.

En outre, si~$\Ss$ est un faisceau de~$\Os_{X}$-modules, alors le morphisme pr\'ec\'edent est un isomorphisme de~$\Os_{Y,y}$-modules.
\end{thm}

\begin{thm}\label{finisuite}
Supposons que l'application~$\varphi$ est un morphisme topologique fini. Soit~$\Ss' \to \Ss \to \Ss''$ une suite exacte de faisceaux en groupes ab\'eliens sur~$X$. Alors la suite des images directes
$$\varphi_{*}\Ss' \to \varphi_{*}\Ss \to \varphi_{*}\Ss''$$
est encore exacte.
\end{thm}

\begin{thm}\label{finicohomologie}
Supposons que l'application~$\varphi$ est un morphisme topologique fini. Soit~$\Ss$ un faisceau en groupes ab\'eliens sur~$X$. Alors, quel que soit~$q\in\N$, il existe un isomorphisme de groupes naturel 
$$H^q(X,\Ss) \simeq H^q(Y,\varphi_{*}\Ss).$$
\end{thm}

\index{Morphisme fini!au sens topologique|)}

\section{Th\'eor\`eme de division de {Weierstra\ss} global}\label{Wglobal}

Soit~$(\As,\|.\|)$ un anneau de Banach uniforme. Nous noterons~$B=\Ms(\As)$. Soient~$b$ un point de~$B$, $U$ un voisinage compact de~$b$ dans~$B$ et~$R$ un nombre r\'eel strictement positif. Le th\'eor\`eme de {Weierstra\ss} classique permet, sous certaines conditions, de diviser une s\'erie \`a coefficients dans~$\Bs(U)$ de rayon de convergence sup\'erieur \`a~$R$ par une autre. Pour pouvoir effectuer cette division, il est cependant n\'ecessaire, en g\'en\'eral, d'autoriser le voisinage~$U$ de~$b$ et de le rayon de convergence~$R$ \`a diminuer. Dans le th\'eor\`eme qui suit, nous montrons que si le diviseur est un polyn\^ome unitaire et que le nombre r\'eel~$R$ est assez grand, ces restrictions sont inutiles.

\begin{thm}[Th\'eor\`eme de division de Weierstra{\ss} global]\label{vraimentglobal}\index{Theoreme de Weierstrass@Th\'eor\`eme de Weierstra{\ss}!division globale}
Soient $p\in\N$ et $G \in \As[T]$ un polyn\^ome unitaire de degr\'e $p$. Alors il existe un nombre r\'eel~$v>0$ v\'erifiant la propri\'et\'e suivante : pour toute partie compacte~$U$ de~$B$, tout nombre r\'eel $w\ge v$ et tout \'el\'ement~$F$ de~$\Bs(U)\of{\la}{|T|\le w}{\ra}$, il existe un unique couple $(Q,R) \in (\Bs(U)\of{\la}{|T|\le w}{\ra})^2$ tel que 
\begin{enumerate}[\it i)]
\item $R$ soit un polyn\^ome de degr\'e strictement inf\'erieur \`a $p$ ;
\item $F=QG+R$.
\end{enumerate}
En outre, il existe une constante $C\in\R_{+}^*$, ind\'ependante de $U$, $w$ et $F$, telle que l'on ait les in\'egalit\'es
$$\left\{{\renewcommand{\arraystretch}{1.3}\begin{array}{rcl}
\|Q\|_{U,w} &\le& C\,\|F\|_{U,w}\ ;\\
\|R\|_{U,w} &\le& C\,\|F\|_{U,w}.
\end{array}}\right.$$
\end{thm}
\begin{proof}
Notons $$G= T^p + \sum_{k=0}^{p-1} g_k\, T^k$$ o\`u, quel que soit \mbox{$k\in\cn{0}{p-1}$}, \mbox{$g_k\in \As$}. Soit~$U$ une partie compacte de~$B$. Soit~$u>0$. Tout \'el\'ement $\varphi$ de $\Bs(U)\of{\la}{|T|\le u}{\ra}$ peut s'\'ecrire de fa\c{c}on unique sous la forme 
$$\varphi = \alpha(\varphi)\, T^p + \beta(\varphi),$$
o\`u $\alpha(\varphi)$ d\'esigne un \'el\'ement de $\Bs(U)\of{\la}{|T|\le u}{\ra}$ et $\beta(\varphi)$ un \'el\'ement de $\Bs(U)[T]$ de degr\'e strictement inf\'erieur \`a $p$. Remarquons, d\`es \`a pr\'esent, que, quel que soit $\varphi \in \Bs(U)\of{\la}{|T|\le u}{\ra}$, nous avons 
$$\|\varphi\|_{U,u} = \|\alpha(\varphi)\|_{U,u}\, u^p + \|\beta(\varphi)\|_{U,u}.$$

Consid\'erons, \`a pr\'esent, l'endomorphisme 
$$ A_{U,u} :
{\renewcommand{\arraystretch}{1.2}\begin{array}{ccc}
\Bs(U)\of{\la}{|T|\le u}{\ra} & \to &\Bs(U)\of{\la}{|T|\le u}{\ra} \\
\varphi & \mapsto & \alpha(\varphi)\, G + \beta(\varphi)
\end{array}}.$$
Remarquons que, quel que soit $\varphi \in \Bs(U)\of{\la}{|T|\le u}{\ra}$, nous avons
$${\renewcommand{\arraystretch}{1.3}\begin{array}{rcl}
\|A_{U,u}(\varphi)-\varphi\|_{U,u} &=& \|\alpha(\varphi)\, (G-T^p)\|_{U,u}\\
&\le&  \|\alpha(\varphi)\|_{U,u}\, \|G-T^p\|_{U,u}\\
&\le& u^{-p}\,\|\varphi\|_{U,u}\, \|G-T^p\|_{U,u}\\
&\le& u^{-p}\,\|\varphi\|_{U,u}\, \left(\disp\sum_{k=0}^{p-1} \|g_{k}\|_{U}\, u^k \right)\\
&\le&\, \left(\disp\sum_{k=0}^{p-1} \|g_{k}\|_{B}\,u^{k-p} \right)\,  \|\varphi\|_{U,u}
\end{array}}$$

Il existe~$v>0$ tel que
$$\sum_{k=0}^{p-1} \|g_{k}\|_{B}\,v^{k-p} \le\frac{1}{2}.$$

Soit $w\ge v$. On dispose alors de l'in\'egalit\'e 
$$\|A_{U,w}-I\|_{U,w} \le\frac{1}{2}.$$
Par cons\'equent, l'endomorphisme $A_{U,w}=I+(A_{U,w}-I)$ est inversible.

Soit $F\in \Bs(U)\of{\la}{|\bT|\le w}{\ra}$. Il existe un unique couple $(Q,R)$, avec \mbox{$Q\in\Bs(U)\of{\la}{|T|\le w}{\ra}$} et $R\in\Bs(U)[T]$ de degr\'e strictement inf\'erieur \`a $p$, tel que 
$$F=QG+R.$$
Avec les notations pr\'ec\'edentes, nous avons $Q=\alpha(A_{U,w}^{-1}(F))$ et $R=\beta(A_{U,w}^{-1}(F))$. Puisque $\|A_{U,w}-I\|_{U,w}\le 1/2$, nous avons 
$$\|A_{U,w}^{-1}\|_{U,w}\le \sum_{i=0}^{+\infty} 2^{-i} = 2.$$ 
On en d\'eduit que 
$$\|Q\|_{U,w}\le 2v^{-p}\, \|F\|_{U,w}$$
et que 
$$\|R\|_{U,w}\le 2\, \|F\|_{U,w}.$$
\end{proof}  

Soit~$G(T)$ un polyn\^ome unitaire \`a coefficients dans~$\As$. Notons~$p\in\N$ son degr\'e. Fixons un nombre r\'eel~$w>0$. Soit~$U$ une partie compacte de~$B$. Munissons l'alg\`ebre quotient $\Bs(U)[T]/(G(T))$ de la semi-norme r\'esiduelle $\|.\|_{U,w,\textrm{r\'es}}$ induite par la norme~$\|.\|_{U,w}$ sur~$\Bs(U)[T]/(G(T))$. Par
d\'efinition, quel que soit~$F$ dans~$\Bs(U)[T]/(G(T))$, nous avons
$$\|F\|_{U,w,\textrm{r\'es}}= \inf\left\{ \left\|\sum_{i\in\N} a_{i}\, T^i \right\|_{U,w},\ \sum_{i\in\N} a_{i}\, T^i=F \mod G\right\}.$$
Notons~$v_{0}>0$ le nombre r\'eel dont l'existence nous est assur\'ee par le th\'eor\`eme pr\'ec\'edent. Nous noterons~$C_{0}$ la constante associ\'ee.
\newcounter{nres}\setcounter{nres}{\thepage}

\begin{lem}\label{fininorme}
Pour tout nombre r\'eel~\mbox{$w\ge v_{0}$} et toute partie compacte~$U$ de~$B$, les propri\'et\'es suivantes sont satisfaites : 
\begin{enumerate}[\it i)]
\item la semi-norme~$\|.\|_{U,w,\textrm{r\'es}}$ d\'efinie sur le quotient~$\Bs(U)[T]/(G(T))$ est une norme ;
\item l'anneau~$\Bs(U)[T]/(G(T))$ est complet pour la norme~$\|.\|_{U,w,\textrm{r\'es}}$.
\end{enumerate}
\end{lem}
\begin{proof}
Soient~$w\ge v_{0}$ et~$U$ une partie compacte de~$B$. Le th\'eor\`eme~\ref{vraimentglobal} assure que le morphisme naturel
$$\Bs(U)[T]/(G(T)) \to \Bs(U)\of{\la}{|T|\le w}{\ra}/(G(T))$$
est un isomorphisme. Pour montrer que la semi-norme~$\|.\|_{U,w,\textrm{r\'es}}$ est une norme sur le quotient~$\Bs(U)\of{\la}{|T|\le w}{\ra}/(G(T))$, il suffit de montrer que l'id\'eal~$(G)$ est ferm\'e dans l'anneau~$\Bs(U)\of{\la}{|T|\le w}{\ra}$ pour la norme~$\|.\|_{U,w}$. Soit~$(H_{n})_{n\ge 0}$ une suite d'\'el\'ements de~$\Bs(U)\of{\la}{|T|\le w}{\ra}$ tel que la suite~$(F_{n} = GH_{n})_{n\ge 0}$ converge vers un \'el\'ement~$F$ de~$\Bs(U)\of{\la}{|T|\le w}{\ra}$. D'apr\`es le th\'eor\`eme~\ref{vraimentglobal}, quels que soient~$n,m\ge 0$, nous avons
$$\|F_{n}-F_{m}\|_{U,w} \le C_{0}\, \|H_{n}-H_{m}\|_{U,w}.$$
En particulier, la suite~$(H_{n})_{n\ge 0}$ est de Cauchy dans~$\Bs(U)\of{\la}{|T|\le w}{\ra}$. Cet espace \'etant complet, elle converge vers un \'el\'ement~$H$. Nous avons alors
$$F = GH \in (G),$$
ce qui montre que l'id\'eal~$(G)$ est ferm\'e.
\end{proof}

Soit~$U$ une partie compacte de~$B$. Puisque le polyn\^ome~$G$ est unitaire et de degr\'e~$p$, l'application
$$n : \begin{array}{ccc}
\Bs(U)^p & \to & \Bs(U)[T]/(G(T))\\
(a_{0},\ldots,a_{p-1}) & \mapsto & \disp \sum_{i=0}^{p-1} a_{i}\, T^i 
\end{array}$$
est bijective. Nous noterons encore~$\|.\|_{U}$ la norme d\'efinie sur~$\Bs(U)^d$ en prenant le maximum des normes des coefficients. Nous pouvons alors d\'efinir une norme~$\|.\|_{U,\textrm{div}}$ sur~$\Bs(U)[T]/(G(T))$ par la formule
$$\|.\|_{U,\textrm{div}} = \|n^{-1}(.)\|_{U}.$$

\begin{lem}\label{equivdiv}
Pour tout nombre r\'eel~\mbox{$w\ge v_{0}$} et toute partie compacte~$U$ de~$B$, les normes~$\|.\|_{U,\textrm{div}}$ et~$\|.\|_{U,w,\textrm{r\'es}}$ d\'efinies sur~$\Bs(U)[T]/(G(T))$ sont \'equi\-va\-lentes. En particulier, quels que soient~$w_{1},w_{2}\ge v_{0}$, les normes~$\|.\|_{U,w_{1},\textrm{r\'es}}$ et~$\|.\|_{U,w_{2},\textrm{r\'es}}$ d\'efinies sur~$\Bs(U)[T]/(G(T))$ sont \'equivalentes.
\end{lem}
\begin{proof}
Soient~$w\ge v_{0}$ et~$U$ une partie compacte de~$B$. Soit~$F$ un \'el\'ement de~$\Bs(U)[T]/(G(T))$. Notons~$(f_{0},\ldots,f_{p-1})=n^{-1}(f)$ et~$F_{0} = \sum_{i=0}^{p-1} f_{i}\, T^i$ dans~$\Bs(U)[T]$. L'image de~$F_{0}$ dans~$\Bs(U)[T]/(G(T))$ n'est autre que~$F$. Nous avons donc
$$\|F\|_{U,w,\textrm{r\'es}} \le \|F_{0}\|_{U,w} \le \left(\sum_{i=0}^{d-1} w^i\right)\, \|F\|_{U,\textrm{div}}.$$
Soit~$\eps>0$. Il existe un \'el\'ement~$F_{1}$ de~$\Bs(U)[T]$ d'image~$F$ dans~$\Bs(U)[T]/(G(T))$ tel que l'on ait
$$\|F_{1}\|_{U,w} \le \|F\|_{U,w,\textrm{r\'es}} + \eps.$$
Observons que le reste de la division euclidienne de~$F_{1}$ par~$G$ est \'egal \`a~$F_{0}$. D'apr\`es le th\'eor\`eme de division de Weierstra{\ss}~\ref{vraimentglobal}, nous avons donc
$$\|F_{0}\|_{U,w} \le C_{0}\, \|F_{1}\|_{U,w} \le C_{0}\, (\|F\|_{U,w,\textrm{r\'es}} +\eps).$$
On en d\'eduit que
$$\|F\|_{U,\textrm{div}} \le \max_{1\le i\le p-1} (r^{-i})\, \|F_{0}\|_{U,w} \le  \max_{1\le i\le p-1} (r^{-i})\, C_{0}\, (\|F\|_{U,w,\textrm{r\'es}} +\eps).$$
On obtient le r\'esultat souhait\'e en faisant tendre~$\eps$ vers~$0$.
\end{proof}

Il existe des \'el\'ements~$g_{0},\ldots,g_{p-1}$ de~$\As$ tels que l'on ait
$$G  = T^p + \sum_{k=0}^{p-1} g_{k}\, T^k \textrm{ dans } \As[T].$$
Posons
$$v_{1} = \max_{1\le k\le d} (\|g_{k}\|^{1/(d-k)}).$$
Soit~$b$ un point de~$B$. Nous noterons~$G(b)[T]$ l'image du polyn\^ome~$G(T)$ dans l'anneau~$\Hs(b)[T]$. Rappelons que, d'apr\`es la proposition~3.1.2.1 de~\cite{BGR} l'ensemble des points de la droite~$\E{1}{\Hs(b)}$ en lesquels le polyn\^ome~$G(b)[T]$ s'annule est contenu dans le disque ferm\'e de centre~$0$ et de rayon~$v_{1}$.

Posons~$v=\max(v_{0},v_{1})$. Soit~$w\ge v$. D'apr\`es le lemme~\ref{fininorme}, la semi-norme $\|.\|_{B,w,\textrm{r\'es}}$ d\'efinie sur le quotient~$\Qs_{w} = \As[T]/(G(T))$ est une norme qui rend cet anneau complet. Sa d\'efinition nous assure que le morphisme naturel
$$\As \to  \Qs_{w}$$
est born\'e. Notons
$$\varphi_{w} : C_{w} \to B$$
le morphisme induit entre les spectres analytiques. Remarquons que le morphisme surjectif $\As[T] \to \Qs_{w}$ induit un plongement 
$$C_{w} \hookrightarrow \E{1}{\As}$$
Nous identifierons dor\'enavant~$C_{w}$ \`a son image par ce plongement. Elle est contenue dans le lieu d'annulation du polyn\^ome~$G$ sur la droite~$\E{1}{\As}$. Puisque~\mbox{$w\ge v_{1}$}, nous avons m\^eme \'egalit\'e :
$$C_{w} = \left\{\left.x\in\E{1}{\As}\, \right|\, G(x)=0\right\}.$$

Soit~$U$ une partie compacte de~$B$. Nous noterons~$\Qs_{U,w}$ l'anneau de Banach~$\Bs(U)[T]/(G(T))$ muni de la norme~$\|.\|_{U,w,\textrm{r\'es}}$. Le morphisme naturel
$$\Qs_{w} \to \Qs_{U,w}$$
est born\'e et l'image de l'anneau total des fractions de~$\Qs_{w}$ est dense dans~$\Qs_{U,w}$. Le morphisme 
$$\Ms(\Qs_{U,w}) \to C_{w}$$
induit entre les spectres analytiques est donc injectif et nous identifierons do\-r\'e\-na\-vant l'espace~$\Ms(\Qs_{U,w})$ \`a son image dans~$C_{w}$.

\begin{lem}\label{spconvexe}
Soit~$w\ge v$. Soit~$U$ une partie compacte et spectralement convexe de~$B$. Alors
$$\Ms(\Qs_{U,w})=\varphi_{w}^{-1}(U) = \left\{\left.x\in\E{1}{\As}\, \right|\, \pi(x)\in U,\, G(x)=0\right\},$$
o\`u~$\pi$ d\'esigne la projection naturelle de~$\E{1}{\As}$ sur~$B$.
\end{lem}
\begin{proof}
L'inclusion~$\Ms(\Qs_{U,w})\supset\varphi_{w}^{-1}(U)$ est \'evidente. R\'e\-ci\-pro\-que\-ment, la partie compacte~$U$ est suppos\'ee spectralement convexe. Par d\'efinition, cela signifie que~$\Ms(\Bs(U))=U$. On en d\'eduit que~$\Ms(\Qs_{U,w})$ est contenu dans~$\pi^{-1}(U)$. En outre, en tout point~$x$ de~$\Ms(\Qs_{U,w})$, nous avons~$G(x)=0$. On en d\'eduit le r\'esultat attendu.
\end{proof}

Nous allons, \`a pr\'esent, d\'emontrer un r\'esultat permettant d'assurer que les normes de la forme~$\|.\|_{U,w,\textrm{r\'es}}$ sont uniformes. \`A cet effet, nous introduisons une condition technique. Si~$P$ et~$Q$ sont deux polyn\^omes \`a coefficients dans un anneau~$A$, nous notons~R\'es$(P,Q)\in A$ le r\'esultant des polyn\^omes~$P$ et~$Q$. \index{Resultant@R\'esultant}

\begin{defi}\label{conditionR}\index{Condition $(R_{G})$}
Soit~$U$ une partie compacte de~$B$. Nous dirons que~$U$ v\'erifie la {\bf condition~$\boldsymbol{(R_{G})}$} si elle est spectralement convexe et s'il existe un sous-ensemble~$\Gamma_{U}$ de~$U$ v\'erifiant les propri\'et\'es suivantes :
\begin{enumerate}[\it i)]
\item tout \'el\'ement de~$\Bs(U)$ atteint son maximum sur~$\Gamma_{U}$ ;
\item la fonction~R\'es$(G,G')$ est born\'ee inf\'erieurement sur~$\Gamma_{U}$ par un nombre r\'eel~$m_{U}>0$.
\end{enumerate}
\end{defi}

En pratique, nous utiliserons cette d\'efinition dans les deux cas suivants :
\begin{enumerate}
\item la fonction R\'es$(G,G')$ ne s'annule pas sur~$U$ ;
\item l'ensemble~$\Gamma_{U}$ peut \^etre choisi fini et hors du lieu d'annulation de R\'es$(G,G')$.
\end{enumerate}

\begin{lem}\label{resultant}
Soient~$(k,|.|)$ un corps valu\'e complet. Choisissons une cl\^oture alg\'ebrique~$\bar{k}$ de~$k$ et notons encore~$|.|$ l'unique valeur absolue sur~$\bar{k}$ qui prolonge celle d\'efinie sur~$k$. Soit~$d\in\N$ un entier, $g$ un polyn\^ome \`a coefficients dans~$k$, de degr\'e~$d$, unitaire et s\'eparable. Notons~$\alpha_{1},\ldots,\alpha_{d}$ les racines de~$g$ dans~$\bar{k}$. Soit un nombre r\'eel~$r$ v\'erifiant
$$r\ge\max_{1\le i\le d} (|\alpha_{i}|).$$
Posons
$$D = \frac{ d (2r)^{d^2-d}}{|\textrm{R\'es}(g,g')|}.$$
Alors, quel que soit~$f = \sum_{i=0}^{d-1} a_{i}\, T^i$ dans~$k[T]$, nous avons
$$\sum_{i=0}^{d-1} |a_{i}|\, r^i \le D\, \max_{1\le i\le d} (|f(\alpha_{i})|).$$
\end{lem}
\begin{proof}
Puisque le polyn\^ome~$g$ est s\'eparable, les \'el\'ements~$\alpha_{i}$, avec~\mbox{$i\in\cn{1}{d}$}, sont deux \`a deux distincts. D'apr\`es la formule d'interpolation de Lagrange, dans l'anneau~$\bar{k}[T]$, nous avons donc
$$\begin{array}{rcl}
f(T) &=& \disp \sum_{j=1}^d f(\alpha_{j})\, \prod_{i\ne j} \frac{T-\alpha_{i}}{\alpha_{j}-\alpha_{i}}\\
&=& \disp \frac{1}{\prod_{j=1}^d \prod_{i\ne j} (\alpha_{j}-\alpha_{i})}\,    \sum_{j=1}^d f(\alpha_{j}) \left(\prod_{k\ne j} \prod_{l\ne k} (\alpha_{k}-\alpha_{l})\right) \prod_{i\ne j} (T-\alpha_{i}).
\end{array}$$
On en d\'eduit le r\'esultat annonc\'e.
\end{proof}

\begin{prop}\label{finiuniforme}\index{Norme!uniforme}
Pour tout nombre r\'eel~$w\ge v$ et toute partie compacte~$U$ de~$B$ qui v\'erifie la condition~$(R_{G})$, la semi-norme~$\|.\|_{U,w,\textrm{r\'es}}$ est une norme uniforme sur~$\Bs(U)[T]/(G(T))$. 
\end{prop}
\begin{proof}
Soit~$w\ge v$. Soit~$U$ une partie compacte de~$B$ qui v\'erifie la condition~$(R_{G})$. Nous reprenons les notations de la d\'efinition~\ref{conditionR}. Le lemme~\ref{fininorme} nous assure que  la semi-norme~$\|.\|_{U,w,\textrm{r\'es}}$ est une norme sur~$\Bs(U)[T]/(G(T))$. Notons~$\|.\|_{\infty}$ la norme spectrale associ\'ee. Le lemme~\ref{spconvexe} nous fournit une description explicite de cette norme en termes de norme uniforme sur une partie de la droite~$\E{1}{\As}$. Pour montrer que les deux normes sont \'equivalentes, il suffit de montrer qu'il existe une constante~$D\in\R$ telle que, pour tout \'el\'ement~$F$ de~$\Bs(U)[T]/(G(T))$, nous avons
$$\|F\|_{U,w,\textrm{r\'es}} \le D\, \|F\|_{\infty}.$$

Soit~$F$ un \'el\'ement de~$\Bs(U)[T]/(G(T))$. Puisque le polyn\^ome~$G$ est unitaire et de degr\'e~$p$, l'\'el\'ement~$F$ poss\`ede un unique repr\'esentant dans~$\Bs(U)[T]$ de la forme
$$F_{0}(T) = \sum_{k=0}^{p-1} a_{k}\, T^k,$$
avec~$a_{0},\ldots,a_{p-1}\in\Bs(U)$. 


Soit~$b$ un point de~$\Gamma_{U}$. Le r\'esultant des polyn\^omes~$G(b)(T)$ et~$G'(b)(T)$ n'est autre que l'image~$\textrm{R\'es}(G,G')(b)$ de~$\textrm{R\'es}(G,G')$ dans~$\Hs(b)$. Il suffit, pour s'en convaincre, d'utiliser la d\'efinition du r\'esultant comme d\'eterminant de la matrice de Sylvester. Par hypoth\`ese, l'\'el\'ement~$\textrm{R\'es}(G,G')(b)$ de~$\Hs(b)$ n'est pas nul et le polyn\^ome~$G(b)(T)$ est donc s\'eparable. Notons~$\alpha_{1},\ldots,\alpha_{d}$ ses racines dans une cl\^oture alg\'ebrique de~$\Hs(b)$. Lorsque l'on immerge naturellement la fibre~$\varphi^{-1}(b)$ dans la droite analytique~$\E{1}{\Hs(b)}$, l'image est exactement compos\'ee des points rigides qui correspondent aux classes de conjugaison sous l'action du groupe de Galois des racines~$\alpha_{1},\ldots,\alpha_{d}$. En particulier, nous avons
$$\max_{1\le k\le d} (|F_{0}(\alpha_{j})|) = \max_{x\in \varphi^{-1}(b)} (|F(x)|) \le \|F\|_{\infty}.$$
Remarquons, \`a pr\'esent, que, d'apr\`es~\cite{BGR}, proposition~3.1.2.1 , nous avons
$$\max_{1\le k\le p} (|\alpha_{k}|) \le v_{1}\le w.$$
D'apr\`es le lemme~\ref{resultant}, nous avons donc
$$\sum_{k=0}^{p-1} |a_{k}(b)|\, w^k\le  \frac{p (2w)^{p^2-p}}{|\textrm{R\'es}(G(b),G(b)')|} \, \|F\|_{\infty} \le \frac{p (2w)^{p^2-p}}{m_{U}}\, \|F\|_{\infty}.$$

Pour tout indice~$j\in\cn{1}{p-1}$, choisissons un point~$b_{j}$ de~$\Gamma_{U}$ tel que
$$|a_{j}(b_{j})| = \max_{b\in U} (|a_{j}(b)|).$$
Nous avons alors
$$\begin{array}{rcl}
\|F\|_{U,w,\textrm{r\'es}} &\le& \|F_{0}\|_{U,w}\\
&\le& \disp \sum_{k=0}^{p-1} \|a_{k}\|_{U}\, w^k\\
&\le& \disp \sum_{j=0}^{p-1} \sum_{k=0}^{p-1} |a_{k}(b_{j})|\, w^k\\
&\le& \disp \frac{p^2 (2w)^{p^2-p}}{m_{U}}\, \|F\|_{\infty}.
\end{array}$$
On en d\'eduit que la norme~$\|.\|_{U,w,\textrm{r\'es}}$ est uniforme.
\end{proof}

\section{Un exemple}\label{unexemple}

\index{Morphisme fini!hypersurface d'une droite|(}


Gardons les notations de la section pr\'ec\'edente. Fixons un nombre r\'eel~$w \ge v$. Nous le conserverons tout au long de cette section et nous autoriserons donc \`a supprimer la lettre~$w$ plac\'ee en indice, lorsque cela ne pr\^ete pas \`a confusion. Nous nous int\'eresserons, ici, au morphisme 
$$\varphi : C \to B$$
induit par le morphisme
$$\As \to \As[T]/(G(T)) = \Qs$$
et, plus particuli\`erement, au faisceau~$\varphi_{*}\Os_{C}$. Nous montrerons que, sous certaines hypoth\`eses, c'est un faisceau de $\Os_{B}$-modules libre, comme dans le cadre classique. 

Commen\c{c}ons par montrer que le morphisme~$\varphi$ est un morphisme topologique fini, au sens de la d\'efinition~\ref{topofini}.

\begin{lem}\label{phifini}\index{Morphisme fini!au sens topologique}
Le morphisme~$\varphi$ est un morphisme topologique fini.
\end{lem}
\begin{proof}
Le fait que le morphisme~$\varphi$ soit continu est imm\'ediat. Puisque l'espace~$C$ est compact, on en d\'eduit aussit\^ot que le morphisme~$\varphi$ est \'egalement ferm\'e.

Pour finir, montrons que les fibres du morphisme~$\varphi$ sont finies. Soit~$b$ un point de~$B$. La fibre~$\varphi^{-1}(b)$ est constitu\'ee de l'ensemble des \'el\'ements du disque de centre~$0$ et de rayon~$w$ de la droite~$\E{1}{\Hs(b)}$ en lesquels le polyn\^ome~$H(b)$ s'annule. Puisque ce polyn\^ome est unitaire, il n'est pas nul et l'ensemble~$\varphi^{-1}(b)$ est fini.
\end{proof}

Soit~$b$ un point de~$B$. Notons~$c_{1},\ldots,c_{r}$, avec~$r\in\N^*$, ses ant\'ec\'edents par le morphisme~$\varphi$. Nous supposerons qu'il existe un syt\`eme fondamental~$\Rs$ de voisinages de~$b$ dans~$B$ form\'e de parties compactes qui v\'erifient la condition~$(R_{G})$. 

Soient~$U$ un \'el\'ement de~$\Rs$. Nous allons construire un morphisme
$$\psi_{U} : \Bs(U)[T]/(G(T)) \to \Bs(\varphi^{-1}(U)).$$
Rappelons que l'anneau~$\Bs(\varphi^{-1}(U))$ est un sous-anneau de l'anneau~$\Cs(\varphi^{-1}(U))$ des fonctions
$$f : \varphi^{-1}(U) \to \bigsqcup_{x\in \varphi^{-1}(U)} \Hs(x)$$
qui v\'erifient~$f(x)\in\Hs(x)$, quel que soit~$x\in \varphi^{-1}(U)$. D'apr\`es le lemme~\ref{spconvexe}, nous avons
$$\Ms(\Bs(U)[T]/(G(T))) = \varphi^{-1}(U).$$
Cette remarque nous permet de construire un morphisme
$$\psi_{U,c} : 
\begin{array}{ccc}
\Bs(U)[T]/(G(T)) &\to& \Cs(\varphi^{-1}(U))\\
F & \mapsto & (x\in\varphi^{-1}(U) \mapsto F(x) \in\Hs(x))
\end{array}.$$

\begin{lem}
L'image du morphisme~$\psi_{U,c}$ est contenue dans~$\Bs(\varphi^{-1}(U))$.
\end{lem}
\begin{proof}
Soit~$F$ un \'el\'ement de~$\Bs(U)[T]/(G(T))$. D'apr\`es le th\'eor\`eme~\ref{vraimentglobal}, il existe un \'el\'ement 
$$F_{0} =\sum_{k=0}^{p-1} f_{i}\, T^k \in\Bs(U)[T]$$ 
v\'erifiant les conditions suivantes :
\begin{enumerate}[\it i)]
\item $F =F_{0}$ dans~$\Bs(U)[T]/(G(T))$ ;
\item $\|F_{0}\|_{U,w} \le C\, \|F\|_{U,w,\textrm{r\'es}}$. 
\end{enumerate} 
Soit~$k\in\cn{0}{p-1}$. Par d\'efinition de~$\Bs(U)$, il existe une suite~$(p_{k,n})_{n\ge 0}$ d'\'el\'ements de~$\As$ et une suite~$(q_{k,n})_{n\ge 0}$ d'\'el\'ements de~$\As$ ne s'annulant pas sur~$U$ telles que la suite~$(p_{k,n}/q_{k,n})_{n\ge 0}$ converge vers~$f_{k}$ dans~$\Bs(U)$ pour la norme~$\|.\|_{U}$.

Pour~$n\in\N$, posons
$$P_{n} = \frac{1}{\disp \prod_{0\le k\le p-1} q_{k,n}}\, \sum_{k=0}^{p-1} p_{k,n} \left(\prod_{l\ne k} q_{l,n}\right) T^k \in \Ks(U)[T].$$
Son image modulo~$G(T)$ d\'efinit un \'el\'ement de~$\Ks(\varphi^{-1}(U))$, que nous noterons~$Q_{n}$. Quel que soit~$n\in\N$ et quel que soit~$x\in\varphi^{-1}(U)$, nous avons
$$\begin{array}{rcl}
|Q_{n}(x)-F_{x}| &=& |P_{n}(x)-F_{0}(x)|\\
&\le& \disp \sum_{k=0}^{p-1} \left|\frac{p_{k,n}(x)}{q_{k,n}(x)} - f_{k}(x)\right|\, |T^k(x)|\\
&\le&  \disp \sum_{k=0}^{p-1} \left\|\frac{p_{k,n}}{q_{k,n}} - f_{k}\right\|_{U}\, \|T\|_{\varphi^{-1}(U)}^k.
\end{array}$$
On en d\'eduit que la suite~$(Q_{n})_{n\ge 0}$ d'\'el\'ements de~$\Ks(\varphi^{-1}(U))$ converge vers l'\'el\'ement~$\psi_{U,c}(F)$ pour la norme~$\|.\|_{\varphi^{-1}(U)}$. Par cons\'equent, l'\'el\'ement~$\psi_{U,c}(F)$ appartient \`a~$\Bs(\varphi^{-1}(U))$.

%
\end{proof}

Notons
$$\psi_{U} : \Bs(U)[T]/(G(T)) \to \Bs(\varphi^{-1}(U))$$
le morphisme d\'eduit de~$\psi_{U,c}$ par corestriction.

\begin{prop}
Le morphisme $\psi_{U}$ est un isomorphisme.
\end{prop}
\begin{proof}
D'apr\`es la proposition~\ref{finiuniforme}, la norme~$\|.\|_{U,w,\textrm{r\'es}}$ d\'efinie sur $\Bs(U)[T]/(G(T))$ est \'equivalente \`a sa norme spectrale et, d'apr\`es le lemme~\ref{spconvexe}, cette norme spectrale n'est autre que la norme uniforme sur~$\varphi^{-1}(U)$. Le caract\`ere injectif du morphisme s'en d\'eduit aussit\^ot.

Soit~$F$ un \'el\'ement de~$\Bs(\varphi^{-1}(U))$. Par d\'efinition, il existe une suite~$(P_{n})_{n\ge 0}$ d'\'el\'ements de~$\Qs$ et une suite~$(Q_{n})_{n\ge 0}$ d'\'el\'ements de~$\Qs$ ne s'annulant pas sur~$\varphi^{-1}(U)$ telles que la suite~$(P_{n}/Q_{n})_{n\ge 0}$ converge vers~$F$ pour~$\|.\|_{U,w,\textrm{r\'es}}$. Soit~$n\in\N$. Notons~$P_{U,n}$ et~$Q_{U,n}$ les images respectives de~$P_{n}$ et~$Q_{n}$ dans~$\Qs_{U}$. Par hypoth\`ese, l'\'el\'ement~$Q_{U,n}$ ne s'annule pas sur~$\varphi^{-1}(U)=\Ms(\Qs_{U})$. D'apr\`es le corollaire~1.2.4 de~\cite{rouge}, il est donc inversible dans~$\Qs_{U}$. La suite~$(P_{U,n} Q_{U,n}^{-1})_{n\ge 0}$ de~$\Qs_{U}$ est de Cauchy dans~$\Qs_{U}$. Elle converge donc vers un \'el\'ement de~$\Qs_{U}$ dont l'image par le morphisme~$\psi_{U}$ est l'\'el\'ement~$F$ dont nous sommes partis.
\end{proof}

Nous disposons donc, \`a pr\'esent, d'un isomorphisme
$$\psi_{b} : \varinjlim_{U \in \Rs} \Bs(U)[T]/(G(T))  \xrightarrow[]{\sim} \varinjlim_{U \in \Rs} \Bs(\varphi^{-1}(U))  \xrightarrow[]{\sim} \varinjlim_{U \in \Us} \Bs(\varphi^{-1}(U)),$$
o\`u~$\Us$ d\'esigne l'ensemble des voisinages du point~$b$ dans~$B$. En effet, la premi\`ere fl\`eche est un isomorphisme en vertu de la proposition qui pr\'ec\`ede et la seconde gr\^ace au fait que l'ensemble~$\Rs$ est, par hypoth\`ese, cofinal dans~$\Us$.

Pour tout \'el\'ement~$U$ de~$\Us$, la partie~$\varphi^{-1}(U)$ est un voisinage de la fibre~$\varphi^{-1}(b)$ dans~$C$. Nous disposons donc d'un morphisme de restriction
$$\chi_{b} : \varinjlim_{U \in \Us} \Bs(\varphi^{-1}(U))  \to \prod_{i=1}^r \Os_{C,c_{i}}.$$

\begin{lem}\label{finiinj}
Le morphisme $\chi_{b}$ est injectif.
\end{lem}
\begin{proof}
Ce r\'esultat d\'ecoule directement du lemme~\ref{voisfini}.
\end{proof}

Int\'eressons-nous, \`a pr\'esent, \`a la surjectivit\'e du morphisme~$\chi_{b}$. Pour cela, il nous faut introduire une nouvelle condition. Rappelons que, pour tout point~$b$ de~$B$, nous notons~$\kappa(b)=\Os_{B,b}/\m_{b}$ le corps r\'esiduel du point~$b$ et que le corps~$\Hs(b)$ est son compl\'et\'e pour la valeur absolue associ\'ee \`a~$b$.

\begin{defi}\label{conditionIG}\index{Condition $(I_{G})$}
Nous dirons qu'un point~$b$ de~$B$ satisfait la {\bf condition~$\boldsymbol{(I_{G})}$} si tout facteur irr\'eductible dans~$\kappa(b)[T]$ du polyn\^ome~$G(T)$ reste irr\'eductible dans~$\Hs(b)[T]$.
\end{defi}

Nous supposerons d\'esormais que le point~$b$ satisfait la condition~$(I_{G})$. Remarquons que tel est toujours le cas si le polyn\^ome~$G(b)(T)$ est irr\'eductible (ou, de mani\`ere \'equivalente, si~$r=1$). Dans l'anneau~$\Hs(b)[T]$, \'ecrivons l'image du polyn\^ome~$G(T)$ sous la forme
$$G(b)(T) = \prod_{i=1}^r h_{i}(T)^{n_{i}},$$
o\`u~$r$ est un entier strictement positif, $h_{1},\ldots,h_{r}$ sont des polyn\^omes irr\'eductibles et unitaires \`a coefficients dans~$\Hs(b)$ et $n_{1},\ldots,n_{r}$ des entiers strictement positifs. Les points~$c_{1},\ldots,c_{r}$ de~$C$ sont donc les points de la droite~$\E{1}{\Hs(b)}$ d\'efinis par l'annulation des polyn\^omes~$h_{1},\ldots,h_{r}$. Quitte \`a changer l'ordre des polyn\^omes, nous pouvons supposer que, quel que soit~$i\in\cn{1}{r}$, le point~$c_{i}$ est d\'efini par l'\'equation~$h_{i}=0$.

La condition~$(I_{G})$ assure que la d\'ecomposition en produits de facteurs ir\-r\'e\-duc\-ti\-bles du polyn\^ome~$G(T)$ dans~$\kappa(b)[T]$ et dans~$\Hs(b)[T]$ est identique. On en d\'eduit que, quel que soit~$i\in\cn{1}{r}$, le polyn\^ome~$h_{i}$ est \`a coefficients dans~$\kappa(b)$. D'apr\`es la proposition~\ref{Hensel}, l'anneau local~$\Os_{B,b}$ est hens\'elien. Par cons\'equent, il existe des polyn\^omes~$H_{1},\ldots,H_{r}$ unitaires \`a coefficients dans~$\Os_{B,b}$ qui v\'erifient les propri\'et\'es suivantes :
\begin{enumerate}
\item $G =\disp  \prod_{i=1}^r H_{i}$ dans~$\Os_{B,b}[T]$ ;
\item quel que soit~$i\in\cn{1}{r}$, nous avons~$H_{i}(b) = h_{i}^{n_{i}}$ dans~$\Hs(b)[T]$.\\
\end{enumerate}

\begin{lem}
Il existe un voisinage~$W_{1}$ de~$c_{1}$ dans~$C$, \ldots, un voisinage~$W_{r}$ de~$c_{r}$ dans~$C$ tels que, quel que soit~$j\in\cn{1}{r}$ et quel que soit~$\eps>0$, il existe une fonction~$F_{j,\eps}\in\Ks(W)$, avec  
$$W=\bigcup_{1\le i\le r} W_{i}$$
v\'erifiant les propri\'et\'es suivantes :
\begin{enumerate}[\it i)]
\item $\|F_{j,\eps}-1\|_{W_{j}}\le \eps$ ;
\item quel que soit~$i\ne j$, $\|F_{j,\eps}\|_{W_{i}}\le \eps$.
\end{enumerate}
\end{lem}
\begin{proof}
Il suffit de d\'emontrer le r\'esultat ind\'ependamment pour chacun des indices~$j\in\cn{1}{r}$. Le r\'esultat attendu s'en d\'eduit en consid\'erant, pour chaque indice~$i\in\cn{1}{r}$, l'intersection des ouverts~$W_{i}$ construits et en restreignant les fonctions~$F_{j,\eps}$. 


Soit~$j\in\cn{1}{r}$. Il existe un voisinage~$V$ de~$b$ dans~$B$ tel que les fonctions $H_{1},\ldots,H_{r}$ appartiennent \`a~$\Bs(V)[T]$. Choisissons des voisinages compacts~$W_{1}$ de~$c_{1}$, \ldots, $W_{r}$ de~$c_{r}$ dans~$\varphi^{-1}(V)$, deux \`a deux disjoints. Quitte \`a restreindre ces voisinages, nous pouvons supposer que, quel que soit~$k\in\cn{1}{r}$, la fonction~$H_{k}$ ne s'annule pas sur la partie compacte~$W_{i}$, pour~$i\ne k$. Il existe alors deux nombres r\'eels~$m,M>0$ tels que, quel que soit~$k\in\cn{1}{r}$ et quel que soit~$i\ne k$, nous ayons
$$m < \min_{x\in W_{i}}(|H_{k}(x)|) \le \|H_{k}\|_{W_{i}} < M.$$
Remarquons que nous pouvons restreindre les voisinages~$W_{i}$, avec~$i\in\cn{1}{r}$, sans changer les valeurs des constantes~$m$ et~$M$. En particulier, nous pouvons supposer que nous avons
$$\|H_{j}\|_{W_{j}} < \frac{1}{2}\, m^{r-1}$$
et, quel que soit~$i\ne j$, 
$$\|H_{i}\|_{W_{i}} < \frac{1}{2}\, m M^{2-r}.$$
Par densit\'e de~$\Ks(W)$ dans~$\Bs(W)\supset\Bs(V)[T]$, nous pouvons supposer qu'il existe des \'el\'ements \mbox{$K_{1},\ldots,K_{r}$} de~$\Ks(W)$ qui v\'erifient les m\^emes in\'egalit\'es que $H_{1},\ldots,H_{r}$. 

Soit~$N\in\N^*$. Montrons que la fonction~$D_{N} = K_{j}^N + \prod_{i\ne j} K_{i}^N$ ne s'annule pas sur~$W$. Sur~$W_{j}$, tout d'abord, nous avons
$$\min_{x\in W_{j}} \left( \prod_{i\ne j} |K_{i}^N(x)|\right) \ge  \prod_{i\ne j} \min_{x\in W_{j}} (|K_{i}(x)|^N) \ge m^{N(r-1)}$$ 
et
$$\|K_{j}^N\|_{W_{j}} \le 2^{-N} m^{N(r-1)}.$$
On en d\'eduit que
$$\min_{x\in W_{j}} (|D_{N}(x)|) \ge (1-2^{-N}) m^{N(r-1)}\ge \frac{m^{N(r-1)}}{2}.$$
Soit~$i\ne j$. Nous avons
$$\min_{x\in W_{i}}(|K_{j}^N(x)|) \ge m^N$$
et
$$\left\| K_{i}^N\, \prod_{k\ne i,j} K_{k}^N\right\|_{W_{i}} \le 2^{-N}\, m^N M^{N(2-r)} M^{N(r-2)} \le 2^{-N}\, m^N.$$
On en d\'eduit que
$$\min_{x\in W_{i}} (|D_{N}(x)|) \ge (1-2^{-N}) m^{N}\ge \frac{m^{N}}{2}.$$
En particulier, l'\'el\'ement~$D_{N}$ de~$\Ks(W)$ est inversible.

Consid\'erons l'\'el\'ement~$F_{N} = D_{N}^{-1}\, \prod_{i\ne j} K_{i}^N$ de~$\Ks(W)$. Il v\'erifie
$$\|F_{N}-1\|_{W_{j}} = \|D_{N}^{-1}\, K_{j}^N\|_{W_{j}} \le 2 m^{-N(r-1)}\, 2^{-N} m^{N(r-1)} \le 2^{1-N}$$
et, quel que soit~$i\ne j$, 
$$\|F_{N}\|_{W_{i}} \le 2m^{-N}\, 2^{-N}m^NM^{N(2-r)}\, M^{N(r-2)} \le 2^{1-N}.$$
Quel que soit~$\eps>0$, quitte \`a choisir un nombre entier~$N$ assez grand, l'\'el\'ement~$F_{N}$ v\'erifie les propri\'et\'es demand\'ees.
\end{proof}

\begin{lem}\label{finisurj}
Le morphisme~$\chi_{b}$ est surjectif.
\end{lem}
\begin{proof}
Il suffit de montrer que, quel que soit~$i\in\cn{1}{r}$ et quel que soit~$f$ dans~$\Os_{C,c_{i}}$, il existe un \'el\'ement~$F$ de~$\varinjlim_{U \in \Us} \Bs(\varphi^{-1}(U))$ dont l'image dans~$\Os_{C,c_{i}}$ est \'egale \`a~$f$ et l'image dans~$\Os_{C,c_{j}}$, pour tout~$j\ne i$, est nulle.

Soient~$i\in\cn{1}{r}$ et~$f\in\Os_{C,c_{i}}$. Il existe un voisinage~$V_{i}$ de~$c_{i}$ dans~$C$ sur lequel la fonction~$f$ est d\'efinie. Quitte \`a restreindre ce voisinage, nous pouvons supposer qu'il existe une suite~$(p_{n})_{n\ge 0}$ d'\'el\'ements de~$\Qs$ et une suite~$(q_{n})_{n\ge 0}$ d'\'el\'ements de~$\Qs$ qui ne s'annulent pas sur~$V_{i}$ tels que la suite~$(p_{n}/q_{n})_{n\ge 0}$ converge vers~$f$ pour la norme~$\|.\|_{V_{i}}$. 

Nous reprenons, \`a pr\'esent, les notations du lemme pr\'ec\'edent. Quitte \`a restreindre le voisinage~$V_{i}$, nous pouvons supposer qu'il est compact et contenu dans~$W_{i}$. Nous noterons
$$V = V_{i} \cup \left(\bigcup_{j\ne i} W_{j}\right).$$
Soit~$n\in\N$. Il existe un nombre r\'eel~$m>0$ tel que, quel que soit~$x$ dans~$V_{i}$, nous ayons~$|q_{n}(x)|>m$ et un nombre r\'eel~$M\ge m$ tel que~$\|p_{n}\|_{V}\le M$ et~$\|q_{n}\|_{V}\le M$. Posons
$$\lambda_{n} =\frac{m}{2^n(r+1)M}\le \frac{1}{r+1} \textrm{ et } \mu_{n} = \frac{\lambda_{n}}{M}\le\frac{1}{M(r+1)}.$$ 
Consid\'erons l'\'el\'ement de~$\Ks(V)$ d\'efini par
$$Q_{n} = F_{i,\mu_{n}}q_{n} + \sum_{j\ne i} F_{j,\lambda_{n}}.$$
Montrons qu'il ne s'annule pas sur~$W$. Quel que soit~$x\in V_{i}$, nous avons
$$|Q_{n}(x)| \ge |F_{i,\mu_{n}}(x)q_{n}(x)| - \sum_{j\ne i} |F_{j,\lambda_{n}}(x)| \ge 1 - \frac{r}{r+1} \ge \frac{1}{r+1}.$$
Soit~$j\ne i$. Quel que soit~$x\in W_{j}$, nous avons
$$|Q_{n}(x)| \ge |F_{j,\lambda_{n}}(x)| - |F_{i,\mu_{n}}(x)q_{n}(x)| - \sum_{k\ne i,j} |F_{k,\lambda_{n}}(x)| \ge \frac{1}{r+1}.$$
On en d\'eduit que~$Q_{n}$ est inversible dans~$\Ks(V)$. Consid\'erons, \`a pr\'esent, l'\'el\'ement de~$\Ks(V)$ d\'efini par
$$R_{n} = F_{i,\mu_{n}}p_{n}Q_{n}^{-1}.$$
Quel que soit~$j\ne i$, nous avons
$$\|R_{n}\|_{W_{j}} \le \frac{m}{2^nM^2(r+1)}\, M(r+1)\le \frac{m}{2^nM} \le \frac{1}{2^n}.$$
Au-dessus de~$V_{i}$, nous avons
$$R_{n} - \frac{p_{n}}{q_{n}} = \frac{p_{n}}{Q_{n}q_{n}}\, (F_{i,\mu_{n}}q_{n} - Q_{n}) =  \frac{p_{n}}{Q_{n}q_{n}}\, \sum_{j\ne i} F_{j,\lambda_{n}}.$$
On en d\'eduit que
$$\left\|R_{n} - \frac{p_{n}}{q_{n}}\right\|_{V_{i}} \le \frac{M(r+1)}{m}\, \frac{(r-1)m}{2^n(r+1)M} \le \frac{r-1}{2^n}.$$

On d\'eduit de ces in\'egalit\'es que la suite~$(R_{n})_{n\ge 0}$ converge pour la norme~$\|.\|_{V}$ vers la fonction qui co\"{\i}ncide avec~$f$ sur~$V_{i}$ et qui est nulle sur~$W_{j}$, quel que soit~$j\ne i$.
\end{proof}

Venons-en, maintenant, \`a la description de l'anneau local~$(\varphi_{*}\Os_{C})_{b}$. Il nous suffit pour cela de regrouper les r\'esultats obtenus pr\'ec\'edemment.

\index{Morphisme fini!image directe du faisceau structural|(}

\begin{thm}\label{casreduit}
Soit~$b$ un point de~$B$. Supposons que le point~$b$  v\'erifie la condition~$(I_{G})$ et poss\`ede un syst\`eme fondamental de voisinages compacts qui satisfont la condition~$(R_{G})$. Alors le morphisme
$$\alpha_{b} : \begin{array}{ccc}
\Os_{B,b}^p & \to & (\varphi_{*}\Os_{C})_{b}\\
(a_{0},\ldots,a_{p-1}) & \mapsto & \disp \sum_{i=0}^p a_{i}\, T^i
\end{array}$$
est un isomorphisme de~$\Os_{B,b}$-modules.
\end{thm}
\begin{proof}
Notons
$$\beta_{b} : \begin{array}{ccc}
\Os_{B,b}^p & \to & \Os_{B,b}[T]/(G(T)) \\
(a_{0},\ldots,a_{p-1}) & \mapsto & \disp \sum_{i=0}^p a_{i}\, T^i
\end{array}.$$
C'est un isomorphisme, car le polyn\^ome~$G(T)$ est unitaire et de degr\'e~$p$.

Notons~$\gamma_{b}$ le morphisme naturel
$$\gamma_{b} :  \Os_{B,b}[T]/(G(T))  \to \varinjlim_{U\in\Rs} \Bs(U)[T]/(G(T)).$$
Il est bien d\'efini car~$\Rs$ est, par hypoth\`ese, un syst\`eme fondamental de voisinages du point~$b$ dans~$B$ et c'est \'egalement un isomorphisme.

Notons~$\delta_{b}$ le morphisme induit par la restriction
$$\delta_{b} : \prod_{i=1}^r \Os_{C,c_{i}}  \to  (\varphi_{*}\Os_{C})_{b}.$$
D'apr\`es le th\'eor\`eme~\ref{finifibres}, c'est encore un isomorphisme.

Avec les notations pr\'ec\'edentes, le morphisme~$\alpha_{b}$ se d\'ecompose de la fa\c{c}on suivante :
$$\alpha_{b} = \delta_{b} \circ \chi_{b} \circ \psi_{b} \circ \gamma_{b} \circ \beta_{b}.$$
Nous avons d\'emontr\'e plus haut que les morphismes~$\chi_{b}$ et~$\psi_{b}$ sont des isomorphismes. On en d\'eduit le r\'esultat attendu.


\end{proof}

Nous tirons imm\'ediatement les cons\'equences de ce r\'esultat en termes globaux.

\begin{cor}\label{casreduitglobal}
Supposons que tout point de~$B$ v\'erifie la condition~$(I_{G})$ et poss\`ede un syst\`eme fondamental de voisinages compacts qui satisfont la condition~$(R_{G})$. Alors, le morphisme
$$\alpha : \begin{array}{ccc}
\Os_{B}^p & \to & \varphi_{*}\Os_{C}\\
(a_{0},\ldots,a_{p-1}) & \mapsto & \disp \sum_{i=0}^p a_{i}\, T^i
\end{array}$$
est un isomorphisme de~$\Os_{B}$-modules. En particulier, pour toute partie~$V$ de~$B$, le morphisme naturel
$$\Os_{B}(V)[T]/(G(T)) \to \Os_{C}(\varphi^{-1}(V))$$
est un isomorphisme.
\end{cor}
\begin{proof}
La premi\`ere partie du r\'esultat d\'ecoule imm\'ediatement du th\'eor\`eme pr\'ec\'edent. On en d\'eduit que, pour toute partie~$V$ de~$B$, le morphisme naturel
$$\Os_{B}(V)[T]/(G(T)) \to (\varphi_{*}\Os_{C})(V)$$
est un isomorphisme. Il nous reste \`a remarquer que le morphisme naturel
$$ (\varphi_{*}\Os_{C})(V) = \varinjlim_{\overset{U \supset V}{U\textrm{ ouvert}}} \Os_{C}(\varphi^{-1}(U)) \to \varinjlim_{\overset{W \supset \varphi^{-1}(V)}{W\textrm{ ouvert}}} \Os_{C}(W) =  \Os_{C}(\varphi^{-1}(V))$$
est un isomorphisme. En effet, d'apr\`es le lemme~\ref{phifini}, le morphisme~$\varphi$ est un morphisme topologique fini. Il suffit alors d'appliquer le corollaire~\ref{corvoisfini}.
\end{proof}

\begin{cor}\label{casreduitStein}
Supposons que tout point de~$B$ v\'erifie la condition~$(I_{G})$ et poss\`ede un syst\`eme fondamental de voisinages compacts qui satisfont la condition~$(R_{G})$. Supposons que le faisceau~$\Os_{B}$ est coh\'erent. Pour toute partie~$V$ de~$B$, nous noterons
$$\varphi_{V} : \varphi^{-1}(V) \to V$$
le morphisme d\'eduit de~$\varphi$ par restriction et corestriction. Alors, pour toute partie~$V$ de~$B$ et tout faisceau coh\'erent~$\Fs$ sur~$\varphi^{-1}(V)$, le faisceau~$(\varphi_{V})_{*}\Fs$ est coh\'erent.
\end{cor}
\begin{proof}
D'apr\`es le corollaire~\ref{casreduitglobal}, le faisceau~$\varphi_{*}\Os_{C}$ est isomorphe au faisceau~$\Os_{B}^p$. C'est donc un faisceau coh\'erent. Soient~$V$ une partie de~$B$ et~$\Fs$ un faisceau coh\'erent sur~$\varphi^{-1}(V)$. Soit~$b$ un point de~$V$. Notons~$c_{1},\ldots,c_{r}$, avec~$r\in\N$, ses ant\'ec\'edents par le morphisme~$\varphi$. Ils sont en nombre fini, d'apr\`es le lemme~\ref{phifini}. Soit~$i\in\cn{1}{r}$. Il existe un voisinage~$U_{i}$ du point~$c_{i}$ dans~$\varphi^{-1}(V)$, des entiers~$p_{i}$ et~$q_{i}$ et une suite exacte
$$0 \to \Os_{U_{i}}^{p_{i}} \to  \Os_{U_{i}}^{q_{i}} \to \Fs_{U_{i}} \to 0.$$
Nous pouvons supposer que les entiers~$p_{i}$, avec~$i\in\cn{1}{r}$, sont \'egaux \`a un m\^eme entier~$p$, que les entiers~$q_{i}$, avec~$i\in\cn{1}{r}$ sont \'egaux \`a m\^eme entier~$q$ et que les voisinages~$U_{i}$, avec~$i\in\cn{1}{r}$, sont deux \`a deux disjoints. Notons~$U$ leur r\'eunion. D'apr\`es le lemme~\ref{voisfini}, quitte \`a restreindre encore les voisinages~$U_{i}$, nous pouvons supposer que la partie~$U$ est de la forme~$\varphi^{-1}(W)$, o\`u~$W$ est un voisinage du point~$b$ dans~$V$. Nous pouvons regrouper les suites pr\'ec\'edentes en une suite exacte
$$0 \to \Os_{U}^{p} \to  \Os_{U}^{q} \to \Fs_{U} \to 0.$$
D'apr\`es le th\'eor\`eme~\ref{finisuite}, la suite
$$0 \to ((\varphi_{W})_{*}\Os_{U})^{p} \to  ((\varphi_{W})_{*}\Os_{U})^{q} \to (\varphi_{W})_{*}\Fs_{U} \to 0$$
est encore exacte. D'apr\`es le corollaire~\ref{corvoisfini}, le faisceau~$(\varphi_{W})_{*}\Os_{U}$ est la restriction \`a~$W$ du faisceau~$\varphi_{*}\Os_{C}$. C'est donc un faisceau coh\'erent. On en d\'eduit que le faisceau~$(\varphi_{W})_{*}\Fs_{U}$, qui n'est autre que la restriction \`a~$W$ du faisceau~$(\varphi_{V})_{*}\Fs$, d'apr\`es le m\^eme corollaire, est coh\'erent. Par cons\'equent, le faisceau~$(\varphi_{V})_{*}\Fs$ est coh\'erent.
%
\end{proof}

\index{Morphisme fini!image directe du faisceau structural|)}

\index{Morphisme fini!hypersurface d'une droite|)}

\section{Th\'eor\`eme de division de {Weierstra\ss} en un point rigide}\label{Wrigide}

Le th\'eor\`eme de division de {Weierstra\ss} \ref{division} que nous avons d\'emontr\'e g\'en\'eralise le th\'eor\`eme classique sur~$\C$ et permet de d\'ecrire l'anneau local au voisinage du point~$0$ d'une fibre. En g\'eom\'etrie analytique complexe, il est toujours possible de se ramener \`a ce cas \`a l'aide d'une translation. Sur une base quelconque, en revanche, un tel artifice est impossible. Nous allons cependant montrer ici que l'\'etude des morphismes finis que nous avons entreprise permet d'obtenir un th\'eor\`eme de division de {Weierstra\ss} pour les points rigides des fibres.

\bigskip

Soit~$(\As,\|.\|)$ un anneau de Banach uniforme. Nous notons~$B=\Ms(\As)$, \mbox{$X=\E{1}{\As}$} (avec variable~$T$) et~$\pi : X\to B$ le morphisme de projection. Soit \mbox{$s>0$}. Consid\'erons l'alg\`ebre $\As\of{\la}{|T|\le s}{\ra}$ munie de la norme~$\|.\|_{s}$. Nous noterons~$\As_{s}$ son compl\'et\'e pour la norme uniforme sur son spectre analytique. Le morphisme $\As[T]\to \As_{s}$ induit une application continue et injective
$$\Ms(\As_{s}) \hookrightarrow \E{1}{\As}$$
dont l'image est le disque ferm\'e~$\overline{D}(s)$. Nous identifierons dor\'enavant le spectre analytique~$\Ms(\As_{s})$ \`a ce disque. 
\newcounter{Ass}\setcounter{Ass}{\thepage}

\begin{defi}\label{conditionS}\index{Condition $(S)$}\index{Condition $(S_{P})$}
Soient~$b$ un point de~$B$ et~$P(T)$ un polyn\^ome \`a coefficients dans~$\As$ unitaire dont l'image dans~$\Hs(b)[T]$ est irr\'eductible. Nous dirons que le point~$b$ satisfait la {\bf condition~$\boldsymbol{(S_{P})}$} s'il existe un nombre r\'eel~$s>0$ et un syst\`eme fondamental~$\Us_{b,P}$ de voisinages compacts et spectralement convexes de~$b$ dans~$B$ tel que, quel que soient~$U\in\Us_{b,P}$ et~$r\in\of{]}{0,s}{]}$, la partie compacte~$\overline{D}_{U}(r)$ de~$\Ms(\As_{s})$ v\'erifie la condition~$(R_{P(S)-T})$.

Nous dirons qu'un point~$b$ de~$B$ satisfait la {\bf condition~$\boldsymbol{(S)}$} si, pour tout polyn\^ome unitaire~$P(T)$ \`a coefficients dans~$\Os_{B,b}$ dont l'image dans~$\Hs(b)[T]$ est irr\'eductible, il existe un voisinage compact et spectralement convexe~$V$ du point~$b$ dans~$B$ sur lequel le polyn\^ome~$P(T)$ est d\'efini et tel que le point~$b$ de~$\Ms(\Bs(V))$ satisfasse la condition~$(S_{P})$.
\end{defi}



\begin{rem}\label{remspconvexe}\index{Compact!spectralement convexe}
D'apr\`es la proposition \ref{stabilitespconvexe}, si~$U$ d\'esigne une partie compacte et spectralement convexe de~$B$ et~$r$ et~$s$ deux nombres r\'eels v\'erifiant~$0<r\le s$, alors la partie compacte~$\overline{D}_{U}(r)$ de~$\Ms(\As_{s})$ est spectralement convexe.
\end{rem}




Soient~$b$ un point de~$B$ et~$P(T)$ un polyn\^ome \`a coefficients dans~$\As$ unitaire et dont l'image dans~$\Hs(b)[T]$ est irr\'eductible. Nous noterons~$y$ l'unique point de la fibre~$\pi^{-1}(b)$ d\'efini par l'\'equation~$P=0$. Nous allons d\'ecrire l'anneau local de la droite analytique en ce point en nous ramenant \`a la situation d\'ecrite dans les sections qui pr\'ec\`edent. Pour cela, nous supposerons que le point~$b$ v\'erifie la condition~$(S_{P})$.\\

La condition~$(S_{P})$ assure qu'il existe un nombre r\'eel strictement positif~$s$ et un syst\`eme fondamental~$\Us_{b,P}$ de voisinages compacts et spectralement convexes de~$b$ dans~$B$ tel que, quel que soient~$U\in\Us_{b,P}$ et~$r\in\of{]}{0,s}{]}$, la partie compacte~$\overline{D}_{U}(r)$ de~$\Ms(\As_{s})$ v\'erifie la condition~$(R_{P(S)-T})$. Soit~$U$ un \'el\'ement de~$\Us_{b,P}$. Consid\'erons l'alg\`ebre de Banach~$\As' = \Bs(\overline{D}_{U}(s))$ et munissons-la de sa norme uniforme~$\|.\|'$. Notons~$B'$ son spectre analytique. Consid\'erons le nombre r\'eel~$v>0$ dont l'existence est d\'emontr\'ee dans la section~\ref{Wglobal} pour l'alg\`ebre~$\As'$ et le polyn\^ome~$P(S)-T$ de~$\As'[S]$. Choisissons un nombre r\'eel~$v'>v$. Notons~$x\in B'$ le point de la fibre au-dessus de~$b$ d\'efini par l'\'equation~$T=0$. Nous avons un isomorphisme $\Hs(b)\xrightarrow[]{\sim} \Hs(x)$. Par hypoth\`ese, la partie~$B'$ v\'erifie la condition~$(R_{P(S)-T})$. Par cons\'equent, d'apr\`es la proposition~\ref{finiuniforme}, la semi-norme~$\|.\|'_{B',v',\textrm{r\'es}}$ d\'efinie sur le quotient~$\As'[S]/(P(S)-T)$ est une norme uniforme. Notons~$C'$ le spectre analytique de~$\As'[S]/(P(S)-T)$ et~$\varphi' : C' \to B'$ le morphisme naturel. Puisque le polyn\^ome~$P(S)$ est irr\'eductible dans~$\Hs(x)[S]$, la fibre~${\varphi'}^{-1}(x)$ ne comporte qu'un seul point. C'est le point rigide de la fibre au-dessus de~$b$ de l'espace affine de dimension~$2$ associ\'e \`a l'id\'eal maximal~$(P(S),T)$. Nous noterons~$z$ ce point. Remarquons que, par choix de~$v'$, la partie~$C'$ est un voisinage du point~$z$ dans le ferm\'e de Zariski de~$\E{2}{\As}$ d\'efini par l'\'equation~$P(S)-T=0$.

Puisque le polyn\^ome~$(P(S)-T)(x)=P(S)(x) \in\Hs(x)[S]$ est irr\'eductible, le point~$x$ de~$X=\E{1}{\As}$ (avec variable~$T$) satisfait la condition~$(I_{P(S)-T})$. D'apr\`es la proposition~\ref{voisdep}, l'ensemble des parties de la forme~$\overline{D}_{V}(r)$, avec~$V\in\Us_{b,P}$ et~\mbox{$r>0$}, est un syst\`eme fondamental de voisinages compacts du point~$x$ dans~$X$. Quitte \`a ne consid\'erer les parties pr\'ec\'edentes que sous les conditions~$V\subset U$ et~$r\le 1$, nous obtenons un syst\`eme fondamental de voisinages compacts et spectralement convexes du point~$x$ dans~$B'$. Par hypoth\`ese, ces parties satisfont la condition~$(R_{P(S)-T})$. Nous pouvons donc appliquer le th\'eor\`eme~\ref{casreduit}. Il assure que le morphisme naturel
$$\Os_{X,x}[S]/(P(S)-T) \xrightarrow[]{\sim} \Os_{B',x}[S]/(P(S)-T) \to \Os_{C',z}$$
est un isomorphisme.

\bigskip

Consid\'erons, \`a pr\'esent, l'alg\`ebre $\Bs(U)\of{\la}{|S|\le v'}{\ra}$ munie de la norme~$\|.\|_{U,v'}$. C'est une alg\`ebre de Banach dont nous noterons~$B''$ le spectre analytique. Notons~$\|.\|''$ la norme uniforme sur~$B''$ et~$\As''$ l'alg\`ebre de Banach obtenue en compl\'etant l'alg\`ebre $\Bs(U)\of{\la}{|S|\le v'}{\ra}$ pour cette norme. Consid\'erons le nombre r\'eel~$v>0$ dont l'existence est d\'emontr\'ee dans la section~\ref{Wglobal} pour l'alg\`ebre~$\As''$ et le polyn\^ome~$T-P(S)$ de~$\As''[S]$. Choisissons un nombre r\'eel~$v''\ge\max(v,1)$. Remarquons que la condition~$(R_{T-P(S)})$ est trivialement v\'erifi\'ee pour tout partie compacte et spectralement convexe de~$B''$ et, en particulier, pour la partie~$B''$ elle-m\^eme. D'apr\`es la proposition~\ref{finiuniforme}, la semi-norme~$\|.\|''_{B'',v'',\textrm{r\'es}}$ d\'efinie sur le quotient~$\As''[T]/(T-P(S))$ est une norme uniforme. Notons~$C''$ le spectre analytique de~$\As''[T]/(T-P(S))$ et~$\varphi'' : C'' \to B''$ le morphisme naturel. Puisque le polyn\^ome~$T-P(S)\in\As''[T]$ est de degr\'e~$1$, la fibre~${\varphi''}^{-1}(y)$ ne comporte qu'un seul point. C'est le point rigide de la fibre au-dessus de~$b$ de l'espace affine de dimension~$2$ associ\'e \`a l'id\'eal maximal~$(P(S),T)$, comme pr\'ec\'edemment. Nous noterons donc encore~$z$ ce point. 

Le point~$y$ de~$B''$ (avec variable~$S$) satisfait \'evidemment la condition~$(I_{T-P(S)})$. La remarque~\ref{remspconvexe} montre que le point~$y$ de~$B''$ poss\`ede un syst\`eme fondamental de voisinages compacts qui satisfont la condition~$(R_{T-P(S)})$. Il suffit, par exemple, de consid\'erer l'ensemble des voisinages compacts rationnels du point~$y$. Nous pouvons donc appliquer le th\'eor\`eme~\ref{casreduit}. On en d\'eduit que le morphisme naturel
$$\Os_{Y,y} \xrightarrow[]{\sim} \Os_{B'',y} \to \Os_{C'',z}$$
est un isomorphisme.

\bigskip

Pour finir, remarquons que les parties~$C'$ et~$C''$ se plongent naturellement dans l'espace affine de dimension~$2$ au-dessus de~$B$. Par choix de~$v''$, une fois identifi\'es les espaces et leur plongement, nous avons l'inclusion~$C'\subset C''$. On en d\'eduit qu'en tout point~$c$ int\'erieur \`a~$C'$, le morphisme de restriction
$$\Os_{C'',c} \to \Os_{C',c}$$
est un isomorphisme. En effet, en un tel point, l'anneau local est form\'e des fonctions qui sont localement limites uniformes de fractions rationnelles sans p\^oles \`a coefficients dans~$\As$. En particulier, nous avons un isomorphisme
$$\Os_{C'',z} \xrightarrow[]{\sim} \Os_{C',z}.$$

Il ne nous reste plus, \`a pr\'esent, qu'\`a combiner ces r\'esultats pour obtenir une description explicite de l'anneau local~$\Os_{Y,y}$. 

\begin{thm}\label{anneaulocalrigide}\index{Anneau local en un point!rigide de $\E{1}{\As}$}
Sous la condition~$(S_{P})$, le morphisme naturel
$$\Os_{X,x}[S]/(P(S)-T) \to \Os_{C',z} \xleftarrow[]{\sim} \Os_{Y,y}$$
est un isomorphisme.
\end{thm}

Forts de cette description, nous pouvons, \`a pr\'esent, d\'emontrer un th\'eor\`eme de division de {Weierstra\ss} au voisinage des points rigides des fibres de la droite analytique. Rappelons les notations : $\As$ est un anneau de Banach muni d'une norme uniforme, $B=\Ms(\As)$ est son spectre analytique, $b$ est un point de~$B$, $P(S)$ est un polyn\^ome unitaire \`a coefficients dans~$\As$ dont l'image dans~$\Hs(b)[S]$ est irr\'eductible, $Y=\E{1}{\As}$ est la droite analytique au-dessus de~$B$ (nous notons~$S$ la variable correspondante) et~$y$ est l'unique point de la fibre au-dessus de~$b$ d\'efini par l'\'equation~$P(y)=0$. 

\begin{thm}\label{divisionrigide}\index{Theoreme de Weierstrass@Th\'eor\`eme de Weierstra{\ss}!division au voisinage d'un point rigide}
Supposons que le point~$b$ de~$B$ satisfait la condition~$(S)$. Soit~$G(S)$ un polyn\^ome \`a coefficients dans~$\Os_{B,b}$. Notons~$n$ la valuation~$P$-adique de l'image de ce polyn\^ome dans~$\Hs(b)[S]$. Alors, pour tout \'el\'ement~$F$ de~$\Os_{Y,y}$, il existe un unique couple~$(Q,R)$ d'\'el\'ements de~$\Os_{Y,y}$ v\'erifiant les propri\'et\'es suivantes :
\begin{enumerate}[\it i)]
\item l'\'el\'ement~$R$ est un polyn\^ome \`a coefficients dans~$\Os_{B,b}$ de degr\'e strictement inf\'erieur \`a~$nd$ ;
\item nous avons l'\'egalit\'e~$F=QG+R$.
\end{enumerate} 

En outre, si l'\'el\'ement~$F$ appartient \`a~$\Os_{B,b}[S]$, il en est de m\^eme pour les \'el\'ements~$Q$ et~$R$.
\end{thm}
\begin{proof}
On se ram\`ene imm\'ediatement \`a traiter le m\^eme probl\`eme dans l'anneau local~$\Os_{C',z}$ et avec le polyn\^ome~$G=P^n$. Le r\'esultat se d\'eduit alors simplement du th\'eor\`eme de division de {Weierstra\ss} classique dans l'anneau~$\Os_{X,x}$. D\'etaillons la preuve de l'existence de la division. Nous disposons d'un isomorphisme
$$\Os_{X,x}[S]/(P(S)-T)  \xrightarrow[]{\sim} \Os_{C',z}.$$
Puisque le polyn\^ome~$P(S)-T$ de~$\Os_{X,x}[S]$ est unitaire, il existe des \'el\'ements $f_{0},\ldots,f_{d-1}$ de~$\Os_{X,x}$ tels que l'on ait
$$F = \sum_{i=0}^{d-1} f_{i}\, S^i \mod (P(S)-T).$$
Soit~$i\in\cn{1}{d-1}$. D'apr\`es le th\'eor\`eme~\ref{division}, il existe un \'el\'ement~$q_{i}$ de~$\Os_{X,x}$ et un polyn\^ome~$r_{i}(T)$ \`a coefficients dans~$\Os_{B,b}$ de degr\'e inf\'erieur \`a~$n-1$ tels que l'on ait l'\'egalit\'e
$$f_{i} = q_{i}T^n + r_{i}.$$
Par cons\'equent, dans l'anneau~$\Os_{X,x}[S]/(P(S)-T)$, nous avons
$$\begin{array}{rcl}
F &=&\disp \left(\sum_{i=0}^{d-1} a_{i}q_{i}S^i\right)\, T^n + \sum_{i=0}^{d-1} r_{i}(T)S^i\\
&=& \disp  \left(\sum_{i=0}^{d-1} a_{i}q_{i}S^i\right)\, T^n + \sum_{i=0}^{d-1} r_{i}(P(S))S^i
\end{array}.$$
Pour conclure, il nous suffit de remarquer que le degr\'e du polyn\^ome 
$$\sum_{i=0}^{d-1} r_{i}(P(S))\, S^i \in \Os_{B,b}[S]$$ 
est inf\'erieur \`a $(n-1)d+(d-1)=nd-1$.

La remarque finale est claire lorsque le point~$y$ est rationnel, en utilisant le fait que l'anneau~$\Os_{B,b}[S]$ se plonge dans~$\Os_{Y,y}$ et l'unicit\'e de la division. Le cas d'un point~$y$ quelconque se ram\`ene \`a celui d'un point rationnel par le m\^eme raisonnement que pr\'ec\'edemment.
\end{proof}

\section{Endomorphismes de la droite}\label{endodroite}

\index{Morphisme fini!endomorphisme d'une droite|(}

Dans cette partie, nous \'etudions les morphismes finis d'une partie de la droite analytique dans elle-m\^eme donn\'es par un polyn\^ome \`a coefficient dominant inversible. Maintenant que nous disposons du th\'eor\`eme de division de Weierstra{\ss} pour les points rigides, nous pouvons suivre pas \`a pas les raisonnements utilis\'es en g\'eom\'etrie analytique complexe.

\bigskip

Soit~$(\As,\|.\|)$ un anneau de Banach uniforme. Nous notons~$B=\Ms(\As)$, $X=\E{1}{\As}$ (avec variable~$T$) et~$\pi : X\to B$ le morphisme de projection. Fixons, d\`es \`a pr\'esent, les notations. Notons~$K$ l'anneau total des fractions de~$\As$. Soit
$$P(T) = \sum_{i=0}^d a_{i}\, T^i,$$
avec~$d\in\N^*$ et~$a_{0},\ldots,a_{d}\in K$, un polyn\^ome non constant \`a coefficients dans~$K$. Pour toute partie~$V$ compacte et spectralement convexe de l'espace~$B$ sur laquelle les coefficients de~$P$ sont d\'efinis et le coefficient~$a_{d}$ inversible, le morphisme naturel
$$\Bs(V)[T] \to \Bs(V)[T,S]/(P(S)-T) \xrightarrow[]{\sim} \Bs(V)[S]$$
induit un morphisme continu de la partie~$\pi^{-1}(V)$ dans elle-m\^eme. Soit~$U$ une partie localement connexe de l'espace~$B$ sur laquelle les coefficients de~$P$ sont d\'efinis et le coefficient~$a_{d}$ inversible. Tout point de~$U$ poss\`ede un syst\`eme fondamental de voisinages compacts et spectralement convexes. Par cons\'equent, nous pouvons construire un morphisme
$$\varphi : \pi^{-1}(U) \to \pi^{-1}(U)$$
en recollant des morphismes du type pr\'ec\'edent. Afin d'\'eviter les confusions, nous noterons respectivement~$Z$ et~$Y$ la source et le but du morphisme~$\varphi$. Nous consid\'ererons donc le morphisme
$$\varphi : Z \to Y.$$

\begin{prop}\label{phifini2}\index{Morphisme fini!au sens topologique}
Le morphisme~$\varphi$ est un morphisme topologique fini.
\end{prop}
\begin{proof}
Le fait que le morphisme~$\varphi$ soit continu est imm\'ediat. Pour montrer qu'il est ferm\'e, nous allons montrer qu'il est topologiquement propre, c'est-\`a-dire que l'image r\'eciproque de toute partie compacte est encore compacte.

Soit~$E$ une partie compacte de~$Y$. Il existe une partie compacte~$C$ de~$B$ et un nombre r\'eel~$r$ tels que la partie~$E$ soit contenue dans le disque compact $\overline{D}_{C}(r)$. La partie $\varphi^{-1}(E)$ est alors une partie ferm\'ee de 
$$\varphi^{-1}(\overline{D}_{C}(r)) = \left\{z\in Z\, \big|\, \pi(z)\in C,\, |P(S)(z)|\le r \right\}.$$
D'apr\`es le corollaire \ref{partiecompacte}, cette derni\`ere partie est compacte. On en d\'eduit que la partie $\varphi^{-1}(E)$ l'est \'egalement.


Pour finir, montrons que les fibres du morphisme~$\varphi$ sont finies. Soit~$y$ un point de~$Y$. L'ensemble de ses ant\'ec\'edents par l'application~$\varphi$ est l'ensemble des points de l'espace analytique~$\E{1}{\Hs(y)}$, dont nous noterons~$S$ la variable, qui annulent le polyn\^ome
$$Q_{y}(S) = P(y)(S)-T(y) =  \sum_{i=0}^d a_{i}(y)\, S^i - T(y) \in \Hs(y)[S].$$
Puisque le polyn\^ome~$P$ n'est pas constant et que son coefficient dominant ne s'annule pas sur~$Y$, le polyn\^ome~$Q_{y}(S)$ n'est pas nul. On en d\'eduit que l'ensemble~$\varphi^{-1}(y)$ est fini.

\end{proof}






Posons
$$G(S) = P(S)-T \in \Os(U)[T][S].$$
Consid\'erons~$\E{2}{\As}$ l'espace affine analytique de dimension~$2$ sur~$\As$ avec variables~$S$ et~$T$. Notons~$Z'$ l'ouvert de~$\E{2}{\As}$ form\'e des points dont la projection sur l'espace~$B$ appartient \`a~$U$. Le polyn\^ome~$G$ d\'efinit une fonction analytique sur l'espace~$Z'$. Nous identifierons l'espace analytique~$Z$ avec le ferm\'e de Zariski de l'espace~$Z'$ d\'efini par l'\'equation~$G=0$. Soit~$y$ un point de~$Y$. Notons~$z_{1},\ldots,z_{t}$, avec~$t\in\N^*$, ses ant\'ec\'edents par le morphisme~$\varphi$. Le th\'eor\`eme qui suit est l'analogue du th\'eor\`eme~2 de~\cite{GR}, I, \S 2. 

\begin{rem}
Les d\'efinitions~\ref{conditionIG} et~\ref{conditionS} des conditions~$(I_{G})$ et~$(S)$ \'etant locales, elles s'adaptent sans peine au cas des points d'un espace analytique qui n'est pas un spectre analytique. Nous nous autoriserons donc \`a les utiliser encore sans plus de pr\'ecautions.
\end{rem}


\begin{thm}
Supposons que le point de~$y$ de~$Y$ satisfait les conditions~$(I_{G})$ et~$(S)$. Soit~$(f_{1},\ldots,f_{t}) \in \prod_{i=1}^t \Os_{Z',z_{i}}$. Alors, il existe un unique \'el\'ement~$(r,q_{1},\ldots,q_{t})$ de~$\Os_{Y,y}[S]\times\prod_{i=1}^t \Os_{Z',z_{i}}$ v\'erifiant les propri\'et\'es suivantes :
\begin{enumerate}[\it i)]
\item le polyn\^ome~$r$ est de degr\'e strictement inf\'erieur \`a~$d$ ;
\item quel que soit~$i\in\cn{1}{t}$, nous avons~$f_{i} = q_{i}\, G + r$ dans~$\Os_{Z',z_{i}}$.
\end{enumerate}
\end{thm}
\begin{proof}
Dans~$\Hs(y)[S]$, \'ecrivons le polyn\^ome~$G$ sous la forme
$$G(y)[S] = a_{d}(y)\, \prod_{i=1}^t p_{i}(S)^{n_{i}},$$
o\`u~$p_{1},\ldots,p_{t}$ sont des polyn\^omes irr\'eductibles et unitaires \`a coefficients dans~$\Hs(y)$ et $n_{1},\ldots,n_{t}$ des \'el\'ements de~$\N^*$. Pour~$i\in\cn{1}{t}$, notons~$d_{i}$ le degr\'e du polyn\^ome~$p_{i}$. Les points~$z_{1},\ldots,z_{t}$ de~$Z'$ sont donc les \'el\'ements de la fibre au-dessus du point~$y$ d\'efinis par l'annulation des polyn\^omes~$p_{1},\ldots,p_{t}$. Quitte \`a changer l'ordre des polyn\^omes, nous pouvons supposer que, quel que soit~$i\in\cn{1}{t}$, le point~$z_{i}$ est d\'efini par l'\'equation~$p_{i}=0$.

D'apr\`es la condition~$(I_{G})$, la d\'ecomposition en produits de facteurs irr\'eductibles du polyn\^ome~$G[S]$ dans~$\kappa(y)[S]$ et dans~$\Hs(y)[S]$ est identique. On en d\'eduit que, quel que soit~$i\in\cn{1}{t}$, le polyn\^ome~$p_{i}$ est \`a coefficients dans~$\kappa(y)$. D'apr\`es la proposition~\ref{Hensel}, l'anneau local~$\Os_{Y,y}$ est hens\'elien. Par cons\'equent, il existe des polyn\^omes~$G_{1},\ldots,G_{t}$ unitaires \`a coefficients dans~$\Os_{Y,y}$ v\'erifiant les propri\'et\'es suivantes :
\begin{enumerate}
\item $\disp G = a_{d}\, \prod_{i=1}^t G_{i}$ dans~$\Os_{Y,y}[S]$ ;
\item quel que soit~$i\in\cn{1}{t}$, nous avons~$G_{i}(y) = p_{i}^{n_{i}}$.\\
\end{enumerate}

D\'emontrons, maintenant, l'existence de l'\'ecriture annonc\'ee. Il suffit de le d\'emontrer pour des $t$-uplets~$(f_{1},\ldots,f_{t})$ comportant un seul terme non nul. Le r\'esultat g\'en\'eral en d\'ecoulera par addition. Soit~$i\in\cn{1}{t}$ et supposons que, quel que soit~$j\ne i$, nous avons~$f_{j}=0$. Posons
$$e_{i} = a_{d}\, \prod_{j\ne i} G_{j}(S).$$
La fonction~$e_{i}$ est inversible au voisinage du point~$z_{i}$. Puisque le point de~$y$ de~$Y$ satisfait la condition~$(S)$, nous pouvons appliquer le th\'eor\`eme de division de {Weierstra\ss}~\ref{divisionrigide}. On en d\'eduit qu'il existe un \'el\'ement~$q_{i}$ de~$\Os_{Z',z_{i}}$ et un polyn\^ome~$r'$ \`a coefficients dans~$\Os_{Y,y}$ et de degr\'e strictement inf\'erieur \`a~$n_{i}d_{i}$ tels que
$$f_{i}\, e_{i}^{-1} = q_{i}\, G_{i} + r' \textrm{ dans } \Os_{Z',z_{i}}.$$
En multipliant l'\'egalit\'e par~$e_{i}$, nous obtenons
$$f_{i} = q_{i}\, G + r  \textrm{ dans } \Os_{Z',z_{i}},$$
o\`u~$r=r'\, e_{i}$ est un polyn\^ome de degr\'e strictement inf\'erieur \`a~$d$.

Soit~$j\ne i$. La fonction~$G_{i}$ est inversible au voisinage du point~$z_{j}$. Posons
$$q_{j} = -r'\, G_{i}^{-1}  \textrm{ dans } \Os_{Z',z_{j}}.$$
Nous avons alors
$$0 = q_{j}\, G + r  \textrm{ dans } \Os_{Z',z_{j}}.$$
\ \\

Pour finir, d\'emontrons l'unicit\'e de l'\'ecriture obtenue. Soit~$(r,q_{1},\ldots,q_{t})$ un \'el\'ement de~$\Os_{Y,y}[S]\times\prod_{i=1}^t \Os_{Z',z_{i}}$ v\'erifiant les propri\'et\'es suivantes :
\begin{enumerate}[\it i)]
\item le polyn\^ome~$r$ est de degr\'e strictement inf\'erieur \`a~$d$ ;
\item quel que soit~$i\in\cn{1}{t}$, nous avons~$f_{i} = q_{i}\, G + r$ dans~$\Os_{Z',z_{i}}$.
\end{enumerate}
Pour montrer que l'\'ecriture est unique, nous pouvons supposer que, quel que soit~$i\in\cn{1}{t}$, nous avons~$f_{i}=0$ et montrer alors que~$r=0$ et que, quel que soit~$i\in\cn{1}{t}$, $q_{i}=0$. Supposons donc que, quel que soit~$i\in\cn{1}{t}$, nous avons~$f_{i}=0$. Soit~$i\in\cn{1}{t}$. Avec les m\^emes notations que pr\'ec\'edemment, nous obtenons l'\'egalit\'e
$$-r = (q_{i}e_{i}) G_{i} \textrm{ dans } \Os_{Z',z_{i}}.$$
D'apr\`es la remarque finale du th\'eor\`eme~\ref{divisionrigide}, cette \'egalit\'e vaut dans~$\Os_{Y,y}[S]$. Puisque les polyn\^omes $G_{1},\ldots,G_{t}$ sont premiers entre eux deux \`a deux, leur produit~$a_{d}^{-1}G$ divise~$r$. Pour des raisons de degr\'e, cela impose au polyn\^ome~$r$ d'\^etre nul. Par unicit\'e de la division euclidienne dans chacun des anneaux~$\Os_{Z',z_{i}}$, avec~$i\in\cn{1}{t}$, nous en d\'eduisons que les fonctions~$q_{1},\ldots,q_{t}$ sont \'egalement nulles. 
\end{proof}

Nous parvenons enfin au r\'esultat attendu.

\index{Morphisme fini!image directe du faisceau structural|(}

\begin{thm}\label{thfini}
Supposons que le point de~$y$ de~$Y$ v\'erifie les conditions~$(I_{G})$ et~$(S)$. Alors, le morphisme
$$\begin{array}{ccc}
\Os_{Y,y}^d & \to & (\varphi_{*}\Os_{Z})_{y}\\
(c_{0},\ldots,c_{d-1}) & \mapsto & \disp \sum_{i=0}^{d-1} c_{i}\, S^i
\end{array}$$
est un isomorphisme. 
\end{thm}
\begin{proof}
Puisque le polyn\^ome~$G$ est de degr\'e~$d$ est que son coefficient dominant est inversible sur~$\pi(Y)$, le morphisme
$$\begin{array}{ccc}
\Os_{Y,y}^d & \to & \Os_{Y,y}[S]/(G(S))\\
(c_{0},\ldots,c_{d-1}) & \mapsto & \disp \sum_{i=0}^{d-1} c_{i}\, S^i
\end{array}$$
est un isomorphisme. Le r\'esultat d\'ecoule alors du th\'eor\`eme pr\'ec\'edent gr\^ace \`a l'\'egalit\'e
$$ (\varphi_{*}\Os_{Z})_{y} \simeq \prod_{i=1}^t \Os_{Z,z_{i}}.$$
\end{proof}

Nous d\'eduisons de ce r\'esultat deux corollaires. Leur d\'emonstration est similaire \`a celle des corollaires~\ref{casreduitglobal} et~\ref{casreduitStein}.

\begin{cor}\label{thfinicor1}
Supposons que tout point de~$Y$ v\'erifie les conditions~$(I_{G})$ et~$(S)$. Alors, le morphisme
$$\begin{array}{ccc}
\Os_{Y}^d & \to & \varphi_{*}\Os_{Z}\\
(c_{0},\ldots,c_{d-1}) & \mapsto & \disp \sum_{i=0}^{d-1} c_{i}\, S^i
\end{array}$$
est un isomorphisme. En particulier, pour toute partie~$V$ de~$Y$, le morphisme naturel
$$\Os_{Y}(V)[S]/(P(S)-T) \to \Os_{Z}(\varphi^{-1}(V))$$
est un isomorphisme.
\end{cor}

\begin{cor}\label{thfinicor2}
Supposons que tout point de~$Y$ v\'erifie les conditions~$(I_{G})$ et~$(S)$. Supposons que le faisceau~$\Os_{Y}$ est coh\'erent. Pour toute partie~$V$ de~$Y$, nous noterons
$$\varphi_{V} : \varphi^{-1}(V) \to V$$
le morphisme d\'eduit de~$\varphi$ par restriction et corestriction. Alors, pour toute partie~$V$ de~$Y$ et tout faisceau de $\Os_{\varphi^{-1}(V)}$-modules coh\'erent~$\Fs$, le faisceau de $\Os_{V}$-modules~$(\varphi_{V})_{*}\Fs$ est coh\'erent.
\end{cor}

\index{Morphisme fini!image directe du faisceau structural|)}

\index{Morphisme fini!endomorphisme d'une droite|)}

\section{Au-dessus d'un anneau d'entiers de corps de nombres}\label{adcdn}

\index{Morphisme fini!au-dessus de A@au-dessus de $A$|(}

Nous souhaitons disposer des r\'esultats \'etablis \`a la section pr\'ec\'edente lorsque la base est le spectre d'un anneau d'entiers de corps de nombres.  Nous nous pla\c{c}ons dans ce cadre et reprenons les notations du chapitre \ref{chapitredroite}. Nous souhaitons montrer que les hypoth\`eses du th\'eor\`eme~\ref{thfini} sont satisfaites. Commen\c{c}ons par nous int\'eresser \`a la condition~$(I_{G})$.
\newcounter{nodfini}\setcounter{nodfini}{\thepage}

\begin{prop}
Soient~$x$ un point de~$X$ et~$P(S)$ un polyn\^ome irr\'eductible de~$\kappa(x)[S]$. Alors l'image de~$P(S)$ dans~$\Hs(x)[S]$ est irr\'eductible.
\end{prop}
\begin{proof}
Supposons, tout d'abord, que la caract\'eristique du corps r\'esiduel~$\kappa(x)$ est un nombre premier. Le point~$x$ appartient alors \`a une fibre extr\^eme. D'apr\`es le th\'eor\`eme \ref{recrigext}, le corollaire \ref{avd3ext} ou la proposition \ref{avd2ext}, les corps $\kappa(x)$ et~$\Hs(x)$ sont naturellement isomorphes et le r\'esultat est tautologique.

Supposons, \`a pr\'esent, que la caract\'eristique du corps r\'esiduel~$\kappa(x)$ est nulle. Dans ce cas, le polyn\^ome~$P$ est s\'eparable. D'apr\`es le th\'eor\`eme \ref{resume}, {\it iv}, le corps~$\kappa(x)$ est hens\'elien. Nous concluons alors par la proposition~2.4.1 de~\cite{bleu}.
\end{proof}

\begin{cor}\label{XIG}\index{Condition $(I_{G})$!sur A1@sur $\AA$}
Soient~$x$ un point de~$X$ et~$G(S)$ un polyn\^ome \`a coefficients dans l'anneau local~$\Os_{X,x}$. Le point~$x$ v\'erifie la condition~$(I_{G})$.
\end{cor}

Int\'eressons-nous, \`a pr\'esent, \`a la condition~$(S)$.

\begin{lem}
Soient~$U$ une partie compacte et spectralement convexe de~$B$ et~$r$ un nombre r\'eel strictement positif. Supposons que les valeurs absolues associ\'ees aux points de~$U$ sont ultram\'etriques. Notons~$\Bs(U)\of{\{}{|T|\le r}{\}}$ l'alg\`ebre constitu\'ee des \'el\'ements de la forme
$$\sum_{i\ge 0} b_{i}\, T^i \textrm{ de } \Bs(U)[\![T]\!]$$
tels que la suite~$(\|b_{i}\|_{U}\, r^i)_{i\ge 0}$ tend vers~$0$. Munissons-la de la norme d\'efinie par
$$\left\|\sum_{i\ge 0} b_{i}\, T^i\right\|_{U,r,\textrm{um}} = \max_{i\ge 0} (\|b_{i}\|_{U}\, r^i).$$
C'est alors une alg\`ebre compl\`ete et le morphisme $A[T] \to \Bs(\overline{D}_{U}(r))$ induit un isomorphisme d'alg\`ebres norm\'ees
$$\Bs(U)\of{\{}{|T|\le r}{\}} \xrightarrow[]{\sim} \Bs(\overline{D}_{U}(r)).$$
\end{lem}
\begin{proof}
Puisque la partie~$U$ de~$B$ est spectralement convexe, le morphisme $A[T] \to \Bs(U)\of{\{}{|T|\le r}{\}}$ induit une injection continue
$$\Ms(\Bs(U)\of{\{}{|T|\le r}{\}}) \to \E{1}{A}$$
dont l'image est le disque~$\overline{D}_{U}(r)$. Soit~$F$ un \'el\'ement de~$A[T]$ qui ne s'annule pas sur le disque~$\overline{D}_{U}(r)$. Alors, d'apr\`es~\cite{rouge}, corollaire~1.2.4, l'\'el\'ement~$F$ est inversible dans~$\Bs(U)\of{\{}{|T|\le r}{\}}$. On en d\'eduit que le morphisme $A[T] \to \Bs(U)\of{\{}{|T|\le r}{\}}$ se prolonge en un morphisme injectif
$$\Ks(\overline{D}_{U}(r)) \hookrightarrow \Bs(U)\of{\{}{|T|\le r}{\}}.$$

Comparons, maintenant, la norme~$\|.\|_{U,r,\textrm{um}}$ \`a la norme uniforme sur le disque $\overline{D}_{U}(r)$. Soit~$F = \sum_{i\ge 0} f_{i}\, T^i \in A[T]$. Soit~$b$ un point de~$U$. La semi-norme associ\'ee au point~$b$ est, par hypoth\`ese, ultram\'etrique. Par cons\'equent, la norme uniforme~$\|.\|_{b,r}$ sur le disque ferm\'e de rayon~$r$ au-dessus du point~$b$ v\'erifie
$$\|F\|_{b,r} = \max_{i\ge 0} (|a_{i}(b)|\, r^i).$$
On en d\'eduit que
$$\|F\|_{\overline{D}_{U}(r)} = \max_{b\in U} (\|F\|_{b,r}) =  \max_{i\ge 0} (\|b_{i}\|_{U}\, r^i) = \|F\|_{U,r,\textrm{um}}.$$
Le morphisme pr\'ec\'edent se prolonge donc au compl\'et\'e de~$\Ks(\overline{D}_{U}(r))$. On en d\'eduit un morphisme injectif 
$$\Bs(\overline{D}_{U}(r)) \hookrightarrow \Bs(U)\of{\{}{|T|\le r}{\}}.$$
L'image de ce morphisme contient tous les polyn\^omes \`a coefficients dans~$\Bs(U)$. L'ensemble de ces polyn\^omes \'etant dense dans~$ \Bs(U)\of{\{}{|T|\le r}{\}}$, le morphisme pr\'ec\'edent est un isomorphisme.
\end{proof}

Rappelons que si~$(\As,\|.\|)$ d\'esigne un anneau de Banach et~$t$ un \'el\'ement de~$\R_{+}^*$, nous notons~$\As_{t}$ le compl\'et\'e de l'alg\`ebre~$\As\of{\la}{|T|\le t}{\ra}$ pour la norme uniforme sur son spectre analytique et que nous identifions ce spectre au disque ferm\'e~$\overline{D}(t)$ contenu dans la droite analytique~$\E{1}{\As}$. 


\begin{lem}\label{condR1}
Soient~$t>0$, $x$ un point de~$\overline{D}(t)$ et~$U$ un voisinage compact et connexe du point~$x$ dans~$\overline{D}(t)$ v\'erifiant les propri\'et\'es suivantes :
\begin{enumerate}[\it i)]
\item les valeurs absolues associ\'ees aux points de~$U$ sont ultram\'etriques ;
\item la partie~$U$ est spectralement convexe ;
\item la partie~$U$ poss\`ede un bord analytique fini et alg\'ebriquement trivial.
\end{enumerate}
Soit~$P(S)$ un polyn\^ome de~$\Os(U)[S]$ tel que l'\'el\'ement~R\'es$(P,P')$ de~$\Os(U)$ n'est pas nul. Alors, quel que soient~$s>0$ et~$r\in{]}{0,s}{]}$, le disque~$\overline{D}_{U}(r)$ (avec variable~$T_{0}$) de~$\Ms((\As_{t})_{s})$ v\'erifie la condition~$R_{P(S)-T_{0}}$.
\end{lem}
\begin{proof}
Soient~$s>0$ et~$r\in{]}{0,s}{]}$. Notons~$V=\overline{D}_{U}(r)$ muni de la variable~$T_{0}$. D'apr\`es la proposition \ref{stabilitespconvexe}, une telle partie de~$\Ms((\As_{t})_{s})$ est spectralement convexe.

Pour tout point~$y$ de~$\Gamma$, notons~$y_{r}$ le point~$\eta_{r}$ de la fibre de~$V$ au-dessus du point~$y$. Remarquons, d\`es \`a pr\'esent, que, pour tout point~$y$ de~$\Gamma$, l'\'el\'ement~$T_{0}(y_{r})$ de~$\Hs(y_{r})$ est transcendant sur~$\Hs(y)$. Notons
$$\Gamma_{V} = \{y_{r},\, y\in\Gamma\}.$$
Tout \'el\'ement de~$\Bs(V)$ atteint son maximum sur~$\Gamma_{V}$. Il suffit, pour s'en convaincre, d'utiliser la description explicite d\'emontr\'ee dans le lemme qui pr\'ec\`ede.

Pour montrer que le disque~$V$ v\'erifie la condition~$R_{P(S)-T_{0}}$, il nous suffit donc de montrer que la fonction~R\'es$(P(S)-T_{0},P'(S))$ ne s'annule pas sur~$\Gamma_{V}$. Soit~$y$ un point de~$\Gamma$. Nous avons~R\'es$(P(S),P'(S))(y)\ne 0$. En effet, dans le cas contraire, puisque l'anneau local~$\Os_{X,y}$ est un corps, la fonction~R\'es$(P(S),P'(S))$ serait nulle au voisinage du point~$y$ et donc nulle sur~$U$, d'apr\`es le principe du prolongement analytique. Consid\'erons l'\'el\'ement~$R_{y}(T_{0}) = \textrm{R\'es}(P(S)-T_{0},P'(S))$ de~$\Hs(y)[T]$. Nous venons de montrer que~$R_{y}(0)\ne 0$. On en d\'eduit que le polyn\^ome~$R_{y}$ n'est pas nul, puis que~$R_{y}(T_{0}(y_{r}))$ n'est pas nul, car~$T_{0}(y_{r})$ est transcendant sur~$\Hs(y)$. C'est le r\'esultat attendu. 
\end{proof}

\begin{lem}\label{condR2}
Soit~$x$ un point de~$X$. Soit~$P(S)$ un polyn\^ome de~$\Os_{X,x}[S]$ dont l'image dans~$\Hs(x)[S]$ est irr\'eductible. Alors l'\'el\'ement~R\'es$(P,P')$ de~$\Os_{X,x}$ n'est pas nul.
\end{lem}
\begin{proof}
Le polyn\^ome~$P(S)$ est \'egalement irr\'eductible dans~$\kappa(x)[S]$ et donc dans~$\Os_{X,x}[S]$, puisque l'anneau local~$\Os_{X,x}$ est hens\'elien. Or l'anneau~$\Os_{X,x}$ est int\`egre et son corps de fractions est parfait, car de caract\'eristique nulle. En effet, c'est une extension du corps des fractions de l'anneau~$\Os_{B,\pi(x)}$, dont la caract\'eristique est nulle. On en d\'eduit le r\'esultat voulu.

\end{proof}

\begin{prop}\label{XS}\index{Condition $(S)$!sur A1@sur $\AA$}
Tout point de~$X$ satisfait la condition~$(S)$. 
\end{prop}
\begin{proof}
Soit~$x$ un point de~$X$. Supposons, tout d'abord, que le corps r\'esiduel compl\'et\'e~$\Hs(x)$ est parfait. Soit~$P(S)$ un \'el\'ement de~$\Os_{X,x}[S]$ dont l'image dans~$\Hs(x)[S]$ est irr\'eductible. Puisque le corps r\'esiduel compl\'et\'e~$\Hs(x)$ est parfait, ce polyn\^ome est s\'eparable et nous avons~$\textrm{R\'es}(P(S),P'(S))(x)>0$. Par cons\'equent, il existe un nombre r\'eel~$m>0$ et un voisinage compact et spectralement convexe~$V$ du point~$x$ dans~$X$ sur lequel le polyn\^ome~$P(S)$ est d\'efini et la fonction~$\textrm{R\'es}(P(S),P'(S))$ minor\'ee par~$m$. Consid\'erons, \`a pr\'esent, le polyn\^ome~$\textrm{R\'es}(P(S)-T_{0},P'(S))$ de~$\Bs(V)[T_{0}]$. Nous venons de montrer que son coefficient constant est minor\'e sur~$V$. Par cons\'equent, il existe~$s>0$ tel que ce polyn\^ome ne s'annule pas sur~$\overline{D}_{V}(s)$. Pour tout voisinage compact et spectralement convexe~$U$ de~$x$ dans~$V$ et tout nombre r\'eel~$r\in\of{]}{0,s}{]}$, la condition~$(R_{P(S)-T_{0}})$ est alors v\'erifi\'ee sur le disque~$\overline{D}_{U}(r)$. 

Supposons, \`a pr\'esent, que le corps r\'esiduel compl\'et\'e~$\Hs(x)$ n'est pas parfait. Le point~$x$ appartient alors n\'ecessairement \`a une fibre extr\^eme de l'espace~$X$. D'apr\`es le th\'eor\`eme \ref{Shilovum}, il poss\`ede un syst\`eme fondamental de voisinages v\'erifiant les conditions du lemme~\ref{condR1}. On conclut alors \`a l'aide du lemme~\ref{condR2}.
\end{proof}



Nous pouvons, \`a pr\'esent, appliquer le th\'eor\`eme~\ref{thfini} et ses corollaires. Nous allons notamment en d\'eduire une expression explicite des anneaux de sections globales au voisinage des lemniscates. Soit~$V$ une partie de l'espace~$B$ et~$P(S)$ un polyn\^ome \`a coefficients dans~$\Os(D)$ dont le coefficient dominant est inversible. Soient~$s$ et~$t$ deux \'el\'ements de~$\R_{+}$ v\'erifiant l'in\'egalit\'e $s\le t$. Posons
$${\renewcommand{\arraystretch}{1.5}\begin{array}{cl}
& C_{0} = \left\{x\in X_{V}\, \big|\, s \le |T(x)| \le t\right\},\\
& L_{0} = \left\{x\in X_{V}\, \big|\, s \le |P(T)(x)| \le t\right\},\\
& C_{1} = \left\{x\in X_{V}\, \big|\, s < |T(x)| \le t\right\},\\
& L_{1} = \left\{x\in X_{V}\, \big|\, s < |P(T)(x)| \le t\right\},\\
& C_{2} = \left\{x\in X_{V}\, \big|\, s \le |T(x)| < t\right\},\\
& L_{2} = \left\{x\in X_{V}\, \big|\, s \le |P(T)(x)| < t\right\},\\
& C_{3} = \left\{x\in X_{V}\, \big|\, s < |T(x)| < t\right\},\\
& L_{3} = \left\{x\in X_{V}\, \big|\, s < |P(T)(x)| < t\right\},\\
& C_{4} = \left\{x\in X_{V}\, \big|\,  |T(x)| \ge s\right\},\\
& L_{4} = \left\{x\in X_{V}\, \big|\, |P(T)(x)| \ge s\right\},\\
& C_{5} = \left\{x\in X_{V}\, \big|\, |T(x)|> s\right\}\\
\textrm{et} & L_{5} = \left\{x\in X_{V}\, \big|\, |P(T)(x)|> s\right\}.
\end{array}}$$
Choisissons un couple~$(C,L)$ parmi les couples~$(C_{i},L_{i})$, avec $i\in\cn{0}{5}$.

%

%

Proc\'edons, \`a pr\'esent, comme dans la section~\ref{endodroite}. Pour toute partie~$V$ compacte et spectralement convexe contenue dans~$D$, le morphisme naturel
$$\Bs(V)[T] \to \Bs(V)[T,S]/(P(S)-T) \xrightarrow[]{\sim} \Bs(V)[S]$$
induit un morphisme continu de la partie~$X_{V}$ dans elle-m\^eme. Ces morphismes se recollent en un morphisme
$$\varphi : X_{V} \to X_{V}.$$
Remarquons que nous avons l'\'egalit\'e
$$L = \varphi^{-1}(C).$$

D'apr\`es le corollaire~\ref{XIG} et la proposition~\ref{XS}, les hypoth\`eses des corollaires~\ref{thfinicor1} et~\ref{thfinicor2} sont v\'erifi\'ees. En appliquant le corollaire~\ref{thfinicor1}, nous obtenons le r\'esultat suivant.

\index{Faisceau structural!sections sur une lemniscate de A1@section sur une lemniscate de $\AA$}

\begin{thm}\label{sectionslemniscates}
Le morphisme naturel
$$\Os(C)[S]/(P(S)-T) \to \Os(L)$$
est un isomorphisme.
\end{thm}

Le th\'eor\`eme \ref{isocouronne} nous permet d'en d\'eduire une description explicite des anneaux de sections globales des lemniscates.



\index{Morphisme fini!au-dessus de A@au-dessus de $A$|)}

%% file: Stein.tex
\chapter{Espaces de Stein}\label{chapitreStein}

Ce chapitre est consacr\'e \`a l'\'etude de quelques sous-espaces de Stein de la droite analytique au-dessus d'un anneau d'entiers de corps de nombres. Nous y utiliserons les notations du chapitre \ref{chapitredroite}. Pr\'ecis\'ement, nous d\'emontrons que certaines parties assez simples, disques, couronnes ou lemniscates relatifs, sont des espaces de Stein.

Le num\'ero \ref{defiStein} contient les d\'efinitions dans un cadre g\'en\'eral : nous appellerons espace de Stein tout espace annel\'e qui satisfait les conclusions des th\'eor\`emes~A et~B de H.~Cartan. Nous rappelons \'egalement quelques propri\'et\'es classiques de ces espaces.

Au num\'ero \ref{cgplc}, nous nous sommes attach\'e \`a d\'egager des conditions sous lesquelles une r\'eunion de deux parties compactes et de Stein est encore un espace de Stein. Les notions introduites peuvent sembler absconses, mais elle ne sont que formalisations des m\'ethodes de la g\'eom\'etrie analytique complexe. 

Nous reprenons ensuite le cadre du chapitre \ref{chapitredroite}, celui de la droite affine analytique au-dessus d'un anneau d'entiers de corps de nombres. Nous utilisons alors les r\'esultats obtenus pour montrer, par r\'ecurrence, que les parties compactes et connexes de l'espace de base sont des espaces de Stein, au num\'ero \ref{pcdlb}, ainsi que les couronnes compactes des fibres, au num\'ero \ref{pcdf}, et les couronnes compactes et connexes de la droite, au num\'ero \ref{ccdld}.

Au num\'ero \ref{ldld}, finalement, nous traiterons le cas des couronnes ouvertes et, plus g\'en\'eralement, de toute partie de la forme
$$\left\{x\in X_{V}\, \big|\, s < |P(T)(x)| < t\right\},$$
o\`u~$V$ est une partie connexe de la base, $P(T)$ un polyn\^ome unitaire \`a coefficients dans~$\Os(V)$ et~$r$ et~$s$ deux nombres r\'eels. Nous reprenons, l\`a encore, les m\'ethodes de la g\'eom\'etrie analytique complexe. Nous indiquons tout d'abord des conditions sous lesquelles une partie qui poss\`ede une exhaustion par des parties compactes et de Stein est elle-m\^eme un espace de Stein. Nous d\'emontrons ensuite un r\'esultat de fermeture pour certains germes de faisceaux, qui nous semble pr\'esenter un int\'er\^et ind\'ependant. Nous concluons finalement en d\'ecrivant explicitement des exhaustions de Stein pour les couronnes ouvertes et en utilisant les r\'esultats sur les morphismes finis d\'emontr\'es au chapitre \ref{chapitrefini}.




\section{D\'efinitions}\label{defiStein}

Soit $(X,\Os_{X})$ un espace localement annel\'e. Avant d'en venir aux d\'e\-mon\-stra\-tions annonc\'ees, nous rappelons quelques propri\'et\'es et d\'efinitions dans un cadre g\'en\'eral. Expliquons, tout d'abord, ce que nous entendons par espace de Stein. Nous utiliserons la d\'efinition cohomologique classique.

\begin{defi}\index{Espace de Stein!theoreme A@th\'eor\`eme A|see{Th\'eor\`eme A}}\index{Theoreme A@Th\'eor\`eme A}\index{Theoreme@Th\'eor\`eme!A|see{Th\'eor\`eme A}}
Soit~$\Fs$ un faisceau de $\Os_{X}$-modules. Nous dirons que le faisceau~$\Fs$ satisfait le {\bf th\'eor\`eme A} si, pour tout point~$x$ de~$X$, le $\Os_{X,x}$-module~$\Fs_{x}$ est engendr\'e par l'ensemble de ses sections globales~$\Fs(X)$.

Soit~$Y$ une partie de~$X$. Nous dirons que le faisceau~$\Fs$ satisfait le th\'eor\`eme~A sur~$Y$ si le faisceau de $\Os_{Y}$-modules~$\Fs_{|Y}$ satisfait le th\'eor\`eme~A.
\end{defi}

Remarquons que, par d\'efinition (\emph{cf.} d\'efinition \ref{defitf}), le th\'eor\`eme~A est satisfait localement pour les faisceaux de type fini. \'Enon\c{c}ons ce r\'esultat sous forme d'un lemme afin de pouvoir nous y faire r\'ef\'erer ult\'erieurement.

\begin{lem}\label{thAlocal}\index{Theoreme A@Th\'eor\`eme A!cas local}
Soit~$\Fs$ un faisceau de $\Os_{X}$-modules de type fini. Soit~$x$ un point de~$X$. Il existe un voisinage~$V$ du point~$x$ dans~$X$ tel que le faisceau~$\Fs$ satisfasse le th\'eor\`eme~A sur~$V$. 
\end{lem}





Signalons \'egalement que lorsque nous consid\'erons des parties compactes, nous pouvons pr\'eciser le r\'esultat du th\'eor\`eme~A.

\begin{lem}\label{thAcompact}\index{Theoreme A@Th\'eor\`eme A!cas d'un compact}
Soit~$\Fs$ un faisceau de $\Os_{X}$-modules de type fini qui satisfait le th\'eor\`eme~A. Si l'espace~$X$ est compact, il existe un ensemble fini d'\'el\'ements de~$\Fs(X)$ dont les images engendrent le $\Os_{X,x}$-module~$\Fs_{x}$ en tout point~$x$ de~$X$.
\end{lem}
\begin{proof}
Soit~$x$ un point de l'espace~$X$. Puisque le faisceau~$\Fs$ est un $\Os_{X}$-module de type fini, il existe un voisinage~$U$ du point~$x$ dans~$X$, un entier~$p$ et des \'el\'ements $F_{1},\ldots,F_{p}$ de~$\Fs(X)$ tels que, pour tout point~$y$ de~$U$, le $\Os_{X,y}$-module~$\Fs_{y}$ soit engendr\'e par les germes $(F_{1})_{y},\ldots,(F_{p})_{y}$. 

D'apr\`es le th\'eor\`eme~A, il existe un entier~$q$ et des \'el\'ements $G_{1},\ldots,G_{q}$ de~$\Fs(X)$ tels le $\Os_{X,x}$-module~$\Fs_{x}$ soit engendr\'e par les germes $(G_{1})_{x},\ldots,(G_{q})_{x}$. En particulier, il existe une famille $(a_{i,j})_{1\le i\le p, 1\le j\le q}$ d'\'el\'ements de~$\Os_{X,x}$ tels que l'on ait 
$$\forall i\in \cn{1}{p},\, (F_{i})_{x} = \sum_{j=1}^q a_{i,j}\, (G_{j})_{x} \textrm{ dans } \Fs_{x}.$$ 
Il existe un voisinage~$V$ du point~$x$ dans~$U$ sur lequel les \'el\'ements~$a_{i,j}$, avec $i\in\cn{1}{p}$ et $j\in\cn{1}{q}$, sont d\'efinis et les \'egalit\'es pr\'ec\'edentes sont valables. On en d\'eduit que, pour tout point~$y$ de~$V$, le $\Os_{X,y}$-module~$\Fs_{y}$ est engendr\'e par les germes $(G_{1})_{y},\ldots,(G_{q})_{y}$. 

On conclut finalement en utilisant la compacit\'e de l'espace~$X$.
\end{proof}

\begin{defi}\index{Espace de Stein!theoreme B@th\'eor\`eme B|see{Th\'eor\`eme B}}\index{Theoreme B@Th\'eor\`eme B}\index{Theoreme@Th\'eor\`eme!B|see{Th\'eor\`eme B}}

Soit $\Fs$ un faisceau de $\Os_{X}$-modules. Nous dirons que le faisceau~$\Fs$ satisfait le {\bf th\'eor\`eme B} si, quel que soit $q\in \N^*$, nous avons
$$H^q(X,\Fs) = 0.$$

Soit~$Y$ une partie de~$X$. Nous dirons que le faisceau~$\Fs$ satisfait le th\'eor\`eme~B sur~$Y$ si le faisceau de $\Os_{Y}$-modules~$\Fs_{|Y}$ satisfait le th\'eor\`eme~B.
\end{defi}


\begin{defi}\index{Espace de Stein}\index{Stein|see{Espace de Stein}}
Nous dirons que l'espace $X$ est {\bf un espace de Stein} si tout faisceau de $\Os_{X}$-modules coh\'erent satisfait les th\'eor\`emes $A$ et $B$.
\end{defi}





\begin{rem}
Attention, cette d\'efinition d'espace de Stein est plus faible que la d\'efinition classique pour les espaces analytiques sur un corps ultram\'etrique (\emph{cf.}~\cite{AB}, d\'efinition~2.3). 
\end{rem}

\bigskip

Soit~$Y$ une partie de~$X$. Lorsque~$Y$ est compacte, les propri\'et\'es de finitude des faisceaux coh\'erents imposent des liens entre les faisceaux coh\'erents sur~$Y$ et les faisceaux coh\'erents d\'efinis sur un voisinage de~$Y$ dans~$X$. Le r\'esultat qui suit est d\'emontr\'e \`a la proposition~$1$ de~\cite{SemCartan419}. La preuve qui y figure est \'ecrite dans le langage de la g\'eom\'etrie analytique complexe, mais elle s'adapte \`a notre cadre, sans la moindre modification.

\begin{prop}\index{Faisceau!coherent@coh\'erent!au voisinage d'un compact}
Supposons que la partie~$Y$ est compacte. Soit~$\Fs$ un faisceau de $\Os_{Y}$-modules coh\'erent. Alors, il existe un voisinage ouvert~$U$ du compact~$Y$ dans~$X$ et un faisceau de $\Os_{X}$-modules coh\'erent~$\Gs$ tel que l'on ait
$$\Fs = \Gs_{|Y}.$$ 
\end{prop}

\begin{cor}\label{Steincompact}\index{Espace de Stein!compact}
Supposons que la partie~$Y$ est compacte et poss\`ede un syst\`eme fondamental de voisinages form\'e d'espaces de Stein. Alors, la partie~$Y$ est de Stein.
\end{cor}




\bigskip

Mentionnons que l'on peut remplacer les conditions qui figurent dans la d\'efinition d'espace de Stein par des conditions plus faibles. En effet, le th\'eor\`eme~A se d\'eduit du th\'eor\`eme~B. La nullit\'e du premier groupe de cohomologie \`a coefficient dans n'importe quel faisceau coh\'erent suffit d'ailleurs \`a assurer le r\'esultat (\emph{cf.}~\cite{GR}, IV, \S 1, th\'eor\`eme~2).

\begin{thm}\label{BimpliqueA}\index{Theoreme A@Th\'eor\`eme A!aaavient de theoreme B@$\Leftarrow$ th\'eor\`eme B}\index{Theoreme B@Th\'eor\`eme B!aaaimplique theoreme A@$\Rightarrow$ th\'eor\`eme A}
Supposons que, pour tout faisceau de $\Os_{X}$-modules coh\'erent~$\Fs$, on ait $H^1(X,\Fs)=0$. Alors tout faisceau de $\Os_{X}$-modules coh\'erent satisfait le th\'eor\`eme~A.
\end{thm}

On d\'eduit de ce r\'esultat une stabilit\'e de la notion d'espace de Stein par morphisme fini.

\begin{thm}\label{finiStein}\index{Espace de Stein!stabilite par morphisme fini@stabilit\'e par morphisme fini}\index{Morphisme fini!applications aux espaces de Stein}
Soit $\varphi : Y \to X$ un morphisme topologique fini. Supposons que, pour tout faisceau de $\Os_{Y}$-modules coh\'erent~$\Fs$, le faisceau de $\Os_{X}$-modules~$\varphi_{*}\Fs$ est coh\'erent. Alors l'espace~$Y$ est un espace de Stein.
\end{thm}
\begin{proof}
Soit~$\Fs$ un faisceau de $\Os_{Y}$-modules coh\'erent. D'apr\`es le th\'eor\`eme~\ref{finicohomologie}, pour tout entier~$q$, nous avons un isomorphisme
$$H^q(Y,\Fs) \simeq H^q(X,\varphi_{*}\Fs).$$
Or le faisceau~$\varphi_{*}\Fs$ est coh\'erent et la partie~$X$ est de Stein. On en d\'eduit que, pour tout entier~$q\ge 1$, nous avons
$$H^q(Y,\Fs)=0.$$
Nous venons de montrer que tout faisceau de $\Os_{Y}$-modules coh\'erent satisfait le th\'eor\`eme~$B$. Le th\'eor\`eme~\ref{BimpliqueA} assure alors que la partie~$Y$ est de Stein.
\end{proof}

\section{Cadre g\'en\'eral pour les compacts}\label{cgplc}

Dans cette premi\`ere partie, nous nous sommes attacher \`a d\'egager un cadre g\'en\'eral pour d\'emontrer que des parties compactes sont des espaces de Stein. Nous consid\'ererons donc ici un espace localement annel\'e $(X,\Os_{X})$ et deux parties compactes~$K^-$ et~$K^+$ de l'espace topologique sous-jacent. Posons $L=K^-\cap K^+$ et $M=K^-\cup K^+$.
\newcounter{KLM}\setcounter{KLM}{\thepage}

\subsection{Lemmes de Cousin et de Cartan}

Il n'est gu\`ere ais\'e de travailler directement avec les anneaux de fonctions au voisinages de compacts. Nous allons donc introduire une d\'efinition qui nous permettra de consid\'erer plut\^ot des anneaux de Banach.

\begin{defi}\index{Systeme@Syst\`eme!de Banach}
Un {\bf syst\`eme de Banach associ\'e au couple $(K^-,K^+)$} est la donn\'ee de
\begin{enumerate}[\it i)]
\item un ensemble ordonn\'e filtrant $(A,\le)$ ;
\item un syst\`eme inductif $((\Bs^-_{\alpha},\|.\|^-_{\alpha}),\varphi^-_{\alpha,\beta})$ sur~$A$ \`a valeurs dans la cat\'egorie des anneaux de Banach et des morphismes born\'es entre iceux ;
\item un syst\`eme inductif $((\Bs^+_{\alpha},\|.\|^+_{\alpha}),\varphi^+_{\alpha,\beta})$ sur~$A$ \`a valeurs dans la cat\'egorie des anneaux de Banach et des morphismes born\'es entre iceux ;
\item un syst\`eme inductif $((\Cs_{\alpha},\|.\|_{\alpha}),\varphi_{\alpha,\beta})$ sur~$A$ \`a valeurs dans la cat\'egorie des anneaux de Banach et des morphismes born\'es entre iceux ;
\item pour tout \'el\'ement~$\alpha$ de~$A$, un morphisme born\'e $\psi_{\alpha}^- : \Bs^-_{\alpha} \to \Cs_{\alpha}$ ;
\item pour tout \'el\'ement~$\alpha$ de~$A$, un morphisme born\'e $\psi_{\alpha}^+ : \Bs^+_{\alpha} \to \Cs_{\alpha}$ ;
\item pour tout \'el\'ement~$\alpha$ de~$A$, un morphisme $\rho_{\alpha}^- : \Bs^-_{\alpha} \to \Os(K^-)$ ;
\item pour tout \'el\'ement~$\alpha$ de~$A$, un morphisme $\rho_{\alpha}^+ : \Bs^-_{\alpha} \to \Os(K^+)$ ;
\item pour tout \'el\'ement~$\alpha$ de~$A$, un morphisme $\rho_{\alpha} : \Cs_{\alpha} \to \Os(L)$
\end{enumerate}
v\'erifiant les propri\'et\'es suivantes :
\begin{enumerate}
\item pour tout \'el\'ement~$\alpha$ de~$A$, le diagramme 
$$\xymatrix{
\Bs^-_{\alpha} \ar[r]^{\rho_{\alpha}^-} \ar[d]^{\psi_{\alpha}^-} & \Os(K^-) \ar[d]^{\cdot_{L}}\\
\Cs_{\alpha} \ar[r]^{\rho_{\alpha}}& \Os(L)
}$$
commute ;
\item pour tout \'el\'ement~$\alpha$ de~$A$, le diagramme 
$$\xymatrix{
\Bs^+_{\alpha} \ar[r]^{\rho_{\alpha}^+} \ar[d]^{\psi_{\alpha}^+} & \Os(K^+) \ar[d]^{\cdot_{L}}\\
\Cs_{\alpha} \ar[r]^{\rho_{\alpha}}& \Os(L)
}$$
commute ;
\item pour tous \'el\'ements~$\alpha$ et~$\beta$ de~$A$ tels que $\alpha\le \beta$, le diagramme
$$\xymatrix{
\Cs_{\alpha} \ar[rd]^{\rho_{\alpha}} \ar[dd]_{\varphi_{\beta,\alpha}}&\\
& \Os(L)\\
\Cs_{\beta} \ar[ru]_{\rho_{\beta}}&
}$$ 
commute ;
\item le morphisme
$$\rho : \varinjlim_{\alpha\in A} \Cs_{\alpha} \to \Os(L)$$
induit par la famille de morphismes $(\rho_{\alpha})_{\alpha\in A}$ est un isomorphisme.
\end{enumerate}
\end{defi}

Dans toute la suite de ce paragraphe, nous consid\'ererons un syst\`eme de Banach~$\Omega$ associ\'e au couple $(K^-,K^+)$. Nous aurons besoin d'une propri\'et\'e suppl\'ementaire, connue sous le nom de lemme de Cousin.
\newcounter{Omega}\setcounter{Omega}{\thepage}

\begin{defi}\index{Systeme@Syst\`eme!de Cousin}\index{Lemme!de Cousin|see{Syst\`eme de Cousin}}
Un {\bf syst\`eme de Cousin associ\'e au couple $(K^-,K^+)$} est un syst\`eme de Banach associ\'e au m\^eme couple et pour lequel il existe un nombre r\'eel~$D$ v\'erifiant la propri\'et\'e suivante : pour tout \'el\'ement~$\alpha$ de~$A$ et tout \'el\'ement~$f$ de~$\Cs_{\alpha}$, il existe des \'el\'ements~$f^-$ de~$\Bs_{\alpha}^-$ et~$f^+$ de~$\Bs^+_{\alpha}$ tels que
\begin{enumerate}[\it i)]
\item $f=\psi_{\alpha}^-(f^-)+\psi_{\alpha}^+(f^+)$ ;
\item $\|f^-\|_{\alpha}^- \le D\,\|f\|_{\alpha}$ ;
\item $\|f^+\|_{\alpha}^+ \le D\,\|f\|_{\alpha}$.
\end{enumerate}
\end{defi}

Nous allons montrer que tout syst\`eme de Cousin v\'erifie le lemme de Cartan, en reprenant essentiellement la m\'ethode mise en {\oe}uvre dans~\cite{GR}, III \S 1. 

Commen\c{c}ons par introduire quelques notations. Soient $p,q \in\N^*$ et $(\Ds,\|.\|)$ un anneau de Banach. Nous d\'efinissons la norme d'une matrice $a=(a_{i,j})\in M_{p,q}(\Ds)$ par 
$$\|a\| = \max_{1\le i\le p} \left( \sum_{j=1}^q \|a_{i,j}\|\right).$$ 
La multiplication des matrices est continue par rapport \`a cette norme. En effet, on v\'erifie facilement que, quel que soient $r\in\N^*$, $a\in M_{p,q}(\Ds)$ et $b\in M_{q,r}(\Ds)$, on a 
$$\|ab\|\le \|a\|\, \|b\|.$$
Nous noterons $I\in M_{q}(\Ds)$ la matrice identit\'e. Nous allons, tout d'abord, d\'emontrer quelques lemmes.

\begin{lem}\label{inverse}
Toute matrice~$a$ de~$M_{q}(\Ds)$ v\'erifiant 
$$\|a- I\| \le \frac{1}{2}$$ 
est inversible et son inverse~$a^{-1}$ v\'erifie l'in\'egalit\'e 
$$\|a^{-1}\|\le 2.$$
\end{lem}

\begin{lem}
Soit $(a_{k})_{k\ge 0}$ une suite de $M_{q}(\Ds)$ v\'erifiant la condition 
$$\sum_{k\ge 0} \|a_{k}-I\| \le \frac{1}{2}.$$ 
Alors, quel que soit $n\in\N$, nous avons
$$\|a_{0}\cdots a_{n}-I\| \le 2 \sum_{k=0}^n \|a_{k}-I\|$$
\end{lem}
\begin{proof}
D\'emontrons ce r\'esultat par r\'ecurrence sur l'entier $n\in\N$. Pour $n=0$, c'est \'evident.

Supposons que la formule est vraie pour $n\in\N$. Nous avons
$$\begin{array}{rcl}
a_{0}\cdots a_{n}a_{n+1}-I &= &(a_{0}\cdots a_{n}-I)(a_{n+1}-I)\\ 
&&+ (a_{0}\cdots a_{n}-I)+(a_{n+1}-I).
\end{array}$$
On en d\'eduit que 
$$\begin{array}{rcl}
\disp\|a_{0}\cdots a_{n}a_{n+1}-I\| &\le&\disp \|a_{0}\cdots a_{n}-I\|\, \|a_{n+1}-I\|\\ 
&&\disp+ \|a_{0}\cdots a_{n}-I\| + \|a_{n+1}-I\|\\
&\le &\disp 2 \sum_{k=0}^n \|a_{k}-I\|\\
&&\disp + ( \|a_{0}\cdots a_{n}-I\| + 1)\,  \|a_{n+1}-I\|.
\end{array}$$
En outre, nous avons 
$$\|a_{0}\cdots a_{n}-I\| \le 2\, \sum_{k=0}^n \|a_{k}-I\| \le 1.$$
On en d\'eduit que 
$$\|a_{0}\cdots a_{n+1}-I\| \le 2 \sum_{k=0}^{n+1} \|a_{k}-I\|.$$
\end{proof}

\begin{lem}\label{convergence}
Soit $(g_{k})_{k\ge 0}$ une suite de $M_{q}(\Ds)$ v\'erifiant 
$$\sum\limits_{k\ge 0} \|g_{k}\| \le \frac{1}{4}.$$
Alors la suite de terme g\'en\'eral  
$$P_{n}=(I+g_{0})\cdots(I+g_{n})$$
converge dans $M_{q}(\Ds)$ vers une matrice inversible $P$ v\'erifiant
$$\|P-I\| \le 2 \sum\limits_{k\ge 0} \|g_{k}\|.$$ 
\end{lem}
\begin{proof}
D'apr\`es le lemme pr\'ec\'edent, quels que soient $j\ge i\ge 0$, nous avons
$$\|(I+g_{i})\cdots(I+g_{j})-I\| \le 2 \sum\limits_{k=i}^j \|g_{k}\|.$$
En particulier, quel que soit $n\in\N$, nous avons 
$$\|P_{n}-I\|\le 2\sum_{k=0}^n \|g_{k}\| \le \frac{1}{2}$$ 
et donc $\|P_{n}\| < 3/2$. On en d\'eduit que, quels que soient $m\ge n\ge 0$, on a
$${\renewcommand{\arraystretch}{1.5}\begin{array}{rcl}
\|P_{m}-P_{n}\| &=& \|P_{n}\, ((I+g_{n+1})\cdots(I+g_{m})-I)\|\\ 
&\le& \disp  \frac{3}{2}\,\sum\limits_{k=n+1}^m \|g_{k}\|.
\end{array}}$$
Par cons\'equent, la suite $(P_{n})_{n\ge 0}$ est une suite de Cauchy de $M_{q}(\Ds)$. Puisque cet anneau est complet, elle converge donc vers un \'el\'ement $P$. Nous avons n\'ecessairement 
$$\|P-I\|\le 2\sum\limits_{k\ge 0} \|g_{k}\| \le \frac{1}{2}.$$ 
Le lemme \ref{inverse} nous permet alors de conclure.
\end{proof}

Pla\c{c}ons-nous de nouveau dans le cadre des syst\`emes de Cousin. Pour~$\alpha\in A$, les morphismes~$\psi_{\alpha}^-$, $\psi_{\alpha}^+$, $\rho_{\alpha}^-$, $\rho_{\alpha}^+$ et~$\rho_{\alpha}$ se prolongent naturellement en des morphismes de groupes entre les espaces de matrices ; nous les noterons identiquement.

\begin{lem}
Supposons que~$\Omega$ soit un syst\`eme de Cousin. Soit $\eps\in\R$ v\'erifiant $0<\eps<1/2\, D^{-1}$ et $\beta=4D^2\eps<1$. Soit $\alpha\in A$. Soit $a=I+b \in M_{q}(\Cs_{\alpha})$ telle que $\|b\|_{\alpha} \le\eps$. Alors, il existe 
$a^-=I+b^-\in M_{q}(\Bs_{\alpha}^-)$, $a^+=I+b^+\in M_{q}(\Bs_{\alpha}^+)$ et $\tilde{a}=I+\tilde{b}\in M_{q}(\Cs_{\alpha})$ v\'erifiant les propri\'et\'es suivantes :
\begin{enumerate}[\it i)]
\item $a = \psi_{\alpha}^-(a^-)\, \tilde{a}\, \psi_{\alpha}^+(a^+)$ ;
\item $\|b^-\|_{\alpha}^- \le D\, \|b\|_{\alpha}$ ;
\item $\|b^+\|_{\alpha}^+ \le D\, \|b\|_{\alpha}$ ;
\item $\|\tilde{b}\|_{\alpha} \le \beta_{\alpha}\, \|b\|_{\alpha}$.
\end{enumerate}
\end{lem}
\begin{proof}
En appliquant la propri\'et\'e des syst\`emes de Cousin \`a chaque coefficient de la matrice~$b$, on montre qu'il existe des matrices $b^-\in M_{q}(\Bs_{\alpha}^-)$ et $b^+\in M_{q}(\Bs_{\alpha}^+)$ v\'erifiant les propri\'et\'es suivantes :
\begin{enumerate}
\item $b=\psi_{\alpha}^-(b^-)+\psi_{\alpha}^+(b^+)$ ;
\item $\|b^-\|_{\alpha}^-\le D\,\|b\|_{\alpha}$ ;
\item $\|b^+\|_{\alpha}^+\le D\,\|b\|_{\alpha}$.
\end{enumerate} 
Posons $a^-=I+b^-\in M_{q}(\Bs_{\alpha}^-)$ et $a^+=I+b^+\in M_{q}(\Bs_{\alpha}^+)$. Nous avons alors l'\'egalit\'e
$$\psi_{\alpha}^-(a^-)\, \psi_{\alpha}^+(a^+) = a + \psi_{\alpha}^-(b^-)\, \psi_{\alpha}^+(b^+).$$
Par choix de~$b$, nous avons $D\|b\|_{\alpha}\le 1/2$. D'apr\`es le lemme \ref{inverse}, la matrice~$a^-$ est inversible dans $M_{q}(\Bs^-_{\alpha})$ et v\'erifie $\|(a^-)^{-1}\|_{\alpha}^-\le 2$. De m\^eme, la matrice $a^+$ est inversible dans $M_{q}(\Bs_{\alpha}^+)$ et v\'erifie $\|(a^+)^{-1}\|_{\alpha}^+\le 2$.

Posons 
$$\tilde{a}=\psi_{\alpha}^-(a^-)^{-1}\, a\, \psi_{\alpha}^+(a^+)^{-1} \textrm{ et } \tilde{b} = \tilde{a}- I \textrm{ dans } M_{q}(\Cs_{\alpha}).$$
Nous avons alors
$${\renewcommand{\arraystretch}{1.3}\begin{array}{rcl}
\tilde{b} &=&  \psi_{\alpha}^-\left((a^-)^{-1}\right)\, a\, \psi_{\alpha}^+\left((a^+)^{-1}\right) - I\\
&=& \psi_{\alpha}^-\left((a^-)^{-1}\right)\, \left(\psi_{\alpha}^-(a^-)\,\psi_{\alpha}^+(a^+)-\psi_{\alpha}^-(b^-)\,\psi_{\alpha}^+(b^+)\right)\, \psi_{\alpha}^+\left((a^+)^{-1}\right) - I\\
&=&  - \psi_{\alpha}^-\left((a^-)^{-1}\right)\ \psi_{\alpha}^-(b^-)\, \psi_{\alpha}^+(b^+)\, \psi_{\alpha}^+\left((a^+)^{-1}\right)
\end{array}}$$
et nous en tirons l'in\'egalit\'e
$$\|\tilde{b}\|_{\alpha} \le 4D^2\|b\|_{\alpha}^2\le \beta\|b\|_{\alpha}.$$ 
\end{proof}

Nous voici enfin pr\^ets \`a d\'emontrer le lemme de Cartan.

\begin{thm}[Lemme de Cartan]\label{Cartan}\index{Lemme!de Cartan}
Supposons que~$\Omega$ soit un syst\`eme de Cousin. Alors, il existe $\eps\in\R$ v\'erifiant la propri\'et\'e suivante : quels que soient $\alpha\in A$ et $a\in M_{q}(\Cs_{\alpha})$ v\'erifiant $\|a-I\|_{\alpha} < \eps$, il existe $c^-\in GL_{q}(\Bs_{\alpha}^-)$ et $c^+\in GL_{q}(\Bs_{\alpha}^+)$ telles que
\begin{enumerate}[\it i)]
\item $a = \psi_{\alpha}^-(c^-)\, \psi_{\alpha}^+(c^+)$ ;
\item $\|c^--I\|_{\alpha}^- \le 4D\,\|a-I\|_{\alpha}$ ;
\item $\|c^+-I\|_{\alpha}^+  \le 4D\,\|a-I\|_{\alpha}$.
\end{enumerate}
\end{thm}
\begin{proof}
Choisissons $\eps\in\R$ v\'erifiant les conditions du lemme pr\'ec\'edent ainsi que $\beta \le 1/2$ et $\eps \le 1/(8D)$. Soient $\alpha\in A$ et $a\in M_{q}(\Cs_{\alpha})$ v\'erifiant $\|a-I\|_{\alpha} < \eps$. Posons $b=a-I$ et $M=\|b\|_{\alpha}$. D\'efinissons, \`a pr\'esent, par r\'ecurrence, trois suites~$(b_{k}^-)_{k\ge 0}$, $(b_{k}^+)_{k\ge 0}$ et~$(\tilde{b}_{k})_{k\ge 0}$ de $M_{q}(\Bs_{\alpha}^-)$, $M_{q}(\Bs_{\alpha}^+)$ et $M_{q}(\Cs_{\alpha})$ v\'erifiant les conditions suivantes : quel que soit $k\ge 0$, nous avons
\begin{enumerate}
\item $\|b^-_{k}\|_{\alpha}^-\le DM\beta^{k}$ ;
\item $\|b^+_{k}\|_{\alpha}^+\le DM\beta^{k}$ ;
\item $\|\tilde{b}_{k}\|_{\alpha}\le M\beta^k$
\end{enumerate} 
et, quel que soit $k\ge 1$, nous avons
\begin{enumerate}
\item[4.] $\psi_{\alpha}^-(I+b_{k}^-)\, (I+\tilde{b}_{k})\, \psi_{\alpha}^+(I+b_{k}^+) = I+\tilde{b}_{k-1}.$
\end{enumerate} 

Initialisons la r\'ecurrence en posant $\tilde{b}_{0}=b$. La troisi\`eme propri\'et\'e est alors v\'erifi\'ee, par la d\'efinition m\^eme de $M$. Posons $b_{0}^-=0$ et $b_{0}^+=0$. Les premi\`ere et deuxi\`eme propri\'et\'es sont alors trivialement v\'erifi\'ees.

Soit $k\ge 0$ tels que $b_{k}^-$, $b_{k}^+$ et $\tilde{b}_{k}$ soient d\'ej\`a construits et v\'erifient les propri\'et\'es demand\'ees. Nous avons alors 
$$\|\tilde{b}_{k}\|_{\alpha}\le M\beta^k\le M\le \eps$$ 
et le lemme pr\'ec\'edent appliqu\'e avec $b=\tilde{b}_{k}$ nous fournit trois matrices $b^-$, $b^+$ et $\tilde{b}$. Posons $b_{k+1}^-=b^-$, $b_{k+1}^+=b^+$ et $\tilde{b}_{k+1}= \tilde{b}$. La quatri\`eme propri\'et\'e est alors v\'erifi\'ee.

Nous disposons, en outre, des in\'egalit\'es suivantes : $\|\tilde{b}_{k+1}\|_{\alpha}\le \beta\|\tilde{b}_{k}\|_{\alpha}$, $\|b^-_{k+1}\|_{\alpha}^-\le D\|\tilde{b}_{k}\|_{\alpha}$ et $\|b^+_{k+1}\|_{\alpha}^+\le D\|\tilde{b}_{k}\|_{\alpha}$. On en d\'eduit que les trois premi\`eres propri\'et\'es sont \'egalement v\'erifi\'ees.\\

Pour $n\in\N^*$, posons 
$$P_{n}=(I+b_{1}^-)\cdots(I+b_{n}^-) \in M_{q}(\Bs_{\alpha}^-)$$
et
$$Q_{n}=(I+b_{n}^+)\cdots(I+b_{1}^+) \in M_{q}(\Bs_{\alpha}^+).$$
De la quatri\`eme propri\'et\'e on d\'eduit que, quel que soit $n\in\N$, nous avons
$$a = \psi_{\alpha}^-(P_{n})\, (I+\tilde{b}_{n})\, \psi_{\alpha}^+(Q_{n}).$$
En utilisant les trois premi\`eres et le fait que $\beta\le 1/2$, nous obtenons
$$\sum_{k\ge 0} \|b_{k}^-\|_{\alpha}^- = DM \sum_{k\ge 0} \beta^k = 2DM \le \frac{1}{4}.$$
D'apr\`es le lemme \ref{convergence}, la suite $(P_{n})_{n\ge 0}$ converge dans $\Bs_{\alpha}^-$ vers une matrice inversible $c^- \in GL_{q}(\Bs_{\alpha}^-)$ v\'erifiant 
$$\|c^--I\|_{\alpha}^- \le 2\, \sum_{k\ge 0} \|b_{k}^-\|_{\alpha}^-\le 4DM \le 4D\|a-I\|_{\alpha}.$$
De m\^eme, la suite $(Q_{n})_{n\ge 0}$ converge dans $\Bs_{\alpha}^+$ vers une matrice inversible $c^+ \in GL_{q}(\Bs_{\alpha}^+)$ v\'erifiant 
$$\|c^+-I\|_{\alpha}^+ \le 2\,\sum_{k\ge 0} \|b_{k}^+\|_{\alpha}^+\le 4DM \le 4D\|a-I\|_{\alpha}.$$
Puisque la suite~$(\tilde{b}_{n})_{n\ge 0}$ converge vers~$0$, la suite~$(\psi_{\alpha}^-(P_{n})\, \psi_{\alpha}^+(Q_{n}))_{n\ge 0}$ converge vers~$a$. On en d\'eduit que 
$$a=\psi_{\alpha}^-(c^-)\, \psi_{\alpha}^+(c^+).$$
\end{proof}

\subsection{Prolongement de sections d'un faisceau}

Soit~$\Omega$ un syst\`eme de Banach associ\'e au couple $(K^-,K^+)$. Pour d\'emontrer les th\'eor\`emes~A et~B, nous chercherons \`a prolonger des sections de faisceaux. Pour ce faire, nous introduisons une nouvelle propri\'et\'e pour les syst\`emes de Cousin ; il s'agit d'une propri\'et\'e d'approximation.

\begin{defi}\index{Systeme@Syst\`eme!de Cousin-Runge}
Un {\bf syst\`eme de Cousin-Runge associ\'e au couple $(K^-,K^+)$} est un syst\`eme de Cousin associ\'e au m\^eme couple et tel que, quels que soient~$\alpha\in A$, $p,q\in\N$, $s_{1},\ldots,s_{p},t_{1},\ldots,t_{q}\in\Cs_{\alpha}$ et $\delta\in\R_{+}^*$, nous nous trouvions dans l'une des deux situations suivantes : soit il existe un \'el\'ement inversible~$f$ de~$\Bs_{\alpha}^+$, des \'el\'ements $s'_{1},\ldots,s'_{p}$ de~$\Bs_{\alpha}^+$ et $t'_{1},\ldots,t'_{q}$ de~$\Bs_{\alpha}^-$ tels que, quel que soient $i\in\cn{1}{p}$ et $j\in\cn{1}{q}$, on ait
\begin{enumerate}[\it i)]
\item $\|\psi_{\alpha}^+(f^{-1}s_{i}-s'_{i})\|_{\alpha}\,\|\psi_{\alpha}^+(f)\, \psi_{\alpha}^-(t_{j})\|_{\alpha} \le \delta$ ;
\item $\|\psi_{\alpha}^+(f^{-1}s_{i})\|_{\alpha}\, \|\psi_{\alpha}^+(f)\, \psi_{\alpha}^-(t_{j})-\psi_{\alpha}^-(t'_{j})\|_{\alpha}\le \delta$,
\end{enumerate}
soit il existe un \'el\'ement inversible~$f$ de~$\Bs_{\alpha}^-$, des \'el\'ements $s'_{1},\ldots,s'_{p}$ de~$\Bs_{\alpha}^+$ et $t'_{1},\ldots,t'_{q}$ de~$\Bs_{\alpha}^-$ tels que, quel que soient $i\in\cn{1}{p}$ et $j\in\cn{1}{q}$, on ait
\begin{enumerate}[\it i)]
\item $\|\psi_{\alpha}^-(f^{-1}t_{j}-t'_{j})\|_{\alpha}\,\|\psi_{\alpha}^-(f)\, \psi_{\alpha}^+(s_{i})\|_{\alpha} \le \delta$ ;
\item $\|\psi_{\alpha}^-(f^{-1}t_{j})\|_{\alpha}\, \|\psi_{\alpha}^-(f)\, \psi_{\alpha}^+(s_{i})-\psi_{\alpha}^+(s'_{i})\|_{\alpha}\le \delta$.
\end{enumerate}
\end{defi}

Nous utiliserons cette propri\'et\'e par le biais du lemme suivant.

\begin{lem}\label{lemattache}
Supposons que~$\Omega$ soit un syst\`eme de Cousin-Runge. Soit~$\eps>0$ le nombre r\'eel dont le th\'eor\`eme \ref{Cartan} assure l'existence. Soient $\Fs$ un faisceau de $\Os_{M}$-modules, $p,q\in\N$, \mbox{$T^-\in \Fs(K^-)^p$}, \mbox{$T^+\in\Fs(K^+)^q$}, $U_{0}=(u_{a,i})\in M_{p,q}(\Os(L))$ et $V_{0}=(v_{b,j})\in M_{q,p}(\Os(L))$ telles que, dans~$\Fs(L)$, on ait 
\begin{enumerate}
\item[a)] $T^-=U_{0}\, T^+$ ;
\item[b)] $T^+=V_{0}\, T^-$.
\end{enumerate}

Supposons qu'il existe~$\alpha\in A$, $U,U_{\delta}\in M_{p,q}(\Bs_{\alpha}^+)$ et $V,V_{\delta}\in M_{q,p}(\Bs_{\alpha}^-)$ tels que
\begin{enumerate}
\item[c)] $\varphi^+_{\alpha}(U) = U_{0}$ ;
\item[d)] $\varphi^-_{\alpha}(V) = V_{0}$ ;
\item[e)] $\|\psi_{\alpha}^+(U_{\delta}-U)\|_{\alpha}\, \|\psi_{\alpha}^-(V)\|_{\alpha}  < \eps$ ;
\item[f)] $\|\psi_{\alpha}^+(U)\|_{\alpha}\, \|\psi_{\alpha}^-(V_{\delta}-V)\|_{\alpha} < \eps$.
\end{enumerate}

Alors il existe $S^-\in \Fs(M)^p$, $S^+\in\Fs(M)^q$, \mbox{$A^-\in GL_{p}(\Os(K^-))$} et \mbox{$A^+\in GL_{q}(\Os(K^+))$} v\'erifiant
\begin{enumerate}[\it i)]
\item $S^-=A^- \, T^- \textrm{ dans } \Fs(K^-)$ ;
\item $S^+=A^+ \, T^+ \textrm{ dans } \Fs(K^+)$.
\end{enumerate}
\end{lem}
\begin{proof}
Supposons qu'il existe~$\alpha\in A$, $U,U_{\delta}\in M_{p,q}(\Bs_{\alpha}^+)$ et $V,V_{\delta}\in M_{q,p}(\Bs_{\alpha}^-)$ tels que
\begin{enumerate}
\item[c)] $\rho^+_{\alpha}(U) = U_{0}$ ;
\item[d)] $\rho^-_{\alpha}(V) = V_{0}$ ;
\item[e)] $\|\psi_{\alpha}^+(U_{\delta}-U)\|_{\alpha}\, \|\psi_{\alpha}^-(V)\|_{\alpha}  < \eps$ ;
\item[f)] $\|\psi_{\alpha}^+(U)\|_{\alpha}\, \|\psi_{\alpha}^-(V_{\delta}-V)\|_{\alpha} < \eps$.
\end{enumerate}

Posons $T_{\delta}^- = \rho_{\alpha}^-(V_{\delta})\, T^-$ dans $\Fs(K^-)$. Dans~$\Fs(L)$, nous avons alors
$$T_{\delta}^- - T^+ = \rho_{\alpha}(\psi_{\alpha}^-(V_{\delta}-V))\, T^- = \rho_{\alpha}(\psi_{\alpha}^-(V_{\delta}-V)\, \psi_{\alpha}^+(U))\, T^+.$$
Posons 
$$A=I + \psi_{\alpha}^-(V_{\delta}-V)\, \psi_{\alpha}^+(U) \in M_{q}(\Cs_{\alpha}).$$
Nous avons alors
$$T_{\delta}^-=\rho_{\alpha}(A)\, T^+ \textrm{ dans } \Fs(L)$$
et
$$\|A-I\|_{\alpha} \le \|\psi_{\alpha}^-(V_{\delta}-V)\|_{\alpha}\, \|\psi_{\alpha}^+(U)\|_{\alpha} < \eps.$$

D'apr\`es le th\'eor\`eme~\ref{Cartan}, il existe deux matrices $C^-\in GL_{q}(\Bs_{\alpha}^-)$ et $C^+\in GL_{q}(\Bs_{\alpha}^+)$ telles que
$$A = \psi_{\alpha}^-(C^-)\, \psi_{\alpha}^+(C^+).$$
Posons
$$A^+=\rho_{\alpha}^+(C^+)\in GL_{q}(\Os(K^+)).$$
Dans $\Fs(L)$, nous avons alors
$$A^+\, T^+ = \rho_{\alpha}\left( \psi_{\alpha}^-\left((C^-)^{-1}\right)A\right)\, T^+ = \rho_{\alpha}\left( \psi_{\alpha}^-\left((C^-)^{-1}\right)\right)\, T_{\delta}^-.$$
Nous pouvons donc d\'efinir un \'el\'ement~$S^+$ de~$\Fs(M)^q$ par 
\begin{enumerate}
\item $S^+_{|K^-} = \rho_{\alpha}^-\left((C^-)^{-1}\right)\, T_{\delta}^-$ ;
\item $S^+_{|K^+} = A^+\, T^+$.
\end{enumerate}
On proc\`ede de m\^eme pour construire la section $S^-$.
\end{proof}

Nous parvenons maintenant au r\'esultat permettant de recoller les sections d'un faisceau.

\begin{thm}\label{attachefaisceau}\index{Faisceau!recollement de sections}
Supposons que~$\Omega$ soit un syst\`eme de Cousin-Runge. Soit $\Fs$ un faisceau de $\Os_{M}$-modules. Supposons qu'il existe deux entiers~$p$ et~$q$, une famille $(t_{1}^-, \ldots,t_{p}^-)$ d'\'el\'ements de~$\Fs(K^-)$ et une famille $(t_{1}^+,\ldots,t_{q}^+)$ d'\'el\'ements de~$\Fs(K^+)$ dont les restrictions \`a~$L$ engendrent le m\^eme sous-$\Os(L)$-module de $\Fs(L)$. Alors il existe $s_{1}^-, \ldots,s_{p}^-,s_{1}^+,\ldots,s_{q}^+ \in \Fs(M)$, $a^-\in GL_{p}(\Os(K^-))$ et $a^+\in GL_{q}(\Os(K^+))$ tels que
$$\begin{pmatrix} s_{1}^-\\ \vdots\\ s_{p}^-\end{pmatrix} = a^- \begin{pmatrix} t_{1}^-\\ \vdots\\ t_{p}^-\end{pmatrix} \textrm{ dans } \Fs(K^-)^p$$
et 
$$\begin{pmatrix} s_{1}^+\\ \vdots\\ s_{q}^+\end{pmatrix} = a^+ \begin{pmatrix} t_{1}^+\\ \vdots\\ t_{q}^+\end{pmatrix} \textrm{ dans } \Fs(K^+)^q.$$
\end{thm}
\begin{proof}
Posons 
$$T^-=\begin{pmatrix} t_{1}^-\\ \vdots\\ t_{p}^- \end{pmatrix} \textrm{ et } T^+=\begin{pmatrix} t_{1}^+\\ \vdots\\ t_{q}^+ \end{pmatrix}.$$
Par hypoth\`ese, il existe~$\alpha\in A$, $U=(u_{a,i})\in M_{p,q}(\Cs_{\alpha})$ et $V=(v_{b,j})\in M_{q,p}(\Cs_{\alpha})$ tels qu'on ait les \'egalit\'es 
$$T^-=\rho_{\alpha}(U)\, T^+ \textrm{ et } T^+=\rho_{\alpha}(V)\, T^- \textrm{ dans } \Fs(L).$$

Consid\'erons le nombre r\'eel $\eps>0$ dont le th\'eor\`eme \ref{Cartan} assure l'existence. Quitte \`a \'echanger les compacts~$K^-$ et~$K^+$, nous pouvons supposer que la premi\`ere propri\'et\'e des syst\`emes de Cousin-Runge est v\'erifi\'ee. Il existe alors un \'el\'ement inversible~$f$ de~$\Bs_{\alpha}^+$, des \'el\'ements $\hat{u}_{a,i}$ de~$\Bs_{\alpha}^+$, pour $(a,i)\in\cn{1}{p}\times\cn{1}{q}$, et $\hat{v}_{b,j}$ de~$\Bs_{\alpha}^-$, pour $(b,j)\in\cn{1}{q}\times\cn{1}{p}$, v\'erifiant les conditions suivantes : quel que soient $(a,i)\in\cn{1}{p}\times\cn{1}{q}$ et $(b,j)\in\cn{1}{q}\times\cn{1}{p}$, nous avons
\begin{enumerate}[\it i)]
\item $\|\psi_{\alpha}^+(f^{-1}u_{a,i}-\hat{u}_{a,i})\|_{\alpha}\, \|\psi_{\alpha}^+(f)\, \psi_{\alpha}^-(v_{b,j})\|_{\alpha}  <  \eps$ ;
\item $\|\psi_{\alpha}^+(f^{-1}u_{a,i})\|_{\alpha}\, \|\psi_{\alpha}^+(f)\, \psi_{\alpha}^-(v_{b,j})-\psi_{\alpha}^-(\hat{v}_{b,j})\|_{\alpha}  <  \eps$.
\end{enumerate}

Les matrices $T^-$, $\rho_{\alpha}^+(f)T^+$, $\rho_{\alpha}(\psi_{\alpha}^+(f^{-1}) U)$ et $\rho_{\alpha}(\psi_{\alpha}^+(f)V)$ v\'erifient donc les hypoth\`eses du lemme \ref{lemattache}. Par cons\'equent, il existe $S^-\in \Fs(M)^p$, $S^+\in\Fs(M)^q$, $A^-\in GL_{p}(\Os(K^-))$ et $A^+\in GL_{q}(\Os(K^+))$ tels que
\begin{enumerate}
\item $S^-_{|K^-}=A^- \, T^-$ ;
\item $S^+_{|K^+}=A^+ \, \rho_{\alpha}^+(f)T^+ $.
\end{enumerate}

En posant $$\begin{pmatrix} s_{1}^-\\ \vdots\\ s_{p}^-\end{pmatrix} = S^-, \, \begin{pmatrix} s_{1}^+\\ \vdots\\ s_{q}^+\end{pmatrix} = S^+,$$
ainsi que $a^-=A^-$ et $a^+  = \rho_{\alpha}^+(f) A^+$, on obtient le r\'esultat souhait\'e. Remarquons que $a^+\in  GL_{q}(\Os(K^+))$ car $f$ est inversible dans $\Bs_{\alpha}^+$. 
\end{proof}

Indiquons, \`a pr\'esent, la fa\c{c}on dont ce r\'esultat permet de d\'emontrer le th\'e\-o\-r\`e\-me~A.

\begin{cor}\label{corA}\index{Theoreme A@Th\'eor\`eme A!pour un syst\`eme de Cousin-Runge}
Supposons qu'il existe un syst\`eme de Cousin-Runge associ\'e au couple $(K^-,K^+)$. Soit~$\Fs$ un faisceau de $\Os_{M}$-modules de type fini qui satisfait le th\'eor\`eme~A sur les compacts~$K^-$ et~$K^+$. Alors il le satisfait encore sur leur r\'eunion~$M$. 
\end{cor}
\begin{proof}
D'apr\`es le lemme \ref{thAcompact} il existe deux entiers~$p$ et~$q$, une famille $(t_{1}^-, \ldots,t_{p}^-)$ d'\'el\'ements de~$\Fs(K^-)$ dont l'image engendre le $\Os_{X,x}$-module~$\Fs_{x}$ en tout point~$x$ de~$K^+$ et une famille $(t_{1}^+,\ldots,t_{q}^+)$ d'\'el\'ements de~$\Fs(K^+)$ dont l'image engendre le $\Os_{X,x}$-module~$\Fs_{x}$ en tout point~$x$ de~$K^+$. En particulier, les restrictions \`a~$L$ de ces deux familles engendrent toutes deux~$\Fs(L)$. Nous pouvons donc appliquer le th\'eor\`eme \ref{attachefaisceau}. Les sections $s_{1}^-,\ldots,s_{p}^-,s_{1}^+,\ldots,s_{q}^+$ de~$\Fs(M)$ dont il assure l'existence engendrent le $\Os_{X,x}$-module~$\Fs_{x}$ en tout point~$x$ de~$M$. On en d\'eduit le r\'esultat annonc\'e.
\end{proof}

Expliquons, \`a pr\'esent, comment d\'eduire le th\'eor\`eme~B du th\'eor\`eme~A. Insistons sur le fait que, dans la proposition qui suit, nous n'avons besoin d'associer aucun syst\`eme de Banach au couple $(K^-,K^+)$.

\begin{prop}\label{propB}\index{Theoreme B@Th\'eor\`eme B!pour un couple de compacts}
Supposons que pour tout \'el\'ement~$f$ de~$\Os(L)$, il existe un \'el\'ement~$f^-$ de~$\Os(K^-)$ et un \'el\'ement~$f^+$ de~$\Os(K^+)$ qui v\'erifient l'\'egalit\'e
$$f = f^- + f^+ \textrm{ dans } \Os_{X,x}.$$
Supposons \'egalement que tout faisceau de $\Os_{L}$-modules coh\'erent satisfait le th\'eor\`eme~B. 

Soit~$\Fs$ un faisceau de $\Os_{M}$-modules coh\'erent qui satisfait le th\'eor\`eme~A sur~$M$. Soit 
$$0\to \Fs \xrightarrow[]{d} \Is_{0} \xrightarrow[]{d} \Is_{1} \xrightarrow[]{d} \cdots$$
une r\'esolution flasque du faisceau~$\Fs$. Soient $q\in\N^*$ et~$\gamma$ un cocycle de degr\'e~$q$ sur~$M$. Si~$\gamma$ est un cobord au voisinage des compacts~$K^-$ et~$K^+$, alors c'est un cobord au voisinage de leur r\'eunion~$M$.
\end{prop}
\begin{proof}
Supposons qu'il existe $\beta^-\in \Is_{q-1}(K^-)$ et $\beta^+\in \Is_{q-1}(K^+)$ tels que 
$$d(\beta^{-})=\gamma \textrm{ dans } \Is_{q}(K^-) \textrm{ et }d(\beta^{+})=\gamma  \textrm{ dans } \Is_{q}(K^+).$$ 

Supposons, tout d'abord que $q\ge 2$. Nous avons $d(\beta^--\beta^+)=0$ dans~$\Is_{q}(L)$. D'apr\`es le th\'eor\`eme~B, nous avons $H^{q-1}(L,\Fs)=0$. Par cons\'equent, il existe $\alpha\in\Is_{q-2}(L)$ telle que $d(\alpha)=\beta^--\beta^+$ dans~$\Is_{q-1}(L)$. Puisque le faisceau~$\Is_{q-2}$ est flasque, $\alpha$ se prolonge en une section sur~$M$ que nous noterons identiquement. D\'efinissons $\beta\in\Is_{q-1}(M)$ par $\beta = \beta^-$ au-dessus de $K^-$ et $\beta=\beta^++d(\alpha)$ au-dessus de $K^+$. Nous avons alors l'\'egalit\'e
$$d(\beta)=\gamma \textrm{ dans } \Is_{q}(M)$$
et~$\gamma$ est un cobord au voisinage de~$M$.

Int\'eressons-nous, \`a pr\'esent, au cas $q=1$. Nous avons alors $d(\beta^--\beta^+)=0$ dans~$\Is_{1}(L)$. On en d\'eduit que $\beta^--\beta^+$ est un \'el\'ement de~$\Fs(L)$. D'apr\`es le th\'eor\`eme~A et le lemme \ref{thAcompact}, il existe un entier positif~$m$ et une famille $(u_{1},\ldots,u_{m})$ d'\'el\'ements de~$\Fs(M)$ dont les images engendrent le $\Os_{X,x}$-module~$\Fs_{x}$ en tout point~$x$ de~$M$. En d'autres termes, l'application 
$$\begin{array}{ccc}
\Os^m &\to& \Fs\\
(a_{1},\ldots,a_{m}) & \mapsto &\disp \sum_{i=1}^m a_{i}\, u_{i}
\end{array}$$ 
est surjective au-dessus de $M$. Son noyau~$\Ns$ est un faisceau de $\Os_{M}$-modules coh\'erent. D'apr\`es le th\'eor\`eme~B, nous avons $H^1(L,\Ns)=0$. On en d\'eduit que la famille $(u_{1},\ldots,u_{m})$ engendre $\Fs(L)$ en tant que $\Os(L)$-module. Par cons\'equent, il existe $\lambda_{1},\ldots,\lambda_{m} \in \Os(L)$ tels que 
$$\beta^--\beta^+ = \sum_{i=1}^m \lambda_{i}\, u_{i} \textrm{ dans } \Fs(L).$$

Par hypoth\`ese, quel que soit $i\in\cn{1}{m}$, il existe $\lambda_{i}^-\in\Os(K^-)$ et $\lambda_{i}^+\in\Os(K^+)$ tels que 
$$\lambda_{i}=\lambda_{i}^--\lambda_{i}^+ \textrm{ dans } \Os(L).$$
Nous avons alors l'\'egalit\'e
$$\beta^- - \sum_{i=1}^m \lambda_{i}^-\, u_{i} = \beta^+ - \sum_{i=1}^m \lambda_{i}^+\, u_{i} \textrm{ dans } \Fs(L).$$ 
On en d\'eduit l'existence d'un \'el\'ement~$\beta$ de~$\Is_{0}(M)$ v\'erifiant $d(\beta)=\gamma$ dans~$\Is_{1}(M)$.
\end{proof}

\section{Parties compactes de la base}\label{pcdlb}

Dans la suite de ce chapitre, nous reprenons les notations du chapitre \ref{chapitredroite}. 
\newcounter{nodStein}\setcounter{nodStein}{\thepage}


\bigskip

Nous allons maintenant appliquer les r\'esultats obtenus au paragraphe pr\'ec\'edent pour d\'emontrer que certaines parties compactes de l'espace de base~$B$ sont de Stein. \`A cet effet, nous allons exhiber des syst\`emes de Cousin-Runge. \'Enon\c{c}ons tout d'abord un r\'esultat de th\'eorie des nombres.

\begin{lem}\label{reseau}
Il existe $C\in\R$ tel que, quel que soit 
$$(x_{\sigma})_{\sigma\in\Sigma_{\infty}} \in \prod_{\sigma\in\Sigma_{\infty}} \hat{K}_{\sigma},$$ 
il existe $y\in A$ v\'erifiant
$$\forall \sigma\in\Sigma_{\infty},\, |y-x_{\sigma}|_{\sigma} \le C.$$
\end{lem}
\begin{proof}
Notons $r_{1}$ le nombre de places r\'eelles de $K$ et $2r_{2}$ le nombre de places complexes de $K$. Le r\'esultat d\'ecoule directement du fait que l'image de l'anneau des entiers $A$ par l'application 
$$\begin{array}{rcl}
K & \to & \R^{r_{1}}\times \C^{r_{2}} \simeq \R^{r_{1}+2r_{2}}\\
x & \mapsto & (\sigma(x))_{\sigma\in\Sigma_{\infty}}
\end{array}$$
est un r\'eseau.
\end{proof}

Le lemme qui suit sera utile pour exhiber des syst\`emes de Cousin.

\begin{figure}[htb]
\begin{center}
\input{K0+-.pstex_t}
\caption{Les compacts $K_{0}^-$ et $K_{0}^+$.}
\end{center}
\end{figure}
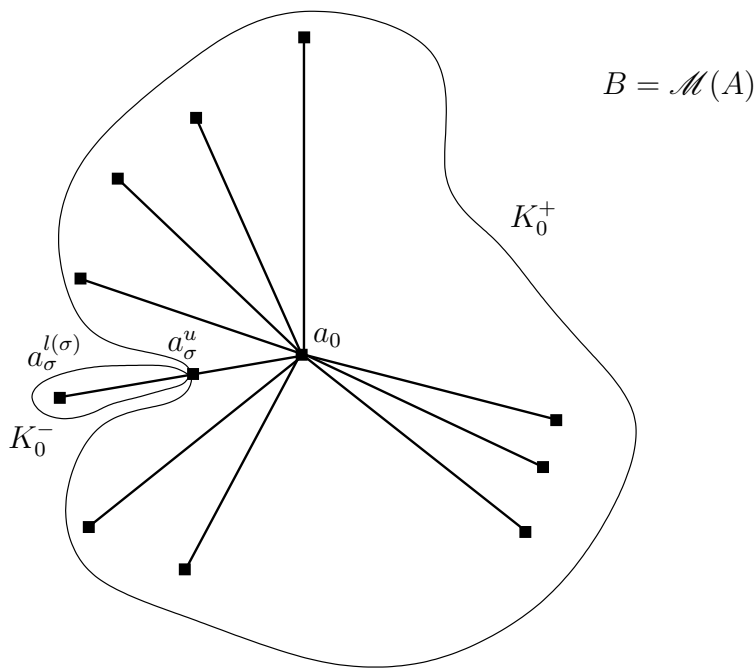

\begin{lem}\label{Cousinbase}
Soient $\sigma\in\Sigma$ et $u\in \of{]}{0,l(\sigma)}{[}$. Posons 
$$K_{0}^- = \of{[}{a_{\sigma}^u,a_{\sigma}^{l(\sigma)}}{]},\ K_{0}^+=B\setminus\of{]}{a_{\sigma}^u,a_{\sigma}^{l(\sigma)}}{]} \textrm{ et } L_{0}=K_{0}^-\cap K_{0}^+=\{a_{\sigma}^u\}.$$ 
\newcounter{KLMzero}\setcounter{KLMzero}{\thepage}

Il existe $D\in\R$ tel que, quel que soit $a\in\Bs(L_{0})$, il existe $a^-\in\Bs(K_{0}^-)$ et $a^+\in\Bs(K_{0}^+)$ v\'erifiant les propri\'et\'es suivantes :
\begin{enumerate}[\it i)]
\item $a=a^--a^+$ dans $\Bs(L_{0})$ ;
\item $\|a^-\|_{K_{0}^-} \le D\,  \|a\|_{L_{0}}$ ;
\item $\|a^+\|_{K_{0}^+} \le D\,\|a\|_{L_{0}}$.
\end{enumerate}
\end{lem}
\begin{proof}
Consid\'erons la constante $C\in\R$ dont le lemme \ref{reseau} assure l'existence. Nous pouvons, sans perdre de g\'en\'eralit\'e, supposer que $C\ge 1$. Soit~$a\in\Bs(L_{0})$. Remarquons que l'anneau~$\Bs(L_{0})$ est isomorphe au corps~$\hat{K}_{\sigma}$ muni de la valeur absolue~$|.|_{\sigma}^u$. Dans le raisonnement qui suit, nous aurons besoin de conna\^{\i}tre le type de $\sigma$.

\bigskip

Supposons, tout d'abord, que $\sigma\in\Sigma_{\infty}$. Dans ce cas, nous avons
$$\Bs(K_{0}^-) = \hat{K}_{\sigma} \textrm{ et } \|.\|_{K_{0}^-} = \max(|.|_{\sigma}^u,|.|_{\sigma})$$
et 
$$\Bs(K_{0}^+) = A \textrm{ et } \|.\|_{K_{0}^+} = \max_{\sigma'\in\Sigma_{\infty}\setminus\{\sigma\}}(|.|_{\sigma'}).$$
Distinguons plusieurs cas. Supposons, tout d'abord, que $|a|_{\sigma}\ge 1$. Puisque~\mbox{$\sigma\in\Sigma_{\infty}$}, le nombre r\'eel~$|a|_{\sigma}^u$ est un \'el\'ement de~$\hat{K}_{\sigma}$. Par d\'efinition de $C$, il existe $b \in A$ v\'erifiant les propri\'et\'es suivantes : 
\begin{enumerate}
\item $\left|b + |a|_{\sigma}^{u}\right|_{\sigma}\le C$ ;
\item $\forall \sigma'\in\Sigma_{\infty}\setminus\{\sigma\}$, $|b|_{\sigma'}\le C$.
\end{enumerate} 
Quel que soit $\sigma'\in \Sigma_{\infty}\setminus\{\sigma\}$, nous avons donc
$$|b|_{\sigma'} \le C\le  C\, |a|_{\sigma}^{u}.$$ 
De nouveau, nous allons distinguer deux cas. Supposons, tout d'abord, que~\mbox{$|b|_{\sigma}\ge 1$}. De la premi\`ere in\'egalit\'e, nous tirons 
$$\begin{array}{rcl}
|b|_{\sigma} &\le& |a|_{\sigma}^{u} +  C\\
&\le & (C+1)\, |a|_{\sigma}^{u}.
\end{array}$$ 
Puisque $|b|_{\sigma}\ge 1$, nous avons \'egalement $|b|_{\sigma}^{u} \le  (C+1)\, |a|_{\sigma}^{u}$. Si $|b|_{\sigma}\le 1$, nous avons encore $|b|_{\sigma} \le  (C+1)\, |a|_{\sigma}^{u}$.

Supposons, \`a pr\'esent, que~$|a|_{\sigma}\le 1$. Nous avons alors l'\'egalit\'e~\mbox{$\|a\|_{K_{0}^-}=|a|_{\sigma}^u$}. Nous posons $b=0$.

Dans tous les cas, il existe $D\in\R$ tel que $\|a+b\|_{K_{0}^-} \le D\, |a|_{\sigma}^{u}$ et \mbox{$\|b\|_{K_{0}^+} \le D\, |a|_{\sigma}^{u}$}. Nous pouvons, par exemple, choisir $D=C+2$. Les \'el\'ements~$a^-=a+b$ de~$\Bs(K_{0}^-)$ et~$a^+=b$ de~$\Bs(K_{0}^+)$ v\'erifient les propri\'et\'es demand\'ees.

\bigskip

Supposons, \`a pr\'esent, que~$\sigma\in\Sigma_{f}$. Nous avons alors
$$\Bs(K_{0}^-) = \hat{A}_{\sigma} \textrm{ et } \|.\|_{K_{0}^-} = |.|_{\sigma}^u$$
et 
$$\Bs(K_{0}^+) = A\left[\frac{1}{\sigma}\right] \textrm{ et } \|.\|_{K_{0}^+} = \max\left(|.|_{\sigma}^u,\max_{\sigma'\in\Sigma_{\infty}}(|.|_{\sigma})\right).$$
Comme pr\'ec\'edemment, nous allons distinguer plusieurs cas. Pour commencer, supposons que $|a|_{\sigma}\ge 1$. D'apr\`es le th\'eor\`eme d'approximation fort, il existe un \'el\'ement~$b$ de $\bigcap_{\sigma'\in\Sigma_{f}\setminus\{\sigma\}} A_{\sigma}$ v\'erifiant
$$|b+a|_{\sigma}\le 1.$$ 
En partculier, $b+a$ appartient \`a~$\Bs(K_{0}^-)$. Par d\'efinition de la constante $C$, il existe $c \in A$ v\'erifiant la propri\'et\'e suivante : quel que soit $\sigma'\in\Sigma_{\infty}$, nous avons $|b+c|_{\sigma'}\le C$. On en d\'eduit que, quel que soit~\mbox{$\sigma'\in \Sigma_{\infty}$}, nous avons
$$|b+c|_{\sigma'} \le C\le C\, |a|_{\sigma}^{u}.$$ 
En outre, nous avons
$$|b+c|_{\sigma}^{u}  \le \max\left(|a|_{\sigma}^{u} , |a+b|_{\sigma}^{u}, |c|_{\sigma}^{u}\right) \le |a|_{\sigma}^{u}.$$

Supposons, \`a pr\'esent, que $|a|_{\sigma}\le 1$. Dans ce cas, $a$ appartient \`a~$\Bs(K^-)$. Nous posons $b=c=0$. 

Dans tous les cas, il existe $D\in\R$ tel que $\|a+b+c\|_{K_{0}^-} \le D\, |a|_{\sigma}^{u}$ et \mbox{$\|b+c\|_{K_{0}^+} \le D\, |a_{k}|_{\sigma}^{u}$}. Nous pouvons, par exemple, choisir $D=C+1$. Les \'el\'ements $a^-=a+b+c$ de~$\Bs(K_{0}^-)$ et~$a^+=b+c$ de~$\Bs(K_{0}^+)$ v\'erifient les propri\'et\'es demand\'ees.
\end{proof}

Int\'eressons-nous, \`a pr\'esent, \`a la propri\'et\'e d'approximation qui intervient dans la d\'efinition des syst\`emes de Cousin-Runge.

\begin{lem}\label{appbase}
Soient $\sigma\in\Sigma$ et $u\in \of{]}{0,l(\sigma)}{[}$. Posons 
$$K_{0}^- = \of{[}{a_{\sigma}^u,a_{\sigma}^{l(\sigma)}}{]},\ K_{0}^+=B\setminus\of{]}{a_{\sigma}^u,a_{\sigma}^{l(\sigma)}}{]} \textrm{ et } L_{0}=K_{0}^-\cap K_{0}^+=\{a_{\sigma}^u\}.$$ 

Soient $p,q\in\N$ et $s_{1},\ldots,s_{p},t_{1},\ldots,t_{q}\in\Bs(L_{0})$. Soit $\delta\in\R_{+}^*$. Si~$\sigma$ appartient \`a~$\Sigma_{f}$, alors il existe un \'el\'ement inversible~$f$ de~$\Bs(K_{0}^+)$ et des \'el\'ements $s'_{1},\ldots,s'_{p}$ de~$\Bs(K_{0}^+)$ et $t'_{1},\ldots,t'_{q}$ de~$\Bs(K_{0}^-)$ tels que, quel que soient $i\in\cn{1}{p}$ et $j\in\cn{1}{q}$, on ait
\begin{enumerate}[\it i)]
\item $\|f^{-1}s_{i}-s'_{i}\|_{L_{0}}\,\|ft_{j}\|_{L_{0}} \le \delta$ ;
\item $\|f^{-1}s_{i}\|_{L_{0}}\, \|ft_{j}-t'_{j}\|_{L_{0}}\le \delta$.
\end{enumerate}

Si~$\sigma$ appartient \`a~$\Sigma_{\infty}$, alors il existe un \'el\'ement inversible~$g$ de~$\Bs(K_{0}^-)$ et des \'el\'ements $s''_{1},\ldots,s''_{p}$ de~$\Bs(K_{0}^+)$ et $t''_{1},\ldots,t''_{q}$ de~$\Bs(K_{0}^-)$ tels que, quel que soient $i\in\cn{1}{p}$ et $j\in\cn{1}{q}$, on ait
\begin{enumerate}[\it i)]
\item $\|gs_{i}-s''_{i}\|_{L_{0}}\, \|g^{-1}t_{j}\|_{L_{0}}\le \delta$ ;
\item $\|gs_{i}\|_{L_{0}}\,\|g^{-1}t_{i}-t''_{i}\|_{L_{0}}\le \delta$.
\end{enumerate}
\end{lem}
\begin{proof}
Posons $M=\max\{\|s_{i}\|_{L_{0}},\|t_{j}\|_{L_{0}},\, 1\le i\le p, 1\le j\le q\}$. L'anneau de Banach~$\Bs(L_{0})$ n'est autre que le corps~$\hat{K}_{\sigma}$ muni de la valeur absolue~$|.|_{\sigma}^u$. Par cons\'equent, pour tout \'el\'ement~$i$ de~$\cn{1}{n}$, il existe un \'el\'ement~$s^*_{i}$ de~$K$ tel que
$$\|s_{i}-s^*_{i}\|_{L_{0}} \le \delta.$$

\bigskip

Distinguons maintenant deux cas. Supposons, tout d'abord, que $\sigma$ est un \'el\'ement de~$\Sigma_{f}$. D'apr\`es le lemme \ref{diviseur}, il existe $h\in A$ telle que 
$$\left\{\begin{array}{l}
|h|_{\sigma}<1\ ;\\
\forall \sigma'\in\Sigma_{f}\setminus\{\sigma\},\, |h|_{\sigma'}=1.
\end{array}\right.$$
Il existe $N\in\N$ tel que, quel que soit $i\in\cn{1}{p}$, on ait 
$$h^N\, s^*_{i} \in \hat{A}_{\sigma}=\Bs(K_{0}^-).$$ 
Posons 
$$f=h^{-N}\in K.$$ 
Son image dans l'anneau~$\Bs(K_{0}^+)$ est inversible. En outre, quel que soit $i\in\cn{1}{p}$, nous avons 
$$\|f^{-1} s_{i} - f^{-1} s^*_{i}\|_{L_{0}} \le \|f^{-1}\|_{L_{0}}\, \|s_{i}-s^*_{i}\|_{L_{0}} \le \delta\, |f^{-1}|_{\sigma}^u.$$

Soit $j\in\cn{1}{q}$. D'apr\`es le th\'eor\`eme d'approximation fort, il existe un \'el\'ement~$t^*_{j}$ de $A[1/\sigma]$ tel que
$$\|ft_{j}-t^*_{j}\|_{L_{0}}\le \delta.$$
En particulier, la fonction~$t^*_{j}$ d\'efinit un \'el\'ement de~$\Bs(K_{0}^+)$ et, quel que soit~$i\in\cn{1}{p}$, nous avons 
$$\|f^{-1} s_{i}\|_{L_{0}}\,  \|ft_{j} - t^*_{j}\|_{L_{0}}\le |f^{-1}|_{\sigma}^u\, \|s_{i}\|_{L_{0}}\, \delta \le M \delta,$$
car $f^{-1}\in A$. Quel que soit $i\in\cn{1}{p}$, nous avons \'egalement
$$\|f^{-1} s_{i} - f^{-1} s^*_{i}\|_{L_{0}}\, \|ft_{j}\|_{L_{0}} \le \delta\, |f^{-1}|_{\sigma}^u\, |f|_{\sigma}^u\, \|t_{j}\|_{L_{0}} \le M\delta.$$

\bigskip

Supposons, \`a pr\'esent, que $\sigma\in\Sigma_{\infty}$. Il existe $g\in A$ tel que, quel que soit \mbox{$i\in\cn{1}{n}$}, on ait $gs^*_{i} \in A$. Remarquons que l'image de~$g$ dans~$\Bs(K_{0}^-)$ est inversible et que~$gs^*_{i}$ est un \'el\'ement de~$\Bs(K_{0}^+)$. En outre, quel que soit~$i\in\cn{1}{n}$, nous avons
$$\|gs_{i} - gs^*_{i}\|_{L_{0}} \le \|g\|_{L_{0}}\, \|s_{i}-s^*_{i}\|_{L_{0}} \le |g|_{\sigma}^u\, \delta.$$
Choisissons un nombre r\'eel~$N>0$ tel que
$$\forall i\in\cn{1}{p},\, \|gs_{i}\|_{L_{0}}\le N.$$

Soit $j\in\cn{1}{q}$. Il existe un \'el\'ement~$t^*_{j}$ de~$K$ tel que
$$\| g^{-1}t_{j} - t^*_{j}\|_{L_{0}}\le \frac{\delta}{N}.$$
La fonction $t^*_{j}$ d\'efinit un \'el\'ement de~$\Bs(K_{0}^-)$ et, quel que soit \mbox{$i\in\cn{1}{p}$}, v\'erifie
$$\|gs_{i}\|_{L_{0}}\, \|g^{-1}t_{j}-t^*_{j}\|_{L_{0}}\le \|gs_{i}\|_{L_{0}}\, \frac{\delta}{N} \le \delta.$$
Quel que soit $i\in\cn{1}{p}$, nous avons encore
$$\|gs_{i} - s^*_{i}\|_{L_{0}}\, \|g^{-1}t_{j}\|_{L_{0}} \le  |g|_{\sigma}^u\, \delta\, |g^{-1}|_{\sigma}^u\, \|t_{j}\|_{L_{0}} \le M\delta.$$
\end{proof}

Les lemmes qui pr\'ec\`edent nous permettent d'exhiber de nombreux syst\`emes de Cousin-Runge.

\begin{prop}\label{CousinRungebase}\index{Systeme@Syst\`eme!de Cousin-Runge!sur MA@sur $\Ms(A)$}
Soient $\sigma\in\Sigma$ et  $u\in \of{]}{0,l(\sigma)}{[}$. Posons 
$$K_{0}^- = \of{[}{a_{\sigma}^u,a_{\sigma}^{l(\sigma)}}{]},\ K_{0}^+=B\setminus\of{]}{a_{\sigma}^u,a_{\sigma}^{l(\sigma)}}{]}.$$ 
Soient~$K^-$ et~$K^+$ deux parties compactes et connexes de l'espace~$B$ dont l'intersection est le singleton~$\{a_{\sigma}^u\}$. Il existe un syst\`eme de Cousin-Runge associ\'e au couple $(K^-,K^+)$.
\end{prop}
\begin{proof}
Ce cas est particuli\`erement simple et nous allons construire un syst\`eme de Cousin-Runge dont l'ensemble~$A$ est r\'eduit \`a un seul \'el\'ement. Quitte \`a \'echanger les compacts~$K^-$ et~$K^+$, nous pouvons supposer que nous avons les inclusions
$$K^-\subset K_{0}^- \textrm{ et } K^+\subset K_{0}^+.$$
Posons
$$(\Bs^{-},\|.\|^-) = (\Bs(K_{0}^-),\|.\|_{K_{0}^-}),$$
$$(\Bs^{+},\|.\|^+) = (\Bs(K_{0}^+),\|.\|_{K_{0}^+})$$
et
$$(\Cs,\|.\|) = (\Hs(a_{\sigma}^u),|.|_{\sigma}^u).$$
On d\'efinit de mani\`ere \'evidente des morphismes born\'es~$\psi^-$ et~$\psi^+$ comme dans la d\'efinition des syst\`emes de Banach. La proposition \ref{isoBO} permet de d\'efinir \'egalement des morphismes $\rho^-$, $\rho^+$ et~$\rho$. L'ensemble de ces donn\'ees forme un syst\`eme de Banach associ\'e au couple $(K^-,K^+)$. Les deux lemmes qui pr\'ec\`edent assurent que c'est un syst\`eme de Cousin-Runge.
\end{proof}

De ce r\'esultat, nous allons d\'eduire que toute partie compacte et connexe de l'espace~$B$ est un espace de Stein.

\begin{thm}\label{Abase}\index{Theoreme A@Th\'eor\`eme A!pour un compact de MA@pour un compact de $\Ms(A)$}
Soit~$M$ une partie compacte et connexe de l'espace~$B$. Tout faisceau de $\Os_{M}$-modules de type fini satisfait le th\'eor\`eme~A.
\end{thm}
\begin{proof}
Soit~$\Fs$ un faisceau de $\Os_{M}$-modules de type fini. Soit~$b$ un point de~$M$. D'apr\`es le lemme \ref{thAlocal}, le faisceau~$\Fs$ v\'erifie le th\'eor\`eme~A sur un voisinage du point~$b$. Par compacit\'e de~$M$, il existe un entier~$p$ et des parties compactes et connexes $V_{0},\ldots,V_{p}$ de~$M$ recouvrant~$M$ telles que, quel que soit $i\in\cn{0}{p}$, le faisceau~$\Fs$ v\'erifie le th\'eor\`eme~A sur~$V_{i}$. Nous pouvons, en outre, supposer que, quel que soit $j\in\cn{0}{p-1}$, les compacts~$W_{j}=\bigcup_{0\le i\le j} V_{i}$ et~$V_{j+1}$ s'intersectent en un ensemble r\'eduit \`a un point de la forme~$a_{\sigma}^u$, avec~$\sigma\in\Sigma$ et $u\in\of{]}{0,l(\sigma)}{[}$. On montre alors, par r\'ecurrence et en utilisant \`a chaque \'etape la proposition \ref{CousinRungebase} et le corollaire \ref{corA}, que, quel que soit $j\in\cn{0}{p}$, le faisceau~$\Fs$ v\'erifie le th\'eor\`eme~A sur~$W_{j}$. On obtient le r\'esultat attendu en consid\'erant le cas~$j=p$. 
\end{proof}

\begin{thm}\label{Bbase}\index{Theoreme B@Th\'eor\`eme B!pour un compact de MA@pour un compact de $\Ms(A)$}
Soit~$M$ une partie compacte et connexe de l'espace~$B$. Tout faisceau de $\Os_{M}$-modules coh\'erent satisfait le th\'eor\`eme~B.
\end{thm}
\begin{proof}
Soit~$\Fs$ un faisceau de $\Os_{M}$-modules coh\'erent. Soit 
$$0\to \Fs \xrightarrow[]{d} \Is_{0} \xrightarrow[]{d} \Is_{1} \xrightarrow[]{d} \cdots$$
une r\'esolution flasque du faisceau~$\Fs$. Soient $q\in\N^*$ et~$\gamma$ un cocycle de degr\'e~$q$ sur~$M$. Soit~$b$ un point de~$M$. Par d\'efinition, le cocycle~$\gamma$ est un cobord au voisinage du point~$b$. En raisonnant comme dans la preuve qui pr\'ec\`ede et en utilisant la proposition \ref{propB}, dont la premi\`ere hypoth\`ese est v\'erifi\'ee d'apr\`es la proposition \ref{CousinRungebase}, au lieu du corollaire \ref{corA}, on montre que le cocycle~$\gamma$ est un cobord sur le compact~$M$. Puisque ce r\'esultat vaut pour tout cocycle, nous avons finalement montr\'e que le faisceau~$\Fs$ v\'erifie le th\'eor\`eme~B.
\end{proof}

\begin{cor}\index{Espace de Stein!compact de MA@compact de $\Ms(A)$}
Toute partie compacte et connexe de l'espace~$B$ est un espace de Stein.
\end{cor}

\section{Parties compactes des fibres}\label{pcdf}

Appliquons, \`a pr\'esent, les r\'esultats obtenus dans le cas des parties compactes des fibres de la droite analytique~$X$. Nous commencerons par d\'emontrer l'existence de syst\`emes de Cousin-Runge.

\begin{lem}\label{Cousinum}
Soient $V$ une partie compacte de~$B$ et $u,v,w$ trois nombres r\'eels v\'erifiant $0<u\le v\le w$. Pour tout \'el\'ement~$f$ de l'anneau $\Bs(V)\of{\la}{u \le|T|\le v}{\ra}$, il existe des \'el\'ements~$f^-$ de $\Bs(V)\of{\la}{|T|\le v}{\ra}$ et~$f^+$ de $\Bs(V)\of{\la}{u\le |T|\ge w}{\ra}$ tels que
\begin{enumerate}[\it i)]
\item $f = f^- + f^+$ dans $\Bs(V)\of{\la}{u \le|T|\le v}{\ra}$ ;
\item $\|f^-\|_{V,v} \le \|f\|_{V,u,v}$ ;
\item $\|f^+\|_{V,u,w} \le \|f\|_{V,u,v}$.
\end{enumerate} 
\end{lem}
\begin{proof}
Il existe une suite $(a_{k})_{k\in\Z}$ d'\'el\'ements de~$\Bs(V)$ telle que
$$f  = \sum_{k\in\Z} a_{k}\, T^k \in \Bs(V)\of{\la}{u \le|T|\le v}{\ra}.$$
Posons
$$f^- = \sum_{k\ge 0} a_{k}\, T^k \in \Bs(V)\of{\la}{|T|\le v}{\ra}$$
et
$$f^+ = \sum_{k\le -1} a_{k}\, T^k \in \Bs(V)\of{\la}{u \le|T|\le w}{\ra}.$$
Ces \'el\'ements v\'erifient l'\'egalit\'e 
$$f=f^-+f^+ \textrm{ dans } \Bs(V)\of{\la}{u \le|T|\le v}{\ra}.$$

Int\'eressons-nous, \`a pr\'esent, aux normes de ces s\'eries. Remarquons que
$${\renewcommand{\arraystretch}{1.3}\begin{array}{rcl}
\disp \left\|\sum_{k\in\Z} a_{k}\, T^k\right\|_{V,u,v} &=& \disp\sum_{k\in \Z} \|a_{k}\|_{V}\, \max(u^k,v^k)\\
&=& \disp\sum_{k\le -1} \|a_{k}\|_{V}\, u^k + \sum_{k\ge 0} \|a_{k}\|_{V}\, v^k.
\end{array}}$$
Nous avons
$$\|f^-\|_{V,v} = \sum_{k\ge 0} \|a_{k}\|_{V}\, v^k \le \|f\|_{V,u,v}$$
et
$$\|f^+\|_{V,u,w} = \sum_{k\le -1} \|a_{k}\|_{V}\, u^k \le \|f\|_{V,u,v}.$$
\end{proof}

\begin{lem}
Soient $V$ une partie compacte de~$B$ et $u,v,w$ trois nombres r\'eels v\'erifiant $0<u\le v\le w$. Soient~$p$ et~$q$ deux entiers et $s_{1},\ldots,s_{p},t_{1},\ldots,t_{q}$ des \'el\'ements de $\Bs(V)\of{\la}{u \le|T|\le v}{\ra}$. Soit $\delta\in\R_{+}^*$. Alors, il existe un \'el\'ement inversible~$f$ de $\Bs(V)\of{\la}{u \le|T|\le w}{\ra}$, des \'el\'ements $s'_{1},\ldots,s'_{p}$ de~$\Bs(V)\of{\la}{u \le|T|\le w}{\ra}$ et $t'_{1},\ldots,t'_{q}$ de~$\Bs(V)\of{\la}{|T|\le v}{\ra}$ tels que, quel que soient $i\in\cn{1}{p}$ et $j\in\cn{1}{q}$, on ait
\begin{enumerate}[\it i)]
\item $\|f^{-1}s_{i}-s'_{i}\|_{V,u,v}\, \|ft_{j}\|_{V,u,v} \le \delta$ ;
\item $\|f^{-1}s_{i}\|_{V,u,v}\, \|ft_{j}-t'_{j}\|_{V,u,v} \le \delta$.
\end{enumerate}
\end{lem}
\begin{proof}
Pour $i\in\cn{1}{p}$ et $j\in\cn{1}{q}$, notons
$$s_{i} = \sum_{k\in\Z} a_{k}^{(i)}\, T^k \in \Bs(V)\of{\la}{u \le|T|\le v}{\ra}$$
et
$$t_{j} = \sum_{k\in\Z} b_{k}^{(j)}\, T^k \in \Bs(V)\of{\la}{u \le|T|\le v}{\ra}.$$
Soit~$M>0$ tel que
$$\max_{1\le i\le p}(\|s_{i}\|_{V,u,v})\le M \textrm{ et } \max_{1\le j\le q} (\|t_{j}\|_{V,u,v}) \le M.$$
Il existe $k_{0}\le 0$ tel que, quels que soit $j\in\cn{1}{q}$, on ait
$$\sum_{k\le k_{0}-1} \|b_{k}^{(j)}\|_{V}\, u^k \le \frac{\delta}{M}.$$
Posons
$$f=T^{-k_{0}} \in\Bs(V)\of{\la}{u \le|T|\le w}{\ra}.$$
C'est un \'el\'ement inversible de~$\Bs(V)\of{\la}{u \le|T|\le w}{\ra}$. Pour $j\in\cn{1}{q}$, posons
$$t'_{j}=f\, \sum_{k\ge k_{0}} b_{k}^{(j)}\, T^k = \sum_{k\ge 0} b_{k+k_{0}}^{(j)}\, T^k \in\Bs(V)\of{\la}{|T|\le v}{\ra}.$$
Quels que soient $i\in\cn{1}{p}$ et $j\in\cn{1}{q}$, nous avons alors
$${\renewcommand{\arraystretch}{1.3}\begin{array}{rcl}
\disp \|f^{-1}s_{i}\|_{V,u,v}\, \|ft_{j}-t'_{j}\|_{V,u,v} & \le & \disp \|T^{k_{0}}\|_{V,u,v}\, \|s_{i}\|_{V,u,v}\, \left\|f\, \sum_{k\le k_{0}-1} b_{k}^{(j)}\, T^k\right\|_{V,u,v}\\
& \le & \disp u^{k_{0}}\, M\, \left\|\sum_{k\le k_{0}-1} b_{k}^{(j)}\, T^{k-k_{0}}\right\|_{V,u,v}\\
& \le & \disp u^{k_{0}}\, M\, \sum_{k\le k_{0}-1} \|b_{k}^{(j)}\|_{V}\, u^{k-k_{0}}\\
& \le & \delta.
\end{array}}$$

Soit $i\in\cn{1}{p}$. Il existe un \'el\'ement~$s''_{i}$ de~$\Bs(V)[T,T^{-1}]$ tel que
$$\|s_{i}-s''_{i}\|_{V,u,v}\le \frac{\delta}{M}\, \left(\frac{v}{u}\right)^{k_{0}}.$$
Posons
$$s'_{i}=f^{-1} s''_{i} = T^{k_{0}} s''_{i} \in \Bs(V)\of{\la}{u \le|T|\le w}{\ra}.$$
Quel que soit $j\in\cn{1}{q}$, nous avons alors
$$\|f^{-1}s_{i} - s'_{i}\|_{V,u,v}\, \|ft_{j}\|_{V,u,v}\le u^{k_{0}}\,  \frac{\delta}{M}\, \left(\frac{v}{u}\right)^{k_{0}}\, v^{-k_{0}}\, M \le \delta.$$
\end{proof}

\begin{prop}\label{propCousinRunge}\index{Systeme@Syst\`eme!de Cousin-Runge!dans les fibres de A1@dans les fibres de $\AA$}
Soit~$b$ un point de~$B$. Soit~$r$ un \'el\'ement de $\R_{+}^*\setminus\sqrt{|\Hs(b)^*|}$. Notons~$x$ le point~$\eta_{r}$ de la fibre~$X_{b}$. Soient~$K^-$ et~$K^+$ deux parties compactes et connexes de la fibre~$X_{b}$ dont l'intersection est \'egale au singleton~$\{x\}$. Alors, il existe un syst\`eme de Cousin-Runge associ\'e au couple $(K^-,K^+)$.
\end{prop}
\begin{proof}
Il existe un nombre r\'eel~$w$ tel que la partie compacte~$K^-\cup K^+$ soit contenue dans le disque ouvert de centre~$0$ et de rayon~$w$ de la fibre~$X_{b}$. Quitte \`a \'echanger les compacts~$K^-$ et~$K^+$, nous pouvons supposer que
$$K^-\subset \left\{y\in X_{b}\, \big|\, |T(y)|\le r\right\}.$$

Soit~$(V_{n})_{n\in\N}$ une suite d\'ecroissante de voisinages compacts du point~$b$ dans~$B$ qui compose un syst\`eme fondamental de voisinages de ce point. On d\'eduit facilement l'existence d'une telle suite de la description explicite de la topologie de l'espace~$B$ pr\'esent\'ee au num\'ero \ref{descriptionMA}. Soit~$(u_{n})_{n\in\N}$ une suite croissante et de limite~$r$ d'\'el\'ements de $\of{]}{0,r}{[}$. Soit~$(v_{n})_{n\in\N}$ une suite d\'ecroissante et de limite~$r$ d'\'el\'ements de $\of{]}{r,w}{[}$. Pour tout \'el\'ement~$n$ de~$\N$, nous posons
$$(\Bs^{-}_{n},\|.\|_{n}^-) = (\Bs(V_{n})\of{\la}{|T|\le v_{n}}{\ra},\|.\|_{V_{n},v_{n}}),$$
$$(\Bs^{+}_{n},\|.\|_{n}^+) = (\Bs(V_{n})\of{\la}{u_{n}\le |T| \le w}{\ra},\|.\|_{V_{n},u_{n},w})$$
et
$$(\Cs,\|.\|_{n}) = (\Bs(V_{n})\of{\la}{u_{n}\le |T| \le v_{n}}{\ra},\|.\|_{V_{n},u_{n},v_{n}}).$$
Quels que soient~$n\in\N$ et~$m\in\N$, on d\'efinit de mani\`ere \'evidente des morphismes $\psi^-_{n}$, $\psi^+_{n}$, $\rho^-_{n}$, $\rho^+_{n}$ et~$\rho_{n}$ comme dans la d\'efinition des syst\`emes de Banach. Le fait que les trois derniers soient born\'es d\'ecoule de la proposition \ref{spectrecouronne}. L'ensemble de ces donn\'ees forme un syst\`eme de Banach associ\'e au couple $(K^-,K^+)$. Les trois premi\`eres propri\'et\'es sont \'evidentes et la derni\`ere d\'ecoule du th\'eor\`eme \ref{anneaulocal}. Les deux lemmes qui pr\'ec\`edent assurent que ce syst\`eme est un syst\`eme de Cousin-Runge.
\end{proof}





\begin{cor}
Soit~$b$ un point de l'espace~$B$. Soit~$P_{b}(T)$ un polyn\^ome \`a coefficients dans~$\Hs(b)$. Soit~$r$ un \'el\'ement de $\R_{+}^*\setminus\sqrt{|\Hs(b)^*|}$. Posons
$$L_{0} = \{z\in X_{b}\, |\,  |P_{b}(T)(z)| = r\}.$$
Soient~$s$ et~$t$ deux \'el\'ements de~$\R_{+}$ tels que $s \le r \le t$. Consid\'erons les parties compactes de~$X$ d\'efinies par
$$K_{0}^-=\{z\in X_{b}\, |\, s\le |P_{b}(T)(z)| \le r\}$$
et 
$$K_{0}^+=\{z\in X_{b}\, |\, r\le |P_{b}(T)(z)| \le t\}.$$
Notons~$M_{0}$ leur r\'eunion. Soit~$\Fs$ un faisceau de $\Os_{M_{0}}$-modules de type fini qui satisfait le th\'eor\`eme~A sur les compacts~$K_{0}^-$ et~$K_{0}^+$. Alors il le satisfait encore sur~$M_{0}$. 
\end{cor}
\begin{proof}
D'apr\`es le lemme \ref{translation}, il existe un voisinage ouvert~$U$ du point~$b$ dans~$B$ et un polyn\^ome~$P(T)$ \`a coefficients dans~$\Os(U)$ dont l'image dans $\Hs(b)[T]$ est~$P_{b}(T)$. Comme expliqu\'e au num\'ero \ref{endodroite}, le morphime
$$\Os(U)[T] \to \Os(U)[T,S]/(P(S)-T) \xrightarrow[]{\sim} \Os(U)[S]$$
induit un morphisme
$$\varphi : Z = X_{U} \to X_{U} = Y.$$
C'est un morphisme topologique fini, d'apr\`es la proposition \ref{phifini2}. Posons
$$K^-=\{z\in X_{b}\, |\, s\le |T(z)| \le r\}$$
et 
$$K^+=\{z\in X_{b}\, |\, r\le |T(z)| \le t\}.$$
Ces deux compacts ont pour intersection l'ensemble r\'eduit au point~$\eta_{r}$ de la fibre~$X_{b}$, point que nous noterons~$y$. Notons~$M$ leur r\'eunion. Un calcul direct montre que, pour tous nombres r\'eels~$u$ et~$v$, nous avons
$$\varphi^{-1}\left(\{z\in X_{b}\, |\, u\le |T(z)| \le v\}\right) = \{z\in X_{b}\, |\, u\le |P_{b}(T)(z)| \le v\}.$$
En particulier, nous avons
$$\varphi^{-1}(K^-) = K_{0}^-,\ \varphi^{-1}(K^+)=K_{0}^+,\  \varphi^{-1}(y) = L_{0} \textrm{ et } \varphi^{-1}(M)=M_{0}.$$
On en d\'eduit que la partie compacte~$L_{0}$ est finie.

D'apr\`es le lemme \ref{thAcompact}, il existe un entier~$p$ et des \'el\'ements $t_{1}^-,\ldots,t_{p}^-$ de~$\Fs(K_{0}^-)$ dont les images engendrent le $\Os_{Z,z}$-module~$\Fs_{z}$, pour tout \'el\'ement~$z$ de~$K_{0}^-$. De m\^eme, il existe un entier~$q$ et des \'el\'ements $t_{1}^+,\ldots,t_{q}^+$ de~$\Fs(K_{0}^+)$ dont les images engendrent le $\Os_{Z,z}$-module~$\Fs_{z}$, pour tout \'el\'ement~$z$ de~$K_{0}^+$. Le corollaire \ref{XIG} et la proposition \ref{XS} nous permettent d'appliquer le th\'eor\`eme \ref{thfini}. Il assure que la famille $(1,S,\ldots,S^{d-1})$, o\`u~$d$ d\'esigne le degr\'e du polyn\^ome~$P$, engendre le $\Os_{Y,y}$-module~$(\varphi_{*}\Os_{Z})_{y}$. Quitte \`a remplacer la famille $(t_{i}^-)_{1\le i\le p}$ de~$\varphi_{*}\Fs(K^-)$ par la famille $(S^k\, t_{i}^-)_{0\le k\le d-1,1\le i\le p}$, nous pouvons supposer que le sous-$\Os_{Y,y}$-module de~$(\varphi_{*}\Fs)_{y}$ qu'elle engendre est identique au sous-$(\varphi_{*}\Os_{Z})_{y}$-module de~$(\varphi_{*}\Fs)_{y}$ qu'elle engendre. D'apr\`es le th\'eor\`eme \ref{finifibres}, les morphismes naturels
$$(\varphi_{*}\Os_{Z})_{y}\to \prod_{z\in L_{0}}\Os_{Z,z} \textrm{ et } (\varphi_{*}\Fs)_{y} \to \prod_{z\in L_{0}} \Fs_{z}$$ 
sont des isomorphismes. Pour tout \'el\'ement~$z$ de~$L_{0}$, la famille $(t_{1}^-,\ldots,t_{p}^-)$ engendre le $\Os_{Z,z}$-module~$\Fs_{z}$. Par cons\'equent, elle engendre \'egalement le $\Os_{Y,y}$-module~$(\varphi_{*}\Fs)_{y}$. De m\^eme, quitte \`a remplacer la famille $(t_{i}^+)_{1\le i\le q}$ de~$\varphi_{*}\Fs(K^+)$ par la famille $(S^k\, t_{i}^+)_{0\le k\le d-1,1\le i\le q}$, nous pouvons supposer qu'elle engendre le m\^eme $\Os_{Y,y}$-module~$(\varphi_{*}\Fs)_{y}$.

D'apr\`es le th\'eor\`eme \ref{attachefaisceau} et la proposition \ref{propCousinRunge}, il existe alors des \'el\'ements $s_{1}^-, \ldots,s_{p}^-,s_{1}^+,\ldots,s_{q}^+$ de~$\Fs(M)$, $a^-$ de~$GL_{p}(\Os(K^-))$ et~$a^+$ de~$GL_{q}(\Os(K^+))$ tels que
$$\begin{pmatrix} s_{1}^-\\ \vdots\\ s_{p}^-\end{pmatrix} = a^- \begin{pmatrix} t_{1}^-\\ \vdots\\ t_{p}^-\end{pmatrix} \textrm{ dans } \Fs(K^-)^p$$
et 
$$\begin{pmatrix} s_{1}^+\\ \vdots\\ s_{q}^+\end{pmatrix} = a^+ \begin{pmatrix} t_{1}^+\\ \vdots\\ t_{q}^+\end{pmatrix} \textrm{ dans } \Fs(K^+)^q.$$

Les matrices~$a^-$ et~$a^+$ induisent des \'el\'ements~$a_{0}^-$ et~$a_{0}^+$ de~$GL_{p}(\Os(K_{0}^-))$ et~$GL(\Os(K_{0}^+))$ tels que
$$\begin{pmatrix} s_{1}^-\\ \vdots\\ s_{p}^-\end{pmatrix} = a_{0}^- \begin{pmatrix} t_{1}^-\\ \vdots\\ t_{p}^-\end{pmatrix} \textrm{ dans } \Fs(K_{0}^-)^p$$
et 
$$\begin{pmatrix} s_{1}^+\\ \vdots\\ s_{q}^+\end{pmatrix} = a_{0}^+ \begin{pmatrix} t_{1}^+\\ \vdots\\ t_{q}^+\end{pmatrix} \textrm{ dans } \Fs(K_{0}^+)^q.$$
Les sections $s_{1}^-,\ldots,s_{p}^-,s_{1}^+,\ldots,s_{q}^+$ de~$\Fs(M_{0})$ engendrent alors le $\Os_{Z,z}$-module~$\Fs_{z}$ en tout point~$z$ de~$M_{0}$. On en d\'eduit le r\'esultat annonc\'e.
\end{proof}

Int\'eressons-nous, maintenant, plus sp\'ecifiquement au cas des fibres centrale et extr\^emes.

\begin{thm}\label{thAfibretriv}\index{Theoreme A@Th\'eor\`eme A!sur une fibre extr\^eme de A1@sur une fibre extr\^eme de $\AA$}\index{Theoreme A@Th\'eor\`eme A!sur la fibre centrale de A1@sur la fibre centrale de $\AA$}
Soit~$b$ un point central ou extr\^eme de l'espace~$B$. Soit~$M$ une partie compacte et connexe de la fibre~$X_{b}$. Tout faisceau de $\Os_{M}$-modules de type fini satisfait le th\'eor\`eme~A.
\end{thm}
\begin{proof}
Soit~$\Fs$ un faisceau de $\Os_{M}$-modules de type fini. D'apr\`es le lemme \ref{thAlocal}, tout point de~$M$ poss\`ede un voisinage sur lequel le faisceau~$\Fs$ v\'erifie le th\'eor\`eme~A. Par compacit\'e de~$M$, il existe un entier~$p$ et des parties compactes et connexes $M_{0},\ldots,M_{p}$ de~$M$ recouvrant~$M$ telles que, quel que soit $i\in\cn{0}{p}$, le faisceau~$\Fs$ v\'erifie le th\'eor\`eme~A sur~$M_{i}$. Nous pouvons, en outre, supposer que, quel que soit $j\in\cn{0}{p-1}$, les compacts~$N_{j}=\bigcup_{0\le i\le j} M_{i}$ et~$M_{j+1}$ s'intersectent en un ensemble r\'eduit \`a un point de type~$3$. Quel que soit $j\in\cn{0}{p-1}$, il existe alors un polyn\^ome irr\'eductible~$P_{b}(T)$ \`a coefficients dans~$\Hs(b)$, un \'el\'ement~$r$ de $\R_{+}^*\setminus\{1\}$ et des \'el\'ements~$s$ et~$t$ de~$\R_{+}$ v\'erifiant $s \le r \le t$ tels que l'on ait soit
$$N_{j} \subset \{z\in X_{b}\, |\, s\le |P_{b}(T)(z)| \le r\} \textrm{ et } M_{j+1} \subset \{z\in X_{b}\, |\, r\le |P_{b}(T)(z)| \le t\},$$
soit 
$$M_{j+1} \subset \{z\in X_{b}\, |\, s\le |P_{b}(T)(z)| \le r\} \textrm{ et } N_{j} \subset \{z\in X_{b}\, |\, r\le |P_{b}(T)(z)| \le t\}.$$
On montre alors, par r\'ecurrence et en utilisant le corollaire pr\'ec\'edent \`a chaque \'etape, que, quel que soit $j\in\cn{0}{p}$, le faisceau~$\Fs$ v\'erifie le th\'eor\`eme~A sur~$N_{j}$. On obtient le r\'esultat attendu en consid\'erant le cas~$j=p$. 
\end{proof}

\begin{rem}\label{topodroiteum}
Nous pouvons en fait d\'emontrer le r\'esultat pr\'ec\'edent pour tous les points de~$B_{\textrm{um}}$. Il suffit de savoir \'ecrire tout compact~$M$ de~$X_{b}$ comme r\'eunion de compacts $M_{0},\ldots,M_{p}$, pour un certain entier~$p$, v\'erifiant les m\^emes propri\'et\'es que ceux de la preuve du th\'eor\`eme : pour tout \'el\'ement~$j$ de~$\cn{0}{p-1}$, il existe un polyn\^ome~$P_{b}(T)$ \`a coefficients dans~$\Hs(b)$, un \'el\'ement~$r$ de $\R_{+}^*\setminus\sqrt{|\Hs(b)^*|}$ et des \'el\'ements~$s$ et~$t$ de~$\R_{+}$ v\'erifiant $s \le r \le t$ tels que l'on ait 
$$N_{j} = \bigcup_{1\le i\le j} M_{i} = \left\{z\in X_{b}\, \big|\, |P_{b}(T)(z)| = r\right\}$$
et soit
$$N_{j} \subset \{z\in X_{b}\, |\, s\le |P_{b}(T)(z)| \le r\} \textrm{ et } M_{j+1} \subset \{z\in X_{b}\, |\, r\le |P_{b}(T)(z)| \le t\},$$
soit 
$$M_{j+1} \subset \{z\in X_{b}\, |\, s\le |P_{b}(T)(z)| \le r\} \textrm{ et } N_{j} \subset \{z\in X_{b}\, |\, r\le |P_{b}(T)(z)| \le t\}.$$

On peut d\'emontrer que, pour tout \'el\'ement~$r$ de $\R_{+}^*\setminus\sqrt{|\Hs(b)^*|}$ et tout polyn\^ome irr\'eductible~$P_{b}$ \`a coefficients dans~$\Hs(b)$, l'ensemble
$$\left\{z\in X_{b}\, \big|\, |P_{b}(T)(z)| = r\right\}$$
est r\'eduit \`a un point. Le r\'esultat concernant le d\'ecoupage des compacts s'obtient alors en utilisant le fait que les points du type pr\'ec\'edent sont denses et la structure d'arbre de l'espace~$X_{b}$.
\end{rem}

\bigskip

\begin{lem}\label{Cousintype3}
Soit~$b$ un point de l'espace~$B$. Soit~$P_{b}(T)$ un polyn\^ome \`a coefficients dans~$\Hs(b)$. Soit $r$ un \'el\'ement de $\R_{+}^*\setminus\sqrt{|\Hs(b)^*|}$. Posons
$$L_{0} = \left\{z\in X_{b}\, \big|\, |P_{b}(T)(z)| = r \right\}.$$
Soit $t\ge r$. Consid\'erons les parties compactes de~$X$ d\'efinies par
$$K_{0}^-=\left\{z\in X_{b}\, \big|\, |P_{b}(T)(z)| \le r \right\}$$
et 
$$K_{0}^+=\left\{z\in X_{b}\, \big|\, r\le |P_{b}(T)(z)| \le t\right\}.$$
Leur intersection est le compact~$L_{0}$. Pour tout \'el\'ement~$f$ de~$\Os(L_{0})$, il existe un \'el\'ement~$f^-$ de~$\Os(K_{0}^-)$ et un \'el\'ement~$f^+$ de~$\Os(K_{0}^+)$ qui v\'erifient l'\'egalit\'e
$$f = f^- + f^+ \textrm{ dans } \Os(L_{0}).$$
\end{lem}
\begin{proof}
Commen\c{c}ons par le m\^eme raisonnement que dans le corollaire qui pr\'ec\`ede. D'apr\`es le lemme \ref{translation}, il existe un voisinage ouvert~$U$ du point~$b$ dans~$B$ et un polyn\^ome~$P(T)$ \`a coefficients dans~$\Os(U)$ dont l'image dans $\Hs(b)[T]$ est~$P_{b}(T)$. Comme expliqu\'e au num\'ero \ref{endodroite}, le morphisme
$$\Os(U)[T] \to \Os(U)[T,S]/(P(S)-T) \xrightarrow[]{\sim} \Os(U)[S]$$
induit un morphisme
$$\varphi : Z = X_{U} \to X_{U} = Y.$$
Posons
$$K^-=\left\{y\in X_{b}\, \big|\, |T(y)| \le r\right\}$$
et 
$$K^+=\left\{y\in X_{b}\, \big|\, r\le |T(y)| \le t\right\}.$$
Ces deux compacts ont pour intersection l'ensemble r\'eduit au point~$\eta_{r}$ de la fibre~$X_{b}$, point que nous noterons~$x$. D'apr\`es le th\'eor\`eme \ref{sectionslemniscates}, les morphismes naturels
$$\Os(K^-)[S]/(P(S)-T) \to \Os(K_{0}^-),$$ 
$$\Os(K^+)[S]/(P(S)-T) \to \Os(K_{0}^+)$$
et 
$$\Os_{Y,x}[S]/(P(S)-T) \to \Os(L_{0})$$
sont des isomorphismes. Par cons\'equent, il suffit de d\'emontrer le r\'esultat pour le point~$x$ et les parties compactes~$K^-$ et~$K^+$. Le r\'esultat d\'ecoule alors de l'existence d'un syst\`eme de Cousin associ\'e au couple $(K^-,K^+)$ (\emph{cf.}~proposition \ref{propCousinRunge}).
\end{proof}

\begin{thm}\label{thBfibretriv}\index{Theoreme B@Th\'eor\`eme B!sur une fibre extr\^eme de A1@sur une fibre extr\^eme de $\AA$}\index{Theoreme B@Th\'eor\`eme B!sur la fibre centrale de A1@sur la fibre centrale de $\AA$}
Soit~$b$ un point central ou extr\^eme de l'espace~$B$. Soit~$M$ une partie compacte et connexe de la fibre~$X_{b}$. Tout faisceau de $\Os_{M}$-modules de type fini satisfait le th\'eor\`eme~B.
\end{thm}
\begin{proof}
Soit~$\Fs$ un faisceau de $\Os_{M}$-modules coh\'erent. Soit 
$$0\to \Fs \xrightarrow[]{d} \Is_{0} \xrightarrow[]{d} \Is_{1} \xrightarrow[]{d} \cdots$$
une r\'esolution flasque du faisceau~$\Fs$. Soient $q\in\N^*$ et~$\gamma$ un cocycle de degr\'e~$q$ sur~$M$. Par d\'efinition, tout point de~$M$ poss\`ede un voisinage sur lequel le cocycle~$\gamma$ est un cobord. Par compacit\'e de~$M$, il existe un entier~$p$ et des parties compactes et connexes $M_{0},\ldots,M_{p}$ de~$M$ recouvrant~$M$ telles que, quel que soit $i\in\cn{0}{p}$, le cocycle~$\gamma$ soit un cobord sur~$M_{i}$. Nous pouvons, en outre, supposer que, quel que soit $j\in\cn{0}{p-1}$, les compacts~$N_{j}=\bigcup_{0\le i\le j} M_{i}$ et~$M_{j+1}$ s'intersectent en un ensemble r\'eduit \`a un point de type~$3$. 

Montrons, par r\'ecurrence, que, quel que soit $j\in\cn{0}{p}$, le cocycle~$\gamma$ est un cobord sur le compact~$N_{j}$. Le cas~$j=0$ est vrai par hypoth\`ese. Soit $j\in\cn{0}{p-1}$ et supposons que le cocycle~$\gamma$ est un cobord sur le compact~$N_{j}$. D'apr\`es le lemme pr\'ec\'edent, pour tout \'el\'ement~$f$ de~$\Os(N_{j}\cap M_{j+1})$, il existe un \'el\'ement~$f^-$ de~$\Os(N_{j})$ et un \'el\'ement~$f^+$ de~$\Os(M_{j+1})$ qui v\'erifient l'\'egalit\'e
$$f = f^- + f^+ \textrm{ dans } \Os(N_{j}\cap M_{j+1}).$$
En outre, puisque l'intersection $N_{j}\cap M_{j+1}$ est r\'eduite \`a un point, tout faisceau de $\Os_{N_{j}\cap M_{j+1}}$-modules v\'erifie le th\'eor\`eme~B. Finalement, tout faisceau de $\Os_{N_{j}\cup M_{j+1}}$-modules coh\'erent satisfait le th\'eor\`eme~A, d'apr\`es le th\'eor\`eme \ref{thAfibretriv}. La proposition \ref{propB} assure alors que le cocycle~$\gamma$ est un cobord sur le compact $N_{j} \cup M_{j+1} =N_{j+1}$, ce qu'il fallait d\'emontrer.

Nous avons en particulier prouv\'e que le cocycle~$\gamma$ est un cobord sur le compact $M = N_{p}$. Puisque ce r\'esultat vaut pour tout cocycle, nous avons finalement montr\'e que le faisceau~$\Fs$ v\'erifie le th\'eor\`eme~B.
\end{proof}

\begin{rem}
Nous pouvons en fait d\'emontrer le r\'esultat pr\'ec\'edent pour tous les points de~$B_{\textrm{um}}$ en proc\'edant de m\^eme et en utilisant le r\'esultat dont il est question \`a la remarque \ref{topodroiteum}.
\end{rem}

\bigskip

Le cas des parties compactes des fibres internes peut se traiter, comme toujours, en se ramenant au cas classique des espaces analytiques sur un corps valu\'e complet. Nous pouvons en d\'eduire des r\'esultats ind\'ependants de la fibre consid\'er\'ee, par exemple dans le cas des couronnes. 



\begin{thm}\label{thABfibre}\index{Espace de Stein!couronne compacte d'une fibre de A1@couronne compacte d'une fibre de $\AA$}\index{Disque!de Stein|see{Espace de Stein}}\index{Couronne!de Stein|see{Espace de Stein}}
Soient $b$ un point de l'espace~$B$ et~$r$ et~$s$ deux \'el\'ements de~$\R_{+}$ v\'erifiant l'in\'egalit\'e $r\le s$. Posons
$$C = \{y\in X_{b}\, |\, r\le |T(y)|\le s\}.$$ 
La couronne~$C$ est un sous-espace de Stein de la droite analytique~$X$.
\end{thm}
\begin{proof}
Si~$b$ est un point central ou extr\^eme de l'espace~$B$, le r\'esultat d\'ecoule des th\'eor\`emes \ref{thAfibretriv} et \ref{thBfibretriv}.

Supposons d\'esormais que le point~$b$ est un point interne de l'espace~$B$. Notons 
$$j_{b} : X_{b} \hookrightarrow X$$
le morphisme d'inclusion. Soit $\Fs$ un faisceau coh\'erent sur $C$. Le faisceau de $\Os_{X_{b}}$-modules $j_{b}^{-1}\Fs$ est encore un faisceau coh\'erent sur $C$. D'apr\`es la proposition~\ref{isointerne}, il nous suffit de montrer que le faisceau~$j_{b}^{-1}\Fs$ v\'erifie les th\'eor\`emes~A et~B.

Distinguons, \`a pr\'esent, deux cas. Si le point $b$ appartient \`a une branche ultram\'etrique, son corps r\'esiduel $\Hs(b)$ est muni d'une valeur absolue ultram\'etrique non triviale. D'apr\`es le th\'eor\`eme~2.4 de~\cite{AB}, la proposition~3.3.4 de~\cite{rouge} et le th\'eor\`eme \ref{BimpliqueA}, pour tous \'el\'ements~$r'$ et~$s'$ de~$\R_{+}$, la partie de l'espace analytique~$\E{1}{\Hs(b)}$ d\'efinie par
$$\{y\in X_{b}\, |\, r'< |T(y)|< s'\}$$
est un espace de Stein (dans notre sens). Or, d'apr\`es le lemme \ref{lemvoisdisquedroite}, l'ensemble des parties de la forme
$$\{y\in X_{b}\, |\, r'< |T(y)|< s'\},$$
o\`u~$r'$ et~$s'$ sont deux \'el\'ements de~$\R_{+}$ v\'erifiant $r'\prec r$ et $s'>s$, est un syst\`eme fondamental de voisinages de~$C$ dans l'espace analytique~$\E{1}{\Hs(b)}$. Le corollaire \ref{Steincompact} assure alors que la partie compacte~$C$ de l'espace analytique~$\E{1}{\Hs(b)}$ est de Stein.


Si le point $b$ appartient \`a une branche archim\'edienne, le faisceau coh\'erent~$j_{b}^{-1}\Fs$ v\'erifie encore les th\'eor\`emes~A et~B. En effet, les couronnes ferm\'ees de $\C$ sont des espaces de Stein.
\end{proof}

\begin{rem}
Lorsque le point~$b$ est un point interne d'une branche ultram\'etrique, nous avons essentiellement red\'emontr\'e un cas particulier de la proposition~3.1 de~\cite{GK}.
\end{rem}

\section{Couronnes compactes de la droite}\label{ccdld}

Dans ce paragraphe, nous d\'emontrons que certaines parties compactes de la droite analytique~$X$ sont de Stein. Comme pr\'ec\'edemment, nous commencerons par exhiber des syst\`emes de Cousin-Runge. 



\begin{lem}\label{Cousin}
Soient $\sigma\in\Sigma$, $u\in \of{]}{0,l(\sigma)}{[}$ et $s,t\in\R_{+}$ tels que $s\le t$. Posons 
$$K_{0}^- = \of{[}{a_{\sigma}^u,a_{\sigma}^{l(\sigma)}}{]},\ K_{0}^+=B\setminus\of{]}{a_{\sigma}^u,a_{\sigma}^{l(\sigma)}}{]} \textrm{ et } L_{0}=K_{0}^-\cap K_{0}^+=\{a_{\sigma}^u\}.$$

Soit~$D\in\R$ la constante dont le lemme \ref{Cousinbase} assure l'existence. Quels que soient $s,t\in{[}{0,+\infty}{[}$, avec $s \le t$, et quel que soit $f\in\Bs(L_{0})\of{\la}{s\le |T|\le t}{\ra}$, il existe $f^-\in\Bs(K_{0}^-)\of{\la}{s\le |T|\le t}{\ra}$ et $f^+\in\Bs(K_{0}^+)\of{\la}{s\le |T|\le t}{\ra}$ v\'erifiant les propri\'et\'es suivantes :
\begin{enumerate}[\it i)]
\item $f=f^--f^+$ dans $\Bs(L_{0})\of{\la}{s\le |T|\le t}{\ra}$ ;
\item $\|f^-\|_{K_{0}^-,s,t} \le D\,  \|f\|_{L_{0},s,t}$ ;
\item $\|f^+\|_{K_{0}^+,s,t} \le D\,\|f\|_{L_{0},s,t}$.
\end{enumerate}
\end{lem}
\begin{proof}
Soient $s,t\in{[}{0,+\infty}{[}$, avec $s \le t$, et $f\in\Bs(L_{0})\of{\la}{s\le |T|\le t}{\ra}$. Par d\'efinition, il existe une famille~$(a_{k})_{k\in\Z}$ de~$\Bs(L_{0})=\hat{K}_{\sigma}$ telle que l'on ait
$$ f = \sum_{k\in\Z} a_{k}\,T^k$$
et que les s\'eries
$$\sum_{k\ge 0} a_{k}\, t^k \textrm{ et } \sum_{k\le 0} a_{k}\, s^k$$ 
convergent. Soit~$k\in\Z$. D'apr\`es le lemme~\ref{Cousinbase}, il existe des \'el\'ements~$a_{k}^-$ de~$\Bs(K_{0}^-)$ et~$a_{k}^+$ de~$\Bs(K_{0}^+)$ v\'erifiant les propri\'et\'es suivantes :
\begin{enumerate}[\it i)]
\item $a_{k}=a_{k}^--a_{k}^+$ dans $\Bs(L_{0})$ ;
\item $\|a_{k}^-\|_{K_{0}^-} \le D\,  \|a_{k}\|_{L_{0}}$ ;
\item $\|a_{k}^+\|_{K_{0}^+} \le D\,\|a_{k}\|_{L_{0}}$.
\end{enumerate}
Posons
$$f^-= \sum\limits_{k\in\Z} a_{k}^-\, T^k$$
et  
$$f^+= \sum\limits_{k\in\Z} a_{k}^+\, T^k.$$
Ces s\'eries v\'erifient les conditions requises.
\end{proof}

Contrairement au pr\'ec\'edent, le r\'esultat d'approximation ne nous semble pas pouvoir se d\'eduire du r\'esultat similaire pour les parties compactes de la base (\emph{cf.} lemme \ref{appbase}).

\begin{lem}\label{app}
Soient $\sigma\in\Sigma$, $u\in \of{]}{0,l(\sigma)}{[}$ et $s,t\in\R_{+}$ tels que $s\le t$. Posons 
$$K_{0}^- = \of{[}{a_{\sigma}^u,a_{\sigma}^{l(\sigma)}}{]},\ K_{0}^+=B\setminus\of{]}{a_{\sigma}^u,a_{\sigma}^{l(\sigma)}}{]} \textrm{ et } L_{0}=K_{0}^-\cap K_{0}^+=\{a_{\sigma}^u\}.$$ 

Soient $p,q\in\N$ et $s_{1},\ldots,s_{p},t_{1},\ldots,t_{q}\in\Bs(L_{0})\of{\la}{s\le |T|\le t}{\ra}$. Soit $\delta\in\R_{+}^*$. Si~$\sigma$ appartient \`a~$\Sigma_{f}$, alors il existe un \'el\'ement inversible $f$ de $\Bs(K_{0}^+)\of{\la}{s\le |T|\le t}{\ra}$ et des \'el\'ements $s'_{1},\ldots,s'_{p}$ de $\Bs(K_{0}^+)\of{\la}{s\le |T|\le t}{\ra}$ et $t'_{1},\ldots,t'_{q}$ de $\Bs(K_{0}^-)\of{\la}{s\le |T|\le t}{\ra}$ tels que, quel que soient $i\in\cn{1}{p}$ et $j\in\cn{1}{q}$, on ait
\begin{enumerate}[\it i)]
\item $\|f^{-1}s_{i}-s'_{i}\|_{L_{0},s,t}\,\|ft_{j}\|_{L_{0},s,t} \le \delta$ ;
\item $\|f^{-1}s_{i}\|_{L_{0},s,t}\, \|ft_{j}-t'_{j}\|_{L_{0},s,t}\le \delta$.
\end{enumerate}

Si~$\sigma$ appartient \`a~$\Sigma_{\infty}$, alors il existe un \'el\'ement inversible~$g$ de $\Bs(K_{0}^-)\of{\la}{s\le |T|\le t}{\ra}$ et des \'el\'ements $s''_{1},\ldots,s''_{p}$ de $\Bs(K_{0}^+)\of{\la}{s\le |T|\le t}{\ra}$ et $t''_{1},\ldots,t''_{q}$ de $\Bs(K_{0}^-)\of{\la}{s\le |T|\le t}{\ra}$ tels que, quel que soient $i\in\cn{1}{p}$ et $j\in\cn{1}{q}$, on ait
\begin{enumerate}[\it i)]
\item $\|gs_{i}-s''_{i}\|_{L_{0},s,t}\, \|g^{-1}t_{j}\|_{L_{0},s,t}\le \delta$ ;
\item $\|gs_{i}\|_{L_{0},s,t}\,\|g^{-1}t_{i}-t''_{i}\|_{L_{0},s,t}\le \delta$.
\end{enumerate}
\end{lem}
\begin{proof}
Posons $M=\max\{\|s_{i}\|_{L_{0},s,t},\|t_{j}\|_{L_{0},s,t},\, 1\le i\le p, 1\le j\le q\}$. Soit $i\in\cn{1}{p}$. La fonction~$s_{i}$ appartient \`a $\Bs(L_{0})\of{\la}{s\le |T|\le t}{\ra}$. Par cons\'equent, il existe une famille~$(a_{k})_{k\in\Z}$ de~$\hat{K}_{\sigma}$ telle que
$$s_{i} = \sum_{k\in\Z} a_{k}\,T^k$$
et les s\'eries
$$\sum_{k\ge 0} |a_{k}|_{\sigma}^u\, t^k \textrm{ et } \sum_{k\le 0} |a_{k}|_{\sigma}^u\, s^k$$
convergent. Il existe $n_{i},n'_{i}\in\Z$ tel que 
$$\left\| s_{i} - \sum\limits_{k=n_{i}}^{n'_{i}} a_{k}\,T^k \right\|_{L_{0},s,t}\le \delta.$$
Il existe \'egalement $s^*_{i} \in K[T,T^{-1}]$ tel que 
$$\left\| \sum\limits_{k=n_{i}}^{n'_{i}} a_{k}\,T^k - s^*_{i} \right\|_{L_{0},s,t}\le \delta.$$

\bigskip

Distinguons deux cas. Supposons, tout d'abord, que $\sigma\in\Sigma_{f}$. D'apr\`es le lemme \ref{diviseur}, il existe $h\in A$ telle que 
$$\left\{\begin{array}{l}
|h|_{\sigma}<1\ ;\\
\forall \sigma'\in\Sigma_{f}\setminus\{\sigma\},\, |h|_{\sigma'}=1.
\end{array}\right.$$
Il existe $N\in\N$ tel que, quel que soit $i\in\cn{1}{p}$, on ait $h^N\, s^*_{i} \in \hat{A}_{\sigma}[T,T^{-1}]$. En particulier, la fonction~$h^N\, s^*_{i}$ d\'efinit un \'el\'ement de $\Bs(K_{0}^-)\of{\la}{s\le |T|\le t}{\ra}$. Posons 
$$f=h^{-N}\in K.$$ 
C'est un \'el\'ement inversible de $\Bs(K_{0}^+)\of{\la}{s\le |T|\le t}{\ra}$. En outre, nous avons 
$$\|f^{-1} s_{i} - f^{-1} s^*_{i}\|_{L_{0},s,t} \le \|f^{-1}\|_{L_{0},s,t}\, \|s_{i}-s^*_{i}\|_{L_{0},s,t} \le 2\delta\, |f^{-1}|_{\sigma}^u.$$

Soit $j\in\cn{1}{q}$. La fonction~$t_{j}$ appartient \`a~$\Bs(L_{0})\of{\la}{s\le |T|\le t}{\ra}$. Par cons\'equent, il existe une famille~$(b_{k})_{k\in\Z}$ de~$\hat{K}_{\sigma}$ telle que
$$ft_{j} = \sum_{k\in\Z} b_{k}\,T^k$$
et les s\'eries
$$\sum_{k\ge 0} |b_{k}|_{\sigma}^u\, t^k \textrm{ et } \sum_{k\le 0} |b_{k}|_{\sigma}^u\, s^k$$
convergent. Il existe $m_{j},m'_{j}\in\Z$ tel que 
$$\left\|f t_{j} - \sum\limits_{k=m_{j}}^{m'_{j}} b_{k}\,T^k \right\|_{L_{0},s,t}\le \delta.$$
Par le th\'eor\`eme d'approximation fort, quel que soit $\eps>0$, il existe des \'el\'ements \mbox{$c_{m_{j}},\ldots,c_{m'_{j}}$} de~$K$ tels que, quel que soit $k\in\cn{m_{j}}{m'_{j}}$, on ait 
\begin{enumerate}
\item $\forall \sigma'\in \Sigma_{f}\setminus\{\sigma\}$, $c_{k} \in \hat{A}_{\sigma'}$ ;
\item $|b_{k}-c_{k}|_{\sigma}^{u}\le \eps$.
\end{enumerate}
On en d\'eduit qu'il existe $t^*_{j} \in A[1/\sigma][T,T^{-1}]$ tel que 
$$\left\| \sum\limits_{k=m_{j}}^{m'_{j}} b_{k}\,T^k - t^*_{j} \right\|_{L_{0},s,t}\le \delta.$$
En particulier, la fonction~$t^*_{j}$ d\'efinit un \'el\'ement de~$\Bs(K_{0}^+)\of{\la}{s\le |T|\le t}{\ra}$ et, quel que soit~$i\in\cn{1}{p}$, nous avons 
$$\|f^{-1} s_{i}\|_{L_{0},s,t}\,  \|ft_{j} - t^*_{j}\|_{L_{0},s,t}\le |f^{-1}|_{\sigma}^u\, \|s_{i}\|_{L_{0},s,t}\, 2\delta \le 2 M \delta,$$
car $f^{-1}\in A$. Quel que soit $i\in\cn{1}{p}$, nous avons \'egalement
$$\|f^{-1} s_{i} - f^{-1} s^*_{i}\|_{L_{0},s,t}\, \|ft_{j}\|_{L_{0},s,t} \le 2\delta\, |f^{-1}|_{\sigma}^u\, |f|_{\sigma}^u\, \|t_{j}\|_{L_{0},s,t} \le 2M\delta.$$

\bigskip

Supposons, \`a pr\'esent, que $\sigma\in\Sigma_{\infty}$. Il existe $g\in A$ tel que, quel que soit \mbox{$i\in\cn{1}{n}$}, on ait $gs^*_{i} \in A[T,T^{-1}]$. Remarquons que l'image de~$g$ l'anneau $\Bs(K_{0}^-)\of{\la}{s\le |T|\le t}{\ra}$ est inversible et que~$gs^*_{i}$ d\'efinit un \'el\'ement de~$\Bs(K_{0}^+)\of{\la}{s\le |T|\le t}{\ra}$. En outre, quel que soit~$i\in\cn{1}{n}$, nous avons
$$\|gs_{i} - gs^*_{i}\|_{L_{0},s,t} \le \|g\|_{L_{0}}\, \|s_{i}-s^*_{i}\|_{L_{0},s,t} \le 2\,|g|_{\sigma}^u\, \delta.$$

Soit $j\in\cn{1}{q}$. Puisque la fonction~$g^{-1}t_{j}$ appartient \`a~$\Bs(L_{0})\of{\la}{s\le |T|\le t}{\ra}$, il existe une famille~$(b_{k})_{k\in\Z}$ de~$\hat{K}_{\sigma}$ telle que
$$g^{-1}t_{j} = \sum_{k\in\Z} b_{k}\,T^k$$
et les s\'eries
$$\sum_{k\ge 0} |b_{k}|_{\sigma}^{u}\, t^k \textrm{ et } \sum_{k\le 0} |b_{k}|_{\sigma}^{u}\, s^k$$
convergent. Il existe $m_{j},m'_{j}\in\Z$ tels que 
$$\left\| g^{-1}t_{j} - \sum\limits_{k=m_{j}}^{m'_{j}} b_{k}\,T^k \right\|_{L_{0},s,t}\le \frac{\delta}{2\,\|gs_{i}\|_{L_{0},s,t}}.$$
En approchant chacun des coefficients $b_{k}$, avec $k\in\cn{m_{j}}{m'_{j}}$, on montre qu'il existe \'egalement $t^*_{j} \in K[T,T^{-1}]$ tel que 
$$\left\| \sum\limits_{k=m_{j}}^{m'_{j}} b_{k}\,T^k - t^*_{j} \right\|_{L_{0},s,t}\le \frac{\delta}{2\,\|gs_{i}\|_{L_{0},s,t}}.$$
La fonction $t^*_{j}$ d\'efinit un \'el\'ement de~$\Bs(K_{0}^-)\of{\la}{s\le |T|\le t}{\ra}$ et, quel que soit \mbox{$i\in\cn{1}{p}$}, v\'erifie
$$\|gs_{i}\|_{L_{0},s,t}\, \|g^{-1}t_{j}-t^*_{j}\|_{L_{0},s,t}\le \|gs_{i}\|_{L_{0},s,t}\, 2\,  \frac{\delta}{2\,\|gs_{i}\|_{L_{0},s,t}} \le \delta.$$
Quel que soit $i\in\cn{1}{p}$, nous avons encore
$$\|gs_{i} - s^*_{i}\|_{L_{0},s,t}\, \|g^{-1}t_{j}\|_{L_{0},s,t} \le 2\, |g|_{\sigma}^u\, \delta\, |g^{-1}|_{\sigma}^u\, \|t_{j}\|_{L_{0},s,t} \le 2M\delta.$$
\end{proof}

\begin{prop}\label{CousinRungecouronne}\index{Systeme@Syst\`eme!de Cousin-Runge!pour les couronnes de A1@pour les couronnes de $\AA$}
Soient $\sigma\in\Sigma$, $u\in \of{]}{0,l(\sigma)}{[}$ et $s,t\in\R_{+}$ tels que $s\le t$. Posons 
$$K_{0}^- = \of{[}{a_{\sigma}^u,a_{\sigma}^{l(\sigma)}}{]},\ K_{0}^+=B\setminus\of{]}{a_{\sigma}^u,a_{\sigma}^{l(\sigma)}}{]},\ L_{0}=K_{0}^-\cap K_{0}^+=\{a_{\sigma}^u\}$$ 
et
$$L=\overline{C}_{L_{0}}(s,t) = \{x\in X\, | s\le |T(x)|\le t\}.$$
Soient~$K^-$ une partie compacte de $\overline{C}_{K_{0}^-}(s,t)$ et~$K^+$ une partie compacte de $\overline{C}_{K_{0}^+}(s,t)$ dont l'intersection est le compact~$L$. Il existe un syst\`eme de Cousin-Runge associ\'e au couple $(K^-,K^+)$.
\end{prop}
\begin{proof}

Soit~$(s_{n})_{n\in\N}$ une suite croissante et de limite~$s$ d'\'el\'ements de $\of{[}{0,s}{[}$. Soit~$(t_{n})_{n\in\N}$ une suite d\'ecroissante et de limite~$t$ d'\'el\'ements de $\of{]}{t,+\infty}{[}$. Pour tout \'el\'ement~$n$ de~$\N$, nous posons
$$(\Bs^{-}_{n},\|.\|_{n}^-) = (\Bs(K_{0}^-)\of{\la}{s_{n}\le |T|\le t_{n}}{\ra},\|.\|_{K_{0}^-,s_{n},t_{n}}),$$
$$(\Bs^{+}_{n},\|.\|_{n}^+) = (\Bs(K_{0}^+)\of{\la}{s_{n}\le |T| \le t_{n}}{\ra},\|.\|_{K_{0}^+,s_{n},t_{n}})$$
et
$$(\Cs,\|.\|_{n}) = (\Bs(L_{0})\of{\la}{s_{n}\le |T| \le t_{n}}{\ra},\|.\|_{L_{0},s_{n},t_{n}}).$$
Quels que soient~$n\in\N$ et~$m\in\N$, on d\'efinit de mani\`ere \'evidente des morphismes born\'es~$\psi^-_{n}$ et~$\psi^+_{n}$ comme dans la d\'efinition des syst\`emes de Banach. Le th\'eor\`eme \ref{isocouronne} permet de d\'efinir \'egalement des morphismes $\rho^-_{n}$, $\rho^+_{n}$ et~$\rho_{n}$. La proposition \ref{spectrecouronne} assure qu'ils sont born\'es. L'ensemble de ces donn\'ees forme un syst\`eme de Banach associ\'e au couple $(K^-,K^+)$. Les trois premi\`eres propri\'et\'es sont \'evidentes et la derni\`ere d\'ecoule de nouveau du th\'eor\`eme \ref{isocouronne}, joint \`a la proposition \ref{isoBO}. Les deux lemmes qui pr\'ec\`edent assurent que ce syst\`eme est un syst\`eme de Cousin-Runge.
\end{proof}

Nous allons d\'eduire de ces r\'esultats le fait que les couronnes compactes et connexes de la droite analytique~$X$ sont des espaces de Stein.

\begin{thm}\label{A}\index{Theoreme A@Th\'eor\`eme A!couronne compacte de A1@couronne compacte de $\AA$}
Soit~$V$ une partie compacte et connexe de l'espace~$B$. Soient~$s$ et~$t$ deux nombres r\'eels tels que $0\le s\le t$. Posons
$$M = \overline{C}_{V}(s,t) = \{x\in X_{V}\, |\, s\le |T(x)|\le t\}.$$
Tout faisceau de $\Os_{M}$-modules de type fini satisfait le th\'eor\`eme~A.
\end{thm}
\begin{proof}
Soit~$\Fs$ un faisceau de $\Os_{M}$-modules de type fini. Soit~$b$ un point de~$V$. D'apr\`es le th\'eor\`eme \ref{thABfibre}, le faisceau~$\Fs$ v\'erifie le th\'eor\`eme~A sur le compact $X_{b}\cap M$ et donc sur un voisinage de ce compact, d'apr\`es le lemme \ref{thAlocal}. En utilisant le lemme qui pr\'ec\`ede, on en d\'eduit qu'il existe un voisinage compact~$V_{b}$ du point~$b$ dans~$V$ tel que le faisceau~$\Fs$ v\'erifie le th\'eor\`eme~A sur le compact~$X_{V_{b}}\cap M$. Par compacit\'e de~$M$, il existe un entier~$p$ et des parties compactes et connexes $V_{0},\ldots,V_{p}$ de~$V$ recouvrant~$V$ telles que, quel que soit $i\in\cn{0}{p}$, le faisceau~$\Fs$ v\'erifie le th\'eor\`eme~A sur~$X_{V_{i}}\cap M$. Nous pouvons, en outre, supposer que, quel que soit $j\in\cn{0}{p-1}$, les compacts~$W_{j}=\bigcup_{0\le i\le j} V_{i}$ et~$V_{j+1}$ s'intersectent en un ensemble r\'eduit \`a un point de type~$3$. On montre alors, par r\'ecurrence et en utilisant \`a chaque \'etape la proposition \ref{CousinRungecouronne} et le corollaire \ref{corA}, que, quel que soit $j\in\cn{0}{p}$, le faisceau~$\Fs$ v\'erifie le th\'eor\`eme~A sur~$X_{W_{j}}\cap M$. On obtient le r\'esultat attendu en consid\'erant le cas~$j=p$. 
\end{proof}

\begin{thm}\label{B}\index{Theoreme B@Th\'eor\`eme B!couronne compacte de A1@couronne compacte de $\AA$}
Soit~$V$ une partie compacte et connexe de l'espace~$B$. Soient~$s$ et~$t$ deux nombres r\'eels tels que $0\le s\le t$. Posons
$$M = \overline{C}_{V}(s,t) = \{x\in X_{V}\, |\, s\le |T(x)|\le t\}.$$
Tout faisceau de $\Os_{M}$-modules coh\'erent satisfait le th\'eor\`eme~B.
\end{thm}
\begin{proof}
Soit~$\Fs$ un faisceau de $\Os_{M}$-modules coh\'erent. Soit 
$$0\to \Fs \xrightarrow[]{d} \Is_{0} \xrightarrow[]{d} \Is_{1} \xrightarrow[]{d} \cdots$$
une r\'esolution flasque du faisceau~$\Fs$. Soient $q\in\N^*$ et~$\gamma$ un cocycle de degr\'e~$q$ sur~$M$. Soit~$b$ un point de~$V$. D'apr\`es le th\'eor\`eme \ref{thABfibre}, le faisceau~$\Fs$ v\'erifie le th\'eor\`eme~B sur le compact $X_{b}\cap M$. Par cons\'equent, le cocycle~$\gamma$ est un cobord au voisinage du compact $X_{b}\cap M$. En utilisant le lemme \ref{lemvoisdisquedroite}, on en d\'eduit qu'il existe un voisinage compact~$V_{b}$ du point~$b$ dans~$V$ tel que le cocycle~$\gamma$ soit sur le compact~$X_{V_{b}}\cap M$. En raisonnant comme dans la preuve qui pr\'ec\`ede et en utilisant la proposition \ref{propB}, dont la premi\`ere hypoth\`ese est v\'erifi\'ee d'apr\`es la proposition \ref{CousinRungecouronne}, au lieu du corollaire \ref{corA}, on montre que le cocycle~$\gamma$ est un cobord sur le compact~$M$. Puisque ce r\'esultat vaut pour tout cocycle, nous avons finalement montr\'e que le faisceau~$\Fs$ v\'erifie le th\'eor\`eme~B.
\end{proof}

\begin{thm}\label{ABcompact}\index{Espace de Stein!couronne compacte de A1@couronne compacte de $\AA$}
Soit~$V$ une partie compacte et connexe de l'espace~$B$. Soient~$s$ et~$t$ deux nombres r\'eels tels que $0\le s\le t$. La couronne compacte
$$\overline{C}_{V}(s,t) = \{x\in X_{V}\, |\, s\le |T(x)|\le t\}$$
est un espace de Stein.
\end{thm}



\section{Lemniscates de la droite}\label{ldld}

Dans cette partie, nous allons montrer que les th\'eor\`emes~A et~B sont satisfaits pour les faisceaux coh\'erents d\'efinies sur les couronnes ouvertes de la droite analytique~$X$ et les lemniscates. Ici encore, nous nous inspirerons des techniques utilis\'ees en g\'eom\'etrie analytique complexe. Pour toute couronne ouverte~$C$, nous consid\`ererons une famille croissante de couronnes ferm\'ees dont la r\'eunion est \'egale \`a~$C$. Nous montrerons alors que cette famille forme une exhaustion de Stein (\emph{cf.} \cite{GR}, IV, \S 1, d\'efinition~6). Il nous restera alors \`a montrer que toute partie poss\'edant une exhaustion de Stein est de Stein.

La preuve que nous proposons ici suit de tr\`es pr\`es l'ouvrage \cite{GR} de H. Grauert et R. Remmert. Plus pr\'ecis\'ement, nous nous sommes inspir\'es de la partie IV,~\S 1 pour les d\'efinition et propri\'et\'es des exhaustions de Stein et de la partie IV,~\S 4, pour montrer que les familles croissantes de couronnes ferm\'ees consid\'er\'ees en satisfont les conditions. 

Nous traiterons finalement le cas des lemniscates en faisant appel au th\'eor\`eme \ref{finiStein} et aux r\'esultats sur les morphismes finis d\'emontr\'es au chapitre \ref{chapitrefini}.

\subsection{Exhaustions de Stein}

Commen\c{c}ons par rappeler la d\'efinition d'une exhaustion.

\begin{defi}\index{Exhaustion}
Soit $S$ un espace topologique. Une suite $(K_{n})_{n\in\N}$ de parties compactes de $S$ est une {\bf exhaustion de $S$} si elle v\'erifie les propri\'et\'es suivantes :
\begin{enumerate}[\it i)]
\item quel que soit $n\in\N$, le compact $K_{n}$ est contenu dans l'int\'erieur de $K_{n+1}$ ;
\item la r\'eunion des compacts $K_{n}$ est \'egale \`a $S$.
\end{enumerate}
\end{defi}

Le r\'esultat qui suit est classique (\emph{cf.}~\cite{GR}, IV, \S 1, th\'eor\`eme~4) et nous permettra de d\'emontrer une partie du th\'eor\`eme~B pour les faisceaux coh\'erents d\'efinis sur les couronnes ouvertes.

\begin{thm}
Soient $S$ un espace topologique et~$(K_{n})_{n\in\N}$ une exhaustion de~$S$. Soient $\Ss$ un faisceau de groupes ab\'eliens sur $S$ et $q\ge 2$ un nombre entier. Supposons que, quel que soit $n\in\N$, on ait 
$$H^{q-1}(K_{n},\Ss)=H^{q}(K_{n},\Ss)=0.$$
Alors on a \'egalement
$$H^q(S,\Ss)=0.$$
\end{thm}

\begin{defi}\index{Morphisme borne@Morphisme born\'e}
Soient $(A,\|.\|_{A})$ et $(B,\|.\|_{B})$ deux anneaux munis de semi-normes. Soit $\varphi : A\to B$ un morphisme d'anneaux. Nous dirons que le {\bf morphisme} $\varphi$ est {\bf born\'e} s'il existe un nombre r\'eel~$M$ tel que, pour tout \'el\'ement~$a$ de~$A$, nous ayons
$$\|\varphi(a)\|_{B}\le M\, \|a\|_{A}.$$
\end{defi}

Venons-en, \`a pr\'esent, aux exhaustions de Stein (\emph{cf.}~\cite{GR}, IV, \S 1, d\'efinition~6). 

\begin{defi}\label{exhStein}\index{Exhaustion!de Stein}\index{Espace de Stein!exhaustion de Stein|see{Exhaustion}}
Soient $(S,\Os_{S})$ un espace localement annel\'e et~$\Ss$ un faisceau de $\Os_{S}$-modules coh\'erent. Une suite $(K_{n})_{n\in\N}$ de parties compactes et de Stein de~$S$ est une {\bf exhaustion de Stein de~$S$ relativement au faisceau~$\Ss$} si c'est une exhaustion de~$S$ et si, quel que soit $n\in\N$, il existe une semi-norme~$\|.\|_{n}$ sur~$\Ss(K_{n})$ telle que, quel que soit $n\in\N$, les propri\'et\'es suivantes soient v\'erifi\'ees : 
\begin{enumerate}[\it i)]
\item la partie $\Ss(S)_{|K_{n}}$ de $\Ss(K_{n})$ est dense pour~$\|.\|_{n}$;
\item l'application de restriction $(\Ss(K_{n+1}),\|.\|_{n+1}) \to (\Ss(K_{n}),\|.\|_{n})$ est born\'ee ;
\item l'application de restriction $(\Ss(K_{n+1}),\|.\|_{n+1}) \to (\Ss(K_{n}),\|.\|_{n})$ envoie toute suite de Cauchy sur une suite convergente ;
\item tout \'el\'ement $s$ de $\Ss(K_{n+1})$ v\'erifiant $\|s\|_{n+1}=0$ est nul sur $K_{n}$.
\end{enumerate}
\end{defi}

Cette notion nous permettra de compl\'eter la d\'emonstration du th\'eor\`eme~B pour les faisceaux coh\'erents d\'efinis sur les couronnes ouvertes, par l'interm\'ediaire du r\'esultat suivant (\emph{cf.}~\cite{GR}, IV, \S 1, th\'eor\`eme~7).

\begin{thm}
Soient $(S,\Os_{S})$ un espace localement annel\'e et~$\Ss$ un faisceau de $\Os_{S}$-modules coh\'erent. Supposons qu'il existe une exhaustion de Stein de~$S$ relativement au faisceau~$\Ss$. Alors nous avons 
$$H^1(S,\Ss)=0.$$
\end{thm}



En regroupant les r\'esultats des deux th\'eor\`emes qui pr\'ec\`edent et celui du th\'eor\`eme \ref{BimpliqueA}, nous obtenons le r\'esultat suivant.

\begin{thm}\label{csqexh}
Soit $(S,\Os_{S})$ un espace localement annel\'e. Supposons que, pour tout faisceau de $\Os_{S}$-modules coh\'erent~$\Ss$, l'espace~$S$ poss\`ede une exhaustion de Stein relativement \`a~$\Ss$. Alors, l'espace~$S$ est de Stein.
\end{thm}

\subsection{Fermeture des modules}

Pour montrer que les exhaustions naturelles des couronnes ouvertes par des couronnes ferm\'ees sont bien des exhaustions de Stein, nous avons besoin de r\'esultats de fermeture sur certains faisceaux de modules. Nous leur consacrons cette partie. Les preuves que nous proposons sont inspir\'ees de~\cite{GuRo}, II, D, th\'eor\`emes~2 et~3. 


Commen\c{c}ons par introduire une notation. Soient~$(Y,\Os_{Y})$ un espace analytique, $y$ un point de~$Y$ et~$V$ un voisinage du point~$y$ dans~$Y$. Soient~$p\in\N$ et~$\Ms$ un sous-module de~$\Os_{Y,y}^p$. Nous noterons~$\Ms(V)$ le~$\Os_{Y}(V)$-module constitu\'e des \'el\'ements~$F$ de~$\Os_{Y}(V)^p$ dont le germe~$F_{y}$ en~$y$ appartient \`a~$\Ms$. D\'efinissons, maintenant, la notion de module fortement engendr\'e. Nous l'utiliserons constamment dans cette partie.

\begin{defi}\index{Fortement engendre@Fortement engendr\'e}
Soient~$(Y,\Os_{Y})$ un espace analytique et~$y$ un point de~$Y$. Soient~$p\in\N$ et~$\Ms$ un sous-module de~$\Os_{Y,y}^p$. Soient~$V$ un voisinage du point~$y$ dans~$Y$ et~$\|.\|$ une norme sur~$\Os_{Y}(V)$. Nous munissons le module produit~$\Os_{Y}(V)^p$ de la norme, que nous noterons encore~$\|.\|$, donn\'ee par le maximum des normes des coefficients. Soient~$q\in\N$ et~$G_{1},\ldots,G_{q}$ des \'el\'ements de~$\Os_{Y}(V)^p$. Nous dirons que la famille~$(G_{1},\ldots,G_{q})$ {\bf engendre fortement le module~$\Ms$ sur~$V$ pour la norme~$\|.\|$} s'il existe une constante~$C\in\R$ telle que, pour tout \'el\'ement~$F$ de~$\Ms(V)$, il existe des fonctions~$f_{1},\ldots,f_{q}$ dans~$\Os_{Y}(V)$ satisfaisant les propri\'et\'es suivantes :
\begin{enumerate}[\it i)]
\item $\disp F = \sum_{i=1}^q f_{i}\, G_{i}$ dans~$\Ms(V)$ ;
\item quel que soit~$i\in\cn{1}{q}$, nous avons~$\|f_{i}\| \le C\, \|F\|$.
\end{enumerate}

Nous dirons que le module~$\Ms$ est {\bf fortement engendr\'e sur~$V$ pour la norme~$\|.\|$} s'il existe une famille finie de~$\Os_{Y}(V)^p$qui engendre fortement le module~$\Ms$ sur~$V$ pour la norme~$\|.\|$.
\end{defi}

Les syst\`emes de g\'en\'erateurs forts jouissent de propri\'et\'es agr\'eables.


\begin{lem}\label{fese}
Soient~$(Y,\Os_{Y})$, $(Y',\Os_{Y'})$ et~$(Y'',\Os_{Y''})$ des espaces analytiques, $y$, $y'$ et~$y''$ des points de~$Y$, $Y'$ et~$Y''$, $p$, $p'$, $p''$ des entiers et~$\Ms$, $\Ms'$, $\Ms''$ des sous-modules de~$\Os_{Y,y}^p$, $\Os_{Y',y'}^{p'}$ et~$\Os_{Y'',y''}^{p''}$. Soient~$V$, $V'$ et~$V''$ des voisinages des points~$y$, $y'$ et~$y''$ dans~$Y$, $Y'$ et~$Y''$ et~$\|.\|$, $\|.\|'$ et~$\|.\|''$ des normes sur~$\Os_{Y}(V)$, $\Os_{Y'}(V')$ et~$\Os_{Y''}(V'')$. Supposons qu'il existe une suite exacte courte de groupes ab\'eliens
$$0\to \Ms'(V') \xrightarrow[]{u} \Ms(V) \xrightarrow[]{v} \Ms''(V'') \to 0$$
v\'erifiant les propri\'et\'es suivantes :
\begin{enumerate}[\it i)]
\item le morphisme~$u$ est une isom\'etrie ;
\item il existe un morphisme born\'e~$u_{0} : \Os_{Y'}(V') \to \Os_{Y}(V)$ qui v\'erifie
$$\forall f'\in \Os_{Y'}(V'),\, \forall F'\in\Ms(V'),\, u(f'\, F') = u_{0}(f')\, u(F')\ ;$$
\item le morphisme~$v$ est born\'e ;
\item il existe un morphisme born\'e~$\tau : \Os_{Y''}(V'') \to \Os_{Y}(V)$ qui v\'erifie
$$\forall f''\in \Os_{Y''}(V''),\, \forall F\in\Ms(V),\, v(\tau(f'')\, F) = f''\, v(F).$$
\end{enumerate}
Si les modules~$\Ms'$ et~$\Ms''$ sont fortement engendr\'es sur~$V'$ et~$V''$ pour les normes~$\|.\|'$ et~$\|.\|''$,  alors le module~$\Ms$ est fortement engendr\'e sur~$V$ pour la norme~$\|.\|$.
\end{lem}
\begin{proof}
Commen\c{c}ons par traduire les hypoth\`eses sur les morphismes born\'es. Il existe des constantes~$D_{u_{0}},D_{v},D_{\tau}\in\R$ telles que, quel que soit $f'\in\Os_{Y'}(V')$, on ait
$$\|u_{0}(f')\| \le D_{u_{0}}\, \|f'\|',$$
quel que soit \mbox{$F\in\Ms(V)$}, on ait
$$\|v(F)\|'' \le D_{v}\, \|F\|$$ 
et, quel que soit~$f''\in\Os_{Y''}(V'')$, on ait
$$\|\tau(f'')\|\le D_{\tau}\, \|f''\|''.$$

Supposons que les modules~$\Ms'$ et~$\Ms''$ sont fortement engendr\'es sur~$V'$ et~$V''$ pour les normes~$\|.\|'$ et~$\|.\|''$. Il existe un entier~$r'\in\N$ et une famille~$(G'_{1},\ldots,G'_{r'})$ de~$\Ms'(V')$ qui engendre fortement le module~$\Ms'$ sur~$V'$ pour la norme~$\|.\|'$, avec une certaine constante~$C'\in\R$. De m\^eme, il existe un entier~$r''\in\N$ et une famille~$(G''_{1},\ldots,G''_{r''})$ de~$\Ms''(V'')$ qui engendre fortement le module~$\Ms''$ sur~$V''$ pour la norme~$\|.\|''$, avec une certaine constante~$C''\in\R$. Quel que soit~$i\in\cn{1}{r'}$, nous posons
$$H'_{i} = u(G'_{i}).$$
Quel que soit~$j\in\cn{1}{r''}$, nous choisissons un \'el\'ement~$H''_{j}$ de~$\Ms(V)$ tel que
$$v(H''_{j})= G''_{j}.$$
Nous allons montrer que la famille $(H'_{1},\ldots,H'_{r'},H''_{1},\ldots, H''_{r''})$ de~$\Ms(V)$ engendre fortement le module~$\Ms$ sur~$V$ pour la norme~$\|.\|$. 

Soit~$F\in\Ms(V)$. Alors~$v(F)\in\Ms'(V')$. Il existe donc~$f''_{1},\ldots,f''_{r''}\in\Os_{Y''}(V'')$ tels que l'on ait
\begin{enumerate}[\it i)]
\item $\disp v(F) = \sum_{j=1}^{r''} f''_{j}\, G''_{j}$ ;
\item $\forall j \in\cn{1}{r''}$, $\|f''_{j}\|'' \le C''\, \|v(F)\|''$.
\end{enumerate}
Posons
$$F_{0} = F - \sum_{j=1}^{r''} \tau(f''_{j})\, H''_{j}.$$
Quel que soit~$j\in\cn{1}{r''}$, nous avons
$$\|\tau(f''_{j})\| \le D_{\tau}\, \|f''_{j}\|'' \le D_{\tau}C''\,\|v(F)\|''\le D_{\tau}C''D_{v}\,\|F\|.$$
Nous en d\'eduisons que
$$\|F_{0}\| \le \left(1+D_{\tau}C''D_{v}\,\sum_{j=1}^{r''} \|H''_{j}\|\right)\,\|F\|.$$
Posons
$$M = 1+D_{\tau}C''D_{v}\,\sum_{j=1}^{r''} \|H''_{j}\|.$$

Nous avons
$$v(F_{0}) = v(F) -  \sum_{j=1}^{r''} v(\tau(f''_{j})\, H''_{j}) = v(F) -  \sum_{j=1}^{r''} f''_{j}\, G''_{j} =0.$$
Par cons\'equent, $F_{0} \in \Ker(v)=\Im(u)$. On en d\'eduit qu'il existe~$F'\in\Ms'(V')$ tel que
$$u(F')=F_{0}.$$ 
Il existe \'egalement~$f'_{1},\ldots,f'_{r'}\in\Os_{Y'}(V')$ tels que l'on ait
\begin{enumerate}[\it i)]
\item $\disp F' = \sum_{i=1}^{r'} f'_{i}\, G'_{i}$ ;
\item $\forall i \in\cn{1}{r'}$, $\|f'_{i}\|' \le C'\, \|F'\|'$.
\end{enumerate}
Nous avons finalement
$$\begin{array}{rcl}
F &=& \disp F_{0} +  \sum_{j=1}^{r''} \tau(f''_{j})\, H''_{j}\\
&=& \disp u\left( \sum_{i=1}^{r'} f'_{i}\, G'_{i} \right) +  \sum_{j=1}^{r''} \tau(f''_{j})\, H''_{j}\\
&=& \disp  \sum_{i=1}^{r'}  u_{0}(f'_{i})\, H'_{i} +  \sum_{j=1}^{r''} \tau(f''_{j})\, H''_{j}.
\end{array}$$
Nous avons vu pr\'ec\'edemment que la norme des coefficients~$\tau(f''_{j})$, avec~\mbox{$j\in\cn{1}{r''}$}, est born\'ee en fonction de celle de~$\|F\|$. En outre, quel que soit~$i\in\cn{1}{r'}$, nous avons
$$\|u_{0}(f'_{i})\| \le D_{u_{0}}\, \|f'_{i}\|' \le D_{u_{0}}C'\, \|F'\|' \le D_{u_{0}}C'\, \|F_{0}\|\le D_{u_{0}}C'M\, \|F\|.$$
On en d\'eduit que la famille $(H'_{1},\ldots,H'_{r'},H''_{1},\ldots, H''_{r''})$ de~$\Ms(V)$ engendre fortement le module~$\Ms$ sur~$V$ pour la norme~$\|.\|$. 
\end{proof}

\begin{cor}\label{corfe}
Soient~$(Y,\Os_{Y})$ un espace analytique et~$y$ un point de~$Y$. Soient~$V$ un voisinage du point~$y$ dans~$Y$ et~$\|.\|$ une norme sur~$\Os_{Y}(V)$. Supposons que tous les id\'eaux de~$\Os_{Y,y}$ sont fortement engendr\'es sur~$V$ pour la norme~$\|.\|$. Alors, quel que soit~$p\in\N^*$, tous les sous-modules de~$\Os_{Y,y}^p$ sont fortement engendr\'es sur~$V$ pour la norme~$\|.\|$.
\end{cor}
\begin{proof}
Nous allons d\'emontrer le r\'esultat par r\'ecurrence. L'initialisation pour~$p=1$ n'est autre que l'hypoth\`ese. Soit~$p\in\N^*$ pour lequel le r\'esultat est vrai. Soit~$\Ms$ un sous-module de~$\Os_{Y,y}^{p+1}$. Notons~$\Ms'$ l'id\'eal de~$\Os_{Y,y}$ compos\'e des \'el\'ements~$f$ de~$\Os_{Y,y}$ tels que~$(0,\ldots,0,f)$ appartient \`a~$\Ms$. Notons~$\Ms''$ le sous-module de~$\Os_{Y,y}^p$ dont les \'el\'ements sont les $p$ premi\`eres composantes des \'el\'ements de~$\Ms$. Les morphismes naturels
$$0\to \Ms'(V) \xrightarrow[]{u} \Ms(V) \xrightarrow[]{v} \Ms''(V) \to 0$$
forment une suite exacte courte de groupes ab\'eliens. Montrons que les propri\'et\'es du lemme~\ref{fese} sont v\'erifi\'ees. Le morphisme~$u$ est bien une isom\'etrie. Choisissons pour~$u_{0}$ le morphisme identit\'e sur~$\Os_{Y}(V)$. Les propri\'et\'es du point~{\it ii)} sont alors v\'erifi\'ees. Le morphisme~$v$ est born\'e (et l'on peut m\^eme choisir la constante~$1$). Nous pouvons choisir pour~$\tau$ le morphisme identit\'e sur~$\Os_{Y}(V)$. L'hypoth\`ese de l'\'enonc\'e nous assure que l'id\'eal~$\Ms'$ est fortement engendr\'e sur~$V$ pour la norme~$\|.\|$. L'hypoth\`ese de r\'ecurrence nous assure que tel est \'egalement le cas pour le module~$\Ms''$. D'apr\`es le lemme~\ref{fese}, le module~$\Ms$ est, lui aussi, fortement engendr\'e sur~$V$ pour la norme~$\|.\|$.
\end{proof}

\'Enon\c{c}ons, \`a pr\'esent, quelques conditions permettant d'assurer que certains modules poss\`edent des syst\`emes de g\'e\-n\'e\-ra\-teurs forts.

\begin{lem}\label{genefortscorps}
Soient~$(Y,\Os_{Y})$ un espace analytique et~$y$ un point de~$Y$. Soit~$V$ un voisinage du point~$y$ dans~$Y$. Munissons l'anneau~$\Os_{Y}(V)$ de la norme uniforme~$\|.\|_{V}$. Supposons que le morphisme de restriction $\Os_{Y}(V) \to \Os_{Y,y}$ est injectif et que l'anneau local~$\Os_{Y,y}$ est un corps. Alors, quel que soit~$p\in\N^*$, tous les sous-modules de~$\Os_{Y,y}^p$ sont fortement engendr\'es sur~$V$ pour la norme~$\|.\|_{V}$.
\end{lem}
\begin{proof}
D'apr\`es le corollaire \ref{corfe}, il suffit de montrer que tous les id\'eaux de~$\Os_{Y,y}$ sont fortement engendr\'es sur~$V$ pour la norme~$\|.\|_{V}$. Puisque l'anneau local~$\Os_{Y,y}$ est un corps, il ne poss\`ede que deux id\'eaux : $\Os_{Y,y}$ et~$(0)$. Il est \'evident que la famille~$(1)$ engendre fortement l'id\'eal~$\Os_{Y,y}$ sur~$V$ pour la norme~$\|.\|_{V}$. L'injectivit\'e du morphisme $\Os_{Y}(V) \to \Os_{Y,y}$ assure que la famille~$(0)$ engendre fortement l'id\'eal~$\Os_{Y,y}$ sur~$V$ pour la norme~$\|.\|_{V}$.
\end{proof}

\begin{lem}
Soient~$(Y,\Os_{Y})$ un espace analytique et~$y$ un point de~$Y$. Supposons que l'anneau local~$\Os_{Y,y}$ est un anneau de valuation discr\`ete. Soit~$V$ un voisinage du point~$y$ dans~$Y$ et~$\pi$ une uniformisante forte de l'anneau~$\Os_{Y,y}$ sur~$V$. Alors, la famille~$(\pi)$ engendre fortement l'id\'eal~$\pi\, \Os_{Y,y}$ sur~$V$ pour la norme~$\|.\|_{V}$.
\end{lem}
\begin{proof}
C'est une simple traduction des d\'efinitions.
\end{proof}

\begin{cor}\label{genefortsavd}
Soient~$(Y,\Os_{Y})$ un espace analytique et~$y$ un point de~$Y$. Supposons que l'anneau local~$\Os_{Y,y}$ est un anneau de valuation discr\`ete. Notons~$\m$ son id\'eal maximal. Soit~$V$ un voisinage du point~$y$ dans~$Y$ et~$\pi$ une uniformisante forte de l'anneau~$\Os_{Y,y}$ sur~$V$. Supposons que le morphisme de restriction $\Os_{Y}(V) \to \Os_{Y,y}$ est injectif. Alors, quel que soit~$p\in\N^*$, tous les sous-modules de~$\Os_{Y,y}^p$ sont fortement engendr\'es sur~$V$ pour la norme~$\|.\|_{V}$.
\end{cor}
\begin{proof}
D'apr\`es le corollaire \ref{corfe}, il suffit de montrer que tous les id\'eaux de~$\Os_{Y,y}$ sont fortement engendr\'es sur~$V$ pour la norme~$\|.\|_{V}$. Puisque l'anneau local~$\Os_{Y,y}$ est un anneau de valuation discr\`ete, ses id\'eaux sont de la forme~$(0)$ ou~$(\pi^n)$ avec~$n\in\N$. L'injectivit\'e du morphisme $\Os_{Y}(V) \to \Os_{Y,y}$ assure que la famille~$(0)$ engendre fortement l'id\'eal~$(0)$ sur~$V$ pour la norme~$\|.\|_{V}$. On constate imm\'ediatement que la famille~$(1)$ engendre fortement~$\Os_{Y,y}$ sur~$V$ pour la norme~$\|.\|_{V}$. Finalement, on montre, par r\'ecurrence, en utilisant le lemme pr\'ec\'edent, que, quel que soit~$n\in\N$, la famille~$(\pi^n)$ engendre fortement l'id\'eal~$\pi^n\,\Os_{Y,y}$ sur~$V$ pour la norme~$\|.\|_{V}$.
\end{proof}

Appliquons ces r\'esultats au cas de la base~$B$ et de l'espace~$X$.

\begin{cor}\label{modulescorpsavd}\index{Fortement engendre@Fortement engendr\'e!corps}\index{Fortement engendre@Fortement engendr\'e!anneau de valuation discrete@anneau de valuation discr\`ete}
Soit~$Y$ l'un des deux espaces~$B$ et~$X$. Soit~$y$ un point de~$Y$ en lequel l'anneau local~$\Os_{Y,y}$ est un corps ou un anneau de valuation discr\`ete. Il existe un syst\`eme fondamental~$\Vs$ de voisinages compacts du point~$y$ dans~$Y$ tel que, pour tout \'el\'ement~$V$ de~$\Vs$ et tout entier~$p\in\N^*$, tout sous-module de~$\Os_{Y,y}^p$ est fortement engendr\'e sur~$V$ pour la norme~$\|.\|_{V}$.


\end{cor}
\begin{proof}
D'apr\`es la proposition \ref{prolanB} et le th\'eor\`eme \ref{prolan}, le principe du prolongement analytique vaut au voisinage de tout point de l'espaces~$Y$. Par cons\'equent, pour toute partie connexe~$V$ contenant le point~$y$, le morphisme de restriction $\Os_{Y}(V) \to \Os_{Y,y}$ est injectif.

Supposons, tout d'abord, que l'anneau local~$\Os_{Y,y}$ est un corps. D'apr\`es le lemme \ref{genefortscorps} et la remarque qui pr\'ec\`ede, il suffit de d\'emontrer que le point~$y$ poss\`ede un syst\`eme fondamental de voisinages compacts et connexes. C'est \'evident pour l'espace~$B$ et c'est encore vrai pour l'espace~$X$, d'apr\`es le th\'eor\`eme \ref{resume}.

Supposons, \`a pr\'esent, que l'anneau local~$\Os_{Y,y}$ est un anneau de valuation discr\`ete. Soit~$\pi$ une uniformisante de l'anneau~$\Os_{Y,y}$ et~$U$ un voisinage du point~$y$ sur lequel elle est d\'efinie. D'apr\`es le corollaire \ref{genefortsavd} et la remarque figurant au d\'ebut de la preuve, il suffit de montrer que le point~$y$ de~$Y$ poss\`ede un syst\`eme fondamental~$\Vs$ de voisinages compacts et connexes tel que, pour tout \'el\'ement~$V$ de~$\Vs$, la fonction~$\pi$ est une uniformisante forte de l'anneau~$\Os_{Y,y}$ sur~$V$. C'est le cas, d'apr\`es le lemme \ref{uniforteB}  et le th\'eor\`eme \ref{uniforte}.
\end{proof}

Il nous reste \`a traiter le cas des points rigides des fibres extr\^emes de l'espace~$X$. Soient~$V$ une partie compacte de~$B$ et~$t$ un nombre r\'eel strictement positif. D'apr\`es la proposition~\ref{imagedisque}, le morphisme naturel $A[T]\to \Os(V)[\![T]\!]$ se prolonge un un morphisme injectif
$$j_{V,t} : \Os(\overline{D}_{V}(t)) \hookrightarrow \Os(V)\of{\la}{|T|\le t}{\ra}^\dag.$$
Soit~$f$ un \'el\'ement de~$\Os(\overline{D}_{V}(t))$. Il existe une suite~$(a_{k})_{k\in\N}$ d'\'el\'ements de~$\Os(V)$ telle que l'on ait l'\'egalit\'e
$$j_{V,t}(f) = \sum_{k\in\N} a_{k}\, T^k \textrm{ dans } \Os(V)\of{\la}{|T|\le t}{\ra}^\dag.$$
Nous posons alors
$$\|f\|_{V,t} = \sum_{k\in\N} \|a_{k}\|_{V}\, t^k \in\R_{+}.$$
La fonction~$\|.\|_{V,t}$ d\'efinit une norme sur l'anneau~$\Os(\overline{D}_{V}(t))$.

\begin{lem}
Soit~$\m$ un \'el\'ement de~$\Sigma_{f}$. Posons~$V_{0}=\of{[}{a_{\m},\tilde{a}_{\m}}{]}$. Soit~$x$ le point de la fibre extr\^eme~$\tilde{X}_{\m}$ d\'efini par l'\'equation $T(x)=0$. Soient~$V$ un voisinage compact et connexe du point~$\tilde{a}_{\m}$ dans~$V_{0}$ et~$t$ un \'el\'ement de l'intervalle~$\of{]}{0,1}{[}$. Pour tout \'el\'ement~$F$ de~$\Os(\overline{D}_{V}(t))$ dont l'image dans l'anneau local~$\Os_{X,x}$ est divisible par~$\pi_{\m}$, il existe un \'el\'ement~$d_{V,t}(F)$ de~$\Os(\overline{D}_{V}(t))$ qui v\'erifie l'\'egalit\'e
$$F = \pi_{\m}\, d_{V,t}(F) \textrm{ dans } \Os(\overline{D}_{V}(t)).$$
\end{lem}
\begin{proof}
Soit un \'el\'ement~$F$ de~$\Os(\overline{D}_{V}(t))$ dont l'image dans l'anneau local~$\Os_{X,x}$ est divisible par~$\pi_{\m}$. Consid\'erons la restriction de la fonction~$F$ \`a la trace~$E$ du disque~$\overline{D}_{V}(t)$ sur la fibre extr\^eme~$\tilde{X_{\m}}$. Cette fonction est nulle au voisinage du point~$x$. Le principe du prolongement analytique sur~$\tilde{X}_{\m}$ assure qu'elle est nulle en tout point de~$E$. Or, d'apr\`es le corollaire \ref{avd3ext}, en tout point~$y$ de~$E$ diff\'erent de~$x$, l'anneau local~$\Os_{X,y}$ est un anneau de valuation discr\`ete d'uniformisante~$\pi_{\m}$. On en d\'eduit que la fonction~$F$ est divisible par~$\pi_{\m}$ au voisinage de tout point de~$E$. Soit~$z$ un point de~$\overline{D}_{V}(t)\setminus E$. La fonction~$\pi_{\m}$ est inversible au voisinage de ce point. Par cons\'equent, la fonction~$F$ est multiple de~$\pi_{\m}$ au voisinage de ce point.

En utilisant le fait que les anneaux locaux sont int\`egres et que le principe du prolongement analytique vaut sur la partie connexe~$\overline{D}_{V}(t)$ de~$X$, nous obtenons l'existence d'un \'el\'ement~$d_{V,t}(F)$ de~$\Os(\overline{D}_{V}(t))$ qui v\'erifie l'\'egalit\'e
$$F = \pi_{\m}\, d_{V,t}(F) \textrm{ dans } \Os(\overline{D}_{V}(t)).$$
\end{proof}

\begin{lem}\label{lemmemod}
Soit~$\m$ un \'el\'ement de~$\Sigma_{f}$. Posons~$V_{0}=\of{[}{a_{\m},\tilde{a}_{\m}}{]}$. Soit~$x$ le point de la fibre extr\^eme~$\tilde{X}_{\m}$ d\'efini par l'\'equation $T(x)=0$. Soit~$\Is$ un id\'eal de~$\Os_{X,x}$. Supposons que pour tout voisinage~$W$ du point~$x$ dans~$X$, il existe un voisinage compact et connexe~$V$ du point~$\tilde{a}_{\m}$ dans~$V_{0}$ et un nombre r\'eel~$t>0$ tels que le disque compact~$\overline{D}_{V}(t)$ soit contenu dans~$W$ et l'id\'eal~$\pi_{\m}\Is$ soit fortement engendr\'e sur~$\overline{D}_{V}(t)$ pour la norme~$\|.\|_{V,t}$. Alors, il en est de m\^eme pour l'id\'eal~$\Is$.
\end{lem}
\begin{proof}
Soit~$W$ un voisinage du point~$x$ dans~$X$. Par hypoth\`ese, il existe un voisinage compact et connexe~$V$ du point~$\tilde{a}_{\m}$ dans~$V_{0}$, un nombre r\'eel~$t>0$, un entier~$q$ et des \'el\'ements $G_{1},\ldots,G_{q}$ de~$\Os(\overline{D}_{V}(t))$ v\'erifiant les propri\'et\'es suivantes :
\begin{enumerate}[\it i)]
\item le disque compact~$\overline{D}_{V}(t)$ est contenu dans~$W$ ;
\item la famille~$(G_{1},\ldots,G_{q})$ engendre fortement le module~$\pi_{\m}\Is$ sur~$\overline{D}_{V}(t)$ pour la norme~$\|.\|_{V,t}$, avec une certaine constante~$C$. 
\end{enumerate}

Nous reprenons les notations du lemme qui pr\'ec\`ede. Soit~$F$ un \'el\'ement de~$\Is(\overline{D}_{V}(t))$. La fonction~$\pi_{\m}\, F$ appartient alors \`a~$\pi_{\m}\Is(\overline{D}_{V}(t))$. Il existe donc des \'el\'ements $f_{1},\ldots,f_{q}$ de~$\Os(\overline{D}_{V}(t))$ satisfaisant les propri\'et\'es suivantes :
\begin{enumerate}[\it i)]
\item $\disp \pi_{\m}\,F = \sum_{i=1}^q f_{i}\, G_{i} = \pi_{\m}\, \sum_{i=1}^q f_{i}\, d_{V,t}(G_{i})$ dans~$\Os(\overline{D}_{V}(t))$ ;
\item quel que soit~$i\in\cn{1}{q}$, nous avons~$\|f_{i}\|_{V,t} \le C\, \|\pi_{\m}\, F\|_{V,t}$.
\end{enumerate}
La partie~$\overline{D}_{V}(t)$ \'etant connexe, l'int\'egrit\'e des anneaux locaux et le principe du prolongement analytique nous assurent que
$$F = \sum_{i=1}^q f_{i}\, d_{V,t}(G_{i}) \textrm{ dans } \Os(\overline{D}_{V}(t)).$$
En outre, pour tout \'el\'ement~$i$ de~$\cn{1}{q}$, nous avons
$$\|f_{i}\|_{V,t} \le C\, \|\pi_{\m}\, F\|_{V,t} \le C\|\pi_{\m}\|_{V}\, \|F\|_{V,t}.$$ 
On en d\'eduit que la famille~$(d(G_{1}),\ldots,d(G_{q}))$ engendre fortement le module~$\Is$ sur~$\overline{D}_{V}(t)$ pour la norme~$\|.\|_{V,t}$.
\end{proof}

\begin{rem}
L'implication r\'eciproque de celle \'enonc\'ee dans le lemme pr\'ec\'edent est valable et sa d\'emonstration est d'ailleurs \'evidente.
\end{rem}

\begin{prop}\label{modulesrationnel}\index{Fortement engendre@Fortement engendr\'e!point rationnel d'une fibre extreme de A1@point rationnel d'une fibre extr\^eme de $\AA$}
Soit~$\m$ un \'el\'ement de~$\Sigma_{f}$. Posons~$V_{0}=\of{[}{a_{\m},\tilde{a}_{\m}}{]}$. Soit~$x$ le point de la fibre extr\^eme~$\tilde{X}_{\m}$ d\'efini par l'\'equation $T(x)=0$. Soient~$p$ un entier non nul et~$\Ms$ un sous-module de~$\Os_{Y,y}^p$. Soit $W$ un voisinage du point~$x$ dans~$X$. Alors, il existe un voisinage compact et connexe~$V$ du point~$\tilde{a}_{\m}$ dans~$V_{0}$ et un nombre r\'eel~$t>0$ tels que le disque compact~$\overline{D}_{V}(t)$ soit contenu dans~$W$ et le module~$\Ms$ soit fortement engendr\'e sur~$\overline{D}_{V}(t)$ pour la norme~$\|.\|_{V,t}$.
\end{prop}
\begin{proof}
D'apr\`es le corollaire~\ref{corfe}, le cas~$p=1$ entra\^ine les autres. Nous pouvons donc supposer que~$p=1$. Dans ce cas, le module~$\Ms$ est un id\'eal de~$\Os_{X,x}$. Dans le cas o\`u l'id\'eal~$\Ms$ est nul, le principe du prolongement analytique nous permet de conclure. Nous supposerons, d\'esormais, que l'id\'eal~$\Ms$ n'est pas nul. Rappelons que, d'apr\`es le th\'eor\`eme~\ref{anneaulocal}, l'anneau local~$\Os_{X,x}$ est un anneau de s\'eries convergentes \`a coefficients dans~$\Os_{B,\tilde{a}_{\m}}$. Plus pr\'ecis\'ement, il est naturellement isomorphe \`a l'anneau~$L_{\tilde{a}_{\m}}$ d\'efini \`a la section~\ref{algloc}. Reprenons, \`a pr\'esent, les notations du lemme~\ref{lemU}. Notons
$$w=\min\{v(F)\,|\, F\in\Ms, G\ne 0\}.$$ 
D'apr\`es ce lemme, il existe un id\'eal~$\Ns$ de~$\Os_{X,x}$ v\'erifiant les propri\'et\'es suivantes :
\begin{enumerate}[\it i)]
\item $\Ms = \pi_{\m}^w\,\Ns$ ;
\item il existe un \'el\'ement~$G$ de~$\Ns$ qui v\'erifie~$G(\tilde{a}_{\m})=0$. 
\end{enumerate}
D'apr\`es le lemme~\ref{lemmemod}, il suffit de d\'emontrer le r\'esultat voulu pour l'id\'eal~$\Ns$.

Il existe un voisinage compact et connexe~$V$ de~$\tilde{a}_{\m}$ dans~$V_{0}$ et un nombre r\'eel~$t>0$ tels que 
$$G \in\Bs(V)\of{\la}{|T|\le t}{\ra}.$$ 
D'apr\`es le th\'eor\`eme de pr\'eparation de Weierstra{\ss} (\emph{cf.} th\'eor\`eme~\ref{preparation}), il existe une fonction inversible $E \in \Os_{X,x}$ et un polyn\^ome~\mbox{$\Omega\in\Os_{B,b}[T]$} distingu\'e de degr\'e $d\in\N$ tel que l'on ait l'\'egalit\'e $G = E\, \Omega$ dans~$\Os_{X,x}$. Quitte \`a restreindre~$V$ et \`a diminuer~$t$, nous pouvons supposer que cette \'egalit\'e vaut dans~$\Bs(V)\of{\la}{|T|\le t}{\ra}$. Remarquons que~$\Omega$ est un \'el\'ement de~$\Ns$. D'apr\`es la proposition~\ref{voisdep}, quitte \`a diminuer encore~$V$ et~$t$, nous pouvons supposer que le disque compact~$\overline{D}_{V}(t)$ est contenu dans~$W$.

D'apr\`es le th\'eor\`eme de division de Weierstra{\ss} semi-local, quitte \`a diminuer encore~$V$ et~$t$, nous pouvons supposer que, quel que soit~$u\in\of{[}{t,1}{]}$, pour tout \'el\'ement~$F$ de~$\Bs(V)\of{\la}{|T|\le u}{\ra}$, il existe un unique \'el\'ement~$(Q,R)$ appartenant \`a $(\Bs(V)\of{\la}{|T|\le u}{\ra})^2$ tel que 
\begin{enumerate}[\it i)]
\item $R$ soit un polyn\^ome de degr\'e strictement inf\'erieur \`a $d$ ;
\item $F=Q\,\Omega+R$.
\end{enumerate}
En outre, il existe une constante $C\in\R_{+}^*$, ind\'ependante de~$u$ et de~$F$, telle que l'on ait les in\'egalit\'es
$$\left\{{\renewcommand{\arraystretch}{1.3}\begin{array}{rcl}
\|Q\|_{V,u} &\le& C\,\|F\|_{V,u}\ ;\\
\|R\|_{V,u} &\le& C\,\|F\|_{V,u}.
\end{array}}\right.$$

Soit~$F$ un \'el\'ement de~$\Os(\overline{D}_{V}(t))$. D'apr\`es la proposition~\ref{imagedisque}, il existe \mbox{$u\in\of{]}{t,1}{]}$} tel que
$$F\in\Os(V)\of{\la}{|T|\le u}{\ra} = \Os(V)\of{\la}{|T|\le u}{\ra}.$$ 
En appliquant le r\'esultat pr\'ec\'edent, nous obtenons deux \'el\'ements~$Q$ et~$R$ de $\Bs(V)\of{\la}{|T|\le u}{\ra}$ et donc de~$\Os(\overline{D}_{V}(t))$, d'apr\`es le th\'eor\`eme~\ref{isodisque}. On en d\'eduit que~$Q\,\Omega$ appartient \`a~$\Ns(\overline{D}_{V}(t))$ et donc que~$R$ appartient \`a~$\Ns(\overline{D}_{V}(t))$. Il existe~$a_{0},\ldots,a_{d-1}\in \Os(V)$ tels que
$$R(T) = \sum_{i=0}^{d-1} a_{i}\, T^i.$$
Nous d\'efinissons un morphisme de groupes~$r$ en associant \`a l'\'el\'ement~$F$ la famille~$(a_{0},\ldots,a_{d-1})$. Les majorations du th\'eor\`eme de division de Weierstra{\ss} et du lemme~\ref{inegcoeff} nous assurent que
$$\|r(F)\|_{V,t} \le C t^{1-d}\, \|F\|_{V,t}.$$

Notons~$\Ns''$ le sous-$\Os_{B,\tilde{a}_{\m}}$-module de~$\Os_{B,\tilde{a}_{\m}}^d$ form\'e par les familles de coefficients des polyn\^omes de~$\Ns$ dont le degr\'e est strictement inf\'erieur \`a~$d$. Notons~$\Ns'$ l'id\'eal de~$\Os_{X,x}$ engendr\'e par~$\Omega$ et 
$$u : \Ns'(\overline{D}_{V}(t)) \to \Ns(\overline{D}_{V}(t))$$ 
l'injection canonique. D'apr\`es le th\'eor\`eme de division de Weierstra{\ss}, nous disposons alors d'une suite exacte
$$0 \to \Ns'(\overline{D}_{V}(t))  \xrightarrow[]{u} \Ns(\overline{D}_{V}(t)) \xrightarrow[]{r} \Ns''(\overline{D}_{V}(t)) \to 0.$$
Montrons qu'elle v\'erifie les conditions du lemme~\ref{fese}. Le morphisme~$u$ est bien une isom\'etrie. Nous pouvons choisir l'identit\'e de~$\Os_{Y}(\overline{D}_{V}(t))$ pour le morphisme~$u_{0}$. Nous avons montr\'e pr\'ec\'edemment que le morphisme~$r$ \'etait born\'e. Pour le morphisme~$\tau$, nous choisissons le morphisme naturel~$\Os(V) \to \Os(\overline{D}_{V}(t))$. Il est \'egalement born\'e.

En outre, la famille~$(G)$ engendre fortement le module~$\Ns'$ sur~$V$ pour la norme~$\|.\|_{V,t}$, toujours d'apr\`es le th\'eor\`eme de division de Weierstra{\ss}. La description explicite de l'espace~$V$ et des fonctions sur cet espace assure que le module~$\Ns''$ est \'egalement fortement engendr\'e sur~$V$ pour la norme~$\|.\|_{V}$. Nous d\'eduisons alors du lemme~\ref{fese} que le module~$\Ns$ est fortement engendr\'e sur~$V$ pour la norme~$\|.\|_{V,t}$.

\end{proof}

\begin{rem}
Ce r\'esultat vaut \'egalement pour les points rationnels des autres fibres. La d\'emonstration en est d'ailleurs plus simple puisque l'anneau local en le point de la base est alors un corps. Il vaut encore pour les points rationnels des fibres des espaces affines de dimension plus grande. Nous pourrions \'egalement l'adapter pour les points rigides, \`a condition de prendre la peine d\'efinir des normes ad\'equates.
\end{rem}

D\'emontrons, \`a pr\'esent, le r\'esultat sur la fermeture des modules que nous avions en vue.

\begin{thm}\label{modfermes}\index{Fermeture des modules}\index{Faisceau structural!fermeture des modules}
Soient $x$ un point de $X$, $p$ un entier non nul et~$\Ms$ un sous-module de~$\Os_{X,x}^p$. Soient $U$ un voisinage de $x$ dans $X$ et $F$ un \'el\'ement de $\Os(U)^p$. Supposons qu'il existe une suite $(F_{k})_{k\in\N}$ de $\Os(U)^p$ qui converge vers uniform\'ement vers $F$ sur $U$ et que, quel que soit $k\in\N$, on ait $(F_{k})_{x} \in\Ms$. Alors, on a
$$F_{x} \in \Ms.$$
\end{thm}
\begin{proof}
Nous devons distinguer plusieurs cas : celui o\`u l'anneau local $\Os_{X,x}$ est un corps, celui o\`u c'est un anneau de valuation discr\`ete et celui o\`u le point $x$ est un point rigide de sa fibre. La d\'emonstration est similaire dans les trois cas. Nous ne traiterons que le dernier qui est le plus difficile, en particulier \`a cause de la diff\'erence, pour les fonctions d\'efinies sur des disques, entre leur norme en tant que s\'erie et leur norme uniforme. Seuls les point rigides des fibres extr\^emes ne sont pas trait\'es dans les autres cas. Nous supposerons donc que $x$ est de ce type. D'apr\`es la proposition~\ref{isoext}, nous pouvons nous ramener au cas d'un point rationnel. Quitte \`a nous placer sur un voisinage assez petit du point $x$, puis \`a effectuer une translation, nous pouvons supposer que le point $x$ est le point de sa fibre d\'efini par l'\'equation $T(x)=0$.

D'apr\`es la proposition~\ref{voisdep}, il existe un voisinage~$W$ de~$b$ dans~$B$ et un nombre r\'eel~$u>0$ tels que la partie~$\overline{D}_{V}(t)$ soit contenue dans~$U$. D'apr\`es le th\'eor\`eme \ref{modulesrationnel}, il existe un voisinage compact et connexe~$V$ de~$b$ dans~$W$, un nombre r\'eel~$t\in\of{]}{0,u}{[}$, un entier~$q\in\N$ et des \'el\'ements $G_{1},\ldots, G_{q}$ de~$\overline{D}_{V}(t)$ tels que la famille~$(G_{1},\ldots,G_{q})$ engendre fortement le module~$\Ms$ sur~$\overline{D}_{V}(t)$ pour la norme~$\|.\|_{V,t}$, avec une certaine constante~$C$. 

Quitte \`a extraire une sous-suite de $(F_{k})_{k\in\N}$, nous pouvons supposer que, quel que soit $k\in\N^*$, nous avons 
$$\|F_{k}-F_{k-1}\|_{\overline{D}_{V}(u)} \le 2^{-k}.$$ 
D'apr\`es la proposition \ref{comparaisoncouronne}, nous avons alors 
$$\|F_{k}-F_{k-1}\|_{V,t} \le \frac{u}{u-t}\, 2^{-k}.$$ 

Construisons, \`a pr\'esent, par r\'ecurrence, des suites $(f_{k,1})_{k\in\N},\ldots,(f_{k,q})_{k\in\N}$ de~$\Os(\overline{D}_{V}(t))$ v\'erifiant les propri\'et\'es suivantes : quel que soit $k\in\N$, nous avons 
$$F_{k} = \sum_{j=1}^q f_{k,j} G_{j}$$
et, quel que soit $k\in\N^*$, nous avons
$$\forall j\in\cn{1}{q},\, \|f_{k,j}-f_{k-1,j}\|_{\overline{D}_{V}(t)} \le \frac{C}{2^{k}}.$$

Initialisons la r\'ecurrence. Pour construire $f_{0,1},\ldots,f_{0,q}$, il suffit d'utiliser le fait que la famille~$(G_{1},\ldots,G_{q})$ engendre fortement le module~$\Ms$ sur~$\overline{D}_{V}(t)$ pour la norme~$\|.\|_{V,t}$ avec la constante~$C$ et de l'appliquer \`a la fonction~$F_{0}$. 

Soit $k\in\N^*$ et supposons avoir construit $f_{k-1,1},\ldots,f_{k-1,q} \in\Os(\overline{D}_{V}(t))$ v\'erifiant les propri\'et\'es demand\'ees. En appliquant la propri\'et\'e de g\'en\'eration forte \`a la fonction $F_{k}-F_{k-1}$, on montre qu'il existe $g_{k,1},\ldots,g_{k,q}\in \Os(\overline{D}_{V}(t))$ v\'erifiant 
$$F_{k}-F_{k-1} = \sum_{j=1}^q g_{k,j} G_{j}$$
et
$$\forall j\in\cn{1}{q},\, \|g_{k,j}\|_{V,t} \le C \|F_{k}-F_{k-1}\|_{V,t} \le \frac{C}{2^{k}}.$$
Pour $j\in\cn{1}{q}$, posons 
$$f_{k,j}=f_{k-1,j}+g_{k,j}.$$
On obtient alors le r\'esultat voulu car, quel que soit~$j\in\cn{1}{q}$, nous avons 
$$\|g_{k,j}\|_{\overline{D}_{V}(t)} \le \|g_{k,j}\|_{V,t}.$$

Soit $j\in\cn{1}{q}$. D'apr\`es les in\'egalit\'es pr\'ec\'edentes, la suite~$(f_{k,j})_{k\in\N}$ est de Cauchy dans $\Os(\overline{D}_{V}(t))$. Soit~$U_{0}$ un voisinage du point~$x$ dans~$X$ contenu dans l'int\'erieur de~$\overline{D}_{V}(t)$. La suite~$(f_{k,j})_{k\in\N}$ converge alors dans $\Os(U_{0})$. Notons \mbox{$f_{j}\in\Os(U_{0})$} sa limite. Nous avons alors 
$$F= \sum_{j=1}^q f_{j} G_{j} \textrm{ dans } \Os(U_{0}).$$
On en d\'eduit finalement que 
$$F_{x} \in\Ms.$$
\end{proof}

\subsection{Conclusion}\label{parexh}

Nous nous int\'eresserons ici \`a l'\'etude des lemniscates au-dessus de n'importe quelle partie connexe de l'espace de base~$B$. Commen\c{c}ons par \'enoncer un r\'esultat topologique. Il se d\'emontre \`a l'aide des descriptions explicites du num\'ero \ref{descriptionMA}  

\begin{lem}\label{exhbase}\index{Exhaustion!sur MA@sur $\Ms(A)$}
Toute partie connexe de l'espace~$B$ poss\`ede une exhaustion par des parties compactes et connexes.
\end{lem}

\index{Exhaustion!de Stein!pour les couronnes de A1@pour les couronnes de $\AA$|(}

Soit~$V$ une partie connexe de l'espace~$B$. Soit $(V_{n})_{n\in\N}$ une exhaustion de~$V$ par des parties compactes et connexes de l'espace~$B$.

Soient $s,t\in\of{[}{0,+\infty}{[}$, avec $s < t$. Soient $(s_{n})_{n\in\N}$ et $(t_{n})_{n\in\N}$ deux suites r\'eelles v\'erifiant les conditions suivantes :
\begin{enumerate}[\it i)]
\item la suite $(s_{n})_{n\in\N}$ est strictement d\'ecroissante et tend vers $s$ ;
\item la suite $(t_{n})_{n\in\N}$ est strictement croissante et tend vers $t$ ;
\item $s_{0} \le t_{0}$.
\end{enumerate}
Soit $(u_{n})_{n\in\N}$ une suite strictement croissante et tendant vers l'infini d'\'el\'ements de $\of{[}{s_{0},+\infty}{[}$.

Posons
$$\begin{array}{cccl}
&V^{(0)}_{s,t} &=& \left\{x\in X_{V}\, \big|\, s \le |T(x)| \le t\right\},\\
&V^{(1)}_{s,t} &=& \left\{x\in X_{V}\, \big|\, s < |T(x)| \le t\right\},\\
&V^{(2)}_{s,t} &=& \left\{x\in X_{V}\, \big|\, s \le |T(x)| < t\right\},\\
&V^{(3)}_{s,t} &=& \left\{x\in X_{V}\, \big|\, s < |T(x)| < t\right\},\\
&V^{(4)}_{s} &=& \left\{x\in X_{V}\, \big|\,  |T(x)| \ge s\right\}\\
\textrm{et} &V^{(5)}_{s} &=& \left\{x\in X_{V}\, \big|\,  |T(x)|>s\right\}.
\end{array}$$


D\'esignons par $C$ l'une de ses six parties de $X$. D\'efinissons alors une exhaustion $(C_{n})_{n\in\N}$ de~$C$ par des parties compactes en posant, pour tout $n\in\N$, 
$$\begin{array}{ccclcl}
&C_{n} &=& X_{V_{n}} \cap \overline{C}(s,t) &\textrm{si} &C=V^{(0)}_{s,t},\\
&C_{n} &=& X_{V_{n}} \cap \overline{C}(s_{n},t) &\textrm{si}& C=V^{(1)}_{s,t},\\
&C_{n} &=& X_{V_{n}} \cap \overline{C}(s,t_{n}) &\textrm{si} &C=V^{(2)}_{s,t},\\
&C_{n} &=& X_{V_{n}} \cap \overline{C}(s_{n},t_{n}) &\textrm{si}& C=V^{(3)}_{s,t},\\
&C_{n} &=& X_{V_{n}} \cap \overline{C}(s,u_{n}) &\textrm{si}& C=V^{(4)}_{s}\\
\textrm{et} & C_{n} &=& X_{V_{n}} \cap \overline{C}(s_{n},u_{n}) &\textrm{si}& C=V^{(5)}_{s}.
\end{array}$$


Nous allons montrer que l'exhaustion $(C_{n})_{n\in\N}$ est une exhaustion de Stein de~$C$ relativement \`a tout faisceau de $\Os_{C}$-modules coh\'erent. Nous savons d\'ej\`a, d'apr\`es le th\'eor\`eme \ref{lemniscateStein} que, pour tout \'el\'ement~$n$ de~$\N$, la partie~$C_{n}$ est de Stein. Fixons un faisceau de $\Os_{C}$-modules coh\'erent $\Ss$.

Il nous faut, \`a pr\'esent, d\'efinir une semi-norme sur chacune des couronnes compactes consid\'er\'ees. Soit $n\in\N$. D'apr\`es le th\'eor\`eme~A et le lemme \ref{thAcompact}, il existe un entier $l_{n}\in\N^*$ et un morphisme de $\Os_{C_{n}}$-modules surjectif
$$\alpha_{n} :\Os_{C_{n}}^{l_{n}} \to \Ss_{C_{n}}.$$
Le th\'eor\`eme~B assure qu'il induit un morphisme de $\Os(C_{n})$-modules surjectif
$$\eps_{n} : \Os(C_{n})^{l_{n}} \to \Ss(C_{n}).$$

Introduisons une notation. Pour toutes parties~$E$ et~$F$ de~$X$ v\'erifiant $E\subset F$ et tout entier positif~$l$, nous noterons~$\|.\|_{\infty,E}$ la semi-norme sur l'anneau $\Os(F)^l$ obtenue en prenant le maximum des normes uniformes sur~$E$ des coefficients.

Nous d\'efinissons alors une semi-norme $\|.\|_{n}$ sur $\Ss(C_{n})$ en posant, pour tou\-te section $s\in\Ss(C_{n})$, 
$$\|s\|_{n} = \inf\{\|t\|_{\infty,C_{n}},\, t\in\eps_{n}^{-1}(s)\}.$$

Il nous reste \`a v\'erifier que les conditions de la d\'efinition \ref{exhStein} sont satisfaites. Soit $n\in\N$. Introduisons, tout d'abord, quelques notations. Nous d\'esignerons par~$r_{n}$ et~$\rho_{n}$ les applications de restriction suivantes :
$$r_{n} : (\Os(C_{n+1})^{l_{n+1}},\|.\|_{\infty,C_{n+1}}) \to (\Os(C_{n})^{l_{n+1}},\|.\|_{\infty,C_{n}})$$
et
$$\rho_{n} : (\Ss(C_{n+1}),\|.\|_{n+1}) \to (\Ss(C_{n}),\|.\|_{n}).$$
Le morphisme~$r_{n}$ est born\'e. 


D'apr\`es le th\'eor\`eme B, le morphisme surjectif $\alpha_{n+1} : \Os_{C_{n+1}}^{l_{n+1}} \to \Ss_{C_{n+1}}$ consid\'er\'e pr\'ec\'edemment induit un morphisme surjectif
$$\eps'_{n} : \Os(C_{n})^{l_{n+1}} \to \Ss(C_{n}).$$
Nous pouvons donc d\'efinir une nouvelle semi-norme $\|.\|'_{n}$ sur $\Ss(C_{n})$ en posant, pour toute section $s\in\Ss(C_{n})$, 
$$\|s\|'_{n} = \inf\{\|t\|_{\infty,C_{n}},\, t\in{\eps'}_{n}^{-1}(s)\}.$$

Nous noterons 
$$\sigma_{n} : \Ss(C_{n}),\|.\|'_{n}) \to (\Ss(C_{n}),\|.\|_{n})$$
le morphisme identit\'e allant de l'anneau $\Ss(C_{n})$ muni de la norme $\|.\|'_{n}$ \`a l'anneau $\Ss(C_{n})$ muni de la norme $\|.\|_{n}$. 


\begin{lem}
Quel que soit $n\in\N$, il existe un morphisme born\'e 
$$\eta_{n} : \Os(C_{n})^{l_{n+1}} \to \Os(C_{n})^{l_{n}}$$ 
qui fait commuter le diagramme suivant :
$$\xymatrix{
\Os(C_{n})^{l_{n+1}} \ar[r]^{\ \eps'_{n}} \ar[d]^{\eta_{n}} & \Ss(C_{n}) \ar[d]^{\sigma_{n}} \\
\Os(C_{n})^{l_{n}} \ar[r]^{\eps_{n}}  & \Ss(C_{n})
}.$$
\end{lem}
\begin{proof}
Soit $(e_{1},\ldots,e_{l_{n+1}})$ la base canonique du $\Os(C_{n+1})$-module $\Os(C_{n+1})^{l_{n+1}}$. Quel que soit $i\in\cn{1}{l_{n+1}}$, on choisit $g_{i}\in \Os(C_{n})^{l_{n}}$ tel que
$$\eps_{n}(g_{i})=(\sigma_{n} \circ \eps'_{n})(e_{i}) \textrm{ dans } \Ss(C_{n}).$$
L'application 
$$\eta_{n} :
\begin{array}{rcl}
\Os(C_{n})^{l_{n+1}} &\to& \Os(C_{n})^{l_{n}}\\
\disp \sum\limits_{i=1}^{l_{n+1}} f_{i}\, e_{i} & \mapsto & \disp \sum\limits_{i=1}^{l_{n+1}} {f_{i}}\, g_{i}
\end{array}$$
convient. Elle fait clairement commuter le diagramme qui pr\'ec\`ede. En outre, pour tous \'el\'ements $f_{1},\ldots,f_{l_{n+1}}$ de~$\Os(C_{n})$, nous avons
$$\begin{array}{rcl}
\disp \left\|  \sum\limits_{i=1}^{l_{n+1}} {f_{i}}\, g_{i}\right\|_{\infty,C_{n}} & \le & \disp \sum\limits_{i=1}^{l_{n+1}} \|f_{i}\|_{C_{n}}\, \|g_{i}\|_{\infty,C_{n}}\\
& \le & \disp \max_{1\le i\le l_{n+1}}(\|f_{i}\|_{C_{n}})\, \sum\limits_{i=1}^{l_{n+1}} \|g_{i}\|_{\infty,C_{n}}\\
& \le & \disp \left\|  \sum\limits_{i=1}^{l_{n+1}} {f_{i}}\, e_{i}\right\|_{\infty,C_{n}}\,  \sum\limits_{i=1}^{l_{n+1}} \|g_{i}\|_{\infty,C_{n}}\\
\end{array}$$ 
\end{proof}

Finalement, nous obtenons le diagramme commutatif suivant :
$$\xymatrix{
(\Os(C_{n+1})^{l_{n+1}},\|.\|_{\infty,C_{n+1}}) \ar[r]^{\ \quad \eps_{n+1}} \ar[d]^{r_{n}} & (\Ss(C_{n+1}),\|.\|_{n+1}) \ar[d]^{\cdot_{|C_{n}}} \ar@/^4pc/[dd]^{\rho_{n}}\\
(\Os(C_{n})^{l_{n+1}},\|.\|_{\infty,C_{n}}) \ar[r]^{\quad \eps'_{n}} \ar[d]^{\eta_{n}} & (\Ss(C_{n}),\|.\|'_{n}) \ar[d]^{\sigma_{n}} \\
(\Os(C_{n})^{l_{n}},\|.\|_{\infty,C_{n}}) \ar[r]^{\quad \eps_{n}}  & (\Ss(C_{n}),\|.\|_{n}).
}$$

D\'emontrons, \`a pr\'esent, que les conditions de la d\'efinition \ref{exhStein} sont satisfaites. 

\begin{lem}\label{lemborne}
Pout tout entier positif~$n$, le morphisme~$\rho_{n}$ est born\'e.
\end{lem}
\begin{proof}
Soit $n\in\N$. Les morphismes $r_{n}$, $\eta_{n}$ et~$\eps_{n}$ sont born\'es. Par cons\'equent, il existe un nombre r\'eel~$M$ tel que, pour tout \'el\'ement~$t$ de $\Os(C_{n+1})^{l_{n+1}}$, nous ayons
$$\|\eps_{n}\circ \eta_{n}\circ r_{n} (t)\|_{n} \le M\, \|t\|_{\infty,C_{n+1}}.$$

Soit~$s$ un \'el\'ement de~$\Ss(C_{n+1})$. Soit~$\delta>0$. Il existe un \'el\'ement~$t_{\delta}$ de $\Os(C_{n+1})^{l_{n+1}}$ tel que 
$$ \|t_{\delta}\|_{\infty,C_{n+1}} \le \|s\|_{n+1} + \delta.$$
On en d\'eduit que
$${\renewcommand{\arraystretch}{1.3}\begin{array}{rcl}
\|\rho_{n}(s)\|_{n} & = & \|\eps_{n}\circ \eta_{n}\circ r_{n} (t_{\delta})\|_{n}\\
& \le & M \|t_{\delta}\|_{\infty,C_{n+1}}\\
& \le & M \|s\|_{n+1} + M\delta.
\end{array}}$$
On obtient le r\'esultat voulu en faisant tendre le nombre r\'eel~$\delta$ vers~$0$.
\end{proof}

\begin{lem}\label{lemdense}
Pout tout entier positif~$n$, l'image du morphisme $\rho_{n}$ est dense dans $\Ss(C_{n})$ pour la norme~$\|.\|_{n}$.
\end{lem}
\begin{proof}
Soit $n\in\N$. Puisque le morphisme~$\sigma_{n}$ est surjectif et born\'e, il suffit de montrer que l'image du morphisme de restriction
$$\Ss(C_{n+1}) \to \Ss(C_{n})$$
est dense pour la norme~$\|.\|'_{n}$.

Soit~$s$ un \'el\'ement de~$\Ss(C_{n})$. Soit $\delta>0$. Il existe un \'el\'ement~$t$ de $\Os(C_{n})^{l_{n+1}}$ tel que $\eps'_{n}(t)=s$. On d\'eduit du lemme \ref{lemvoisdisquedroite} et de la proposition \ref{imagecouronne} (respectivement \ref{imagedisque}) que l'anneau $A[T,T^{-1}]$ (respectivement $A[T]$) est dense dans~$\Os(C_{n})$ si \mbox{$s_{n}>0$} (respectivement $s_{n}=0$). En particulier, l'anneau $\Os(C_{n+1})$ est dense dans l'anneau $\Os(C_{n})$ et il existe un \'el\'ement~$t'$ de $\Os(C_{n+1})^{l_{n+1}}$ tel que
$$\|r_{n}(t')-t\|_{\infty,C_{n}} \le \delta.$$
Posons $s'=\eps_{n+1}(t') \in \Ss(C_{n+1})$. Nous avons alors 
$$\left\|s'_{|_{C_{n}}}-s\right\|'_{n}\le \delta.$$
\end{proof}



\begin{lem}\label{nulinterieur}
Soit $n\in\N$. Soit $s\in\Ss(C_{n+1})$ telle que $\|s\|_{n+1}=0$. Alors la section $s$ est nulle sur l'ouvert $(C_{n+1})^\circ$. En particulier, elle est nulle sur $C_{n}$.
\end{lem}
\begin{proof}
Par hypoth\`ese, il existe $t\in \eps_{n+1}^{-1}(s)$ et une suite $(t_{j})_{j\in\N}$ de $\Ker(\eps_{n+1})$ v\'erifiant 
$$\lim\limits_{j\to +\infty} \|t-t_{j}\|_{\infty,C_{n+1}}= 0.$$
En d'autres termes, la suite $(t_{j})_{j\in\N}$ converge uniform\'ement vers $t$ sur $(C_{n+1})^\circ$.

Soit $x\in(C_{n+1})^\circ$. La suite des germes $((t_{j})_{x})_{j\in\N}$ converge vers $t_{x}$ dans $\Os_{Y,x}^{l_{n+1}}$. D'apr\`es le th\'eor\`eme \ref{modfermes}, nous avons 
$$t_{x} \in \mathscr{K}er(\eps_{n+1})_{x}.$$
Par cons\'equent, $t\in\mathscr{K}er(\eps_{n+1})((C_{n+1})^\circ)$ et la section $s$ est nulle sur $(C_{n+1})^\circ$.
\end{proof}

\begin{lem}\label{Cauchyconverge}
Soit $n\in\N$. Soit $(s_{k})_{k\in\N}$ une suite d'\'el\'ements de~$\Ss(C_{n+1})$ qui est de Cauchy pour la semi-norme~$\|.\|_{n+1}$. Il existe un \'el\'ement~$s$ de~$\Ss(C_{n})$ tel que la suite $(\rho_{n}(s_{k}))_{k\in\N}$ converge vers~$s$ pour la semi-norme~$\|.\|_{n}$.

Si~$s'$ est une limite de la suite $(\rho_{n}(s_{k}))_{k\in\N}$ dans~$\Ss(C_{n})$, alors elle co\"{\i}ncide avec l'\'el\'ement~$s$ sur l'ouvert $(C_{n})^\circ$.
\end{lem}
\begin{proof}
Il existe une application $\alpha : \N \to \N$ strictement croissante telle que
$$\forall k\in\N,\, \|s_{\alpha(k)}-s_{\alpha(k+1)}\|_{n+1} \le \frac{1}{2^{k+1}}.$$
Pour tout entier positif~$k$, choisissons un \'el\'ement~$d_{k}$ de $\Os(C_{n+1})^{l_{n+1}}$ qui rel\`eve l'\'el\'ement $s_{\alpha(k)}-s_{\alpha(k+1)}$ de~$\Ss(C_{n+1})$ et v\'erifie
$$\|d_{k}\|_{\infty,C_{n+1}} \le \frac{1}{2^k}.$$
Choisissons \'egalement un \'el\'ement~$t_{0}$ de $\Os(C_{n+1})^{l_{n+1}}$ qui rel\`eve~$s_{0}$. Pour tout entier positif~$k$, posons
$$t_{k} = t_{0} + \sum_{l=0}^k d_{l}.$$
C'est un \'el\'ement de $\Os(C_{n+1})^{l_{n+1}}$ qui rel\`eve~$s_{k}$. On v\'erifie ais\'ement que la suite $(t_{k})_{k\in\N}$ est une suite de Cauchy de $\Os(C_{n+1})^{l_{n+1}}$. Puisque la couronne~$C_{n}$ est contenue dans l'int\'erieur de~$C_{n+1}$, la suite $(r_{n}(t_{k}))_{k\in\N}$ converge dans $\Os(C_{n})^{l_{n+1}}$. Notons~$t$ sa limite. 

Puisque les morphismes~$\eta_{n}$ et~$\eps_{n}$ sont born\'es, la suite $(\eps_{n}(\eta_{n}(r_{n}(t_{k}))))_{k\in\N}$ de~$\Ss(C_{n})$ converge vers $s=\eps_{n}(\eta_{n}(t))$. Or, pour tout entier positif~$k$, nous avons
$$\eps_{n}(\eta_{n}(r_{n}(t_{k}))) = \rho_{n}(s_{\alpha(k)}).$$
Par cons\'equent, la suite $(\rho_{n}(s_{k}))_{k\in\N}$ de $\Ss(C_{n})$ poss\`ede une valeur d'adh\'erence. Puisque le morphisme~$\rho_{n}$ est born\'e et que la suite $(s_{k})_{k\in\N}$ est de Cauchy, la suite $(\rho_{n}(s_{k}))_{k\in\N}$ l'est encore. On en d\'eduit qu'elle converge vers~$s$.

Soit~$s'$ une limite de la suite $(\rho_{n}(s_{k}))_{k\in\N}$ dans~$\Ss(C_{n})$. Nous avons $\|s'-s\|_{n}=0$. D'apr\`es le lemme \ref{nulinterieur}, les \'el\'ements~$s$ et~$s'$ co\"{\i}ncident sur l'ouvert $(C_{n})^\circ$.
\end{proof}

\begin{lem}
Soit $n\in\N$. L'image du morphisme
$$\Ss(C) \to  \Ss(C_{n})$$ 
est dense pour la semi-norme~$\|.\|_{n}$.
\end{lem}
\begin{proof}
D'apr\`es le lemme \ref{lemborne}, pour tout entier $k\ge n$, il existe $M_{k}\ge 1$ tel que, pour tout \'el\'ement~$t$ de~$\Ss(C_{k+1})$, nous ayons
$$\left\|t_{|C_{k}}\right\|_{k} \le M_{k}\, \|t\|_{k+1}.$$


Soit~$s$ un \'el\'ement de~$\Ss(C_{n})$. Soit $\delta>0$. Choisissons une suite $(\delta_{k})_{k\ge n}$ d'\'el\'ements de $\R_{+}^*$ telle que
$$\sum_{k\ge n} \left(\prod_{i=n}^{k-1} M_{i}\right) \delta_{k} \le \delta.$$
En utilisant le lemme \ref{lemdense}, on montre, par r\'ecurrence, qu'il existe un \'el\'ement 
$$(s_{k})_{k\ge n} \in \prod_{k\ge n} \Ss(C_{k})$$
tel que $s_{n}=s$ et, pour tout entier $k\ge n$,
$$\left\|{s_{k+1}}_{|C_{k}} - s_{k}\right\|_{k} \le \delta_{k}.$$

Soit~$k\ge n$. Pour tout entier $l\ge k$, nous avons
$$\|{s_{l+1}}_{|C_{k}} - {s_{l}}_{|C_{k}}\|_{k} = \|({s_{l+1}}_{|C_{l}}-s_{l})_{|C_{k}}\|_{k} \le \left(\prod_{i=k}^{l-1} M_{i}\right) \delta_{l} \le \left(\prod_{i=n}^{l-1} M_{i}\right) \delta_{l}.$$
On en d\'eduit que la suite $({s_{l}}_{|C_{k}})_{l\ge k}$ de $\Ss(C_{k})$ est de Cauchy. D'apr\`es le lemme \ref{Cauchyconverge}, elle poss\`ede une limite~$t_{k}$ dans~$\Ss(C_{k})$.

Soient~$k_{1}$ et~$k_{2}$ deux entiers v\'erifiant $k_{1}\ge k_{2}\ge n$. Puisque le morphisme de restriction de $\Ss(C_{k_{1}})$ \`a $\Ss(C_{k_{2}})$ est born\'e, l'\'el\'ement ${t_{k_{2}}}_{|C_{k_{1}}}$ de $\Ss(C_{k_{1}}$ est une limite de la suite $({s_{l}}_{|C_{k_{1}}})_{l\ge k_{2}}$. D'apr\`es le lemme \ref{Cauchyconverge}, les \'el\'ements ${t_{k_{2}}}_{|C_{k_{1}}}$ et $t_{k_{1}}$ co\"{\i}ncident sur $(C_{k_{1}})^\circ$. Puisque $(C_{k})_{k\ge n}$ est une exhaustion de~$C$, la famille $(t_{k})_{k\ge n}$ d\'etermine une section~$t$ de~$\Ss(C)$. 

Pour tout entier $k\ge n$, nous avons
$$t_{|C_{n}} - s = t_{|C_{n}} - {s_{k+1}}_{|C_{n}} + \sum_{l=n}^k ({s_{l+1}}_{|C_{n}} - {s_{l}}_{|C_{n}}).$$
Par cons\'equent, pour tout entier $k\ge n$, nous avons
$$\|t_{|C_{n}}-s\|_{n} \le \|t_{|C_{n}} - {s_{k+1}}_{|C_{n}}\|_{n} +  \sum_{l=n}^k \left(\prod_{i=n}^{l-1} M_{i}\right) \delta_{l}.$$
En faisant tendre $k$ vers l'infini, nous obtenons
$$\|t_{|C_{n}}-s\|_{n} \le \delta.$$
Cela termine la d\'emonstration.
\end{proof}

Les r\'esultats des quatre lemmes qui pr\'ec\`edent correspondent aux quatre conditions requises pour que l'exhaustion $(C_{n})_{n\in\N}$ soit une exhaustion de Stein relativement au faisceau~$\Ss$ (\emph{cf.} d\'efinition \ref{exhStein}). Nous avons donc d\'emontr\'e le r\'esultat suivant.

\begin{thm}
La suite $(C_{n})_{n\in\N}$ est une exhaustion de Stein de la couronne~$C$, relativement \`a tout faisceau de $\Os_{C}$-modules coh\'erent.
\end{thm}

\index{Exhaustion!de Stein!pour les couronnes de A1@pour les couronnes de $\AA$|)}

Le th\'eor\`eme~\ref{csqexh} nous permet alors d'en d\'eduire le r\'esultat voulu.

\begin{thm}\label{ABouvert}\index{Espace de Stein!couronne quelconque de A1@couronne quelconque de $\AA$}
La couronne $C$ est une partie de Stein de la droite analytique~$X$.
\end{thm}

\`A l'aide des r\'esultats sur les morphismes finis que nous avons obtenus, nous pouvons d\'eduire que d'autres parties de la droite analytique~$X$ sont de Stein. Regroupons ces r\'esultats dans le th\'eor\`eme qui suit.

\begin{thm}\label{lemniscateStein}\index{Espace de Stein!lemniscate de A1@lemniscate de $\AA$}
Soit~$V$ une partie connexe de l'espace~$B$. Soient~$s$ et~$t$ deux nombres r\'eels tels que $0\le s\le t$. Soit~$P$ un polyn\^ome \`a coefficients dans~$\Os(V)$ dont le coefficient dominant est inversible. Les parties suivantes de la droite analytique~$X$ sont des espaces de Stein :
\begin{enumerate}[\it i)]
\item $\left\{x\in X_{V}\, \big|\, s \le |P(T)(x)| \le t\right\}$ ;
\item $\left\{x\in X_{V}\, \big|\, s \le |P(T)(x)| < t\right\}$ ;
\item $\left\{x\in X_{V}\, \big|\, s < |P(T)(x)| \le t\right\}$ ;
\item $\left\{x\in X_{V}\, \big|\, s < |P(T)(x)| < t\right\}$ ;
\item $\left\{x\in X_{V}\, \big|\, |P(T)(x)| \ge s\right\}$ ;
\item $\left\{x\in X_{V}\, \big|\, |P(T)(x)| > s\right\}$.
\end{enumerate}
\end{thm}
\begin{proof}
Comme expliqu\'e au num\'ero \ref{endodroite}, le morphime
$$\Os(V)[T] \to \Os(V)[T,S]/(P(S)-T) \xrightarrow[]{\sim} \Os(U)[S]$$
induit un morphisme
$$\varphi : X_{U} \to X_{U}.$$
Chacune des parties qui figure dans l'\'enonc\'e est l'image r\'eciproque d'une couronne par ce morphisme. D'apr\`es le corollaire \ref{ABouvert}, les couronnes sont des espaces de Stein. Nous pouvons donc conclure en utilisant le th\'eor\`eme \ref{finiStein}. Les hypoth\`eses en sont v\'erifi\'ees d'apr\`es la proposition \ref{phifini2}, le corollaire \ref{thfinicor2}, le corollaire \ref{XIG}, la proposition \ref{XS} et le th\'eor\`eme \ref{coherence}.

\end{proof}

%% file: K0+-.pstex_t
\begin{picture}(0,0)%
\includegraphics{K0+-.pstex}%
\end{picture}%
\setlength{\unitlength}{3947sp}%
\begingroup\makeatletter\ifx\SetFigFont\undefined%
\gdef\SetFigFont#1#2#3#4#5{%
  \reset@font\fontsize{#1}{#2pt}%
  \fontfamily{#3}\fontseries{#4}\fontshape{#5}%
  \selectfont}%
\fi\endgroup%
\begin{picture}(3969,4146)(4130,-6056)
\put(7875,-2439){\makebox(0,0)[lb]{\smash{{\SetFigFont{12}{14.4}{\rmdefault}{\mddefault}{\updefault}{\color[rgb]{0,0,0}$B=\Ms(A)$}%
}}}}
\put(7295,-3279){\makebox(0,0)[lb]{\smash{{\SetFigFont{12}{14.4}{\rmdefault}{\mddefault}{\updefault}{\color[rgb]{0,0,0}$K_0^+$}%
}}}}
\put(4145,-4669){\makebox(0,0)[lb]{\smash{{\SetFigFont{12}{14.4}{\rmdefault}{\mddefault}{\updefault}{\color[rgb]{0,0,0}$K_0^-$}%
}}}}
\put(4265,-4154){\makebox(0,0)[lb]{\smash{{\SetFigFont{12}{14.4}{\rmdefault}{\mddefault}{\updefault}{\color[rgb]{0,0,0}$a_\sigma^{l(\sigma)}$}%
}}}}
\put(5145,-4059){\makebox(0,0)[lb]{\smash{{\SetFigFont{12}{14.4}{\rmdefault}{\mddefault}{\updefault}{\color[rgb]{0,0,0}$a_\sigma^u$}%
}}}}
\put(6061,-4001){\makebox(0,0)[lb]{\smash{{\SetFigFont{12}{14.4}{\rmdefault}{\mddefault}{\updefault}{\color[rgb]{0,0,0}$a_0$}%
}}}}
\end{picture}%

%% file: applications.tex
\chapter{Applications}\label{chapitreapplications}

Dans ce chapitre, nous exposons quelques r\'esultats sur les s\'eries arithm\'etiques convergentes. Rappelons que nous d\'esignons par cette expression les s\'eries \`a coefficients dans un anneau d'entiers de corps de nombres, \'eventuellement localis\'e par une partie multiplicative finiment engendr\'ee, qui poss\`edent un rayon de convergence strictement positif en toute place. Nous allons montrer que les th\'eor\`emes g\'eom\'etriques que nous avons obtenus jusqu'ici peuvent \^etre appliqu\'es \`a leur \'etude. 

Nous consacrons le num\'ero \ref{pdca} aux probl\`emes de Cousin. Rappelons que le probl\`eme de Cousin multiplicatif consiste \`a prescrire l'ordre des z\'eros et des p\^oles d'une fonction m\'eromorphe et que le probl\`eme de Cousin additif consiste \`a prescrire ses parties principales (c'est-\`a-dire ses parties non holomorphes). En g\'eom\'etrie analytique complexe, l'origine de ces questions remonte au $\textrm{XIX}^\textrm{\`eme}$ si\`ecle. Elle sont, d\'esormais, bien comprises et la th\'eorie des espaces de Stein permet de leur apporter une solution \'el\'egante. Pour plus de pr\'ecisions, l'on consultera avec profit le deuxi\`eme paragraphe du chapitre~V de l'ouvrage~\cite{GR} de H. Grauert et R. Remmert.

Au num\'ero \ref{ndadsa}, nous nous int\'eresserons \`a la noeth\'erianit\'e de certains anneaux de s\'eries arithm\'etiques convergentes. Pour tout nombre r\'eel positif~$r$, notons
$\Zr$
l'anneau form\'e des s\'eries en une variable \`a coefficients entiers dont le rayon de convergence complexe est strictement sup\'erieur \`a~$r$. Dans l'article \cite{Harbater}, D.~Harbater d\'emontre, par une preuve purement alg\'ebrique, que, pour tout nombre r\'eel positif~$r$, l'anneau $\Zr$ est noeth\'erien (\emph{cf.}~th\'eor\`eme~1.8). En g\'eom\'etrie analytique complexe, on trouve un r\'esultat analogue dans l'article~\cite{Frisch} de J.~Frisch, qui sera ensuite pr\'ecis\'e par~Y.-T.~Siu, dans~\cite{Siu}. Nous adapterons leur m\'ethode, tr\`es g\'eom\'etrique, dans le cadre de la droite analytique au-dessus d'un anneau d'entiers de corps de nombres afin d'\'etendre le r\'esultat de D.~Harbater.

Finalement, au num\'ero \ref{pdgi}, nous proposons une nouvelle d\'emonstration d'un r\'esultat de D.~Harbater li\'e au probl\`eme de Galois inverse. Notons
$\Zun$
l'anneau form\'e des s\'eries en une variable \`a coefficients entiers dont le rayon de convergence complexe est sup\'erieur ou \'egal \`a~$1$. Le corollaire~3.8 de l'article \cite{galoiscovers} assure que tout groupe fini est le groupe de Galois d'une extension du corps $\Frac(\Z_{1^-}[\![T]\!])$. Nous proposons une d\'emonstration g\'eom\'etrique et conceptuellement tr\`es simple de ce r\'esultat. 

\bigskip

De nouveau, nous reprenons les notations du chapitre \ref{chapitredroite}.
\newcounter{nodapp}\setcounter{nodapp}{\thepage}

\section{Probl\`emes de Cousin arithm\'etiques}\label{pdca}

Dans cette partie, nous nous int\'eresserons aux probl\`emes de Cousin pour les anneaux de s\'eries arithm\'etiques.

Nous allons nous int\'eresser \`a ces probl\`emes sur la droite analytique~\mbox{$X=\E{1}{A}$} au-dessus de~$B=\Ms(A)$. Puisque les seules fonctions m\'eromorphes sur~$X$ sont les fractions rationnelles (\emph{cf}.~corollaire~\ref{fnmero}), nous n'\'etudierons pas v\'eritablement les probl\`emes de Cousin sur l'espace~$X$, mais nous restreindrons au disque unit\'e ouvert de rayon~$1$. \`A cet effet, nous utiliserons les r\'esultats obtenus au chapitre pr\'ec\'edent sur les sous-espaces de Stein de $X$. Signalons que les d\'emonstrations que nous proposons pr\'esentent encore des similitudes frappantes avec celles de la g\'eom\'etrie analytique complexe. 

Fixons quelques notations. Posons
$$\D = \mathring{D}(0,1) = \left\{x\in X\, \big|\, |T(x)|<1\right\}$$
et, quel que soit~$\sigma\in\Sigma$,
$$\D_{a_{\sigma}} = \left\{ x\in X_{a_{\sigma}}\, \big|\, |T(x)|<1\right\}.$$
\newcounter{D}\setcounter{D}{\thepage}

\subsection{Probl\`eme de Cousin multiplicatif}

Annon\c{c}ons tout de suite un r\'esultat n\'egatif : le probl\`eme de Cousin multiplicatif n'admet pas toujours de solution sur le disque~$\D$, c'est-\`a-dire qu'il existe un diviseur qui ne provient d'aucune fonction m\'eromorphe. En fait, tel est d\'ej\`a le cas sur un corps ultram\'etrique, d\`es que celui-ci n'est pas maximalement complet. Ce r\'esultat est d\^u \`a M. Lazard (\emph{cf.} \cite{Lazard}, proposition~6).\index{Probleme de Cousin@Probl\`eme de Cousin!multiplicatif!sur un corps ultram\'etrique} Fixer les ordres des z\'eros est donc impossible, mais nous allons montrer que nous pouvons les minorer.

\begin{defi}\index{Ordre d'annulation en un point rigide}
Soit~$x$ un point rigide de~$\D_{a_{\sigma}}$. Notons~$p_{x}\in\hat{K}_{\m}[T]$ le polyn\^ome irr\'eductible et unitaire qui lui est associ\'e. L'anneau local~$\Os_{\D,x}$ est alors un anneau de valuation discr\`ete dont~$p_{x}$ est une uniformisante. Soient~$f$ une fonction d\'efinie sur un voisinage du point~$x$ et~$n$ un entier. Nous dirons que la fonction~$f$ {\bf s'annule \`a l'ordre~$n$ en~$x$} si~$p_{x}^n$ divise~$f$ dans l'anneau local~$\Os_{X,x}$.
\end{defi}

Introduisons une autre d\'efinition afin de pr\'eciser sous quelles conditions nous entendons prescrire les ordres d'annulation. 


\begin{defi}\index{Distribution!d'ordres}\newcounter{o}\setcounter{o}{\thepage}
Une {\bf distribution d'ordres~$o$ sur $\D$} est la donn\'ee de 
\begin{enumerate}[\it i)]
\item un sous-ensemble fini $\Sigma_{o}$ de $\Sigma$ ; 
\item pour tout $\sigma\in\Sigma_{o}$, un sous-ensemble $E_{\sigma}$ de points rigides de $\D_{a_{\sigma}}$ ; 
\item pour tout $\sigma\in\Sigma_{o}$ et tout point~$e\in E_{\sigma}$, un nombre entier~$n_{e}$ 
\end{enumerate}
v\'erifiant la condition suivante : quel que soit~$\sigma\in\Sigma_{o}$, l'ensemble~$E_{\sigma}$ est ferm\'e, discret et ne contient pas le point~$0$.
\end{defi}

\`A toute distribution d'ordres est donc associ\'e un diviseur de Cartier sur le disque ouvert analytique~$\D_{a_{\sigma}}$. Il est presque imm\'ediat que ce diviseur s'\'etend en un diviseur de Cartier sur~$\D \cap X'_{\sigma}$. Pour l'\'etendre \'egalement \`a la fibre centrale, nous utiliserons le r\'esultat topologique qui suit.
\begin{lem}\label{topodiscret}
Soient~$\sigma\in\Sigma$, $I$ un ensemble, $\Pi=(P_{i})_{i\in I}$ une famille de polyn\^omes \`a coefficients dans~$\hat{K}_{\sigma}$, deux \`a deux distincts, irr\'eductibles et unitaires et~$(x_{i})_{i\in I}$ la famille de points rigides de~$X_{a_{\sigma}}$ associ\'ee. Supposons que l'ensemble~$E$ des points~$x_{i}$, avec~$i\in I$, soit contenu dans~$\D_{a_{\sigma}}$, ferm\'e et discret dans~$\D_{a_{\sigma}}$ et \'evite le point~$0$. Alors la partie
$$V_{\Pi} = \bigcup_{i\in I} \{ y\in X'_{\sigma}\, |\, P_{i}(y)=0\}$$
est ferm\'ee dans~$(X'_{\sigma}\cup X_{0})\cap \D$. 
\end{lem}
\begin{proof}
Nous allons montrer que le compl\'ementaire~$U$ de~$V_{\Pi}$ dans la partie~$(X'_{\sigma}\cup X_{0})\cap \D$ est ouvert. Par hypoth\`ese, la partie~$U\cap \D_{a_{\sigma}}$ est ouverte. La structure de produit de~$X'_{\sigma}$ (\emph{cf.}~propositions~\ref{produitum} et~\ref{produitarc}) nous permet d'en d\'eduire que la partie~$U\cap X'_{\sigma}$ est encore ouverte. 

Soit~$y$ un point de~$U\cap X_{0}=\D\cap X_{0}$. Il existe un \'el\'ement~$r$ de~$\of{[}{0,1}{[}$ tel que~$y$ soit le point~$\eta_{r}$ de la fibre centrale~$X_{0}$. Puisque la partie~$E$ du disque~$\D_{a_{\sigma}}$ est ferm\'ee et ne contient pas $0$, il existe $t>0$ v\'erifiant 
$$\left\{z\in E\, \big|\, |T(z)|<t\right\} = \emptyset.$$
Par cons\'equent, la partie 
$$\bigcup_{0<\eps\le 1} \{z\in X_{a_{\sigma}^\eps}\,|\, |T(z)|<t^\eps\}$$
ne coupe pas $V_{\Pi}$. 

Soit $s\in\of{]}{r,1}{[}$. Il existe~$\alpha\in\of{]}{0,1}{]}$ tel que $t^\alpha>s$. La partie d\'efinie par 
$$V=\{z\in \pi^{-1}(\of{[}{a_{0},a_{\sigma}^\alpha}{[})\,|\, |T(z)|<s\}$$
est un voisinage de~$y$ dans~$X_{\sigma}$. Observons qu'elle ne coupe pas~$V_{\Pi}$. En effet, la partie~$V_{\Pi}$ ne coupe pas la fibre centrale~$X_{0}$ et ne coupe pas non plus $V\cap X'_{\sigma}$, par choix de $s$. Finalement, nous avons bien montr\'e que la partie~$V_{\Pi}$ est ferm\'ee dans~$(X'_{\sigma}\cup X_{0})\cap \D$.
\end{proof}

Soit~$o$ une distribution d'ordres sur~$\D$. Pour montrer qu'il existe une fonction analytique qui poss\`ede des z\'eros d'ordre sup\'erieur \`a ceux prescrits par~$o$, nous allons commencer par interpr\'eter une telle fonction comme une section d'un faisceau. \`A cet effet, construisons explicitement le diviseur de Cartier mentionn\'e plus haut. Plus pr\'ecis\'ement, nous allons associer \`a la distribution d'ordres~$o$ un sous-faisceau inversible~$\Ss_{o}$ de~$\Os$ sur l'espace
$$\D_{o}=\D\setminus\left(\bigcup_{\m\in\Sigma_{o}\cap\Sigma_{f}} \tilde{X}_{\m}\right).$$
Soient~$\sigma\in\Sigma_{o}$. Pour chaque \'el\'ement~$e$ de~$E_{\sigma}$, choisissons un voisinage ouvert~$U_{e}$ du point~$e$ dans~$\D_{a_{\sigma}}$ et \'evitant le point~$0$. Quitte \`a restreindre ces ouverts, nous pouvons supposer qu'ils sont deux \`a deux disjoints. Soit~$e\in E_{\sigma}$.  Notons~$p_{e}$ le polyn\^ome \`a coefficients dans~$\hat{K}_{\sigma}$, irr\'eductible et unitaire associ\'e \`a ce point. L'image de l'ouvert~$U_{e}$ par le flot,
$$V_{e} = \bigcup_{y\in U_{e}} T_{X}(y),$$
est un voisinage ouvert dans~$\D_{o}$ du ferm\'e de Zariski
$$Z_{e} = \{y\in X'_{\sigma}\, |\, p_{e}(y)=0\}.$$
Pour~$f\in E_{\sigma}\setminus\{e\}$, les ouverts~$V_{e}$ et~$V_{f}$ sont disjoints. D\'efinissons le faisceau~$\Ss_{o}$ sur l'ouvert~$V_{e}$ par
$${\Ss_{o}}_{|V_{e}} = p_{e}^{n_{e}}\, \Os_{|V_{e}}.$$

D'apr\`es le lemme~\ref{topodiscret}, la partie
$$U=\D_{o}\setminus\left(\bigcup_{\sigma\in\Sigma_{o}, e\in E_{\sigma}} Z_{e}\right)$$
est ouverte. Nous y d\'efinissons le faisceau $\Ss_{o}$ par
$${\Ss_{o}}_{|U} = \Os_{|U}.$$
On v\'erifie sans peine que cette d\'efinition est coh\'erente avec les pr\'ec\'edentes et que le faisceau~$\Ss_{0}$ ainsi construit est un sous-faisceau inversible de~$\Os_{|\D_{o}}$.

\begin{thm}\label{zeromin}\index{Probleme de Cousin@Probl\`eme de Cousin!multiplicatif!sur A1@sur $\AA$}\index{Espace de Stein!couronne quelconque de A1@couronne quelconque de $\AA$}
Soit~$o$ une distribution d'ordres sur~$\D$. Alors il existe une fonction~$\varphi$ holomorphe sur~$\D_{o}$ et non nulle v\'erifiant la condition suivante : quel que soient~$\sigma\in\Sigma_{o}$ et~$e\in E_{\sigma}$, la fonction~$\varphi$ s'annule au point~$e$ \`a un ordre sup\'erieur \`a~$n_{e}$.
\end{thm}
\begin{proof}
Le faisceau~$\Ss_{o}$ construit pr\'ec\'edemment est inversible et donc coh\'erent. D'apr\`es le th\'eor\`eme \ref{ABouvert}, ce faisceau satisfait le th\'eor\`eme~A sur~$\D_{o}$. On en d\'eduit qu'il existe une section globale non nulle~$\varphi$ du faisceau~$\Ss_{o}$ sur~$\D_{o}$. Cette fonction convient.
\end{proof}


\subsection{Probl\`eme de Cousin additif}

\index{Probleme de Cousin@Probl\`eme de Cousin!additif!sur C@sur $\C$}\index{Theoreme@Th\'eor\`eme!de Mittag-Leffler}\index{Mittag-Leffler|see{Th\'eor\`eme de Mittag-Leffler}}

Soient~$F$ un ensemble ferm\'e et discret de points de~$\C$ et~$(R_{f})_{f\in F}$ une famille de polyn\^omes \`a coefficients dans~$\C$ sans terme constant. En g\'eom\'etrie analytique complexe, la r\'esolution du probl\`eme de Cousin additif sur~$\C$, appel\'e encore th\'eor\`eme de Mittag-Leffler, nous assure qu'il existe une fonction m\'eromorphe~$\varphi$ sur~$\C$ v\'erifiant les propri\'et\'es suivantes :
\begin{enumerate}[\it i)]
\item la fonction $\varphi$ est holomorphe sur $\C\setminus F$ ;
\item pour tout point~$f$ de~$F$, nous avons $\varphi(z) - R\left(\frac{1}{z-f}\right)$ dans~$\Os_{\C,f}$.
\end{enumerate}

Comme pr\'ec\'edemment, nous allons chercher \`a adapter ce r\'esultat pour des fonctions m\'eromorphes sur le disque unit\'e ouvert~$\D$. Rappelons que nous avons introduit le faisceau des fonctions m\'eromorphes \`a la d\'efinition \ref{defimero}.
Commen\c{c}ons par une nouvelle d\'efinition.

\begin{defi}\index{Parties principales}\index{Faisceau!des parties principales}\newcounter{Ps}\setcounter{Ps}{\thepage}
Le faisceau quotient 
$$\Ps=\Ms/\Os$$
est appel\'e {\bf faiceau des parties principales} sur $X$.
\end{defi}

Par construction, nous disposons de la suite exacte courte 
$$0\to \Os \to \Ms \to \Ps\to 0.$$
Soit $U$ un ouvert de $X$. La suite exacte longue de cohomologie associ\'ee commence comme suit :
$$0\to \Os(U) \to \Ms(U) \to \Ps(U) \to H^1(U,\Os) \to \cdots$$
En particulier, si le groupe $H^1(U,\Os)$ est nul, alors l'application canonique 
$$\Ms(U) \to \Ps(U)$$
est surjective. Cette simple remarque permet de d\'emontrer le th\'eor\`eme de Mittag-Leffler en l'appliquant avec $U=\C$. Nous allons adopter la m\^eme d\'emarche pour apporter une solution au probl\`eme de Cousin additif sur l'espace analytique~$X$.

Soit $\sigma\in\Sigma$. Fixons une cl\^oture alg\'ebrique~$L_{\sigma}$ de~$\hat{K}_{\sigma}$. Soit~$x$ un point rigide de~$\D_{a_{\sigma}}$. Le th\'eor\`eme~\ref{isorig} assure qu'il existe un \'el\'ement~$\alpha(x)$ de~$L_{\sigma}$ tel que l'on ait un isomorphisme
$$\hat{K}_{\sigma}(\alpha(x)) \simeq \Hs(x)$$
et un voisinage~$U'_{x}$ du point rationnel~$\alpha(x)$ de $\E{1}{\Hs(x)}$ tel que le morphisme naturel
$$u_{x} : U'_{x} \to \E{1}{\hat{K}_{\sigma}}$$  
induise un isomorphisme sur son image~$U_{x}$. En particulier, nous avons un isomorphisme
$$v_{x} : \Os_{\E{1}{\hat{K}_{\sigma}},x} \xrightarrow[]{\sim} \Os_{\E{1}{\Hs(x)},\alpha(x)}.$$

\begin{defi}\index{Distribution!de parties principales}\newcounter{p}\setcounter{p}{\thepage}
Une {\bf distribution $p$ de parties principales sur~$\D$} est la donn\'ee de 
\begin{enumerate}[\it i)]
\item un sous-ensemble fini $\Sigma_{p}$ de $\Sigma$ ;
\item pour tout $\sigma\in\Sigma_{p}$, un sous-ensemble $F_{\sigma}$ de points rigides de $\D_{a_{\sigma}}$ ;
\item pour tout $\sigma\in\Sigma_{\Delta}$ et tout point $f\in F_{\sigma}$, un \'el\'ement $R_{f}$ de $\Hs(f)[T]$ sans terme constant
\end{enumerate}
v\'erifiant la condition suivante :
quel que soit~$\sigma\in\Sigma_{p}$, l'ensemble~$F_{\sigma}$ est ferm\'e, discret et ne contient pas le point~$0$.
\end{defi}


Si~$p$ d\'esigne une distribution de parties principales sur~$\D$, nous posons
$$\D_{p}=\D\setminus\left(\bigcup_{\m\in\Sigma_{p}\cap\Sigma_{f}} \tilde{X}_{\m}\right).$$

\begin{thm}\label{thpp}\index{Probleme de Cousin@Probl\`eme de Cousin!additif!sur A1@sur $\AA$}\index{Espace de Stein!couronne quelconque de A1@couronne quelconque de $\AA$}
Soit $p$ une distribution de parties principales sur $\D$. Alors, il existe une fonction $\varphi$ m\'eromorphe sur $\D_{p}$ v\'erifiant les conditions suivantes :
\begin{enumerate}[\it i)]
\item quel que soit $\sigma\notin\Sigma_{p}$, la s\'erie $\varphi$ d\'efinit une fonction holomorphe sur $\D_{a_{\sigma}}$ ;
\item quel que soit $\sigma\in\Sigma_{p}$ la fonction $\varphi$ d\'efinit une fonction m\'eromorphe sur $\D_{a_{\sigma}}$, holomorphe sur le compl\'ementaire de $F_{\sigma}$ ;
\item quel que soient~$\sigma\in\Sigma_{p}$ et~$f \in F_{\sigma}$, nous avons
$${v_{f}}^{*}\varphi - R_{f}\left(\frac{1}{T-\alpha(f)}\right) \in \Os_{\E{1}{\Hs(f)},\alpha(f)}\ ;$$
\item $\disp \varphi \in \left( A\left[\frac{1}{\Sigma_{p}}\right] [\![T]\!]\right) \cap \Os_{\D_{a_{\sigma}},0}$. 
\end{enumerate}
\end{thm}
\begin{proof}
Nous allons associer \`a la distribution de parties principales~$p$ une section~$s_{p}$ du faisceau~$\Ps$ sur~$\D_{p}$. Soit~$\sigma\in\Sigma_{p}$. Pour chaque \'el\'ement~$f$ de~$F_{\sigma}$, nous avons d\'efini pr\'ec\'edemment un voisinage ouvert~$U_{f}$ du point~$f$ dans~$\D_{a_{\sigma}}$. Puisque la partie~$F_{\sigma}$ est discr\`ete et ne contient pas~$0$, quitte \`a restreindre ces ouverts, nous pouvons supposer qu'ils sont deux \`a deux disjoints et \'evitent le point~$0$. Soit~$f\in F_{\sigma}$. En utilisant les propositions~\ref{isolocal} et~\ref{flot}, on montre que l'isomorphisme~$u_{f}^{-1}$, d\'efini sur~$U_{f}$, se prolonge \`a l'image de l'ouvert~$U_{f}$ par le flot,
$$V_{f} = \bigcup_{y\in U_{f}} T_{X}(y).$$
C'est un voisinage ouvert dans~$\D_{p}$ du ferm\'e de Zariski
$$Z_{f} = \{y\in X'_{\sigma}\, |\, p_{f}(y)=0\}.$$
Pour~$g\in F_{\sigma}\setminus\{f\}$, les ouverts~$V_{f}$ et~$V_{g}$ sont disjoints. D\'efinissons la section~$s_{p}$ du faisceau~$\Ps$ sur l'ouvert~$V_{f}$ par
$${s_{p}}_{|V_{f}} = (u_{f}^{-1})^*\left(R_{f}\left[\frac{1}{T-\alpha(f)}\right]\right).$$

D'apr\`es le lemme~\ref{topodiscret}, la partie
$$U=\D_{p}\setminus\left(\bigcup_{\sigma\in\Sigma_{p}, f\in F_{\sigma}} Z_{f}\right)$$
est ouverte. Nous y d\'efinissons la section~$s_{p}$ par
$${s_{p}}_{|U} = 0.$$
On v\'erifie sans peine que cette d\'efinition est coh\'erente avec les pr\'ec\'edentes et que nous avons bien construit ainsi une section~$s_{p}$ de~$\Ps$ sur l'ouvert~$\D_{p}$.

D'apr\`es le th\'eor\`eme~\ref{ABouvert}, nous avons $H^1(\D_{p},\Os)=0$. On en d\'eduit que le morphisme canonique
$$\Ms(\D_{p}) \to \Ps(\D_{p})$$
est surjectif. Par cons\'equent, la section $s_{p}$ poss\`ede un ant\'ec\'edent $\varphi$ par ce morphisme. Quel que soit $\sigma\in\Sigma$, la fonction $\varphi$ d\'efinit une fonction m\'eromorphe sur~$\D_{a_{\sigma}}$ qui poss\`ede les propri\'et\'es prescrites par l'\'enonc\'e.

Remarquons \'egalement que la fonction $\varphi$ est holomorphe au voisinage de la section nulle de $\D_{p}$. On en d\'eduit que le d\'eveloppement en $0$ de $\varphi$ est \`a coefficients dans $A[1/\Sigma_{p}]$, par la proposition~\ref{imagedisque}.
\end{proof}

Sous cette forme, le r\'esultat du th\'eor\`eme peut \^etre obtenu \`a partir du r\'esultat analogue de g\'eom\'etrie analytique complexe et d'un argument d'approximation. Nous en proposons, \`a pr\'esent, un raffinement qui, \`a notre connaissance, ne peut se d\'emontrer ainsi.

\begin{thm}\label{thppo}\index{Probleme de Cousin@Probl\`eme de Cousin!multiplicatif!sur A1@sur $\AA$}\index{Probleme de Cousin@Probl\`eme de Cousin!additif!sur A1@sur $\AA$}
Soient~$o$ une distribution d'ordres sur~$\D$ et~$p$ une distribution de parties principales sur~$\D$. Supposons que, quel que soit~$\sigma\in\Sigma_{o}\cap \Sigma_{p}$, les ensembles~$E_{\sigma}$ et~$F_{\sigma}$ soient disjoints. Alors, il existe une fonction $\varphi$ m\'eromorphe sur $\D'=\D_{o}\cap \D_{p}$ v\'erifiant les conditions suivantes :
\begin{enumerate}[\it i)]
\item quel que soit $\sigma\notin\Sigma_{p}$, la s\'erie $\varphi$ d\'efinit une fonction holomorphe sur $\D_{a_{\sigma}}$ ;
\item quel que soit $\sigma\in\Sigma_{p}$ la fonction $\varphi$ d\'efinit une fonction m\'eromorphe sur $\D_{a_{\sigma}}$, holomorphe sur le compl\'ementaire de $F_{\sigma}$ ;
\item quel que soient~$\sigma\in\Sigma_{p}$ et~$f \in F_{\sigma}$, nous avons
$${v_{f}}^{*}\varphi - R_{f}\left(\frac{1}{T-\alpha(f)}\right) \in \Os_{\E{1}{\Hs(f)},\alpha(f)}\ ;$$
\item quel que soient~$\sigma\in\Sigma_{o}$ et~$e\in E_{\sigma}$, la fonction~$\varphi$ s'annule au point~$e$ \`a un ordre sup\'erieur \`a~$n_{e}$ ;
\item $\disp \varphi \in \left(A\left[\frac{1}{\Sigma_{o}\cup\Sigma_{p}}\right][\![T]\!]\right) \cap \Os_{\D_{a_{\sigma}},0}$. 
\end{enumerate}
\end{thm} 
\begin{proof}
Il suffit de reprendre la preuve du th\'eor\`eme pr\'ec\'edent en l'appliquant \`a d'autres faisceaux. Juste avant le th\'eor\`eme~\ref{zeromin}, nous avons construit un sous-faisceau~$\Ss_{o}$ de~$\Os_{|\D_{o}}$. Construisons un sous-faisceau~$\Ts_{o}$ de~$\Ms_{|\D_{o}}$ par la m\^eme m\'ethode. Reprenons les notations utilis\'ees lors de la d\'efinition du faisceau~$\Ss_{o}$. Nous pouvons, en outre, supposer que les ouverts~$U_{e}$, et donc~$V_{e}$, sont connexes. Soient~$\sigma\in\Sigma_{o}$ et~$e\in E_{\sigma}$. Notons~$S_{e}$ l'ensemble des \'el\'ements de~$\Os_{|V_{e}}$ qui ne sont pas identiquement nuls sur~$Z_{e}$. C'est une partie multiplicative de~$\Os_{|V_{e}}$. Nous posons
$${\Ts_{o}}_{|V_{e}} = p_{e}^{n_{e}}\, S_{e}^{-1} \Os_{|V_{e}}.$$
Nous posons \'egalement
$${\Ts_{o}}_{|U} = \Ms_{|U}.$$
Nous avons bien construit ainsi un sous-faisceau de~$\Ms_{|\D_{o}}$.

Le faisceau~$\Ss_{o}$ s'injecte dans ce faisceau. Nous allons, \`a pr\'esent, construire une section~$s_{p}$ du faisceau quotient~$\Ts_{o}/\Ss_{o}$ sur l'ouvert~$\D'=\D_{o}\cap \D_{p}$. Nous pouvons proc\'eder exactement comme dans la preuve du th\'eor\`eme pr\'ec\'edent. Il suffit de prendre garde \`a choisir des ouverts~$U_{f}$ qui \'evitent les points des ensembles~$E_{\sigma}$. 

Consid\'erons la suite exacte courte
$$0 \to \Ss_{o} \to \Ts_{o} \to \Ts_{o}/\Ss_{o} \to 0.$$
Le faisceau~$\Ss_{o}$ est inversible et donc coh\'erent. D'apr\`es le th\'eor\`eme~\ref{ABouvert}, nous avons donc~$H^1(\D',\Ss_{o})=0$. On en d\'eduit que le morphisme canonique
$$\Ts(\D') \to (\Ts_{o}/\Ss_{o})(\D')$$
est surjectif. Par cons\'equent, la section $s_{p}$ poss\`ede un ant\'ec\'edent $\varphi$ par ce morphisme. Cette fonction poss\`ede les propri\'et\'es requises.
\end{proof}

Nous donnerons \`a la fin de la partie suivante (\emph{cf.} corollaire \ref{corppo}) une interpr\'etation en termes de s\'eries de ce th\'eor\`eme.

\subsection{Th\'eor\`eme de Poincar\'e}

Dans la lign\'ee des probl\`emes de Cousin, le th\'eor\`eme de Poincar\'e sur~$\C$ nous assure que toute fonction m\'eromorphe s'\'ecrit {\it globalement} comme un quotient de deux fonctions holomorphes. Ici encore, les techniques des espaces de Stein s'av\`ereront utiles. 

\begin{thm}\label{thPoincare}\index{Poincare@Poincar\'e|see{Th\'eor\`eme de Poincar\'e}}\index{Fonctions meromorphes@Fonctions m\'eromorphes!theoreme de Poincare@th\'eor\`eme de Poincar\'e}\index{Theoreme@Th\'eor\`eme!de Poincare@de Poincar\'e}\index{Espace de Stein!fonctions meromorphes@fonctions m\'eromorphes}
Soit~$M$ une partie connexe et de Stein de la droite~$X$. L'anneau~$\Os(M)$ est int\`egre et le morphisme naturel
$$\Frac(\Os(M)) \to \Ms(M)$$
est un isomorphisme.
\end{thm}
\begin{proof}
Le corollaire \ref{anneauintegre} assure donc que l'anneau~$\Os(M)$ est int\`egre. Il suffit de d\'emontrer que le morphisme naturel
$$\Frac(\Os(M)) \to \Ms(M)$$
est surjectif. Soit~$h$ un \'el\'ement de~$\Ms(M)$. Le faisceau de $\Os_{M}$-modules $\Os_{M} \cap h\Os_{M}$ est coh\'erent. Puisque la fonction nulle appartient \'evidemment \`a l'image du morphisme pr\'ec\'edent, nous pouvons supposer que~$h$ n'est pas nulle. Le faisceau $\Os_{M} \cap h\Os_{M}$ n'est alors pas nul. D'apr\`es le th\'eor\`eme~A, il poss\`ede une section globale non-nulle~$f$ sur~$M$. On en d\'eduit le r\'esultat voulu.
\end{proof}

Ce th\'eor\`eme nous permet, par exemple, de d\'ecrire les fonctions m\'eromorphes sur le disque ouvert de centre~$0$ et de rayon~$1$ comme quotient de fonctions holomorphes sur ce disque. Nous allons utiliser ce r\'esultat pour donner une version explicite, c'est-\`a-dire en termes de s\'eries convergentes, du th\'eor\`eme~\ref{thppo}.

Soit $\sigma\in\Sigma$. Soient~$L_{\sigma}$ une cl\^oture alg\'ebrique de~$\hat{K}_{\sigma}$ et~$\hat{L}_{\sigma}$ son compl\'et\'e pour la valeur absolue~$|.|_{\sigma}$. Remarquons que le groupe de Galois Gal($L_{\sigma}/\hat{K}_{\sigma}$) agit sur~$\hat{L}_{\sigma}$. Pour tout \'el\'ement~$x$ de~$L_{\sigma}$, nous noterons~$p_{x}$ le polyn\^ome minimal unitaire de~$x$ sur~$\hat{K}_{\sigma}$. Nous noterons \'egalement 
$$L_{\sigma}^{\circ\circ} = \left\{ x \in L_{\sigma}\, \big|\, |x|_{\sigma} < 1 \right\}.$$


Rappelons, finalement, que l'on peut interpr\'eter les fonctions holomorphes sur~$\E{1}{\hat{K}_{\sigma}}$ comme des fonctions holomorphes sur $\E{1}{\hat{L}_{\sigma}}$ invariantes par le groupe de Galois Gal($L_{\sigma}/\hat{K}_{\sigma}$). 

\begin{cor}\label{corppo}\index{Probleme de Cousin@Probl\`eme de Cousin!multiplicatif!sur A1@sur $\AA$}\index{Probleme de Cousin@Probl\`eme de Cousin!additif!sur A1@sur $\AA$}
Soit~$\Sigma_{\Delta}$ une partie finie de~$\Sigma$. Pour $\sigma\in\Sigma_{\Delta}$, soient~$E_{\sigma}$ et~$F_{\sigma}$ deux sous-ensembles de~$L_{\sigma}^{\circ\circ}$ disjoints, ferm\'es, discrets et \'evitant~$0$. Pour~$\sigma\in\Sigma_{\Delta}$ et~$e\in E_{\sigma}$, soit~$n_{e}$ un entier. Pour~$\sigma\in\Sigma_{\Delta}$ et~$f\in F_{\sigma}$, soit~$R_{f}$ un polyn\^ome \`a coefficients dans~$\Hs(f)$ sans terme constant. Supposons que 
\begin{enumerate}[\it i)]
\item quel que soient~$\sigma\in\Sigma_{\Delta}$, $e\in E_{\sigma}$ et~$\tau \in \textrm{Gal}(L_{\sigma}/\hat{K}_{\sigma})$, nous avons
$$\tau(e)\in E_{\sigma} \textrm{ et } n_{\tau(e)} = n_{e}\ ;$$
\item quel que soient~$\sigma\in\Sigma_{\Delta}$, $f\in F_{\sigma}$ et~$\tau \in \textrm{Gal}(L_{\sigma}/\hat{K}_{\sigma})$, nous avons
$$\tau(f)\in F_{\sigma} \textrm{ et } R_{\tau(f)} = \tau(R_{f}).$$
\end{enumerate}
Alors, il existe deux s\'eries $u,v \in A[1/\Sigma_{\Delta}][\![T]\!]$ v\'erifiant les propri\'et\'es suivantes :
\begin{enumerate}[a)]
\item quel que soit $\sigma\notin\Sigma_{\Delta}$, la s\'erie $u/v$, vue comme fonction analytique sur $\hat{L}_{\sigma}$, est d\'eveloppable en $0$ en une s\'erie enti\`ere de rayon de convergence sup\'erieur \`a $1$  ;
\item quel que soit $\sigma\in\Sigma_{\Sigma}$ et $z\notin F_{\sigma}$, la s\'erie $u/v$, vue comme fonction analytique sur $\hat{L}_{\sigma}$, est d\'eveloppable en $z$ en une s\'erie enti\`ere de rayon de convergence strictement positif ;
\item quel que soit $\sigma\in\Sigma_{\Sigma}$ et $e\in E_{\sigma}$, la s\'erie $u/v$, vue comme fonction analytique sur $\hat{L}_{\sigma}$, s'annule en~$e$ \`a un ordre sup\'erieur \`a~$n_{e}$ ;
\item quel que soient $\sigma\in\Sigma_{\Delta}$ et $f \in F_{\sigma}$, la s\'erie $u/v$, vue comme fonction analytique sur $\hat{L}_{\sigma}$, est d\'eveloppable en $f$ en une s\'erie de Laurent de partie principale $R_{f}\left(\frac{1}{T-f}\right)$ et de rayon de convergence strictement positif.
\end{enumerate}
\end{cor}

\section{Noeth\'erianit\'e d'anneaux de s\'eries arithm\'etiques}\label{ndadsa}

\subsection{Sous-vari\'et\'es analytiques}\label{cint}

\index{Droite affine analytique!sur A@sur $A$!sous-variete analytique@sous-vari\'et\'e analytique|(}

Jusqu'ici, nous avons \'etudi\'e les propri\'et\'es de la droite analytique~$X$ ou de certaines de ces parties, comme les disques et les couronnes relatifs. Il est \'egalement naturel de s'int\'eresser aux ferm\'es analytiques de la droite~$X$, c'est-\`a-dire aux parties d\'efinies localement par l'annulation de fonctions analytiques. Nous en proposons ici une br\`eve \'etude.



\begin{defi}
Soit $U$ un ouvert de $X$. On appelle {\bf sous-vari\'et\'e analytique} de $U$ tout espace localement annel\'e de la forme 
$$(V(\Is),\Os_{U}/\Is),$$ 
o\`u $\Is$ est un faisceau d'id\'eaux de type fini de $\Os_{U}$. 
\end{defi}

\begin{rem}
Soient~$U$ un ouvert de~$X$ et~$\Is$ un faisceau d'id\'eaux de type fini de~$\Os_{U}$. L'espace topologique~$V(\Is)$ est donc ferm\'e dans~$U$. Puisque le faisceau~$\Os_{U}$ est coh\'erent, le faisceau d'id\'eaux de type fini~$\Is$ l'est \'egalement. Nous en d\'eduisons que le faisceau~$\Os_{U}/\Is$ l'est encore.
\end{rem}

\begin{defi}
Soient $U$ un ouvert de $X$ et $(Z,\Os_{Z})$ une sous-vari\'et\'e analytique de~$U$. Soit $x$ un point de $Z$. On dit que la sous-vari\'et\'e $(Z,\Os_{Z})$ est {\bf int\`egre en $x$} si l'anneau local $\Os_{Z,x}$ est int\`egre. On dit que la sous-vari\'et\'e $(Z,\Os_{Z})$ est {\bf int\`egre} si elle est int\`egre en chacun de ses points.
\end{defi}

Nous allons, \`a pr\'esent, d\'ecrire les germes de sous-vari\'et\'es analytiques int\`egres en un point. Soit $x$ un point de $X$. Soient $U$ un voisinage ouvert de $x$ dans $X$ et $\Is$ un faisceau d'id\'eaux de $\Os_{U}$ tel que la sous-vari\'et\'e analytique 
$$(Z,\Os_{Z}) = (V(\Is),\Os_{U}/\Is)$$
soit int\`egre en $x$. L'id\'eal $\Is_{x}$ est donc un id\'eal premier de $\Os_{X,x}$. Nous allons distinguer plusieurs cas.

\bigskip
 
Supposons tout d'abord, que l'anneau local $\Os_{X,x}$ est un corps. L'id\'eal $\Is_{x}$ ne peut alors \^etre que l'id\'eal nul. Par le principe du prolongement analytique (\emph{cf.} th\'eor\`eme~\ref{prolan}), au voisinage du point $x$, l'id\'eal $\Is$ est nul et la sous-vari\'et\'e $(Z,\Os_{Z})$ co\"{\i}ncide avec $(X,\Os_{X})$.\\

Supposons, \`a pr\'esent, que l'anneau local $\Os_{X,x}$ est un anneau de valuation discr\`ete d'uniformisante $\tau$. L'id\'eal $\Is_{x}$ est alors soit l'id\'eal nul, soit l'id\'eal $(\tau)$. Si $\Is_{x}=(0)$, localement, la sous-vari\'et\'e $(Z,\Os_{Z})$ n'est autre que l'espace total, comme pr\'ec\'edemment. Supposons donc que $\Is_{x}=(\tau)$. D'apr\`es \ref{modulescorpsavd}, l'id\'eal $\Is$ est localement engendr\'e par~$\tau$. Distinguons de nouveau plusieurs cas. 

Supposons, tout d'abord, que le point $b=\pi(x)$ est un point interne de $B$. Il existe $\sigma \in\Sigma$ tel que ce point appartienne \`a la branche $\sigma$-adique ouverte. Il existe donc un polyn\^ome $P(T) \in \Hs(b)[T] = K_{\sigma}[T]$ irr\'eductible et unitaire tel que le point $x$ soit le point de la fibre $X_{b}$ d\'efini par l'\'equation $P(T)(x)=0$. En outre, nous pouvons supposer que $\tau=P(T)$. Notons $V$ un voisinage ouvert et connexe de $b$ dans $\pi(U)$ au-dessus duquel l'id\'eal $\Is$ est engendr\'e par $P(T)$. Nous pouvons supposer que $V$ est contenu dans la branche $\sigma$-adique ouverte. Alors l'application qui \`a tout point~$c$ de $V$ associe l'unique point $y$ de la fibre $X_{c}$ d\'efini par l'\'equation $P(T)(y)=0$ r\'ealise un hom\'eomorphisme de $V$ sur $X_{V}\cap Z$. On en d\'eduit que $X_{V}\cap Z$ est connexe et localement connexe par arcs. En outre, en tout point $y$ de $X_{V}\cap Z$, l'anneau local $\Os_{Z,y}$ est un corps. Par cons\'equent, les parties ouvertes et connexes de la sous-vari\'et\'e $X_{V}\cap Z$ v\'erifient le principe du prolongement analytique.

Supposons, \`a pr\'esent, que $b=\pi(x)$ soit le point central $a_{0}$ de $B$. Il existe encore un polyn\^ome $P(T) \in \Hs(b)[T] = K[T]$, irr\'eductible et unitaire, tel que le point~$x$ soit le point de la fibre~$X_{b}$ d\'efini par l'\'equation $P(T)(x)=0$. Nous pouvons \'egalement supposer que $\tau=P(T)$. Au voisinage de $x$, la sous-vari\'et\'e d\'efinie par l'\'equation \mbox{$P(T)=0$} est un rev\^etement topologique de $B$, ramifi\'e au point $x$. Il suffit de choisir pour voisinage de $x$ un ouvert de $X$ sur lequel~$\Is$ est engendr\'e par $P(T)$ et qui \'evite les fibres extr\^emes $\tilde{X}_{\m}$ correspondant \`a un id\'eal~$\m$ tel que le polyn\^ome $P(T)$ ait des racines multiples dans $k_{\m}$ (il n'existe qu'un nombre fini de tels id\'eaux). Comme pr\'ec\'edemment, il existe un voisinage $W$ de~$x$ dans~$U$ tel que la sous-vari\'et\'e $W\cap Z$ soit connexe, localement connexe par arcs et que ses parties ouvertes et connexes v\'erifient le principe du prolongement analytique.

Supposons, pour finir, que $b=\pi(x)$ soit un point extr\^eme de $B$. Il existe alors $\m\in\Sigma_{f}$ tel que $b=\tilde{a}_{\m}$. L'anneau local $\Os_{X,x}$ est un anneau de valuation discr\`ete si, et seulement si, le point $x$ est de type 2 ou 3. Nous pouvons alors choisir l'uniformisante $\tau=\pi_{\m}$. Par cons\'equent, au voisinage du point $x$, la sous-vari\'et\'e~$Z$ n'est autre que la fibre $\tilde{X}_{\m}$. De nouveau, nous en d\'eduisons qu'il existe un voisinage $W$ de $x$ dans $U$ tel que la sous-vari\'et\'e $W\cap Z$ soit connexe, localement connexe par arcs et que ses parties ouvertes et connexes v\'erifient le principe du prolongement analytique.

\bigskip

Il nous reste \`a traiter le cas o\`u l'anneau local $\Os_{X,x}$ n'est ni un corps, ni un anneau local. Le point $x$ est alors n\'ecessairement un point rigide d'une fibre extr\^eme : il existe $\m\in\Sigma_{f}$ et un polyn\^ome irr\'eductible et unitaire $P(T)\in k_{\m}[T]$ tel que $x$ soit l'unique point de la fibre $\tilde{X}_{\m}$ d\'efini par l'\'equation $P(T)(x)=0$. L'id\'eal maximal de $\Os_{X,x}$ est $(\pi_{\sigma},P(T))$. L'id\'eal premier $\Is_{x}$ peut \^etre de plusieurs sortes. Tout d'abord, comme dans les cas pr\'ec\'edents, nous pouvons avoir $\Is_{x}=(0)$. La sous-vari\'et\'e $Z$ co\"{\i}ncide alors localement avec l'espace $X$ tout entier. Si l'id\'eal $\Is_{x}$ est de hauteur 2, c'est l'id\'eal maximal $\m_{x}$ et la sous-vari\'et\'e~$Z$ est, localement, r\'eduite au point $x$. Si l'id\'eal $\Is_{x}$ est de hauteur 1, alors nous pouvons avoir $\Is_{x}=(\pi_{\m})$, auquel cas la sous-vari\'et\'e $Z$ co\"{\i}ncide localement avec la fibre $\tilde{X}_{\m}$, ou bien $\Is_{x}=(Q(T))$, o\`u $Q(T)$ est un polyn\^ome irr\'eductible de~$\hat{A}_{\m}[T]$ qui rel\`eve $P(T)$. Dans ce dernier cas, il est encore possible de construire une section de $\pi$ qui soit un hom\'eomorphisme d'un voisinage de $\tilde{a}_{\m}$ dans $B$ vers un voisinage de $x$ dans $Z$. Dans tous les cas, il existe un voisinage $W$ de $x$ dans $U$ tel que la sous-vari\'et\'e $W\cap Z$ soit connexe, localement connexe par arcs et que ses parties ouvertes et connexes v\'erifient le principe du prolongement analytique.

\bigskip

\`A l'aide de ces descriptions explicites, nous obtenons les r\'esultats suivants.

\begin{prop}\label{locintegre}
Soit $x$ un point de $X$. Soient $U$ un voisinage ouvert de~$x$ dans $X$ et $\Is$ un faisceau d'id\'eaux de $\Os_{U}$ tel que la sous-vari\'et\'e analytique 
$$(Z,\Os_{Z}) = (V(\Is),\Os_{U}/\Is)$$
soit int\`egre en $x$. Alors il existe un voisinage ouvert $V$ de $x$ dans $X$ tel que la sous-vari\'et\'e $Z \cap V$ de $V$ soit int\`egre. 
\end{prop}

\begin{prop}\label{propintegre}\index{Prolongement analytique!pour une sous-variete analytique integre de A1@pour une sous-vari\'et\'e analytique int\`egre de $\AA$}
Soient $U$ un ouvert de $X$ et $(Z,\Os_{Z})$ une sous-vari\'et\'e analytique int\`egre de $U$. Alors $Z$ est localement connexe par arcs et ses parties ouvertes et connexes satisfont au principe du prolongement analytique.
\end{prop}

\index{Droite affine analytique!sur A@sur $A$!sous-variete analytique@sous-vari\'et\'e analytique|)}

\subsection{Th\'eor\`eme de Frisch}

Dans ce paragraphe, nous d\'emontrons que l'anneau des germes de fonctions analytiques au voisinage de certains compacts est noeth\'erien. Le premier r\'esultat de ce type a \'et\'e obtenu par J. Frisch dans le cadre de la g\'eom\'etrie analytique complexe (\emph{cf.} \cite{Frisch}, th\'eor\`eme I, 9) :

\begin{thm*}[J.~Frisch]\index{Theoreme@Th\'eor\`eme!de Frisch!sur un corps archimedien@sur un corps archim\'edien}\index{Frisch|see{Th\'eor\`eme de Frisch}}
Soit $X$ une vari\'et\'e analytique r\'eelle ou complexe. Soit~$K$ une partie compacte de~$X$, semi-analytique et de Stein. Alors l'anneau des fonctions analytiques au voisinage de $K$ est noeth\'erien.
\end{thm*}

\begin{defi}\label{morcelable}\index{Morcelable}
Soient $E$ une partie de $X$ et $x$ un point de $E$. La partie~$E$ est dite {\bf morcelable au voisinage du point~$x$} si, pour tout voisinage ouvert~$U$ de~$x$ dans~$X$ et toute sous-vari\'et\'e analytique $Z$ de $U$ int\`egre en $x$, il existe un voisinage $V$ de $x$ dans $E\cap U$ qui poss\`ede un syst\`eme fondamental de voisinages ouverts dans $U$ dont les traces sur $Z$ sont connexes. 

La partie~$E$ est dite {\bf morcelable} si elle est morcelable au voisinage de chacun de ses points.
\end{defi}

\begin{prop}\label{stationne}\index{Faisceau!coherent@coh\'erent!sur une partie morcelable de A1@sur une partie morcelable de $\AA$}
Soit $E$ une partie morcelable de $X$. Soient $\Fs$ un faisceau coh\'erent sur $E$ et $(\Fs_{n})_{n\in\N}$ une suite croissante de sous-faisceaux coh\'erents de $\Fs$. Alors la suite $(\Fs_{n})_{n\in\N}$ est localement stationnaire dans $E$ au sens o\`u, quel que soit $x\in E$, il existe un entier $n_{0}\in\N$ et un voisinage $U$ de $x$ dans $E$ tels que 
$$\forall n\ge n_{0},\, \forall z\in U,\, (\Fs_{n_{0}})_{z} \xrightarrow{\sim} (\Fs_{n})_{z}.$$
\end{prop}
\begin{proof}
Soit $x\in E$. Il existe $n_0\in \N$ tel que, quel que soit $n\ge n_0$, on ait 
$$(\Fs_{n_0})_x\xrightarrow{\sim} (\Fs_n)_x.$$
Quitte \`a remplacer $\Fs$ par $\Fs/\Fs_{n_0}$ et $\Fs_n$ par $\Fs_n/\Fs_{n_0}$, pour $n\ge n_0$, puis \`a d\'ecaler les indices, nous pouvons supposer que $$(\Fs_n)_x = 0,$$ quel que soit $n\in\N$. Puisque $\Fs_{x}$ est un module de type fini  sur $\Os_{X,x}$, il existe un entier $r\in\N$ et une filtration
$$0=M^{0} \subset M^{1} \subset \cdots \subset M^{r} = \Fs_{x}$$
de $\Fs_{x}$ par des sous-modules de type fini et des id\'eaux premiers $\p_{0},\ldots,\p_{r}$ de~$\Os_{X,x}$ v\'erifiant la condition suivante : quel que soit $i\in\cn{0}{r-1}$, on dispose d'un isomorphisme
$$M^{(i+1)}/M^{(i)} \simeq \Os_{X,x}/\p_{i}.$$
Cette filtration et ces isomorphismes se prolongent au niveau des faisceaux. Il existe une filtration de $\Fs$
$$0=\Fs^{(0)}\subset \Fs^{(1)}\subset \cdots \subset \Fs^{(r)}=\Fs$$
par des sous-faisceaux coh\'erents d\'efinis au voisinage de $a$ et $r$ sous-vari\'et\'es analytiques $Z_0,\ldots,Z_{r-1}$ d\'efinies au voisinage de $x$, int\`egres en $x$ et v\'erifiant la condition suivante : quel que soit $i\in\cn{0}{r-1}$, on dispose d'un isomorphisme de faisceaux
$$\Fs^{(i+1)}/\Fs^{(i)} \simeq \Os_{Z_i}.$$

Il nous suffit, \`a pr\'esent, de montrer que, pour chaque $i\in\cn{0}{r-1}$, la sous-suite $(\Gs_{i,n})_{n\in\N}$ de $\Fs^{(i)}/\Fs^{(i+1)} \simeq \Os_{Z_i}$ induite par $(\Fs_n)_{n\in\N}$ stationne au voisinage de $x$ dans $E$ et m\^eme au voisinage de $x$ dans $E \cap Z_i$. Soit $U$ un voisinage ouvert de $x$ dans $X$ sur lequel $Z_{i}$ est d\'efinie. D'apr\`es la proposition \ref{locintegre}, nous pouvons supposer que $Z_{i}\cap U$ est une sous-vari\'et\'e int\`egre de $U$. Par hypoth\`ese, la partie $E$ de $X$ est morcelable au voisinage du point $x$. Il existe donc un voisinage~$V$ de $x$ dans $E\cap U$ qui poss\`ede un syst\`eme fondamental de voisinages ouverts dans $X$ dont les intersections avec $Z_{i}$ sont connexes. 

Soient $n\in\N$ et $f\in\Gs_{i,n}$. Il existe un voisinage ouvert $W$ de $V$ dans $X$ sur lequel la fonction $f$ est d\'efinie et tel que $W\cap Z_{i}$ soit une sous-vari\'et\'e int\`egre et connexe de $W$. Puisque $(\Gs_{i,n})_{x}=0$, la fonction $f$ est nulle au voisinage de $x$ dans $Z_{i}$. D'apr\`es \ref{propintegre}, $W\cap Z_{i}$ v\'erifie le principe du prolongement analytique. On en d\'eduit que $f$ est nulle sur $W\cap Z_{i}$. Finalement, le faisceau $\Gs_{i,n}$ est nul sur~$V\cap Z_{i}$, et donc sur $V$.
\end{proof}

\begin{cor}
Soient $E$ une partie de $X$ morcelable et de Stein, $\Fs$ un faisceau coh\'erent sur $E$ et $(f_{i})_{i\in I}$ une famille de sections de $\Fs$ sur $E$. Le sous-faisceau de $\Fs$ engendr\'e par la famille $(f_{i})_{i\in I}$ est coh\'erent.
\end{cor}

\begin{thm}\label{Frisch}\index{Theoreme@Th\'eor\`eme!de Frisch!sur A1@sur $\AA$}\index{Faisceau structural!germes au voisinage d'un compact}\index{Noetherianite@Noeth\'erianit\'e!des anneaux de sections sur un compact de A1@des anneaux de sections sur un compact de $\AA$}\index{Noetherianite@Noeth\'erianit\'e!des anneaux locaux|see{Anneau local en un point}}
Soient $E$ une partie compacte morcelable et de Stein de~$X$. L'anneau $\Os(E)$ des germes de fonctions analytiques au voisinage de $E$ est noeth\'erien.
\end{thm}
\begin{proof}
Soit $(I_n)_{n\in\N}$ une suite croissante d'id\'eaux de type fini de~$\Os(E)$. Pour $n\in\N$, notons $\Is_n$ le faisceau d'id\'eaux coh\'erents de $\Os_X$ engendr\'e par $I_n$. D'apr\`es la proposition \ref{stationne} et la compacit\'e de $E$, il existe un rang~\mbox{$n_{0}\in\N$} \`a partir duquel la suite $(\Is_n)_{n\in\N}$ stationne.

Puisque l'id\'eal $I_{n_0}$ est de type fini, il poss\`ede un syst\`eme g\'en\'erateur fini \mbox{$(f_1,\ldots,f_p)$}, avec $p\in\N$ et, quel que soit $i\in\cn{1}{p}$, $f_i\in \Os(E)$. Le morphisme de faisceaux 
$$\varphi : \begin{array}{ccc}
\Os_X^p & \to & \Is_{n_0}\\
(a_1,\ldots,a_p) & \mapsto & a_1 f_1 + \ldots + a_p f_p
\end{array}$$
est alors surjectif.

Soit $n\ge n_0$. Notons $\Gs$ le noyau du morphisme de faisceaux $\varphi$. C'est encore un faisceau coh\'erent sur $E$. Nous disposons de la suite exacte 
$$0\to \Gs \to \Os^p \to \Is_n \to 0.$$ 
Puisque $H^1(E,\Gs)=0$, le morphisme 
$$\begin{array}{ccc}
\Os(E)^p & \to & \Is_n(E)\\
(a_1,\ldots,a_p) & \mapsto & a_1 f_1 + \ldots + a_p f_p
\end{array},$$ 
est surjectif. Par cons\'equent, nous avons
$$I_{n} \subset \Is_{n}(E) = (f_1,\ldots,f_p)\, \Os(E) \subset I_{n_0}.$$
On en d\'eduit que $I_n=I_{n_0}$.
\end{proof}

\subsection{S\'eries arithm\'etiques}

Dans ce paragraphe, nous appliquons le th\'eor\`eme obtenu afin de d\'emontrer la noeth\'erianit\'e de certains anneaux de s\'eries arithm\'etiques. Il est vraisemblable que l'analogue du th\'eor\`eme de Frisch vaut pour toute partie semi-analytique de la droite $X=\E{1}{A}$. Cependant, pour le d\'emontrer par la m\'ethode pr\'esent\'ee ci-dessus, il nous faudrait savoir que les parties semi-analytiques de~$X$ sont localement connexes. Nous ne nous lancerons pas dans la d\'emonstration de ce r\'esultat et nous contenterons d'adapter le th\'eor\`eme de Frisch au cas des couronnes ferm\'ees au-dessus de certaines parties compactes de l'espace~$B$.

Soit~$V$ une partie compacte et connexe de l'espace~$B$. Soit~$s$ et~$t$ deux nombres r\'eels v\'erifiant $0\le s\le t$. Posons
$$C = \overline{C}_{V}(s,t) =  \{x\in X_{V}\, |\, s\le |T(x)|\le t\}..$$

\begin{prop}\index{Couronne!connexite par arcs locale@connexit\'e par arcs locale}\index{Disque!connexite par arcs locale@connexit\'e par arcs locale}
La couronnne~$C$ de $X$ est localement connexe par arcs.
\end{prop}
\begin{proof}
Si $x$ est un point int\'erieur \`a $C$, le r\'esultat est vrai car il l'est pour l'espace $X$ lui-m\^eme, d'apr\`es le th\'eor\`eme~\ref{resume}. Nous supposerons donc, d\'esormais, que le point $x$ est situ\'e sur le bord de la couronne $C$. En particulier, nous avons n\'ecessairement $|T(x)|=s$ ou $|T(x)|=t$. Nous supposerons que~\mbox{$|T(x)|=t$}. L'autre cas se traite de m\^eme. Nous allons distinguer selon le type du point $x$ et de son projet\'e $\pi(x)$ sur la base.

\bigskip

Supposons, tout d'abord, que le point $\pi(x)$ soit un point extr\^eme : il existe $\m\in\Sigma_{f}$ tel que $\pi(x)=\tilde{a}_{\m}$. Si le point $x$ est le point $\eta_{s}$, alors le r\'esultat provient du corollaire \ref{voissectiondep}, si $t\ne 1$, et de la proposition \ref{sectionGaussext}, si $t=1$. Il faut plus pr\'ecis\'ement revenir \`a la description explicite des sections donn\'ee dans la preuve de ces propositions. Il nous reste \`a traiter le cas o\`u $x$ v\'erifie $|T(x)|=1$, mais n'est pas le point $\eta_{1}$. Un tel point appartient n\'ecessairement \`a l'int\'erieur de la couronne $C$. En effet, il existe $\tilde{\alpha} \in \tilde{k}_{\m}^*$ et $u\in\of{[}{0,1}{[}$ tels que $x=\eta_{\tilde{\alpha},u}$. Choisissons un relev\'e $\alpha$ de $\tilde{\alpha}$ dans $\hat{A}_{\m}$. Soit $v\in\of{]}{u,1}{[}$. Alors le voisinage de $x$ dans $X$ d\'efini par 
$$U = \left\{y\in \pi^{-1}(\of{]}{a_{0},\tilde{a}_{\m}}{]})\, \big|\, |(T-\alpha)(y)| < v \right\}$$
est contenu dans $\overline{C}(s,1)$. En effet, soient $\eps\in\of{]}{0,+\infty}{]}$ et $y\in U\cap X_{a_{\m}^\eps}$. Nous avons $ |(T-\alpha)(y)| < v < 1$. Puisque $|\alpha(y)| = |\alpha|_{\m}^\eps=1$, cela impose que $|T(y)|=1$. 

Lorsque le point $\pi(x)$ est le point central $a_{0}$ de $B$, le r\'esultat se d\'emontre de fa\c{c}on identique.

\bigskip

Venons-en, \`a pr\'esent, au cas de la partie archim\'edienne de $X$. Soit $\sigma\in\Sigma_{\infty}$. Rappelons que, d'apr\`es la proposition \ref{produitarc}, l'application 
$$\varphi : {\renewcommand{\arraystretch}{1.2}\begin{array}{ccc}
X_{a_{\sigma}} \times \of{]}{0,1}{]} & \to & X'_{\sigma}\\
(x,\eps) & \mapsto & x^\eps 
\end{array}}$$
est un hom\'eomorphisme. Nous supposerons que $K_{\sigma}=\C$. Le cas $K_{\sigma}=\R$ se traite de m\^eme. Nous avons 
$$\varphi^{-1}(X'_{\sigma}\cap \overline{C}(s,t)) = \left\{(u,z)\in \of{]}{0,1}{]} \times \C\, \big|\, s^{1/u} \le |z| \le t^{1/u}\right\}.$$
Cette partie est localement connexe par arcs et il en est de m\^eme de son intersection avec la couronne~$C$.

\bigskip

Il nous reste \`a traiter le cas o\`u le point $\pi(x)$ est de la forme $a_{\sigma}^\lambda$, avec $\sigma\in\Sigma_{f}$ et $\lambda\in\of{]}{0,+\infty}{[}$. Nous pouvons supposer que $\lambda=1$. Comme dans le cas des fibres au-dessus d'un corps trivialement valu\'e, il nous suffit de traiter le cas o\`u $x$ est le point $\eta_{t}$ de sa fibre. Nous supposerons que $t\in\of{]}{0,1}{[}$. Les autres cas se traitent de m\^eme. Soit $U$ un voisinage de $x$ dans $X$. Il existe un voisinage connexe par arcs $V$ de $x$ dans $X_{a_{\sigma}} \cap U$. Il existe $\beta\in\of{]}{0,1}{[}$ tel que, quel que soit $u \in \of{]}{t^{1/\beta},t^\beta}{[}$, on ait $\eta_{u} \in X_{a_{\sigma}}\cap V$. D'apr\`es la proposition \ref{produitum}, quitte \`a augmenter $\beta$, nous pouvons supposer que la partie
$$W = \left\{ x^\eps, \, x\in X_{a_{\sigma}} \cap V, \eps \in \of{]}{\beta,1/\beta}{[}\right\}$$ 
est un voisinage de $x$ dans $U$. La trace de~$W$ sur chaque fibre est connexe par arcs en tant qu'intersection sur un arbre de deux parties connexes par arcs (l'une \'etant hom\'eomorphe \`a $V$, l'autre  \'etant une couronne). En outre, ces fibres sont jointes par une section depuis la base : l'application qui au point $a_{\sigma,\eps}$, avec $\eps\in \of{]}{\beta,1/\beta}{[}$, associe le point $\eta_{t}$ de sa fibre. On en d\'eduit que la trace de la partie~$W$ sur la couronne~$C$ est connexe par arcs.
\end{proof}

\begin{cor}\index{Couronne!morcelable}\index{Disque!morcelable}\index{Morcelable!couronne}
La couronne~$C$ de $X$ est morcelable.
\end{cor}
\begin{proof}
Soit $x$ un point de $C$. Soient $U$ un voisinage ouvert de $x$ dans~$X$ et $Z$ une sous-vari\'et\'e analytique de $U$ int\`egre en $x$. Nous devons montrer qu'il existe un voisinage $V$ de $x$ dans $E\cap U$ qui poss\`ede un syst\`eme fondamental de voisinages ouverts dans $U$ dont les traces sur $Z$ sont connexes. 

Supposons, tout d'abord, que $Z=U$ au voisinage de $x$. Dans ce cas, la proposition pr\'ec\'edente nous permet de conclure. Si, maintenant, $Z$ est une sous-vari\'et\'e analytique stricte de $U$, nous en connaissons pr\'ecis\'ement la forme gr\^ace aux descriptions donn\'ees dans la partie \ref{cint}. En particulier, au voisinage du point $x$, la sous-vari\'et\'e $Z$ est soit un point, soit hom\'eomorphe \`a un intervalle, soit une fibre extr\^eme. Le r\'esultat est imm\'ediat dans chacun de ces cas. 
\end{proof}

\begin{thm}\index{Couronne!noetherianite de l'anneau des sections@noeth\'erianit\'e de l'anneau des sections}\index{Disque!noetherianite de l'anneau des sections@noeth\'erianit\'e de l'anneau des sections}\index{Noetherianite@Noeth\'erianit\'e!des anneaux de sections sur une couronne de A1@des anneaux de sections sur une couronne de $\AA$}\index{Espace de Stein!couronne compacte de A1@couronne compacte de $\AA$}
L'anneau $\Os(C)$ des germes de fonctions analytiques au voisinage de la couronne $C$ de $X$ est noeth\'erien.
\end{thm}
\begin{proof}
Une telle partie est morcelable en vertu du corollaire pr\'e\-c\'e\-dent. Nous savons \'egalement qu'elle est de Stein, d'apr\`es le th\'eor\`eme \ref{ABcompact}. Le th\'eor\`eme \ref{Frisch} nous permet donc de conclure.
\end{proof}

\begin{cor}\label{Harbater}\index{Noetherianite@Noeth\'erianit\'e!d'anneaux de series arithmetiques@d'anneaux de s\'eries arithm\'etiques}
Soient $\Sigma'$ un sous-ensemble fini de $\Sigma$ contenant $\Sigma_{\infty}$ et~$(r_{\sigma})_{\sigma\in\Sigma'}$ une famille d'\'el\'ements de $\of{]}{0,1}{[}$. Il existe un \'el\'ement $N\in A^*$ tel que
$$\bigcap_{\sigma\in\Sigma'} A_{\sigma} = A\left[\frac{1}{N}\right].$$ 
Le sous-anneau de $K(\!(T)\!)$ constitu\'e des s\'eries de la forme $\sum_{k\ge k_{0}} a_{k}\, T^k$ v\'erifiant les conditions
\begin{enumerate}[\it i)]
\item $k_{0}\in\Z$,
\item $\forall k\ge k_{0}$, $a_{k}\in A[1/N]$,
\item $\forall \sigma\in\Sigma'$, $\exists r > r_{\sigma}$, $\disp \lim_{k\to +\infty} |a_{k}|_{\sigma}\, r^k =0$ 
\end{enumerate}
est noeth\'erien.

Le sous-anneau de $K[\![T]\!]$ constitu\'e des s\'eries de la forme $\sum_{k\ge 0} a_{k}\, T^k$ v\'erifiant les conditions
\begin{enumerate}[\it i)]
\item $\forall k\ge 0$, $a_{k}\in A[1/N]$,
\item $\forall \sigma\in\Sigma'$, $\exists r > r_{\sigma}$, $\disp \lim_{k\to +\infty} |a_{k}|_{\sigma}\, r^k =0$ 
\end{enumerate}
est noeth\'erien.
\end{cor}
\begin{proof}
Il suffit d'appliquer le th\'eor\`eme pr\'ec\'edent \`a une couronne bien choisie. Posons
$$t = \max_{\sigma\in\Sigma'}(r_{\sigma}) \in\of{]}{0,1}{[}.$$
Quel que soit $\sigma\in\Sigma'$, il existe $\eps_{\sigma} \in\of{]}{0,1}{]}$ tel que 
$$t^{1/\eps_{\sigma}} = r_{\sigma}.$$
D\'efinissons une partie compacte $V$ de $B$ par 
$$V = \left( \bigcup_{\sigma\in\Sigma'} \of{[}{a_{0},a_{\sigma}^{\eps_{\sigma}}}{]} \right) \cup \left( \bigcup_{\sigma\notin \Sigma'} B_{\sigma} \right).$$
Soit~$s\in\of{]}{0,t}{]}$. D'apr\`es la proposition \ref{descriptionsectionscouronnes}, le premier anneau consid\'er\'e n'est autre que l'anneau~$\Os(\overline{C}_{V}(s,t))$. Il est noeth\'erien, en vertu du th\'eor\`eme pr\'ec\'edent.

Le second \'enonc\'e s'obtient de m\^eme en consid\'erant le disque~$\overline{D}_{V}(t)$, au lieu de la couronne~$\overline{C}_{V}(s,t)$. 

\end{proof}

Comme cas particulier du th\'eor\`eme, nous retrouvons un r\'esultat de D. Harbater (\emph{cf.} \cite{Harbater}, th\'eor\`eme 1.8). Signalons que notre d\'emonstration se distingue tr\`es nettement de la sienne, qui passe par une description explicite de tous les id\'eaux premiers de l'anneau \'etudi\'e.

\begin{cor}\index{Theoreme@Th\'eor\`eme!d'Harbater}\index{Harbater|see{Th\'eor\`eme d'Harbater}}\newcounter{Zr}\setcounter{Zr}{\thepage}
Soit~$r_{\infty}\in\of{]}{0,1}{[}$. Consid\'erons le sous-anneau $\Z_{r^+}[\![T]\!]$ de~$\Z[\![T]\!]$ constitu\'e des s\'eries de la forme $\sum_{k\ge 0} a_{k}\, T^k$ v\'erifiant la condition
$$\exists r > r_{\infty},\ \disp \lim_{k\to +\infty} |a_{k}|_{\infty}\, r^k =0.$$
C'est l'anneau des fonctions holomorphes au voisinage du disque de centre $0$ et de rayon $r_{\infty}$ de $\C$ dont le d\'eveloppement en s\'erie enti\`ere en $0$ est \`a coefficients entiers. L'anneau~$\Z_{r^+}[\![T]\!]$
est noeth\'erien.
\end{cor}
\begin{proof}
Il suffit d'appliquer le second r\'esultat du th\'eor\`eme pr\'ec\'edent avec~$K=\Q$ et~$\Sigma'=\Sigma_{\infty}$.
\end{proof}


\section{Probl\`eme de Galois inverse}\label{pdgi}

\index{Probleme de Galois inverse@Probl\`eme de Galois inverse|(}

Nous exposons ici une application de notre th\'eorie au probl\`eme de Galois inverse. Pr\'ecis\'ement, nous nous proposons de d\'emontrer que tout groupe groupe fini est le groupe de Galois d'une extension finie et galoisienne du corps~$\Ms(\D)$, o\`u~$\D$ d\'esigne le disque relatif ouvert de rayon~$1$ centr\'e en la section nulle.. Signalons que, dans le cas o\`u le corps de nombres~$K$ consid\'er\'e n'est autre que le corps~$\Q$, nous red\'emontrons un r\'esultat de D.~Harbater (\emph{cf.} \cite{galoiscovers}, corollaire~3.8). Nous souhaitons insister sur le fait que la d\'emonstration que nous proposons est purement g\'eom\'etrique, ce qui la distingue de celle de D.~Harbater, tr\`es alg\'ebrique. 

Nous utiliserons un proc\'ed\'e classique : construction de rev\^etements galoisiens cycliques, puis recollement de ces rev\^etements afin d'en obtenir un nouveau ayant pour groupe de Galois un groupe fini prescrit. L'on trouvera une introduction tr\`es agr\'eable \`a ces techniques dans l'article \cite{LiuGalois}, o\`u Q.~Liu d\'emontre -- d'apr\`es D.~Harbater et en suivant une id\'ee de J.-P.~Serre -- que, pour tout nombre premier~$p$, tout groupe fini est le groupe de Galois d'une extension finie et galoisienne du corps~$\Q_{p}(T)$. La mise en {\oe}uvre de ces deux \'etapes que nous proposons nous semble particuli\`erement simple, une fois connues les propri\'et\'es de la droite analytique sur un anneau d'entiers de corps de nombres. Les m\'ethodes utilis\'ees par D.~Harbater nous paraissent d'une difficult\'e technique bien sup\'erieure (\emph{cf.}~\cite{mockcovers}, proposition~2.2 pour la construction des rev\^etements cycliques et~\cite{galoiscovers}, th\'eor\`eme~3.6, dont la preuve fait appel aux r\'esultats de l'article~\cite{Harbater},  pour le recollement).

Mentionnons pour finir que nous allons en fait construire des faisceaux d'al\-g\`e\-bres coh\'erents ayant pour groupe d'automorphismes un groupe fini prescrit. Bien entendu, ces faisceaux sont les images directes de faisceaux structuraux de rev\^etements ramifi\'es de la droite analytique sur un anneau d'entiers de corps de nombres et il y aurait tout int\'er\^et \`a mener plut\^ot nos constructions dans ce langage. Nous nous en abstenons uniquement parce qu'aucune r\'ef\'erence concernant ces espaces n'est disponible. Nous indiquerons cependant en remarque les traductions dans ce cadre ; elles sont imm\'ediates pour qui dispose d'une bonne th\'eorie des courbes analytiques sur un anneau d'entiers de corps de nombres.

\bigskip

Introduisons quelques notations. Rappelons que nous notons
$$\D = \mathring{D}(0,1) = \left\{ x\in X\, \big|\, |T(x)| <1\right\}.$$
Pour tout \'el\'ement~$\m$ de~$\Sigma_{f}$, nous posons
$$\begin{array}{cccl}
&\D_{\m} &=& \D \cap X_{\m},\\ 
&\D'_{\m} &=& \D_{\m} \setminus X_{0}\\ 
\textrm{et} & \D''_{\m} &=& \D'_{m} \setminus \tilde{X}_{\m}.
\end{array}$$
\newcounter{Dm}\setcounter{Dm}{\thepage}

\subsection{Construction locale de rev\^etements cycliques}

\index{Probleme de Galois inverse@Probl\`eme de Galois inverse!construction locale|(}

Soient~$V$ une partie de~$X$ et~$P$ un polyn\^ome unitaire \`a coefficients dans~$\Os(V)$. Notons~$n$ son degr\'e. On d\'efinit un pr\'e\-fais\-ceau~$\Fs_{P}$ sur~$V$ en posant, pour toute partie ouverte~$U$ de~$V$,
$$\Fs_{P}(U) = \Os(U)[S]/(P(S))$$
et en utilisant les morphismes de restriction induits par ceux du faisceau~$\Os$.

\begin{lem}\label{FPcoherent}\index{Faisceau!coherent@coh\'erent}
Le pr\'efaisceau~$\Fs_{P}$ est un faisceau de $\Os_{V}$-alg\`ebres coh\'erent. 
\end{lem}
\begin{proof}
On constate imm\'ediatement que le pr\'efaisceau~$\Fs_{P}$ est un pr\'efaisceau de $\Os_{V}$-alg\`ebres. Il nous suffit donc de montrer que c'est un faisceau et un faisceau de $\Os_{V}$-modules coh\'erent. Puisque le polyn\^ome~$P$ est unitaire, le morphisme de $\Os_{V}$-modules
$$\begin{array}{ccc}
\Os_{V}^n & \to & \Fs\\
(a_{0},\ldots,a_{n-1}) & \mapsto & \disp \sum_{i=0}^{n-1} a_{i}\, S^i
\end{array}$$
est un isomorphisme. On en d\'eduit que le pr\'efaisceau~$\Fs$ est un faisceau, puis qu'il est coh\'erent, car le faisceau structural~$\Os$ l'est, en vertu du th\'eor\`eme \ref{coherence}.
\end{proof}

\begin{rem}
Le faisceau~$\Fs_{P}$ est l'image directe du faisceau structural d'une courbe analytique sur~$A$. Celle-ci nous est donn\'ee comme un rev\^etement ramifi\'e, de degr\'e inf\'erieur \`a~$n$, de la partie~$V$ de la droite analytique~$\E{1}{A}$ .
\end{rem}

Nous allons, \`a pr\'esent, restreindre notre \'etude aux faisceaux~$\Fs_{P}$ pour une classe de polyn\^omes~$P$ particuliers. Soient~$n$ un entier sup\'erieur \`a~$1$, $p$ un nombre premier congru \`a~$1$ modulo~$n$ et~$\m$ un id\'eal maximal de l'anneau~$A$ contenant~$p$. Posons
$$Q(S) = S^n-\pi_{\m}^n-T \in \Os(\D''_{\m})[S].$$
Le r\'esultat du lemme suivant donne la raison du choix des entiers~$n$ et~$p$.

\begin{lem}\index{Racines de l'unite@Racines de l'unit\'e}
L'anneau~$A_{\m}$ contient~$n$ racines $n^{\textrm{\`emes}}$ de l'unit\'e.
\end{lem}
\begin{proof}
Puisque l'anneau~$A_{\m}$ contient l'anneau~$\Z_{p}$, il suffit de montrer que le polyn\^ome~$U^n-1$ poss\`ede~$n$ racines dans~$\Z_{p}$. Le groupe multiplicatif~$\F_{p}^*$ du corps r\'esiduel~$\F_{p}$ de~$\Z_{p}$ est cyclique et d'ordre~$p-1$. Puisque~$n$ divise~$p-1$, le groupe~$\F_{p}^*$ contient un \'el\'ement d'ordre exactement~$n$ et le polyn\^ome~$U^n-1$ est scind\'e \`a racines simples sur~$\F_{p}$. Le lemme de Hensel assure qu'il l'est encore sur~$\Z_{p}$.
\end{proof}

Pour tout entier positif~$i$ et tout nombre rationnel~$k$, posons
$$C_{k}^i = \frac{k\, (k-1)\, \cdots (k-i+1)}{i!} \in \Q.$$
Rappelons que nous avons l'\'egalit\'e
$$\left(\sum_{i\ge 0} C_{\frac{1}{n}}^i\, Z^i\right)^n = 1+Z \textrm{ dans } \Q[\![Z]\!].$$

\begin{lem}
Pour tout \'el\'ement~$x$ de~$\D'_{\m}$ et tout entier positif~$i$, nous avons
$$\left|C_{\frac{1}{n}}^i(x)\right| \le 1.$$
\end{lem}
\begin{proof}
Soient~$x$ un \'el\'ement de~$\D'_{\m}$ et~$i$ un entier positif. Notons~$|.|_{x}$ la valeur absolue sur le corps~$\Hs(x)$. Remarquons que l'application
$$C_{k}^{i} : \Q \to \Q$$
est polynomiale, et donc continue lorsque l'on munit le corps~$\Q$ de la valeur absolue~$|.|_{x}$. Nous savons que, pour tout entier~$l$, l'\'el\'ement~$C_{l}^i$ est entier. Il v\'erifie donc l'in\'egalit\'e
$$\left|C_{l}^i\right|_{x} \le 1,$$ 
puisque la valeur absolue~$|.|_{x}$ est ultram\'etrique. En outre, le nombre premier~$p$ ne divise pas l'entier~$n$. Par cons\'equent, le nombre rationnel~$\frac{1}{n}$ appartient \`a~$\Z_{p}$ et il est donc limite d'\'el\'ements de~$\Z$ pour la valeur absolue~$|.|_{x}$. On en d\'eduit le r\'esultat voulu. 
\end{proof}

Fixons~$\zeta$ une racine primitive $n^\textrm{\`eme}$ de l'unit\'e. Notons~$\tau$ la permutation cyclique $(1\ 2\ \cdots\ n)$ de l'ensemble~$\cn{1}{n}$. Posons
$$U =  \left\{x\in D'_{\m}\, \big|\, |T(x)|< |\pi_{m}(x)|^n  \right\}.$$

\begin{prop}\label{nonramifie}
Il existe un isomorphisme de $\Os_{U}$-alg\`ebres
$$\varphi : \Fs_{Q} \to \Os^n$$
tel que, pour tout ouvert~$V$ de~$U$ et tout \'el\'ement~$s$ de~$\Fs(V)$, nous ayons
$$\varphi(\zeta s) = \tau(\varphi(s)).$$
\end{prop}
\begin{proof}
Consid\'erons la fonction
$$f = \pi_{\m}^{-n}\, T$$
d\'efinie sur~$X'_{\m}$. Pour tout \'el\'ement~$r$ de~$\of{]}{0,1}{[}$, consid\'erons la partie~$V_{r}$ de~$\D'_{\m}$ d\'efinie par
$$V_{r} = \left\{x\in \D'_{\m}\, \big|\, |T(x)| \le r |\pi_{\m}(x)|^n \right\}.$$
Pour tout \'el\'ement~$x$ de~$V_{r}$ et tous entiers~$a\ge 0$ et~$b\ge a$, nous avons
$$\left| \sum_{i=a}^b C_{\frac{1}{n}}^i(x)\, f(x)^i \right| \le r^a.$$
On en d\'eduit que la s\'erie $\sum_{i\ge 0} C_{\frac{1}{n}}^i\, f^i$ converge uniform\'ement sur~$V_{r}$. Puisque tout point de~$U$ poss\`ede un voisinage de la forme~$V_{r}$, pour un certain \'el\'ement~$r$ de~$\of{]}{0,1}{[}$, la s\'erie $\sum_{i\ge 0} C_{\frac{1}{n}}^i\, f^i$ d\'efinit une fonction~$g$ sur~$U$. Cette fonction v\'erifie l'\'egalit\'e
$$g^n = 1+f = 1+\pi_{\m}^{-n}\, T \textrm{ dans } \Os(U).$$
On en d\'eduit que nous avons l'\'egalit\'e
$$Q(S) = S^n-\pi_{\m}^n-T = \prod_{j=0}^n (S - \pi_{\m}\, \zeta^{j}\, g) \textrm{ dans } \Os(U)[S].$$
Par cons\'equent, le morphisme
$$\begin{array}{ccc}
\Fs_{Q} & \to & \Os^n\\
F(S) & \mapsto & \left(F(\pi_{\m}\,g), F(\pi_{\m}\, \zeta^{-1}\, g),\ldots, F(\pi_{\m}\, \zeta^{-(n-1)}\, g)\right) 
\end{array}$$
est un isomorphisme. On v\'erifie imm\'ediatement qu'il satisfait la condition requise.
\end{proof}

\begin{rem}
La premi\`ere partie du r\'esultat signifie que le rev\^etement associ\'e au faisceau~$\Fs_{Q}$ est trivial au-dessus de l'ouvert~$U$. La seconde assure que le groupe $\of{\la}{\zeta}{\ra} \simeq \Z/n\Z$ agit sur le rev\^etement par une permutation cyclique des feuillets du lieu trivial.  
\end{rem}

\begin{lem}
Le polyn\^ome~$Q(S)=S^n-\pi_{\m}^n-T$ est irr\'eductible sur le corps Frac$(\Os(\D''_{\m}))$. En particulier, l'anneau $\Fs_{Q}(\D''_{\m})$ est int\`egre.
\end{lem}
\begin{proof}
Notons~$x$ le point~$0$ de la fibre extr\^eme~$\tilde{X}_{\m}$. D'apr\`es le corollaire \ref{descratext}, l'anneau local en ce point est isomorphe \`a l'anneau~$\hat{A}_{\m}[\![T]\!]$. Remarquons que le polyn\^ome~$Q(S)$ est irr\'eductible sur le corps~$\hat{K}_{\m}(\!(T)\!)$. En effet, il n'y a aucune racine, pour des raisons de valuation $T$-adique, est s\'eparable et le groupe de Galois de son extension de d\'ecomposition agit transitivement sur ses racines.

D'apr\`es le principe du prolongement analytique (\emph{cf.} th\'eor\`eme~\ref{prolan}), le morphisme naturel $\Os(\D''_{\m})\to \Os_{X,x}$ est injectif. Par cons\'equent, le corps $\Frac(\Os(\D''_{m}))$ est un sous-corps de~$\Frac(\Os_{X,x})$, et donc de~$\hat{K}_{\m}(\!(T)\!)$. On en d\'eduit que le polyn\^ome~$Q(S)$ est irr\'eductible sur le corps Frac($\Os(\D''_{\m})$).

Puisque le polyn\^ome~$Q(S)$ est unitaire, l'unicit\'e de la division euclidienne assure que le morphisme 
$$\Os(\D''_{\m})[S]/(Q(S)) \to \Frac(\Os(\D''_{\m}))[S]/(Q(S))$$
est injectif. Puisque l'anneau au but est int\`egre, celui \`a la source, qui n'est autre que l'anneau~$\Fs_{Q}(\D''_{\m})$, l'est \'egalement.


\end{proof}

\begin{rem}
Ce r\'esultat signifie que la courbe associ\'ee au faisceau~$\Fs_{Q}$ est int\`egre, c'est-\`a-dire r\'eduite et irr\'eductible.
\end{rem}

Nous pouvons \^etre encore plus pr\'ecis. 

\begin{lem}\label{integreunebranche}
Soient~$x$ un point de~$U$ et~$i$ un \'el\'ement de~$\cn{1}{n}$. Le morphisme
$$\rho : \Fs_{Q}(\D''_{\m}) \to \Fs_{x} \xrightarrow[\sim]{\varphi_{x}} \Os_{X,x}^n \xrightarrow[]{p_{i}} \Os_{X,x},$$
o\`u~$p_{i}$ est la projection sur le $i^\textrm{\`eme}$ facteur, est injectif.
\end{lem}
\begin{proof}
Soit~$s$ un \'el\'ement de l'anneau $\Fs_{Q}(\D''_{\m})=\Os(\D''_{\m})/(Q(S))$ dont l'image par le morphisme~$\rho$ est nulle. Choisissons un \'el\'ement~$F(S)$ de $\Os(\D''_{\m})[S]$ qui repr\'esente la section~$s$. Reprenons les notations de la preuve de la proposition \ref{nonramifie}. Par hypoth\`ese, nous avons
$$R(\pi_{\m}\, \zeta^{-i}\, g)=0 \textrm{ dans } \Os_{X,x}.$$
Pour montrer que l'\'el\'ement~$s$ est nul, il suffit de montrer que le polyn\^ome~$Q(S)$ est le polyn\^ome minimal de l'\'el\'ement $\pi_{\m}\, \zeta^{-i}\, g$ sur le corps~$\Frac(\Os(\D''_{\m}))$. C'est bien le cas, puisque le lemme pr\'ec\'edent assure que le polyn\^ome~$Q$ est irr\'eductible sur le corps~$\Frac(\Os(\D''_{\m}))$.
\end{proof}

\begin{rem}
Ce r\'esultat est une sorte de principe du prolongement analytique sur la courbe associ\'ee au faisceau~$\Fs_{Q}$ : si une fonction holomorphe sur la courbe est nulle au voisinage d'un point de l'un des feuillets du rev\^etement, alors elle est nulle partout. On attend que ce principe vaille pour toute courbe irr\'eductible.
\end{rem}

Terminons par un r\'esultat topologique. 

\begin{lem}\label{Festferme}
La partie 
$$F = \D''_{\m} \setminus U = \left\{x\in \D''_{\m}\, \big|\, |T(x)|\ge |\pi_{m}(x)|^n \right\}$$
est ferm\'ee dans le disque~$\D$.
\end{lem}
\begin{proof}
Il suffit de montrer que~$F$ est ferm\'ee dans~$\D_{\m}$ puisque cette derni\`ere partie est elle-m\^eme ferm\'ee dans~$\D$. En d'autres termes, nous souhaitons montrer que la partie 
$$V = U \cup (\D\cap X_{0})$$
est ouverte dans~$\D_{\m}$. Puisque~$U$ est une partie ouverte de~$\D_{\m}$, il suffit de montrer que~$V$ est voisinage dans~$\D_{\m}$ de chacun des points de~$\D\cap X_{0}$.

Soit~$x$ un point de~$\D\cap X_{0}$. 
Posons $r=|T(x)|$. C'est un \'el\'ement de l'intervalle~$\of{]}{0,1}{[}$. Soient~$s$ un \'el\'ement de~$\of{]}{r,1}{[}$ et~$\eps$ un \'el\'ement de~$\of{]}{0,1}{[}$ tels que l'on ait $|\pi_{\m}|_{\m}^{n\eps} > s$. La partie 
$$\left\{y\in \pi^{-1}({[}{a_{0},a_{\m}^\eps}{[}) \, \big|\, |T(y)|<s \right\}$$
est un voisinage ouvert du point~$x$ dans~$\D_{\m}$ qui est contenu dans~$V$. 
\end{proof}



\index{Probleme de Galois inverse@Probl\`eme de Galois inverse!construction locale|)}

\subsection{Recollement}

\index{Probleme de Galois inverse@Probl\`eme de Galois inverse!recollement|(}

Soit~$G$ un groupe fini. Notons~$n$ son ordre et $g_{1},\ldots,g_{n}$ ses \'el\'ements. Chacun de ces \'el\'ements engendre un sous-groupe cyclique de~$G$. Nous allons construire, par la m\'ethode mise en place au num\'ero pr\'ec\'edent, un rev\^etement galoisien cyclique associ\'e \`a chacun des \'el\'ements du groupe~$G$. Il ne nous restera plus ensuite qu'\`a les recoller convenablement. 

En termes g\'eom\'etriques, nous allons recoller les rev\^etements au-dessus de leur lieu de trivialit\'e en tenant compte des relations entre les \'el\'ements du groupe~$G$. Ce proc\'ed\'e est simple et naturel et l'on ne doit pas se laisser rebuter par la technicit\'e apparente de la construction qui suit.

\bigskip

Soit~$i$ un \'el\'ement de~$\cn{1}{n}$. Notons~$n_{i}$ l'ordre de l'\'el\'ement~$g_{i}$ dans le groupe~$G$. C'est un diviseur de~$n$ et nous noterons~$d_{i}$ le quotient. Soient~$p_{i}$ un nombre premier congru \`a~$1$ modulo~$n_{i}$ et~$\m_{i}$ un id\'eal maximal de l'anneau~$A$ qui contient~$p_{i}$. Soit~$\zeta_{i}$ une racine primitive $n_{i}^\textrm{\`eme}$ de l'unit\'e dans~$\hat{A}_{\m_{i}}$. Notons~$\Fs_{i}$ le faisceau $\Fs_{S^{n_{i}}-\pi_{\m_{i}}^{n_{i}}-T}$ sur~$\D''_{\m_{i}}$ et~$\Gs_{i}$ le faisceau~$\Fs_{i}^{d_{i}}$. Posons
$$U_{i} = \left\{x\in \D''_{\m_{i}}\, \big|\, |T(x)| < |\pi_{m_{i}}(x)|^{n_{i}} \right\} \textrm{ et } F_{i} = \D''_{\m_{i}}\setminus U_{i}.$$
Notons~$\tau_{i}$ la permutation cyclique $(1\ 2\ \cdots\ n_{i})$ de l'ensemble~$\cn{1}{n_{i}}$. D'apr\`es la proposition \ref{nonramifie}, il existe un isomorphisme de $\Os_{U_{i}}$-alg\`ebres
$$\varphi_{i} : \Fs_{i} \xrightarrow[]{\sim} \Os^{n_{i}}$$
qui v\'erifie la condition suivante : pour tout ouvert~$V$ de~$U_{i}$ et toute section~$s$ de~$\Fs_{i}$ sur~$V$, nous ayons
$$\varphi_{i}(\zeta_{i}\, s) = \tau_{i}(\varphi_{i}(s)).$$

Choisissons des \'el\'ements $a_{i,0},\ldots,a_{i,d_{i}-1}$ de~$\cn{1}{n}$ de sorte que tout \'el\'ement du quotient $G/\la g_{i} \ra$ poss\`ede un repr\'esentant et un seul parmi les \'el\'ements $g_{a_{i,0}},\ldots,g_{a_{i,d_{i}-1}}$. Notons~$\sigma_{i}$ la permutation de l'ensemble~$\cn{1}{n}$ telle que
$$\forall u\in\cn{0}{d_{i}-1},\, \forall v\in\cn{1}{n_{i}},\, g_{a_{i,u}}\, g_{i}^{v-1} = g_{\sigma_{i}(un_{i}+v)}.$$

\bigskip

D'apr\`es le lemme \ref{Festferme}, pour tout \'el\'ement~$i$ de~$\cn{1}{n}$, la partie~$F_{i}$ est ferm\'ee dans~$\D$. D\'efinissons une partie ouverte de~$\D$ par 
$$U_{0} = \D\setminus \bigcup_{1\le i\le n} F_{i}.$$
Notons~$\Gs_{0}$ le faisceau~$\Os^{n}$ sur~$U_{0}$.

\begin{lem}\label{U0connexe}
L'ouvert~$U_{0}$ de la droite~$X$ est connexe.
\end{lem}
\begin{proof}
Notons $\Sigma_{0}=\{\m_{1},\ldots,\m_{n}\}$. Par d\'efinition, nous avons
$$U_{0} = \bigcup_{1\le i\le n} \left\{x\in \D''_{\m_{i}}\, \big|\, |T(x)| < |\pi_{m_{i}}(x)| \right\} \cup \bigcup_{\sigma\in\Sigma\setminus\Sigma_{0}} B_{\sigma}.$$
La projection de cette partie est 
$$B_{0} = B\setminus \bigcup_{1\le i\le n} \{\tilde{a}_{\m_{i}}\},$$
qui est connexe. En outre, la section nulle d\'efinie une section continue de l'application $\pi : U_{0} \to B_{0}$ et, pour tout \'el\'ement~$b$ de~$B_{0}$, la partie~$X_{b}\cap U_{0}$ est connexe. On en d\'eduit que la partie~$U_{0}$ est connexe. 
\end{proof}

La famille $(U_{0},\D''_{\m_{1}},\ldots,\D''_{\m_{n}})$ d\'efinit un recouvrement ouvert du disque~$\D$. Les seules intersections de deux \'el\'ements de cette famille \`a n'\^etre pas vides sont celles de la forme $\D''_{\m_{i}}\cap \D''_{\m_{i}}$, pour $i\in\cn{0}{n}$, et $\D''_{\m_{i}}\cap U_{0}=U_{i}$, pour $i\in\cn{1}{n}$. Pour d\'efinir un faisceau d'alg\`ebres~$\Gs$ en recollant les faisceaux $\Gs_{0},\ldots,\Gs_{n}$, il nous suffit de choisir un isomorphisme de $\Os$-alg\`ebres entre~$\Gs_{i}$ et~$\Gs_{0}$ au-dessus de l'ouvert $U_{i}$, pour tout \'el\'ement~$i$ de~$\cn{1}{n}$. Nous utiliserons l'isomorphisme
$$\psi_{i} : \Gs_{i} \xrightarrow[\sim]{(\varphi_{i},\ldots,\varphi_{i})} \Os^n \xrightarrow[\sim]{\sigma_{i}^{-1}} \Os^n = \Gs_{0}.$$ 
Remarquons que le faisceau~$\Gs$ est coh\'erent, car nous l'avons construit en recollant des faisceaux coh\'erents (\emph{cf.} lemme \ref{FPcoherent}).

\begin{prop}\label{GinjAutGs}
Il existe un morphisme de groupes injectif du groupe~$G$ dans le groupe des automorphismes de $\Os$-alg\`ebres du faisceau~$\Gs$.
\end{prop}
\begin{proof}
Soit~$h$ un \'el\'ement du groupe~$G$. Notons~$\alpha_{h}$ la permutation de l'ensemble~$\cn{1}{n}$ telle que
$$\forall j\in\cn{1}{n},\, h g_{j} = g_{\alpha_{h}(j)}.$$
On d\'efinit \`a l'aide de cette permutation un automorphisme de $\Os_{U_{0}}$-alg\`ebres~$\mu_{h}$ de~$\Gs_{0}$ :
$$\mu_{h} : \Gs_{0} = \Os^n \xrightarrow[\sim]{\alpha_{h}} \Os^n=\Gs_{0}.$$

Soit~$i$ un \'el\'ement de~$\cn{1}{n}$. Notons~$\beta_{h}$ la permutation de l'ensemble~$\cn{0}{d_{i}-1}$ telle que
$$\forall u\in\cn{0}{d_{i}-1},\, hg_{a_{i,u}} = g_{a_{i,\beta_{h}(u)}} \textrm{ dans } G/\la g_{i} \ra.$$
Notons~$\beta'_{h}$ la permutation de l'ensemble~$\cn{1}{d_{i}}$ d\'efinie par
$$\forall u\in\cn{1}{d_{i}},\ \beta'_{h}(u) = \beta_{h}(u-1)+1.$$
Elle induit un automorphisme de $\Os_{U_{i}}$-alg\`ebres de~$\Gs_{i}$ :
$$\Gs_{i} = \Fs_{i}^{d_{i}} \xrightarrow[\sim]{\beta'_{h}} \Fs_{i}^{d_{i}} =\Gs_{i}.$$
Soit~$u$ un \'el\'ement de~$\cn{0}{d_{i}-1}$. Il existe un \'el\'ement~$m_{i,u}$ de~$\cn{0}{n_{i}-1}$ tel que
$$hg_{a_{i,u}} = g_{a_{i,\beta_{h}(u)}}\, g_{i}^{m_{i,u}} \textrm{ dans } G.$$
On d\'efinit alors un automorphisme~$\gamma_{h}$ de $\Os_{U_{i}}$-alg\`ebres de~$\Gs_{i}$ :
$$\Gs_{i} = \Fs_{i}^{d_{i}}  \xrightarrow[\sim]{(\varphi_{i}^{m_{i,1}},\ldots,\varphi_{i}^{m_{i},{d_{i}}})} \Fs_{i}^{d_{i}} = \Gs_{i}.$$
Un simple calcul montre qu'au-dessus de l'ouvert~$U_{i}$, les automorphismes~$\mu_{h}$ et~$\gamma_{h}\circ\beta'_{h}$ de~$\Gs$ co\"incident.

Nous avons donc construit une application
$$\mu : \begin{array}{ccc}
G & \to & \textrm{Aut}_{\Os}(\Gs)\\
h & \mapsto & \mu_{h}
\end{array}.$$
Montrons que c'est un morphisme de groupes. Soient~$h_{1}$ et~$h_{2}$ deux \'el\'ements de~$G$. Pour tout \'el\'ement~$j$ de~$\cn{1}{n}$, nous avons
$$\begin{array}{rcl}
g_{\alpha_{h_{1}h_{2}}(j)} &=& h_{1}h_{2}g_{j}\\
&=& h_{1}g_{\alpha_{h_{2}}(j)}\\
&=& g_{\alpha_{h_{1}}(\alpha_{h_{2}}(j))}.
\end{array}$$
Par cons\'equent, nous avons $\alpha_{h_{1}h_{2}} = \alpha_{h_{1}} \circ \alpha_{h_{2}}$ et donc $\mu_{h_{1}h_{2}} = \mu_{h_{1}} \circ \mu_{h_{2}}$. Par cons\'equent, l'application~$\mu$ est un morphisme de groupes.

Montrons finalement que le morphisme~$\mu$ est injectif. Soient~$h_{1}$ et~$h_{2}$ deux \'el\'ements de~$G$ tels que $\alpha_{h_{1}}=\alpha_{h_{2}}$. Nous avons alors
$$h_{1}g_{1} = g_{\alpha_{h_{1}}(1)} = g_{\alpha_{h_{2}}(1)} = h_{2}g_{1}.$$
On en d\'eduit que $h_{1}=h_{2}$. Par cons\'equent, le morphisme~$\mu$ est injectif.
\end{proof}

\begin{rem}
Il n'est gu\`ere difficile de montrer que le morphisme~$\mu$ construit pr\'ec\'edemment est en fait un isomorphisme de groupes. 
\end{rem}

\begin{cor}\label{GinjAutGsD}\index{Espace de Stein!couronne quelconque de A1@couronne quelconque de $\AA$}
Il existe un morphisme de groupes injectif du groupe~$G$ dans le groupe des automorphismes de $\Os(\D)$-alg\`ebres du faisceau~$\Gs(\D)$.
\end{cor}
\begin{proof}
Soient~$\As$ et~$\Bs$ deux faisceaux de $\Os_{\D}$-alg\`ebres coh\'erents. Consid\'erons l'application surjective
$$\textrm{Mor}_{\Os}(\As,\Bs) \to \textrm{Mor}_{\Os(\D)}(\As(\D),\Bs(\D)).$$
Elle est injective car les faisceaux~$\As$ et~$\Bs$ satisfont le th\'eor\`eme~A (\emph{cf.} corollaire \ref{ABouvert}).

On en d\'eduit que le morphisme de groupes injectif
$$\mu : G \to \textrm{Aut}_{\Os}(\Gs)$$
construit pr\'ec\'edemment induit un morphisme de groupes injectif
$$\mu_{\D} : G \to \textrm{Aut}_{\Os(\D)}(\Gs(\D)).$$
\end{proof}

\begin{lem}\label{racinedegreinfn}
Tout \'el\'ement de~$\Gs(\D)$ annule un polyn\^ome unitaire \`a coefficients dans~$\Ms(\D)$ de degr\'e inf\'erieur \`a~$n$.
\end{lem}
\begin{proof}
Soit~$s$ un \'el\'ement de~$\Gs(\D)$. Nous supposerons, tout d'abord, qu'il existe un point~$x_{0}$ de~$U_{0}$ tel que toutes les coordonn\'ees de son image~$s_{x_{0}}$ dans $\Gs_{x_{0}}=\Os_{X,x_{0}}^n$ soient distinctes. Puisque l'ouvert~$U_{0}$ est connexe, le principe du prolongement analytique (\emph{cf.} th\'eor\`eme \ref{prolan}) assure qu'en tout point~$x$ de~$U_{0}$, toutes les coordonn\'ees du germe~$s_{x}$ sont distinctes. Notons $a_{1},\ldots,a_{n}$ les coordonn\'ees de l'image de~$s$ dans $\Gs(U_{0})=\Os(U_{0})^n$. Posons
$$M(Z) = \prod_{j=1}^n (Z-a_{i}) \in \Os(U_{0})[Z].$$
En tout point~$x$ de~$U_{0}$, l'image du polyn\^ome~$M$ est l'unique polyn\^ome unitaire de degr\'e inf\'erieur \`a~$n$ \`a coefficients dans~$\Ms_{x}$ qui annule le germe~$s_{x}$. 

Pour tout \'el\'ement~$j$ de~$\cn{0}{n}$, posons $V_{j} = U_{0} \cup \bigcup_{1\le i\le j} \D''_{\m_{i}}$. Montrons, par r\'ecurrence, que pour tout \'el\'ement~$j$ de~$\cn{0}{n}$, il existe un polyn\^ome unitaire~$N_{j}$ de degr\'e~$n$ \`a coefficients dans~$\Ms(V_{j})$ qui annule l'\'el\'ement~$s_{|V_{j}}$ de~$\Gs(V_{j})$. Nous avons d\'ej\`a trait\'e le cas~$j=0$. Soit maintenant un \'el\'ement~$j$ de~$\cn{0}{n-1}$ pour lequel l'hypoth\`ese de r\'ecurrence est v\'erifi\'ee. Puisque l'ouvert~$\D''_{\m_{j+1}}$ est connexe, l'anneau~$\Ms(\D''_{\m_{i}})$ est un corps, d'apr\`es le corollaire \ref{anneauintegre}, et tout \'el\'ement de l'anneau $\Os(\D''_{\m_{i}})/(S^{n_{j+1}}-\pi_{\m_{j+1}}^{n_{j+1}}-T)$ est annul\'e par un polyn\^ome unitaire de degr\'e inf\'erieur \`a~$n_{j+1}$ \`a coefficients dans le corps~$\Ms(\D''_{\m_{i}})$. On en d\'eduit que l'\'el\'ement~$s_{|U_{j+1}}$ de~$\Gs(\D''_{\m_{i}})$ est annul\'e par un polyn\^ome unitaire~$M_{j+1}$ de degr\'e inf\'erieur \`a~$n$ \`a coefficients dans le corps~$\Ms(\D''_{\m_{i}})$. Soit~$x$ un \'el\'ement de $\D''_{\m_{i}}\cap U_{0} = U_{j+1}$. Nous avons d\'emontr\'e qu'il existe un unique polyn\^ome unitaire de degr\'e inf\'erieur \`a~$n$ \`a coefficients dans~$\Ms_{x}$ qui annule le germe~$s_{x}$. On en d\'eduit que les images les images des polyn\^omes~$N_{j}$ et~$M_{j+1}$ dans~$\Ms_{x}[Z]$ co\"incident. Par cons\'equent, les images de ces polyn\^omes dans $\Ms(U_{j+1})[Z]$ co\"incident. On en d\'eduit que le polyn\^ome~$N_{j}$ se prolonge en un polyn\^ome unitaire~$N_{j+1}$ de degr\'e inf\'erieur \`a~$n$ \`a coefficients dans~$\Ms(V_{j+1})$ qui annule l'\'el\'ement~$s_{|V_{j+1}}$ de~$\Gs(V_{j+1})$. 

On d\'eduit finalement le r\'esultat attendu du cas~$j=n$.

\bigskip

Soit~$x_{0}$ un point de l'ouvert~$U_{0}$. La fibre du faisceau~$\Gs$ au point~$x_{0}$ est isomorphe \`a l'alg\`ebre $\Os_{X,x_{0}}^n$. D'apr\`es le th\'eor\`eme \ref{lemniscateStein}, le faisceau~$\Gs$ v\'erifie le th\'eor\`eme~A sur le disque~$\D$. On en d\'eduit qu'il existe un \'el\'ement~$s_{0}$ de~$\Gs(\D)$ dont toutes les coordonn\'ees de l'image dans la fibre $\Gs_{x_{0}}=\Os_{X,x_{0}}^n$ sont distinctes. 

Soit~$s$ un \'el\'ement de~$\Os(\D)$. Il existe un \'el\'ement~$\lambda$ de~$\Os(\D)$ tel que toutes les coordonn\'ees du germe de la section $s_{1} = s+\lambda s_{0}$ au point~$x_{0}$ soient distinctes. Le raisonnement qui pr\'ec\`ede montre qu'il existe deux polyn\^omes unitaires~$P_{0}$ et~$P_{1}$ de degr\'e inf\'erieur \`a~$n$ \`a coefficients dans~$\Ms(\D)$ qui annulent respectivement les sections~$s_{0}$ et~$s_{1}$. D'apr\`es le corollaire \ref{anneauintegre}, l'anneau~$\Ms(\D)$ est un corps. On en d\'eduit qu'il existe un polyn\^ome unitaire~$P$ de degr\'e inf\'erieur \`a~$n$ \`a coefficients dans~$\Ms(\D)$ qui annule la section~$s$.
\end{proof}

\begin{lem}\label{GsDintegre}
L'alg\`ebre~$\Gs(\D)$ est int\`egre.
\end{lem}
\begin{proof}
Remarquons, tout d'abord, que l'alg\`ebre~$\Gs(\D)$ n'est pas nulle. En effet, les \'el\'ements~$0$ et~$1$ sont distincts. Il nous reste \`a montrer qu'elle ne contient aucun diviseur de z\'ero.

Soient~$s$ et~$t$ deux \'el\'ements de~$\Gs(\D)$ dont le produit est nul. Au-dessus de l'ouvert~$U_{0}$, le faisceau~$\Gs$ n'est autre que le faisceau~$\Os^n$. D'apr\`es le lemme \ref{U0connexe} et le corollaire \ref{anneauintegre}, l'anneau~$\Os(U_{0})$ est int\`egre. Par cons\'equent, la premi\`ere coordonn\'ee de l'une des deux sections doit \^etre nulle. Supposons que ce soit celle de~$s$. Notons 
$$s = (s_{1},\ldots,s_{n}) \textrm{ dans }\Os(U_{0})^n.$$
Nous supposons donc que~$s_{1}=0$.

Soit~$i$ un \'el\'ement de~$\cn{1}{n}$. Il existe un \'el\'ement~$j$ de~$\cn{1}{n}$ tel que~$g_{j}=g_{i}\, g_{1}^{-1}$ dans le groupe~$G$. Notons
$$s = (t_{1},\ldots,t_{d_{j}}) \textrm{ dans } \Fs_{j}(\D''_{\m_{j}})^{d_{j}}.$$
Il existe des \'el\'ements~$u$ de~$\cn{0}{d_{j}-1}$ et~$v$ de~$\cn{1}{n_{j}}$ tels que
$$g_{1} = g_{a_{j,u}}\, g_{j}^{v-1}.$$
Par d\'efinition du morphisme~$\psi_{j}$, nous avons alors
$$\varphi_{j}(t_{u+1}) = 
\left(\begin{array}{c}
s_{\sigma_{j}(un_{j}+1)}\\
\vdots\\
s_{\sigma_{j}(un_{j}+n_{j})}
\end{array}\right) \textrm{ dans } \Os(U_{j})^{n_{j}}.$$
Par d\'efinition de~$\sigma_{j}$, nous avons $\sigma_{j}(un_{j}+v)=1$. Par cons\'equent, l'\'el\'ement~$s_{\sigma_{j}(un_{j}+v)}$ de~$\Os(U_{j})$ est nul. Le lemme \ref{integreunebranche} assure que l'\'el\'ement~$t_{u+1}$ de~$\Fs_{j}(\D''_{\m_{j}})$ est \'egalement nul. Nous avons choisi l'\'el\'ement~$j$ de fa\c{c}on \`a avoir l'\'egalit\'e $g_{1}=g_{i}$ dans $G/\la g_{j} \ra$. On en d\'eduit qu'il existe un \'el\'ement~$w$ de~$\cn{1}{n_{j}}$ tel que $s_{i}=s_{\sigma_{j}(un_{j}+w)}$ dans~$\Os(U_{j})$. Par cons\'equent, l'\'el\'ement~$s_{i}$ est nul dans~$\Os(U_{j})$ et donc dans~$\Os(U_{0})$, par le principe du prolongement analytique.

Nous avons montr\'e que l'\'el\'ement~$s_{|U_{0}}$ de~$\Gs_{0}(U_{0})$ est nul. En utilisant de fa\c{c}on r\'ep\'et\'ee le lemme \ref{integreunebranche}, on en d\'eduit que l'\'el\'ement~$s$ de~$\Gs(\D)$ est nul. Par cons\'equent, l'alg\`ebre~$\Gs(\D)$ est int\`egre.
\end{proof}

\begin{lem}\label{AalgfermeGsD}
L'anneau~$A$ est alg\'ebriquement ferm\'e dans l'anneau~$\Gs(\D)$.
\end{lem}
\begin{proof}
Soit~$P$ un polyn\^ome unitaire \`a coefficients dans~$A$ sans racines dans~$A$. Supposons, par l'absurde, qu'il existe une section~$s$ de~$\Gs(\D)$ qui est racine du polyn\^ome~$P$. Notons~$x$ le point~$0$ de la fibre centrale~$X_{0}$ de l'espace~$X$. C'est un point de l'ouvert~$U_{0}$. Notons~$a$ la premi\`ere coordonn\'ee de l'image du germe~$s_{x}$ par l'isomorphisme $\Gs_{x} \xrightarrow[]{\sim} \Os_{X,x}^n$. C'est un \'el\'ement de~$\Os_{X,x}$ qui v\'erifie l'\'egalit\'e $P(a)=0$. D'apr\`es le corollaire \ref{descriptioncentral13}, l'anneau local~$\Os_{X,x}$ se plonge dans l'anneau~$K[\![T]\!]$. On en d\'eduit que le polyn\^ome~$P$ poss\`ede une racine dans l'anneau~$K[\![T]\!]$ et donc dans le corps~$K$. Puisque l'anneau~$A$ est alg\'ebriquement ferm\'e dans le corps~$K$, cette racine doit appartenir \`a~$A$. Nous avons abouti \`a une contradiction. On en d\'eduit le r\'esultat annonc\'e. 
\end{proof}

Introduisons une d\'efinition correspondant \`a cette propri\'et\'e.

\begin{defi}\index{Extension reguliere@Extension r\'eguli\`ere}
Une extension~$L$ du corps~$\Ms(\D)$ est dite r\'eguli\`ere si le corps~$K$ est alg\'ebriquement ferm\'e dans~$L$.
\end{defi}

Regroupons, \`a pr\'esent, les r\'esultats obtenus.

\begin{prop}\index{Fonctions meromorphes@Fonctions m\'eromorphes!extension galoisienne de MD@extension galoisienne de $\Ms(\D)$}
L'extension de corps
$$\Ms(\D) \to \Frac(\Gs(\D))$$
est finie de degr\'e~$n$, r\'eguli\`ere et galoisienne de groupe de Galois~$G$. 
\end{prop}
\begin{proof}
L'extension $\Ms(\D) \to \Frac(\Gs(\D))$ est finie et de degr\'e inf\'erieur \`a~$n$ d'apr\`es le lemme \ref{racinedegreinfn}. Elle est r\'eguli\`ere d'apr\`es le lemme \ref{AalgfermeGsD}. On d\'eduit du corollaire \ref{GinjAutGsD} que le groupe~$G$ s'injecte dans le groupe des $\Ms(\D)$-automorphismes du corps Frac($\Gs(\D)$). Or le groupe~$G$ a pour cardinal~$n$. On en d\'eduit que l'extension $\Ms(\D) \to \Frac(\Gs(\D))$ est exactement de degr\'e~$n$, qu'elle est galoisienne et que son groupe de Galois est isomorphe au groupe~$G$.
\end{proof}

Puisque nous sommes partis d'un groupe fini~$G$ arbitraire, nous avons finalement d\'emontr\'e le r\'esultat suivant.

\begin{thm}\index{Theoreme@Th\'eor\`eme!d'Harbater}\index{Probleme de Galois inverse@Probl\`eme de Galois inverse!Theoreme d'Harbater@Th\'eor\`eme d'Harbater}
Tout groupe fini est le groupe de Galois d'une extension finie, galoisienne et r\'eguli\`ere  du corps~$\Ms(\D)$. 
\end{thm}

Pour finir, donnons une description explicite du corps~$\Ms(\D)$. Rappelons que l'anneau~$A$ est munie de la norme~$\|.\|$ d\'efinie de la fa\c{c}on suivante :
$$\forall f\in A,\, \|f\| = \max_{\sigma : K\hookrightarrow \C}(|\sigma(f)|_{\infty}).$$

\begin{prop}\newcounter{Aun}\setcounter{Aun}{\thepage}
Notons~$A_{1^-}[\![T]\!]$ le sous-anneau de~$A[\![T]\!]$ form\'e des s\'eries 
$$\sum_{k\ge 0} a_{k}\, T^k$$
qui v\'erifient la condition suivante :
$$\forall r<1,\ \lim_{k \to +\infty}\|a_{k}\|\, r^k = 0.$$

Le morphisme naturel $A[T] \to \Ms(\D)$ induit un isomorphisme
$$\Frac(A_{1^-}[\![T]\!]) \xrightarrow[]{\sim} \Ms(\D).$$
\end{prop}
\begin{proof}
D'apr\`es les th\'eor\`emes \ref{lemniscateStein} et \ref{thPoincare}, le morphisme naturel
$$\Frac(\Os(\D)) \to \Ms(\D)$$
est un isomorphisme.

On montre que le morphisme $A[T] \to \Os(\D)$ induit un isomorphisme
$$A_{1^-}[\![T]\!] \xrightarrow[]{\sim} \Os(\D)$$
en \'ecrivant
$$\Os(\D) = \varinjlim_{r<1} \Os(\overline{D}(r))$$
et en utilisant la description explicite de l'anneau~$\Os(\overline{D}(r))$ fournie par le th\'eor\`eme \ref{imagedisque}. On en d\'eduit le r\'esultat annonc\'e.
\end{proof}

\index{Probleme de Galois inverse@Probl\`eme de Galois inverse!recollement|)}

\index{Probleme de Galois inverse@Probl\`eme de Galois inverse|)}

%% file: notations.tex
\chapter*{Glossaire des notations}

{\noindent\bf Espaces analytiques}

\bigskip

\begin{supertabular}{ll}

$\As$ & anneau de Banach, p.\theAs\\

$\Ms(\As)$ & spectre analytique de l'anneau de Banach~$\As$, p.\theMAs\\

$\E{n}{\As}$ & espace affine analytique de dimension~$n$ sur l'anneau de Banach~$\As$, p.\theAnAs\\

$\p_{x}$ & noyau de la semi-norme associ\'ee au point~$x$, p.\thepx, \thepxbis\\

$\Hs(x)$ & corps r\'esiduel compl\'et\'e du point~$x$, p.\theHsx, \theHsxbis\\

$f(x)$ & valeur de la fonction~$f$ au point~$x$, p.\thefdex, \thefdexbis\\

$\alphab$ & point rationnel, p.\thealphab\\

$\Os$ & faisceau structural, p.\theOs\\

$\m_{x}$ & id\'eal maximal de l'anneau local~$\Os_{x}$, p.\themx\\

$\kappa(x)$ & corps r\'esiduel du point~$x$, p.\thekappax\\

\end{supertabular}

\bigskip\bigskip

{\noindent\bf Points de la droite affine analytique au-dessus d'un corps trivialement valu\'e}

\bigskip

\begin{supertabular}{ll}

$\eta_{1}$ & point de Gau{\ss}, associ\'e \`a la valeur absolue triviale, p.\theetauntriv\\

$\eta_{P,0}$ & (pour~$P$ irr\'eductible) unique point d\'efini par l'\'equation $P=0$, p.\theetaPzero\\

$\eta_{P,r}$ & (pour~$P$ irr\'eductible et $r\in\R_{+}^*\setminus\{1\}$) unique point d\'efini par l'\'equation $|P|=r$, p.\theetaPr\\

$\eta_{\alpha,r}$ & ($r\in\R_{+}^*\setminus\{1\}$) autre notation pour le point $\eta_{T-\alpha,r}$, p.\theetaalphar\\

$\eta_{r}$ & (pour $r\in\R_{+}^*\setminus\{1\}$) autre notation pour le point $\eta_{0,r}$, p.\theetar\\

$\alpha$ & point rationnel, p.\thealpha\\

\end{supertabular}

\bigskip\bigskip

{\noindent\bf Points de la droite affine analytique au-dessus d'un corps ultram\'etrique complet}

\bigskip

\begin{supertabular}{ll}

$\alpha$ & point rationnel, p.\thealpha, \thealphapac\\

$\eta_{\alpha,r}$ & point de Shilov du disque de centre~$\alpha$ et de rayon~$r$, p.\theetaalpharac, \theetaalpharpac\\

$\eta_{r}$ & autre notation pour le point $\eta_{0,r}$, p.\theetarac, \theetarpac\\

$\eta_{\alphab,\br}$ & point de type~$4$, p.\theetaalphabbr\\

$\eta_{P,0}$ & (pour~$P$ irr\'eductible) unique point d\'efini par l'\'equation $P=0$, p.\theetaPzeropac\\

\end{supertabular}

\bigskip\bigskip


{\noindent\bf Faisceaux et fonctions}

\bigskip

\begin{supertabular}{ll}

$\Ks$ & pr\'efaisceau des fractions rationnelles sans p\^oles, p.\theKs\\

$\Os$ & faisceau structural, p.\theOs\\

$\Bs(V)$ & (pour $V$ compact) compl\'et\'e de l'anneau $\Ks(V)$ pour la norme uniforme, p.\theBsV\\

$\Ps$ & faisceau des parties principales, p.\thePs\\

$o$ & distribution d'ordres, p.\theo\\

$p$ & distribution de parties principales, p.\thep\\

$\Z_{r^+}[\![T]\!]$ & s\'eries \`a coefficients entiers de rayon de convergence strictement sup\'erieur \`a~$r$, p.\theZr\\

$A_{1^-}[\![T]\!]$ & s\'eries \`a coefficients dans~$A$ de rayon de convergence sup\'erieur ou \'egal \`a~$1$, p.\theAun\\

\end{supertabular}

\bigskip\bigskip

{\noindent\bf Flot}

\bigskip

\begin{supertabular}{ll}

$I_{b}$, $I_{x}$ & intervalle de d\'efinition du flot, p.\theIb\\

$|.|_{x}^\eps$ & (pour $\eps\in I_{x}$) image d'un point par le flot, p.\thevaeps\\

$|.|_{x}^0$ & prolongement du flot en~$0$, p.\thevazero\\

$D_{Y}$ & graphe du flot dans~$Y$, p.\theDY\\

$I_{Y}(x)$ & intervalle de d\'efinition du flot dans~$Y$, p.\theIYx\\

$T_{Y}(x)$ & trajectoire du point~$x$ dans~$Y$, p.\theTYx\\

$D_{Y}(x)$ & graphe du flot de~$x$ dans~$Y$, p.\theDYx\\

\end{supertabular}

\bigskip\bigskip

{\noindent\bf Disques et couronnes}

\bigskip

\begin{supertabular}{ll}

$\br$, $\bs$, $\bt$ & polyrayons, p.\thebs\\

$\bk$ & $n$-uplet d'entiers, p.\thebk\\

$\bs^\bk$ & puissance d'un polyrayon, p.\thebsbk, \thebsbkbis\\

$\bT$ & multivariable, p.\thebT\\

$\bT^\bk$ & puissance d'une multivariable, p.\thebTbk\\

$\bs\le \bt$, $\bs<\bt$, $\bs\prec\bt$ & relations d'ordre entre polyrayons, p.\therelo\\ 

$\bmax(\bs^\bk,\bt^\bk)$ & maximum pour les polyrayons, p.\thebmax, \thebmaxbis\\

$\bmin(\bs^\bk,\bt^\bk)$ & minimum pour les polyrayons, p.\thebmin\\

$\overline{D} (\bt)$ & disque ferm\'e, p.\theDbt\\

$\overline{C} ( \bs,\bt )$ & couronne ferm\'ee, p.\theCbsbt, \theCbsbtbis\\

$\overline{D}_{V} (\bt)$ & disque ferm\'e au-dessus de~$V$, p.\thecrel\\

$\mathring{D}_{V} (\bt)$ & disque ouvert au-dessus de~$V$, p.\thecrel\\

$\overline{C}_{V} (\bs,\bt)$ & couronne ferm\'ee au-dessus de~$V$, p.\thecrel\\

$\mathring{C}_{V} (\bs,\bt)$ & couronne ferm\'ee  au-dessus de~$V$, p.\thecrel\\

$\D$ & disque ouvert de rayon~$1$ de~$\E{1}{A}$, p.\theD\\

$\D_{a_{\sigma}}$ & trace de~$\D$ sur~$X_{a_{\sigma}}$, p.\theD\\

$\D_{\m}$ & trace de~$\D$ sur~$X_{\m}$, p.\theDm\\

$\D'_{\m}$ & trace de~$\D$ sur~$X'_{\m}$, p.\theDm\\

$\D''_{\m}$ & trace de~$\D$ sur~$X''_{\m}$, p.\theDm\\

$\As\of{\la}{|\bT| \le \bt}{\ra}$ & alg\`ebre de s\'eries associ\'ee au disque $\overline{D} (\bt)$, p.\thealgd\\

$\As\of{\la}{\bs \le |\bT| \le \bt}{\ra}$ & alg\`ebre de s\'eries associ\'ee \`a la couronne $\overline{C} (\bs,\bt)$, p.\thealgc, \thealgcbis\\

$\As_{s}$ & compl\'et\'e de $\As\of{\la}{|T|\le s}{\ra}$ pour la norme uniforme sur son spectre, p.\theAss\\

$\Os(V)\of{\la}{|\bT| \le \bt}{\ra}^\dag$ & alg\`ebre de s\'eries surconvergentes associ\'ee au disque $\overline{D}_{V}(\bt)$, p.\thesurcvd\\

$\Os(V)\of{\la}{\bs \le |\bT| \le \bt}{\ra}^\dag$ & alg\`ebre de s\'eries surconvergentes associ\'ee \`a la couronne $\overline{C}_{V}(\bs,\bt)$, p.\thesurcvc\\

$\|.\|_{\bt}$ & norme sur l'alg\`ebre $\As\of{\la}{|\bT| \le \bt}{\ra}$, p.\thenormed\\

$\|.\|_{V,\bt}$ & norme sur l'alg\`ebre $\Bs(V)\of{\la}{|\bT|\le \bt}{\ra}$, p.\thenVt\\

$\|.\|_{U,w,\textrm{r\'es}}$ & semi-norme induite par $\|.\|_{U,w}$\\ 
& sur l'anneau quotient $\Bs(U)[T]/(G(T))$, p.\thenres\\

$\|.\|_{\bs,\bt}$ & norme sur l'alg\`ebre $\As\of{\la}{\bs \le |\bT| \le \bt}{\ra}$, p.\thenormec, \thenormecbis\\

$\|.\|_{V,\bs,\bt}$ & norme sur l'alg\`ebre $\Bs(V)\of{\la}{\bs \le |\bT|\le \bt}{\ra}$, p.\thenVst\\

$L_{b}$ & anneau local limite d'alg\`ebres de disques, p.\theLb\\

$L_{b,\br}$ & anneau local limite d'alg\`ebres de couronnes, p.\theLbbr\\

\end{supertabular}

\bigskip\bigskip

{\noindent\bf Corps de nombres}

\bigskip

\begin{supertabular}{ll}

$K$ & corps de nombres, p.\theK\\

$A$ & anneau des entiers du corps~$K$, p.\theK\\

$|.|_{\infty}$ & valeur absolue usuelle, p.\thevainfinie\\

$|.|_{0}$ & valeur absolue triviale, p.\thevatriviale\\

$v_{p}$ & valuation $p$-adique, p.\thevapadique\\

$|.|_{p}$ & valeur absolue $p$-adique normalis\'ee par $|p|_{p} = \frac{1}{p}$, p.\thevapadique\\

$\m$ & id\'eal maximal de~$A$, p.\them\\

$A_{\m}$ & localis\'e de~$A$ en~$\m$, p.\them\\

$\pi_{\m}$ & uniformisante de $A_{\m}$, p.\them\\

$k_{\m}$ & corps r\'esiduel de~$A_{\m}$, p.\them\\

$\hat{A}_{\m}$ & compl\'et\'e de~$A_{\m}$ pour la topologie $\m$-adique, p.\them\\

$\hat{K}_{\m}$ & corps des fractions de $\hat{A}_{\m}$, p.\them\\

$p_{\m}$ & nombre premier v\'erifiant $\m \cap \Z = p_{\m}\, \Z$, p.\them\\

$n_{\m}$ & degr\'e de l'extension $\hat{K}_{\m}/Q_{p_{\m}}$, p.\them\\

$|.|_{\m}$ & valeur absolue sur~$K$ prolongeant $|.|_{p_{\m}}$, p.\them\\

$\sigma$ & plongement du corps~$K$ dans~$\C$ (et parfois id\'eal maximal de $A$), p.\thesigma\\

$\hat{K}_{\sigma}$ & $\R$ ou~$\C$ selon que le plongement est r\'eel ou complexe, p.\thesigma\\

$|.|_{\sigma}$ & valeur absolue associ\'ee au plongement, p.\thesigma\\

$n_{\sigma}$ & degr\'e de l'extension $\hat{K}_{\sigma}/\R$, p.\thesigma\\

$r_{1}$ & nombre de plongements r\'eels de~$K$, p.\therun\\

$r_{2}$ & moiti\'e du nombre de plongements complexes non r\'eels de~$K$, p.\therun\\

$\Sigma_{f}$ & ensemble des id\'eaux maximaux de~$A$, p.\therun\\

$\Sigma_{\infty}$ & ensemble des plongements complexes de~$K$ \`a conjugaison pr\`es, p.\therun\\

$\Sigma$ & r\'eunion des deux ensembles pr\'ec\'edents, p.\therun\\

$l(\sigma)$ & $1$ si $\sigma\in\Sigma_{f}$, $+\infty$ si $\sigma\in\Sigma_{\infty}$, p.\therun\\

\end{supertabular}

\bigskip\bigskip

{\noindent\bf Spectre analytique d'un anneau d'entiers de corps de nombres}

\bigskip

\begin{supertabular}{ll}

$\Ms(A)$, $B$ & spectre analytique de l'anneau~$A$, p.\theMA\\

$\Ms(A)_{\sigma}$, $B_{\sigma}$ & branche $\sigma$-adique, p.\thebranche\\

$\Ms(A)'_{\sigma}$, $B'_{\sigma}$ & branche $\sigma$-adique ouverte, p.\thebranche\\

$\Ms(A)''_{\sigma}$, $B''_{\sigma}$ & branche $\sigma$-adique semi-ouverte, p.\thebranche\\

$a_{0}$ & point associ\'e \`a la valeur absolue triviale sur~$A$, p.\theazero\\

$a_{\m}$ & point associ\'e \`a la valeur absolue $|.|_{\m}$, p.\theame\\

$a_{\m}^\eps$ & point associ\'e \`a la valeur absolue $|.|_{\m}^\eps$, p.\theame\\

$a_{\m}^0$ & autre notation pour le point~$a_{0}$, p.\theame\\

$\tilde{a}_{\m}$ & point associ\'e \`a la valeur absolue triviale sur $k_{\m}$, p.\thetam\\

$a_{\m}^{+\infty}$ & autre notation pour le point~$\tilde{a}_{\m}$, p.\thetam\\

$a_{\sigma}$ & point associ\'e \`a la valeur absolue $|.|_{\sigma}$, p.\thesigma\\

$a_{\sigma}^\eps$ & point associ\'e \`a la valeur absolue $|.|_{\sigma}^\eps$, p.\thesigma\\

$a_{\sigma}^0$ & autre notation pour le point $a_{0}$, p.\thesigma\\

\end{supertabular}

\bigskip\bigskip

{\noindent\bf Espace affine analytique au-dessus d'un anneau d'entiers de corps de nombres}

\bigskip

\begin{supertabular}{ll}

$\E{n}{A}$ & espace affine analytique de dimension~$n$ au-dessus de~$A$, p.\theMA\\

$X$ & $\E{n}{A}$ au chapitre~\ref{chapitreespace}, p.\thepi\\
& $\E{1}{A}$ aux chapitres~\ref{chapitredroite}, \ref{chapitreStein}, \ref{chapitreapplications} et au num\'ero \ref{adcdn}, p.\thenod, \thenodfini, \thenodStein, \thenodapp\\

$\pi$ & morphisme de projection naturel $X\to B$, p.\thepi\\

$X_{V}$ & image r\'eciproque de la partie~$V$ par le morphisme~$\pi$, p.\thepi\\

$X_{\sigma}$ & partie $\sigma$-adique, p.\thepartie\\

$X'_{\sigma}$ & partie $\sigma$-adique ouverte, p.\thepartie\\

$X''_{\sigma}$ & partie $\sigma$-adique semi-ouverte, p.\thepartie\\

$X_{b}$ & fibre du morphisme~$\pi$ au-dessus du point~$b$, p.\thepi\\

$X_{0}$ & fibre du morphisme~$\pi$ au-dessus du point~$a_{0}$, p.\thepartie\\

$\tilde{X}_{\m}$ & fibre du morphisme~$\pi$ au-dessus du point~$\tilde{a}_{\m}$, p.\thepartie\\

\end{supertabular}

\bigskip\bigskip

{\noindent\bf Espaces de Stein}

\bigskip

\begin{supertabular}{ll}

$K^-$, $K^+$ & compacts, p.\theKLM\\ 

$L$ & intersection de~$K^-$ et~$K^+$, p.\theKLM\\

$M$ & r\'eunion de~$K^-$ et~$K^+$, p.\theKLM\\

$\Omega$ & syst\`eme de Banach associ\'e au couple $(K^-,K^+)$, p.\theOmega\\

$K_{0}^-$, $K_{0}^+$ & compacts de~$B$, p.\theKLMzero\\

$L_{0}$ & intersection de~$K_{0}^-$ et~$K_{0}^+$, p.\theKLMzero\\

\end{supertabular}

%% file: memoire.bbl
\begin{thebibliography}{10}

\bibitem{rouge}
Vladimir~G. Berkovich.
\newblock {\em Spectral theory and analytic geometry over non-{A}rchimedean
  fields}, volume~33 of {\em Mathematical Surveys and Monographs}.
\newblock American Mathematical Society, Providence, RI, 1990.

\bibitem{bleu}
Vladimir~G. Berkovich.
\newblock \'{E}tale cohomology for non-{A}rchimedean analytic spaces.
\newblock {\em Inst. Hautes \'Etudes Sci. Publ. Math.}, (78):5--161 (1994),
  1993.

\bibitem{BCR}
J.~Bochnak, M.~Coste, and M.-F. Roy.
\newblock {\em G\'eom\'etrie alg\'ebrique r\'eelle}, volume~12 of {\em
  Ergebnisse der Mathematik und ihrer Grenzgebiete (3) [Results in Mathematics
  and Related Areas (3)]}.
\newblock Springer-Verlag, Berlin, 1987.

\bibitem{BGR}
Siegfried Bosch, Ulrich G{\"u}ntzer, and Reinhold Remmert.
\newblock {\em Non-{A}rchimedean analysis}, volume 261 of {\em Grundlehren der
  Mathematischen Wissenschaften}.
\newblock Springer-Verlag, Berlin, 1984.
\newblock A systematic approach to rigid analytic geometry.

\bibitem{algcom56}
N.~Bourbaki.
\newblock {\em \'{E}l\'ements de math\'ematique. {F}asc. {XXX}. {A}lg\`ebre
  commutative. {C}hapitre 5: {E}ntiers. {C}hapitre 6: {V}aluations}.
\newblock Actualit\'es Scientifiques et Industrielles, No. 1308. Hermann,
  Paris, 1964.

\bibitem{diviseurs}
N.~Bourbaki.
\newblock {\em \'{E}l\'ements de math\'ematique. {F}asc. {XXXI}. {A}lg\`ebre
  commutative. {C}hapitre 7: {D}iviseurs}.
\newblock Actualit\'es Scientifiques et Industrielles, No. 1314. Hermann,
  Paris, 1965.

\bibitem{SemCartan419}
Henri Cartan.
\newblock Faisceaux analytiques sur les vari{\'e}t{\'e}s de {S}tein :
  d{\'e}monstration des th{\'e}or{\`e}mes fondamentaux.
\newblock In {\em S{\'e}minaire Henri Cartan}, volume 4, expos{\'e} \no 19,
  pages 1--15, 1951--1952.

\bibitem{AntoineBourbaki}
Antoine Ducros.
\newblock Espaces analytiques $p$-adiques au sens de {B}erkovich.
\newblock {\em S\'em. Bourbaki}, (Exp. 958), 2006.

\bibitem{AntoineGazette}
Antoine Ducros.
\newblock G\'eom\'etrie analytique {$p$}-adique : la th\'eorie de {B}erkovich.
\newblock {\em Gaz. Math.}, (111):12--27, 2007.

\bibitem{Shilov}
Alain Escassut and Nicolas Ma\"iNetti.
\newblock Shilov boundary for normed algebras.
\newblock In {\em Topics in analysis and its applications}, volume 147 of {\em
  NATO Sci. Ser. II Math. Phys. Chem.}, pages 1--10. Kluwer Acad. Publ.,
  Dordrecht, 2004.

\bibitem{Frisch}
Jacques Frisch.
\newblock Points de platitude d'un morphisme d'espaces analytiques complexes.
\newblock {\em Invent. Math.}, 4:118--138, 1967.

\bibitem{TF}
Roger Godement.
\newblock {\em Topologie alg\'ebrique et th\'eorie des faisceaux}.
\newblock Hermann, Paris, 1973.
\newblock Troisi\`eme \'edition revue et corrig\'ee, Publications de l'Institut
  de Math\'ematique de l'Universit\'e de Strasbourg, XIII, Actualit\'es
  Scientifiques et Industrielles, No. 1252.

\bibitem{GR}
Hans Grauert and Reinhold Remmert.
\newblock {\em Theory of {S}tein spaces}.
\newblock Classics in Mathematics. Springer-Verlag, Berlin, 2004.
\newblock Translated from the German by Alan Huckleberry, Reprint of the 1979
  translation.

\bibitem{GK}
Elmar Grosse-Kl{\"o}nne.
\newblock Rigid analytic spaces with overconvergent structure sheaf.
\newblock {\em J. Reine Angew. Math.}, 519:73--95, 2000.

\bibitem{EGAISpringer}
Alexandre Grothendieck and Jean Dieudonn\'e.
\newblock {\em \'{E}l\'ements de g\'eom\'etrie alg\'ebrique. {I}. {L}e langage
  des sch\'emas.}, volume~8 of {\em Grundlehren der Mathematischen
  Wissenschaften}.
\newblock Springer-Verlag, Berlin-Heidelberg-New York, 1971.

\bibitem{GuRo}
Robert~C. Gunning and Hugo Rossi.
\newblock {\em Analytic functions of several complex variables}.
\newblock Prentice-Hall Inc., Englewood Cliffs, N.J., 1965.

\bibitem{Harbater}
David Harbater.
\newblock Convergent arithmetic power series.
\newblock {\em Amer. J. Math.}, 106(4):801--846, 1984.

\bibitem{mockcovers}
David Harbater.
\newblock Mock covers and {G}alois extensions.
\newblock {\em J. Algebra}, 91(2):281--293, 1984.

\bibitem{galoiscovers}
David Harbater.
\newblock Galois covers of an arithmetic surface.
\newblock {\em Amer. J. Math.}, 110(5):849--885, 1988.

\bibitem{AB}
Reinhardt Kiehl.
\newblock Theorem {A} und {T}heorem {B} in der nichtarchimedischen
  {F}unktionentheorie.
\newblock {\em Invent. Math.}, 2:256--273, 1967.

\bibitem{Lazard}
Michel Lazard.
\newblock Les z\'eros des fonctions analytiques d'une variable sur un corps
  valu\'e complet.
\newblock {\em Inst. Hautes \'Etudes Sci. Publ. Math.}, (14):47--75, 1962.

\bibitem{LiuGalois}
Qing Liu.
\newblock Tout groupe fini est un groupe de {G}alois sur {${\bf Q}\sb p(T)$},
  d'apr\`es {H}arbater.
\newblock In {\em Recent developments in the inverse Galois problem (Seattle,
  WA, 1993)}, volume 186 of {\em Contemp. Math.}, pages 261--265. Amer. Math.
  Soc., Providence, RI, 1995.

\bibitem{dimension}
Jun-iti Nagata.
\newblock {\em Modern dimension theory}.
\newblock Bibliotheca Mathematica, Vol. VI. Edited with the cooperation of the
  ``Mathematisch Centrum'' and the ``Wiskundig Genootschap'' at Amsterdam.
  Interscience Publishers John Wiley \& Sons, Inc., New York, 1965.

\bibitem{alh}
Michel Raynaud.
\newblock {\em Anneaux locaux hens\'eliens}.
\newblock Lecture Notes in Mathematics, Vol. 169. Springer-Verlag, Berlin,
  1970.

\bibitem{CL}
Jean-Pierre Serre.
\newblock {\em Corps locaux}.
\newblock Hermann, Paris, 1968.
\newblock Deuxi\`eme \'edition, Publications de l'Universit\'e de Nancago, No.
  VIII.

\bibitem{Siu}
Yum-Tong Siu.
\newblock Noetherianness of rings of holomorphic functions on {S}tein compact
  series.
\newblock {\em Proc. Amer. Math. Soc.}, 21:483--489, 1969.

\bibitem{Tate}
John Tate.
\newblock Rigid analytic spaces.
\newblock {\em Invent. Math.}, 12:257--289, 1971.

\end{thebibliography}
